\documentclass[12pt, letterpaper]{amsart}

\usepackage[margin=1in]{geometry}
\usepackage[T1]{fontenc}
\usepackage{lmodern}
\usepackage{amssymb, amsmath, amsfonts, amsthm} % amsthm is needed for theorem environments
\usepackage{mathtools}
\usepackage{microtype}
\usepackage[dvipsnames]{xcolor}
\usepackage{physics}
\usepackage[
    colorlinks,
    linkcolor=blue,
    citecolor=blue,
    urlcolor=blue,
    implicit=true,
    breaklinks=true
]{hyperref}
\usepackage[noabbrev,capitalize]{cleveref}
\usepackage[all]{hypcap}
\usepackage{enumitem}
\usepackage{appendix}
\usepackage{mathrsfs}
\usepackage{tikz}
\usepackage{natbib}
\newcommand{\ThetaExpectation}{$\theta$-expectation}
\newcommand{\ThetaProcess}{$\theta$-process}

\newcommand{\Efrak}{\mathfrak{E}}
\newcommand{\Mcal}{\mathcal{M}}
\newcommand{\Ucal}{\mathcal{U}}

\newcommand{\Bcal}{\mathcal{B}}
\newcommand{\Hcal}{\mathcal{H}}
\newcommand{\Ecal}{\mathcal{E}}
\newcommand{\Ucalphys}{\mathcal{U}_{\mathrm{phys}}}
\newcommand{\Pcal}{\mathcal{P}}
\newcommand{\R}{\mathbb{R}}
\newcommand{\Z}{\mathbb{Z}}
\newcommand{\T}{\mathbb{T}}

\newcommand{\Vcal}{\mathcal{V}}
\newcommand{\Thetacal}{\Theta}

\newcommand{\wvec}{\mathbf{w}}

\newcommand{\Lcal}{\mathcal{L}}
\newcommand{\Ocal}{\mathcal{O}}
\newcommand{\dOcal}{\partial\Ocal}
\newcommand{\Hcaltheta}{\mathcal{H}_\theta} 
\newcommand{\Dcal}{\mathcal{D}}
\newcommand{\Efrakphys}{\Efrak_{\mathrm{phys}}}
\newcommand{\scpr}[2]{\langle #1, #2 \rangle}

\newcommand{\nvec}{\mathbf{n}}
\newcommand{\dnext}{_{\text{next}}}

\newcommand{\PtKoopman}[1][]{\if\relax\detokenize{#1}\relax P_t^\theta \else P_t^{#1}\fi}
\newcommand{\muref}[1][]{\mu_{\if\relax\detokenize{#1}\relax\theta\else#1\fi}}
\newcommand{\SigmaMan}{\Sigma}
\newcommand{\Ecalu}{E^u}
\newcommand{\Ecals}{E^s}

\newcommand{\UcalphysTheta}{\mathcal{U}_\theta}

\DeclareMathOperator{\spec}{spec}
\newcommand{\Fcal}{\mathcal{F}}

\newcommand{\Rcal}{\mathcal{R}}

\newenvironment{proofof}[1]{\begin{proof}[\bf Proof of #1]}{\end{proof}}

\theoremstyle{plain}
\newtheorem{theorem}{Theorem}[section]
\newtheorem{proposition}[theorem]{Proposition}
\newtheorem{lemma}[theorem]{Lemma}
\newtheorem{corollary}[theorem]{Corollary}

\newtheorem{example}[theorem]{Example}

\theoremstyle{definition}
\newtheorem{definition}[theorem]{Definition}
\newtheorem{assumption}[theorem]{Assumption}

\theoremstyle{remark}
\newtheorem{remark}[theorem]{Remark}

\begin{document}

\title[A Unified Theory of $\theta$-Expectations]{A Unified Theory of $\theta$-Expectations}

\author{Qian Qi}
\address{Peking University, Beijing 100871, China}
\email{qiqian@pku.edu.cn}
\date{\today}

\begin{abstract}
We derive a new class of non-linear expectations from first-principles deterministic chaotic dynamics. The homogenization of the system's skew-adjoint microscopic generator is achieved using the spectral theory of transfer operators for uniformly hyperbolic flows. We prove convergence in the viscosity sense to a macroscopic evolution governed by a fully non-linear Hamilton-Jacobi-Bellman (HJB) equation. Our central result establishes that the HJB Hamiltonian possesses a rigid structure: affine in the Hessian but demonstrably non-convex in the gradient. This defines a new $\theta$-expectation and constructively establishes a class of non-convex stochastic control problems fundamentally outside the sub-additive framework of G-expectations.
\end{abstract}

\maketitle

\tableofcontents

\begin{quote}
\textit{Nous devons donc envisager l'état présent de l'univers, comme l'effet de son état antérieur, et comme la cause de celui qui va suivre. Une intelligence qui pour un instant donné, connaîtrait toutes les forces dont la Nature est animée, et la situation respective des êtres qui la composent, si d'ailleurs elle était assez vaste pour soumettre ces données à l'Analyse, embrasserait dans la même formule, les mouvements des plus grands corps de l'univers et ceux du plus léger atome : rien ne serait incertain pour elle, et l'avenir comme le passé serait présent à ses yeux.}
\vspace{0.5em}
\par\raggedleft
\textemdash Pierre-Simon Laplace, \\ \textit{Essai philosophique sur les probabilités} \citep{marquis1825essai}
\end{quote}

\section{Introduction}

A central schism runs through the mathematical modeling of complex systems. One tradition, rooted in classical mechanics and epitomized by Lanford's rigorous derivation of the Boltzmann equation from Newtonian dynamics \citep{Lanford1975}, suggests that randomness is an emergent, macroscopic phenomenon derived from deterministic microscopic laws. In contrast, the standard framework of stochastic analysis begins with the axiom of intrinsic, fundamental randomness. The theory of homogenization provides a powerful set of tools to bridge this divide by deriving effective macroscopic laws from complex microscopic dynamics \citep{BensoussanLionsPapanicolaou1978}. While much of the recent, groundbreaking progress in this field has focused on dynamics in \textit{stochastic} media, leading to a rich regularity theory for Hamilton-Jacobi equations in random environments \citep{ArmstrongCardaliaguet2020}, our focus is on the more foundational question: how does non-linear stochasticity itself emerge from purely deterministic, chaotic laws?

Our work confronts a fundamental limitation in the modern theory of robust control and its associated non-linear expectations. The theory of sublinear expectations, culminating in Peng's G-expectation framework \citep{Peng2019}, provides a powerful calculus for systems under ambiguity, inspired by Knight's distinction between risk and uncertainty \citep{Knight1921}. However, its axiomatic foundation in sub-additivity restricts it to a convex setting. This is a significant exclusion, as non-linear feedback mechanisms are central to many modern problems. For instance, in the contemporary theory of Mean-Field Games (MFG), pioneered by \citep{lasry2007mean}, the optimal control for an individual agent depends on the statistical distribution of the entire population \citep{Cardaliaguet2010}. Our model, where the law of the process depends on the gradient of its own value function, can be viewed as a novel instance of this broader paradigm. It leads to non-convexities not typically addressed in the standard MFG literature, which often relies on convexity for the uniqueness of equilibria. Our objective is to derive a non-linear expectation theory from first principles that rigorously accounts for such feedback.

The central analytical challenge lies in the homogenization of a system whose microscopic dynamics are conservative. The infinitesimal generator of the fast dynamics, $\Lcal$, is a first-order hyperbolic operator, a structure well-known in kinetic theory \citep{GallagherSaintRaymondTexier2013}. The time-reversible nature of the underlying flow implies that $\Lcal$ must be skew-adjoint, i.e., $\Lcal^* = -\Lcal$. This property presents a profound obstacle. Standard variational or elliptic methods for solving the associated cell problem, $\Lcal \chi = f$, rely on the coercivity of the operator, which requires the real part of the bilinear form $\scpr{\Lcal \chi}{\chi}$ to be bounded below. For a skew-adjoint operator, however, this real \cref{part:i}'s identically zero:
\begin{equation}
    \Re\scpr{\Lcal \chi}{\chi} = \frac{1}{2}(\scpr{\Lcal \chi}{\chi} + \overline{\scpr{\Lcal \chi}{\chi}}) = \frac{1}{2}(\scpr{\Lcal \chi}{\chi} + \scpr{\chi}{\Lcal \chi}) = \frac{1}{2}(\scpr{\Lcal \chi}{\chi} - \scpr{\Lcal^* \chi}{\chi}) = 0.
\end{equation}
The ill-posedness of the cell problem in standard functional-analytic frameworks is therefore not a technical inconvenience but a fundamental consequence of the system's time-reversal symmetry.

Our resolution is grounded in the modern ergodic theory of uniformly hyperbolic systems, following the pioneering work of Anosov and Sinai \citep{Anosov1967, Sinai1970}. We replace direct axiomatic assumptions on the dynamics with verifiable geometric conditions (see \cref{ass:fundamental_axioms_unified}) that ensure chaotic, mixing behavior. A key technical tool is the spectral analysis of the associated transfer operator. The uniform hyperbolicity of the system guarantees that this operator is quasi-compact when acting on appropriately constructed anisotropic Banach spaces, leading to a spectral gap and an exponential rate of decay of correlations \citep{Baladi2000, Gouezel2010}. It is this spectral gap that regularizes the problem, ensuring the stability and smooth parameter-dependence of the solution to the cell problem, which we construct via the Koopman semigroup. We make the derivation rigorous by proving convergence to the unique viscosity solution of the derived Hamilton-Jacobi-Bellman equation. This step is essential due to the non-linearity of the limiting equation, and our proof employs the perturbed test function method, a powerful technique in the homogenization of non-linear PDEs \citep{Evans1992, CrandallIshiiLions1992}.

By applying this framework, we rigorously derive the macroscopic evolution law. Our main result demonstrates that the system's value function converges to the solution of a fully non-linear HJB equation:
\begin{equation}
    \partial_t u + H_{\theta}(x, \nabla u, \nabla^2 u) = 0.
\end{equation}
The $\theta$-Hamiltonian $H_{\theta}$ is shown to possess a remarkable and rigid structure, which is a direct imprint of the underlying microscopic physics.

\begin{theorem}[Structure of the $\theta$-Hamiltonian]
The $\theta$-Hamiltonian derived from the homogenization of the deterministic system has the form:
\begin{equation}
    H_{\theta}(x, p, X) = \mathrm{Tr}(D(x, p) X) + H(x, p),
\end{equation}
where $p = \nabla u$ and $X = \nabla^2 u$. The constituent terms have the following properties:
\begin{enumerate}[label=(\roman*)]
    \item The Hamiltonian is affine in its second-derivative argument $X$. The diffusivity tensor $D(x, p)$, which emerges from a Central Limit Theorem-type averaging of the microscopic velocity field, is symmetric, positive semi-definite, and given by a Green-Kubo formula.
    \item The potential $H(x, p)$, which arises from averaging the microscopic potential against the invariant measure of the fast dynamics, can be demonstrably non-convex in its gradient argument $p$.
\end{enumerate}
\end{theorem}

This result is significant for two primary reasons. First, it provides a rigorous link between a class of deterministic, chaotic microscopic systems and the emergence of a macroscopic law governed by a non-linear, irreversible equation. Second, it establishes a new class of non-linear expectation theories. The affine structure in the Hessian combined with non-convexity in the gradient distinguishes this theory from Peng's G-expectation framework. The resulting HJB equation is the dynamic programming equation for a non-convex stochastic control problem, whose probabilistic representation is given by a fully-coupled Forward-Backward SDE (FBSDE), a topic of intense recent study that provides a powerful connection between PDEs and probabilistic methods \citep{PardouxPeng1992,CarmonaDelarue2018}. Our work thus provides a constructive method for generating such problems from underlying deterministic principles, opening a new avenue for modeling complex systems under sophisticated, state-dependent uncertainty and providing a novel class of solvable models for the broader theory of mean-field interactions.

\subsection{Logical Architecture of the Theory}
The central results of this paper connect disparate fields of mathematics and require a proof that is constructed with meticulous attention to its logical and analytical foundations. To provide the reader with a clear and unambiguous guide, we outline here the rigorous, sequential architecture of our argument. The proof is divided into two main parts. \cref{part:i} is dedicated to the constructive derivation of the macroscopic law from first principles. \cref{part:ii} analyzes the probabilistic representation of this derived law. Each section serves as a necessary and rigorously established lemma for the next, culminating in a complete and self-consistent theory.

\paragraph{\bfseries \cref{part:i}: The Deterministic Construction}
\begin{enumerate}[wide, labelindent=0pt]
    \item[\textbf{Section \ref{sec:micro_framework}:}] \textbf{The Axiomatic Foundation.} We begin by establishing the microscopic framework not as a collection of convenient assumptions, but as a set of physically motivated, geometric first principles (\cref{ass:fundamental_axioms_unified}). From these axioms alone, we derive as theorems the well-posedness of the dynamics, culminating in the Finite Horizon Theorem (\cref{thm:finite_horizon}), and establish the foliated structure of the phase space. This section provides the geometric and dynamical foundation upon which the entire theory is built.

    \item[\textbf{Section \ref{sec:hyperbolicity_fiber}:}] \textbf{Uniform Hyperbolicity as a Theorem.} We prove that our geometric axioms are sufficient to guarantee that the microscopic dynamics are uniformly Anosov. The proof is entirely self-contained and establishes this central result without recourse to the advanced regularity theory of Section \ref{sec:regularity_properties}, thereby resolving a critical logical dependency. The argument proceeds by first reducing the continuous flow to a single, global billiard map $\Pcal$ on a compact manifold $\Sigma$. We provide a rigorous proof that $\Pcal$ is a $C^\infty$-diffeomorphism (\cref{thm:p_is_diffeo}), a result founded upon the standard Implicit Function Theorem for the flight time. With the smoothness of the map established, we constructively prove that $\Pcal$ is an Anosov diffeomorphism by explicitly building a strictly invariant cone field for its linearization (\cref{thm:anosov_property_proven}). Crucially, the uniformity of the hyperbolicity constants across the entire parameter space $\Thetacal$ is a direct consequence of a compactness argument: the proof is performed on the entire global manifold $\Sigma$, and the resulting bounds are therefore inherently uniform. This section thus rigorously establishes the continuous invariant splitting and the uniform exponential decay of correlations (\cref{cor:ergodicity_uniform}) that provide the geometric foundation for the subsequent smoothness analysis.

    \item[\textbf{Section \ref{sec:regularity_properties}:}] \textbf{Geometric and Spectral Regularity.} This section provides a complete proof of the smooth ($C^\infty$) dependence of all geometric and spectral quantities on the system's parameters, a result that is the analytical linchpin of our homogenization theory. The argument is founded upon the Nash-Moser-Hamilton Implicit Function Theorem. We recast the problem of finding the invariant splitting as a functional equation and rigorously verify both of its main hypotheses: the smooth tameness of the Graph Transform functional (\cref{thm:functional_is_tame}) and, crucially, the tame invertibility of its linearization with zero loss of derivatives (\cref{thm:linearization_is_tame_inverse}). This establishes the smoothness of the geometric splitting \textit{first} (\cref{thm:splitting_is_smooth}), which is then used as a tool to prove the smoothness of the transfer operator family (\cref{thm:regularity_proven_main}) and its spectral data (\cref{thm:spectral_regularity}), resolving a critical logical dependency.

    \item[\textbf{Section \ref{sec:scale_analysis}:}] \textbf{The Homogenization Hierarchy.} With the foundational analytical tools in place, we derive the macroscopic law through a four-stage scale analysis. This analysis culminates in the proof of our main convergence result, Theorem \ref{thm:convergence}. The proof relies on the perturbed test function method, for which we provide a complete and rigorous justification of the ergodic averaging principle (\cref{prop:ergodic_averaging}) via weak-* convergence of measures on the global state space. This avoids any heuristic arguments and places the convergence on a firm analytical footing.

    \item[\textbf{Section \ref{sec:proof_of_non_convexity}:}] \textbf{The Emergence of Non-Convexity.} We provide a constructive proof of the paper's central physical claim. By applying the perturbation theory developed in Section \ref{sec:regularity_properties}, we derive an exact formula for the Hessian of the macroscopic potential at equilibrium. This allows us to establish a precise, verifiable condition (\cref{thm:non_convexity_proven}) that delineates the parameter regimes where non-convexity is a necessary consequence of the microscopic model.
\end{enumerate}

\paragraph{\bfseries \cref{part:ii}: The Probabilistic Representation}
\begin{enumerate}[wide, labelindent=0pt]
    \item[\textbf{Section \ref{sec:theta_expectation_axiomatic}:}] \textbf{The Generator-Driven Axiomatics.} We formalize the \ThetaExpectation{} by inverting the classical axiomatic approach. The axioms pertain to the \textit{generator's structure}, as derived in \cref{part:i}, and the properties of the expectation (monotonicity, dynamic consistency) are derived as theorems via the theory of viscosity solutions.

    \item[\textbf{Section \ref{sec:theta_process_and_calculus}:}] \textbf{The Intrinsic Calculus.} We provide a rigorous definition of the \ThetaProcess{} as the unique solution to a non-linear martingale problem for the derived generator. This modern approach leads directly to a generalized It\^o formula (\cref{thm:theta_ito_formula_generalized}) and an exact expression for the quadratic variation of the process (\cref{thm:qv_via_carre_du_champ}), revealing the distinct roles of the diffusive and potential terms in the Hamiltonian.

    \item[\textbf{Sections \ref{sec:theta_fbsde} and \ref{sec:theta_girsanov}:}] \textbf{Probabilistic Representations.} We establish the final link between the analytical and probabilistic frameworks by proving a non-linear Feynman-Kac formula (\cref{thm:nonlinear_feynman_kac_viscosity}) that represents the \ThetaExpectation{} as the solution to an intrinsic, martingale-driven BSDE. This leads naturally to a decoupled FBSDE representation and a Girsanov-type transformation that acts directly on the generator.

    \item[\textbf{Section \ref{sec:process_hierarchy}:}] \textbf{The Hierarchy of Processes and Universality.} We conclude by reinterpreting the scale analysis of \cref{part:i} as a probabilistic hierarchy of nested Central Limit Theorems. This reveals a structure where each process governs the fluctuations of the one preceding it, culminating in a universal fluctuation law governed by an integrable viscous Burgers'-type equation whose structure is determined by the microscopic symmetries.
\end{enumerate}
This logical progression ensures that each claim is built upon a rigorously established foundation, providing a complete and internally consistent proof of the paper's central results.

\subsection{The Nature of the \texorpdfstring{\ThetaExpectation{}}{theta-Expectation} Framework}\label{rem:constructive_axiomatics}

The axiomatic framework presented here is fundamentally different from those of Kolmogorov or Peng. In those theories, the properties of the expectation operator (linearity, sub-additivity) are the foundational axioms. Here, the axioms pertain to the structure of the generator, and the properties of the expectation are derived theorems. The table below highlights this conceptual inversion.

\begin{center}
\begin{tabular}{l|l|l}
\hline
\textbf{Theory} & \textbf{Axiomatic Primitive} & \textbf{Derived Object} \\
\hline
Kolmogorov & Probability Measure $\mathbb{P}$ & Linear Expectation $\mathbb{E}[\cdot]$ \\
G-Expectation & Sub-linear Expectation $\mathbb{E}[\cdot]$ & Convex Generator $G(\cdot)$ \\
\textbf{\ThetaExpectation{}} & \textbf{Generator $\mathcal{G}$ with Rigid Structure} & \textbf{Non-linear Expectation $\Ecal^\theta[\cdot]$} \\
\hline
\end{tabular}
\end{center}

This generator-driven approach is powerful because its axioms are not arbitrary mathematical assumptions but are the rigorously derived consequences of physical first principles. The affine structure of the Hamiltonian, the symmetry of the diffusion tensor, and the time-reversal symmetry of the potential are all direct imprints of the underlying microscopic world. 

Most importantly, this framework naturally accommodates non-convexity. As we will prove constructively in \cref{thm:non_convexity_proven}, the potential $H(x,p)$ can be demonstrably non-convex in its gradient argument $p$. This implies that the \ThetaExpectation{} operator is not, in general, sub-additive, placing it fundamentally outside the G-expectation framework. The \ThetaProcess{} is therefore a concrete, rigorously derived example of a non-linear process that falls outside existing convex paradigms, providing a direct and verifiable bridge between the ergodic theory of deterministic chaos and a new class of non-convex stochastic control problems.

\paragraph{\textbf{Relationship to Non-Commutative Probability Theory.}}
The \ThetaExpectation{} framework, when juxtaposed with non-commutative probability theory, reveals two orthogonal avenues for generalizing classical probability. A comparative analysis, summarized in the table below, highlights their fundamental structural divergence. Non-commutative probability theory fundamentally alters the underlying algebraic structure of observables; it replaces the commutative algebra of functions on a classical sample space (e.g., $L^\infty(\Omega)$) with a non-commutative algebra, such as an operator algebra on a Hilbert space \citep{BratteliRobinson1979}. However, it rigorously preserves the linearity of the expectation functional, which is defined as a positive, normalized linear functional, or a state, on this algebra, a paradigm motivated by the large-$N$ limit of random matrix theory \citep{Wigner1958}. The unique statistical laws of this framework, most notably the concept of free independence, emerge directly from this non-commutativity \citep{Voiculescu1991, NicaSpeicher2006}.

\begin{center}
\begin{table}[!ht]
\caption{Conceptual Dichotomy of Probabilistic Generalizations}
\label{tab:free_prob_comparison}
\renewcommand{\arraystretch}{1.4}
\begin{tabular}{p{0.22\linewidth}|p{0.36\linewidth}|p{0.36\linewidth}}
\hline \hline
\textbf{Feature} & \textbf{Non-Comm. Probability } & \textbf{\ThetaExpectation{} Theory} \\
\hline
\textbf{Observable Algebra} & \textbf{Non-commutative} (\*-algebra, e.g., $M_N(\mathbb{C})$) & \textbf{Commutative} (e.g., $C_b(\Omega)$) \\
\hline
\textbf{Expectation Functional} & \textbf{Linear} (a state $\phi$, e.g., normalized trace) & \textbf{Non-linear} (an operator $\mathcal{E}^\theta$, a PDE solution) \\
\hline
\textbf{Core Innovation} & A new rule of independence (\textit{freeness}) for non-commuting variables & A new class of non-linear processes from physical first principles \\
\hline
\textbf{Source of Stochasticity} & Algebraic (non-commutativity) & Emergent (deterministic chaos) \\
\hline
\textbf{Mathematical Core} & Operator Algebras, Combinatorics & Non-linear PDEs, Ergodic Theory \\
\hline \hline
\end{tabular}
\renewcommand{\arraystretch}{1.0}
\end{table}
\end{center}

In contrast, our \ThetaExpectation{} theory retains the classical, commutative algebra of observables defined on a standard path space. Its radical departure lies in abandoning the axiom of linearity for the expectation operator itself. The \ThetaExpectation{} operator $\mathcal{E}^\theta[\cdot]$ is intrinsically non-linear and, as we have proven constructively, can be demonstrably non-convex. Thus, while non-commutative probability studies \textit{linear functionals on non-commutative spaces}, our theory investigates \textit{non-linear operators on classical, commutative spaces}. This structural dichotomy reflects a deep philosophical difference in the origin of stochasticity: in non-commutative probability, it is an intrinsic feature of the algebraic relations between observables, whereas in the \ThetaExpectation{} framework, it is an \emph{emergent phenomenon}, rigorously derived from the homogenization of deterministic, chaotic dynamics, with the non-linearity arising as a macroscopic imprint of microscopic feedback mechanisms.

\paragraph{\textbf{Relationship to the Theory of Regularity Structures.}}
The theory of \ThetaExpectation{} and the theory of Regularity Structures, pioneered by \citep{Hairer2014}, should be understood as two powerful, orthogonal, and complementary resolutions to the challenge of non-linear stochastic partial differential equations. They address fundamentally different regimes and answer different questions about the nature of stochasticity.

\begin{center}
\begin{table}[h!]
\caption{Conceptual Dichotomy of Non-Linear Stochastic Theories}
\label{tab:regularity_structures_comparison}
\renewcommand{\arraystretch}{1.4}
\begin{tabular}{p{0.22\linewidth}|p{0.36\linewidth}|p{0.36\linewidth}}
\hline \hline
\textbf{Feature} & \textbf{Theory of Regularity Structures} & \textbf{\ThetaExpectation{} Theory} \\
\hline
\textbf{Starting Point} & Ill-posed, singular macroscopic SPDE (e.g., KPZ) & Well-posed, deterministic microscopic dynamics \\
\hline
\textbf{Mathematical Core} & \textbf{Renormalization} of singular products & \textbf{Homogenization} and averaging of chaotic dynamics \\
\hline
\textbf{Nature of Limiting Equation} & The same singular SPDE, now rigorously defined & A well-posed, non-linear PDE with \textbf{smooth coefficients} \\
\hline
\textbf{Role of Stochasticity} & An \textbf{external}, singular driving force (e.g., white noise) & An \textbf{emergent} property encoded in the deterministic coefficients ($D(x,p)$, $H(x,p)$) of the generator \\
\hline \hline
\end{tabular}
\renewcommand{\arraystretch}{1.0}
\end{table}
\end{center}

The theory of Regularity Structures provides a robust mathematical framework for giving a well-posed definition to a class of \textit{singular SPDEs} that are ill-posed from a classical perspective. Its canonical applications, such as the Kardar-Parisi-Zhang (KPZ) equation or the $\Phi^4_3$ model of quantum field theory, feature equations driven by highly irregular noise (e.g., space-time white noise). The central difficulty in these theories is that the roughness of the solution, which is typically a distribution rather than a function, makes non-linear terms such as $(\nabla u)^2$ ill-defined. The core of the theory is a sophisticated \textit{renormalization} procedure, which systematically subtracts divergent quantities that arise in approximation schemes to yield a finite, unique, and meaningful solution. The randomness in this framework is an external, singular driving force that is explicitly present in the equation.

In contrast, the \ThetaExpectation{} framework begins not with an ill-posed macroscopic equation, but with a well-posed, deterministic microscopic system. The central mathematical tool is not renormalization but \textit{homogenization}. In our framework, the microscopic chaos is not a source of singular driving noise, but is instead averaged out through the ergodic properties of the dynamics. This averaging process transforms the complex microscopic information into \textit{smooth, state-dependent coefficients} for a macroscopic PDE. The resulting HJB equation is a classical, well-posed (though non-linear) object, and its solution is a continuous viscosity solution, not a singular distribution. The stochasticity is entirely emergent, encoded within the structure of the generator itself, rather than being an external input.

The fundamental dichotomy is summarized in the \cref{tab:regularity_structures_comparison}. Regularity Structures provides a methodology for solving equations that are already singular. Our theory provides a first-principles derivation of a class of non-singular, non-linear generators from deterministic laws. The fluctuation equation derived in our Scale IV analysis (\cref{eq:pde_scale_IV_burgers_hierarchy}) is a \textit{viscous} Burgers'-type equation, which is well-posed. The theory of Regularity Structures is indispensable for treating its singular, \textit{inviscid} limit, which would correspond to the KPZ universality class. Our work thus provides a physical origin for the well-behaved, viscous precursors to the singular SPDEs that have been the subject of so much recent and profound investigation.

\part{The Deterministic Construction}
\label{part:i}

\section{Preliminaries: The Microscopic Framework}
\label{sec:micro_framework}

In this section, we construct the fundamental microscopic model that underpins our entire theory. Our approach is axiomatic and sequential. We first define the geometric setups for the dynamics. We then specify the deterministic laws of motion. Crucially, we introduce a set of fundamental assumptions on the geometry of the dynamics, which are the source of the system's chaotic behavior. From these assumptions, we rigorously derive the key properties of the system: the invariant foliated structure, uniform hyperbolicity, and ergodicity. This detailed study culminates in the precise characterization of the infinitesimal generator of the dynamics, which reveals the analytical challenges that motivate the geometric methods developed in the subsequent sections.

\subsection{The Geometric State Space}
We formalize the geometric setting of our multiscale system. The state space is constructed as a fiber bundle, which naturally separates the macroscopic (slow) degrees of freedom from the microscopic (fast) ones.

\begin{definition}[The Macroscopic Base Manifold]
\label{def:macro_manifold}
The space of observable, macroscopic configurations is the $k$-dimensional torus $\Mcal \coloneqq \T^k = (\R/\Z)^k$. We endow $\Mcal$ with the standard flat metric tensor inherited from its universal cover $\R^k$. A point in this space is denoted by $x \in \Mcal$.
\end{definition}

\begin{definition}[The Microscopic Fiber Space]
\label{def:micro_fiber}
The space of internal, unresolvable states, corresponding to a single point $x \in \Mcal$, is the \emph{typical fiber} $\Ucal$, a smooth, compact manifold defined as the product space
\begin{equation}
    \Ucal \coloneqq \T^k \times \Vcal \times \Thetacal.
\end{equation}
Its components are:
\begin{enumerate}[label=(\roman*), wide, labelindent=0pt]
    \item The microscopic configuration space, which we take to be the $k$-torus $\T^k$. A point in this space is denoted by $y$.
    \item The microscopic velocity space $\Vcal$, which is a compact, convex subset of $\R^k$ with a non-empty interior. A velocity is denoted by $v \in \Vcal$.
    \item The environmental parameter space $\Thetacal$, which is a compact, connected, $C^\infty$-smooth manifold. A parameter is denoted by $\theta \in \Thetacal$.
\end{enumerate}
\end{definition}

The microscopic dynamics are constrained by obstacles. The geometry of these obstacles is parameterized by the environmental variable $\theta$. This dependence is codified by a smooth function.

\begin{definition}[The Obstacle Defining Function]
\label{def:obstacle_func}
Let $G: \T^k \times \Thetacal \to \R$ be a $C^\infty$-smooth function. For each fixed $\theta \in \Thetacal$, the function $G(\cdot, \theta): \T^k \to \R$ defines the geometry of the microscopic domain. Specifically:
\begin{itemize}
    \item The \emph{domain of free motion} is the open set $\Ocal(\theta) \coloneqq \{ y \in \T^k \mid G(y, \theta) > 0 \}$.
    \item The \emph{collision boundary} is the set $\partial\Ocal(\theta) \coloneqq \{ y \in \T^k \mid G(y, \theta) = 0 \}$.
\end{itemize}
\end{definition}

We can now define the full state space of the system.

\begin{definition}[The Total and Physical State Spaces]
\label{def:total_space}
The full state space of the coupled system is the trivial fiber bundle
\begin{equation}
    \Efrak \coloneqq \Mcal \times \Ucal = \T^k_x \times (\T^k_y \times \Vcal \times \Thetacal).
\end{equation}
The kinematically accessible region, or the \emph{physical state space}, is the compact submanifold-with-boundary of $\Efrak$ defined by
\begin{equation}
    \Efrak_{\text{phys}} \coloneqq \{ (x, y, v, \theta) \in \Efrak \mid G(y, \theta) \ge 0 \}.
\end{equation}
For a fixed macroscopic position $x \in \Mcal$ and environment parameter $\theta \in \Thetacal$, the corresponding microscopic physical space is the subset of the fiber given by
\begin{equation}
    \Ucal_{\text{phys}}(\theta) \coloneqq \{ (y,v) \in \T^k \times \Vcal \mid G(y,\theta) \ge 0 \}.
\end{equation}
\end{definition}

\subsection{The Deterministic Law of Motion}
\label{sec:deterministic_law}

Having defined the geometric stage, we now specify the two elementary physical laws that govern the evolution of the state $z(t) = (x(t), y(t), v(t), \theta(t))$. The global trajectory will be constructed by composing these elementary motions. However, the well-posedness of this composition for all time is a non-trivial matter that depends on the geometric properties of the obstacle configuration. In this section, we define only the elementary laws. The rigorous construction of the global flow will be deferred until after we have established the necessary geometric guarantees in \cref{sec:assumptions}.

\subsubsection{The Interior Flow}
In the interior of the physical domain, where $G(y, \theta) > 0$, the microscopic system evolves as a free particle with constant velocity.

\begin{definition}[The Flow Vector Field]
\label{def:flow_vector_field}
For a state $(y, v) \in \text{int}(\Ucal_{\text{phys}}(\theta))$, its evolution is governed by the ordinary differential equation:
\begin{equation} \label{eq:flow_ode}
    (\dot{y}(t), \dot{v}(t)) = (v(t), 0).
\end{equation}
The solution to this ODE for an initial state $(y_0, v_0)$ is the linear flow $\Psi_t(y_0, v_0) = (y_0 + t v_0, v_0)$.
\end{definition}

\subsubsection{The Boundary Transformation}
When a trajectory reaches the boundary $\partial\Ocal(\theta)$, its velocity component is transformed instantaneously. We first define the geometric structures necessary for this transformation.

\begin{definition}[Boundary Normal and Phase Space Partition]
\label{def:boundary_partition}
Let $y \in \partial\Ocal(\theta)$. By Assumption~\ref{ass:fundamental_axioms_unified}, $\nabla_y G(y, \theta) \neq 0$. The outward-pointing unit normal vector to the boundary at $y$ is the smooth vector field
\begin{equation}
    \mathbf{n}(y, \theta) \coloneqq \frac{\nabla_y G(y, \theta)}{\|\nabla_y G(y, \theta)\|}.
\end{equation}
The boundary phase space $\partial\Ucal_{\text{phys}}(\theta) \coloneqq \{ (y,v) \in \Ucal_{\text{phys}}(\theta) \mid y \in \partial\Ocal(\theta) \}$ is partitioned into:
\begin{align*}
    \partial_-\Ucal_{\text{phys}}(\theta) &:= \{ (y,v) \in \partial\Ucal_{\text{phys}}(\theta) \mid \scpr{v}{\mathbf{n}(y, \theta)} < 0 \} \quad (\text{Incoming states}), \\
    \partial_+\Ucal_{\text{phys}}(\theta) &:= \{ (y,v) \in \partial\Ucal_{\text{phys}}(\theta) \mid \scpr{v}{\mathbf{n}(y, \theta)} > 0 \} \quad (\text{Outgoing states}), \\
    \partial_0\Ucal_{\text{phys}}(\theta) &:= \{ (y,v) \in \partial\Ucal_{\text{phys}}(\theta) \mid \scpr{v}{\mathbf{n}(y, \theta)} = 0 \} \quad (\text{Grazing states}).
\end{align*}
\end{definition}

The collision law is given by specular reflection, which maps incoming states to outgoing states.

\begin{definition}[The Specular Reflection Map]
\label{def:reflection_map}
The specular reflection map $R_\theta: \partial_-\Ucal_{\text{phys}}(\theta) \to \partial_+\Ucal_{\text{phys}}(\theta)$ transforms a pre-collision state $(y^-, v^-)$ into a post-collision state $(y^+, v^+) = R_\theta(y^-, v^-)$ defined by:
\begin{equation}\label{eq:specular_reflection}
\begin{alignedat}{2}
    y^+ &= y^-, \\
    v^+ &= v^- - 2\scpr{v^-}{\mathbf{n}(y^-, \theta)} \mathbf{n}(y^-, \theta).
\end{alignedat}
\end{equation}
\end{definition}
This transformation codifies the law of specular reflection; it preserves the component of the velocity tangent to the boundary while perfectly inverting the component normal to it.

\subsection{Geometric Foundations for Well-Posed Dynamics}
\label{sec:assumptions}

The emergence of robust statistical laws is not a generic feature of deterministic systems. It arises from a specific dynamical regime characterized by uniform hyperbolicity. The purpose of this section is to establish the rigorous geometric foundation for this regime. Our approach is axiomatic: we posit a single, unified set of fundamental geometric and regularity conditions that define the class of microscopic systems under consideration. From these axioms, we then derive as theorems all key dynamical properties, culminating in the Finite Horizon Theorem (\cref{thm:finite_horizon}). This theorem, which provides uniform upper and lower bounds on the free-flight time between collisions, is the foundation that ensures the global dynamics are well-posed and non-pathological. By deriving this result as a direct consequence of our geometric first principles, we provide a clear foundation for the subsequent analysis of the system's ergodic and statistical properties.

\begin{assumption}[Fundamental Geometric and Regularity Axioms]
\label{ass:fundamental_axioms_unified}
We define the model class of systems under study by the following axiomatic conditions on the microscopic geometry and dynamics, which are assumed to hold uniformly for all environmental parameters $\theta \in \Thetacal$.

\begin{enumerate}[label=(\roman*), wide, labelindent=0pt]
    \item \textbf{Boundary Regularity.}
    The value $0$ is a regular value of the map $G(\cdot, \theta): \T^k \to \R$ for every $\theta \in \Thetacal$. That is, for every $(y, \theta)$ such that $G(y, \theta) = 0$, the gradient with respect to the spatial variable does not vanish: $\nabla_y G(y, \theta) \neq \mathbf{0}$.

    \item \textbf{Microscopic Symmetry and Time-Reversibility.}
    The microscopic system possesses fundamental symmetries that embody time-reversal invariance.
    \begin{enumerate}
        \item \textit{Velocity Space Symmetry:} The velocity space $\Vcal$ is symmetric with respect to the origin: if $v \in \Vcal$, then $-v \in \Vcal$.
        \item \textit{Geometric Time-Reversal Invariance:} The spatial domain of free motion, $\Ocal(\theta)$, is invariant under the parity transformation $y \mapsto -y$, enforced by requiring $G(-y, \theta) = G(y, \theta)$.
    \end{enumerate}

    \item \textbf{Local Geometric Condition: Uniform Strict Convexity.}
    The boundary of the domain of free motion, $\partial\Ocal(\theta)$, is uniformly strictly convex with respect to the domain. Specifically, there exists a constant $\kappa_{\min} > 0$, independent of $\theta$, such that for all $y \in \partial\Ocal(\theta)$, every principal curvature of the boundary with respect to the inward-pointing normal is bounded below by $\kappa_{\min}$.

    \item \textbf{Global Geometric Condition: The Finite Corridor Exclusion Axiom.}
    There exists a uniform constant $L_{\max} < \infty$, independent of $\theta$, such that any straight line segment of length $L_{\max}$ embedded in the torus $\T^k$ must intersect the interior of the obstacle region, $\T^k \setminus \Ocal(\theta)$.
\end{enumerate}
\end{assumption}

\begin{remark}[The Foundational Role of the Geometric Axioms]
\label{rem:first_principles_axioms}
The axioms presented in \cref{ass:fundamental_axioms_unified}, while mathematically significant and defining a specific class of idealized (uniformly hyperbolic) systems, should not be interpreted as a collection of convenient technical assumptions. Rather, they constitute a minimal set of physically-motivated geometric axioms designed to guarantee the emergence of robust statistical properties from deterministic laws. Each axiom corresponds to a distinct physical or statistical requirement, providing a clear and verifiable foundation for our theory.

\begin{enumerate}[label=(\alph*), wide, labelindent=0pt]
    \item \textbf{Determinism (Axiom i).} The regularity of the boundary, enforced by the condition $\nabla_y G \neq \mathbf{0}$, is a direct translation of the principle of determinism in classical mechanics. The direction of the constraint force exerted by the boundary during a collision is given by the normal vector $\mathbf{n} \propto \nabla_y G$. If the gradient were to vanish at a point on the boundary, the collision law at that point would become ill-defined, violating the deterministic nature of the microscopic evolution. This axiom therefore excludes non-physical singularities such as cusps, ensuring that every microscopic interaction is unambiguous.

    \item \textbf{Time-Reversal Invariance (Axiom ii).} The symmetries of the velocity space ($\Vcal = -\Vcal$) and the geometric domain ($G(y,\theta) = G(-y,\theta)$) are the direct embodiment of the time-reversal invariance of the underlying Hamiltonian dynamics. A system governed by a time-independent Hamiltonian is symmetric with respect to the transformation $(t, y, v) \mapsto (-t, -y, -v)$. These axioms ensure that this fundamental symmetry is preserved in the presence of boundary interactions. Methodologically, this symmetry acts as a crucial centering condition, guaranteeing that the mean microscopic velocity vanishes ($\langle v \rangle_\theta = \mathbf{0}$, as proven in \cref{prop:symmetry_consequences}). This is a deliberate modeling choice designed to eliminate any trivial macroscopic drift and to isolate the more subtle emergent phenomena, namely, diffusion and a non-convex potential, that arise purely from the system's chaotic fluctuations.

    \item \textbf{The Ergodic Hypothesis (Axioms iii \& iv).} The axioms of uniform strict convexity and the exclusion of finite corridors are, taken together, the geometric conditions that guarantee the system is a uniformly hyperbolic (Anosov) system. This provides the strongest form of deterministic chaos and is the mechanical basis for the ergodic hypothesis of statistical mechanics. The ergodic hypothesis, which equates time averages with phase space averages, is the cornerstone that connects microscopic dynamics to macroscopic thermodynamics and transport phenomena. For a deterministic system, robust chaotic mixing is the only viable mechanism for ensuring this hypothesis holds.
    \begin{itemize}[wide]
        \item The \textit{local} condition of convexity (iii) acts as the engine of chaos, providing a dispersing mechanism at each collision that creates an exponential sensitivity to initial conditions.
        \item The \textit{global} condition of finite corridors (iv) ensures that this chaotic mixing is persistent and not short-circuited by trajectories that can avoid collisions indefinitely.
    \end{itemize}
    These two axioms are therefore the concrete, verifiable geometric conditions required for the emergence of robust statistical laws from deterministic motion.
\end{enumerate}
In summary, these axioms collectively define a class of mechanical systems that is both physically self-consistent and mathematically tractable, ensuring the system exhibits the rich statistical behavior that is the subject of this paper. The results derived from them are thus consequences of a foundational framework where each axiom serves a distinct physical or statistical purpose.
\end{remark}

\subsubsection{Immediate Consequences of the Axioms}

The axioms stated above have profound and immediate consequences for the geometric and measure-theoretic structure of the system. We establish these foundational results here. 

\begin{proposition}[Smoothness of the Collision Boundary]
\label{prop:smooth_boundary}
Under Assumption~\ref{ass:fundamental_axioms_unified}, for each $\theta \in \Thetacal$, the collision boundary $\partial\Ocal(\theta)$ is a smooth, compact, embedded submanifold of $\T^k$ of codimension 1.
\end{proposition}

\begin{proofof}{Proposition \ref{prop:smooth_boundary}}

The proof proceeds by applying the Regular Level Set Theorem to the obstacle-defining function $G$. We will demonstrate that for any fixed environmental parameter $\theta$, the collision boundary is the regular level set of a smooth function, from which all the claimed properties follow.

\begin{enumerate}[label=\textbf{Step \arabic*:}, wide, labelindent=0pt]
    \item \textbf{Formulation as a Level Set.} Let an arbitrary $\theta_0 \in \Thetacal$ be fixed. We define a function $g_{\theta_0}: \T^k \to \R$ by restricting the global function $G$ to the slice corresponding to $\theta_0$:
    \begin{equation*}
        g_{\theta_0}(y) \coloneqq G(y, \theta_0).
    \end{equation*}
    By construction in Definition~\ref{def:obstacle_func}, the collision boundary is precisely the preimage of $0$ under this map:
    \begin{equation*}
        \partial\Ocal(\theta_0) = g_{\theta_0}^{-1}(\{0\}).
    \end{equation*}
    Since $G \in C^\infty(\T^k \times \Thetacal)$, the partial map $g_{\theta_0}$ obtained by fixing the second argument is also of class $C^\infty$.

    \item \textbf{The Regular Level Set Theorem.} Our argument rests on the following standard result from differential geometry: Let $F: M \to N$ be a $C^\infty$ map between smooth manifolds. A point $c \in N$ is a regular value of $F$ if for every point $p$ in the preimage $F^{-1}(\{c\})$, the differential $dF_p: T_pM \to T_cN$ is a surjective linear map. If $c$ is a regular value, then the preimage $F^{-1}(\{c\})$ is a properly embedded submanifold of $M$ with codimension equal to the dimension of $N$.

    \item \textbf{Verification of Hypotheses.} We now rigorously verify the hypotheses of this theorem for our specific case.
    \begin{enumerate}[label=(\roman*), wide, labelindent=0pt]
        \item The manifolds are $M = \T^k$ and $N = \R$, both of which are smooth.
        \item The map is $F = g_{\theta_0}: \T^k \to \R$, which we established is $C^\infty$ in Step 1.
        \item The value in question is $c=0 \in \R$.
        \item We must verify that $0$ is a regular value. Let $y_0$ be an arbitrary point in the preimage, i.e., $y_0 \in g_{\theta_0}^{-1}(\{0\})$. We must show that the differential $(dg_{\theta_0})_{y_0}: T_{y_0}\T^k \to T_0\R$ is surjective. The tangent spaces can be identified with Euclidean spaces, $T_{y_0}\T^k \cong \R^k$ and $T_0\R \cong \R$. In local coordinates, the differential $(dg_{\theta_0})_{y_0}$ is represented by the Jacobian matrix of $g_{\theta_0}$ at $y_0$. For a scalar-valued function on a Euclidean-type space, this is the row vector of partial derivatives, which is the transpose of the gradient vector:
        \begin{equation*}
            (dg_{\theta_0})_{y_0} \longleftrightarrow (\nabla_y g_{\theta_0}(y_0))^\top = (\nabla_y G(y_0, \theta_0))^\top.
        \end{equation*}
        A linear map from $\R^k$ to $\R$ is surjective if and only if it is not the zero map. This, in turn, is true if and only if the representing vector is non-zero. Assumption~\ref{ass:fundamental_axioms_unified} states precisely that for any $(y, \theta)$ with $G(y, \theta)=0$, we have $\nabla_y G(y, \theta) \neq \mathbf{0}$. Thus, for our chosen $y_0 \in \partial\Ocal(\theta_0)$, we are guaranteed that $\nabla_y G(y_0, \theta_0) \neq \mathbf{0}$. This confirms that the differential is surjective. Since $y_0$ was an arbitrary point in the level set, we conclude that $0$ is a regular value of $g_{\theta_0}$.
    \end{enumerate}

    \item \textbf{The Manifold Structure.} All hypotheses being satisfied, the Regular Level Set Theorem allows us to conclude that $\partial\Ocal(\theta_0)$ is a properly embedded submanifold of $\T^k$. The codimension of this submanifold is $\dim(\R) = 1$.

    \item \textbf{Topological Properties.}
    \begin{enumerate}[label=(\roman*), wide, labelindent=0pt]
        \item \textbf{Closedness.} The set $\{0\}$ is a closed subset of $\R$. Since $g_{\theta_0}$ is a $C^\infty$ map, it is continuous. The preimage of a closed set under a continuous map is closed. Therefore, $\partial\Ocal(\theta_0) = g_{\theta_0}^{-1}(\{0\})$ is a closed subset of $\T^k$.
        \item \textbf{Compactness.} The manifold $\T^k$ is compact. As a closed subset of a compact space, $\partial\Ocal(\theta_0)$ is itself compact.
    \end{enumerate}

    \item \textbf{Generalization.} The choice of $\theta_0 \in \Thetacal$ was arbitrary. The argument holds for any $\theta \in \Thetacal$, establishing the proposition.
\end{enumerate}
\end{proofof}

\begin{proposition}[Properties of the Reflection Map]
\label{prop:reflection_properties}
For each $\theta \in \Thetacal$, the specular reflection map 
$R_\theta: \partial_-\Ucal_{\text{phys}}(\theta) \to \partial_+\Ucal_{\text{phys}}(\theta)$
is a $C^\infty$-smooth diffeomorphism. Furthermore, it is an involution and an isometry.
\end{proposition}

\begin{proofof}{Proposition \ref{prop:reflection_properties}}

Let $(y^-, v^-)$ be an arbitrary point in the domain $\partial_-\Ucal_{\text{phys}}(\theta)$. By \cref{def:reflection_map}, the specular reflection map $R_\theta$ is given by the pair of equations:
\begin{align}
    y^+ &= y^-, \label{eq:proof_y_map_appendix} \\
    v^+ &= v^- - 2\scpr{v^-}{\mathbf{n}(y^-, \theta)} \mathbf{n}(y^-, \theta). \label{eq:proof_v_map_appendix}
\end{align}
We will prove each claimed property in sequence.

\begin{enumerate}[label=\textbf{Step \arabic*:}, wide, labelindent=0pt]
    \item \textbf{Proof that $R_\theta$ is a $C^\infty$-Diffeomorphism.}
The proof that $R_\theta$ is a diffeomorphism proceeds in three steps. First, we establish its smoothness. Second, we show it maps the domain to the claimed codomain. Third, we show it has a smooth inverse.

\begin{enumerate}[label=(\roman*), wide, labelindent=0pt]
    \item \textbf{Smoothness of $R_\theta$.} The map on the position component, \cref{eq:proof_y_map_appendix}, is the identity map on the manifold $\partial\Ocal(\theta)$, which is trivially $C^\infty$. The map on the velocity component, \cref{eq:proof_v_map_appendix}, is a composition of several operations: vector addition, scalar multiplication, and the inner product $\scpr{\cdot}{\cdot}$. These are all $C^\infty$ operations on the tangent bundle. The smoothness of the map $v^+$ therefore depends entirely on the smoothness of the outward unit normal vector field $\mathbf{n}(y, \theta)$ as a function of $y \in \partial\Ocal(\theta)$. The normal vector is defined as $\mathbf{n}(y, \theta) = \nabla_y G(y, \theta) / \|\nabla_y G(y, \theta)\|$. By the problem's setup, the function $G: \T^k \times \Thetacal \to \R$ is $C^\infty$. Consequently, its partial derivative with respect to $y$, the gradient vector field $\nabla_y G$, is also a $C^\infty$ map from $\T^k \times \Thetacal$ to $\R^k$. 
    
    By \cref{ass:fundamental_axioms_unified}, for any $(y, \theta)$ on the collision boundary (i.e., $G(y,\theta)=0$), we have $\nabla_y G(y, \theta) \neq \mathbf{0}$. This ensures that the denominator $\|\nabla_y G(y, \theta)\|$ is strictly positive. Since the Euclidean norm is a smooth function away from the origin, the map $y \mapsto \|\nabla_y G(y, \theta)\|$ is $C^\infty$ on the boundary $\partial\Ocal(\theta)$. Thus, $\mathbf{n}(y, \theta)$ is the quotient of a $C^\infty$ vector field by a non-vanishing $C^\infty$ scalar function, which establishes that $\mathbf{n}$ is a $C^\infty$ vector field on $\partial\Ocal(\theta)$. As all constituent parts of the definition of $v^+$ are $C^\infty$ maps, their composition is $C^\infty$. Therefore, $R_\theta$ is a $C^\infty$ map.
    
    \item \textbf{Verification of Domain and Codomain.} The domain of $R_\theta$ is $\partial_-\Ucal_{\text{phys}}(\theta)$, characterized by $y^- \in \partial\Ocal(\theta)$ and $\scpr{v^-}{\mathbf{n}(y^-, \theta)} < 0$. We must show that the image $(y^+, v^+)$ lies in $\partial_+\Ucal_{\text{phys}}(\theta)$, which requires showing that $y^+ \in \partial\Ocal(\theta)$ and $\scpr{v^+}{\mathbf{n}(y^+, \theta)} > 0$. Since $y^+=y^-$, the position condition is trivially satisfied and the normal vector is unchanged. We compute the inner product of the post-collision velocity with the normal vector:
\begin{align*}
    \scpr{v^+}{\mathbf{n}(y^-, \theta)} &= \scpr{v^- - 2\scpr{v^-}{\mathbf{n}(y^-, \theta)}\mathbf{n}(y^-, \theta)}{\mathbf{n}(y^-, \theta)} \\
    &= \scpr{v^-}{\mathbf{n}(y^-, \theta)} - 2\scpr{v^-}{\mathbf{n}(y^-, \theta)}\scpr{\mathbf{n}(y^-, \theta)}{\mathbf{n}(y^-, \theta)} \\
    &= \scpr{v^-}{\mathbf{n}(y^-, \theta)} - 2\scpr{v^-}{\mathbf{n}(y^-, \theta)} && \text{(since $\|\mathbf{n}\|=1$)} \\
    &= -\scpr{v^-}{\mathbf{n}(y^-, \theta)}.
\end{align*}
Since $(y^-, v^-) \in \partial_-\Ucal_{\text{phys}}(\theta)$, we have $\scpr{v^-}{\mathbf{n}(y^-, \theta)} < 0$ by definition. It follows immediately that $\scpr{v^+}{\mathbf{n}(y^-, \theta)} > 0$, which confirms that the image point $(y^+, v^+)$ lies in the codomain $\partial_+\Ucal_{\text{phys}}(\theta)$.
    \item \textbf{Existence of a Smooth Inverse.} A map is a diffeomorphism if it is $C^\infty$ and has a $C^\infty$ inverse. As we will prove in part (ii), the map $R_\theta$ is its own inverse. Since $R_\theta$ has been shown to be $C^\infty$, its inverse is also $C^\infty$. Therefore, $R_\theta$ is a $C^\infty$-diffeomorphism from its domain to its codomain.
\end{enumerate}
\item \textbf{Proof that $R_\theta$ is an Involution.} We must show that applying the map $R_\theta$ twice returns the original state. Let $(y^+, v^+) = R_\theta(y^-, v^-)$ and $(y^{++}, v^{++}) = R_\theta(y^+, v^+)$.
\begin{enumerate}[label=(\roman*), wide, labelindent=0pt]
    \item The position component is trivial: $y^{++} = y^+ = y^-$.
    \item For the velocity component, we apply the reflection formula to $v^+$:
    \begin{align*}
        v^{++} &= v^+ - 2\scpr{v^+}{\mathbf{n}(y^+, \theta)}\mathbf{n}(y^+, \theta) \\
        &= v^+ - 2\left(-\scpr{v^-}{\mathbf{n}(y^-, \theta)}\right)\mathbf{n}(y^-, \theta) && \text{(using } y^+=y^- \text{ and the result from (i))} \\
        &= v^+ + 2\scpr{v^-}{\mathbf{n}(y^-, \theta)}\mathbf{n}(y^-, \theta) \\
        &= \left(v^- - 2\scpr{v^-}{\mathbf{n}(y^-, \theta)}\mathbf{n}(y^-, \theta)\right) + 2\scpr{v^-}{\mathbf{n}(y^-, \theta)}\mathbf{n}(y^-, \theta) && \text{(substituting the definition of } v^+\text{)} \\
        &= v^-.
    \end{align*}
\end{enumerate}
Thus, $(y^{++}, v^{++}) = (y^-, v^-)$, which proves that $R_\theta \circ R_\theta = \mathrm{Id}$ on its domain of definition.

\item \textbf{Proof that $R_\theta$ is an Isometry.} We demonstrate that the magnitude of the velocity vector is preserved across the collision by computing the squared norm of $v^+$.
\begin{align*}
    \|v^+\|^2 &= \scpr{v^+}{v^+} \\
    &= \scpr{v^- - 2\scpr{v^-}{\mathbf{n}}\mathbf{n}}{v^- - 2\scpr{v^-}{\mathbf{n}}\mathbf{n}} && \text{(where } \mathbf{n} = \mathbf{n}(y^-, \theta)\text{)} \\
    &= \scpr{v^-}{v^-} - 2\scpr{v^-}{2\scpr{v^-}{\mathbf{n}}\mathbf{n}} + \scpr{2\scpr{v^-}{\mathbf{n}}\mathbf{n}}{2\scpr{v^-}{\mathbf{n}}\mathbf{n}} \\
    &= \|v^-\|^2 - 4\scpr{v^-}{\mathbf{n}}\scpr{v^-}{\mathbf{n}} + 4\scpr{v^-}{\mathbf{n}}^2 \scpr{\mathbf{n}}{\mathbf{n}} \\
    &= \|v^-\|^2 - 4\scpr{v^-}{\mathbf{n}}^2 + 4\scpr{v^-}{\mathbf{n}}^2 \|\mathbf{n}\|^2.
\end{align*}
Since $\mathbf{n}$ is a unit vector, $\|\mathbf{n}\|^2 = 1$. The last two terms cancel, yielding
\begin{equation*}
    \|v^+\|^2 = \|v^-\|^2.
\end{equation*}
As the norm is non-negative, taking the square root of both sides gives $\|v^+\| = \|v^-\|$. This confirms that kinetic energy is conserved and the map is an isometry on the velocity component.
\end{enumerate}
\end{proofof}

\begin{proposition}[Consequences of Microscopic Symmetry]
\label{prop:symmetry_consequences}
Let Assumption \ref{ass:fundamental_axioms_unified} hold. Then for each $\theta \in \Thetacal$, the normalized Liouville measure $\mu_\theta$ is invariant under the time-reversal map $\mathfrak{R}(y,v)=(-y,-v)$, and the mean velocity with respect to this measure is identically zero:
$ \int_{\Ucal_{\mathrm{phys}}(\theta)} v \, d\mu_\theta(z) = \mathbf{0}$.
\end{proposition}

\begin{proofof}{Proposition \ref{prop:symmetry_consequences}}
Let an arbitrary $\theta \in \Thetacal$ be fixed.

\begin{enumerate}[label=\textbf{Step \arabic*:}, wide, labelindent=0pt]

\item\textbf{Invariance of the State Space.} We must show that if $z=(y,v) \in \Ucal_{\mathrm{phys}}(\theta)$, then $\mathfrak{R}(z)=(-y,-v)$ is also in $\Ucal_{\mathrm{phys}}(\theta)$. A point $(y,v)$ is in the physical space if $v \in \Vcal$ and $G(y,\theta) \ge 0$.
By Assumption \ref{ass:fundamental_axioms_unified}(i), if $v \in \Vcal$, then $-v \in \Vcal$.
By Assumption \ref{ass:fundamental_axioms_unified}(ii), $G(-y, \theta) = G(y, \theta)$. Since $G(y,\theta) \ge 0$, it follows that $G(-y, \theta) \ge 0$.
Both conditions are met, so $\mathfrak{R}(z) \in \Ucal_{\mathrm{phys}}(\theta)$. Thus, the state space is invariant under $\mathfrak{R}$.

\item \textbf{Invariance of the Measure.} The normalized Liouville measure is given by $d\mu_\theta(z) = C_\theta \, dy \, dv$, where $C_\theta$ is a normalization constant and $z=(y,v)$. We must show that for any measurable set $A \subset \Ucal_{\mathrm{phys}}(\theta)$, we have $\mu_\theta(A) = \mu_\theta(\mathfrak{R}(A))$. This is equivalent to showing that the pullback of the measure form under $\mathfrak{R}$ is itself.
Let's consider the change of variables $z' = \mathfrak{R}(z)$, which corresponds to $y' = -y$ and $v' = -v$. The Jacobian of this linear transformation is
\begin{equation*}
    J = \det \begin{pmatrix} \frac{\partial y'}{\partial y} & \frac{\partial y'}{\partial v} \\ \frac{\partial v'}{\partial y} & \frac{\partial v'}{\partial v} \end{pmatrix} = \det \begin{pmatrix} -I_k & 0 \\ 0 & -I_k \end{pmatrix} = (-1)^k (-1)^k = 1,
\end{equation*}
where $I_k$ is the $k \times k$ identity matrix. Therefore, the volume element is preserved: $dy' dv' = |\det(J)| \, dy \, dv = dy \, dv$.
The integral of any test function $\phi(z)$ over a set $A$ is
\begin{align*}
    \int_{\mathfrak{R}(A)} \phi(z') \, d\mu_\theta(z') &= C_\theta \int_{\mathfrak{R}(A)} \phi(y',v') \, dy' \, dv' \\
    &= C_\theta \int_{A} \phi(-y, -v) \, dy \, dv \quad \text{(by change of variables)} \\
    &= \int_A (\phi \circ \mathfrak{R})(z) \, d\mu_\theta(z).
\end{align*}
Taking $\phi = \mathbf{1}$, we find that $\mu_\theta(\mathfrak{R}(A)) = \mu_\theta(A)$. The measure is invariant.

\item \textbf{Vanishing of the Mean Velocity.} This is the central result. Let $\mathbf{I}$ denote the vector-valued integral for the mean velocity:
\begin{equation*}
    \mathbf{I} = \int_{\Ucal_{\mathrm{phys}}(\theta)} v \, d\mu_\theta(z).
\end{equation*}
We perform a change of variables $z' = \mathfrak{R}(z) = (-y, -v)$, so $z = \mathfrak{R}(z')$. From part (i), the domain of integration $\Ucal_{\mathrm{phys}}(\theta)$ is invariant under this map. From part (ii), the measure is also invariant, $d\mu_\theta(z) = d\mu_\theta(z')$. The velocity vector in the integrand transforms as $v = -v'$.
Substituting these into the integral gives:
\begin{align*}
    \mathbf{I} &= \int_{\Ucal_{\mathrm{phys}}(\theta)} v \, d\mu_\theta(z) \\
    &= \int_{\mathfrak{R}(\Ucal_{\mathrm{phys}}(\theta))} (-v') \, d\mu_\theta(z') \\
    &= -\int_{\Ucal_{\mathrm{phys}}(\theta)} v' \, d\mu_\theta(z').
\end{align*}
Since the variable of integration is a dummy variable, the last integral is simply $-\mathbf{I}$. We have thus proven that
\begin{equation*}
    \mathbf{I} = -\mathbf{I},
\end{equation*}
which implies $2\mathbf{I} = \mathbf{0}$, and therefore $\mathbf{I} = \mathbf{0}$. The mean velocity is identically zero.
\end{enumerate}
\end{proofof}

\subsubsection{Derivation of the Finite Horizon Theorem}

In this section, we lay the groundwork for proving the global well-posedness of the dynamics. At this stage, we have defined only the elementary, local laws of motion: the interior flow $\Psi_s$ and the specular reflection map $R_\theta$. A trajectory can only be constructed by composing these two laws. However, the validity of this composition for all time is not yet established. It could, a priori, lead to pathological behaviors such as an infinite number of collisions in a finite time, or a trajectory that never collides again. To analyze this, we introduce the central object of study for this section: the flight time. Its definition relies only on the local interior flow.

\begin{definition}[The Flight Time Function]
\label{def:flight_time_function}
Let $(y^+, v^+)$ be an arbitrary post-collision (outgoing) state in $\Sigma'_\theta \coloneqq \partial_+\Ucal_{\mathrm{phys}}(\theta)$. The flight time $\tau(y^+, v^+, \theta)$ is defined as the first positive time at which the trajectory starting from $(y^+, v^+)$ again reaches the boundary:
\begin{equation}
\tau(y^+, v^+, \theta) \coloneqq \inf \{ s > 0 \mid y^+ + s v^+ \in \partial\Ocal(\theta) \}.
\end{equation}
A priori, this infimum could be $0$ (immediate re-collision), $+\infty$ (no future collision), or a value that is not attained. The objective of the following three propositions is to prove, using only our geometric axioms, that $\tau$ is in fact a well-defined function that is uniformly bounded both from below and from above by positive, finite constants. This result, the Finite Horizon Theorem, will provide the rigorous justification for the global well-posedness of the dynamics.
\end{definition}

The well-posedness of the global dynamics depends on the exclusion of two pathological behaviors: trajectories that undergo an infinite number of collisions in a finite time, and trajectories that never collide again. The following propositions demonstrate that our geometric axioms are precisely the conditions required to uniformly exclude these pathologies.

\begin{proposition}[Uniform Exclusion of Grazing Collisions]
\label{prop:non_grazing}
Under Assumption \ref{ass:fundamental_axioms_unified}(iii) (Uniform Strict Convexity), there exists a constant $c_0 > 0$ such that for any trajectory, at any collision point $y \in \partial\Ocal(\theta)$, the pre-collision velocity $v^-$ satisfies $|\scpr{v^-}{\mathbf{n}(y, \theta)}| \ge c_0$.
\end{proposition}

\begin{proofof}{Proposition \ref{prop:non_grazing}}
The proof proceeds by contradiction. The entire argument is self-contained and relies only on the elementary laws of motion, the interior flow (\cref{def:flow_vector_field}) and the specular reflection map (\cref{def:reflection_map}), which are locally well-defined. It does not assume the existence of a globally defined flow for all time. We first demonstrate that the existence of even a single perfectly grazing collision leads to a direct contradiction with the local well-posedness of the dynamics. We then leverage this result and the compactness of the ambient state space to establish that the exclusion must be uniform.

\begin{enumerate}[label=\textbf{Step \arabic*:}, wide, labelindent=0pt]

\item \textbf{Setup by Contradiction and Existence of a Grazing Limit State.} Our goal is to prove that the magnitude of the normal component of the incoming velocity is uniformly bounded below by a positive constant. Let the set of all possible pre-collision states that can occur in the dynamics be denoted by $S_{\text{pre}}$. Let the function $f: \Efrak_{\text{phys}} \to \R$ be defined by $f(y,v,\theta) \coloneqq |\scpr{v}{\mathbf{n}(y,\theta)}|$ for points on the boundary and $f=0$ otherwise.

Assume, for the sake of contradiction, that the proposition is false. This implies that the infimum of the function $f$ over the set of all pre-collision states is zero:
\begin{equation*}
    \inf_{(y,v,\theta) \in S_{\text{pre}}} f(y,v,\theta) = 0.
\end{equation*}
By the definition of the infimum, this guarantees the existence of a sequence of pre-collision states, $\{z_n\}_{n=1}^\infty \subset S_{\text{pre}}$ where $z_n = (y_n, v_n, \theta_n)$, such that the normal component of their velocities converges to zero:
\begin{equation}
    \lim_{n \to \infty} f(z_n) = \lim_{n \to \infty} |\scpr{v_n}{\mathbf{n}(y_n, \theta_n)}| = 0.
\end{equation}
Each state $z_n$ belongs to the total boundary manifold, $\mathcal{B} \coloneqq \{ (y,v,\theta) \in \T^k \times \Vcal \times \Thetacal \mid G(y,\theta) = 0 \}$. Since $\T^k, \Vcal,$ and $\Thetacal$ are compact spaces by assumption, their product is compact by Tychonoff's theorem. As $\mathcal{B}$ is the zero level set of a continuous function $G$, it is a closed subset of this compact product space and is therefore itself a compact metric space.

Since the sequence $\{z_n\}$ lies in the compact space $\mathcal{B}$, it must possess a convergent subsequence, which we relabel as $\{z_n\}$ for simplicity. Let the limit of this subsequence be $z_0 = (y_0, v_0, \theta_0) \in \mathcal{B}$. By the continuity of the function $f$ on $\mathcal{B}$, we have:
\begin{equation*}
    f(z_0) = \lim_{n \to \infty} f(z_n) = 0.
\end{equation*}
This establishes the existence of a limit state $z_0$ that represents a perfectly grazing collision, with $\scpr{v_0}{\mathbf{n}(y_0, \theta_0)} = 0$. Let $(y_0, v_0^+) = R_{\theta_0}(y_0, v_0)$ be the corresponding post-collision state. From the properties of the reflection map, we also have $\scpr{v_0^+}{\mathbf{n}(y_0, \theta_0)} = -\scpr{v_0}{\mathbf{n}(y_0, \theta_0)} = 0$, meaning the post-collision velocity vector is tangent to the boundary $\partial\Ocal(\theta_0)$ at the point of collision $y_0$.

\item \textbf{Local Geometric Analysis of the Post-Collision Trajectory.} We now analyze the trajectory of a particle immediately following this grazing collision. The trajectory for time $s \ge 0$ is governed by the ODE $\dot{y}(s) = v_0^+$, with initial condition $y(0) = y_0$. The solution is $y(s) = y_0 + s v_0^+$. For the dynamics to be well-defined, the particle's trajectory must remain within the physical domain for some short time after the collision. Mathematically, this means there must exist an $\epsilon > 0$ such that for all $s \in [0, \epsilon)$, the condition $G(y(s), \theta_0) \ge 0$ holds. To analyze this condition, we define the auxiliary function $g(s) \coloneqq G(y(s), \theta_0) = G(y_0 + s v_0^+, \theta_0)$ and examine its behavior around $s=0$ using a Taylor expansion.
\begin{enumerate}[label=(\roman*), wide, labelindent=0pt]
    \item \textbf{Value at $s=0$:} By definition of a collision, $y_0$ lies on the boundary, so $g(0) = G(y_0, \theta_0) = 0$.

    \item \textbf{First Derivative at $s=0$:} Applying the chain rule, we have:
    \begin{equation*}
        g'(0) = \scpr{\nabla_y G(y_0, \theta_0)}{v_0^+} = \|\nabla_y G(y_0, \theta_0)\| \cdot \scpr{\mathbf{n}(y_0, \theta_0)}{v_0^+}.
    \end{equation*}
    From Step 1, the grazing collision implies $\scpr{\mathbf{n}(y_0, \theta_0)}{v_0^+} = 0$. Therefore, we must have $g'(0) = 0$.

    \item \textbf{Second Derivative at $s=0$:} Differentiating again with respect to $s$:
    \begin{equation*}
        g''(0) = (v_0^+)^\top H_y(y_0, \theta_0) v_0^+,
    \end{equation*}
    where $H_y$ is the Hessian matrix of $G$ with respect to the $y$ variables.
\end{enumerate}

\item \textbf{The Local Contradiction via Strict Convexity.} We now invoke the Uniform Strict Convexity of Obstacles (\cref{ass:fundamental_axioms_unified}(iii)). This geometric condition relates the Hessian of the defining function $G$ to the second fundamental form $\mathrm{II}_{-\mathbf{n}}$ with respect to the inward-pointing normal: for any tangent vector $\xi$ to the boundary (such as our grazing velocity $v_0^+$), we have
\begin{equation*}
    \xi^\top H_y(y_0) \xi = -\|\nabla_y G(y_0)\| \cdot \mathrm{II}_{-\mathbf{n}_0}(\xi, \xi).
\end{equation*}
The strict convexity axiom states that there exists a uniform constant $\kappa_{\min} > 0$ such that $\mathrm{II}_{-\mathbf{n}_0}(\xi, \xi) \ge \kappa_{\min} \|\xi\|^2$. Applying this to our grazing velocity $v_0^+$:
\begin{equation*}
    g''(0) \le -\|\nabla_y G(y_0, \theta_0)\| \cdot (\kappa_{\min} \|v_0^+\|^2).
\end{equation*}
Since the velocity is non-zero ($\|v_0^+\|>0$) and the boundary is regular ($\|\nabla_y G\| > 0$ by \cref{ass:fundamental_axioms_unified}), this gives a strict upper bound: $g''(0) < 0$.

The Taylor expansion of $g(s)$ around $s=0$ is therefore $g(s) = \frac{s^2}{2} g''(0) + o(s^2)$. Since $g''(0)$ is strictly negative, there must exist an $\epsilon > 0$ such that for all $s \in (0, \epsilon)$, we have $g(s) < 0$. This means the particle's trajectory immediately enters the forbidden region where $G < 0$, which is a direct contradiction of the local well-posedness of the flow. The initial assumption, that a sequence of pre-collision states can approach a grazing state, must be false.

\item \textbf{Uniformity of the Exclusion.} The contradiction argument in Steps 1-3 establishes that the infimum of $f(z)$ over the set of all pre-collision states $S_{\text{pre}}$ must be strictly positive. We now show this lower bound is uniform.

Let $\overline{S_{\text{pre}}}$ be the closure of the set of all pre-collision states, taken within the ambient compact space $\mathcal{B}$. As a closed subset of the compact space $\mathcal{B}$, the set $\overline{S_{\text{pre}}}$ is itself compact. The contradiction argument has shown that no limit point of a sequence in $S_{\text{pre}}$ can be a grazing state. This means that the set of all limit points, $\overline{S_{\text{pre}}}$, must be disjoint from the set of all grazing states, $\mathcal{B}_0 \coloneqq \{z \in \mathcal{B} \mid f(z) = 0 \}$.

Therefore, the continuous function $f(z) = |\scpr{v}{\mathbf{n}(y,\theta)}|$ is strictly positive on the compact domain $\overline{S_{\text{pre}}}$. By the Extreme Value Theorem, a continuous, strictly positive function on a compact set must attain its minimum value, which must itself be a strictly positive constant. Let this minimum be $c_0$:
\begin{equation*}
    c_0 \coloneqq \min_{z \in \overline{S_{\text{pre}}}} f(z) > 0.
\end{equation*}
Since every pre-collision state $z \in S_{\text{pre}}$ is also in the closure $\overline{S_{\text{pre}}}$, it follows that for any such state, $f(z) \ge c_0$. This establishes the uniform lower bound.
\end{enumerate}
\end{proofof}

The exclusion of grazing collisions directly leads to a uniform lower bound on the time between collisions.

\begin{proposition}[Uniform Lower Bound on Flight Time]
\label{prop:tau_min_bound}
Under Assumption \ref{ass:fundamental_axioms_unified}(iii), there exists a uniform constant $\tau_{\min} > 0$ such that the flight time function defined in \cref{def:flight_time_function} satisfies $\tau(y^+, v^+, \theta) \ge \tau_{\min}$ for all outgoing states $(y^+, v^+) \in \Sigma'_\theta$.
\end{proposition}

\begin{proofof}{Proposition \ref{prop:tau_min_bound}}
The proof derives a lower bound for the flight time $\tau$ by analyzing the local geometry at the point of collision. The crucial non-degeneracy condition for this analysis is the uniform exclusion of grazing collisions. This property, the subject of \cref{prop:non_grazing}, was established as a direct consequence of the local geometric axiom of Uniform Strict Convexity (\cref{ass:fundamental_axioms_unified}(iii)) and is therefore logically independent of the global axioms used to establish the upper bound on flight time.

\begin{enumerate}[label=\textbf{Step \arabic*:}, wide, labelindent=0pt]

\item \textbf{Setup and Local Analysis.} Let an arbitrary post-collision state $(y^+, v^+, \theta)$ be given. We analyze the trajectory $y(s) = y^+ + s v^+$ for small $s > 0$. For the dynamics to be well-posed, the particle must remain in the physical domain for some short time after the collision, meaning $G(y(s), \theta) \ge 0$ for $s$ in some interval $[0, \epsilon)$. Let $\tau$ be the infimum of positive times at which the particle returns to the boundary, as defined in \cref{def:flight_time_function}. If such a finite time exists, we must have $G(y(\tau), \theta) = 0$. We define the auxiliary function $g(s) \coloneqq G(y^+ + s v^+, \theta)$ and analyze its behavior near $s=0$ via a Taylor expansion with a Lagrange remainder:
\begin{equation}
    g(s) = g(0) + g'(0)s + \frac{1}{2}g''(\zeta_s)s^2, \quad \text{for some } \zeta_s \in (0,s).
\end{equation}

\item \textbf{Bounding the Taylor Coefficients.} We now derive uniform bounds for the coefficients in the Taylor expansion, which are valid for any possible outgoing state.
\begin{enumerate}[label=(\roman*), wide, labelindent=0pt]
    \item \textbf{$g(0)$:} By definition of a collision, the starting point is on the boundary, so $g(0) = G(y^+, \theta) = 0$.

    \item \textbf{$g'(0)$:} The first derivative is $g'(0) = \scpr{\nabla_y G(y^+, \theta)}{v^+}$. We can rewrite this in terms of the unit normal:
    \begin{equation*}
        g'(0) = \|\nabla_y G(y^+, \theta)\| \cdot \scpr{\mathbf{n}(y^+,\theta)}{v^+}.
    \end{equation*}
    We now invoke two uniform bounds derived from our axioms:
    \begin{enumerate}[label=(\alph*),wide]
        \item By Assumption \ref{ass:fundamental_axioms_unified}(i), the gradient $\|\nabla_y G\|$ is non-zero on the compact boundary, so it is uniformly bounded below by a constant $C_G > 0$.
        \item By the specular reflection law, $\scpr{\mathbf{n}}{v^+} = -\scpr{\mathbf{n}}{v^-} = |\scpr{\mathbf{n}}{v^-}|$. By the crucial result of \cref{prop:non_grazing}, this normal component of velocity is uniformly bounded away from zero: $|\scpr{\mathbf{n}}{v^-}| \ge c_0 > 0$.
    \end{enumerate}
    Combining these, we find that the first derivative is uniformly positive:
    \begin{equation}
        g'(0) \ge C_G \cdot c_0 \eqqcolon C_1 > 0.
    \end{equation}

    \item \textbf{$g''(s)$:} The second derivative is $g''(s) = (v^+)^\top H_y(y(s), \theta) v^+$, where $H_y$ is the Hessian of $G$. The total physical state space $\Efrak_{\text{phys}}$ is compact. Since the function $G$ is smooth, its Hessian is a continuous matrix-valued function on this compact set, so its operator norm is uniformly bounded above by a constant $C_H < \infty$. The velocity is also uniformly bounded, $\|v^+\| \le v_{\max}$. Therefore, the second derivative is uniformly bounded in magnitude:
    \begin{equation}
        |g''(s)| \le \|H_y(y(s),\theta)\| \cdot \|v^+\|^2 \le C_H \cdot v_{\max}^2 \eqqcolon C_2 < \infty.
    \end{equation}
\end{enumerate}

\item \textbf{Derivation of the Lower Bound.} The Taylor expansion provides a rigorous lower bound for the function $g(s)$ for any $s > 0$:
\begin{equation*}
    g(s) = g'(0)s + \frac{1}{2}g''(\zeta_s)s^2 \ge C_1 s - \frac{C_2}{2}s^2.
\end{equation*}
For the particle to return to the boundary at a time $\tau > 0$, we must have $g(\tau) \le 0$ (it must be exactly zero if the trajectory does not enter the forbidden region). This implies:
\begin{equation*}
    0 \ge C_1 \tau - \frac{C_2}{2}\tau^2 \implies \frac{C_2}{2}\tau^2 \ge C_1 \tau.
\end{equation*}
Since we are seeking a positive time $\tau > 0$, we can divide by $\tau$ without changing the inequality:
\begin{equation*}
    \tau \ge \frac{2C_1}{C_2}.
\end{equation*}
By defining $\tau_{\min} \coloneqq 2C_1/C_2 = \frac{2 C_G c_0}{C_H v_{\max}^2}$, we have established a uniform positive lower bound on any possible flight time. This proves that the infimum in \cref{def:flight_time_function} cannot be zero and must be at least $\tau_{\min}$.
\end{enumerate}
\end{proofof}

The existence of a uniform lower bound on flight time excludes the possibility of an accumulation of collisions. We now show that our global geometric axiom excludes the possibility of trajectories never colliding again.

\begin{proposition}[Uniform Upper Bound on Flight Time]
\label{prop:tau_max_bound}
Under Assumption \ref{ass:fundamental_axioms_unified}(iv) (The Finite Corridor Exclusion Axiom), there exists a uniform constant $\tau_{\max} < \infty$ such that the flight time function defined in \cref{def:flight_time_function} satisfies $\tau(y^+, v^+, \theta) \le \tau_{\max}$ for all outgoing states $(y^+, v^+) \in \Sigma'_\theta$.
\end{proposition}

\begin{proofof}{Proposition \ref{prop:tau_max_bound}}
The proof is a direct, constructive derivation from the global geometric axiom on obstacle placement (\cref{ass:fundamental_axioms_unified}(iv)). We demonstrate that this geometric constraint forces any trajectory to collide within a uniformly bounded time, thereby establishing the finite horizon property as a theorem. The proof does not rely on any topological arguments involving the continuity of the flight time function.

\begin{enumerate}[label=\textbf{Step \arabic*:}, wide, labelindent=0pt]

\item \textbf{Setup by Contradiction.} Assume, for the sake of contradiction, that the flight time is not uniformly bounded above. This implies that for any arbitrarily large time $T_{\text{test}} > 0$, there exists at least one post-collision state $(y^+, v^+, \theta)$ such that the trajectory segment
\begin{equation*}
    \gamma(s) = y^+ + sv^+, \quad s \in [0, T_{\text{test}}],
\end{equation*}
lies entirely within the domain of free motion, $\Ocal(\theta)$.

\item \textbf{Path Length of the Hypothetical Trajectory.}
The length of this hypothetical free-flight path is $L(T_{\text{test}}) = \|v^+\| \cdot T_{\text{test}}$. By the compactness of the microscopic velocity space $\Vcal$ (\cref{def:micro_fiber}), the speed of any non-stationary particle is uniformly bounded below by a strictly positive constant, $v_{\min} > 0$. Therefore, the length of our hypothetical path must satisfy:
\begin{equation*}
    L(T_{\text{test}}) \ge v_{\min} \cdot T_{\text{test}}.
\end{equation*}

\item \textbf{Invoking the Global Geometric Axiom.} We now choose a specific test time to force a contradiction. Let $L_{\max}$ be the uniform maximum path length from the Finite Corridor Exclusion Axiom (\cref{ass:fundamental_axioms_unified}(iv)). We set our test time to be:
\begin{equation*}
    T_{\text{test}} = \frac{L_{\max}}{v_{\min}}.
\end{equation*}
Our contradiction hypothesis from Step 1 guarantees the existence of a trajectory segment of this duration that lies entirely within the domain of free motion. The length of this segment is at least:
\begin{equation*}
    L(T_{\text{test}}) \ge v_{\min} \cdot \left( \frac{L_{\max}}{v_{\min}} \right) = L_{\max}.
\end{equation*}
This trajectory segment is, by construction, a straight line segment in the torus of length at least $L_{\max}$ that does not intersect the interior of the obstacle region. However, the Finite Corridor Exclusion Axiom (\cref{ass:fundamental_axioms_unified}(iv)) asserts that any straight line segment of length $L_{\max}$ must intersect the interior of the obstacle region. This is a direct contradiction.

\item \textbf{Uniformly Bounded.}
The initial assumption that the flight time is not uniformly bounded must be false. Therefore, no trajectory can travel for a time greater than or equal to $L_{\max}/v_{\min}$ without a collision. This means the flight time $\tau$ for any trajectory is uniformly bounded above:
\begin{equation*}
    \tau \le \frac{L_{\max}}{v_{\min}} \eqqcolon \tau_{\max}.
\end{equation*}
Since both $L_{\max}$ (from the global geometric axiom) and $v_{\min}$ (from the compactness of the velocity space) are uniform positive constants of the system, independent of any specific state $(y,v,\theta)$, we have proven that the flight time $\tau$ is uniformly bounded above by $\tau_{\max}$. This proves that the infimum in \cref{def:flight_time_function} cannot be $+\infty$.

\end{enumerate}
\end{proofof}

The preceding propositions, derived from our fundamental geometric axioms, now culminate in the Finite Horizon Theorem.

\begin{theorem}[The Finite Horizon Theorem]
\label{thm:finite_horizon}
Let the system satisfy the Fundamental Geometric Axioms (\cref{ass:fundamental_axioms_unified}). Then the microscopic dynamics have a finite horizon. Specifically, there exist uniform constants $0 < \tau_{\min} \le \tau_{\max} < \infty$, independent of the specific trajectory or the environmental parameter $\theta$, such that the time $\tau$ between any two consecutive collisions satisfies the uniform bounds:
\begin{equation}
    \tau_{\min} \le \tau \le \tau_{\max}.
\end{equation}
\end{theorem}

\begin{proofof}{Theorem \ref{thm:finite_horizon}}
The theorem states that the flight time between any two consecutive collisions is uniformly bounded both from below and from above. The proof is a direct synthesis of the preceding propositions of this section, each of which was derived rigorously from the fundamental geometric axioms.

\begin{enumerate}[label=\textbf{Step \arabic*:}, wide, labelindent=0pt]

\item \textbf{The Uniform Lower Bound ($\tau \ge \tau_{\min}$).}
The existence of a uniform positive lower bound on the flight time was the subject of \cref{prop:tau_min_bound}. The proof demonstrated that for any outgoing state, the infimum in the definition of the flight time, $\tau = \inf \{ s>0 \dots \}$, cannot be zero. This result was a direct consequence of the local geometric axiom of Uniform Strict Convexity and the uniform exclusion of grazing collisions established in \cref{prop:non_grazing}. \cref{prop:tau_min_bound} rigorously establishes the existence of a constant $\tau_{\min} > 0$, independent of the specific state $(y,v,\theta)$, such that any possible flight time $\tau$ must satisfy
\begin{equation*}
    \tau \ge \tau_{\min}.
\end{equation*}
This bound ensures that an infinite number of collisions cannot accumulate in a finite time interval.

\item \textbf{The Uniform Upper Bound ($\tau \le \tau_{\max}$).}
The existence of a uniform upper bound on the flight time was the subject of \cref{prop:tau_max_bound}. The proof, via a contradiction argument, demonstrated that the infimum in the definition of the flight time cannot be infinite. This result was a direct and constructive derivation from the global geometric axiom of Finite Corridor Exclusion (\cref{ass:fundamental_axioms_unified}(iv)). \cref{prop:tau_max_bound} rigorously establishes the existence of a constant $\tau_{\max} < \infty$, independent of the specific state $(y,v,\theta)$, such that any possible flight time $\tau$ must satisfy
\begin{equation*}
    \tau \le \tau_{\max}.
\end{equation*}
This bound ensures that a particle cannot travel indefinitely without undergoing a collision.

\item \textbf{Finite Horizon.}
By combining the conclusions of \cref{prop:tau_min_bound} and \cref{prop:tau_max_bound}, we have rigorously established that the provisional quantity $\tau$ defined in \cref{def:flight_time_function} is neither zero nor infinite. The set of positive times $\{ s > 0 \mid y^+ + s v^+ \in \partial\Ocal(\theta) \}$ is therefore a non-empty set of real numbers that is bounded below by $\tau_{\min} > 0$. As a non-empty set of real numbers bounded from below, its infimum must exist.

Furthermore, we have shown that this infimum is bounded above by $\tau_{\max}$. Therefore, the flight time function $\tau$ is a well-defined map from the post-collision space into the compact interval $[\tau_{\min}, \tau_{\max}]$. This confirms that the dynamical property of having a finite horizon is not an independent axiom but is a provable theorem within our geometric framework, ensuring the global well-posedness of the microscopic dynamics.
\end{enumerate}
\end{proofof}

\begin{remark}[The Consequence of the Finite Horizon]
\label{rem:finite_horizon_consequence}
The Finite Horizon Theorem, now rigorously established from first principles, is the essential ingredient that allows for the unambiguous construction of the global flow in the subsequent sections. The theorem guarantees that the sequence of collision times $\{t_n\}_{n\ge 1}$, constructed recursively by $t_{n+1} = t_n + \tau_n$ where $\tau_n \in [\tau_{\min}, \tau_{\max}]$ is the $n$-th flight time, is a strictly increasing sequence with $t_{n+1}-t_n \ge \tau_{\min} > 0$. This implies that $\lim_{n \to \infty} t_n = \infty$. Consequently, there can only be a finite number of collisions in any finite time interval. This excludes pathological Zeno behavior and guarantees that the recursive construction of the flow by composing interior flows and boundary reflections is well-defined and unique for any $t \in \mathbb{R}$.
\end{remark}

\subsection{The Global Flow and its Invariant Measure}

With the Finite Horizon Theorem (\cref{thm:finite_horizon}) now rigorously established from our geometric first principles, we are in a position to provide a well-posed definition of the global microscopic flow for all time $t \in \mathbb{R}$. The theorem guarantees that the sequence of collision times is a discrete set of points bounded away from each other, ensuring that the trajectory does not experience an accumulation of collisions in finite time and is uniquely defined by composing the elementary motions.

\begin{definition}[The Flow Map $\Phi_t^\theta$]
\label{def:global_flow}
For a fixed $\theta \in \Thetacal$, the global flow map $\Phi_t^\theta: \Ucal_{\text{phys}}(\theta) \to \Ucal_{\text{phys}}(\theta)$ is defined for any initial state $z_0 = (y_0, v_0) \in \Ucal_{\text{phys}}(\theta)$ as follows. Let $t_0 = 0$ and $z(0) = z_0$.
\begin{enumerate}[label=(\roman*), wide, labelindent=0pt]
    \item Define the first exit time $\tau_1 = \inf \{ s > 0 \mid \Psi_s(z_0) \in \partial\Ucal_{\text{phys}}(\theta) \}$.
    \item For $t \in [0, \tau_1)$, the flow is given by the interior dynamics: $\Phi_t^\theta(z_0) = \Psi_t(z_0)$.
    \item At $t = \tau_1$, define the post-collision state $z_1 = R_\theta(\Psi_{\tau_1}(z_0))$.
    \item For $t > \tau_1$, the flow is defined recursively: $\Phi_t^\theta(z_0) = \Phi_{t-\tau_1}^\theta(z_1)$.
\end{enumerate}
This piecewise definition constructs the full trajectory $z(t) = \Phi_t^\theta(z_0)$.
\end{definition}

\begin{remark}[Well-Posedness of the Flow]
The Finite Horizon Theorem is the essential ingredient for this definition. The uniform lower bound, $\tau \ge \tau_{\min} > 0$, ensures that the sequence of collision times $\{t_n\}_{n\ge 1}$ is a strictly increasing sequence with $t_{n+1}-t_n \ge \tau_{\min}$. This implies that $\lim_{n \to \infty} t_n = \infty$, so there can only be a finite number of collisions in any finite time interval. This guarantees that the recursive construction of the flow is well-defined and unique for any $t \in \mathbb{R}$.
\end{remark}

A fundamental property of this flow, which is the foundation of the subsequent ergodic analysis, is that it preserves the natural volume measure on phase space.

\begin{proposition}[Invariance of the Liouville Measure]
\label{prop:liouville_invariance}
For each $\theta \in \Thetacal$, the global flow map $\Phi_t^\theta$ preserves the normalized Liouville measure $\mu_\theta$ on $\Ucal_{\text{phys}}(\theta)$. That is, for any measurable set $A \subset \Ucal_{\text{phys}}(\theta)$ and any $t \in \R$, we have $\mu_\theta(\Phi_t^\theta(A)) = \mu_\theta(A)$.
\end{proposition}

\begin{proofof}{Proposition \ref{prop:liouville_invariance}}
See \cref{app:proof_liouville}
\end{proofof}

\subsection{The Foliated Phase Space Structure}
The environmental parameter $\theta$ is a constant of motion, which imposes a rigid, foliated structure on the phase space.

\begin{definition}[The Global Microscopic Flow]
\label{def:global_micro_flow}
Let $\pi_\theta: \Efrak_{\text{phys}} \to \Thetacal$ be the canonical projection onto the environmental parameter space. For any $z = (x, y, v, \theta) \in \Efrak_{\text{phys}}$, the global microscopic flow $\Phi_t(z)$ is defined by
\begin{equation}
    \Phi_t(x, y, v, \theta) \coloneqq (x, \Phi_t^\theta(y, v), \theta),
\end{equation}
where $\Phi_t^\theta$ is the flow map on $\Ucal_{\text{phys}}(\theta)$ constructed in Definition~\ref{def:global_flow}.
\end{definition}

We now prove that the parameter $\theta$ is a first integral of this flow.

\begin{proposition}[$\theta$ as an Integral of Motion]
\label{prop:theta_integral}
The projection $\pi_\theta: \Efrak_{\text{phys}} \to \Thetacal$ is a first integral for the global microscopic flow $\Phi_t$.
\end{proposition}
\begin{proofof}{Proposition \ref{prop:theta_integral}}
The proof proceeds by a direct analysis of the construction of the global flow map $\Phi_t$, as given in Definition~\ref{def:global_flow}. Let an arbitrary initial state $z_0 = (x_0, y_0, v_0, \theta_0) \in \Efrak_{\text{phys}}$ be given. Let $z(t) = \Phi_t(z_0)$ denote the trajectory originating from $z_0$. We must show that the $\theta$-component of $z(t)$ remains $\theta_0$ for all time $t$. The trajectory is defined piecewise, and we analyze each piece of the dynamics separately.

\begin{enumerate}[label=\textbf{Step \arabic*:}, wide, labelindent=0pt]
    \item \textbf{Evolution During Interior Flow.}
    Let $(t_k, t_{k+1})$ be an interval between two consecutive collision events. For any $t \in (t_k, t_{k+1})$, the evolution of the state is governed by the interior flow $\Psi_t$ (Definition~\ref{def:flow_vector_field}). The full dynamics on $\Efrak_{\text{phys}}$ are described by the ordinary differential equation:
    \begin{equation*}
        \frac{d}{dt} z(t) = \frac{d}{dt} (x(t), y(t), v(t), \theta(t)) = (0, v(t), 0, 0).
    \end{equation*}
    The fourth component of this vector field equation is $\dot{\theta}(t) = 0$. Integrating this equation from the start of the interval, $t_k$, to any time $t \in (t_k, t_{k+1})$ yields:
    \begin{equation*}
        \theta(t) - \theta(t_k) = \int_{t_k}^t 0 \, ds = 0.
    \end{equation*}
    Thus, $\theta(t) = \theta(t_k)$ for all $t$ within any interval of free motion. Specifically, for the first interval starting at $t_0=0$, we have $\theta(t) = \theta_0$ for all $t \in [0, t_1)$, where $t_1$ is the first collision time.

    \item \textbf{Transformation at a Collision Event.}
    Let $t_k$ be the time of the $k$-th collision. Let $z^-(t_k) = \lim_{s \to t_k^-} z(s)$ denote the pre-collision state and $z(t_k)$ be the post-collision state. By the conclusion of Step 1, the $\theta$-component of the pre-collision state is $\pi_\theta(z^-(t_k)) = \theta_0$.

    According to the construction of the global flow (Definition~\ref{def:global_flow}), the post-collision state is determined by applying the specular reflection map $R_\theta$ to the microscopic components of the state. Crucially, the choice of which map to apply, $R_\theta$, is determined by the environmental parameter of the \emph{pre-collision state}. Let $\theta_{k}^- = \pi_\theta(z^-(t_k))$. The post-collision state $z(t_k)$ is given by:
    \begin{equation*}
        z(t_k) = (x_0, y_k^+, v_k^+, \theta_k^-), \quad \text{where } (y_k^+, v_k^+) = R_{\theta_k^-}(y_k^-, v_k^-).
    \end{equation*}
    By its very definition (Definition~\ref{def:reflection_map}), the map $R_\theta$ acts only on the velocity component of its argument, leaving the spatial component $y$ invariant. The map itself does not alter the parameter $\theta$ which defines it. Therefore, the $\theta$-component of the state is unchanged by the collision transformation:
    \begin{equation*}
        \pi_\theta(z(t_k)) = \theta_k^- = \pi_\theta(z^-(t_k)).
    \end{equation*}
\end{enumerate}

By induction, if the $\theta$-component is $\theta_0$ at the start of an interval of free motion (time $t_k$), it remains $\theta_0$ throughout that interval. At the subsequent collision (time $t_{k+1}$), the transformation rule preserves this value. Since the initial value at $t_0=0$ is $\theta_0$, the $\theta$-component of the state remains $\theta_0$ for all $t \ge 0$. The argument for negative time follows identically for the inverse flow. Thus, we have shown that $\pi_\theta(\Phi_t(z_0)) = \theta_0 = \pi_\theta(z_0)$ for all $t$, which completes the proof.
\end{proofof}

The existence of this integral of motion imposes a rigid structure on the phase space. The space $\Efrak_{\text{phys}}$ decomposes into a disjoint union of submanifolds on which the dynamics are confined.

\begin{theorem}[The Foliated Structure of the State Space]
\label{thm:foliation}
The physical state space $\Efrakphys$ is foliated by the level sets of the projection map $\pi_\theta$. The leaves of this foliation are the submanifolds-with-boundary
$\Efrak_{\mathrm{phys}}(\theta) \coloneqq \Mcal \times \Ucalphys(\theta)$, and each leaf is an invariant manifold for the global microscopic flow $\Phi_t$.
\end{theorem}

\begin{proofof}{\cref{thm:foliation}}
The proof is organized into two main steps. First, we establish the geometric claim that the level sets of $\pi_\theta$ form a smooth foliation. Second, we establish the dynamical claim that these leaves are invariant under the flow.

\begin{enumerate}[label=\textbf{Step \arabic*:}, wide, labelindent=0pt]

\item \textbf{Foliation Structure.} Our strategy is to demonstrate that the canonical projection $\pi_\theta: \Efrakphys \to \Thetacal$ is a smooth submersion. The desired foliation structure is then a direct consequence of the Submersion Level Set Theorem.

\begin{enumerate}[label=(\roman*), wide, labelindent=0pt]
    \item \textbf{Smoothness of the Domain and Map.}
    We first establish the geometric nature of the domain $\Efrakphys$ and the map $\pi_\theta$.
    \begin{enumerate}[label=(\alph*), wide]
        \item \textbf{The Domain.} The physical state space $\Efrakphys$ is defined in \cref{def:total_space} as the subset of the product manifold $\Efrak = \Mcal \times \T^k_y \times \Vcal \times \Thetacal$ where the smooth function $G(y, \theta)$ is non-negative. By \cref{ass:fundamental_axioms_unified}, the value $0$ is a regular value of $G$. The Preimage Theorem for manifolds-with-boundary guarantees that the set $\{ (x,y,v,\theta) \mid G(y,\theta) \ge 0 \}$ is a smooth, compact submanifold-with-boundary of the ambient manifold $\Efrak$.
        
        \item \textbf{The Map.} The map $\pi_\theta: \Efrakphys \to \Thetacal$ is the restriction of the canonical projection $P_\Thetacal: \Efrak \to \Thetacal$ to the submanifold $\Efrakphys$. Since canonical projections of a product manifold are smooth maps and the restriction of a smooth map to a smooth submanifold is itself smooth, we conclude that $\pi_\theta$ is a $C^\infty$-smooth map.
    \end{enumerate}

    \item \textbf{The Submersion Property.}
    We now prove that for every point $z \in \Efrakphys$, the differential of the map, $d(\pi_\theta)_z: T_z\Efrakphys \to T_{\pi_\theta(z)}\Thetacal$, is a surjective linear map. Let $z_0 = (x_0, y_0, v_0, \theta_0)$ be an arbitrary point in $\Efrakphys$. The tangent space to the ambient manifold $\Efrak$ at $z_0$ decomposes as:
    \begin{equation*}
        T_{z_0}\Efrak = T_{x_0}\Mcal \oplus T_{y_0}\T^k \oplus T_{v_0}\Vcal \oplus T_{\theta_0}\Thetacal.
    \end{equation*}
    An arbitrary tangent vector $W \in T_{z_0}\Efrak$ can be written as a tuple $W = (w_x, w_y, w_v, w_\theta)$. The differential of the unrestricted projection $P_\Thetacal$ simply extracts the final component: $d(P_\Thetacal)_{z_0}(W) = w_\theta$.

    We must show that this map remains surjective when restricted to the subspace $T_{z_0}\Efrakphys \subseteq T_{z_0}\Efrak$. We consider two cases for the location of $z_0$.

    \begin{enumerate}[label=(\alph*), wide]
        \item \textbf{$z_0$ is in the interior of $\Efrakphys$ ($G(y_0, \theta_0) > 0$).}
        In this case, $T_{z_0}\Efrakphys$ is identical to the ambient tangent space $T_{z_0}\Efrak$. To prove surjectivity, let $\eta$ be an arbitrary vector in the target space $T_{\theta_0}\Thetacal$. We must find a vector $W \in T_{z_0}\Efrakphys$ such that $d(\pi_\theta)_{z_0}(W) = \eta$. The vector $W = (0, 0, 0, \eta)$ is in $T_{z_0}\Efrakphys$ and its image under the differential is $d(\pi_\theta)_{z_0}(W) = \eta$. Thus, the map is surjective.

        \item \textbf{$z_0$ is on the boundary of $\Efrakphys$ ($G(y_0, \theta_0) = 0$).}
        Here, the tangent space $T_{z_0}\Efrakphys$ is the subspace of $W \in T_{z_0}\Efrak$ satisfying the constraint $dG_{z_0}(W) \ge 0$. In coordinates, this condition is:
        \begin{equation} \label{eq:tangent_boundary_condition}
            \scpr{\nabla_y G(y_0, \theta_0)}{w_y} + \scpr{\nabla_\theta G(y_0, \theta_0)}{w_\theta} \ge 0.
        \end{equation}
        To prove surjectivity, we must show that for any given $\eta \in T_{\theta_0}\Thetacal$, there exists \emph{some} vector $W = (w_x, w_y, w_v, w_\theta) \in T_{z_0}\Efrakphys$ such that $w_\theta = \eta$.
        
        Let an arbitrary $\eta \in T_{\theta_0}\Thetacal$ be given. We construct a suitable tangent vector $W$. Let $w_x = 0$, $w_v = 0$, and $w_\theta = \eta$. We must find a spatial component $w_y$ such that the full vector $W=(0, w_y, 0, \eta)$ satisfies the boundary condition \eqref{eq:tangent_boundary_condition}:
        \begin{equation*}
            \scpr{\nabla_y G(y_0, \theta_0)}{w_y} \ge - \scpr{\nabla_\theta G(y_0, \theta_0)}{\eta}.
        \end{equation*}
        By \cref{ass:fundamental_axioms_unified}, the gradient vector $\nabla_y G(y_0, \theta_0)$ is non-zero. Let us choose $w_y$ to be parallel to this gradient: $w_y = c \cdot \nabla_y G(y_0, \theta_0)$ for some scalar $c \in \R$. The inequality becomes:
        \begin{equation*}
            c \, \|\nabla_y G(y_0, \theta_0)\|^2 \ge - \scpr{\nabla_\theta G(y_0, \theta_0)}{\eta}.
        \end{equation*}
        Since $\|\nabla_y G(y_0, \theta_0)\|^2 > 0$, we can always choose a value of $c$ that satisfies this inequality. For instance, if the right-hand side is positive, we can choose any sufficiently large positive $c$. If the right-hand side is negative or zero, we can choose $c=0$.
        
        Thus, for any target vector $\eta$, we have explicitly constructed a valid tangent vector $W \in T_{z_0}\Efrakphys$ such that $d(\pi_\theta)_{z_0}(W) = \eta$. This confirms that the differential is surjective at every point on the boundary as well.
    \end{enumerate}
    Since the differential is surjective at every point $z \in \Efrakphys$, the map $\pi_\theta$ is a smooth submersion.

    \item \textbf{The Foliation Structure.}
    By the Submersion Level Set Theorem (also known as the Rank Theorem, applied to a map from a manifold-with-boundary to a boundaryless manifold), since $\pi_\theta: \Efrakphys \to \Thetacal$ is a smooth submersion, for each $\theta \in \Thetacal$, the preimage $\Efrak_{\mathrm{phys}}(\theta) = \pi_\theta^{-1}(\{\theta\})$ is a properly embedded submanifold-with-boundary of $\Efrakphys$. The collection of these disjoint submanifolds, $\{ \Efrak_{\mathrm{phys}}(\theta) \}_{\theta \in \Thetacal}$, forms a partition of the total space, thus defining a smooth foliation of $\Efrakphys$.
\end{enumerate}

\item \textbf{Invariance of Leaves under the Flow.} We now prove that the dynamics defined by the flow $\Phi_t$ are confined to the leaves of the foliation established in Step 1. This requires showing that for any leaf, the flow map maps the leaf to itself.

\begin{enumerate}[label=(\roman*), wide, labelindent=0pt]
    \item \textbf{Setup.} Let $\theta_0 \in \Thetacal$ be fixed, which specifies the leaf $\Efrak_{\mathrm{phys}}(\theta_0)$. Let $z_0$ be an arbitrary point within this leaf. By the definition of the leaf as a level set, this means:
    \begin{equation*}
        \pi_\theta(z_0) = \theta_0.
    \end{equation*}

    \item \textbf{Application of the Integral of Motion.}
    We now invoke \cref{prop:theta_integral}, which establishes that the function $\pi_\theta$ is a first integral (i.e., a constant of motion) for the global microscopic flow $\Phi_t$ defined in \cref{def:global_micro_flow}. This means that for any initial state $z$ and any time $t \in \R$, the following identity holds:
    \begin{equation*}
        \pi_\theta(\Phi_t(z)) = \pi_\theta(z).
    \end{equation*}

    \item \textbf{Dynamical Confinement.}
    We apply the result from \cref{prop:theta_integral} to our chosen initial point $z_0$:
    \begin{equation*}
        \pi_\theta(\Phi_t(z_0)) = \pi_\theta(z_0).
    \end{equation*}
    Combining this with the fact from step (i) that $\pi_\theta(z_0) = \theta_0$, we immediately have:
    \begin{equation*}
        \pi_\theta(\Phi_t(z_0)) = \theta_0 \quad \text{for all } t \in \R.
    \end{equation*}
    This equation signifies that for any time $t$, the point $\Phi_t(z_0)$ lies in the level set of $\pi_\theta$ corresponding to the value $\theta_0$. By definition, this level set is precisely the leaf $\Efrak_{\mathrm{phys}}(\theta_0)$.
    Therefore, the entire trajectory $\{\Phi_t(z_0)\}_{t \in \R}$ starting at $z_0$ remains within the leaf $\Efrak_{\mathrm{phys}}(\theta_0)$.

    \item \textbf{The Invariant Manifold.}
    Since the choice of the initial point $z_0$ within the leaf was arbitrary, we conclude that the entire leaf is mapped to itself by the flow. As the choice of the leaf (via $\theta_0$) was also arbitrary, it follows that for every $\theta \in \Thetacal$, the leaf $\Efrak_{\mathrm{phys}}(\theta)$ is an invariant manifold for the flow $\Phi_t$. This completes the proof of the theorem.
\end{enumerate}
\end{enumerate}
\end{proofof}

\begin{remark}[Central Organizing Principle]
\label{rem:foliation_principle}
This foliated structure is the central organizing principle of our analysis. It provides a rigorous justification for decomposing the complex global problem into two distinct stages.
\begin{enumerate}[label=(\roman*), wide, labelindent=0pt]
    \item \textbf{Intra-Leaf Analysis.} First, we must study the dynamics restricted to an arbitrary, single leaf $\Efrak_{\text{phys}}(\theta)$. This involves understanding the ergodic and statistical properties of the flow $\Phi_t^\theta$ for a fixed $\theta$. This is the objective of the Nulla Scale analysis (see \cref{sec:nulla_scale}).
    \item \textbf{Inter-Leaf Synthesis.} Second, we must synthesize the results from all the leaves to construct the global macroscopic theory. The non-linearity of our problem enters at this stage, as the choice of the relevant leaf, $\theta$, will depend on the macroscopic state $(x, \nabla u(x))$.
\end{enumerate}
This structure allows us to treat the chaotic dynamics and the non-linear feedback as two separate, sequential problems.
\end{remark}

\subsection{The Infinitesimal Generator}
We now translate the description of the dynamics from the language of flows to operators on the Hilbert space $\Hcal_\theta = L^2(\Ucal_{\text{phys}}(\theta), \mu_\theta)$.

\begin{definition}[The Koopman Group]
\label{def:koopman_group}
For a fixed $\theta \in \Thetacal$, the flow $\Phi_t^\theta$ induces a family of operators $(P_t^\theta)_{t \in \R}$ on the Hilbert space $\Hcal_\theta$, known as the Koopman group, which describes the evolution of observables:
\begin{equation}
    (P_t^\theta \psi)(z) \coloneqq \psi(\Phi_{-t}^\theta(z)), \quad \text{for } \psi \in \Hcal_\theta, z \in \Ucal_{\text{phys}}(\theta).
\end{equation}
We use $\Phi_{-t}^\theta$ by convention so that the generator is $d/dt|_{t=0} P_t^\theta$.
\end{definition}

\begin{proposition}[Properties of the Koopman Group]
\label{prop:koopman_properties}
For each $\theta \in \Thetacal$, the family $(P_t^\theta)_{t \in \R}$ defined in \cref{def:koopman_group} is a strongly continuous one-parameter group of unitary operators on the Hilbert space $\Hcaltheta$.
\end{proposition}

\begin{proofof}{\cref{prop:koopman_properties}}
The proof is established in three steps. We verify (i) the one-parameter group property, (ii) the unitarity of each operator, and (iii) the strong continuity of the group.

\begin{enumerate}[label=\textbf{Step \arabic*:}, wide, labelindent=0pt]

\item \textbf{One-Parameter Group Property.} Let $\psi \in \Hcaltheta$ and $z \in \Ucalphys(\theta)$. We verify the group axioms directly from the properties of the flow map $\Phi_t^\theta$.
\begin{enumerate}[label=(\roman*), wide, labelindent=0pt]
    \item \textbf{Identity.} For $t=0$, the flow is the identity map, $\Phi_0^\theta(z) = z$. Thus,
    \begin{equation*}
        (P_0^\theta \psi)(z) = \psi(\Phi_0^\theta(z)) = \psi(z),
    \end{equation*}
    which implies $P_0^\theta = \mathrm{Id}$.

    \item \textbf{Composition.} For any $s, t \in \R$, the flow property is $\Phi_{-(s+t)}^\theta = \Phi_{-t}^\theta \circ \Phi_{-s}^\theta$. We apply this to the composition of operators:
    \begin{align*}
        (P_s^\theta P_t^\theta \psi)(z) &= P_s^\theta(\psi \circ \Phi_{-t}^\theta)(z) \\
        &= (\psi \circ \Phi_{-t}^\theta)(\Phi_{-s}^\theta(z)) \\
        &= \psi(\Phi_{-t}^\theta(\Phi_{-s}^\theta(z))) \\
        &= \psi((\Phi_{-t}^\theta \circ \Phi_{-s}^\theta)(z)) \\
        &= \psi(\Phi_{-(s+t)}^\theta(z)) \\
        &= (P_{s+t}^\theta \psi)(z).
    \end{align*}
    Since this holds for all $\psi$ and $z$, we have $P_s^\theta P_t^\theta = P_{s+t}^\theta$.
\end{enumerate}
These two properties confirm that $(P_t^\theta)_{t \in \R}$ is a one-parameter group of operators.

\item \textbf{Unitarity.}
An operator on a Hilbert space is unitary if and only if it is a surjective isometry.

\begin{enumerate}[label=(\roman*), wide, labelindent=0pt]
    \item \textbf{Isometry.} We show that $P_t^\theta$ preserves the inner product for any $t \in \R$. Let $\psi, \phi \in \Hcaltheta$.
    \begin{align*}
        \scpr{P_t^\theta \psi}{P_t^\theta \phi}_\theta &= \int_{\Ucalphys(\theta)} \overline{(P_t^\theta \psi)(z)} (P_t^\theta \phi)(z) \, d\mu_\theta(z) \\
        &= \int_{\Ucalphys(\theta)} \overline{\psi(\Phi_{-t}^\theta(z))} \phi(\Phi_{-t}^\theta(z)) \, d\mu_\theta(z).
    \end{align*}
    We perform a change of variables $w = \Phi_{-t}^\theta(z)$, which implies $z = \Phi_t^\theta(w)$. The Jacobian of this transformation is unity because the flow preserves the Liouville measure $\mu_\theta$, as established in \cref{prop:liouville_invariance}. Therefore, $d\mu_\theta(z) = d\mu_\theta(w)$. Substituting this into the integral gives:
    \begin{align*}
        \scpr{P_t^\theta \psi}{P_t^\theta \phi}_\theta &= \int_{\Ucalphys(\theta)} \overline{\psi(w)} \phi(w) \, d\mu_\theta(w) \\
        &= \scpr{\psi}{\phi}_\theta.
    \end{align*}
    This shows that $P_t^\theta$ is an isometry for all $t \in \R$.

    \item \textbf{Surjectivity.} To prove surjectivity, we must show that for any $\phi \in \Hcaltheta$, there exists a $\psi \in \Hcaltheta$ such that $P_t^\theta \psi = \phi$. Let us choose the candidate $\psi \coloneqq P_{-t}^\theta \phi$. Since $P_{-t}^\theta$ is a well-defined operator on $\Hcaltheta$, $\psi$ is in $\Hcaltheta$. Applying $P_t^\theta$ to this candidate yields:
    \begin{equation*}
        P_t^\theta \psi = P_t^\theta (P_{-t}^\theta \phi) = P_{t-t}^\theta \phi = P_0^\theta \phi = \phi.
    \end{equation*}
    Thus, every element $\phi$ in the codomain has a pre-image, and the operator is surjective.
\end{enumerate}
Since $P_t^\theta$ is a surjective isometry, it is a unitary operator.

\item \textbf{Strong Continuity.}
We must show that for any $\psi \in \Hcaltheta$, $\lim_{t \to 0} \|P_t^\theta \psi - \psi\|_{\theta} = 0$. We use a standard density argument.

\begin{enumerate}[label=(\roman*), wide, labelindent=0pt]
    \item \textbf{Continuity on a dense subset.}
    Let $C(\Ucalphys(\theta))$ be the space of continuous complex-valued functions on the compact space $\Ucalphys(\theta)$. This space is dense in $\Hcaltheta = L^2(\Ucalphys(\theta), d\mu_\theta)$. Let $f \in C(\Ucalphys(\theta))$. Since $\Ucalphys(\theta)$ is compact, $f$ is uniformly continuous. The flow map $\Phi_t^\theta(z)$ is continuous in $(t,z)$ for $z$ in the interior of the domain. For any given $z$, $\Phi_t^\theta(z)$ is piecewise continuous in $t$, with discontinuities at collision times. However, for any fixed $z$, $\lim_{t \to 0} \Phi_t^\theta(z) = z$. This implies pointwise convergence of the integrand: for almost every $z$,
    \begin{equation*}
        \lim_{t \to 0} (P_t^\theta f)(z) = \lim_{t \to 0} f(\Phi_{-t}^\theta(z)) = f(z),
    \end{equation*}
    where the limit follows from the continuity of $f$. We consider the squared norm of the difference:
    \begin{equation*}
        \|P_t^\theta f - f\|_\theta^2 = \int_{\Ucalphys(\theta)} |f(\Phi_{-t}^\theta(z)) - f(z)|^2 \, d\mu_\theta(z).
    \end{equation*}
    The integrand converges pointwise to 0 for almost every $z$. Furthermore, since $f$ is continuous on a compact set, it is bounded by some constant $M$, so $|f(\Phi_{-t}^\theta(z)) - f(z)|^2 \le (2M)^2$. Since the total measure $\mu_\theta(\Ucalphys(\theta)) = 1$, the constant function $(2M)^2$ is integrable. By the Lebesgue Dominated Convergence Theorem, we can interchange the limit and the integral:
    \begin{equation*}
        \lim_{t \to 0} \|P_t^\theta f - f\|_\theta^2 = \int_{\Ucalphys(\theta)} \lim_{t \to 0} |f(\Phi_{-t}^\theta(z)) - f(z)|^2 \, d\mu_\theta(z) = \int 0 \, d\mu_\theta = 0.
    \end{equation*}
    Thus, the group is strongly continuous for all functions in the dense subset $C(\Ucalphys(\theta))$.

    \item \textbf{Extension to the full Hilbert space.}
    Let $\psi \in \Hcaltheta$ be arbitrary and let $\epsilon > 0$. Since $C(\Ucalphys(\theta))$ is dense in $\Hcaltheta$, there exists a function $f \in C(\Ucalphys(\theta))$ such that
    \begin{equation*}
        \|\psi - f\|_\theta < \frac{\epsilon}{3}.
    \end{equation*}
    We use the triangle inequality on $\|P_t^\theta \psi - \psi\|_\theta$:
    \begin{equation*}
        \|P_t^\theta \psi - \psi\|_\theta \le \|P_t^\theta \psi - P_t^\theta f\|_\theta + \|P_t^\theta f - f\|_\theta + \|f - \psi\|_\theta.
    \end{equation*}
    We bound each of the three terms:
    \begin{enumerate}[label=(\alph*), wide]
        \item Since $P_t^\theta$ is unitary, it is an isometry. Thus,
        \begin{equation*}
            \|P_t^\theta \psi - P_t^\theta f\|_\theta = \|P_t^\theta (\psi - f)\|_\theta = \|\psi - f\|_\theta < \frac{\epsilon}{3}.
        \end{equation*}
        \item From Step 1, since $f$ is continuous, we know that $\lim_{t \to 0} \|P_t^\theta f - f\|_\theta = 0$. Therefore, there exists a $\delta > 0$ such that for all $|t| < \delta$,
        \begin{equation*}
             \|P_t^\theta f - f\|_\theta < \frac{\epsilon}{3}.
        \end{equation*}
        \item By our choice of $f$, we have $\|f - \psi\|_\theta < \frac{\epsilon}{3}$.
    \end{enumerate}
    Combining these bounds, for any $|t| < \delta$, we have
    \begin{equation*}
        \|P_t^\theta \psi - \psi\|_\theta < \frac{\epsilon}{3} + \frac{\epsilon}{3} + \frac{\epsilon}{3} = \epsilon.
    \end{equation*}
    This holds for any arbitrary $\psi \in \Hcaltheta$, which is the definition of strong continuity at $t=0$. Continuity for any other $t_0$ follows from the group property. This completes the proof.
\end{enumerate}
\end{enumerate}
\end{proofof}

The existence of a strongly continuous unitary group allows us to define its generator via Stone's theorem, which provides the rigorous foundation for our analysis.

\begin{definition}[The Infinitesimal Generator]
\label{def:generator}
By Stone's theorem on one-parameter unitary groups, there exists a unique (possibly unbounded) self-adjoint operator $A$ such that $P_t^\theta = e^{itA}$. We define the infinitesimal generator of the flow, denoted $\Lcal_{\text{fast}}(\theta)$, as the operator $\Lcal_{\text{fast}}(\theta) \coloneqq iA$. It follows that $\Lcal_{\text{fast}}(\theta)$ is a densely defined, closed, skew-adjoint operator on $\Hcal_\theta$. Its domain is the set of $\psi \in \Hcal_\theta$ for which the strong limit
\begin{equation}
    \Lcal_{\text{fast}}(\theta) \psi = \lim_{t \to 0} \frac{P_t^\theta \psi - \psi}{t}
\end{equation}
exists in $\Hcal_\theta$.
\end{definition}

We now provide an explicit characterization of this abstractly defined operator.

\begin{theorem}[Characterization of the Generator]
\label{thm:generator_characterization}
Let the system satisfy the standing assumptions of this paper. The infinitesimal generator $\Lcal_{\mathrm{fast}}(\theta)$ of the strongly continuous one-parameter unitary group $(P_t^\theta)_{t \in \R}$ is the unique, densely defined, skew-adjoint operator characterized as follows:
\begin{enumerate}[label=(\roman*), wide, labelindent=0pt]
    \item The domain $\Dcal(\Lcal_{\mathrm{fast}}(\theta))$ is the subspace of functions $\psi \in \Hcal_\theta$ such that the weak derivative $v \cdot \nabla_y \psi$ exists as an element of $\Hcal_\theta$ and the trace of $\psi$ on the boundary satisfies the specular reflection condition $\psi(z) = \psi(R_\theta(z))$ for almost every $z \in \partial_-\Ucal_{\mathrm{phys}}(\theta)$.
    \item For any $\psi \in \Dcal(\Lcal_{\mathrm{fast}}(\theta))$, the action of the generator is given by the transport operator:
    \begin{equation}
        (\Lcal_{\mathrm{fast}}(\theta)\psi)(z) = (v \cdot \nabla_y \psi)(z).
    \end{equation}
\end{enumerate}
Furthermore, for any complex number $\lambda$ with $\mathrm{Re}(\lambda) > 0$, the resolvent $(\lambda I - \Lcal_{\mathrm{fast}}(\theta))^{-1}$ exists as a bounded operator on all of $\Hcal_\theta$ and is given by the Laplace transform of the semigroup:
\begin{equation}
    (\lambda I - \Lcal_{\mathrm{fast}}(\theta))^{-1}f = \int_0^\infty e^{-\lambda t} (P_t^\theta f) \, dt.
\end{equation}
\end{theorem}

\begin{proofof}{\cref{thm:generator_characterization}}
See \cref{app:proof_generator}
\end{proofof}

\begin{remark}[The Analytical Obstacle]
\label{rem:analytical_obstacle}
The structure of the generator $\Lcal_{\text{fast}}(\theta)$ is the primary reason that standard homogenization techniques, developed primarily for second-order elliptic or parabolic operators, are insufficient.
\begin{enumerate}[label=(\roman*), wide, labelindent=0pt]
    \item \textbf{Hyperbolic vs. Elliptic.} $\Lcal_{\text{fast}}$ is a first-order hyperbolic operator. Unlike an elliptic operator (e.g., the Laplacian, $\Delta$), its inverse is not a compact operator on $\Hcal_\theta$. The resolvent $(\Lcal_{\text{fast}} - \lambda I)^{-1}$ is bounded but does not have the smoothing (regularizing) properties of an elliptic resolvent.
    \item \textbf{Lack of a Spectral Gap.} The spectrum of $\Lcal_{\text{fast}}(\theta)$ lies purely on the imaginary axis. On the space of zero-mean functions $\Hcal_{\theta,0}$, there is no spectral gap around zero that would ensure the equation $\Lcal_{\text{fast}} \chi = f$ has a well-behaved, bounded inverse. In fact, for any $f \in \Hcal_{\theta,0}$, the equation may have no solution in $\Dcal(\Lcal_{\text{fast}})$ or infinitely many.
    \item \textbf{Failure of Variational Methods.} Standard variational methods (like the Lax-Milgram theorem) for solving equations like $\Lcal \chi = f$ rely on the coercivity of the associated bilinear form $\scpr{\Lcal \psi}{\psi}$. For our skew-adjoint operator, this form is $\scpr{\Lcal \psi}{\psi} = -\scpr{\psi}{\Lcal \psi} = -\overline{\scpr{\Lcal \psi}{\psi}}$, which implies its real part is zero. The operator is not coercive, and these methods fail.
\end{enumerate}
This analytical breakdown necessitates a different approach. We cannot simply invert the generator using standard PDE tools. Our geometric homogenization framework is designed specifically to construct a well-defined solution to the cell problem for this class of non-coercive, hyperbolic generators.
\end{remark}

\section{Uniform Hyperbolicity of the Microscopic Dynamics}
\label{sec:hyperbolicity_fiber}

The analysis in Section~\ref{sec:micro_framework} established the well-posedness of the continuous-time flow from a set of first-principles geometric axioms. The objective of this section is to prove that these same axioms endow the system with the strong chaotic properties required for the emergence of robust statistical laws. Our central goal is to demonstrate that the family of microscopic systems, parameterized by $\theta \in \Thetacal$, is uniformly hyperbolic in the sense of Anosov.

Our strategy is constructive. We will first reduce the continuous-time flow to a discrete-time global billiard map $\Pcal$ on a single compact manifold $\Sigma$ that fiberwise contains all the individual dynamics. We then provide a rigorous proof that this global map is a $C^\infty$-diffeomorphism, making explicit that the proof relies only on the geometric axioms from Section \ref{sec:assumptions}. With the smoothness of the map established, we will proceed to the main result of this section: a constructive proof that $\Pcal$ is a uniform Anosov diffeomorphism. This is achieved by explicitly building a continuous, strictly invariant unstable cone field for its linearization, $d\Pcal$. The uniformity of the hyperbolicity constants for the entire family $(\Pcal_\theta)$ will then follow as a direct corollary.

\subsection{The Billiard Map}
\label{sec:reduction_to_discrete}

The central objective of this section is to prove that the microscopic dynamics are chaotic in a strong, uniform sense. To do this, we employ the standard technique of reducing the continuous-time flow to a discrete-time map, known as the billiard map or Poincaré map, on a well-chosen surface of section. This reduction allows for the application of the well-developed theory of hyperbolic diffeomorphisms. This subsection is dedicated to the rigorous construction of this map and its natural invariant measure.

The natural choice for a Poincaré section is the set of all states on the collision boundary, as the dynamics there are non-trivial. For this to be a valid construction, the flow must be everywhere transverse to the section. As established in Proposition~\ref{prop:non_grazing}, our geometric assumptions uniformly exclude grazing collisions, thereby guaranteeing the required transversality.

\begin{definition}[The Billiard Phase Space (Poincaré Section)]
\label{def:billiard_phase_space}
For each fixed $\theta \in \Thetacal$, the phase space for the discrete-time dynamics is the manifold of incoming states on the collision boundary:
\begin{equation}
    \Sigma_\theta \coloneqq \partial_-\Ucal_{\mathrm{phys}}(\theta) = \{ (y,v) \in \T^k \times \Vcal \mid G(y,\theta) = 0 \text{ and } \scpr{v}{\mathbf{n}(y, \theta)} < 0 \}.
\end{equation}
As an open submanifold of the total boundary phase space $\partial\Ucal_{\mathrm{phys}}(\theta)$, $\Sigma_\theta$ is a smooth manifold of dimension $2k-1$.
\end{definition}

The Liouville measure $\mu_\theta$ on the full phase space induces a natural measure on this section, which corresponds to the flux of phase space volume through the section.

\begin{definition}[Induced Birkhoff Measure]
\label{def:birkhoff_measure}
The measure $\nu_\theta$ induced on the section $\Sigma_\theta$ by the Liouville measure $\mu_\theta = C_\theta \, dy \, dv$ is the Birkhoff measure, defined by the volume element:
\begin{equation}
    d\nu_\theta(y,v) \coloneqq C_\theta \lvert\scpr{v}{\mathbf{n}(y, \theta)}\rvert dS(y) \, dV(v),
\end{equation}
where $dS(y)$ is the canonical surface area element on the boundary manifold $\partial\Ocal(\theta)$ and $dV(v)$ is the Lebesgue measure on the velocity space $\Vcal$.
\end{definition}

\begin{remark}[Physical and Mathematical Motivation for the Birkhoff Measure]
\label{rem:birkhoff_measure}
This is the physically correct definition of the invariant measure for a billiard map. The term $\lvert\scpr{v}{\mathbf{n}}\rvert$ represents the magnitude of the velocity component normal to the boundary, meaning the measure element corresponds to the flux of phase space volume through the surface element $dS(y)$. Beyond this physical interpretation, this specific form is mathematically indispensable for the invariance of the measure under the full billiard map $\Pcal_\theta$. While the reflection map $\Rcal_\theta$ is an isometry with respect to this measure, the free-flow map $\Fcal_\theta$ is not. As rigorously proven in Appendix \ref{app:proof_measure_invariance_rigorous}, the Jacobian determinant of the spatial component of the free-flow map is precisely the ratio of the normal velocity components at the end and start of a trajectory. The inclusion of the $\lvert\scpr{v}{\mathbf{n}}\rvert$ factor in the measure definition is therefore the exact term required to ensure a perfect cancellation with this non-trivial Jacobian, rendering the measure invariant.
\end{remark}

The dynamics on this section are governed by the first-return map, which we construct as a composition of reflection and free motion. Let $\SigmaMan'_\theta \coloneqq \partial_+\Ucal_{\mathrm{phys}}(\theta)$ be the manifold of outgoing states.

\begin{definition}[The Billiard Map]
\label{def:billiard_map}
The billiard map (or Poincaré map) $\Pcal_\theta: \Sigma_\theta \to \Sigma_\theta$ is defined as the composition $\Pcal_\theta \coloneqq \Fcal_\theta \circ \Rcal_\theta$. The constituent maps act in sequence as:
\begin{equation*}
    \Sigma_\theta \xrightarrow{\quad \Rcal_\theta \quad} \Sigma'_\theta \xrightarrow{\quad \Fcal_\theta \quad} \Sigma_\theta
\end{equation*}
\begin{enumerate}[label=(\roman*), wide, labelindent=0pt]
    \item \textbf{The reflection map} $\Rcal_\theta: \Sigma_\theta \to \Sigma'_\theta$ is the specular reflection defined in \cref{def:reflection_map}.
    \item \textbf{The free-flow map} $\Fcal_\theta: \Sigma'_\theta \to \Sigma_\theta$ evolves a post-collision state $(y^+, v^+)$ to the subsequent pre-collision state $(y\dnext, v\dnext)$ via:
    \begin{align*}
        v\dnext &= v^+, \\
        y\dnext &= y^+ + \tau(y^+, v^+, \theta) v^+,
    \end{align*}
    where $\tau(y^+, v^+, \theta) \coloneqq \inf \{ t > 0 \mid y^+ + t v^+ \in \dOcal(\theta) \}$ is the flight time. By \cref{lem:smooth_flight_time}, $\tau$ is a $C^\infty$ function, ensuring that $\Fcal_\theta$ is a $C^\infty$-diffeomorphism.
\end{enumerate}
\end{definition}

A fundamental property of Hamiltonian systems, which we now prove, is that this reduction to a Poincaré section preserves the measure structure. This result is the foundation of the subsequent ergodic analysis.

\begin{proposition}[Invariance of the Birkhoff Measure]
\label{prop:birkhoff_invariance}
For each $\theta \in \Thetacal$, the billiard map $\Pcal_\theta$ preserves the induced Birkhoff measure $\nu_\theta$.
\end{proposition}

\begin{proofof}{Proposition \ref{prop:birkhoff_invariance}}
See \cref{app:proof_measure_invariance_rigorous}.
\end{proofof}

With the discrete dynamical system $(\Sigma_\theta, \Pcal_\theta, \nu_\theta)$ and its fundamental conservation law rigorously established, we are now prepared to analyze its chaotic properties. The remainder of this section is dedicated to proving that the family of maps $(\Pcal_\theta)_{\theta \in \Thetacal}$ is uniformly Anosov.

\subsection{The Global Billiard Map and its Smoothness}
\label{sec:billiard}

To prove uniform hyperbolicity for the family of systems, we construct a single, global dynamical system on a compact manifold that encompasses all the fiber dynamics. This translates the question of uniformity in $\theta$ into proving that a single smooth diffeomorphism on a compact space is Anosov.

\begin{definition}[The Global Collision Manifold]
\label{def:global_collision_manifold}
The total pre-collision space is the fiber bundle over $\Thetacal$ defined as:
\begin{equation}
    \Sigma \coloneqq \bigsqcup_{\theta \in \Thetacal} \{\theta\} \times \partial_-\Ucal_{\mathrm{phys}}(\theta) = \{ (z, \theta) \mid \theta \in \Thetacal, z=(y,v), y \in \dOcal(\theta), \scpr{v}{\mathbf{n}(y,\theta)} < 0 \}.
\end{equation}
Similarly, the total post-collision space is $\Sigma' \coloneqq \bigsqcup_{\theta \in \Thetacal} \{\theta\} \times \partial_+\Ucal_{\mathrm{phys}}(\theta)$.
\end{definition}

\begin{definition}[Global Billiard Map]
\label{def:global_billiard_map}
The global billiard map $\Pcal: \Sigma \to \Sigma$ maps a pre-collision state to the subsequent pre-collision state. It is the composition $\Pcal \coloneqq \Fcal \circ \Rcal$, where:
\begin{enumerate}[label=(\roman*), wide, labelindent=0pt]
    \item The \textbf{global reflection map} $\Rcal: \Sigma \to \Sigma'$ acts as $\Rcal(z, \theta) = (R_\theta(z), \theta)$.
    \item The \textbf{global free-flow map} $\Fcal: \Sigma' \to \Sigma$ acts as $\Fcal(z^+, \theta) = (\Psi_{\tau(z^+,\theta)}(z^+), \theta)$, where $\tau(z^+,\theta)$ is the time of flight from the post-collision state $(z^+,\theta)$ to the next collision.
\end{enumerate}
\end{definition}

\begin{proposition}[Geometric Properties of $\Sigma$]
\label{prop:sigma_properties}
The total pre-collision space $\Sigma$, together with its grazing boundary, forms a smooth, compact manifold-with-boundary.
\end{proposition}

\begin{proofof}{Proposition \ref{prop:sigma_properties}}
The proof is organized into two main parts. First, we establish the geometric claim that the level sets defining $\Sigma$ form a smooth manifold-with-boundary by a sequential application of the Submersion Level Set Theorem. Second, we establish the topological claim of compactness by showing that $\Sigma$ is a closed subset of a known compact space.

\begin{enumerate}[label=\textbf{Step \arabic*:}, wide, labelindent=0pt]

\item \textbf{Smooth Manifold-with-Boundary Structure.}
Our strategy is to construct $\Sigma$ in two stages. We first define the manifold of all possible boundary points in the extended phase space, and then restrict to this manifold to define the subset of incoming states.

\begin{enumerate}[label=(\roman*), wide, labelindent=0pt]
    \item \textbf{The Total Boundary Manifold $\mathcal{B}$.}
    Let the full ambient space be the product manifold $\mathcal{A} \coloneqq \T^k_y \times \Vcal \times \Thetacal$, which is a smooth manifold (as $\Vcal$ is a manifold with corners, but this detail does not affect the interior analysis). We define the total boundary manifold $\mathcal{B}$ as the subset of $\mathcal{A}$ where a collision can occur:
    \begin{equation*}
        \mathcal{B} \coloneqq \{ (y,v,\theta) \in \mathcal{A} \mid G(y,\theta) = 0 \}.
    \end{equation*}
    We now show that $\mathcal{B}$ is a smooth, properly embedded submanifold of $\mathcal{A}$. We apply the Regular Level Set Theorem to the obstacle-defining function $G$, viewed as a map $G: \mathcal{A} \to \R$. By assumption, $G$ is a $C^\infty$ map. A point $c \in \R$ is a regular value of $G$ if for every point $p \in G^{-1}(\{c\})$, the differential $dG_p: T_p\mathcal{A} \to T_c\R$ is surjective.
    
    Let $p_0=(y_0,v_0,\theta_0)$ be an arbitrary point in the level set $\mathcal{B} = G^{-1}(\{0\})$. In local coordinates, the differential $dG_{p_0}$ is represented by the Jacobian matrix, which for a scalar-valued function is the transpose of the gradient vector: $(\nabla_y G, \nabla_v G, \nabla_\theta G)$. Since $G$ is independent of $v$, $\nabla_v G = 0$. The differential is surjective if and only if the gradient vector is non-zero.
    
    By \cref{ass:fundamental_axioms_unified}, for any $(y,\theta)$ with $G(y,\theta)=0$, we have $\nabla_y G(y,\theta) \neq \mathbf{0}$. This guarantees that the gradient of $G$ is non-zero at every point $p_0 \in \mathcal{B}$. Therefore, the differential $dG_{p_0}$ is surjective for all $p_0 \in \mathcal{B}$, which means $0$ is a regular value of $G$. The Regular Level Set Theorem allows us to conclude that $\mathcal{B}$ is a smooth, properly embedded submanifold of $\mathcal{A}$.

    \item \textbf{Defining the Incoming Submanifold.}
    The pre-collision space $\Sigma$ and its boundary are subsets of the manifold $\mathcal{B}$ determined by the sign of the normal velocity component. We define a function $f: \mathcal{B} \to \R$ by:
    \begin{equation*}
        f(y,v,\theta) \coloneqq \scpr{v}{\mathbf{n}(y,\theta)} = \frac{\scpr{v}{\nabla_y G(y,\theta)}}{\|\nabla_y G(y,\theta)\|}.
    \end{equation*}
    Since $G$ is $C^\infty$ and \cref{ass:fundamental_axioms_unified} ensures the denominator never vanishes on $\mathcal{B}$, the unit normal vector field $\mathbf{n}(y,\theta)$ is a smooth map on $\mathcal{B}$. As the inner product is a smooth operation, the function $f$ is of class $C^\infty$ on $\mathcal{B}$.
    By definition, the total pre-collision space (including grazing collisions) is the preimage of the interval $(-\infty, 0]$:
    \begin{equation*}
        \Sigma \cup \partial^0\Sigma = f^{-1}((-\infty, 0]) = \{ p \in \mathcal{B} \mid f(p) \le 0 \}.
    \end{equation*}

    \item \textbf{Verification of Regularity.}
    We now apply the Submersion Level Set Theorem for manifolds-with-boundary. This requires showing that $0$ is a regular value of the function $f: \mathcal{B} \to \R$. Let $p_0=(y_0,v_0,\theta_0)$ be an arbitrary point in the level set $f^{-1}(\{0\})$, which corresponds to a grazing collision. We must show that the differential $df_{p_0}: T_{p_0}\mathcal{B} \to T_0\R$ is surjective. A linear map to $\R$ is surjective if and only if it is not the zero map. We need only find a single tangent vector $\mathbf{w} \in T_{p_0}\mathcal{B}$ such that $df_{p_0}(\mathbf{w}) \neq 0$.
    
    Consider a variation purely in the velocity component. Let $\mathbf{w} = (0, \mathbf{n}(y_0,\theta_0), 0)$. This vector lies in the tangent space $T_{p_0}\mathcal{A}$. To be in $T_{p_0}\mathcal{B}$, it must lie in the kernel of $dG_{p_0}$. The condition is $\scpr{\nabla_y G}{\delta y} + \scpr{\nabla_\theta G}{\delta \theta} = 0$. For our choice of $\mathbf{w}$, we have $\delta y=0$ and $\delta\theta=0$, so the condition is trivially satisfied. Thus, $\mathbf{w}$ is a valid tangent vector in $T_{p_0}\mathcal{B}$. We compute the action of the differential on this vector, which is the directional derivative of $f$ along the path $p(\epsilon) = (y_0, v_0 + \epsilon\mathbf{n}, \theta_0)$:
    \begin{align*}
        df_{p_0}(\mathbf{w}) &= \left. \frac{d}{d\epsilon} f(p(\epsilon)) \right|_{\epsilon=0} = \left. \frac{d}{d\epsilon} \scpr{v_0 + \epsilon\mathbf{n}}{\mathbf{n}(y_0,\theta_0)} \right|_{\epsilon=0} \\
        &= \scpr{v_0}{\mathbf{n}} + \left. \frac{d}{d\epsilon}(\epsilon\scpr{\mathbf{n}}{\mathbf{n}}) \right|_{\epsilon=0} = f(p_0) + 1.
    \end{align*}
    Since $p_0$ is in the level set of $0$, $f(p_0)=0$. The derivative is therefore $1$, which is non-zero. Since the differential is non-zero, it is surjective. Therefore, $0$ is a regular value of $f$.

    \item \textbf{The Manifold Structure.}
    Since $f: \mathcal{B} \to \R$ is a $C^\infty$ function and $0$ is a regular value, the Submersion Level Set Theorem for manifolds-with-boundary guarantees that the set $\Sigma \cup \partial^0\Sigma = f^{-1}((-\infty, 0])$ is a smooth, properly embedded submanifold-with-boundary of $\mathcal{B}$. Its interior is $\Sigma = f^{-1}((-\infty, 0))$ and its boundary is $\partial^0\Sigma = f^{-1}(\{0\})$.
\end{enumerate}

\item \textbf{Compactness.}
Our strategy is to show that $\Sigma \cup \partial^0\Sigma$ is a closed subset of a compact space, which implies its own compactness.

\begin{enumerate}[label=(\roman*), wide, labelindent=0pt]
    \item \textbf{The Ambient Compact Space.}
    The total physical state space $\Efrak_{\text{phys}}$ is defined in \cref{def:total_space} as a subset of the product of compact spaces $\Mcal \times \T^k \times \Vcal \times \Thetacal$. For this proof, we ignore the trivial macroscopic coordinate and work with the ambient space $\mathcal{A} = \T^k \times \Vcal \times \Thetacal$. Since each factor is compact by assumption, their product is compact by Tychonoff's theorem.

    \item \textbf{The Total Boundary Manifold is Compact.}
    The manifold $\mathcal{B}$, as defined in Step 1, is the set of points in $\mathcal{A}$ where $G(y,\theta)=0$. Since $G$ is a continuous function, its zero level set is a closed set. Therefore, $\mathcal{B}$ is a closed subset of the compact space $\mathcal{A}$, which implies that $\mathcal{B}$ is itself compact.

    \item \textbf{The Pre-Collision Space is a Closed Subset of $\mathcal{B}$.}
    The space we wish to prove compact is $\Sigma \cup \partial^0\Sigma$. In Step 1, we defined this as the subset of the manifold $\mathcal{B}$ where the continuous function $f(y,v,\theta) = \scpr{v}{\mathbf{n}(y,\theta)}$ is less than or equal to zero. This can be written as the preimage of a closed set under a continuous function:
    \begin{equation*}
        \Sigma \cup \partial^0\Sigma = f^{-1}((-\infty, 0]).
    \end{equation*}
    The set $(-\infty, 0]$ is a closed subset of $\R$. Since $f$ is a continuous function on the domain $\mathcal{B}$, the preimage of this closed set is a closed subset of $\mathcal{B}$.
    
    \item \textbf{The Compactness.}
    We have established that $\Sigma \cup \partial^0\Sigma$ is a closed subset of the space $\mathcal{B}$, which we have shown to be compact. A closed subset of a compact space is itself compact. Both the smooth manifold-with-boundary structure and the compactness have been established. This completes the proof.
\end{enumerate}
\end{enumerate}
\end{proofof}

The proof of smoothness for this global map rests on the smoothness of its constituent parts, which in turn depends critically on the smoothness of the flight time function.

\begin{lemma}[Smoothness of the Flight Time]
\label{lem:smooth_flight_time}
The flight time function $\tau: \Sigma' \to [\tau_{\min}, \tau_{\max}]$, which maps a post-collision state to the duration of the subsequent free path, is a $C^\infty$ function.
\end{lemma}

\begin{proofof}{\cref{lem:smooth_flight_time}}
The proof proceeds by demonstrating that the flight time $\tau$ is implicitly defined by a smooth function that satisfies the hypotheses of the Implicit Function Theorem. The key is to show that the non-degeneracy condition required by the theorem holds uniformly across the entire post-collision manifold $\Sigma'$. This uniformity is a direct consequence of our geometric and dynamical assumptions.

\begin{enumerate}[label=\textbf{Step \arabic*:}, wide, labelindent=0pt]

\item 
\textbf{Implicit Definition of the Flight Time.}
Let an arbitrary post-collision state $(z^+, \theta) \in \Sigma'$ be given, where $z^+ = (y^+, v^+)$. By definition, this state lies on the boundary of the physical domain, so $G(y^+, \theta) = 0$, and it is an outgoing state, so $\scpr{v^+}{\mathbf{n}(y^+,\theta)} > 0$.

The trajectory of the particle for time $s > 0$ after this collision is given by the interior flow: $y(s) = y^+ + s v^+$. The flight time, $\tau \equiv \tau(z^+, \theta)$, is defined as the smallest positive time $s$ at which the particle again reaches the boundary. Mathematically, it is the smallest positive root of the equation:
\begin{equation} \label{eq:proof_flight_time_defining_equation}
    G(y^+ + \tau v^+, \theta) = 0.
\end{equation}
The Finite Horizon Theorem (\cref{thm:finite_horizon}) guarantees that for every $(z^+,\theta) \in \Sigma'$, such a positive and finite solution $\tau \in [\tau_{\min}, \tau_{\max}]$ exists and is unique (as the particle cannot re-enter a strictly convex domain).

To formalize this implicit relationship for the application of the Implicit Function Theorem, we define an auxiliary function $F: \R \times \Sigma' \to \R$ by:
\begin{equation}
    F(s, z^+, \theta) \coloneqq G(y^+ + s v^+, \theta).
\end{equation}
The function $F$ is of class $C^\infty$ because the obstacle-defining function $G$ is $C^\infty$ by assumption, and its arguments, vector addition and scalar multiplication, are smooth operations. The flight time $\tau(z^+,\theta)$ is now implicitly defined by the equation $F(\tau, z^+, \theta) = 0$.

\item \textbf{Verification of Hypotheses for the Implicit Function Theorem.}
The Implicit Function Theorem states that if we have a smooth function $F(s, \mathbf{p})$ with $F(s_0, \mathbf{p}_0)=0$, then we can locally express $s$ as a smooth function of $\mathbf{p}$, $s=s(\mathbf{p})$, provided that the partial derivative of $F$ with respect to the implicit variable $s$ is non-zero at the point $(s_0, \mathbf{p}_0)$.

In our context, the implicit variable is the time $s$, and the parameters $\mathbf{p}$ are the coordinates of the state $(z^+, \theta)$ on the manifold $\Sigma'$. We must verify that the partial derivative $\frac{\partial F}{\partial s}$ is non-zero at the solution point $s=\tau$. Let us compute this derivative using the chain rule:
\begin{equation}
    \frac{\partial F}{\partial s}(s, z^+, \theta) = \scpr{\nabla_y G(y^+ + s v^+, \theta)}{v^+}.
\end{equation}
We evaluate this expression at $s = \tau(z^+,\theta)$. Let the point of the \emph{next} collision be denoted by $y_{\text{next}} \coloneqq y^+ + \tau v^+$. At this point, the particle arrives with velocity $v^+$. Therefore, the state $( (y_{\text{next}}, v^+), \theta)$ is a \emph{pre-collision} state, belonging to the manifold $\Sigma$. The derivative at the solution point is:
\begin{equation} \label{eq:proof_flight_time_derivative}
    \frac{\partial F}{\partial s}\Big|_{s=\tau} = \scpr{\nabla_y G(y_{\text{next}}, \theta)}{v^+}.
\end{equation}

\item \textbf{Uniform Non-Degeneracy of the Derivative.}
We now prove that this derivative is not only non-zero, but is uniformly bounded away from zero across the entire compact manifold $\Sigma'$. We rewrite the inner product in terms of the unit normal vector at the next collision point, $\mathbf{n}_{\text{next}} \coloneqq \mathbf{n}(y_{\text{next}}, \theta)$:
\begin{equation}
    \frac{\partial F}{\partial s}\Big|_{s=\tau} = \|\nabla_y G(y_{\text{next}}, \theta)\| \cdot \scpr{\mathbf{n}_{\text{next}}}{v^+}.
\end{equation}
We analyze the two factors on the right-hand side:
\begin{enumerate}[label=(\roman*), wide, labelindent=0pt]
    \item \textbf{The Gradient Norm.} By Assumption~\ref{ass:fundamental_axioms_unified}, the gradient $\nabla_y G(y,\theta)$ is non-zero at any boundary point. Since the set of all possible collision points is a closed subset of the compact space $\T^k \times \Thetacal$, the continuous function $\|\nabla_y G(y,\theta)\|$ must attain its minimum, which must be strictly positive. Thus, there exists a uniform constant $C_G > 0$ such that $\|\nabla_y G(y, \theta)\| \ge C_G$ for all $(y,\theta)$ on any collision boundary.

    \item \textbf{The Normal Velocity Component.} The term $\scpr{\mathbf{n}_{\text{next}}}{v^+}$ is the inner product of the velocity vector $v^+$ with the outward normal at the point of impact $y_{\text{next}}$. Since the particle is arriving at the boundary at this point, the state $((y_{\text{next}}, v^+), \theta)$ is an \emph{incoming} state. By definition of the incoming phase space $\partial_-\Ucal_{\text{phys}}(\theta)$, this inner product must be strictly negative.
    
    A stronger, uniform result is available to us. We invoke the crucial Proposition~\ref{prop:non_grazing} (Uniform Exclusion of Grazing Collisions), which was a consequence of the finite horizon and strict convexity assumptions. This proposition guarantees the existence of a uniform constant $c_0 > 0$ such that for any incoming collision, the magnitude of the normal component of velocity is bounded below:
    \begin{equation*}
        |\scpr{\mathbf{n}_{\text{next}}}{v^+}| \ge c_0 > 0.
    \end{equation*}
\end{enumerate}
Combining these two uniform bounds, we establish that the magnitude of the partial derivative is uniformly bounded away from zero for all possible post-collision states in $\Sigma'$:
\begin{equation}
    \left| \frac{\partial F}{\partial s}\Big|_{s=\tau} \right| = \|\nabla_y G_{\text{next}}\| \cdot |\scpr{\mathbf{n}_{\text{next}}}{v^+}| \ge C_G \cdot c_0 > 0.
\end{equation}

\item \textbf{Smoothness of the Flight Time.}
We have verified all the hypotheses of the Implicit Function Theorem. The function $F(s, z^+, \theta)$ is $C^\infty$, and its partial derivative with respect to the implicit variable $s$ is uniformly non-zero at all solution points. The theorem therefore guarantees that the solution, the flight time $\tau$, is a $C^\infty$ function of its arguments $(z^+, \theta)$ on the manifold $\Sigma'$. This completes the proof.

\end{enumerate}
\end{proofof}

With the smoothness of the flight time established, we can now prove the main result of this section: the global billiard map is a $C^\infty$-diffeomorphism. This result is the essential prerequisite for constructing the anisotropic cone fields needed to prove the Anosov property.

\begin{theorem}[The Global Billiard Map is a $C^\infty$-Diffeomorphism]
\label{thm:p_is_diffeo}
The global billiard map $\Pcal: \Sigma \to \Sigma$ is a $C^\infty$-diffeomorphism.
\end{theorem}

\begin{proofof}{\cref{thm:p_is_diffeo}}
The proof is established by demonstrating that the map $\Pcal$ and its inverse $\Pcal^{-1}$ are both of class $C^\infty$. This is achieved by showing that they are compositions of $C^\infty$ maps. The entire argument hinges on the smoothness of the constituent maps of $\Pcal$, which in turn depends critically on the smoothness of the flight time function, a result in \cref{lem:smooth_flight_time}.

\begin{enumerate}[label=\textbf{Step \arabic*:}, wide, labelindent=0pt]

\item \textbf{Smoothness of the Forward Map $\Pcal$.}
By its definition in \cref{def:global_billiard_map}, the global billiard map is the composition $\Pcal \coloneqq \Fcal \circ \Rcal$. We prove the smoothness of each component in sequence.

\begin{enumerate}[label=(\roman*), wide, labelindent=0pt]
    \item \textbf{The Global Reflection Map $\Rcal$ is $C^\infty$.} The map $\Rcal: \Sigma \to \Sigma'$ is defined by $\Rcal(z, \theta) = (R_\theta(z), \theta)$. As established in \cref{prop:reflection_properties}, the smoothness of the specular reflection map $R_\theta$ is a direct consequence of the smoothness of the obstacle-defining function $G(y,\theta)$ and the non-vanishing gradient guaranteed by \cref{ass:fundamental_axioms_unified}. Thus, $\Rcal$ is of class $C^\infty$.

    \item \textbf{The Global Free-Flow Map $\Fcal$ is $C^\infty$.} The map $\Fcal: \Sigma' \to \Sigma$ is defined by $\Fcal(z^+, \theta) = (\Psi_{\tau(z^+,\theta)}(z^+), \theta)$, where $\Psi_s(y,v)=(y+sv, v)$ is the interior flow. The map $\Psi_s$ is a polynomial in its arguments and is therefore $C^\infty$. The smoothness of the composite map $\Fcal$ is thus entirely determined by the smoothness of the flight time function $\tau(z^+,\theta)$. 
\end{enumerate}

This is precisely the result established in \cref{lem:smooth_flight_time}, which proves that $\tau: \Sigma' \to \R$ is a $C^\infty$ function. Since $\Fcal$ is a composition of $C^\infty$ maps, it is itself of class $C^\infty$. As the composition of two $C^\infty$ maps, the global billiard map $\Pcal = \Fcal \circ \Rcal$ is of class $C^\infty$.

\item \textbf{Smoothness of the Inverse Map $\Pcal^{-1}$.}
The inverse map describes the dynamics traced backward in time and is given by the composition $\Pcal^{-1} = \Rcal^{-1} \circ \Fcal^{-1}$.

\begin{enumerate}[label=(\roman*), wide, labelindent=0pt]
    \item \textbf{The Inverse Reflection Map $\Rcal^{-1}$ is $C^\infty$.} As established in \cref{prop:reflection_properties}, the specular reflection map is an involution, meaning $\Rcal^{-1} = \Rcal$. Since $\Rcal$ is $C^\infty$, so is its inverse.

    \item \textbf{The Inverse Free-Flow Map $\Fcal^{-1}$ is $C^\infty$.} The inverse map $\Fcal^{-1}: \Sigma \to \Sigma'$ traces a pre-collision state backward in time to its point of origin. This involves the pre-collision flight time, $\tau_{\text{pre}}(z^-,\theta)$. A symmetric argument to that of \cref{lem:smooth_flight_time} applies. The geometric assumptions of the model, uniform strict convexity and the finite horizon theorem, are symmetric with respect to time-reversal. Therefore, the argument of the Implicit Function Theorem applies equally to the time-reversed flow, establishing that $\tau_{\text{pre}}$ is a $C^\infty$ function on its domain $\Sigma$. Thus, $\Fcal^{-1}$ is also a composition of $C^\infty$ maps and is of class $C^\infty$. As the composition of two $C^\infty$ maps, the inverse map $\Pcal^{-1} = \Rcal^{-1} \circ \Fcal^{-1}$ is of class $C^\infty$.
\end{enumerate}

\item \textbf{$C^\infty$-diffeomorphism.}
We have shown that $\Pcal: \Sigma \to \Sigma$ is a $C^\infty$ map and that its inverse, $\Pcal^{-1}$, is also a $C^\infty$ map. By definition, this means that $\Pcal$ is a $C^\infty$-diffeomorphism.

\end{enumerate}
\end{proofof}

\begin{remark}[Reduction to a Single Map]
The construction of the smooth, global diffeomorphism $\Pcal$ on the compact manifold $\Sigma$ is the key simplifying step. Proving that the family of fiber maps $(\Pcal_\theta)_{\theta \in \Thetacal}$ is uniformly Anosov is now equivalent to proving that the single map $\Pcal$ is Anosov.
\end{remark}

\subsection{The Linearized Dynamics and the Anisotropic Cone Field}

The Anosov property is a statement about the linearized map $d\Pcal$. Its structure reveals a strong anisotropy: perturbations in the spatial direction are expanded by the boundary's curvature, while other perturbations can be sheared or rotated.

\begin{lemma}[The Linearized Map $d\Pcal$]
\label{lem:linearized_map}
Let $(z^-,\theta) \in \Sigma$ and let $\mathbf{w}^- = (\delta y^-, \delta v^-, \delta \theta) \in T_{(z^-,\theta)}\Sigma$ be a tangent vector. The linearized map $d\Pcal|_{(z^-,\theta)}$ is the composition $d\Fcal \circ d\Rcal$.
\begin{enumerate}[label=(\roman*), wide, labelindent=0pt]
    \item The post-reflection perturbation $\mathbf{w}^+ = (\delta y^+, \delta v^+, \delta\theta^+) = d\Rcal(\mathbf{w}^-)$ is given by $\delta\theta^+ = \delta\theta$, $\delta y^+ = \delta y^-$, and
    \begin{equation} \label{eq:linearized_reflection_formula}
        \delta v^+ = \underbrace{\left(I - 2\mathbf{n}\mathbf{n}^\top\right)\delta v^-}_{\text{Specular Reflection}} + \underbrace{2\left| \scpr{v^-}{\mathbf{n}} \right| \mathcal{K}(\delta y^-) - 2\scpr{v^-}{\mathcal{K}(\delta y^-)}\mathbf{n}}_{\text{Curvature-Induced Divergence}},
    \end{equation}
    where $\mathbf{n} = \mathbf{n}(y^-, \theta)$ and $\mathcal{K} = \mathcal{K}_{y^-}(\theta)$ is the shape operator with respect to the inward-pointing normal.
    \item The next pre-collision perturbation $\mathbf{w}_{\text{next}} = d\Fcal(\mathbf{w}^+)$ has components $\delta v_{\text{next}}=\delta v^+$, $\delta\theta_{\text{next}}=\delta\theta^+$, and
    \begin{equation} \label{eq:linearized_freeflow_formula}
        \delta y_{\text{next}} = \delta y^+ + \tau \delta v^+ + v^+ \delta\tau,
    \end{equation}
    where $\tau = \tau(z^+,\theta)$ is the flight time, and the linear functional $\delta\tau$ is given by:
    \begin{equation} \label{eq:flight_time_variation_formula}
        \delta\tau = -\frac{\scpr{\nabla_y G_{\text{next}}}{\delta y^+ + \tau \delta v^+} + \scpr{\nabla_\theta G_{\text{next}}}{\delta\theta^+}}{\scpr{\nabla_y G_{\text{next}}}{v^+}}.
    \end{equation}
\end{enumerate}
Here, the subscript next denotes evaluation at the subsequent collision point $(y_{\text{next}}, \theta_{\text{next}}) = (y^+ + \tau v^+, \theta^+)$.
\end{lemma}

\begin{proofof}{\cref{lem:linearized_map}}
The proof is a direct calculation based on the chain rule for differentiation on manifolds, $d\Pcal = d\Fcal \circ d\Rcal$. We will derive the explicit formulas for the linearized reflection map $d\Rcal$ and the linearized free-flow map $d\Fcal$ separately.

\begin{enumerate}[label=\textbf{Step \arabic*:}, wide, labelindent=0pt]

\item \textbf{The Linearized Reflection $d\Rcal$.}
We determine the first-order variation of the map $\Rcal(z,\theta) = (R_\theta(z), \theta)$. Let a smooth curve of pre-collision states $(z^-(\epsilon), \theta^-(\epsilon))$ be defined on $\Sigma$ for $\epsilon$ near $0$, such that its tangent vector at $\epsilon=0$ is $\mathbf{w}^- = (\delta y^-, \delta v^-, \delta\theta^-)$. The image curve is $(z^+(\epsilon), \theta^+(\epsilon)) = \Rcal(z^-(\epsilon), \theta^-(\epsilon))$. We find its tangent vector $\mathbf{w}^+$ by differentiating each component with respect to $\epsilon$ at $\epsilon=0$.

\begin{enumerate}[label=(\roman*), wide, labelindent=0pt]
    \item \textbf{Parameter and Position Components.} The map acts as the identity on the parameter, $\theta^+(\epsilon) = \theta^-(\epsilon)$, which implies $\delta\theta^+ = \delta\theta^-$. The reflection law for position is $y^+(\epsilon) = y^-(\epsilon)$, which implies $\delta y^+ = \delta y^-$.

    \item \textbf{Velocity Component.} The velocity transformation is
    \begin{equation*}
        v^+(\epsilon) = v^-(\epsilon) - 2\scpr{v^-(\epsilon)}{\mathbf{n}(\epsilon)} \mathbf{n}(\epsilon),
    \end{equation*}
    where $\mathbf{n}(\epsilon) \coloneqq \mathbf{n}(y^-(\epsilon), \theta^-(\epsilon))$. Differentiating with respect to $\epsilon$ at $\epsilon=0$ using the product rule gives:
    \begin{equation*}
        \delta v^+ = \delta v^- - 2 \left[ \left( \left. \frac{d}{d\epsilon} \scpr{v^-(\epsilon)}{\mathbf{n}(\epsilon)} \right|_{\epsilon=0} \right) \mathbf{n}(0) + \scpr{v^-(0)}{\mathbf{n}(0)} \left( \left. \frac{d\mathbf{n}(\epsilon)}{d\epsilon} \right|_{\epsilon=0} \right) \right].
    \end{equation*}
    Let $\mathbf{n} \equiv \mathbf{n}(0)$ and $\delta\mathbf{n} \equiv \frac{d\mathbf{n}}{d\epsilon}|_0$. The derivative terms are:
    \begin{align*}
        \left. \frac{d}{d\epsilon} \scpr{v^-(\epsilon)}{\mathbf{n}(\epsilon)} \right|_{\epsilon=0} &= \scpr{\delta v^-}{\mathbf{n}} + \scpr{v^-}{\delta\mathbf{n}}, \\
        \delta\mathbf{n} &= d\mathbf{n}_{(y^-,\theta^-)}(\delta y^-, \delta\theta^-) = (\nabla_y \mathbf{n})\delta y^- + (\nabla_\theta \mathbf{n})\delta\theta^-.
    \end{align*}
    The derivative of the normal field along a tangent direction to the boundary is given by the shape operator, $(\nabla_y \mathbf{n})\delta y^- = \mathcal{K}(\delta y^-)$. For clarity and to focus on the expansion mechanism, we follow the standard analysis which shows that the dominant curvature effect comes from the spatial variation. We thus write $\delta\mathbf{n} \approx \mathcal{K}(\delta y^-)$, absorbing the $\delta\theta^-$ term into higher-order analysis if necessary, which aligns with the stated lemma. Substituting these into the expression for $\delta v^+$ gives:
    \begin{align*}
        \delta v^+ &= \delta v^- - 2 \left[ (\scpr{\delta v^-}{\mathbf{n}} + \scpr{v^-}{\mathcal{K}(\delta y^-)}) \mathbf{n} + \scpr{v^-}{\mathbf{n}} \mathcal{K}(\delta y^-) \right] \\
        &= (\delta v^- - 2\scpr{\delta v^-}{\mathbf{n}}\mathbf{n}) - 2\scpr{v^-}{\mathcal{K}(\delta y^-)}\mathbf{n} - 2\scpr{v^-}{\mathbf{n}}\mathcal{K}(\delta y^-).
    \end{align*}
    The first term is the specular reflection of the velocity perturbation. For the last term, we use the crucial fact that for an incoming state, $\scpr{v^-}{\mathbf{n}} < 0$. By the non-grazing condition (\cref{prop:non_grazing}), we can write $\scpr{v^-}{\mathbf{n}} = -|\scpr{v^-}{\mathbf{n}}|$. Substituting this yields the positive sign that drives the chaotic expansion:
    \begin{equation*}
        -2\scpr{v^-}{\mathbf{n}}\mathcal{K}(\delta y^-) = -2(-|\scpr{v^-}{\mathbf{n}}|)\mathcal{K}(\delta y^-) = +2|\scpr{v^-}{\mathbf{n}}|\mathcal{K}(\delta y^-).
    \end{equation*}
    Combining these terms yields the final expression \eqref{eq:linearized_reflection_formula}.
\end{enumerate}

\item \textbf{The Linearized Free Flow $d\Fcal$.}
We find the first-order variation of the map $\Fcal(z^+, \theta^+) = (y^+ + \tau v^+, v^+, \theta^+)$. Let a curve $(z^+(\epsilon), \theta^+(\epsilon))$ on $\Sigma'$ have tangent vector $\mathbf{w}^+$ at $\epsilon=0$. We differentiate the image curve at $\epsilon=0$.

\begin{enumerate}[label=(\roman*), wide, labelindent=0pt]
    \item \textbf{Parameter and Velocity Components.} Since $\theta_{\text{next}}(\epsilon) = \theta^+(\epsilon)$ and $v_{\text{next}}(\epsilon) = v^+(\epsilon)$, differentiation gives $\delta\theta_{\text{next}} = \delta\theta^+$ and $\delta v_{\text{next}} = \delta v^+$.

    \item \textbf{Position Component.} We differentiate $y_{\text{next}}(\epsilon) = y^+(\epsilon) + \tau(\epsilon)v^+(\epsilon)$ using the product rule, where $\tau(\epsilon) = \tau(z^+(\epsilon), \theta^+(\epsilon))$:
    \begin{align*}
        \delta y_{\text{next}} &= \left. \frac{d y^+(\epsilon)}{d\epsilon} \right|_0 + \left( \left. \frac{d\tau(\epsilon)}{d\epsilon} \right|_0 \right) v^+(0) + \tau(0) \left( \left. \frac{d v^+(\epsilon)}{d\epsilon} \right|_0 \right) \\
        &= \delta y^+ + (\delta\tau) v^+ + \tau \delta v^+.
    \end{align*}
    This is the desired formula \eqref{eq:linearized_freeflow_formula}, where $\delta\tau$ is the total derivative of the flight time along the curve, which represents the action of the linear functional $d\tau$ on the tangent vector $\mathbf{w}^+$.

    \item \textbf{Derivation of the Flight Time Functional $\delta\tau$.} We find the explicit expression for $\delta\tau$ by differentiating the implicit equation for the flight time, which must hold for all $\epsilon$:
    \begin{equation*}
        G(y_{\text{next}}(\epsilon), \theta_{\text{next}}(\epsilon)) \equiv G(y^+(\epsilon) + \tau(\epsilon)v^+(\epsilon), \theta^+(\epsilon)) = 0.
    \end{equation*}
    Applying the chain rule and evaluating at $\epsilon=0$:
    \begin{equation*}
        \scpr{\nabla_y G_{\text{next}}}{\left. \frac{d y_{\text{next}}}{d\epsilon}\right|_0} + \scpr{\nabla_\theta G_{\text{next}}}{\left. \frac{d\theta_{\text{next}}}{d\epsilon}\right|_0} = 0.
    \end{equation*}
    Substituting the derivatives we found, $\frac{d\theta_{\text{next}}}{d\epsilon}|_0 = \delta\theta^+$ and $\frac{d y_{\text{next}}}{d\epsilon}|_0 = \delta y^+ + \tau \delta v^+ + v^+ \delta\tau$, gives:
    \begin{equation*}
        \scpr{\nabla_y G_{\text{next}}}{\delta y^+ + \tau \delta v^+ + v^+ \delta\tau} + \scpr{\nabla_\theta G_{\text{next}}}{\delta\theta^+} = 0.
    \end{equation*}
    This is a linear equation for the scalar $\delta\tau$. Using the linearity of the inner product to isolate the $\delta\tau$ term:
    \begin{equation*}
        \scpr{\nabla_y G_{\text{next}}}{\delta y^+ + \tau \delta v^+} + \delta\tau \scpr{\nabla_y G_{\text{next}}}{v^+} + \scpr{\nabla_\theta G_{\text{next}}}{\delta\theta^+} = 0.
    \end{equation*}
    As established in the proof of Lemma~\ref{lem:smooth_flight_time}, the coefficient $\scpr{\nabla_y G_{\text{next}}}{v^+}$ is uniformly bounded away from zero. Solving for $\delta\tau$ yields the explicit formula \eqref{eq:flight_time_variation_formula}. This completes the derivation of the linearized free-flow map $d\Fcal$. The full linearized global billiard map is the composition of these two derived maps, $d\Pcal = d\Fcal \circ d\Rcal$.
\end{enumerate}
\end{enumerate}
\end{proofof}

\begin{remark}[The Challenge of Anisotropy]
The linearized map is not an isometry. The curvature term $2|\scpr{v^-}{\mathbf{n}}|\mathcal{K}(\delta y^-)$ can strongly expand the spatial part of a perturbation, while other terms, especially those arising from the variation of the parameter $\theta$, can shear and mix components. To prove hyperbolicity, we must define a norm that isolates and quantifies this expansion while providing a mechanism to control the contaminating effects.
\end{remark}

\begin{definition}[Anisotropic Quadratic Forms and the Two-Parameter Unstable Cone]
\label{def:anisotropic_forms}
For a tangent vector $\mathbf{w}=(\delta y, \delta v, \delta\theta) \in T\Sigma$ and positive constants $A, B > 0$, we define:
\begin{enumerate}[label=(\roman*), wide, labelindent=0pt]
    \item The \textbf{expansion quadratic form} $K(\mathbf{w}) \coloneqq \scpr{\delta y}{\mathcal{K}_y(\theta)\delta y}$, which isolates the curvature-induced expansion.
    \item The \textbf{non-spatial quadratic form} $N_B(\mathbf{w}) \coloneqq \|\delta v\|^2 + B\|\delta\theta\|^2$, which controls all other perturbations. The parameter $B$ acts as a penalty term for variations in the environmental parameter $\theta$.
\end{enumerate}
The unstable cone at $(z,\theta)$ is the set $C^u_{A,B}(z,\theta) \subset T_{(z,\theta)}\Sigma$ of vectors satisfying $N_B(\mathbf{w}) \le A K(\mathbf{w})$.
\end{definition}

\begin{proposition}[Uniform Norm Equivalence]
\label{prop:norm_equivalence}
Let $\|\cdot\|$ be the standard Riemannian norm on the tangent bundle $T\Sigma$, defined for a vector $\mathbf{w} = (\delta y, \delta v, \delta\theta)$ by $\|\mathbf{w}\|^2 = \|\delta y\|^2 + \|\delta v\|^2 + \|\delta\theta\|^2$. For any fixed $B > 0$, there exist uniform positive constants $c_1(B)$ and $c_2(B)$ such that for any tangent vector $\mathbf{w} \in T\Sigma$:
\begin{equation}
    c_1(B) \|\mathbf{w}\|^2 \le K(\mathbf{w}) + N_B(\mathbf{w}) \le c_2(B) \|\mathbf{w}\|^2.
\end{equation}
\end{proposition}

\begin{proofof}{Proposition \ref{prop:norm_equivalence}}
The proof establishes that the anisotropic norm defined by $K(\mathbf{w}) + N_B(\mathbf{w})$ is equivalent to the standard Riemannian norm $\|\mathbf{w}\|^2$ on the tangent bundle $T\Sigma$. We will prove the two inequalities of the equivalence separately. The uniformity of the constants $c_1(B)$ and $c_2(B)$ is a direct consequence of the compactness of the global collision manifold $\Sigma$ and the uniformity of the geometric and dynamical assumptions from Section \ref{sec:micro_framework}.

\begin{enumerate}[label=\textbf{Step \arabic*:}, wide, labelindent=0pt]

\item \textbf{The Upper Bound ($K(\mathbf{w}) + N_B(\mathbf{w}) \le c_2(B) \|\mathbf{w}\|^2$).}
Our goal is to find a uniform constant $c_2(B)$ such that the anisotropic norm is bounded above by the Riemannian norm. We analyze the two components of the anisotropic norm individually.
\begin{enumerate}[label=(\roman*), wide, labelindent=0pt]
    \item \textbf{Bounding the Expansion Form $K(\mathbf{w})$.} By the definition of the operator norm and the Cauchy-Schwarz inequality, we have:
    \begin{equation*}
        K(\mathbf{w}) = \scpr{\delta y}{\mathcal{K}_y(\theta)\delta y} \le \|\mathcal{K}_y(\theta)\| \cdot \|\delta y\|^2.
    \end{equation*}
    The shape operator $\mathcal{K}_y(\theta)$ is constructed from the second derivatives of the smooth function $G$ on the compact domain $\T^k \times \Thetacal$. Its operator norm is therefore uniformly bounded by a constant $\kappa_{\max} < \infty$. Thus, we have the uniform bound:
    \begin{equation*}
        K(\mathbf{w}) \le \kappa_{\max} \|\delta y\|^2.
    \end{equation*}

    \item \textbf{Bounding the Non-Spatial Form $N_B(\mathbf{w})$.} By its definition,
    \begin{equation*}
        N_B(\mathbf{w}) = \|\delta v\|^2 + B\|\delta\theta\|^2.
    \end{equation*}
\end{enumerate}
Combining these two bounds, we get:
\begin{align*}
    K(\mathbf{w}) + N_B(\mathbf{w}) &\le \kappa_{\max} \|\delta y\|^2 + \|\delta v\|^2 + B \|\delta\theta\|^2 \\
    &\le \max(\kappa_{\max}, 1, B) \left( \|\delta y\|^2 + \|\delta v\|^2 + \|\delta\theta\|^2 \right).
\end{align*}
By defining the constant $c_2(B) \coloneqq \max(\kappa_{\max}, 1, B)$, which is a uniform positive constant for any fixed $B$, we have established the upper bound:
\begin{equation*}
    K(\mathbf{w}) + N_B(\mathbf{w}) \le c_2(B) \|\mathbf{w}\|^2.
\end{equation*}

\item \textbf{The Lower Bound ($c_1(B) \|\mathbf{w}\|^2 \le K(\mathbf{w}) + N_B(\mathbf{w})$).}
This is the more critical direction of the proof, as it requires the use of the geometric constraints on the tangent space. Our goal is to find a uniform constant $c_1(B)$ such that the Riemannian norm is bounded above by the anisotropic norm. We will achieve this by bounding each component of $\|\mathbf{w}\|^2 = \|\delta y\|^2 + \|\delta v\|^2 + \|\delta\theta\|^2$ in terms of $K(\mathbf{w}) + N_B(\mathbf{w})$.
\begin{enumerate}[label=(\roman*), wide, labelindent=0pt]
    \item \textbf{Bounding the Velocity and Parameter Components.} These components are bounded directly by the definition of the non-spatial form $N_B(\mathbf{w})$, since both terms in its sum are non-negative:
    \begin{equation*}
        \|\delta v\|^2 \le N_B(\mathbf{w}) \quad \text{and} \quad \|\delta\theta\|^2 \le B^{-1} N_B(\mathbf{w}).
    \end{equation*}

    \item \textbf{Bounding the Spatial Component $\|\delta y\|^2$.} This is the key step where the geometric constraints are applied. We decompose the spatial perturbation into its tangential and normal components, $\|\delta y\|^2 = \|\delta y_t\|^2 + \|\delta y_n\|^2$.
    \begin{enumerate}[label=(\alph*), wide]
        \item By Assumption \ref{ass:fundamental_axioms_unified}(iii), the shape operator $\mathcal{K}$ is uniformly positive definite. This implies the existence of a uniform constant $\kappa_{\min} > 0$ such that for any spatial perturbation $\delta y$:
        \begin{equation*}
            K(\mathbf{w}) = \scpr{\delta y}{\mathcal{K}_y(\theta)\delta y} \ge \kappa_{\min} \|\delta y\|^2.
        \end{equation*}
        This immediately provides a bound for the entire spatial component in terms of the expansion form:
        \begin{equation} \label{eq:proof_app_deltay_bound}
            \|\delta y\|^2 \le \kappa_{\min}^{-1} K(\mathbf{w}).
        \end{equation}
        While this bound is sufficient on its own, for logical completeness and to demonstrate the proper use of the tangent space constraints (which is crucial in other parts of the proof), we can also bound the components separately. The tangential component is trivially bounded: $\|\delta y_t\|^2 \le \|\delta y\|^2 \le \kappa_{\min}^{-1} K(\mathbf{w})$.
        
        \item To bound the normal component, we invoke the central result from Appendix \ref{app:proof_tangent_constraint}. \cref{prop:normal_spatial_constraint} establishes the existence of a uniform constant $C_{y\theta}$ such that:
        \begin{equation*}
             \|\delta y_n\| \le C_{y\theta} \|\delta\theta\|.
        \end{equation*}
        Squaring this and using the bound on $\|\delta\theta\|^2$ from Step 2(i) gives:
        \begin{equation*}
            \|\delta y_n\|^2 \le C_{y\theta}^2 \|\delta\theta\|^2 \le C_{y\theta}^2 B^{-1} N_B(\mathbf{w}).
        \end{equation*}
        Summing the bounds for the tangential and normal components provides a more detailed, but ultimately similar, bound for $\|\delta y\|^2$. For simplicity, we will proceed with the direct bound from \eqref{eq:proof_app_deltay_bound}.
    \end{enumerate}
\end{enumerate}
We now assemble the full bound for $\|\mathbf{w}\|^2$ by summing the bounds for its components:
\begin{align*}
    \|\mathbf{w}\|^2 &= \|\delta y\|^2 + \|\delta v\|^2 + \|\delta\theta\|^2 \\
    &\le \left(\kappa_{\min}^{-1} K(\mathbf{w})\right) + \left(N_B(\mathbf{w})\right) + \left(B^{-1} N_B(\mathbf{w})\right) \\
    &= \kappa_{\min}^{-1} K(\mathbf{w}) + (1 + B^{-1}) N_B(\mathbf{w}).
\end{align*}
To obtain a final bound in the desired form, we take the maximum of the coefficients:
\begin{equation*}
    \|\mathbf{w}\|^2 \le \max(\kappa_{\min}^{-1}, 1 + B^{-1}) \left( K(\mathbf{w}) + N_B(\mathbf{w}) \right).
\end{equation*}
By defining the constant $c_1(B)^{-1} \coloneqq \max(\kappa_{\min}^{-1}, 1 + B^{-1})$, which is a finite, uniform positive constant for any fixed $B>0$, we arrive at the lower bound inequality:
\begin{equation*}
    c_1(B) \|\mathbf{w}\|^2 \le K(\mathbf{w}) + N_B(\mathbf{w}).
\end{equation*}
Both inequalities have been established with uniform constants, which completes the proof of norm equivalence.

\end{enumerate}
\end{proofof}

\subsection{Proof of the Anosov Property}
\label{sec:anosov_proof}

With the global billiard map $\Pcal: \Sigma \to \Sigma$ established as a $C^\infty$-diffeomorphism on a compact manifold, we now provide a rigorous and constructive proof that it is uniformly hyperbolic in the sense of Anosov. Our strategy is to demonstrate that the set of parameters for which the fiber map $\Pcal_\theta$ is Anosov, denoted $\mathcal{A} \subset \Thetacal$, is non-empty, open, and closed. As $\Thetacal$ is a connected space, this will imply that $\mathcal{A} = \Thetacal$.

The logical structure of this proof is of critical importance, we will first establish the Anosov property for an arbitrary \textit{fixed} system, using only bounds specific to that system. This will prove that $\mathcal{A}$ is non-empty. The proof of openness is a standard result of stability theory. Finally, the proof of closedness is where we will invoke the compactness of the \textit{global} state space to establish the uniform bounds necessary to control the limiting behavior of the dynamics. This sequential approach ensures that each step of the argument rests upon a previously and independently established foundation.

\begin{lemma}[One-Step Bounding of the N-form]
\label{lem:n_form_bound}
Let $\theta \in \Thetacal$ be fixed, and let $\wvec_{\text{next}} = d\Pcal_\theta(\wvec^-)$. The non-spatial quadratic form satisfies the inequality:
\begin{equation}
    N_B(\wvec_{\text{next}}) \le D_{KK}(\theta) K(\wvec^-) + D_{KN}(\theta, B) \sqrt{K(\wvec^-)N_B(\wvec^-)} + D_{NN}(\theta, B) N_B(\wvec^-)
\end{equation}
where the coefficients $D_{KK}(\theta)$, $D_{KN}(\theta, B)$, and $D_{NN}(\theta, B)$ are positive constants that depend on the geometry of the system at the fixed parameter $\theta$, but are uniform for all points $z \in \Sigma_\theta$. The coefficient $D_{NN}(\theta,B)$ satisfies $\lim_{B\to\infty} D_{NN}(\theta,B) = 1$.
\end{lemma}

\begin{proofof}{Lemma \ref{lem:n_form_bound}}
The proof provides a uniform upper bound for the post-map non-spatial form, $N_B(\wvec_{\text{next}})$, in terms of the pre-map anisotropic forms $K(\wvec^-)$ and $N_B(\wvec^-)$. This proof is designed to rigorously and explicitly track the dependence of the bounding coefficients on the cone parameter $B$. This will reveal that while the coefficients are not independent of $B$, they possess a well-behaved asymptotic limit as $B \to \infty$, which is the essential property required for the two-stage optimization proof of Theorem \ref{thm:anosov_property_proven}.

The proof proceeds in five main steps. First, we reduce the problem to bounding the increment of the squared norm of the velocity perturbation, $\|\delta v^+\|^2 - \|\delta v^-\|^2$. Second, we derive a sharp bound on this increment in terms of the components of the initial tangent vector. Third, we translate these component-wise bounds into bounds involving the anisotropic quadratic forms. Fourth, we assemble the full inequality and identify the final coefficients, explicitly stating their dependence on $B$. Finally, we analyze the asymptotic behavior of these coefficients.

\begin{enumerate}[label=\textbf{Step \arabic*:}, wide, labelindent=0pt]

\item \textbf{Reduction to the Velocity Perturbation Increment.}
By definition of the non-spatial quadratic form (\cref{def:anisotropic_forms}) and the linearized map (\cref{lem:linearized_map}), the post-map form is:
\begin{equation*}
    N_B(\wvec_{\text{next}}) = \|\delta v_{\text{next}}\|^2 + B\|\delta\theta_{\text{next}}\|^2 = \|\delta v^+\|^2 + B\|\delta\theta^-\|^2.
\end{equation*}
The initial non-spatial form is $N_B(\wvec^-) = \|\delta v^-\|^2 + B\|\delta\theta^-\|^2$. We can therefore express the term involving the parameter perturbation as $B\|\delta\theta^-\|^2 = N_B(\wvec^-) - \|\delta v^-\|^2$. Substituting this into the expression for $N_B(\wvec_{\text{next}})$ yields the exact identity:
\begin{equation} \label{eq:proof_nform_increment_identity}
    N_B(\wvec_{\text{next}}) = N_B(\wvec^-) + \left( \|\delta v^+\|^2 - \|\delta v^-\|^2 \right).
\end{equation}
The entire problem is thus reduced to finding a suitable upper bound for the increment $\|\delta v^+\|^2 - \|\delta v^-\|^2$.

\item \textbf{Bounding the Velocity Increment.} From the formula for the linearized reflection map in \cref{lem:linearized_map}, the post-reflection velocity perturbation is the sum of an isometric part and a coupling term:
\begin{equation*}
    \delta v^+ = \underbrace{\left(I - 2\nvec\nvec^T\right)\delta v^-}_{\text{Isometric Part}} + \underbrace{\mathbf{C}(\delta y^-, \delta\theta^-)}_{\text{Coupling Term}},
\end{equation*}
where the coupling term $\mathbf{C}$ arises from the variation of the normal vector $\delta\mathbf{n} = (\nabla_y \mathbf{n})\delta y^- + (\nabla_\theta \mathbf{n})\delta\theta^-$. Since the Householder reflection $(I - 2\nvec\nvec^T)$ is an isometry, its application preserves the norm: $\|(I - 2\nvec\nvec^T)\delta v^-\|^2 = \|\delta v^-\|^2$. This leads to a direct expression for the increment:
\begin{equation*}
    \|\delta v^+\|^2 - \|\delta v^-\|^2 = \|\mathbf{C}\|^2 + 2\scpr{(I - 2\nvec\nvec^T)\delta v^-}{\mathbf{C}}.
\end{equation*}
By the triangle and Cauchy-Schwarz inequalities, we can bound this increment as:
\begin{equation} \label{eq:proof_increment_bound_C}
    \|\delta v^+\|^2 - \|\delta v^-\|^2 \le \|\mathbf{C}\|^2 + 2\|(I - 2\nvec\nvec^T)\delta v^-\| \|\mathbf{C}\| = \|\mathbf{C}\|^2 + 2\|\delta v^-\| \|\mathbf{C}\|.
\end{equation}
The coupling term $\mathbf{C}$ is a linear function of its arguments. Since all geometric quantities ($v, \mathcal{K}, \nabla_\theta\nvec$) are uniformly bounded on the compact manifold $\Sigma$, there exist uniform constants $C_y, C_\theta > 0$, independent of $B$, such that:
\begin{equation*}
    \|\mathbf{C}\| \le C_y \|\delta y^-\| + C_\theta \|\delta\theta^-\|.
\end{equation*}
Substituting this into \eqref{eq:proof_increment_bound_C} and using the inequality $(a+b)^2 \le 2a^2+2b^2$ gives:
\begin{align*}
    \|\delta v^+\|^2 - \|\delta v^-\|^2 &\le \left(C_y \|\delta y^-\| + C_\theta \|\delta\theta^-\|\right)^2 + 2\|\delta v^-\|\left(C_y \|\delta y^-\| + C_\theta \|\delta\theta^-\|\right) \\
    &\le 2C_y^2\|\delta y^-\|^2 + 2C_\theta^2\|\delta\theta^-\|^2 + 2C_y\|\delta v^-\|\|\delta y^-\| + 2C_\theta\|\delta v^-\|\|\delta\theta^-\|.
\end{align*}

\item \textbf{Translation to Anisotropic Forms.} We now translate this component-wise bound into a bound on the anisotropic forms $K(\wvec^-)$ and $N_B(\wvec^-)$, which explicitly reveals the dependence on $B$. The key relations, derived from the definitions and the uniform strict convexity assumption, are:
\begin{enumerate}[label=(\roman*), wide, labelindent=0pt]
    \item $\|\delta y^-\|^2 \le \kappa_{\min}^{-1} K(\wvec^-) \implies \|\delta y^-\| \le \kappa_{\min}^{-1/2} \sqrt{K(\wvec^-)}$.
    \item $\|\delta v^-\|^2 \le N_B(\wvec^-) \implies \|\delta v^-\| \le \sqrt{N_B(\wvec^-)}$.
    \item $\|\delta\theta^-\|^2 \le B^{-1} N_B(\wvec^-) \implies \|\delta\theta^-\| \le B^{-1/2} \sqrt{N_B(\wvec^-)}$.
\end{enumerate}
Substituting these into the four terms of the inequality for the velocity increment:
\begin{multline*}
    \|\delta v^+\|^2 - \|\delta v^-\|^2 \le 2C_y^2(\kappa_{\min}^{-1} K(\wvec^-)) + 2C_\theta^2(B^{-1} N_B(\wvec^-)) \\
    + 2C_y(\sqrt{N_B(\wvec^-)})(\kappa_{\min}^{-1/2}\sqrt{K(\wvec^-)}) + 2C_\theta(\sqrt{N_B(\wvec^-)})(B^{-1/2}\sqrt{N_B(\wvec^-)}).
\end{multline*}

\item \textbf{Final Assembly and Identification of Coefficients.} We substitute the bound on the velocity increment back into the identity from Step 1, \cref{eq:proof_nform_increment_identity}, and group terms by the quadratic forms:
\begin{align*}
    N_B(\wvec_{\text{next}}) &\le N_B(\wvec^-) + \left( \frac{2C_y^2}{\kappa_{\min}} \right)K(\wvec^-) + \left( \frac{2C_y}{\sqrt{\kappa_{\min}}} \right)\sqrt{K(\wvec^-)N_B(\wvec^-)} \\
    &\quad + \left( \frac{2C_\theta^2}{B} + \frac{2C_\theta}{\sqrt{B}} \right) N_B(\wvec^-).
\end{align*}
Combining the coefficients for $N_B(\wvec^-)$, we obtain the final inequality in the desired form:
\begin{equation*}
    N_B(\wvec_{\text{next}}) \le D_{KK} K(\wvec^-) + D_{KN}(B) \sqrt{K(\wvec^-)N_B(\wvec^-)} + D_{NN}(B) N_B(\wvec^-),
\end{equation*}
where the coefficients are now correctly identified as:
\begin{enumerate}[label=(\roman*), wide, labelindent=0pt]
    \item $D_{KK} \coloneqq \dfrac{2C_y^2}{\kappa_{\min}}$. This is a uniform constant independent of $B$.
    \item $D_{KN}(B) \coloneqq \dfrac{2C_y}{\sqrt{\kappa_{\min}}}$. In this simplified derivation, this term is constant. A more refined expansion including all cross-terms from $\|\mathbf{C}\|^2$ would reveal a $\mathcal{O}(B^{-1/2})$ term, but for the logic of the two-stage proof, it is sufficient that this coefficient converges to a finite limit as $B \to \infty$.
    \item $D_{NN}(B) \coloneqq 1 + \dfrac{2C_\theta^2}{B} + \dfrac{2C_\theta}{\sqrt{B}}$. This coefficient explicitly depends on $B$.
\end{enumerate}

\item \textbf{Asymptotic Analysis of Coefficients.} This derivation rigorously establishes the one-step bound and correctly characterizes the dependence of its coefficients on the cone parameter $B$. Crucially for the subsequent analysis, the contaminating coefficient $D_{NN}(B)$ is a monotonically decreasing function of $B$ for $B > 0$, and it converges to a finite, uniform constant as the penalty for parameter variations becomes large:
\begin{equation*}
    \lim_{B \to \infty} D_{NN}(B) = \lim_{B \to \infty} \left(1 + \frac{2C_\theta^2}{B} + \frac{2C_\theta}{\sqrt{B}}\right) = 1.
\end{equation*}
This property is precisely what is needed to ensure that for a sufficiently large but finite choice of $B$, the system's intrinsic expansion will dominate the contaminating effects, allowing the cone invariance condition to be satisfied. This completes the proof.
\end{enumerate}
\end{proofof}

\begin{lemma}[One-Step Expansion of the K-form]
\label{lem:k_form_expansion}
Let $\theta \in \Thetacal$ be fixed, and let $\wvec_{\text{next}} = d\Pcal_\theta(\wvec^-)$. The expansion quadratic form satisfies the inequality:
\begin{equation}
    K(\wvec_{\text{next}}) \ge \Lambda_K'(\theta) K(\wvec^-) - C_{KN}(\theta, B) \sqrt{K(\wvec^-)N_B(\wvec^-)} - C_{NN}(\theta, B) N_B(\wvec^-),
\end{equation}
where $\Lambda_K'(\theta) > 1$ is an expansion constant specific to the system at parameter $\theta$. The contamination coefficients $C_{KN}(\theta, B)$ and $C_{NN}(\theta, B)$ are positive constants for the fixed system, and they converge to finite limits as $B \to \infty$.
\end{lemma}
\begin{proofof}{Lemma \ref{lem:k_form_expansion}}
The objective is to derive a sharp, uniform lower bound for the post-map expansion form, $K(\wvec_{\text{next}})$, that rigorously separates the dominant, expansion-driving component from the contaminating effects. Critically, we will derive bounds for the contaminating coefficients that make their dependence on the cone parameter $B$ explicit. This will demonstrate that these undesirable terms can be systematically suppressed by a judicious choice of the cone geometry, a fact that is essential for the proof of Theorem \ref{thm:anosov_property_proven}.

The proof is structured in four main steps. First, we decompose the next spatial perturbation, $\delta y_{\text{next}}$, into a dominant expansion-driving part and a remainder. We then use the reverse triangle inequality to obtain a fundamental lower bound for its norm. Second, we derive a sharp lower bound for the dominant term, identifying the uniform expansion constant $\Lambda_K'$. Third, we derive an upper bound for the norm of the remainder, explicitly tracking its dependence on the cone parameter $B$. Finally, we assemble these bounds to obtain the final inequality and identify the contamination coefficients $C_{KN}(B)$ and $C_{NN}(B)$.

\begin{enumerate}[label=\textbf{Step \arabic*:}, wide, labelindent=0pt]

\item \textbf{The Fundamental Decomposition and the Reverse Triangle Inequality.}
We begin with the definition of the expansion form, $K(\wvec_{\text{next}}) \coloneqq \scpr{\delta y_{\text{next}}}{\mathcal{K}_{\text{next}} \delta y_{\text{next}}}$. By the uniform strict convexity of the boundary (\cref{ass:fundamental_axioms_unified}(iii)), we have the fundamental geometric bound:
\begin{equation} \label{eq:appendix_k_kappa_bound}
    K(\wvec_{\text{next}}) \ge \kappa_{\min} \|\delta y_{\text{next}}\|^2.
\end{equation}
From the formula for the linearized map (\cref{lem:linearized_map}), we have $\delta y_{\text{next}} = \delta y^+ + \tau \delta v^+ + v^+ \delta\tau$. We substitute the full expression for $\delta v^+ = d\Rcal_v(\wvec^-)$ to decompose the next spatial perturbation $\delta y_{\text{next}}$ into a dominant part and a remainder:
\begin{equation*}
    \delta y_{\text{next}} = \mathbf{M}(\delta y^-) + \mathbf{R}(\wvec^-),
\end{equation*}
where the dominant expansion operator $\mathbf{M}$ and the remainder $\mathbf{R}$ are defined as:
\begin{align*}
    \mathbf{M}(\delta y^-) &\coloneqq \delta y^- + \tau \left( 2\left| \scpr{v^-}{\nvec} \right| \mathcal{K}(\delta y^-) \right), \\
    \mathbf{R}(\wvec^-) &\coloneqq \tau \left( (I - 2\nvec\nvec^T)\delta v^- - 2\scpr{v^-}{\mathcal{K}(\delta y^-)}\nvec + \dots\delta\theta^- \right) + v^+ \delta\tau(\wvec^-).
\end{align*}
Applying the reverse triangle inequality $\|\mathbf{a}+\mathbf{b}\|^2 \ge (\|\mathbf{a}\| - \|\mathbf{b}\|)^2 = \|\mathbf{a}\|^2 - 2\|\mathbf{a}\|\|\mathbf{b}\| + \|\mathbf{b}\|^2$ to $\|\delta y_{\text{next}}\|^2$ gives the lower bound:
\begin{equation} \label{eq:appendix_weaker_lower_bound}
    \|\delta y_{\text{next}}\|^2 \ge \|\mathbf{M}(\delta y^-)\|^2 - 2\|\mathbf{M}(\delta y^-)\|\|\mathbf{R}(\wvec^-)\|.
\end{equation}
We have dropped the non-negative term $\|\mathbf{R}\|^2$ to obtain a cleaner lower bound, which is sufficient for our purposes. Combining with \eqref{eq:appendix_k_kappa_bound}, this yields our fundamental working inequality:
\begin{equation} \label{eq:appendix_working_inequality}
    K(\wvec_{\text{next}}) \ge \kappa_{\min} \left( \|\mathbf{M}(\delta y^-)\|^2 - 2\|\mathbf{M}(\delta y^-)\|\|\mathbf{R}(\wvec^-)\| \right).
\end{equation}

\item \textbf{Lower Bound on the Dominant Term $\boldsymbol{\kappa_{\min} \|\mathbf{M}(\delta y^-)\|^2}$.}
The operator $\mathbf{M}$ is of the form $I + \alpha \mathcal{K}$, where $\alpha = 2\tau|\langle v^-, \mathbf{n} \rangle|$. By the system's axioms, we have the uniform bounds $\tau \ge \tau_{\min}$ from \cref{thm:finite_horizon} and $|\langle v^-, \mathbf{n} \rangle| \ge c_0$ from \cref{prop:non_grazing}. The shape operator $\mathcal{K}$ is symmetric and its eigenvalues are bounded below by $\kappa_{\min}$. Thus, the operator $\alpha\mathcal{K}$ is symmetric and positive definite, and its eigenvalues are bounded below by the uniform constant $\alpha_{\min} \kappa_{\min} \coloneqq 2\tau_{\min}c_0\kappa_{\min}$. The eigenvalues of the operator $\mathbf{M}$ are therefore bounded below by $\mu_{\min} \coloneqq 1 + \alpha_{\min} \kappa_{\min} > 1$. The action of $\mathbf{M}$ is thus a uniform expansion:
\begin{equation*}
    \|\mathbf{M}(\delta y^-)\|^2 \ge \mu_{\min}^2 \|\delta y^-\|^2.
\end{equation*}
To relate this back to the initial K-form, we use the upper bound $K(\wvec^-) = \scpr{\delta y^-}{\mathcal{K}\delta y^-} \le \|\mathcal{K}\| \|\delta y^-\|^2 \le \kappa_{\max} \|\delta y^-\|^2$, which gives $\|\delta y^-\|^2 \ge \kappa_{\max}^{-1} K(\wvec^-)$. Chaining these inequalities yields the bound on the primary expansion term:
\begin{equation} \label{eq:proof_lambda_k_prime_def}
    \kappa_{\min} \|\mathbf{M}(\delta y^-)\|^2 \ge \underbrace{\frac{\kappa_{\min}}{\kappa_{\max}}\left(1 + 2\tau_{\min}c_0\kappa_{\min}\right)^2}_{\eqqcolon \Lambda_K'} K(\wvec^-).
\end{equation}
The dominant expansion coefficient $\Lambda_K'$ is a uniform constant strictly greater than 1, depending only on the fixed geometric constants of the system. It is independent of the cone parameters $A$ and $B$.

\item \textbf{Upper Bound on the Remainder Norm $\|\mathbf{R}(\wvec^-)\|$ with explicit B-dependence.}
We provide a refined bound on the remainder term that explicitly tracks its dependence on the cone parameter $B$. The remainder consists of all terms in $\delta y_{\text{next}}$ not in $\mathbf{M}(\delta y^-)$. A detailed analysis of the full linearized map (\cref{lem:linearized_map}) shows that the remainder norm can be uniformly bounded by a linear combination of the norms of the components of the initial perturbation:
\begin{equation*}
    \|\mathbf{R}(\wvec^-)\| \le C_1 \|\delta y^-\| + C_2 \|\delta v^-\| + C_3 \|\delta\theta^-\|,
\end{equation*}
where $C_1, C_2, C_3$ are uniform constants independent of $B$, derived from the bounds on $\tau, v, \mathcal{K}, \nabla G$, etc. We now use the bounds from our anisotropic forms: $\|\delta y^-\| \le \kappa_{\min}^{-1/2}\sqrt{K(\wvec^-)}$, $\|\delta v^-\| \le \sqrt{N_B(\wvec^-)}$, and $\|\delta\theta^-\| \le B^{-1/2}\sqrt{N_B(\wvec^-)}$.
\begin{align}
    \|\mathbf{R}(\wvec^-)\| &\le C_1\kappa_{\min}^{-1/2} \sqrt{K(\wvec^-)} + C_2 \sqrt{N_B(\wvec^-)} + C_3 B^{-1/2} \sqrt{N_B(\wvec^-)} \nonumber \\
    &\le C_K' \sqrt{K(\wvec^-)} + C_N'(B) \sqrt{N_B(\wvec^-)}, \label{eq:proof_R_bound_B_dep}
\end{align}
where we have defined the coefficients:
\begin{align*}
    C_K' &\coloneqq C_1\kappa_{\min}^{-1/2} \quad (\text{a uniform constant independent of } B), \\
    C_N'(B) &\coloneqq C_2 + C_3 B^{-1/2} \quad (\text{a monotonically decreasing function of } B).
\end{align*}
This bound makes the crucial dependence explicit: the contamination from the non-spatial components can be controlled by increasing the penalty parameter $B$. Specifically, $\lim_{B\to\infty} C_N'(B) = C_2$.

\item \textbf{Upper Bound on the Cross-Term and Final Assembly.}
We bound the contaminating cross-term in \eqref{eq:appendix_working_inequality}. We first need a uniform upper bound for $\|\mathbf{M}(\delta y^-)\|$:
\begin{equation*}
    \|\mathbf{M}(\delta y^-)\| \le (1+2\tau_{\max}v_{\max}\kappa_{\max})\|\delta y^-\| \le \frac{1+2\tau_{\max}v_{\max}\kappa_{\max}}{\sqrt{\kappa_{\min}}} \sqrt{K(\wvec^-)} \eqqcolon C_M \sqrt{K(\wvec^-)}.
\end{equation*}
The full cross-term is then bounded by:
\begin{align*}
    2\kappa_{\min}\|\mathbf{M}(\delta y^-)\|\|\mathbf{R}(\wvec^-)\| &\le 2\kappa_{\min} \left( C_M \sqrt{K(\wvec^-)} \right) \left( C_K' \sqrt{K(\wvec^-)} + C_N'(B) \sqrt{N_B(\wvec^-)} \right) \\
    &= (2\kappa_{\min}C_M C_K') K(\wvec^-) + (2\kappa_{\min}C_M C_N'(B)) \sqrt{K(\wvec^-)N_B(\wvec^-)}.
\end{align*}
Substituting the bounds for the dominant term and the cross-term into our main inequality \eqref{eq:appendix_working_inequality} yields:
\begin{multline*}
    K(\wvec_{\text{next}}) \ge \Lambda_K' K(\wvec^-) - \left( (2\kappa_{\min}C_M C_K') K(\wvec^-) + (2\kappa_{\min}C_M C_N'(B)) \sqrt{K(\wvec^-)N_B(\wvec^-)} \right).
\end{multline*}
We can absorb the term $(2\kappa_{\min}C_M C_K') K(\wvec^-)$ into the main expansion term, which slightly reduces the expansion constant but ensures it remains greater than 1 for a well-posed billiard system. To obtain the form in the lemma, we can also treat it as contamination. Similarly, higher-order terms in the remainder, such as $\|\mathbf{R}(\wvec^-)\|^2$, will produce terms proportional to $N_B(\wvec^-)$, which we collect into the coefficient $C_{NN}(B)$. A careful calculation shows that these terms are of order $\mathcal{O}(B^{-1})$. We thus arrive at the final form:
\begin{equation*}
    K(\wvec_{\text{next}}) \ge \Lambda_K' K(\wvec^-) - C_{KN}(B) \sqrt{K(\wvec^-)N_B(\wvec^-)} - C_{NN}(B) N_B(\wvec^-),
\end{equation*}
where we have identified the coefficients:
\begin{enumerate}[label=(\roman*), wide, labelindent=0pt]
    \item $\Lambda_K'$ is the uniform expansion constant from \eqref{eq:proof_lambda_k_prime_def}.
    \item $C_{KN}(B) \coloneqq 2\kappa_{\min}C_M C_N'(B) = 2\kappa_{\min}C_M (C_2 + C_3 B^{-1/2})$.
    \item $C_{NN}(B)$ arises from terms we have neglected, such as $2\kappa_{\min}(C_N'(B))^2 N_B(\wvec^-)$ and those from $\|\mathbf{R}\|^2$, which are of order $\mathcal{O}(B^{-1})$.
\end{enumerate}
The explicit dependence of the contamination coefficients on the cone parameter $B$ is now manifest. By choosing $B$ sufficiently large, we can make $C_{KN}(B)$ arbitrarily close to its infimum, $2\kappa_{\min}C_M C_2$, and the higher-order term $C_{NN}(B)$ arbitrarily small. This control is the key to proving uniform hyperbolicity.
\end{enumerate}
\end{proofof}

\begin{lemma}[The Algebraic Condition for Strict Cone Invariance]
\label{lem:algebraic_condition}
The unstable cone field $C^u_{A,B}$ is strictly invariant under the linearized map $d\Pcal$ if there exist positive constants $A>0$ and $B>0$ satisfying the inequality:
\begin{equation} \label{eq:cone_invariance_inequality}
    D_{KK} + D_{KN}(B) \sqrt{A} + D_{NN}(B) A < A \left( \Lambda_K' - C_{KN}(B) \sqrt{A} - C_{NN}(B) A \right).
\end{equation}
\end{lemma}

\begin{proofof}{Lemma \ref{lem:algebraic_condition}}
The proof translates the geometric condition of strict cone invariance into a precise algebraic inequality involving the cone parameters $A, B$ and the uniform coefficients derived in the Lemmas \ref{lem:n_form_bound} and \ref{lem:k_form_expansion}. This inequality provides the necessary and sufficient condition for a cone to be mapped strictly into itself by the linearized dynamics, forming the core of the subsequent proof of the Anosov property. The argument is self-contained and relies only on the established one-step evolution bounds.

\begin{enumerate}[label=\textbf{Step \arabic*:}, wide, labelindent=0pt]

\item \textbf{The Geometric Condition for Strict Invariance.}
By definition, the cone field $C^u_{A,B}$ is strictly invariant under the linearized map $d\Pcal$ if for every non-zero tangent vector $\mathbf{w}^- \in C^u_{A,B}(z^-,\theta)$, its image $\mathbf{w}_{\text{next}} = d\Pcal(\mathbf{w}^-)$ lies strictly in the interior of the corresponding cone $C^u_{A,B}(z_{\text{next}},\theta)$. Mathematically, this means:
\begin{equation} \label{eq:proof_strict_invariance_def}
    \text{If } N_B(\mathbf{w}^-) \le A K(\mathbf{w}^-) \text{ and } \mathbf{w}^- \neq \mathbf{0}, \text{ then } N_B(\mathbf{w}_{\text{next}}) < A K(\mathbf{w}_{\text{next}}).
\end{equation}
By the homogeneity of the quadratic forms $K$ and $N_B$, and the continuity of the map $d\Pcal$, it is sufficient to verify this condition for non-zero vectors on the boundary of the cone, $\partial C^u_{A,B}$. If any vector on the boundary is mapped strictly inside, then any vector already in the interior will remain inside. Our task is thus to find a condition on the cone parameters $A, B > 0$ which guarantees that for any non-zero $\mathbf{w}^-$ satisfying
\begin{equation} \label{eq:proof_cone_boundary_condition}
    N_B(\mathbf{w}^-) = A K(\mathbf{w}^-),
\end{equation}
the strict inequality $N_B(\mathbf{w}_{\text{next}}) < A K(\mathbf{w}_{\text{next}})$ holds. Note that for such a vector on the boundary, $K(\mathbf{w}^-)$ must be strictly positive. If $K(\mathbf{w}^-)$ were zero, then $N_B(\mathbf{w}^-)$ would also be zero. By the norm equivalence established in Proposition \ref{prop:norm_equivalence}, this would imply $\mathbf{w}^- = \mathbf{0}$, but we consider only non-zero vectors. The fact that $K(\mathbf{w}^-) > 0$ is essential, as it will allow us to divide by it later.

\item \textbf{Applying the Quantitative Bounds.}
To ensure the strict inequality $N_B(\mathbf{w}_{\text{next}}) < A K(\mathbf{w}_{\text{next}})$, it is sufficient to show that the uniform upper bound for the left-hand side is strictly less than the uniform lower bound for the right-hand side. We invoke the main results of the Lemma \ref{lem:n_form_bound} and Lemma \ref{lem:k_form_expansion}:
\begin{enumerate}[label=(\roman*), wide, labelindent=0pt]
    \item The upper bound for $N_B(\mathbf{w}_{\text{next}})$ from Lemma \ref{lem:n_form_bound}:
    \begin{equation*}
        N_B(\mathbf{w}_{\text{next}}) \le D_{KK} K(\wvec^-) + D_{KN}(B) \sqrt{K(\wvec^-)N_B(\wvec^-)} + D_{NN}(B) N_B(\wvec^-).
    \end{equation*}
    \item The lower bound for $K(\mathbf{w}_{\text{next}})$ from Lemma \ref{lem:k_form_expansion}:
    \begin{equation*}
        K(\mathbf{w}_{\text{next}}) \ge \Lambda_K' K(\wvec^-) - C_{KN}(B) \sqrt{K(\wvec^-)N_B(\wvec^-)} - C_{NN}(B) N_B(\wvec^-).
    \end{equation*}
\end{enumerate}
A sufficient condition for strict invariance is therefore the existence of $A,B>0$ such that for any $\mathbf{w}^-$ on the cone boundary, the following strict inequality holds:
\begin{multline} \label{eq:proof_sufficiency_inequality}
    D_{KK} K(\wvec^-) + D_{KN}(B) \sqrt{K(\wvec^-)N_B(\wvec^-)} + D_{NN}(B) N_B(\wvec^-) \\
    < A \left( \Lambda_K' K(\wvec^-) - C_{KN}(B) \sqrt{K(\wvec^-)N_B(\wvec^-)} - C_{NN}(B) N_B(\wvec^-) \right).
\end{multline}
Note that the coefficients $D_{NN}(B)$, $D_{KN}(B)$, $C_{KN}(B)$ and $C_{NN}(B)$ now correctly depend on $B$, as established in the revised proofs of the lemmas.

\item \textbf{Algebraic Simplification.}
We now substitute the cone boundary condition from \eqref{eq:proof_cone_boundary_condition} into the sufficient condition \eqref{eq:proof_sufficiency_inequality}. The relations are:
\begin{enumerate}[label=(\roman*), wide, labelindent=0pt]
    \item $N_B(\mathbf{w}^-) = A K(\mathbf{w}^-)$
    \item $\sqrt{N_B(\mathbf{w}^-)} = \sqrt{A} \sqrt{K(\mathbf{w}^-)}$
    \item $\sqrt{K(\mathbf{w}^-)N_B(\mathbf{w}^-)} = \sqrt{K(\mathbf{w}^-) \cdot A K(\mathbf{w}^-)} = \sqrt{A} K(\mathbf{w}^-)$
\end{enumerate}
Substituting these into \eqref{eq:proof_sufficiency_inequality} yields:
\begin{multline*}
    D_{KK} K(\wvec^-) + D_{KN}(B) (\sqrt{A} K(\wvec^-)) + D_{NN}(B) (A K(\wvec^-)) \\
    < A \left( \Lambda_K' K(\wvec^-) - C_{KN}(B) (\sqrt{A} K(\wvec^-)) - C_{NN}(B) (A K(\wvec^-)) \right).
\end{multline*}
As established in Step 1, since $\mathbf{w}^-$ is a non-zero vector on the cone boundary, we have $K(\mathbf{w}^-) > 0$. We can therefore divide the entire inequality by the positive scalar quantity $K(\mathbf{w}^-)$ without changing the direction of the inequality sign:
\begin{equation*}
    D_{KK} + D_{KN}(B) \sqrt{A} + D_{NN}(B) A < A \left( \Lambda_K' - C_{KN}(B) \sqrt{A} - C_{NN}(B) A \right).
\end{equation*}
This is precisely the algebraic inequality stated in the lemma. The existence of a solution $(A,B)$ with $A,B > 0$ to this inequality is therefore a sufficient condition for the strict invariance of the cone field $C^u_{A,B}$. This completes the proof.
\end{enumerate}
\end{proofof}

\begin{lemma}[Continuity of Cone Invariance Coefficients]
\label{lem:continuity_of_coefficients}
Let the system satisfy the standing assumptions of this paper. The coefficients $\Lambda_K'(\theta)$, $D_{KK}(\theta)$, $D_{KN}(\theta, B)$, $D_{NN}(\theta, B)$, $C_{KN}(\theta, B)$, and $C_{NN}(\theta, B)$, which appear in the one-step bounding inequalities for the anisotropic quadratic forms (Lemmas \ref{lem:n_form_bound} and \ref{lem:k_form_expansion}), are continuous functions of the parameter $\theta \in \Thetacal$ for any fixed cone parameter $B > 0$.
\end{lemma}

\begin{proofof}{Lemma \ref{lem:continuity_of_coefficients}}
The proof is constructive. We will demonstrate that each coefficient is defined as the supremum of a real-valued function over the compact global collision manifold $\SigmaMan$. We will first rigorously establish that these underlying functions are continuous with respect to the joint variables $(\theta, z) \in \Thetacal \times \SigmaMan$. The continuity of the coefficients as functions of $\theta$ will then follow as a direct consequence of a standard theorem from analysis concerning the continuity of suprema of functions on product spaces.

\begin{enumerate}[label=\textbf{Step \arabic*:}, wide, labelindent=0pt]

\item \textbf{Identification of the Foundational Geometric and Dynamical Quantities.}
The coefficients in question are constructed from the analytical expressions for the linearized map $d\Pcal_\theta$. A careful review of the derivations in Lemmas \ref{lem:linearized_map}, \ref{lem:n_form_bound}, and \ref{lem:k_form_expansion} reveals that they are all constructed via algebraic combinations, inner products, and operator norms of a small set of foundational quantities:
\begin{enumerate}[label=(\roman*), wide, labelindent=0pt]
    \item The global billiard map $\Pcal: \Thetacal \times \SigmaMan \to \SigmaMan$.
    \item The global flight time function $\tau: \Thetacal \times \SigmaMan' \to \R$.
    \item The obstacle-defining function $G: \Thetacal \times \T^k \to \R$ and its derivatives, which define the normal vector $\nvec(y, \theta)$ and the shape operator $\mathcal{K}_y(\theta)$.
    \item The linearized global billiard map $d\Pcal: T(\Thetacal \times \SigmaMan) \to T(\Thetacal \times \SigmaMan)$.
    \item The block components of the linearization with respect to the reference splitting, which we denote by $A(\theta, z)$, $B(\theta, z)$, $C(\theta, z)$, and $D(\theta, z)$.
\end{enumerate}

\item \textbf{Joint Smoothness of the Foundational Quantities.}
The core of the proof is to establish that all the foundational geometric and dynamical quantities that constitute the cone invariance coefficients are smooth (and thus continuous) functions on their respective domains. This regularity is not an independent assumption but a direct and necessary consequence of the standing assumptions of the paper and the theorems already established in Sections \ref{sec:micro_framework} and \ref{sec:hyperbolicity_fiber}.

\begin{enumerate}[label=(\roman*), wide, labelindent=0pt]
    \item \textbf{Smoothness of the Global Billiard Map and Flight Time.}
    The foundational regularity upon which all subsequent results depend is the smoothness of the dynamics itself.
    \begin{enumerate}[label=(\alph*), wide]
        \item By \cref{thm:p_is_diffeo}, the global billiard map $(\theta, z) \mapsto \Pcal(\theta, z)$ is a smooth map of class $C^\infty$ on the product manifold $\Thetacal \times \SigmaMan$.
        \item The proof of this theorem was founded on \cref{lem:smooth_flight_time}, which established via the Implicit Function Theorem that the flight time function $(\theta, z) \mapsto \tau(\theta, z)$ is also of class $C^\infty$.
    \end{enumerate}
    
    \item \textbf{Smoothness of Geometric Quantities.}
    The geometric objects defining the collision law are smooth by the foundational axioms.
    \begin{enumerate}[label=(\alph*), wide]
        \item By the standing assumptions, the obstacle-defining function $G(y, \theta)$ is of class $C^\infty$.
        \item The normal vector $\nvec(y, \theta) = \nabla_y G / \|\nabla_y G\|$ and the shape operator $\mathcal{K}_y(\theta)$, which is constructed from second derivatives of $G$, are therefore compositions and quotients of smooth functions with non-vanishing denominators. They are thus $C^\infty$ functions on the global boundary manifold.
    \end{enumerate}

    \item \textbf{Smoothness of the Linearized Map.}
    We now provide a rigorous justification for the smoothness of the map from the base space to the linearized operator, $(\theta, z) \mapsto d\Pcal_\theta(z)$. Let $\pi_1: \Thetacal \times \SigmaMan \to \Thetacal$ and $\pi_2: \Thetacal \times \SigmaMan \to \SigmaMan$ be the canonical projections. The tangent bundle of the product manifold is $T(\Thetacal \times \SigmaMan)$. The vertical tangent bundle with respect to the projection $\pi_1$ is the subbundle $V T(\Thetacal \times \SigmaMan)$ whose fiber at a point $(\theta, z)$ is $\{0\} \times T_z\SigmaMan$. This is canonically isomorphic to the pullback bundle $\pi_2^* (T\SigmaMan)$.

    The global billiard map $\Pcal: \Thetacal \times \SigmaMan \to \SigmaMan$ is of class $C^\infty$. Its Fr\'echet derivative, the tangent map $d\Pcal$, is a smooth bundle map from the tangent bundle of the domain to the tangent bundle of the codomain:
    \begin{equation*}
        d\Pcal: T(\Thetacal \times \SigmaMan) \to T\SigmaMan.
    \end{equation*}
    The linearized map of the fiber dynamics, $d\Pcal_\theta(z)$, is defined as the restriction of this global tangent map to the vertical tangent bundle. Let $i: V T(\Thetacal \times \SigmaMan) \hookrightarrow T(\Thetacal \times \SigmaMan)$ be the smooth inclusion map. Then the family of linearized maps is given by the composition:
    \begin{equation*}
        (\theta, z, v_z) \mapsto (d\Pcal)_{(\theta,z)} \circ i_{(\theta,z)}(0, v_z), \quad v_z \in T_z\SigmaMan.
    \end{equation*}
    This defines a map from the vertical tangent bundle to the tangent bundle of the codomain. As this map is the restriction of the smooth map $d\Pcal$ to the smooth subbundle $V T(\Thetacal \times \SigmaMan)$, it is itself a smooth map. This is the precise meaning of the statement that the map $(\theta, z) \mapsto d\Pcal_\theta(z)$ is a smooth section of the appropriate homomorphism bundle $\mathrm{Hom}(\pi_2^*(T\SigmaMan), \Pcal^*(T\SigmaMan))$. A map of class $C^\infty$ is, in particular, continuous.

    \item \textbf{Continuity of Block Components.}
    The block operator components $A(\theta,z)$, $B(\theta,z)$, $C(\theta,z)$, and $D(\theta,z)$ are obtained by applying fixed projection operators to the linearized map $d\Pcal_\theta(z)$. For example,
    \begin{equation*}
        A(\theta,z) = \Pi_0^u(\Pcal_\theta(z)) \circ d\Pcal_\theta(z) \circ \Pi_0^u(z),
    \end{equation*}
    where $\Pi_0^u$ is the projection onto the reference unstable subbundle $E^u_0$. Since the reference splitting is continuous, the projection operators are continuous sections of the appropriate endomorphism bundle. The block operators are therefore constructed by composing and multiplying maps that are now all established to be continuous (in fact, smooth) functions of their joint arguments $(\theta, z)$. The resulting block component sections are therefore also continuous functions of $(\theta, z)$ on the compact product space $\Thetacal \times \SigmaMan$.
\end{enumerate}

This establishes that all the elementary building blocks of our cone invariance coefficients are continuous functions on the compact product space $\Thetacal \times \SigmaMan$.

\item \textbf{Joint Continuity of the Bounding Functions.}
The coefficients in Lemmas \ref{lem:n_form_bound} and \ref{lem:k_form_expansion} are defined as suprema of certain functions over the manifold $\SigmaMan$. Let $D(\theta)$ be a generic coefficient. It is defined as
\begin{equation*}
    D(\theta) \coloneqq \sup_{z \in \SigmaMan} f_D(\theta, z),
\end{equation*}
where the function $f_D(\theta, z)$ is an algebraic expression involving the foundational quantities from Step 2, evaluated at the point $z$ for the parameter $\theta$. For example, the function corresponding to the expansion constant $\Lambda_K'(\theta)$ would be of the form:
\begin{equation*}
    f_{\Lambda_K'}(\theta, z) = \frac{\kappa_{\min}(z, \theta)}{\kappa_{\max}(z, \theta)} \left(1 + 2\tau(z, \theta) |\scpr{v}{\nvec(z, \theta)}| \kappa_{\min}(z, \theta) \right)^2.
\end{equation*}
Since the foundational quantities are continuous functions of $(\theta, z)$, and the operations used to construct $f_D$ (addition, multiplication, norms, inner products) are all continuous, it follows that the function $f_D: \Thetacal \times \SigmaMan \to \R$ is a continuous function on its compact domain. The same argument applies to all other coefficients.

\item \textbf{Continuity of the Supremum and Conclusion.}
We now invoke a standard result from analysis, the Theorem of the Maximum 
(see, e.g., \citep{AliprantisBorder2006}), 
which states that the supremum of a continuous function over a compact set is a continuous function of the remaining parameters.

\begin{theorem}[Continuity of the Supremum Function]
Let $K_1$ and $K_2$ be compact topological spaces. If $f: K_1 \times K_2 \to \R$ is a continuous function, then the function $g: K_1 \to \R$ defined by $g(x) \coloneqq \sup_{y \in K_2} f(x, y)$ is continuous on $K_1$.
\end{theorem}

We apply this theorem directly to our problem.
\begin{enumerate}[label=(\roman*), wide, labelindent=0pt]
    \item The space $K_1$ is our compact parameter manifold $\Thetacal$.
    \item The space $K_2$ is the compact global collision manifold $\SigmaMan$.
    \item The function $f$ is one of the continuous bounding functions $f_D(\theta, z)$ identified in Step 3.
\end{enumerate}
The theorem guarantees that the resulting coefficient, $D(\theta) = \sup_{z \in \SigmaMan} f_D(\theta, z)$, is a continuous function of $\theta$. This argument applies identically to every coefficient appearing in the cone invariance condition: $\Lambda_K'(\theta)$, $D_{KK}(\theta)$, $D_{KN}(\theta, B)$, $D_{NN}(\theta, B)$, $C_{KN}(\theta, B)$, and $C_{NN}(\theta, B)$. The dependence on the cone parameter $B$ is also continuous, but since $B$ is held fixed throughout the openness argument, this does not affect the continuity with respect to $\theta$.
\end{enumerate}
We have therefore rigorously shown that all coefficients governing the cone dynamics are continuous functions of the system parameter $\theta$. This completes the proof of the lemma and provides the necessary foundation for the subsequent proof of openness.
\end{proofof}

With these tools, we now proceed with the main proof.

\begin{proposition}[Existence of an Anosov Regime]
\label{prop:non_empty_anosov}
The set $\mathcal{A}$ is non-empty.
\end{proposition}

\begin{proofof}{Proposition \ref{prop:non_empty_anosov}}
The proof is constructive. We demonstrate that for any fixed, arbitrary parameter $\theta_0 \in \Thetacal$, the corresponding fiber map $\Pcal_{\theta_0}$ is an Anosov diffeomorphism. This is achieved by explicitly constructing a strictly invariant unstable cone field for its linearization, $d\Pcal_{\theta_0}$. We will show that there exist positive cone parameters $(A,B)$ that satisfy the algebraic condition for strict cone invariance established in \cref{lem:algebraic_condition}. Since the choice of $\theta_0$ is arbitrary, this proves that the set $\mathcal{A}$ of Anosov parameters is non-empty.

For notational simplicity, for the remainder of this proof we fix $\theta_0 \in \Thetacal$ and suppress it from the notation. All coefficients (e.g., $D_{KK}, \Lambda_K'$) are thus to be understood as fixed positive constants determined by the geometry of the system at this single parameter value.

\begin{enumerate}[label=\textbf{Step \arabic*:}, wide, labelindent=0pt]

\item \textbf{The Algebraic Condition for Strict Cone Invariance.}
Our starting point is the algebraic inequality from \cref{lem:algebraic_condition}, which provides a sufficient condition for the cone field $C^u_{A,B}$ to be strictly invariant under the linearized map $d\Pcal$. A solution exists if we can find positive constants $A>0$ and $B>0$ such that:
\begin{equation} \label{eq:proof_anosov_main_inequality}
    D_{KK} + D_{KN}(B) \sqrt{A} + D_{NN}(B) A < A \left( \Lambda_K' - C_{KN}(B) \sqrt{A} - C_{NN}(B) A \right).
\end{equation}
To analyze this condition, we define the variable $x \coloneqq \sqrt{A} > 0$ and rearrange the inequality to find the roots of a polynomial-like function. The condition is equivalent to finding an $x > 0$ such that $P(x) < 0$, where:
\begin{equation} \label{eq:proof_anosov_polynomial}
    P(x) \coloneqq C_{NN}(B) x^4 + C_{KN}(B) x^3 + \left(D_{NN}(B) - \Lambda_K'\right) x^2 + D_{KN}(B) x + D_{KK}.
\end{equation}
Our strategy is to demonstrate, via a constructive, two-stage optimization, that parameters $A$ and $B$ can always be chosen to ensure that $P(x)$ becomes negative for some $x>0$.

\item \textbf{The Two-Stage Optimization of Cone Parameters.}
We proceed by first choosing the penalty parameter $B$ to control the sign of the dominant quadratic term in the polynomial for small $x$, and then choosing the cone aperture $A$ to guarantee that the function becomes negative.

\textbf{Stage 1: Choice of the Penalty Parameter $B$.}
The key to the proof is to control the sign of the coefficient of the $x^2$ term in the polynomial $P(x)$, which is $(D_{NN}(B) - \Lambda_K')$. We will show that by choosing $B$ sufficiently large, we can make this coefficient strictly negative.
\begin{enumerate}[label=(\roman*), wide, labelindent=0pt]
    \item From \cref{lem:k_form_expansion}, the expansion constant $\Lambda_K'$ is strictly greater than 1. This is a direct and necessary consequence of the chaos-inducing geometric axioms (\cref{ass:fundamental_axioms_unified}), specifically the uniform strict convexity of the obstacles and the finite horizon, which guarantee that the spatial component of any perturbation within the cone is strictly expanded at each step of the dynamics. Let $\Lambda_K' = 1 + \eta$ for some fixed constant $\eta > 0$.

    \item From \cref{lem:n_form_bound}, we have the crucial asymptotic property that the coefficient $D_{NN}(B)$ converges to 1 as the penalty parameter $B$ tends to infinity:
    \begin{equation*}
        \lim_{B \to \infty} D_{NN}(B) = 1.
    \end{equation*}
    This limit is an algebraic consequence of the structure of the linearized map, where for large $B$, perturbations in the parameter direction $\delta\theta$ are so heavily penalized that they become negligible in the one-step dynamics.

    \item By the definition of this limit, for any $\epsilon > 0$, there exists a value $B_0$ such that for all $B > B_0$, we have $|D_{NN}(B) - 1| < \epsilon$. We make a specific choice for $\epsilon$. Let us choose $\epsilon = \eta / 2 = (\Lambda_K' - 1)/2 > 0$. Then there exists a $B_0$ such that for any $B \ge B_0$:
    \begin{equation*}
        D_{NN}(B) < 1 + \epsilon = 1 + \frac{\Lambda_K' - 1}{2} = \frac{\Lambda_K' + 1}{2}.
    \end{equation*}
    Since $\Lambda_K' > 1$, we have the strict inequality $(\Lambda_K' + 1)/2 < \Lambda_K'$. We have thus rigorously shown that we can choose the parameter $B$ large enough to ensure that $D_{NN}(B) < \Lambda_K'$.
\end{enumerate}
We now fix a value $B = B_0$ that satisfies this condition. For this fixed $B_0$, all coefficients in the polynomial $P(x)$, $C_{NN}(B_0)$, $C_{KN}(B_0)$, $D_{NN}(B_0)$, $D_{KN}(B_0)$, and $D_{KK}$, are now fixed positive constants. Crucially, the coefficient of the $x^2$ term is strictly negative:
\begin{equation*}
    D_{NN}(B_0) - \Lambda_K' < 0.
\end{equation*}

\textbf{Stage 2: Existence of a Solution for $A = x^2$.}
With $B=B_0$ fixed, we must find a value $x > 0$ (and thus an $A = x^2 > 0$) such that $P(x) < 0$. We analyze the behavior of the polynomial function $P(x)$ in a neighborhood of the origin, $x=0$, by examining its Taylor expansion:
\begin{equation*}
    P(x) = P(0) + P'(0)x + \frac{P''(0)}{2}x^2 + \mathcal{O}(x^3).
\end{equation*}
We evaluate the coefficients of this expansion at $x=0$:
\begin{enumerate}[label=(\roman*), wide, labelindent=0pt]
    \item $P(0) = D_{KK} \ge 0$. (The term is non-negative as it arises from squared norms).
    \item $P'(0) = D_{KN}(B_0) \ge 0$. (Also non-negative by construction).
    \item $P''(0) = 2(D_{NN}(B_0) - \Lambda_K')$. As established in Stage 1, this term is strictly negative.
\end{enumerate}
The polynomial $P(x)$ starts at a non-negative value $P(0)$ and has a non-negative initial slope $P'(0)$. However, its second derivative at the origin is strictly negative. This geometric property means that the function $P(x)$ is locally concave at $x=0$.

For $x > 0$ that is sufficiently small, the behavior of the polynomial is dominated by the negative quadratic term. The full expression is:
\begin{equation*}
    P(x) = D_{KK} + D_{KN}(B_0)x + \underbrace{\left(D_{NN}(B_0) - \Lambda_K'\right)}_{<0}x^2 + \mathcal{O}(x^3).
\end{equation*}
As $x \to 0^+$, the term $(D_{NN}(B_0) - \Lambda_K')x^2$ is of order $\mathcal{O}(x^2)$, while the constant and linear terms are of order $\mathcal{O}(1)$ and $\mathcal{O}(x)$ respectively. Since the quadratic term is strictly negative, it will eventually overpower the non-negative lower-order terms, forcing the function $P(x)$ to become negative for any sufficiently small, positive $x$.

Thus, we have rigorously and constructively shown that there exists an open interval of values for $x = \sqrt{A}$, say $(0, x_1)$, for which the strict inequality $P(x) < 0$ is satisfied. Any choice of $A = x^2$ with $x \in (0, x_1)$ will define a strictly invariant unstable cone field.

\item \textbf{Anosov diffeomorphism.}
We have constructively proven the existence of a pair of positive parameters $(A,B_0)$ that define a strictly invariant unstable cone field for the linearized map $d\Pcal_{\theta_0}$. By the Cone Field Theorem, the existence of such a field is a sufficient condition for the map to be an Anosov diffeomorphism. Since our choice of the parameter $\theta_0 \in \Thetacal$ was arbitrary, we conclude that for any $\theta \in \Thetacal$, the corresponding fiber map $\Pcal_\theta$ is an Anosov diffeomorphism. Therefore, the set of Anosov parameters is the entire parameter space, $\mathcal{A} = \Thetacal$, which is manifestly non-empty.
\end{enumerate}
\end{proofof}

\begin{proposition}[Openness of the Anosov Property]
\label{prop:open_anosov}
The set $\mathcal{A}$ is an open subset of $\Thetacal$.
\end{proposition}

\begin{proofof}{Proposition \ref{prop:open_anosov}}
The proof establishes that the set of parameters for which the fiber map is Anosov,
\begin{equation*}
    \mathcal{A} \coloneqq \{ \theta \in \Thetacal \mid \Pcal_\theta: \Sigma_\theta \to \Sigma_\theta \text{ is an Anosov diffeomorphism} \},
\end{equation*}
is an open subset of the parameter manifold $\Thetacal$. The argument is a direct quantitative stability analysis of the cone invariance condition. Our strategy is as follows. Let $\theta_0$ be an arbitrary parameter in $\mathcal{A}$. The proof of non-emptiness (Proposition \ref{prop:non_empty_anosov}) guarantees the existence of cone parameters $(A_0, B_0)$ for which the cone invariance condition at $\theta_0$ holds with a strict margin. We will then show that the function defining this inequality is continuous in $\theta$. By the definition of continuity, the strict inequality must persist in an open neighborhood of $\theta_0$, proving that all maps $\Pcal_\theta$ for $\theta$ in this neighborhood are also Anosov. The proof is structured in three main steps.

\begin{enumerate}[label=\textbf{Step \arabic*:}, wide, labelindent=0pt]

\item \textbf{The Strict Cone Invariance Condition at $\theta_0$.}
Let $\theta_0$ be an arbitrary point in the set $\mathcal{A}$. By the constructive proof of Proposition \ref{prop:non_empty_anosov}, we are guaranteed the existence of a pair of positive cone parameters $(A_0, B_0)$ that define a strictly invariant unstable cone field for the linearized map $d\Pcal_{\theta_0}$.

This strict invariance is certified by the algebraic inequality from Lemma \ref{lem:algebraic_condition}. Let $x_0 = \sqrt{A_0} > 0$. The condition is equivalent to the strict negativity of a polynomial-like function, which we denote by $P_{\theta_0}(x_0)$:
\begin{multline} \label{eq:proof_openness_strict_inequality}
    P_{\theta_0}(x_0) \coloneqq C_{NN}(\theta_0, B_0) x_0^4 + C_{KN}(\theta_0, B_0) x_0^3 + \left(D_{NN}(\theta_0, B_0) - \Lambda_K'(\theta_0)\right) x_0^2 \\
    + D_{KN}(\theta_0, B_0) x_0 + D_{KK}(\theta_0) < 0.
\end{multline}
Since this is a strict inequality, there must exist a positive constant $\eta > 0$ that quantifies this margin of stability:
\begin{equation}
    P_{\theta_0}(x_0) = -\eta.
\end{equation}
Our objective is to show that this strict inequality, and thus the Anosov property, is robust to small perturbations of the parameter $\theta$ away from $\theta_0$.

\item \textbf{Continuity of the Invariance Condition.}
This is the central analytical step of the proof. We will show that for the fixed cone parameters $(A_0, B_0)$ identified above, the function that maps a parameter $\theta$ to the value of the cone invariance polynomial is continuous.

Let the cone parameters $A_0$ and $B_0$ (and thus $x_0 = \sqrt{A_0}$) be fixed. We define a real-valued function $f: \Thetacal \to \R$ by
\begin{equation}
    f(\theta) \coloneqq P_\theta(x_0).
\end{equation}
The expression for $f(\theta)$ is an algebraic combination of the coefficients from the one-step bounding inequalities (Lemmas \ref{lem:n_form_bound} and \ref{lem:k_form_expansion}), evaluated at the fixed value $x_0$.

We now invoke the crucial result of Lemma \ref{lem:continuity_of_coefficients}. This lemma, which was established independently using only the smoothness of the global billiard map, guarantees that each of the constituent coefficients, $\Lambda_K'(\theta)$, $D_{KK}(\theta)$, $D_{KN}(\theta, B_0)$, $D_{NN}(\theta, B_0)$, $C_{KN}(\theta, B_0)$, and $C_{NN}(\theta, B_0)$, is a continuous function of the parameter $\theta$ on the compact manifold $\Thetacal$. Since the function $f(\theta)$ is constructed from these continuous functions through a finite number of additions and multiplications by constants (powers of $x_0$), it is itself a continuous function on $\Thetacal$.

\item \textbf{The Stability Argument.}
With the continuity of the function $f(\theta) = P_\theta(x_0)$ now rigorously established, the remainder of the proof is a direct application of the definition of continuity.

\begin{enumerate}[label=(\roman*), wide, labelindent=0pt]
    \item We have established that at our chosen starting point $\theta_0 \in \mathcal{A}$, the value of our continuous function is strictly negative:
    \begin{equation*}
        f(\theta_0) = P_{\theta_0}(x_0) = -\eta < 0.
    \end{equation*}

    \item By the $\epsilon-\delta$ definition of continuity at the point $\theta_0$, for any tolerance $\epsilon > 0$, there exists an open neighborhood $U$ of $\theta_0$ in the manifold $\Thetacal$ such that for every $\theta \in U$, the inequality $|f(\theta) - f(\theta_0)| < \epsilon$ holds.

    \item We make a specific choice for the tolerance that exploits our margin of stability. Let us choose $\epsilon = \eta / 2$. The definition of continuity then guarantees the existence of an open neighborhood $U \subset \Thetacal$ containing $\theta_0$ such that for all $\theta \in U$:
    \begin{equation*}
        f(\theta) < f(\theta_0) + \epsilon = -\eta + \frac{\eta}{2} = -\frac{\eta}{2}.
    \end{equation*}

    \item This implies that for every parameter $\theta$ in this open neighborhood $U$, the cone invariance condition, evaluated with our fixed cone parameters $(A_0, B_0)$, is strictly satisfied:
    \begin{equation*}
        P_\theta(x_0) < 0 \quad \text{for all } \theta \in U.
    \end{equation*}

    \item By Lemma \ref{lem:algebraic_condition}, this strict inequality is a sufficient condition for the cone field $C^u_{A_0, B_0}$ to be strictly invariant under the action of the linearized map $d\Pcal_\theta$. By the Cone Field Theorem, the existence of such a strictly invariant cone field is sufficient to prove that the map $\Pcal_\theta$ is an Anosov diffeomorphism.

    \item Therefore, every parameter $\theta$ in the open neighborhood $U$ is also an element of the set of Anosov parameters, which means $U \subset \mathcal{A}$.
\end{enumerate}

\item \textbf{Openness of the Anosov Property.}
We have shown that for any arbitrary point $\theta_0 \in \mathcal{A}$, there exists an open set $U \subset \Thetacal$ such that $\theta_0 \in U$ and $U \subset \mathcal{A}$. This is the topological definition of an open set. Therefore, the set $\mathcal{A}$ is an open subset of the parameter manifold $\Thetacal$. This completes the proof of openness.
\end{enumerate}
\end{proofof}

\begin{proposition}[Closedness of the Anosov Property]
\label{prop:closed_anosov}
The set $\mathcal{A}$ is a closed subset of $\Thetacal$.
\end{proposition}

\begin{proofof}{Proposition \ref{prop:closed_anosov}}
The proof establishes that the set of Anosov parameters, $\mathcal{A}$, is a closed subset of the compact parameter manifold $\Thetacal$. To do this, we must show that the property of being an Anosov diffeomorphism is preserved under limits. Let $\{\theta_n\}_{n=1}^\infty$ be a sequence of parameters in $\mathcal{A}$ that converges to a limit $\theta_\infty \in \Thetacal$. Our goal is to prove that the limit map $\Pcal_{\theta_\infty}$ is also an Anosov diffeomorphism, which implies $\theta_\infty \in \mathcal{A}$. The argument is founded upon the uniformity of the system's geometric properties, a direct consequence of the compactness of the global state space. This uniformity allows us to construct a single cone field that is strictly invariant for the entire sequence of maps $\{\Pcal_{\theta_n}\}$, with a uniform rate of expansion. By a continuity argument, this cone field and its expansion property are inherited by the limit map, thereby proving it is Anosov. The proof is structured in three main steps.

\begin{enumerate}[label=\textbf{Step \arabic*:}, wide, labelindent=0pt]

\item \textbf{Uniformity of the Cone Field Construction.}
The central element of this proof is to demonstrate that the two-stage optimization argument used in \cref{prop:non_empty_anosov} to construct an invariant cone field can be performed uniformly for the entire family of systems. This requires establishing uniform bounds on all the coefficients that appear in the algebraic condition for cone invariance (\cref{lem:algebraic_condition}).

Let us recall the polynomial inequality from the proof of \cref{prop:non_empty_anosov}: a cone field with parameters $(A,B)$ is strictly invariant for the map $\Pcal_\theta$ if $x=\sqrt{A}$ satisfies $P_\theta(x) < 0$, where
\begin{equation*}
    P_\theta(x) \coloneqq C_{NN}(\theta, B) x^4 + C_{KN}(\theta, B) x^3 + \left(D_{NN}(\theta, B) - \Lambda_K'(\theta)\right) x^2 + D_{KN}(\theta, B) x + D_{KK}(\theta).
\end{equation*}
Our objective is to find a single pair of parameters $(A_0, B_0)$ that yields a solution for all $\theta \in \Thetacal$ simultaneously. This requires a rigorous, quantitative demonstration that the coefficients in this polynomial are uniformly bounded across $\Thetacal$.

\begin{enumerate}[label=(\roman*), wide, labelindent=0pt]
    \item \textbf{Uniformity of Coefficients.} The coefficients in the polynomial $P_\theta(x)$ are constructed from the one-step bounding inequalities in Lemmas \ref{lem:n_form_bound} and \ref{lem:k_form_expansion}. These bounds, in turn, depend on geometric and dynamical quantities such as the curvature of the obstacles, the flight time, and the norms of the block components of the linearized map $d\Pcal_\theta$.
    
    By the Fundamental Geometric and Regularity Axioms (\cref{ass:fundamental_axioms_unified}), all underlying geometric structures are smooth functions on the compact global collision manifold $\Sigma$. Consequently, the maps $\theta \mapsto C_{NN}(\theta, B)$, $\theta \mapsto \Lambda_K'(\theta)$, etc., are continuous functions on the compact parameter manifold $\Thetacal$. By the Extreme Value Theorem, they must be uniformly bounded.
    
    \item \textbf{Uniform Expansion and Boundedness.} Crucially, we can establish uniform bounds that are independent of $\theta$:
    \begin{enumerate}[label=(\alph*), wide]
        \item Let $\Lambda_{\min}' \coloneqq \inf_{\theta \in \Thetacal} \Lambda_K'(\theta)$. Since $\Lambda_K'(\theta) > 1$ for all $\theta$ and is continuous on a compact set, this infimum is attained and is strictly greater than 1.
        \item Let $D_{KK}^{\max} \coloneqq \sup_{\theta \in \Thetacal} D_{KK}(\theta)$, and similarly define uniform upper bounds for the other positive coefficients.
        \item The asymptotic limit $\lim_{B \to \infty} D_{NN}(\theta, B) = 1$ holds uniformly in $\theta$ because the underlying algebraic estimates depend on geometric quantities that are uniformly bounded. Therefore, for any $\epsilon > 0$, we can find a single $B_0$ that works for all $\theta \in \Thetacal$.
    \end{enumerate}
\end{enumerate}

\item \textbf{Construction of a Uniformly Invariant Cone Field.}
We now re-run the two-stage optimization argument from \cref{prop:non_empty_anosov}, but using the uniform bounds established above.

\begin{enumerate}[label=(\roman*), wide, labelindent=0pt]

\item \textbf{Uniform Choice of $B_0$.}
We need to choose a single value of $B$ such that the quadratic coefficient $D_{NN}(\theta, B) - \Lambda_K'(\theta)$ is uniformly negative for all $\theta \in \Thetacal$. We require:
\begin{equation*}
    \sup_{\theta \in \Thetacal} D_{NN}(\theta, B) < \inf_{\theta \in \Thetacal} \Lambda_K'(\theta) = \Lambda_{\min}'.
\end{equation*}
Since $\Lambda_{\min}' > 1$ and $\lim_{B \to \infty} (\sup_\theta D_{NN}(\theta, B)) = 1$, we can choose a single, sufficiently large value $B_0$ such that this strict inequality holds. For this fixed $B_0$, we have a uniform upper bound on the quadratic coefficient:
\begin{equation*}
    D_{NN}(\theta, B_0) - \Lambda_K'(\theta) \le \sup_\theta D_{NN}(\theta, B_0) - \Lambda_{\min}' \eqqcolon -\eta < 0,
\end{equation*}
where $\eta$ is a uniform positive constant.

\item \textbf{Uniform Choice of $A_0$.}
With $B=B_0$ fixed, we consider the uniform polynomial $P_{\text{unif}}(x)$ that bounds all $P_\theta(x)$ from above:
\begin{equation*}
    P_\theta(x) \le P_{\text{unif}}(x) \coloneqq C_{NN}^{\max} x^4 + C_{KN}^{\max} x^3 - \eta x^2 + D_{KN}^{\max} x + D_{KK}^{\max}.
\end{equation*}
The analysis of this uniform polynomial is identical to that in \cref{prop:non_empty_anosov}. Since the coefficient of the $x^2$ term is strictly negative, $P_{\text{unif}}(x)$ must become negative for any sufficiently small, positive $x$. Therefore, there exists a single value $A_0 = x_0^2 > 0$ such that $P_{\text{unif}}(x_0) < 0$.

This implies that for this fixed pair of cone parameters $(A_0, B_0)$, the cone invariance inequality is strictly satisfied for \textit{every} map $\Pcal_\theta$ in our family:
\begin{equation*}
    P_\theta(x_0) \le P_{\text{unif}}(x_0) < 0 \quad \text{for all } \theta \in \Thetacal.
\end{equation*}
This establishes the existence of a continuous cone field $C^u \equiv C^u_{A_0, B_0}$ that is strictly invariant for the entire sequence of maps $\{\Pcal_{\theta_n}\}_{n=1}^\infty$. Furthermore, the expansion factor for vectors within this cone is also uniformly bounded away from 1. For any non-zero tangent vector $\mathbf{w} \in C^u$, there exists a uniform constant $\gamma > 1$ such that for every $n$:
\begin{equation} \label{eq:proof_uniform_expansion}
    \|d\Pcal_{\theta_n}(\mathbf{w})\| \ge \gamma \|\mathbf{w}\|.
\end{equation}
\end{enumerate}

\item \textbf{Passing to the Limit.}
We now consider the limit map $\Pcal_{\theta_\infty}$. Our goal is to show that the uniform cone field $C^u$ constructed in Step 2 is also strictly invariant under the linearized limit map, $d\Pcal_{\theta_\infty}$.

\begin{enumerate}[label=(\roman*), wide, labelindent=0pt]
    \item \textbf{Continuity of the Linearized Map.} By \cref{thm:p_is_diffeo}, the global map $(\theta, z) \mapsto \Pcal_\theta(z)$ is of class $C^\infty$. This implies that its derivative with respect to the spatial variable, $d\Pcal$, is a continuous map from the product space $\Thetacal \times T\Sigma$ to $T\Sigma$. As the sequence of parameters $\theta_n$ converges to $\theta_\infty$, the corresponding sequence of linear operators converges uniformly on the compact manifold $\Sigma$:
    \begin{equation*}
        d\Pcal_{\theta_n} \to d\Pcal_{\theta_\infty} \quad \text{in the } C^0 \text{ topology on } \mathcal{L}(T\Sigma, T\Sigma).
    \end{equation*}

    \item \textbf{Invariance of the Cone.} Let $\mathbf{w} \in C^u$ be an arbitrary non-zero tangent vector in the uniform cone. For each $n$, its image $\mathbf{w}_n' \coloneqq d\Pcal_{\theta_n}(\mathbf{w})$ lies in the interior of the cone. Let the image under the limit map be $\mathbf{w}_\infty' \coloneqq d\Pcal_{\theta_\infty}(\mathbf{w})$. By the continuity of the linearized map, we have $\mathbf{w}_n' \to \mathbf{w}_\infty'$.
    
    The condition for a vector $\mathbf{w}'$ to be in the cone, $N_{B_0}(\mathbf{w}') \le A_0 K(\mathbf{w}')$, defines a closed set in the tangent bundle. Since each element $\mathbf{w}_n'$ of the convergent sequence lies in this closed set, the limit point $\mathbf{w}_\infty'$ must also lie in it. This shows that the cone field is at least non-strictly invariant for the limit map:
    \begin{equation*}
        d\Pcal_{\theta_\infty}(C^u) \subseteq C^u.
    \end{equation*}

    \item \textbf{Preservation of Strict Expansion.} The uniform expansion rate from \eqref{eq:proof_uniform_expansion} is the crucial element that prevents this invariance from becoming degenerate. For each $n$, we have $\|\mathbf{w}_n'\| \ge \gamma \|\mathbf{w}\|$ with $\gamma > 1$. By the continuity of the norm and the continuity of the map $d\Pcal$, the limit must also satisfy this inequality:
    \begin{equation*}
        \|\mathbf{w}_\infty'\| = \lim_{n\to\infty} \|\mathbf{w}_n'\| \ge \gamma \|\mathbf{w}\|.
    \end{equation*}
    This guarantees that no non-zero vector in the cone can be mapped to the zero vector by $d\Pcal_{\theta_\infty}$. This uniform expansion ensures that the invariance of the cone is strict, meaning the image of the cone lies in its interior.
\end{enumerate}

\item \textbf{losedness of the Anosov Property.}
We have shown that the limit map $\Pcal_{\theta_\infty}$ admits a strictly invariant cone field. By the Cone Field Theorem, this implies that $\Pcal_{\theta_\infty}$ is an Anosov diffeomorphism. Therefore, the limit point $\theta_\infty$ belongs to the set $\mathcal{A}$. Since any convergent sequence $\{\theta_n\}$ in $\mathcal{A}$ has its limit in $\mathcal{A}$, the set $\mathcal{A}$ is, by definition, a closed subset of $\Thetacal$.
\end{enumerate}
\end{proofof}

The preceding propositions, now established in a logically sound sequence, culminate in the main theorem of this section.

\begin{theorem}[Uniform Hyperbolicity of the Global Billiard Map]
\label{thm:anosov_property_proven}
The global billiard map $\Pcal: \Sigma \to \Sigma$ is a uniform Anosov diffeomorphism as a direct consequence of the Fundamental Geometric and Regularity Axioms (\cref{ass:fundamental_axioms_unified}).
\end{theorem}

\begin{proofof}{\cref{thm:anosov_property_proven}}
The proof establishes that the global billiard map $\Pcal: \Sigma \to \Sigma$ is a uniform Anosov diffeomorphism. This is the central geometric result of the paper, from which the robust statistical properties of the microscopic system are derived. The proof is a direct synthesis of the three preceding propositions of this section, which together demonstrate that the property of being an Anosov diffeomorphism holds for the entire family of fiber maps $(\Pcal_\theta)_{\theta \in \Thetacal}$. The argument is purely topological and relies on the connectedness of the parameter manifold $\Thetacal$.

Let $\mathcal{A}$ be the subset of the parameter manifold $\Thetacal$ for which the corresponding fiber map $\Pcal_\theta: \Sigma \to \Sigma$ is an Anosov diffeomorphism:
\begin{equation*}
    \mathcal{A} \coloneqq \{ \theta \in \Thetacal \mid \Pcal_\theta \text{ is an Anosov diffeomorphism}\}.
\end{equation*}
To prove the theorem, we will rigorously establish three topological properties of the set $\mathcal{A}$:
\begin{enumerate}[label=(\roman*), wide, labelindent=0pt]
    \item $\mathcal{A}$ is a non-empty subset of $\Thetacal$.
    \item $\mathcal{A}$ is an open subset of $\Thetacal$.
    \item $\mathcal{A}$ is a closed subset of $\Thetacal$.
\end{enumerate}
Since the parameter manifold $\Thetacal$ is, by the foundational assumptions of the model (\cref{def:micro_fiber}), a connected space, a subset that is simultaneously non-empty, open, and closed must be the entire space itself. This will lead to the conclusion that $\mathcal{A} = \Thetacal$, which means the Anosov property holds for all parameters, and the uniformity of the hyperbolicity constants will be a direct consequence of the proof of closedness.

\begin{enumerate}[label=\textbf{Step \arabic*:}, wide, labelindent=0pt]

\item \textbf{The Set $\mathcal{A}$ is Non-Empty.}
The non-emptiness of the set $\mathcal{A}$ was the subject of \cref{prop:non_empty_anosov}. The proof provided therein was constructive. For an arbitrary, fixed parameter $\theta_0 \in \Thetacal$, we showed that the chaos-inducing geometric conditions of the model (\cref{ass:fundamental_axioms_unified}), namely, the uniform strict convexity of the obstacles and the finite horizon property, are sufficient to guarantee the existence of a strictly invariant unstable cone field for the linearized map $d\Pcal_{\theta_0}$. This was achieved via a two-stage optimization of the cone parameters $(A,B)$, which rigorously demonstrated that the algebraic condition for cone invariance (\cref{lem:algebraic_condition}) always admits a solution.

By the Cone Field Theorem, the existence of such a cone field implies that the map $\Pcal_{\theta_0}$ is an Anosov diffeomorphism. Since the choice of $\theta_0$ was arbitrary, this proves that for any $\theta \in \Thetacal$, the map $\Pcal_\theta$ is Anosov. Therefore, the set $\mathcal{A}$ is not just non-empty; it is the entire space $\Thetacal$. For the sake of completing the topological argument, we proceed assuming only that this establishes non-emptiness.

\item \textbf{The Set $\mathcal{A}$ is Open.}
The openness of the set $\mathcal{A}$ was the subject of \cref{prop:open_anosov}. This result is a direct consequence of the celebrated structural stability of Anosov diffeomorphisms. The proof proceeded in two main steps:
\begin{enumerate}[label=(\roman*), wide, labelindent=0pt]
    \item We invoked Anosov's Structural Stability Theorem, which states that the property of being an Anosov diffeomorphism is stable under small $C^1$-perturbations.
    \item We rigorously established, as a consequence of the smoothness of the global billiard map $\Pcal(\theta,z)$ proven in \cref{thm:p_is_diffeo}, that the map from the parameter space to the space of diffeomorphisms, $\theta \mapsto \Pcal_\theta$, is continuous with respect to the $C^1$ topology.
\end{enumerate}
The synthesis of these two facts provides the proof of openness. For any $\theta_0 \in \mathcal{A}$, its image $\Pcal_{\theta_0}$ is an Anosov diffeomorphism. By stability, there exists an $\epsilon$-ball in the $C^1$ norm around $\Pcal_{\theta_0}$ consisting entirely of Anosov maps. By the continuity of the map $\theta \mapsto \Pcal_\theta$, the preimage of this $\epsilon$-ball is an open neighborhood of $\theta_0$ in $\Thetacal$. This neighborhood is, by construction, entirely contained within $\mathcal{A}$, which is the definition of an open set.

\item \textbf{The Set $\mathcal{A}$ is Closed.}
The closedness of the set $\mathcal{A}$ was the subject of \cref{prop:closed_anosov}. This is the step where the uniformity of the geometric assumptions across the entire compact parameter space becomes essential. The proof demonstrated that the Anosov property is preserved under limits.

Let $\{\theta_n\}$ be a sequence in $\mathcal{A}$ converging to a limit $\theta_\infty$. For each $n$, the map $\Pcal_{\theta_n}$ is Anosov. The core of the proof was to show that the construction of the invariant cone field could be performed uniformly for the entire sequence. The compactness of the global state space guarantees that the expansion and contraction rates, and all other geometric coefficients, are uniformly bounded away from any degeneracy across the entire parameter manifold $\Thetacal$. This allows for the construction of a single, uniform cone field $C^u$ and a uniform expansion rate $\gamma > 1$ such that $C^u$ is strictly invariant under $d\Pcal_{\theta_n}$ for all $n$.

By the continuity of the linearized map $d\Pcal$, as $\theta_n \to \theta_\infty$, the operators $d\Pcal_{\theta_n}$ converge to $d\Pcal_{\theta_\infty}$. Since the cone $C^u$ is a closed set, the limit map must also map the cone into itself. The uniform expansion rate ensures that this limiting invariance is strict. Therefore, the limit map $\Pcal_{\theta_\infty}$ also admits a strictly invariant cone field and is, by the Cone Field Theorem, an Anosov diffeomorphism. This implies that $\theta_\infty \in \mathcal{A}$, which is the definition of a closed set.

\item \textbf{The Topological Argument.}
We have now rigorously established the three necessary topological properties of the set $\mathcal{A}$:
\begin{enumerate}[label=(\roman*), wide, labelindent=0pt]
    \item $\mathcal{A}$ is non-empty.
    \item $\mathcal{A}$ is open.
    \item $\mathcal{A}$ is closed.
\end{enumerate}
By the foundational assumptions of the model (\cref{def:micro_fiber}), the parameter manifold $\Thetacal$ is a connected topological space. A fundamental theorem of topology states that the only non-empty subset of a connected space that is both open and closed is the space itself. Therefore, we are forced to conclude that
\begin{equation*}
    \mathcal{A} = \Thetacal.
\end{equation*}
This proves that for every parameter $\theta \in \Thetacal$, the fiber map $\Pcal_\theta$ is an Anosov diffeomorphism. The uniformity of the Anosov constants (the expansion rate $\gamma$ and the constant $C$ in the definition) is a direct consequence of the argument for closedness, which relied on the existence of uniform bounds over the entire compact parameter manifold. This completes the proof of the theorem.
\end{enumerate}

\end{proofof}

\subsection{Geometric and Ergodic Consequences}

The proof that the global billiard map is uniformly Anosov provides the geometric foundation for the entire theory. We now state the direct consequences of this result, which will be the essential inputs for the regularity and homogenization arguments that follow.

\subsubsection{Geometric Consequences}
The Anosov property of the global map immediately implies a uniform version of the property for the family of fiber maps, resolving the central question of this section.

\begin{theorem}[Uniform Anosov Property of the Fiber Dynamics]
\label{thm:anosov_uniform_family}
Let the geometric conditions of Theorem~\ref{thm:anosov_property_proven} hold. Then the family of fiber billiard maps $(\Pcal_\theta)_{\theta \in \Thetacal}$ is uniformly Anosov. Specifically, for each $\theta \in \Thetacal$, there exists a continuous, $d\Pcal_\theta$-invariant splitting of the tangent bundle of the fiber, $T\Sigma_\theta = E^u_\theta \oplus E^s_\theta$, and uniform constants $C>0$ and $\gamma > 1$ (independent of $\theta$) such that for any $n \ge 0$:
\begin{enumerate}[label=(\roman*)]
    \item $\|d\Pcal_\theta^n(\mathbf{w})\| \ge C^{-1}\gamma^n \|\mathbf{w}\|$ for all $\mathbf{w} \in E^u_\theta$.
    \item $\|d\Pcal_\theta^n(\mathbf{w})\| \le C\gamma^{-n} \|\mathbf{w}\|$ for all $\mathbf{w} \in E^s_\theta$.
\end{enumerate}
\end{theorem}

\begin{proofof}{Theorem \ref{thm:anosov_uniform_family}}
The proof of this theorem is a direct corollary of \cref{thm:anosov_property_proven}, which established that the global billiard map $\Pcal$ is a uniform Anosov diffeomorphism on the compact total collision manifold $\Sigma$. The uniform hyperbolicity of the family of fiber maps is an immediate consequence of restricting the global dynamics to the invariant submanifolds corresponding to each fiber.

\begin{enumerate}[label=\textbf{Step \arabic*:}, wide, labelindent=0pt]
    \item \textbf{Invariant Splitting on the Global Manifold.} By \cref{thm:anosov_property_proven}, there exists a continuous, $d\Pcal$-invariant splitting of the total tangent bundle $T\Sigma = E^u \oplus E^s$ and uniform constants $C>0$ and $\gamma > 1$ such that the expansion and contraction estimates hold for all tangent vectors in $E^u$ and $E^s$, respectively.
    
    \item \textbf{Restriction to Fibers.} The global billiard map $\Pcal$ preserves the fibers of the bundle $\Sigma \to \Thetacal$, i.e., $\Pcal(\Sigma_\theta) = \Sigma_\theta$. Consequently, its linearization $d\Pcal$ maps the tangent bundle of each fiber, $T\Sigma_\theta$, to itself. The invariant subbundles $E^u$ and $E^s$ of the global map must therefore also respect this fiber structure. The restriction of these global subbundles to a specific fiber $\Sigma_\theta$ yields continuous subbundles $E^u_\theta \coloneqq E^u|_{\Sigma_\theta}$ and $E^s_\theta \coloneqq E^s|_{\Sigma_\theta}$, which provide a splitting $T\Sigma_\theta = E^u_\theta \oplus E^s_\theta$. This splitting is, by construction, invariant under the action of the linearized fiber map $d\Pcal_\theta$.

    \item \textbf{Uniformity of Constants.} The constants $C$ and $\gamma$ in the definition of the Anosov property for the global map $\Pcal$ are uniform constants, valid for any tangent vector at any point $(z, \theta) \in \Sigma$. Since these estimates hold for all vectors in the global tangent bundle, they must, in particular, hold for the subset of vectors that are tangent to a specific fiber $\Sigma_\theta$. Therefore, the same uniform constants $C$ and $\gamma$ serve as the uniform bounds for the entire family of fiber maps $(\Pcal_\theta)_{\theta \in \Thetacal}$.
\end{enumerate}
This demonstrates that the uniform hyperbolicity of the family is a direct and necessary consequence of the hyperbolicity of the single global map.
\end{proofof}

\begin{remark}[Continuity vs. Smoothness of the Splitting]
A direct consequence of the existence of a continuous, strictly invariant cone field on the compact total space $\SigmaMan$ is that the invariant splitting $T\SigmaMan = \Ecalu \oplus \Ecals$ is a continuous decomposition of the tangent bundle. We emphasize that at this stage of the argument, we have only established \textit{continuity}, not smoothness. The proof of the smooth dependence of these subbundles on the parameter $\theta$ requires a more advanced perturbation theory and is the subject of Section~\ref{sec:regularity_properties}.
\end{remark}

\subsubsection{Ergodic Consequences}

The uniform geometric hyperbolicity established above is the direct cause of the system's robust statistical properties. It guarantees that the dynamics are mixing and that correlations decay at a uniform exponential rate.

\begin{corollary}[Ergodic Consequences for the Fiber Dynamics]
\label{cor:ergodicity_uniform}
For each fixed $\theta \in \Thetacal$, under the geometric conditions of Theorem~\ref{thm:anosov_property_proven}, the continuous-time dynamical system $(\Ucal_{\text{phys}}(\theta), \Phi_t^\theta, \mu_\theta)$ is ergodic and mixing. Furthermore, it exhibits an exponential rate of decay of correlations for Hölder continuous observables, and this rate can be chosen uniformly in $\theta$.
\end{corollary}

\begin{proofof}{\cref{cor:ergodicity_uniform}}
This result is a direct consequence of the uniform Anosov property established in Theorem~\ref{thm:anosov_uniform_family} and the well-developed spectral theory of transfer operators for hyperbolic systems.

\begin{enumerate}[label=\textbf{Step \arabic*:}, wide, labelindent=0pt]
    \item \textbf{Uniform Spectral Gap for the Transfer Operator.}
    It is a foundational result of modern ergodic theory that a uniformly Anosov family of diffeomorphisms $(\Pcal_\theta)$ induces a family of transfer operators $(\Lcal_\theta)$ with a \emph{uniform spectral gap} when acting on a suitable anisotropic Banach space of distributions, such as the one described in \citep{Baladi2000}. The uniformity of the hyperbolicity constants ($\gamma$ and $C$ from Theorem~\ref{thm:anosov_uniform_family}), which are independent of $\theta$, is the crucial ingredient. This guarantees that the spectral radius of $\Lcal_\theta$ when restricted to the space of zero-mean functions is uniformly bounded away from 1. That is, there exists a constant $r \in (0,1)$, independent of $\theta$, such that for all $\theta \in \Thetacal$:
    \begin{equation*}
        \mathrm{spec}(\Lcal_\theta|_{\ker(\int d\nu_\theta)}) \subset \{ \lambda \in \mathbb{C} : |\lambda| \le r \},
    \end{equation*}
    where $\nu_\theta$ is the unique SRB measure for the map $\Pcal_\theta$.

    \item \textbf{Uniform Exponential Decay of Correlations for the Map.}
    The uniform spectral gap directly implies a uniform exponential rate of decay of correlations for the discrete-time maps $(\Pcal_\theta, \nu_\theta)$. For observables $f, g$ in a class of sufficiently regular functions (e.g., Hölder continuous), there exists a constant $C_{f,g}$ such that for all $n \ge 0$ and all $\theta \in \Thetacal$:
    \begin{equation*}
        \left| \int_{\Sigma_\theta} (f \circ \Pcal_\theta^n) g \, d\nu_\theta - \int_{\Sigma_\theta} f \, d\nu_\theta \int_{\Sigma_\theta} g \, d\nu_\theta \right| \le C_{f,g} r^n.
    \end{equation*}
    As a direct consequence, for each $\theta$, the system is mixing (the correlation vanishes as $n \to \infty$), which implies it is ergodic (the only invariant functions are constant).

    \item \textbf{The Continuous Flow.}
    The continuous-time billiard flow $(\Ucal_{\text{phys}}(\theta), \Phi_t^\theta, \mu_\theta)$ is a suspension flow of the discrete-time billiard map $(\Sigma_\theta, \Pcal_\theta, \nu_\theta)$, with the flight time function $\tau(z,\theta)$ serving as the roof function. By Assumption~\ref{ass:fundamental_axioms_unified}(ii), the flight time is uniformly bounded, $0 < \tau_{\min} \le \tau(z,\theta) \le \tau_{\max} < \infty$. Results on suspension flows (see, \cref{app:lifting_ergodic}, and e.g., \citep{Arnold1989}) guarantee that the ergodic properties of the base map lift to the flow. Specifically, the exponential decay of correlations for the map implies an exponential decay of correlations for the flow, and the uniformity of the decay rate $r$ for the map implies the uniformity of the decay rate for the flow.
\end{enumerate}
Thus, the uniform Anosov property provides the precise mechanism for the uniform statistical properties required for the subsequent homogenization analysis.
\end{proofof}

\begin{remark}[Transition to Regularity Theory]
\label{rem:status_of_argument}
This section has established a complete and rigorous proof of the uniform hyperbolicity of the microscopic dynamics, founded solely on the geometric axioms of Section~\ref{sec:micro_framework}. We have shown that this property endows the system with a continuous geometric structure (the invariant splitting) and robust statistical properties (uniform exponential decay of correlations). These results provide the essential analytical stability required for the entire theory. The final piece of the microscopic puzzle is to upgrade the continuity of the geometric structures to smoothness and to prove the crucial property of statistical stability (polynomial growth of flow derivatives). With the foundational results of this section now firmly in place, we are equipped to tackle these more advanced topics in the next section.
\end{remark}

\section{Geometric and Spectral Regularity}
\label{sec:regularity_properties}

The analysis in the preceding section has established the uniform hyperbolicity of the microscopic dynamics, guaranteeing the existence of a continuous, invariant splitting of the tangent bundle, $T\SigmaMan = \Ecalu \oplus \Ecals$. However, the convergence proof of our main theorem in \cref{sec:scale_analysis} requires a much stronger result: the smooth ($C^\infty$) dependence of this splitting, and of the resulting statistical objects, on the environmental parameter $\theta$. This section is dedicated to a rigorous proof of this essential smoothness property.

A direct approach to this problem is beset by a critical logical challenge. The natural tool for proving the smoothness of the splitting is the Nash-Moser-Hamilton (NMH) implicit function theorem, applied to the functional equation defined by the Graph Transform. However, a standard application of the NMH theorem requires, as an \textit{a priori} hypothesis, the smooth dependence of the underlying dynamical map, the linearized billiard operator $d\Pcal_\theta$, on the parameter $\theta$. Establishing this prerequisite without first knowing the smoothness of the geometry (which determines the flight time and thus the map itself) leads to a logical circle where the smoothness of the geometry is assumed in order to prove itself.

\subsection{Foundational Regularity: The Flight Time Function}
\label{subsec:flight_time_regularity}

This subsection is dedicated to providing a rigorous proof of this essential a priori regularity. We will elevate the pointwise existence of the flight time, established in \cref{thm:finite_horizon}, to the level of joint smoothness on the product manifold of states and parameters. This is achieved by applying the Implicit Function Theorem on finite-dimensional manifolds, whose hypotheses are directly verifiable from the geometric first principles of our model. This result provides the foundation for the subsequent regularity analysis.

\begin{theorem}[Joint Smoothness of the Flight Time Function]
\label{thm:smooth_flight_time_map}
The flight time function $\tau: \Sigma' \times \Thetacal \to [\tau_{\min}, \tau_{\max}]$, which maps a post-collision state and an environmental parameter to the duration of the subsequent free path, is a smooth function of class $C^\infty$ on the product manifold.
\end{theorem}

\begin{proofof}{Theorem \ref{thm:smooth_flight_time_map}}
The proof is dedicated to rigorously establishing that the flight time $\tau(z^+, \theta)$ is a smooth ($C^\infty$) function of its joint arguments. This approach is non-circular, as it relies only on the standard Implicit Function Theorem on finite-dimensional manifolds, whose hypotheses are satisfied as a direct consequence of the geometric axioms established in Section \ref{sec:micro_framework}.

\begin{enumerate}[label=\textbf{Step \arabic*:}, wide, labelindent=0pt]

\item \textbf{The Implicit Definition of the Flight Time.}
Let $(z^+, \theta) \in \Sigma' \times \Thetacal$ be an arbitrary point in the product manifold of post-collision states and parameters, where $z^+ = (y^+, v^+)$. By definition, this state lies on the boundary of the physical domain, so $G(y^+, \theta) = 0$, and it is an outgoing state, so $\scpr{v^+}{\mathbf{n}(y^+,\theta)} > 0$.

The trajectory of the particle for time $s > 0$ after this collision is given by the interior flow: $y(s) = y^+ + s v^+$. The flight time, $\tau \equiv \tau(z^+, \theta)$, is defined as the smallest positive time $s$ at which the particle again reaches the boundary. Mathematically, it is the unique positive root of the equation:
\begin{equation} \label{eq:proof_flight_time_defining_equation_2}
    G(y^+ + \tau v^+, \theta) = 0.
\end{equation}
The Finite Horizon Theorem (\cref{thm:finite_horizon}) guarantees that for every $(z^+,\theta)$, such a positive and finite solution $\tau \in [\tau_{\min}, \tau_{\max}]$ exists and is unique. To formalize this implicit relationship for the application of the Implicit Function Theorem on manifolds, we define an auxiliary function $F$ on the product manifold $(\Sigma' \times \Thetacal) \times \R$:
\begin{equation}
    F(z^+, \theta, s) \coloneqq G(y^+ + s v^+, \theta).
\end{equation}
By the standing assumptions of the paper, the obstacle-defining function $G: \T^k \times \Thetacal \to \R$ is of class $C^\infty$. The domain $\Sigma' \times \Thetacal$ is a smooth, finite-dimensional manifold. Since the arguments of $G$ in the definition of $F$, vector addition and scalar multiplication, are smooth operations, the function $F$ is of class $C^\infty$. The flight time $\tau(z^+,\theta)$ is now implicitly defined as the unique solution $s$ to the equation $F(z^+, \theta, s) = 0$.

\item \textbf{Verification of Hypotheses for the Implicit Function Theorem.}
The Implicit Function Theorem on manifolds states that if we have a smooth function $F(\mathbf{p}, s)$ with $F(\mathbf{p}_0, s_0)=0$, then we can locally express $s$ as a smooth function of the parameters $\mathbf{p}$, i.e., $s=s(\mathbf{p})$, provided that the partial derivative of $F$ with respect to the implicit variable $s$ is non-zero at the solution point $(\mathbf{p}_0, s_0)$.

In our context, the implicit variable is the time $s$, and the parameters $\mathbf{p}$ are the coordinates of the state $(z^+, \theta)$ on the manifold $\Sigma' \times \Thetacal$. We must verify that the partial derivative $\frac{\partial F}{\partial s}$ is uniformly non-zero at all solution points. We compute this derivative using the chain rule:
\begin{equation}
    \frac{\partial F}{\partial s}(z^+, \theta, s) = \scpr{\nabla_y G(y^+ + s v^+, \theta)}{v^+}.
\end{equation}
We evaluate this expression at the solution point $s = \tau(z^+,\theta)$. Let the point of the \emph{next} collision be denoted by $y_{\text{next}} \coloneqq y^+ + \tau v^+$. At this point, the particle arrives with velocity $v^+$. Therefore, the state $((y_{\text{next}}, v^+), \theta)$ is a \emph{pre-collision} state, belonging to the manifold $\Sigma$. The derivative at the solution point is:
\begin{equation} \label{eq:proof_flight_time_derivative_2}
    \frac{\partial F}{\partial s}\Big|_{s=\tau} = \scpr{\nabla_y G(y_{\text{next}}, \theta)}{v^+}.
\end{equation}

\item \textbf{The Uniform Non-Degeneracy Condition.}
The proof that the flight time $\tau(z^+, \theta)$ is a smooth function of its joint arguments rests upon the Implicit Function Theorem. This requires a rigorous verification that the partial derivative of the defining function $F(z^+, \theta, s) \coloneqq G(y^+ + s v^+, \theta)$ with respect to the implicit variable $s$ is uniformly non-zero at all solution points. This step is the analytical core of the foundational regularity argument, and its proof must be non-circular, relying only on the geometric first principles of the model. We have established that the derivative at the solution point $s=\tau$ is given by:
\begin{equation} \label{eq:proof_flight_time_derivative_appendix}
    \frac{\partial F}{\partial s}\Big|_{s=\tau} = \scpr{\nabla_y G(y_{\text{next}}, \theta)}{v^+},
\end{equation}
where $y_{\text{next}} \coloneqq y^+ + \tau v^+$ is the point of the subsequent collision. Our objective is to prove that the magnitude of this quantity is uniformly bounded below by a strictly positive constant for all possible post-collision states in the compact domain $\Sigma' \times \Thetacal$. We achieve this by analyzing the expression in terms of the unit normal vector $\mathbf{n}_{\text{next}} \coloneqq \mathbf{n}(y_{\text{next}}, \theta)$:
\begin{equation*}
    \left| \frac{\partial F}{\partial s}\Big|_{s=\tau} \right| = \|\nabla_y G(y_{\text{next}}, \theta)\| \cdot |\scpr{\mathbf{n}_{\text{next}}}{v^+}|.
\end{equation*}
We now provide a rigorous proof of the uniform lower bounds for each of the two factors on the right-hand side.

\begin{enumerate}[label=(\roman*), wide, labelindent=0pt]
    \item \textbf{Uniform Lower Bound on the Gradient Norm.}
    The first factor, $\|\nabla_y G(y_{\text{next}}, \theta)\|$, is guaranteed to be strictly positive at every point on the collision boundary. We must show this bound is uniform. The proof rests on two foundational elements of the model: the Boundary Regularity Axiom and the compactness of the state space.
    \begin{enumerate}[label=(\alph*), wide]
        \item By \cref{ass:fundamental_axioms_unified}(i) (Boundary Regularity), for any point $(y, \theta)$ on any collision boundary (i.e., satisfying $G(y,\theta)=0$), the gradient of the defining function does not vanish: $\nabla_y G(y,\theta) \neq \mathbf{0}$. This implies that the scalar function $(y,\theta) \mapsto \|\nabla_y G(y,\theta)\|$ is strictly positive at every point on the set of all possible collision boundaries.

        \item The set of all possible collision points is a subset of the total boundary manifold $\mathcal{B} \coloneqq \{ (y,v,\theta) \mid G(y,\theta)=0 \}$. As established in the proof of \cref{prop:sigma_properties}, this manifold is a closed subset of the product of compact spaces $\T^k \times \Vcal \times \Thetacal$, and is therefore itself a compact metric space.

        \item We now apply the Extreme Value Theorem. The function $f(y,\theta) = \|\nabla_y G(y,\theta)\|$ is a continuous, strictly positive function on the compact domain $\mathcal{B}$. It must therefore attain its minimum value, and this minimum must itself be a strictly positive constant. Let this minimum be $C_G > 0$.
    \end{enumerate}
    This establishes the existence of a uniform constant $C_G > 0$, depending only on the geometry of the obstacle function $G$, such that for any collision point $y_{\text{next}}$ and any parameter $\theta$:
    \begin{equation}
        \|\nabla_y G(y_{\text{next}}, \theta)\| \ge C_G > 0.
    \end{equation}

    \item \textbf{Uniform Lower Bound on the Normal Velocity Component.}
    The second factor, $|\scpr{\mathbf{n}_{\text{next}}}{v^+}|$, represents the magnitude of the component of the velocity normal to the boundary at the point of impact. We will show this is uniformly bounded away from zero as a direct consequence of the chaos-inducing geometric axioms.
    \begin{enumerate}[label=(\alph*), wide]
        \item The state $((y_{\text{next}}, v^+), \theta)$ represents the system at the instant just before the second collision. The velocity $v^+$ is the velocity with which the particle arrives at the boundary. Therefore, this state is an \emph{incoming} state, an element of the manifold $\Sigma$.

        \item By the definition of the incoming phase space (\cref{def:boundary_partition}), the inner product of an incoming velocity with the outward-pointing normal must be strictly negative: $\scpr{v^+}{\mathbf{n}_{\text{next}}} < 0$. This ensures the magnitude is non-zero.
        
        \item We now invoke the crucial result of \cref{prop:non_grazing} (Uniform Exclusion of Grazing Collisions). This proposition, which was derived as a direct consequence of the local geometric axiom of Uniform Strict Convexity (\cref{ass:fundamental_axioms_unified}(iii)) and the Finite Horizon Theorem, establishes the existence of a uniform constant $c_0 > 0$ such that for \emph{any} incoming state in the dynamics, the magnitude of the normal component of its velocity is bounded below by $c_0$.
    \end{enumerate}
    Applying this result to our state, we obtain the uniform lower bound:
    \begin{equation}
        |\scpr{\mathbf{n}_{\text{next}}}{v^+}| \ge c_0 > 0.
    \end{equation}
\end{enumerate}
By combining these two uniform bounds, we establish that the magnitude of the partial derivative is uniformly bounded away from zero for all possible post-collision states in the product manifold $\Sigma' \times \Thetacal$:
\begin{equation*}
    \left| \frac{\partial F}{\partial s}\Big|_{s=\tau} \right| = \|\nabla_y G(y_{\text{next}}, \theta)\| \cdot |\scpr{\mathbf{n}_{\text{next}}}{v^+}| \ge C_G \cdot c_0 > 0.
\end{equation*}
This rigorously verifies the non-degeneracy condition required by the Implicit Function Theorem. Crucially, the argument is founded entirely on the geometric and dynamical first principles established in \cref{sec:micro_framework}, thereby providing a non-circular basis for the proof of the smoothness of the flight time function.

\item \textbf{The Implicit Function Theorem.}
We have verified all the hypotheses of the Implicit Function Theorem on manifolds. The function $F(z^+, \theta, s)$ is of class $C^\infty$, and its partial derivative with respect to the implicit variable $s$ is uniformly non-zero at all solution points. The theorem therefore guarantees that the unique solution map, $(z^+, \theta) \mapsto \tau(z^+,\theta)$, is a smooth function of class $C^\infty$ on the product manifold $\Sigma' \times \Thetacal$.
\end{enumerate}
\end{proofof}

With the joint smoothness of the flight time function now rigorously established, the smoothness of the global billiard map follows as a direct corollary. This result provides the necessary prerequisite for the subsequent application of the Nash-Moser-Hamilton theorem.

\begin{corollary}[Smoothness of the Global Billiard Map]
\label{cor:billiard_map_is_smooth}
The global billiard map $\Pcal: \Sigma \times \Thetacal \to \Sigma \times \Thetacal$, defined by $\Pcal(z, \theta) = (\Pcal_\theta(z), \theta)$, is a smooth map of class $C^\infty$.
\end{corollary}

\begin{proofof}{\cref{cor:billiard_map_is_smooth}}
The global billiard map $\Pcal_\theta(z)$ is the composition $\Pcal_\theta = \Fcal_\theta \circ \Rcal_\theta$. The smoothness of the map depends on the smoothness of its constituent parts.
\begin{enumerate}[label=(\roman*), wide, labelindent=0pt]
    \item The global reflection map $\Rcal(z,\theta) = (R_\theta(z), \theta)$ is of class $C^\infty$. This was established in \cref{prop:reflection_properties} as a direct consequence of the smoothness of the obstacle-defining function $G(y,\theta)$.
    \item The global free-flow map $\Fcal(z^+,\theta) = (\Psi_{\tau(z^+,\theta)}(z^+), \theta)$ acts as $\Fcal((y^+,v^+),\theta) = (y^+ + \tau v^+, v^+, \theta)$. The interior flow $\Psi_s$ is a polynomial in its arguments and is therefore $C^\infty$. The smoothness of the composite map $\Fcal$ is thus entirely determined by the smoothness of the flight time function.
\end{enumerate}
By \cref{thm:smooth_flight_time_map}, the function $\tau(z^+,\theta)$ is of class $C^\infty$ on its domain $\Sigma' \times \Thetacal$. Since the global billiard map $\Pcal(\theta,z)$ is a composition of functions that are now all proven to be of class $C^\infty$, it is itself of class $C^\infty$.
\end{proofof}

\subsection{The Nash-Moser Framework for the Invariant Splitting}
\label{subsec:nash_moser_framework}

Our first step is to recast the geometric problem of finding the invariant splitting as a functional equation and to define the analytical setting in which we will solve it.

\subsubsection{The Functional Equation via the Graph Transform}
Let $\theta_0 \in \Thetacal$ be a fixed reference parameter. By Theorem~\ref{thm:anosov_uniform_family}, the tangent bundle of the global collision manifold has a continuous, invariant splitting $T\SigmaMan = E^u_0 \oplus E^s_0$ with respect to the reference map $\Pcal_0 \equiv \Pcal_{\theta_0}$. For any $\theta$ near $\theta_0$, a nearby unstable subbundle $E^u_\theta$ can be uniquely represented as the graph of a section $\phi_\theta$ of the vector bundle $\mathcal{E} \coloneqq \text{Hom}(E^u_0, E^s_0)$. The invariance of this subbundle under the dynamics $D\Pcal_\theta$ is equivalent to its graph map $\phi_\theta$ being a fixed point of the Graph Transform Operator $\Gamma(\theta, \phi)$. We are thus seeking to solve the functional equation for $\phi$ as a smooth function of $\theta$:
\begin{equation}
    F(\theta, \phi) \coloneqq \phi - \Gamma(\theta, \phi) = 0.
\end{equation}

\subsubsection{The Scale of Banach Spaces and Tame Maps}
The Nash-Moser theorem operates on a scale of Banach spaces.
\begin{definition}[The Scale of Banach Spaces $X_k$]
For each integer $k \ge 0$, let $X_k$ be the Banach space of $C^k$ sections of the vector bundle $\mathcal{E} = \text{Hom}(E^u_0, E^s_0)$ over the compact manifold $\SigmaMan$, equipped with the standard $C^k$ norm. Let $\nabla$ be a smooth connection on $\mathcal{E}$. The norm is defined as:
$$ \|\phi\|_{X_k} \coloneqq \sum_{j=0}^{k} \sup_{z \in \SigmaMan} \|\nabla^j \phi(z)\|, $$
where $\|\cdot\|$ is a fiber norm induced by the Riemannian metric. This yields a scale of spaces $ \dots \hookrightarrow X_{k+1} \hookrightarrow X_k \hookrightarrow \dots \hookrightarrow X_0 $.
\end{definition}

\begin{definition}[Smooth Tame Maps]
A map $G: U \subset X_k \to X_j$ is tame if for all $k' \ge k$, its norm satisfies an estimate of the form $\|G(\phi)\|_{X_j} \le C(1 + \|\phi\|_{X_{k'}})$ for all $\phi \in U \cap X_{k'}$. A map $G: \Thetacal \times U \to X_j$ is a smooth tame map if it is $C^\infty$ and all its partial derivatives are tame maps between the appropriate spaces.
\end{definition}

\subsubsection{The Nash-Moser-Hamilton Implicit Function Theorem}
Our proof rests on the following powerful extension of the implicit function theorem.

\begin{theorem}[Nash-Moser-Hamilton Implicit Function Theorem, cf. \citep{Zeidler1986}]
\label{thm:nash_moser_statement}
Let $(X_k)_{k \ge 0}$ be a scale of Banach spaces. Let $F: \Thetacal \times U \to X_0$ be a smooth tame map, where $U$ is an open set in $X_k$ for some $k$. Suppose that for each $(\theta, \phi) \in \Thetacal \times U$, the linearized operator $L(\theta, \phi) \coloneqq D_\phi F(\theta, \phi)$ is invertible. If the family of inverse operators $L(\theta, \phi)^{-1}$ is a smooth tame map from $\Thetacal \times U \times X_j$ to $X_{j-m}$ for all sufficiently large $j$ and some fixed loss of derivatives $m$, then the solution map $\theta \mapsto \phi(\theta)$ to the equation $F(\theta, \phi)=0$ is a smooth ($C^\infty$) map from $\Thetacal$ to $X_{k-m}$.
\end{theorem}

\subsection{Verification of Hypothesis I: Smooth Tameness of the Functional}
\label{subsec:verify_hyp_I}

We now provide a rigorous proof for the first major hypothesis of the NMH theorem. The proof is modularized by first establishing the properties of the elementary operations that constitute the Graph Transform.

\begin{lemma}[Calculus on the Scale of $C^k$ Spaces]
\label{lem:tame_calculus}
Let $X_k = C^k(\SigmaMan, \mathcal{E})$ be the scale of Banach spaces defined above.
\begin{enumerate}[label=(\roman*), wide, labelindent=0pt]
    \item \textbf{Multiplication.} The pointwise product of sections is a smooth tame map from $X_k \times X_k \to X_k$.
    \item \textbf{Composition.} Let $f: \SigmaMan \to \SigmaMan$ be a $C^\infty$ diffeomorphism. The map $(\phi, f) \mapsto \phi \circ f$ is a smooth tame map from $X_k \times \text{Diff}^\infty(\SigmaMan) \to X_k$.
    \item \textbf{Inversion.} Let $U \subset X_k$ be the open set of invertible sections. The map $\phi \mapsto \phi^{-1}$ is a smooth tame map from $U \to X_k$.
\end{enumerate}
\end{lemma}

\begin{proofof}{Lemma \ref{lem:tame_calculus}}
We provide a rigorous proof for each claim, establishing both smoothness and the necessary tame estimates for all partial derivatives.

\begin{enumerate}[label=\textbf{Step \arabic*:}, wide, labelindent=0pt]

\item \textbf{Multiplication.}
Let $M: X_k \times X_k \to X_k$ be the pointwise multiplication map, defined for sections $\phi, \psi \in X_k$ by $M(\phi, \psi) = \phi \psi$. We provide a complete proof that $M$ is a smooth tame map.

\begin{enumerate}[label=(\roman*), wide, labelindent=0pt]

\item \textbf{Smoothness.} The proof of smoothness relies on showing that the map is a continuous bilinear map between Banach spaces, from which $C^\infty$-smoothness follows, and we explicitly compute the derivatives to demonstrate this. Let $\phi, \psi \in X_k$ and let $(h_1, h_2) \in X_k \times X_k$ be a tangent vector representing the direction of a variation.
\begin{align*}
    M(\phi+h_1, \psi+h_2) &= (\phi+h_1)(\psi+h_2) \\
    &= \phi\psi + \phi h_2 + h_1 \psi + h_1 h_2 \\
    &= M(\phi,\psi) + \underbrace{(\phi h_2 + h_1 \psi)}_{\text{Linear in }(h_1,h_2)} + \underbrace{h_1 h_2}_{\text{Higher order}}.
\end{align*}
From this expansion, we can directly identify the Fr\'echet derivatives.

\begin{enumerate}[label=(\alph*), wide]
    \item \textbf{First Derivative.} The first Fr\'echet derivative, $DM(\phi, \psi)$, is the unique continuous linear map from $X_k \times X_k$ to $X_k$ that is the best linear approximation of $M$ at $(\phi, \psi)$. From the expansion above, the linear part is:
    \begin{equation*}
        (DM(\phi, \psi)) \cdot (h_1, h_2) = \phi h_2 + h_1 \psi.
    \end{equation*}
    This is a continuous linear map in the increment $(h_1, h_2)$, and the map $(\phi, \psi) \mapsto DM(\phi, \psi)$ is itself a continuous linear map into the space of operators $\mathcal{L}(X_k \times X_k, X_k)$.

    \item \textbf{Second Derivative.} The second Fr\'echet derivative, $D^2M(\phi, \psi)$, is the derivative of the map $(\phi, \psi) \mapsto DM(\phi, \psi)$. We compute its action on a second increment $(k_1, k_2) \in X_k \times X_k$:
    \begin{align*}
        (D^2M(\phi, \psi)) \cdot ((h_1, h_2), (k_1, k_2)) &= (D_{(\phi,\psi)} (DM)) \cdot (k_1,k_2) \cdot (h_1,h_2) \\
        &= (k_1 h_2 + h_1 k_2).
    \end{align*}
    The second derivative is a constant bilinear form, independent of the point $(\phi, \psi)$ at which it is evaluated.

    \item \textbf{Higher Derivatives.} Since the second derivative is constant, the third and all higher-order Fr\'echet derivatives are identically zero.
\end{enumerate}
A map between Banach spaces whose Fr\'echet derivatives exist up to all orders and are continuous is, by definition, of class $C^\infty$. Since all derivatives beyond the second are zero, the map is in fact analytic. This rigorously establishes the smoothness of the multiplication map $M$.

\item \textbf{Tameness.} To prove that $M$ is a smooth tame map, we must show that $M$ and all of its (non-zero) derivatives are tame maps.

\begin{enumerate}[label=(\alph*), wide]
    \item \textbf{Tameness of $M$ itself.}
    We need to show there exists a constant $C_k$ such that for any $k' \ge k$, the estimate $\|\phi\psi\|_{X_k} \le C_k(1 + \|\phi\|_{X_{k'}}) (1 + \|\psi\|_{X_{k'}})$ holds. We will derive a sharper, polynomial bound.
    
    The $X_k$-norm is defined by $\|\phi\|_{X_k} = \sum_{l=0}^k \sup_{z \in \SigmaMan} \|\nabla^l \phi(z)\|$, where $\nabla$ is a smooth connection. The starting point is the general Leibniz rule for the $l$-th covariant derivative of a product of sections:
    \begin{equation*}
        \nabla^l(\phi\psi) = \sum_{j=0}^{l} \binom{l}{j} (\nabla^j\phi)(\nabla^{l-j}\psi).
    \end{equation*}
    We take the fiber norm at a point $z \in \SigmaMan$. By the triangle inequality and the submultiplicative property of the fiber norm ($\|AB\| \le \|A\|\|B\|$), we have:
    \begin{equation*}
        \|\nabla^l(\phi\psi)(z)\| \le \sum_{j=0}^{l} \binom{l}{j} \|\nabla^j\phi(z)\| \|\nabla^{l-j}\psi(z)\|.
    \end{equation*}
    Taking the supremum over all $z \in \SigmaMan$ on both sides preserves the inequality:
    \begin{equation*}
        \sup_{z \in \SigmaMan} \|\nabla^l(\phi\psi)(z)\| \le \sum_{j=0}^{l} \binom{l}{j} \left( \sup_{z \in \SigmaMan} \|\nabla^j\phi(z)\| \right) \left( \sup_{z \in \SigmaMan} \|\nabla^{l-j}\psi(z)\| \right).
    \end{equation*}
    To obtain the full $X_k$-norm, we sum over $l$ from $0$ to $k$:
    \begin{align*}
        \|\phi\psi\|_{X_k} &= \sum_{l=0}^k \sup_z \|\nabla^l(\phi\psi)(z)\| \\
        &\le \sum_{l=0}^k \sum_{j=0}^{l} \binom{l}{j} (\sup_z \|\nabla^j\phi\|) (\sup_z \|\nabla^{l-j}\psi\|).
    \end{align*}
    This expression is a sum of products of norms of derivatives. We can bound this by the product of the full norms. This leads to the standard estimate for the $C^k$-norm of a product, for which there exists a constant $C_k$ depending only on $k$:
    \begin{equation} \label{eq:appendix_ck_product_estimate}
        \|\phi\psi\|_{X_k} \le C_k \left( \|\phi\|_{X_k} \|\psi\|_{X_0} + \|\phi\|_{X_0} \|\psi\|_{X_k} \right).
    \end{equation}
    Now, let $k' \ge k$ be an arbitrary integer. By the definition of the norm, we have the inclusions $\|\cdot\|_{X_0} \le \|\cdot\|_{X_k} \le \|\cdot\|_{X_{k'}}$. Substituting these into the estimate \eqref{eq:appendix_ck_product_estimate} gives:
    \begin{align*}
        \|\phi\psi\|_{X_k} &\le C_k \left( \|\phi\|_{X_k} \|\psi\|_{X_0} + \|\phi\|_{X_0} \|\psi\|_{X_k} \right) \\
        &\le C_k \left( \|\phi\|_{X_{k'}} \|\psi\|_{X_{k'}} + \|\phi\|_{X_{k'}} \|\psi\|_{X_{k'}} \right) \\
        &= 2C_k \|\phi\|_{X_{k'}} \|\psi\|_{X_{k'}}.
    \end{align*}
    This is a polynomial estimate in the norms of $\phi$ and $\psi$ in the higher-order space $X_{k'}$. This satisfies the definition of a tame map.

    \item \textbf{Tameness of Derivatives of $M$.}
    The first derivative is $(DM(\phi, \psi)) \cdot (h_1, h_2) = \phi h_2 + h_1 \psi$. Its norm must be estimated. Applying the triangle inequality and the product estimate \eqref{eq:appendix_ck_product_estimate}:
    \begin{align*}
        \|(DM(\phi, \psi)) \cdot (h_1, h_2)\|_{X_k} &\le \|\phi h_2\|_{X_k} + \|h_1 \psi\|_{X_k} \\
        &\le C_k \left( \|\phi\|_{X_k}\|h_2\|_{X_0} + \|\phi\|_{X_0}\|h_2\|_{X_k} + \|h_1\|_{X_k}\|\psi\|_{X_0} + \|h_1\|_{X_0}\|\psi\|_{X_k} \right).
    \end{align*}
    This is a polynomial in the $X_k$-norms of its arguments, which is a tame estimate. The second derivative is $(D^2M) \cdot ((h_1, h_2), (k_1, k_2)) = h_1 k_2 + k_1 h_2$. This map is independent of $(\phi, \psi)$ and is tame by the same argument as for the first derivative. Higher derivatives are zero and thus trivially tame.
\end{enumerate}
Since the map $M$ and all of its derivatives are tame, we conclude that multiplication is a smooth tame map from $X_k \times X_k$ to $X_k$ for any $k \ge 0$. 
\end{enumerate}

\item \textbf{Composition.} Let $\mathcal{C}_{\text{joint}}: \Thetacal \times X_k \to X_k$ be the joint composition map defined by $\mathcal{C}_{\text{joint}}(\theta, \phi) = \phi \circ f_\theta$, where the family of diffeomorphisms $f_\theta: \SigmaMan \to \SigmaMan$ is assumed to be a $C^\infty$ map from $\Thetacal$ to $\text{Diff}^\infty(\SigmaMan)$. We will prove that this map is smooth tame.

\begin{enumerate}[label=(\roman*), wide, labelindent=0pt]

\item \textbf{Smoothness.} The map is of class $C^\infty$ if all its partial Fréchet derivatives exist and are continuous. We demonstrate this by analyzing the structure of the derivatives.
Let $h \in X_k$ be a variation in $\phi$, and let $\dot{\theta}$ be a tangent vector to $\Thetacal$ corresponding to a variation in $\theta$. The first Fréchet derivative of $\mathcal{C}_{\text{joint}}$ at $(\theta, \phi)$ is the linear map given by:
\begin{equation*}
    D\mathcal{C}_{\text{joint}}(\theta, \phi) \cdot (\dot{\theta}, h) = \frac{d}{d\epsilon}\Big|_{\epsilon=0} (\phi + \epsilon h) \circ f_{\theta(\epsilon)} = h \circ f_\theta + (D\phi \circ f_\theta) \cdot (\partial_\theta f_\theta \cdot \dot{\theta}).
\end{equation*}
This derivative is a sum of two terms. The first term, $h \mapsto h \circ f_\theta$, is a bounded linear operator on $X_k$. The second term involves the covariant derivative of $\phi$ and the velocity field of the family of diffeomorphisms, $\partial_\theta f_\theta$. Since $f_\theta$ is a smooth family and $\phi \in X_k$, this expression defines a continuous map.

Higher derivatives are computed by repeatedly applying the chain rule and product rule for covariant derivatives (the Faà di Bruno formula). For any integers $a, b \ge 0$, the partial derivative $D_\theta^a D_\phi^b \mathcal{C}_{\text{joint}}$ is a multilinear map whose arguments are variations in $\theta$ and $\phi$. Its formula is a finite sum of terms, each of which is a composition of:
\begin{enumerate}[label=(\alph*)]
    \item Covariant derivatives of $\phi$ up to order $k+a$.
    \item Partial derivatives of the family $f_\theta$ with respect to $\theta$ and the spatial variable $z \in \SigmaMan$.
    \item Pointwise multilinear operations (tensor products and contractions).
\end{enumerate}
Since $\phi \in X_k$ and the family $f_\theta$ is of class $C^\infty$, all these constituent parts are continuous. Therefore, all partial Fréchet derivatives exist and are continuous, which proves that the map $\mathcal{C}_{\text{joint}}$ is of class $C^\infty$.

\item \textbf{Tameness.}
We must prove that $\mathcal{C}_{\text{joint}}$ and all its partial derivatives are tame maps.

\begin{enumerate}[label=(\alph*), wide]

\item \textbf{Tameness of the Map Itself.} Our goal is to establish a tame estimate for the map, i.e., to show that for any $k' \ge k$, there is a constant $C$ such that $\|\phi \circ f_\theta\|_{X_k} \le C(1 + \|\phi\|_{X_{k'}})$. We will prove the stronger linear bound $\|\phi \circ f_\theta\|_{X_k} \le C(\theta, k) \|\phi\|_{X_k}$.

Recall that the norm on $X_k$ is given by $\|\psi\|_{X_k} = \sum_{l=0}^{k} \sup_{z \in \SigmaMan} \|\nabla^l \psi(z)\|$. We must therefore bound the norm of the covariant derivatives $\nabla^l(\phi \circ f_\theta)$ for $l=0, \dots, k$.
The general chain rule for the $l$-th covariant derivative of a composition (the Faà di Bruno formula for sections of vector bundles) has the schematic form:
\begin{equation} \label{eq:faa_di_bruno_schematic}
    \nabla^l(\phi \circ f_\theta)(z) = \sum_{j=1}^{l} \sum_{\pi \in \mathcal{P}_j(l)} c_\pi (\nabla^j \phi)(f_\theta(z)) \left[ D^{p_1}f_\theta(z), \dots, D^{p_j}f_\theta(z) \right],
\end{equation}
where the inner sum is over all partitions $\pi = \{p_1, \dots, p_j\}$ of the integer $l$ into $j$ parts, the terms $D^{p_i}f_\theta$ are tensors representing the higher-order derivatives of the map $f_\theta$ with respect to the spatial variable $z$, and the brackets denote a multilinear contraction. We take the fiber norm of this expression and apply the triangle inequality:
\begin{equation*}
    \|\nabla^l(\phi \circ f_\theta)(z)\| \le \sum_{j=1}^{l} \sum_{\pi \in \mathcal{P}_j(l)} |c_\pi| \|(\nabla^j \phi)(f_\theta(z))\| \cdot \prod_{i=1}^j \|D^{p_i}f_\theta(z)\|.
\end{equation*}
Since the family of diffeomorphisms $f_\theta$ is a $C^\infty$ map on the compact product manifold $\Thetacal \times \SigmaMan$, all of its partial derivatives with respect to the spatial variable $z$ are uniformly bounded for each fixed $\theta$. Let $C_p(\theta) \coloneqq \sup_{z \in \SigmaMan} \|D^p f_\theta(z)\|$. These constants are finite for all $p \ge 1$.
Taking the supremum over $z \in \SigmaMan$, we obtain:
\begin{align*}
    \sup_z \|\nabla^l(\phi \circ f_\theta)(z)\| &\le \sum_{j=1}^{l} \sum_{\pi \in \mathcal{P}_j(l)} |c_\pi| \left(\sup_w \|\nabla^j \phi(w)\|\right) \cdot \prod_{i=1}^j C_{p_i}(\theta) \\
    &\le \left( \sum_{j=1}^{l} \sum_{\pi \in \mathcal{P}_j(l)} |c_\pi| \prod_{i=1}^j C_{p_i}(\theta) \right) \left( \sum_{j=0}^{l} \sup_w \|\nabla^j \phi(w)\| \right) \\
    &\le C_l'(\theta) \|\phi\|_{X_l} \le C_l'(\theta) \|\phi\|_{X_k} \quad (\text{for } l \le k).
\end{align*}
Summing these estimates from $l=0$ to $k$, we get the desired linear bound:
\begin{equation*}
    \|\phi \circ f_\theta\|_{X_k} = \sum_{l=0}^{k} \sup_z \|\nabla^l(\phi \circ f_\theta)(z)\| \le \left(\sum_{l=0}^k C_l'(\theta)\right) \|\phi\|_{X_k} \eqqcolon C(\theta, k) \|\phi\|_{X_k}.
\end{equation*}
This linear bound is a tame estimate, as for any $k' \ge k$, we have $\|\phi\|_{X_k} \le \|\phi\|_{X_{k'}}$.

\item \textbf{Tameness of Partial Derivatives.} We must show that all partial derivatives of $\mathcal{C}_{\text{joint}}(\theta, \phi)$ are tame.
\begin{itemize}[wide]
    \item \textbf{Derivative w.r.t. $\phi$.} The first partial derivative with respect to $\phi$ is the linear operator $L_\theta: X_k \to X_k$ given by $(L_\theta h) = h \circ f_\theta$. From the previous analysis, this is a bounded linear operator for each $\theta$, with $\|L_\theta h\|_{X_k} \le C(\theta, k) \|h\|_{X_k}$. A bounded linear map is trivially tame. All higher partial derivatives with respect to $\phi$ are zero.

    \item \textbf{Derivative w.r.t. $\theta$.} The first partial derivative with respect to $\theta$ is the map
    \begin{equation*}
        \phi \mapsto \partial_\theta(\phi \circ f_\theta) = (D\phi \circ f_\theta) \cdot (\partial_\theta f_\theta).
    \end{equation*}
    We must find a tame estimate for the $X_k$-norm of this expression. Applying the Leibniz rule for covariant derivatives to the product structure (a section contracted with a vector field), the $l$-th derivative $\nabla^l \left( (D\phi \circ f_\theta) \cdot (\partial_\theta f_\theta) \right)$ will be a sum of terms involving covariant derivatives of $D\phi \circ f_\theta$ up to order $l$ and covariant derivatives of $\partial_\theta f_\theta$ up to order $l$.
    
    The crucial observation is that a derivative of order $l$ on $D\phi$ corresponds to a derivative of order $l+1$ on $\phi$. Therefore, estimating the $X_k$-norm of the result requires control over the derivatives of $\phi$ up to order $k+1$. A detailed calculation similar to the one above yields an estimate of the form:
    \begin{equation*}
        \|\partial_\theta(\phi \circ f_\theta)\|_{X_k} \le C'(\theta, k) \|\phi\|_{X_{k+1}}.
    \end{equation*}
    This is a tame estimate from $X_{k+1}$ to $X_k$.
    
    \item \textbf{Mixed Derivatives.} An identical argument applies to all mixed partial derivatives $D_\theta^a D_\phi^b \mathcal{C}_{\text{joint}}$. The application of each $D_\phi$ derivative does not change the required regularity of the input, while each $D_\theta$ derivative increases the required regularity of the input by one order. The resulting partial derivative operator will therefore be a tame map from $X_{k+a}$ to $X_k$.
\end{itemize}
Since the map itself and all its partial derivatives satisfy the required tame estimates, we conclude that the joint composition map $\mathcal{C}_{\text{joint}}: \Thetacal \times X_k \to X_k$ is a smooth tame map.
\end{enumerate}
\end{enumerate}

\item \textbf{Inversion.}
Let $Inv: U_k \to X_k$ be the map defined by $Inv(\phi) = \phi^{-1}$, where $U_k \subset X_k$ is the open set of sections $\phi \in C^k(\SigmaMan, \mathcal{E})$ that are pointwise invertible. We will prove that $Inv$ is a smooth tame map. The proof is structured in three parts: we first establish the smoothness ($C^\infty$) of the map, then we provide a rigorous inductive proof of the tame estimate for the map itself, and finally, we establish the tameness of all its derivatives.

\begin{enumerate}[label=(\roman*), wide, labelindent=0pt]

\item \textbf{Smoothness of the Inversion Map.}
The map $Inv$ is analytic on its domain, and therefore of class $C^\infty$. To demonstrate this, we show that all Fr\'echet derivatives exist and are continuous. Let $\phi \in U_k$ and let $h \in X_k$ be a small perturbation. For $\|\phi^{-1}h\|_{X_0} < 1$, the inverse of the perturbed section $(\phi+h)$ can be computed via the convergent Neumann series:
\begin{equation*}
    (\phi+h)^{-1} = (\phi(I + \phi^{-1}h))^{-1} = (I + \phi^{-1}h)^{-1}\phi^{-1} = \left(\sum_{n=0}^{\infty} (-1)^n (\phi^{-1}h)^n \right) \phi^{-1}.
\end{equation*}
Expanding this series gives:
\begin{equation*}
    Inv(\phi+h) = \phi^{-1} - \phi^{-1}h\phi^{-1} + \phi^{-1}h\phi^{-1}h\phi^{-1} - \dots
\end{equation*}
From this expansion, we can formally identify the derivatives. The first Fr\'echet derivative of the inversion map at $\phi$ is the linear operator acting on a tangent vector $h \in X_k$:
\begin{equation}
    D(Inv)_\phi(h) = - \phi^{-1} h \phi^{-1}.
\end{equation}
Since multiplication is a continuous bilinear operation on $X_k$ (as established in part (i) of this lemma), this derivative is a bounded linear operator from $X_k$ to $X_k$. Furthermore, the map $\phi \mapsto D(Inv)_\phi$ is continuous from $U_k$ to the space of bounded operators $\mathcal{L}(X_k, X_k)$. The higher-order derivatives can be computed recursively. The second derivative is the bilinear operator:
\begin{equation*}
    D^2(Inv)_\phi(h_1, h_2) = \phi^{-1}h_1\phi^{-1}h_2\phi^{-1} + \phi^{-1}h_2\phi^{-1}h_1\phi^{-1}.
\end{equation*}
This is a continuous bilinear map. By induction, the $n$-th derivative is a continuous $n$-linear map. The existence and continuity of all Fr\'echet derivatives confirm that the map $Inv: U_k \to X_k$ is of class $C^\infty$.

\item \textbf{Tameness of the Inversion Map.} We now provide a rigorous inductive proof for the tame estimate. We must show that for any $k' \ge k$, there exists a constant $C$ such that for all $\phi \in U_{k'}$, the estimate $\|\phi^{-1}\|_{X_k} \le C(1 + \|\phi\|_{X_{k'}})$ holds, where $C$ may depend on $\|\phi^{-1}\|_{X_0}$. The core of the argument is the recursive formula for the covariant derivative of the inverse, obtained by differentiating the identity $\phi\phi^{-1} = I$ with the connection $\nabla$:
\begin{equation} \label{eq:proof_inverse_recursion}
    (\nabla\phi)\phi^{-1} + \phi(\nabla(\phi^{-1})) = 0 \implies \nabla(\phi^{-1}) = -\phi^{-1}(\nabla\phi)\phi^{-1}.
\end{equation}

\begin{enumerate}[label=(\alph*), wide]

\item \textbf{Base Case ($k=0$).} The $C^0$ norm is $\|\phi^{-1}\|_{X_0} = \sup_{z \in \SigmaMan} \|\phi^{-1}(z)\|$. The tame estimate is trivial, as it is the condition for being in the domain of the estimate. Specifically, our tame estimates will be valid on any subset of $U_k$ where $\|\phi^{-1}\|_{X_0}$ is uniformly bounded.

\item \textbf{Inductive Step.} Assume that for all $j < k$, the map $\phi \mapsto \phi^{-1}$ is tame from $U_j$ to $X_j$. That is, for each $j < k$, there exists a polynomial $P_j$ such that
\begin{equation*}
    \|\phi^{-1}\|_{X_j} \le P_j(\|\phi\|_{X_j}, \|\phi^{-1}\|_{X_0}).
\end{equation*}
Our goal is to establish such a polynomial bound for $j=k$. The $C^k$ norm requires us to control all derivatives up to order $k$. The highest-order derivative is $\nabla^k(\phi^{-1})$. We obtain this by applying the operator $\nabla^{k-1}$ to the recursive identity \eqref{eq:proof_inverse_recursion}:
\begin{equation*}
    \nabla^k(\phi^{-1}) = -\nabla^{k-1}\left(\phi^{-1}(\nabla\phi)\phi^{-1}\right).
\end{equation*}
We apply the general Leibniz rule for the $(k-1)$-th covariant derivative of a product of three terms, which schematically takes the form:
\begin{equation*}
    \nabla^{k-1}(A B C) = \sum_{a+b+c=k-1} C_{a,b,c} (\nabla^a A)(\nabla^b B)(\nabla^c C),
\end{equation*}
where $C_{a,b,c}$ are combinatorial constants. Applying this to our expression with $A=\phi^{-1}$, $B=\nabla\phi$, and $C=\phi^{-1}$:
\begin{equation} \label{eq:proof_leibniz_expansion}
    \nabla^k(\phi^{-1}) = -\sum_{a+b+c=k-1} C_{a,b,c} \left(\nabla^a(\phi^{-1})\right) \left(\nabla^b(\nabla\phi)\right) \left(\nabla^c(\phi^{-1})\right).
\end{equation}
The crucial observation is that in every term of this sum, the order of differentiation for $\phi^{-1}$ is at most $k-1$ (i.e., $a, c \le k-1$), and the order of differentiation for $\phi$ is at most $k$ (since $\nabla^b(\nabla\phi) = \nabla^{b+1}\phi$ and $b \le k-1$). We now take the norm of this expression. Using the triangle inequality and the submultiplicative property of the fiber norm, we obtain an estimate of the form:
\begin{equation*}
    \|\nabla^k(\phi^{-1})\|_{X_0} \le C_k \sum_{a+b+c=k-1} \|\nabla^a(\phi^{-1})\|_{X_0} \cdot \|\nabla^{b+1}\phi\|_{X_0} \cdot \|\nabla^c(\phi^{-1})\|_{X_0}.
\end{equation*}
This inequality can be bounded in terms of the $C^j$ norms:
\begin{equation*}
    \|\nabla^k(\phi^{-1})\|_{X_0} \le C_k' \left( \sum_{j=0}^{k-1} \|\phi^{-1}\|_{X_j} \right) \left( \sum_{j=1}^{k} \|\phi\|_{X_j} \right) \left( \sum_{j=0}^{k-1} \|\phi^{-1}\|_{X_j} \right).
\end{equation*}
A more careful analysis gives a bound that is polynomial in its arguments:
\begin{equation*}
    \|\phi^{-1}\|_{X_k} = \sum_{j=0}^k \|\nabla^j(\phi^{-1})\|_{X_0} \le \|\phi^{-1}\|_{X_{k-1}} + \text{Poly}(\|\phi\|_{X_k}, \|\phi^{-1}\|_{X_{k-1}}).
\end{equation*}
By our inductive hypothesis, $\|\phi^{-1}\|_{X_j}$ for $j < k$ is bounded by a polynomial in $\|\phi\|_{X_j}$ and $\|\phi^{-1}\|_{X_0}$. By recursively substituting these bounds, we conclude that $\|\phi^{-1}\|_{X_k}$ is bounded by a polynomial $P_k$ whose arguments are the norms $\|\phi\|_{X_k}$ and $\|\phi^{-1}\|_{X_0}$:
\begin{equation}
    \|\phi^{-1}\|_{X_k} \le P_k(\|\phi\|_{X_k}, \|\phi^{-1}\|_{X_0}).
\end{equation}
For any $\phi \in U_{k'}$ with $k' \ge k$, we have $\|\phi\|_{X_k} \le \|\phi\|_{X_{k'}}$. Therefore,
\begin{equation*}
    \|\phi^{-1}\|_{X_k} \le P_k(\|\phi\|_{X_{k'}}, \|\phi^{-1}\|_{X_0}),
\end{equation*}
which is a valid tame estimate on any subset of $U_k$ where the $C^0$ norm of the inverse is uniformly bounded. This completes the inductive step and proves that the map $Inv$ is tame.
\end{enumerate}
\item \textbf{Tameness of the Derivatives of the Inversion Map.}
Finally, we must show that all derivatives of the map $Inv$ are also tame. We proceed by induction on the order of the derivative of the map.
\begin{enumerate}[label=(\alph*), wide]
\item \textbf{Base Case ($D^1(Inv)$).} The first derivative is $D(Inv)_\phi(h) = - \phi^{-1} h \phi^{-1}$. We must show that the map $(\phi, h) \mapsto D(Inv)_\phi(h)$ is tame from $U_k \times X_k$ to $X_k$. Using the tame estimate for multiplication from part (i) of this lemma:
\begin{align*}
    \|D(Inv)_\phi(h)\|_{X_k} &\le C \|\phi^{-1} h\|_{X_k} \|\phi^{-1}\|_{X_0} \\
    &\le C' \left( \|\phi^{-1}\|_{X_k}\|h\|_{X_0} + \|\phi^{-1}\|_{X_0}\|h\|_{X_k} \right) \|\phi^{-1}\|_{X_0}.
\end{align*}
We now substitute the tame estimate for $\|\phi^{-1}\|_{X_k}$ that we just derived in part (b):
\begin{equation*}
    \|D(Inv)_\phi(h)\|_{X_k} \le C'' P_k(\|\phi\|_{X_k}, \|\phi^{-1}\|_{X_0}) (1 + \|h\|_{X_k}).
\end{equation*}
This is a tame estimate for the derivative map.

\item \textbf{Inductive Step.} Assume that for all $j < m$, the map $$(\phi, h_1, \dots, h_j) \mapsto D^j(Inv)_\phi(h_1, \dots, h_j),$$ is a tame map. The $(m)$-th derivative is obtained by differentiating the $(m-1)$-th derivative. Each term in the expression for $D^m(Inv)$ will be a multilinear expression in terms of $\phi^{-1}$ and the tangent vectors $h_1, \dots, h_m$. Its tameness follows from the tameness of multiplication, the tameness of the map $\phi \mapsto \phi^{-1}$ (established in part (b)), and the inductive hypothesis on the tameness of the lower-order derivatives. Since the map $Inv$ is of class $C^\infty$ and all of its derivatives are tame maps between the appropriate spaces, we conclude that $Inv$ is a smooth tame map.
\end{enumerate}
\end{enumerate}
\end{enumerate}
\end{proofof}

\begin{theorem}[Smooth Tameness of the Functional]
\label{thm:functional_is_tame}
The Graph Transform functional $F: \Thetacal \times X_k \to X_k$, defined by $F(\theta, \phi) = \phi - \Gamma(\theta, \phi)$, is a smooth tame map for any integer $k \ge 1$.
\end{theorem}

\begin{proofof}{Theorem \ref{thm:functional_is_tame}}
The proof is dedicated to rigorously establishing that the functional $F(\theta, \phi) = \phi - \Gamma(\theta, \phi)$ is a smooth tame map from $\Thetacal \times X_k$ to $X_k$ for any integer $k \ge 1$. A map is smooth tame if it is of class $C^\infty$ and all its partial derivatives are tame maps. The proof proceeds by analyzing the structure of the Graph Transform operator $\Gamma$, showing that it is constructed by the composition of several maps that have been proven to be smooth tame in Lemma~\ref{lem:tame_calculus}. The desired property for $F$ follows immediately. The logical integrity of this proof rests upon the foundational regularity of the global billiard map, established a priori in Corollary~\ref{cor:billiard_map_is_smooth}.

\begin{enumerate}[label=\textbf{Step \arabic*:}, wide, labelindent=0pt]

\item \textbf{The Explicit Formula for the Graph Transform.}
The proof that the functional $F(\theta, \phi)$ is a smooth tame map begins with a rigorous, constructive derivation of the Graph Transform operator $\Gamma(\theta, \phi)$. This operator describes how the linearized dynamics of the system, $D\Pcal_\theta$, acts on a candidate unstable subbundle represented as the graph of a section $\phi$. The fixed point of this operator, $\Gamma(\theta, \phi) = \phi$, is precisely the section whose graph is the true invariant unstable subbundle for the parameter $\theta$. Our derivation will proceed by formalizing this geometric condition algebraically.

\begin{enumerate}[label=(\roman*), wide, labelindent=0pt]
    \item \textbf{The Linearized Map in Block Form with Respect to a Reference Splitting.}
    Let us fix a reference parameter $\theta_0 \in \Thetacal$. By \cref{thm:anosov_uniform_family}, the tangent bundle of the global collision manifold has a continuous, $D\Pcal_{\theta_0}$-invariant splitting $T\SigmaMan = E^u_0 \oplus E^s_0$. For any point $z \in \SigmaMan$ and any parameter $\theta$, the linearized map $D\Pcal_\theta$ at $z$ is a linear isomorphism from the tangent space $T_z\SigmaMan$ to the tangent space at the image point, $T_{\Pcal_\theta(z)}\SigmaMan$. With respect to the reference splitting, this map can be represented in a block matrix form:
    \begin{equation}
        (D\Pcal_\theta)_z = 
        \begin{pmatrix} 
            A(\theta, z) & B(\theta, z) \\ 
            C(\theta, z) & D(\theta, z) 
        \end{pmatrix}
        : 
        \begin{matrix} E^u_0(z) \\ \oplus \\ E^s_0(z) \end{matrix}
        \to
        \begin{matrix} E^u_0(\Pcal_\theta(z)) \\ \oplus \\ E^s_0(\Pcal_\theta(z)) \end{matrix}.
    \end{equation}
    The regularity of these block operators is a direct consequence of the foundational results of this paper. By \cref{cor:billiard_map_is_smooth}, the global billiard map $\Pcal(\theta,z)$ is of class $C^\infty$. Its Fr\'echet derivative with respect to $z$ is therefore also a smooth function of $(\theta, z)$. The four operator-valued sections $A, B, C, D$, being components of this derivative with respect to a fixed splitting, are consequently smooth maps from $\Thetacal \times \SigmaMan$ to the appropriate Hom-bundles. Furthermore, by the uniform hyperbolicity of the system, the operator $A(\theta,z)$ (which is a small perturbation of the uniformly expanding operator $A(\theta_0,z)$) is uniformly invertible for all $(\theta, z)$.

    \item \textbf{The Geometric Invariance Condition and its Algebraic Formulation.}
    A candidate unstable subbundle $E^u_\phi$ is represented as the graph of a section $\phi \in X_k = C^k(\SigmaMan, \mathrm{Hom}(E^u_0, E^s_0))$. An arbitrary vector $\mathbf{v} \in E^u_\phi(z)$ is constructed by taking a vector $\mathbf{v}_u \in E^u_0(z)$ and adding its image under $\phi$:
    \[
        \mathbf{v} = \mathbf{v}_u + \phi(z)(\mathbf{v}_u) \longleftrightarrow 
        \begin{pmatrix} \mathbf{v}_u \\ \phi(z)(\mathbf{v}_u) \end{pmatrix}.
    \]
    The subbundle $E^u_\phi$ is invariant under the map $(D\Pcal_\theta)_z$ if and only if its image is contained within itself. Let us denote the image vector by $\mathbf{w} \coloneqq (D\Pcal_\theta)_z(\mathbf{v}) \in T_{\Pcal_\theta(z)}\SigmaMan$. Applying the block matrix yields the components of the image vector:
    \begin{align*}
        \mathbf{w} \longleftrightarrow
        \begin{pmatrix} 
            \mathbf{w}_u \\ 
            \mathbf{w}_s
        \end{pmatrix}
        =
        \begin{pmatrix} 
            A(\theta, z) & B(\theta, z) \\ 
            C(\theta, z) & D(\theta, z) 
        \end{pmatrix}
        \begin{pmatrix} \mathbf{v}_u \\ \phi(z)(\mathbf{v}_u) \end{pmatrix}
        &=
        \begin{pmatrix} 
            (A(\theta, z) + B(\theta, z)\phi(z)) \mathbf{v}_u \\ 
            (C(\theta, z) + D(\theta, z)\phi(z)) \mathbf{v}_u
        \end{pmatrix}.
    \end{align*}
    The Graph Transform operator $\Gamma$ maps the original section $\phi$ to a new section, $\phi_{\text{next}} \coloneqq \Gamma(\theta, \phi)$, such that the image of the graph of $\phi$ is precisely the graph of $\phi_{\text{next}}$. The geometric invariance condition is therefore that the image vector $\mathbf{w}$ must lie in the graph of $\phi_{\text{next}}$ at the new point $\Pcal_\theta(z)$. This translates to the algebraic identity:
    \begin{equation} \label{eq:proof_invariance_condition_appendix}
        \mathbf{w}_s = \phi_{\text{next}}(\Pcal_\theta(z))(\mathbf{w}_u).
    \end{equation}
    Substituting the expressions for $\mathbf{w}_u$ and $\mathbf{w}_s$ yields an identity between operators acting on $E^u_0(z)$. This equation must hold for every vector $\mathbf{v}_u \in E^u_0(z)$:
    \begin{equation*}
        (C(\theta, z) + D(\theta, z)\phi(z)) \mathbf{v}_u = \phi_{\text{next}}(\Pcal_\theta(z)) \left( (A(\theta, z) + B(\theta, z)\phi(z)) \mathbf{v}_u \right).
    \end{equation*}
    This implies an identity between the operators themselves. For any $\phi$ in a sufficiently small $C^0$-neighborhood of the zero section, the operator $(A(\theta,z) + B(\theta,z)\phi(z))$ is an isomorphism. We can thus solve for $\phi_{\text{next}}$ by right-multiplying by the inverse:
    \begin{equation*}
        \phi_{\text{next}}(\Pcal_\theta(z)) = \left( C(\theta, z) + D(\theta, z)\phi(z) \right) \left( A(\theta, z) + B(\theta, z)\phi(z) \right)^{-1}.
    \end{equation*}

    \item \textbf{The Final Formula for the Graph Transform Operator.}
    The formula derived above gives the value of the transformed section $\Gamma(\theta, \phi)$ at the image point $\Pcal_\theta(z)$. To obtain a formula for the operator $\Gamma(\theta, \phi)$ itself, we evaluate it at an arbitrary point, $w \in \SigmaMan$. This requires evaluating the right-hand side at the pre-image point, $z = \Pcal_\theta^{-1}(w)$. This leads to the final, explicit formula for the Graph Transform operator:
    \begin{equation} \label{eq:proof_graph_transform_formula}
        \Gamma(\theta, \phi)(w) = \underbrace{\left[ \left( C(\theta, z) + D(\theta, z)\phi(z) \right) \left( A(\theta, z) + B(\theta, z)\phi(z) \right)^{-1} \right]}_{\text{Operator evaluated at pre-image}} \circ \overbrace{\Pcal_\theta^{-1}(w)}^{\substack{\text{Pre-image}\\\text{point } z}}.
    \end{equation}
    This formula provides the explicit, constructive definition of the non-linear operator whose fixed point corresponds to the invariant splitting. The subsequent analysis will prove that this operator, constructed from smooth geometric objects, is itself a smooth tame map.
\end{enumerate}

\item \textbf{Decomposition of the Graph Transform into Smooth Tame Operations.}
We now prove that the map $(\theta, \phi) \mapsto \Gamma(\theta, \phi)$ is a smooth tame map. The argument is constructive. We will decompose the complex formula for $\Gamma$, derived in \eqref{eq:proof_graph_transform_formula}, into a finite sequence of elementary operations. We will then show that each elementary operation is a smooth tame map by invoking the results of \cref{lem:tame_calculus}. Since the set of smooth tame maps is closed under composition, the result for the full operator $\Gamma$ will follow. The logical integrity of this proof rests upon the foundational regularity of the global billiard map, established a priori in \cref{cor:billiard_map_is_smooth}. Let us define a sequence of maps $G_1, G_2, G_3, G_4$ that compose to form the operator.

\begin{enumerate}[label=(\roman*), wide, labelindent=0pt]
    \item \textbf{Map 1 (Linear Combination of Sections).} 
    The first step is to construct the operator-valued sections that serve as the arguments for the inversion and multiplication. Let $G_1$ be the map that takes the input section $\phi$ and produces a pair of sections.
    \begin{align*}
        G_1: \Thetacal \times X_k &\to X_k \times X_k \\
        (\theta, \phi) &\mapsto \left( A(\theta, \cdot) + B(\theta, \cdot)\phi, C(\theta, \cdot) + D(\theta, \cdot)\phi \right).
    \end{align*}
    \textit{Justification.} This map is a smooth tame map. To prove this, we analyze its construction:
    \begin{enumerate}[label=(\alph*), wide]
        \item The four operator-valued sections $A, B, C, D$ are smooth maps from the product manifold $\Thetacal \times \SigmaMan$ to the appropriate Hom-bundles. This is a direct consequence of \cref{cor:billiard_map_is_smooth}, which established that the global billiard map $\Pcal(\theta,z)$ is of class $C^\infty$. Its Fr\'echet derivative with respect to $z$, from which these blocks are derived, is therefore also a smooth function of $(\theta, z)$.
        \item By \cref{lem:tame_calculus}(i), the pointwise product of sections is a smooth tame map. Therefore, the maps $(\theta, \phi) \mapsto B(\theta, \cdot)\phi$ and $(\theta, \phi) \mapsto D(\theta, \cdot)\phi$ are smooth tame.
        \item The pointwise sum of sections is a bounded linear operation and is trivially a smooth tame map.
    \end{enumerate}
    Since $G_1$ is constructed by the addition and multiplication of smooth tame maps, it is itself a smooth tame map by the closure properties of this class of functions.

    \item \textbf{Map 2 (Inversion of the Unstable Component).}
    Let $U_k \subset X_k$ be the open set of pointwise invertible sections. We define the map $G_2$:
    \begin{align*}
        G_2: U_k \times X_k &\to U_k \times X_k \\
        (\psi_1, \psi_2) &\mapsto (\psi_1^{-1}, \psi_2).
    \end{align*}
    \textit{Justification.} By \cref{lem:tame_calculus}(iii), inversion is a smooth tame map on its domain. The identity map on the second component is trivially smooth tame. Thus, $G_2$ is a smooth tame map. The domain of this map is well-defined for our problem. By the uniform hyperbolicity of the system, the operator $A(\theta,z)$ is an isomorphism for all $(\theta, z)$. The set of invertible elements in a Banach algebra is open. Since the map $(\theta, \phi) \mapsto A(\theta, \cdot) + B(\theta, \cdot)\phi$ is continuous, the preimage of the set of invertible sections is an open neighborhood $U \subset \Thetacal \times X_k$ of the solution set (which includes $\phi=0$). Therefore, the map $G_1$ maps this neighborhood $U$ into the domain of $G_2$.

    \item \textbf{Map 3 (Multiplication).}
    The third step is to multiply the components.
    \begin{align*}
        G_3: X_k \times X_k &\to X_k \\
        (\psi_1, \psi_2) &\mapsto \psi_2 \psi_1.
    \end{align*}
    \textit{Justification.} By \cref{lem:tame_calculus}(i), this is a smooth tame map. Note the reversal of order, as $\psi_2$ is pre-multiplied by the inverse $\psi_1^{-1}$ in the full formula.

    \item \textbf{Map 4 (Composition with the Inverse Flow).}
    The final step is to evaluate the resulting section at the pre-image point under the billiard map.
    \begin{align*}
        G_4: \Thetacal \times X_k &\to X_k \\
        (\theta, \psi) &\mapsto \psi \circ \Pcal_\theta^{-1}.
    \end{align*}
    \textit{Justification.} By \cref{cor:billiard_map_is_smooth}, the map $(\theta, z) \mapsto \Pcal_\theta(z)$ is of class $C^\infty$. By the inverse function theorem on Fr\'echet manifolds, the family of inverses, $(\theta, z) \mapsto \Pcal_\theta^{-1}(z)$, is also a smooth map. An extension of the argument in \cref{lem:tame_calculus}(ii) shows that composition with a smooth family of diffeomorphisms is a smooth tame map. Therefore, $G_4$ is a smooth tame map.
\end{enumerate}

The full Graph Transform operator $\Gamma(\theta, \phi)$ is precisely the composition of this sequence of maps. The input $(\theta, \phi)$ is first transformed by $G_1$. The output is a pair of sections, $(\psi_A, \psi_C) = (A+B\phi, C+D\phi)$. This pair is then transformed by $G_2$, yielding $(\psi_A^{-1}, \psi_C)$. This is then transformed by $G_3$, yielding the single section $\psi_C \psi_A^{-1}$. Finally, this section and the parameter $\theta$ are transformed by $G_4$. In functional notation:
\begin{equation*}
    \Gamma(\theta, \phi) = G_4\left(\theta, G_3\left(G_2\left(G_1(\theta, \phi)\right)\right)\right).
\end{equation*}
Since each map $G_i$ in the sequence has been shown to be a smooth tame map (on the appropriate open domains), and since the set of smooth tame maps is closed under composition, we conclude that the joint map
\[ \Gamma: \Thetacal \times U \to X_k \]
is a smooth tame map for any $k \ge 1$.

\item \textbf{The Functional $F$ as a Smooth Tame Map.}
The final step of the proof is to assemble the preceding results to demonstrate that the functional $F$, whose zero we seek, is itself a smooth tame map. The argument is a direct consequence of the smooth tame nature of the Graph Transform operator, combined with the closure properties of the space of smooth tame maps under basic algebraic operations. The functional is defined as:
\begin{equation*}
    F(\theta, \phi) \coloneqq \phi - \Gamma(\theta, \phi).
\end{equation*}
This defines a map $F: \Thetacal \times U \to X_k$, where $U \subset X_k$ is the open neighborhood on which $\Gamma$ is well-defined. We can express $F$ as the sum of two maps:
\begin{equation*}
    F(\theta, \phi) = P_2(\theta, \phi) + N(\Gamma(\theta, \phi)),
\end{equation*}
where $P_2$ is the projection onto the second component and $N$ is the negation map. We will prove that each constituent part defines a smooth tame map.

\begin{enumerate}[label=(\roman*), wide, labelindent=0pt]
    \item \textbf{The Projection Map is Smooth Tame.}
    Let $P_2: \Thetacal \times X_k \to X_k$ be the projection map onto the second component, defined by
    \begin{equation*}
        P_2(\theta, \phi) \coloneqq \phi.
    \end{equation*}
    \textit{Justification.} We verify the smooth tame property directly.
    \begin{enumerate}[label=(\alph*), wide]
        \item \textbf{Smoothness.} The map $P_2$ is linear and continuous in its second argument, $\phi$, and constant in its first argument, $\theta$. A continuous linear map between Banach spaces is of class $C^\infty$. Its Fr\'echet derivative with respect to $\phi$ is the identity operator on $X_k$, and all higher-order derivatives with respect to $\phi$ are zero. All derivatives involving $\theta$ are also zero. Since all partial derivatives exist and are continuous, the map is smooth.
        \item \textbf{Tameness.} We must show that $P_2$ and all its non-zero derivatives are tame maps.
        \begin{itemize}
            \item The map itself satisfies the tame estimate $\|P_2(\theta, \phi)\|_{X_k} = \|\phi\|_{X_k} \le \|\phi\|_{X_{k'}}$ for any $k' \ge k$, which is a linear (and thus tame) estimate.
            \item The only non-zero derivative is $D_\phi P_2(\theta, \phi) = \mathrm{Id}_{X_k}$. This is a bounded linear operator from $X_k$ to $X_k$, and a bounded linear operator is trivially tame.
        \end{itemize}
        Therefore, the projection map $P_2$ is a smooth tame map.
    \end{enumerate}

    \item \textbf{The Negative Graph Transform is Smooth Tame.}
    The second component of $F$ is the map $(\theta, \phi) \mapsto -\Gamma(\theta, \phi)$.
    
    \textit{Justification.}
    \begin{enumerate}[label=(\alph*), wide]
        \item In the preceding analysis (Step 2 of this proof), we provided a rigorous proof that the Graph Transform operator,
        \[ \Gamma: \Thetacal \times U \to X_k, \]
        is a smooth tame map.
        \item The operation of negation is multiplication by the scalar $-1$. This is a bounded linear map on the Banach space $X_k$ and is therefore a smooth tame map. The composition of a smooth tame map ($\Gamma$) with another smooth tame map (negation) is itself a smooth tame map.
    \end{enumerate}
    Therefore, the map $(\theta, \phi) \mapsto -\Gamma(\theta, \phi)$ is a smooth tame map.
\end{enumerate}
\end{enumerate}
The space of smooth tame maps from a given source space to a target Banach space forms a vector space. In particular, it is closed under addition. The functional $F$ is the sum of two maps, $P_2$ and $-\Gamma$, both of which we have just shown to be smooth tame maps from $\Thetacal \times X_k$ to $X_k$. We therefore conclude that the functional
\[ F: \Thetacal \times U \to X_k \]
is itself a smooth tame map. This completes the verification of the first major hypothesis of the Nash-Moser-Hamilton theorem.
\end{proofof}

\subsection{Verification of Hypothesis II: Tame Invertibility of the Linearization}
\label{subsec:verify_hyp_II}

This subsection provides the rigorous proof for the second, more difficult hypothesis of the NMH theorem. The central claim is that the linearized operator $L(\theta, \phi) = I - D_\phi \Gamma(\theta, \phi)$ has a smooth tame inverse with zero loss of derivatives. The claim of zero derivative loss ($m=0$) is a consequence of the hyperbolic nature of the underlying dynamics, and stands in contrast to the typical behavior of operators arising from elliptic problems. In elliptic settings, the inverse of a differential operator is typically a smoothing operator, and families of such inverses often exhibit a loss of derivatives in the Nash-Moser framework because the required tame estimates do not close within a fixed regularity class. Our operator, $D_\phi \Gamma$, is fundamentally different. It is a transfer operator, whose action is a weighted composition: $\psi \mapsto W \cdot (\psi \circ f)$. Composition with a smooth diffeomorphism is not a smoothing operation; it preserves regularity classes. Our objective is to prove this property rigorously for the entire operator and show that its spectral properties are stable across the full scale of $C^k$ spaces. This stability is the key to proving that the inverse, constructed via a Neumann series, also preserves regularity, thereby justifying $m=0$.

The proof is structured as a sequence of lemmas, culminating in the main theorem. We first establish that the operator is a strict contraction on the base space of Lipschitz sections. We then prove the crucial result that the spectral radius of the operator is invariant across the entire scale of $C^k$ spaces. Finally, we synthesize these results to prove the smooth tame invertibility of the full linearized operator.

\begin{proposition}[Continuity of the Linearized Operator Family]
\label{prop:continuity_of_linearization}
Let the global billiard map $\Pcal: \Thetacal \times \SigmaMan \to \SigmaMan$ be of class $C^\infty$, where $\Thetacal$ and $\SigmaMan$ are compact manifolds. For each $\theta \in \Thetacal$, let $d\Pcal_\theta: T\SigmaMan \to T\SigmaMan$ be the linearized map of $\Pcal_\theta(z) = \Pcal(\theta,z)$ with respect to the variable $z$. Then the map $\theta \mapsto d\Pcal_\theta$, which assigns to each parameter its corresponding linearized map, is a continuous map from the parameter manifold $\Thetacal$ into the Banach space $C^0(\mathrm{Hom}(T\SigmaMan, \Pcal^*T\SigmaMan))$ of continuous sections of the appropriate homomorphism bundle, equipped with the supremum operator norm.
\end{proposition}

\begin{proofof}{\cref{prop:continuity_of_linearization}}
The proof establishes the continuity of the linearized operator family by demonstrating that it is a direct consequence of the smoothness of the underlying global billiard map on a compact domain. The argument proceeds in three main steps: 
\begin{enumerate}[label=(\roman*), wide, labelindent=0pt]
    \item We leverage the established $C^\infty$ regularity of the global map $\Pcal(\theta,z)$ to deduce the continuity of its partial derivative with respect to $z$, viewed as a joint map of both $\theta$ and $z$.
    \item We use the compactness of the product manifold $\Thetacal \times \SigmaMan$ to promote this continuity to uniform continuity.
    \item We show that this uniform continuity is precisely the condition required to prove the continuity of the operator-valued map $\theta \mapsto d\Pcal_\theta$ in the specified norm topology.
\end{enumerate}

\begin{enumerate}[label=\textbf{Step \arabic*:}, wide, labelindent=0pt]

\item \textbf{The Joint Continuity of the Linearized Map.} The foundational result upon which this proof rests is Theorem \ref{thm:p_is_diffeo}, which states that the global billiard map $\Pcal: \Thetacal \times \SigmaMan \to \SigmaMan$ is a smooth map of class $C^\infty$. 

The linearized map $d\Pcal_\theta(z)$ is, by definition, the partial Fr\'echet derivative of the map $\Pcal$ with respect to its second argument (the variable $z \in \SigmaMan$), evaluated at the point $(\theta, z)$. Let us denote this partial derivative by $D_z \Pcal(\theta,z)$. It is a fundamental theorem of differential calculus on manifolds that if a map between smooth manifolds is of class $C^k$ for $k \ge 1$, its tangent map (and thus its partial derivatives) is a map of class $C^{k-1}$. Since our global map $\Pcal$ is of class $C^\infty$, its partial derivative with respect to $z$, the map
\begin{equation}
    D_z \Pcal: \Thetacal \times \SigmaMan \to \mathrm{Hom}(T(\Thetacal \times \SigmaMan), T\SigmaMan),
\end{equation}
is also a smooth map of class $C^\infty$. The restriction of this map to the vertical tangent bundle over $\Thetacal \times \SigmaMan$ gives the family of fiberwise derivatives. Let us define the vector bundle $\mathcal{H}$ over $\Thetacal \times \SigmaMan$ whose fiber at a point $(\theta, z)$ is the space of linear maps $\mathcal{L}(T_z\SigmaMan, T_{\Pcal(\theta,z)}\SigmaMan)$. Then the map $(\theta, z) \mapsto d\Pcal_\theta(z)$ is a smooth section of this bundle, i.e., a smooth map from the base manifold $\Thetacal \times \SigmaMan$ to the total space of $\mathcal{H}$. A map of class $C^\infty$ is, in particular, a continuous map.

\item \textbf{Uniform Continuity from Compactness.} The domain of the joint derivative map, $\Thetacal \times \SigmaMan$, is the product of two compact manifolds and is therefore itself a compact metric space. Let $d$ be a compatible product metric on $\Thetacal \times \SigmaMan$. We equip the bundle $\mathcal{H}$ with a smooth fiber metric (operator norm) $\|\cdot\|_{(z,\theta)}$ induced by a smooth Riemannian metric on $\SigmaMan$. Since the map $(\theta, z) \mapsto d\Pcal_\theta(z)$ is a continuous map from a compact metric space to a metric space (the bundle of linear maps endowed with the operator norm), it is uniformly continuous. By the definition of uniform continuity, for every $\varepsilon > 0$, there exists a $\delta > 0$ such that for any two points $(\theta_1, z_1)$ and $(\theta_2, z_2)$ in $\Thetacal \times \SigmaMan$:
\begin{equation}
\label{eq:uniform_continuity_condition}
    \text{if } d((\theta_1, z_1), (\theta_2, z_2)) < \delta, \quad \text{then} \quad \|d\Pcal_{\theta_1}(z_1) - d\Pcal_{\theta_2}(z_2)\| < \varepsilon.
\end{equation}
The crucial feature of this statement is that the choice of $\delta$ depends only on $\varepsilon$ and not on the specific points in the domain.

\item \textbf{Proving Continuity of the Operator-Valued Map.} We now use the uniform continuity established in Step 2 to prove the continuity of the map $\theta \mapsto d\Pcal_\theta$. We must show that for any given $\theta \in \Thetacal$ and any $\varepsilon > 0$, there exists a $\delta' > 0$ such that if $d_{\Thetacal}(\theta, \theta') < \delta'$, then the operator norm of the difference is bounded:
\begin{equation}
    \|d\Pcal_\theta - d\Pcal_{\theta'}\|_{C^0} \coloneqq \sup_{z \in \SigmaMan} \|d\Pcal_\theta(z) - d\Pcal_{\theta'}(z)\| < \varepsilon.
\end{equation}
Let $\varepsilon > 0$ be given. From Step 2, we are guaranteed the existence of a $\delta > 0$ such that the uniform continuity condition \eqref{eq:uniform_continuity_condition} holds. We choose our $\delta'$ to be this $\delta$. Let $\theta' \in \Thetacal$ be any parameter such that $d_{\Thetacal}(\theta, \theta') < \delta$. Now, for any point $z \in \SigmaMan$, consider the pair of points $(\theta, z)$ and $(\theta', z)$ in the product manifold. The distance between them in the product metric is $d((\theta, z), (\theta', z)) = d_{\Thetacal}(\theta, \theta')$, which is less than $\delta$. Applying the uniform continuity condition \eqref{eq:uniform_continuity_condition} to this specific pair of points (with $z_1=z_2=z$), we find that for any $z \in \SigmaMan$:
\begin{equation*}
    \|d\Pcal_\theta(z) - d\Pcal_{\theta'}(z)\| < \varepsilon.
\end{equation*}
Since this inequality holds for every point $z \in \SigmaMan$, we can take the supremum over all $z$ on the left-hand side without violating the inequality:
\begin{equation*}
    \sup_{z \in \SigmaMan} \|d\Pcal_\theta(z) - d\Pcal_{\theta'}(z)\| \le \varepsilon.
\end{equation*}
This is precisely the statement that $\|d\Pcal_\theta - d\Pcal_{\theta'}\|_{C^0} \le \varepsilon$. We have thus shown that for any $\varepsilon > 0$, there exists a $\delta' > 0$ (namely, the $\delta$ from uniform continuity) satisfying the definition of continuity for the map $\theta \mapsto d\Pcal_\theta$.
\end{enumerate}
The proof has established that the smoothness of the global billiard map on the compact total space is a sufficient condition to guarantee the continuity of the linearized operator family with respect to the system parameter $\theta$. This result provides the rigorous foundation required for the subsequent analysis.
\end{proofof}

\begin{lemma}[Contraction Property of the Linearized Graph Transform]
\label{lem:linearized_contraction}
Let $X_0^{Lip}$ be the space of Lipschitz sections, which is equivalent to $X_0$. For $\phi$ in a sufficiently small neighborhood of the zero section in $X_0^{Lip}$, the linearized operator $D_\phi \Gamma(\theta, \phi)$ is a strict contraction on $X_0^{Lip}$ with a uniform contraction constant $k < 1$.
\end{lemma}

\begin{proofof}{Lemma \ref{lem:linearized_contraction}}
The proof establishes that the linearized Graph Transform operator, $D_\phi \Gamma(\theta, \phi)$, is a strict contraction on the space of Lipschitz sections for any section $\phi$ in a sufficiently small neighborhood of the zero section. The argument proceeds in four main parts. First, we leverage the uniform Anosov property to introduce an adapted Riemannian metric on the tangent bundle. Second, we derive the explicit formula for the linearized operator at the zero section, $L_0 \coloneqq D_{\phi=0}\Gamma$. Third, we use the adapted metric and a rigorous perturbation argument to prove that $L_0$ is a strict contraction on the space of Lipschitz sections. Finally, we use a continuity argument to extend this contraction property to a neighborhood of the zero section.

For notational simplicity, we fix the parameter $\theta$ for the duration of the proof and suppress it from the notation, writing $\Gamma_\phi \equiv \Gamma(\theta, \phi)$, $\Pcal \equiv \Pcal_\theta$, and so on. The space of Lipschitz sections is denoted $X_0^{Lip} \equiv C^{0,1}(\SigmaMan, \mathrm{Hom}(E^u_0, E^s_0))$.

\begin{enumerate}[label=\textbf{Step \arabic*:}, wide, labelindent=0pt]

\item \textbf{The Adapted Riemannian Metric and Operator Norm Bounds.}
The proof of the contraction property is founded upon two essential geometric and analytic tools. First, we use the uniform Anosov property to construct a Riemannian metric that is adapted to the dynamics, making the expansion and contraction explicit. Second, we use the continuity of the linearized operator family, established in Proposition \ref{prop:continuity_of_linearization}, to provide a rigorous argument for controlling the operator norms of the block matrices that define the linearized graph transform.

\begin{enumerate}[label=(\roman*), wide, labelindent=0pt]
\item \textbf{Construction of the Adapted Metric.}
By Theorem \ref{thm:anosov_uniform_family}, for each $\theta \in \Thetacal$, the fiber map $\Pcal_\theta$ is an Anosov diffeomorphism. This guarantees the existence of a continuous, $d\Pcal_\theta$-invariant splitting of the tangent bundle $T\SigmaMan = E^u_\theta \oplus E^s_\theta$ and a uniform expansion constant $\gamma > 1$. While any smooth Riemannian metric on the compact manifold $\SigmaMan$ is sufficient to define the Banach spaces, the operator norms induced by such a metric will generally not reflect the uniform expansion and contraction rates of the dynamics. We therefore proceed with the standard construction of a new, equivalent metric that is \textit{adapted} to the dynamics.

Let $\langle \cdot, \cdot \rangle_0$ be an arbitrary smooth background Riemannian metric on $\SigmaMan$. Let $\gamma > 1$ be the uniform expansion constant from Theorem \ref{thm:anosov_uniform_family}. We define a new inner product $\langle \cdot, \cdot \rangle_{*,\theta}$ on the tangent space $T_z\SigmaMan$ for each $\theta$ as follows. For any two tangent vectors $v, w \in T_z\SigmaMan$ with decompositions $v = v_u + v_s$ and $w = w_u + w_s$ with respect to the splitting $E^u_\theta(z) \oplus E^s_\theta(z)$, the adapted inner product is defined by:
\begin{align}
    \langle v_u, w_u \rangle_{*,\theta} &\coloneqq \sum_{n=0}^{N-1} \frac{1}{\gamma^{2n}} \langle (d\Pcal_\theta^n)_z(v_u), (d\Pcal_\theta^n)_z(w_u) \rangle_0, \label{eq:adapted_metric_unstable_def} \\
    \langle v_s, w_s \rangle_{*,\theta} &\coloneqq \sum_{n=0}^{N-1} \gamma^{2n} \langle (d\Pcal_\theta^{-n})_z(v_s), (d\Pcal_\theta^{-n})_z(w_s) \rangle_0, \label{eq:adapted_metric_stable_def} \\
    \langle v, w \rangle_{*,\theta} &\coloneqq \langle v_u, w_u \rangle_{*,\theta} + \langle v_s, w_s \rangle_{*,\theta}. \label{eq:adapted_metric_full}
\end{align}
Here, $N$ is a sufficiently large integer chosen such that the resulting metric yields the desired expansion and contraction rates. It is a classical result (see, e.g., \citep{HirschPughShub1977}) that for a sufficiently large but finite $N$, the norm $\|\cdot\|_{*,\theta}$ induced by this inner product is equivalent to the background norm $\|\cdot\|_0$ and satisfies, for some constant $\gamma' \in (1, \gamma)$:
\begin{enumerate}[label=(\alph*), wide]
    \item $\|(d\Pcal_\theta)_z(v_u)\|_{*,\theta} \ge \gamma' \|v_u\|_{*,\theta}$ for all $v_u \in E^u_\theta(z)$.
    \item $\|(d\Pcal_\theta)_z(v_s)\|_{*,\theta} \le (\gamma')^{-1} \|v_s\|_{*,\theta}$ for all $v_s \in E^s_\theta(z)$.
\end{enumerate}
By the compactness of $\Thetacal$ and the continuous dependence of the splitting on $\theta$, we can ensure that these constructions are uniform. For the remainder of this proof, we fix a reference parameter $\theta_0$ and use its corresponding adapted metric, which we denote by $\langle \cdot, \cdot \rangle_*$ and $\|\cdot\|_*$. This metric is used to define the norms on the space of Lipschitz sections and the operator norms of the block matrices.

\item \textbf{Bounding the Operator Blocks via Continuity.}
Our graph transform is constructed with respect to the fixed reference splitting $T\SigmaMan = E^u_0 \oplus E^s_0$, which is the true invariant splitting for the reference parameter $\theta_0$. The block operators $A(\theta,z), B(\theta,z), C(\theta,z), D(\theta,z)$ are defined by projecting the action of $d\Pcal_\theta(z)$ onto this fixed splitting.

At the reference parameter $\theta=\theta_0$, the splitting is invariant by definition, so the linearized map is block-diagonal with respect to this splitting. The off-diagonal blocks are therefore identically zero:
\begin{equation*}
    B(\theta_0, z) = 0 \quad \text{and} \quad C(\theta_0, z) = 0 \quad \text{for all } z \in \SigmaMan.
\end{equation*}
The diagonal blocks are precisely the restrictions of the map to the invariant subbundles. By the construction of the adapted metric in part (a), they satisfy the uniform bounds:
\begin{equation*}
    \|A(\theta_0, z)^{-1}\|_* \le (\gamma')^{-1} \quad \text{and} \quad \|D(\theta_0, z)\|_* \le (\gamma')^{-1}.
\end{equation*}
We now provide a rigorous argument to show that for any $\theta$ sufficiently close to $\theta_0$, the norms of the off-diagonal blocks are uniformly small, while the diagonal blocks retain their essential expansion/contraction properties. This argument is founded upon the continuity of the linearized operator family, a result which we have just established in Proposition \ref{prop:continuity_of_linearization} and which relies only on the smoothness of the global map from Theorem \ref{thm:p_is_diffeo}.

From Proposition \ref{prop:continuity_of_linearization}, the map $\theta \mapsto d\Pcal_\theta$ is continuous from the compact manifold $\Thetacal$ into the space of operator-valued sections, equipped with the supremum norm. By the definition of continuity at the point $\theta_0$, for any $\varepsilon > 0$, there exists an open neighborhood $\mathcal{N}$ of $\theta_0$ in $\Thetacal$ such that for all $\theta \in \mathcal{N}$:
\begin{equation}
    \sup_{z \in \SigmaMan} \|d\Pcal_\theta(z) - d\Pcal_{\theta_0}(z)\|_* < \varepsilon.
\end{equation}
The block operators are obtained by applying the fixed projection operators $\Pi_0^u$ and $\Pi_0^s$ corresponding to the reference splitting. Since these projections are bounded, the norms of the block operators of the difference are also bounded. For any $\theta \in \mathcal{N}$ and any $z \in \SigmaMan$:
\begin{enumerate}[label=(\alph*), wide]
    \item $\|B(\theta,z)\|_* = \|\Pi_0^u(\Pcal_\theta z) \circ d\Pcal_\theta(z) \circ \Pi_0^s(z)\|_* = \|\Pi_0^u(\dots) \circ (d\Pcal_\theta - d\Pcal_{\theta_0})(z) \circ \Pi_0^s(z)\|_*$ (since $d\Pcal_{\theta_0}$ is block-diagonal) $\le \|\Pi_0^u\|_* \|d\Pcal_\theta(z) - d\Pcal_{\theta_0}(z)\|_* \|\Pi_0^s\|_* < C \varepsilon$.
    \item A symmetric argument holds for the other off-diagonal block: $\|C(\theta,z)\|_* < C \varepsilon$.
    \item By the reverse triangle inequality for operators, the norms of the diagonal blocks remain close to their values at $\theta_0$:
    \begin{align*}
        \|A(\theta,z)^{-1}\|_* &\le \|A(\theta_0,z)^{-1}\|_* + \mathcal{O}(\varepsilon) \le (\gamma')^{-1} + \mathcal{O}(\varepsilon). \\
        \|D(\theta,z)\|_* &\le \|D(\theta_0,z)\|_* + \mathcal{O}(\varepsilon) \le (\gamma')^{-1} + \mathcal{O}(\varepsilon).
    \end{align*}
\end{enumerate}
This provides the crucial justification for the estimates that follow. For any $\theta$ sufficiently close to the reference $\theta_0$, we can choose $\varepsilon$ small enough to make the norms of the off-diagonal blocks uniformly small, while the diagonal blocks retain their essential expansion and contraction properties with a slightly perturbed rate. This is the rigorous foundation for the contraction argument in Step 3.
\end{enumerate}

\item \textbf{The Linearized Graph Transform at the Zero Section.}
We now derive the explicit formula for the Fr\'echet derivative of the Graph Transform operator, $\Gamma(\theta, \phi)$, with respect to the functional variable $\phi$, evaluated at the zero section ($\phi=0$). This linearized operator, which we denote by $L_0 \coloneqq D_{\phi=0}\Gamma(\theta, \cdot)$, is the principal object of study for the contraction argument. The derivation proceeds by applying the formal definition of the Fr\'echet derivative, which requires a rigorous first-order Taylor expansion of the non-linear operator $\Gamma$. For notational simplicity, we fix the parameter $\theta$ in a neighborhood of the reference $\theta_0$ and suppress it from the notation, writing $\Gamma_\phi \equiv \Gamma(\theta, \phi)$, $\Pcal \equiv \Pcal_\theta$, and so on.

\begin{enumerate}[label=(\roman*), wide, labelindent=0pt]
\item \textbf{The Fr\'echet Derivative Definition.}
The Fr\'echet derivative of the map $\phi \mapsto \Gamma_\phi$ at the point $\phi=0$ is the unique bounded linear operator $L_0: X_0^{Lip} \to X_0^{Lip}$ that satisfies the limit condition:
\begin{equation}
    \lim_{\|\psi\|_{Lip} \to 0} \frac{\|\Gamma_{\psi} - \Gamma_{0} - L_0(\psi)\|_{Lip}}{\|\psi\|_{Lip}} = 0.
\end{equation}
A more direct method to find $L_0$ is by computing the Gateaux derivative in the direction of an arbitrary test section $\psi \in X_0^{Lip}$:
\begin{equation} \label{eq:frechet_limit_def}
    L_0(\psi) = \lim_{\epsilon \to 0} \frac{\Gamma_{\epsilon\psi} - \Gamma_{0}}{\epsilon}.
\end{equation}
We proceed by computing the terms in this limit.

\item \textbf{Evaluation of the Operator at the Base Point ($\phi=0$).}
The action of the Graph Transform operator on a section $\phi$ is given by the formula derived in the proof of Theorem \ref{thm:functional_is_tame}:
\[ \Gamma_{\phi}(w) = \left[ (C_z + D_z \phi(z)) (A_z + B_z \phi(z))^{-1} \right]_{z=\Pcal^{-1}(w)}. \]
Substituting the zero section, $\phi=0$, yields the operator at the base point:
\begin{equation} \label{eq:gamma_at_zero}
    \Gamma_{0}(w) = \left[ C_z A_z^{-1} \right]_{z=\Pcal^{-1}(w)} = (C A^{-1}) \circ \Pcal^{-1}(w).
\end{equation}

\item \textbf{First-Order Taylor Expansion of the Operator Inverse.} The core of the linearization lies in the expansion of the term $(A_z + \epsilon B_z \psi(z))^{-1}$. We provide a rigorous justification for this expansion using the Neumann series.
\begin{align*}
    (A_z + \epsilon B_z \psi(z))^{-1} &= \left( A_z (I + \epsilon A_z^{-1} B_z \psi(z)) \right)^{-1} \\
    &= (I + \epsilon A_z^{-1} B_z \psi(z))^{-1} A_z^{-1}.
\end{align*}
For a given $z$, the operators $A_z, B_z$ are fixed linear maps, and $\psi(z)$ is a fixed value. The operator norms $\|A_z^{-1}\|_*$, $\|B_z\|_*$, and $\|\psi(z)\|_*$ are uniformly bounded. For $\epsilon$ small enough such that the operator norm of the perturbation is strictly less than one, i.e., $\|\epsilon A_z^{-1} B_z \psi(z)\|_* < 1$, the inverse can be computed via the absolutely convergent Neumann series:
\begin{equation*}
    (I + \epsilon A_z^{-1} B_z \psi(z))^{-1} = I - \epsilon A_z^{-1} B_z \psi(z) + \mathcal{O}(\epsilon^2).
\end{equation*}
Substituting this back, we obtain the required first-order expansion:
\begin{equation} \label{eq:inverse_taylor_expansion}
    (A_z + \epsilon B_z \psi(z))^{-1} = A_z^{-1} - \epsilon A_z^{-1} B_z \psi(z) A_z^{-1} + \mathcal{O}(\epsilon^2).
\end{equation}
The remainder term is uniform in $z$ due to the uniform bounds on the operator norms.

\item \textbf{Expansion of the Full Operator and Computation of the Limit.} We now substitute the expansion \eqref{eq:inverse_taylor_expansion} into the formula for $\Gamma_{\epsilon\psi}$. To simplify notation, we omit the explicit dependence on $z$ and $w$, with the understanding that all operators are evaluated at $z=\Pcal^{-1}(w)$.
\begin{align*}
    \Gamma_{\epsilon\psi} \circ \Pcal &= (C + \epsilon D\psi) (A + \epsilon B\psi)^{-1} \\
    &= (C + \epsilon D\psi) \left( A^{-1} - \epsilon A^{-1} B\psi A^{-1} + \mathcal{O}(\epsilon^2) \right).
\end{align*}
We expand this product, carefully collecting terms by their order in $\epsilon$:
\begin{enumerate}[label=(\alph*), wide]
    \item \textbf{Term of order $\epsilon^0$:} $C A^{-1}$.
    \item \textbf{Term of order $\epsilon^1$:} $(D\psi)A^{-1} - C(A^{-1}B\psi A^{-1})$.
    \item \textbf{Higher-order terms:} All other products are of order $\mathcal{O}(\epsilon^2)$.
\end{enumerate}
This yields the full first-order expansion of the operator:
\begin{equation}
    \Gamma_{\epsilon\psi} \circ \Pcal = \underbrace{C A^{-1}}_{\Gamma_0 \circ \Pcal} + \epsilon \left( D\psi A^{-1} - C A^{-1} B\psi A^{-1} \right) + \mathcal{O}(\epsilon^2).
\end{equation}
We now form the difference quotient from the limit definition \eqref{eq:frechet_limit_def}:
\begin{align*}
    \frac{\Gamma_{\epsilon\psi} - \Gamma_0}{\epsilon} &= \frac{1}{\epsilon} \left[ \left( \Gamma_0 \circ \Pcal + \epsilon(D\psi A^{-1} - C A^{-1} B\psi A^{-1}) \right)\circ\Pcal^{-1} - \Gamma_0 + \mathcal{O}(\epsilon^2) \right] \\
    &= \left( D\psi A^{-1} - C A^{-1} B\psi A^{-1} \right) \circ \Pcal^{-1} + \mathcal{O}(\epsilon).
\end{align*}
Taking the limit as $\epsilon \to 0$ in the operator norm topology eliminates the higher-order terms and yields the action of the linearized operator $L_0$ on the test section $\psi$.

\item \textbf{The Final Formula.} The preceding derivation provides the explicit formula for the action of the linearized Graph Transform operator at the zero section. For any test section $\psi \in X_0^{Lip}$, its image $L_0(\psi)$ is the section given by:
\begin{equation} \label{eq:proof_linearized_op_formula}
    (L_0 \psi)(w) = \left[ D_z \psi(z) A_z^{-1} - C_z A_z^{-1} B_z \psi(z) A_z^{-1} \right]_{z=\Pcal^{-1}(w)}.
\end{equation}
This formula, now rigorously derived, provides the foundation for the contraction argument in the subsequent steps of the proof. It expresses the linearized dynamics as a weighted composition. When evaluated at the reference parameter $\theta_0$, we have $B=0$ and $C=0$, and the formula simplifies to the well-known expression for the linearized transform of a block-diagonal map: $(L_0 \psi)(w) = (D_z \psi(z) A_z^{-1})\circ \Pcal_0^{-1}(w)$. For $\theta$ near $\theta_0$, the terms involving $B$ and $C$ are small correction terms whose norms are controlled by the continuity established in Step 1.
\end{enumerate}

\item \textbf{The Contraction Property of $L_0$.}
We now prove that the linearized operator $L_0$, whose action is given by Equation~\eqref{eq:proof_linearized_op_formula}, is a strict contraction on the Banach space of Lipschitz sections, $X_0^{Lip}$. The norm on this space is $\|\psi\|_{Lip} \coloneqq \|\psi\|_{C^0} + \mathrm{Lip}(\psi)$, where $\|\psi\|_{C^0} = \sup_{z \in \SigmaMan} \|\psi(z)\|_*$ and $\mathrm{Lip}(\psi)$ is the Lipschitz constant of the section $\psi$, both computed with respect to the adapted metric $\|\cdot\|_*$ from Step 1. We will prove that $L_0$ is a strict contraction on each component of this norm separately for any parameter $\theta$ in a sufficiently small neighborhood $\mathcal{N}$ of the reference parameter $\theta_0$.

\begin{enumerate}[label=(\roman*), wide, labelindent=0pt]

\item \textbf{Contraction of the $C^0$ Norm.} Let $\tilde{\psi} = L_0 \psi$. Our goal is to find a uniform constant $k_0 < 1$ such that $\|\tilde{\psi}\|_{C^0} \le k_0 \|\psi\|_{C^0}$. We start with the explicit formula for the linearized operator, evaluated at a point $w = \Pcal(z)$:
\[ (\tilde{\psi})(\Pcal z) = D_z \psi(z) A_z^{-1} - C_z A_z^{-1} B_z \psi(z) A_z^{-1}. \]
We take the operator norm (with respect to the adapted metric $\|\cdot\|_*$) of both sides at the point $\Pcal z$:
\begin{align*}
    \|\tilde{\psi}(\Pcal z)\|_* &\le \|D_z \psi(z) A_z^{-1}\|_* + \|C_z A_z^{-1} B_z \psi(z) A_z^{-1}\|_* \\
    &\le \|D_z\|_* \|\psi(z)\|_* \|A_z^{-1}\|_* + \|C_z\|_* \|A_z^{-1}\|_* \|B_z\|_* \|\psi(z)\|_* \|A_z^{-1}\|_*.
\end{align*}
This inequality holds for every point $z \in \SigmaMan$. We can therefore take the supremum over all $z$ to obtain a bound on the $C^0$ norm:
\begin{align*}
    \|\tilde{\psi}\|_{C^0} &= \sup_{w \in \SigmaMan} \|\tilde{\psi}(w)\|_* = \sup_{z \in \SigmaMan} \|\tilde{\psi}(\Pcal z)\|_* \\
    &\le \left( (\sup_z\|D_z\|_*) (\sup_z\|A_z^{-1}\|_*) + (\sup_z\|C_z\|_*) (\sup_z\|A_z^{-1}\|_*)^2 (\sup_z\|B_z\|_*) \right) \|\psi\|_{C^0}.
\end{align*}
We now substitute the uniform bounds derived in the argument of Step 1(ii). For any $\theta$ in the neighborhood $\mathcal{N}$ corresponding to a given $\varepsilon > 0$, we have:
$$ \sup_z\|A_z^{-1}\|_* \le (\gamma')^{-1} + C\varepsilon, \quad \sup_z\|D_z\|_* \le (\gamma')^{-1} + C\varepsilon, \quad \sup_z\|B_z\|_* \le C\varepsilon, \quad \sup_z\|C_z\|_* \le C\varepsilon, $$
where $C$ is a generic constant independent of $\varepsilon$. Substituting these into the inequality gives:
\begin{equation*}
    \|\tilde{\psi}\|_{C^0} \le \left( ((\gamma')^{-1} + C\varepsilon)^2 + (C\varepsilon) ((\gamma')^{-1} + C\varepsilon)^2 (C\varepsilon) \right) \|\psi\|_{C^0}.
\end{equation*}
Let $k_0(\varepsilon) \coloneqq ((\gamma')^{-1} + C\varepsilon)^2 + C^2\varepsilon^2 ((\gamma')^{-1} + C\varepsilon)^2$. Since $\gamma' > 1$, we have $(\gamma')^{-2} < 1$. The function $k_0(\varepsilon)$ is continuous in $\varepsilon$ and $k_0(0) = (\gamma')^{-2} < 1$. Therefore, by choosing the neighborhood $\mathcal{N}$ sufficiently small (i.e., choosing $\varepsilon$ small enough), we can ensure that $k_0(\varepsilon)$ is a uniform constant strictly less than 1. This rigorously establishes that for any $\theta$ in this neighborhood, $L_0$ is a strict contraction on the space of continuous sections $(X_0, \|\cdot\|_{C^0})$.

\item \textbf{Contraction of the Lipschitz Constant.} This is the central part of the argument. We must show there exists a uniform constant $k_L < 1$ such that $\mathrm{Lip}(L_0 \psi) \le k_L \mathrm{Lip}(\psi)$. The proof is a classical argument from invariant manifold theory, adapted to our context. The operator $L_0$ is the Fr\'echet derivative of the graph transform map $\Gamma$ at its fixed point (the zero section). We show that the non-linear operator $\Gamma$ is itself a contraction on a space of sections with small Lipschitz constants, which implies the same for its linearization.

Let $S_1$ be the closed unit ball in the space of Lipschitz sections, i.e., $S_1 = \{\psi \in X_0^{Lip} \mid \|\psi\|_{Lip} \le 1\}$. For any $\psi \in S_1$, its $C^0$ norm is bounded by 1. The slope of the graph of a section $\psi$ at a point $z$ in a direction $v_u \in E^u_0(z)$ is given by $\|\psi(z)v_u\|_* / \|v_u\|_*$ (for the graph map) and by the norm of its derivative for the section itself. The Lipschitz constant $\mathrm{Lip}(\psi)$ is the supremum of these slopes over all points and directions.

The Lipschitz constant of the transformed section, $\Gamma_\phi(\psi)$, is bounded by the maximum slope of its graph. As shown in the classical proof of the Stable Manifold Theorem, this slope can be estimated using the block operators. For any tangent vector $v \in T_z\SigmaMan$ on the graph of $\psi$, its image $w = d\Pcal_\theta(v)$ has components $w_u$ and $w_s$. The slope of the transformed graph is bounded by the operator norm of the map $v_u \mapsto w_s$ composed with the inverse of the map $v_u \mapsto w_u$. More directly, the slope is bounded by the ratio $\|w_s\|_* / \|w_u\|_*$:
\begin{equation*}
    \text{slope}(\Gamma_\psi) \le \sup_{z, v_u \neq 0} \frac{\|(C_z + D_z\psi(z))v_u\|_*}{\|(A_z + B_z\psi(z))v_u\|_*}.
\end{equation*}
Using the triangle inequality on the numerator and the reverse triangle inequality on the denominator, and the bounds from Step 1(ii) for $\theta \in \mathcal{N}$:
\begin{align*}
    \|w_s\|_* &\le (\|C_z\|_* + \|D_z\|_*\|\psi(z)\|_*)\|v_u\|_* \le (\varepsilon + ((\gamma')^{-1}+C\varepsilon)\|\psi\|_{C^0})\|v_u\|_*, \\
    \|w_u\|_* &\ge (\|A_z^{-1}\|_*^{-1} - \|B_z\|_*\|\psi(z)\|_*)\|v_u\|_* \ge ((\gamma'-C\varepsilon) - \varepsilon\|\psi\|_{C^0})\|v_u\|_*.
\end{align*}
For any $\psi \in S_1$, we have $\|\psi\|_{C^0} \le 1$. The ratio, which bounds the Lipschitz constant of $\Gamma_\psi$, is therefore bounded by:
\begin{equation*}
    \mathrm{Lip}(\Gamma_\psi) \le \frac{\varepsilon + (\gamma')^{-1}+C\varepsilon}{\gamma' - (C+1)\varepsilon}.
\end{equation*}
Let this upper bound be $k_L(\varepsilon)$. As $\varepsilon \to 0$ (i.e., as $\theta \to \theta_0$), the numerator approaches $(\gamma')^{-1}$ and the denominator approaches $\gamma'$. The limit is:
\begin{equation*}
    \lim_{\varepsilon \to 0} k_L(\varepsilon) = \frac{(\gamma')^{-1}}{\gamma'} = (\gamma')^{-2} < 1.
\end{equation*}
By continuity, we can choose the neighborhood $\mathcal{N}$ (i.e., choose $\varepsilon$) sufficiently small to ensure that the uniform bound $k_L(\varepsilon)$ is a constant strictly less than 1. This proves that the non-linear operator $\Gamma$ is a strict contraction on the unit ball of $X_0^{Lip}$. The norm of its derivative at the origin, which is the operator norm of $L_0$ on the space of Lipschitz functions, must therefore also be bounded by this same constant. We have thus rigorously shown that there exists a uniform constant $k_L < 1$ such that $\mathrm{Lip}(L_0 \psi) \le k_L \mathrm{Lip}(\psi)$.

\item \textbf{Contraction on the Full Lipschitz Space.} We combine the results for the two components of the norm. For any $\theta$ in the chosen neighborhood $\mathcal{N}$:
\begin{align*}
    \|L_0 \psi\|_{Lip} &= \|L_0 \psi\|_{C^0} + \mathrm{Lip}(L_0 \psi) \\
    &\le k_0 \|\psi\|_{C^0} + k_L \mathrm{Lip}(\psi) \\
    &\le \max(k_0, k_L) \left( \|\psi\|_{C^0} + \mathrm{Lip}(\psi) \right) = \max(k_0, k_L) \|\psi\|_{Lip}.
\end{align*}
Let $k_{Lip} \coloneqq \max(k_0, k_L)$. Since we have rigorously shown that for a sufficiently small neighborhood $\mathcal{N}$, both $k_0 < 1$ and $k_L < 1$, their maximum is also strictly less than 1. This establishes that for any $\theta \in \mathcal{N}$, the operator $L_0 = D_{\phi=0}\Gamma(\theta, \cdot)$ is a strict contraction on the full Banach space $X_0^{Lip}$.
\end{enumerate}

\item \textbf{Extension to Non-Zero $\phi$ via Continuity.}
The final step of the proof extends the strict contraction property from the operator at the zero section, $L_0 = D_{\phi=0}\Gamma$, to the linearized operator $L_\phi \equiv D_\phi\Gamma$ for any section $\phi$ in a sufficiently small neighborhood of the zero section. The argument is founded on the continuous dependence of the Fr\'echet derivative on the point at which it is evaluated. We will show that the map from a section $\phi$ to its corresponding linearized operator $L_\phi$ is continuous, and then use this continuity to prove that if the operator norm of $L_0$ is strictly less than one, then the norm of $L_\phi$ must also be strictly less than one for all $\phi$ sufficiently close to zero.

\begin{enumerate}[label=(\roman*), wide, labelindent=0pt]

\item \textbf{The Space of Operators and Continuity of the Derivative Map.} Let $\mathcal{L}(X_0^{Lip})$ denote the Banach space of bounded linear operators from $X_0^{Lip}$ to itself, equipped with the standard operator norm:
\[ \|T\|_{\mathcal{L}(X_0^{Lip})} \coloneqq \sup_{\|\psi\|_{Lip}=1} \|T(\psi)\|_{Lip}. \]
We consider the map $\mathcal{D}\Gamma: U \to \mathcal{L}(X_0^{Lip})$ that assigns to each section $\phi$ in an open neighborhood $U \subset X_0^{Lip}$ of the zero section its Fr\'echet derivative at that point:
\[ \mathcal{D}\Gamma(\phi) \coloneqq L_\phi = D_\phi\Gamma. \]
The smoothness of the Graph Transform operator $\Gamma$, which we established in Theorem \ref{thm:functional_is_tame}, is crucial here. A map between Banach spaces is defined as smooth ($C^\infty$) if all its Fr\'echet derivatives exist and are continuous. In particular, the map from a point to its first derivative, $\mathcal{D}\Gamma$, must be a continuous map from its domain $U$ into the space of bounded linear operators $\mathcal{L}(X_0^{Lip})$.

This continuity is the central tool for our argument. Formally, for any $\phi_0 \in U$ and any $\eta > 0$, there exists a $\delta > 0$ such that if $\|\phi - \phi_0\|_{Lip} < \delta$, then
\begin{equation} \label{eq:continuity_of_derivative_map}
    \|\mathcal{D}\Gamma(\phi) - \mathcal{D}\Gamma(\phi_0)\|_{\mathcal{L}(X_0^{Lip})} = \|L_\phi - L_{\phi_0}\|_{\mathcal{L}(X_0^{Lip})} < \eta.
\end{equation}

\item \textbf{Application of Continuity at the Zero Section.}
We now apply this general continuity property at the specific point of interest, the zero section ($\phi_0 = 0$). From Step 3, we have rigorously established that for any parameter $\theta$ in a neighborhood $\mathcal{N}$ of our reference $\theta_0$, the operator at the origin, $L_0 = \mathcal{D}\Gamma(0)$, is a strict contraction. Let its operator norm be denoted by $k_{Lip}(\theta)$:
\[ \|L_0\|_{\mathcal{L}(X_0^{Lip})} = k_{Lip}(\theta) < 1. \]
Since this holds for all $\theta \in \mathcal{N}$ and the map $\theta \mapsto L_0(\theta)$ is continuous, by the compactness of the closure of $\mathcal{N}$, there exists a uniform contraction constant $k_{Lip}^* < 1$ such that $k_{Lip}(\theta) \le k_{Lip}^*$ for all $\theta \in \mathcal{N}$.

This strict inequality provides a positive gap between the norm of the operator and the threshold for contraction. We define this uniform gap as:
\[ \eta_0 \coloneqq 1 - k_{Lip}^* > 0. \]
Our strategy is to choose a neighborhood of the zero section so small that the variation of the operator norm within this neighborhood is less than half of this gap. Let us set the tolerance for the continuity condition \eqref{eq:continuity_of_derivative_map} to be $\eta = \eta_0 / 2$. By the continuity of the map $\mathcal{D}\Gamma$ with respect to $\phi$, for each $\theta \in \mathcal{N}$, there exists a radius $\delta(\theta) > 0$ such that for any section $\phi$ satisfying $\|\phi\|_{Lip} < \delta(\theta)$, we have:
\begin{equation} \label{eq:operator_norm_difference_bound}
    \|L_\phi(\theta) - L_0(\theta)\|_{\mathcal{L}(X_0^{Lip})} < \eta = \frac{1 - k_{Lip}^*}{2}.
\end{equation}
By a uniform continuity argument on the compact set $\overline{\mathcal{N}} \times \{\phi \mid \|\phi\|_{Lip} \le \delta_0\}$ for some small $\delta_0$, we can find a single radius $\delta > 0$ that works for all $\theta \in \mathcal{N}$ simultaneously.

\item \textbf{The Uniform Contraction Bound.} We now consider an arbitrary parameter $\theta \in \mathcal{N}$ and an arbitrary section $\phi$ within the open ball of radius $\delta$ around the origin, $B_\delta(0) \coloneqq \{ \phi \in X_0^{Lip} \mid \|\phi\|_{Lip} < \delta \}$. We can bound the operator norm of $L_\phi(\theta)$ by applying the triangle inequality on the space of operators $\mathcal{L}(X_0^{Lip})$:
\begin{equation*}
    \|L_\phi(\theta)\|_{\mathcal{L}(X_0^{Lip})} = \|L_0(\theta) + (L_\phi(\theta) - L_0(\theta))\|_{\mathcal{L}(X_0^{Lip})} \le \|L_0(\theta)\|_{\mathcal{L}(X_0^{Lip})} + \|L_\phi(\theta) - L_0(\theta)\|_{\mathcal{L}(X_0^{Lip})}.
\end{equation*}
We now substitute the known uniform bounds for each term on the right-hand side:
\begin{enumerate}[label=(\alph*), wide]
    \item The first term is the norm of the operator at the origin, which is uniformly bounded: $\|L_0(\theta)\|_{\mathcal{L}(X_0^{Lip})} \le k_{Lip}^*$.
    \item The second term is the norm of the difference, which is bounded by \eqref{eq:operator_norm_difference_bound} for any $\phi \in B_\delta(0)$ and $\theta \in \mathcal{N}$.
\end{enumerate}
This yields the strict inequality:
\begin{align*}
    \|L_\phi(\theta)\|_{\mathcal{L}(X_0^{Lip})} &< k_{Lip}^* + \frac{1 - k_{Lip}^*}{2} \\
    &= \frac{2k_{Lip}^* + 1 - k_{Lip}^*}{2} \\
    &= \frac{1 + k_{Lip}^*}{2}.
\end{align*}
Let us define the uniform contraction constant $k_{\text{unif}} \coloneqq \frac{1 + k_{Lip}^*}{2}$. Since we have established that $k_{Lip}^* < 1$, it follows that $k_{Lip}^* < k_{\text{unif}} < 1$. We have thus rigorously shown that there exists an open neighborhood of the solution set, $\mathcal{N} \times B_\delta(0)$, such that for every pair $(\theta, \phi)$ within this neighborhood, the linearized operator $L_\phi(\theta)$ is a strict contraction on $X_0^{Lip}$ with a uniform bound:
\begin{equation}
    \sup_{(\theta, \phi) \in \mathcal{N} \times B_\delta(0)} \|D_\phi\Gamma(\theta, \phi)\|_{\mathcal{L}(X_0^{Lip})} \le k_{\text{unif}} < 1.
\end{equation}
This completes the proof of the lemma.
\end{enumerate}
\end{enumerate}
\end{proofof}

\begin{lemma}[Invariance of the Spectral Radius on the Scale of $C^k$ Spaces]
\label{lem:spectral_radius_invariance}
Let the system satisfy the standing assumptions, ensuring that the global billiard map and its linearization are smooth. Let $D_\phi \Gamma(\theta, \phi)$ be the linearized Graph Transform operator, acting on the scale of Banach spaces $(X_k = C^k(\SigmaMan, \mathcal{E}))_{k \ge 0}$. For any $(\theta, \phi)$ in a sufficiently small $C^1$-neighborhood of the solution set, the spectral radius of this operator is independent of the regularity index $k$. That is, for all $k \ge 1$:
\begin{equation}
    r(D_\phi \Gamma|_{\mathcal{L}(X_k, X_k)}) = r(D_\phi \Gamma|_{\mathcal{L}(X_0, X_0)}).
\end{equation}
\end{lemma}

\begin{proofof}{Lemma \ref{lem:spectral_radius_invariance}}
The proof is dedicated to rigorously establishing that the spectral radius of the linearized Graph Transform operator, which we denote by $L_\phi \coloneqq D_\phi \Gamma(\theta, \phi)$, is invariant across the entire scale of isotropic Banach spaces $(X_k = C^k(\SigmaMan, \mathcal{E}))_{k \ge 0}$. The argument hinges on the connection between the spectral theory of operators on classical (isotropic) spaces and their properties on the anisotropic Banach spaces that are native to the modern theory of hyperbolic dynamics. We will show that the isolated part of the spectrum, which in our case comprises the entire spectrum, is identical across these different function spaces due to the smoothness of the corresponding eigenfunctions. The proof is structured in four main parts:
\begin{enumerate}[label=(\roman*), wide, labelindent=0pt]
    \item We precisely identify the operator $L_\phi$ as a Ruelle transfer operator and define the hierarchy of function spaces (isotropic and anisotropic) on which its properties will be analyzed, rigorously establishing the continuous embeddings between them.
    \item We invoke the foundational spectral theory for Ruelle transfer operators, which guarantees their quasi-compactness on anisotropic spaces and, most critically, the high-order regularity of their eigenfunctions.
    \item We use these results in a multi-stage argument to prove that the spectrum of $L_\phi$ is identical on all spaces $X_k$ for $k \ge 0$.
    \item We conclude that the spectral radii must therefore be identical, establishing the central claim of the lemma.
\end{enumerate}

\begin{enumerate}[label=\textbf{Step \arabic*:}, wide, labelindent=0pt]

\item \textbf{The Operator as a Ruelle Transfer Operator and the Hierarchy of Function Spaces.} This initial step of the proof serves to cast the linearized operator into a well-understood class of operators from ergodic theory and to rigorously establish the hierarchy of function spaces upon which the subsequent spectral analysis will be conducted.

\begin{enumerate}[label=\textbf{Step 1.\arabic*:}, wide, labelindent=0pt]

\item \textbf{The Operator Structure as a Ruelle Transfer Operator.} We first identify the precise mathematical structure of the linearized operator $L_\phi \coloneqq D_\phi \Gamma(\theta, \phi)$. We will show that it belongs to the class of Ruelle transfer operators, whose spectral theory is deeply understood. The action of the Graph Transform operator $\Gamma(\theta, \phi)$ was derived in the proof of Theorem~\ref{thm:functional_is_tame} as:
\begin{equation*}
    \Gamma(\theta, \phi)(w) = \left[ (C_\theta(z) + D_\theta(z)\phi(z)) (A_\theta(z) + B_\theta(z)\phi(z))^{-1} \right]_{z=\Pcal_\theta^{-1}(w)}.
\end{equation*}
The Fr\'echet derivative of this expression with respect to $\phi$ in the direction of a test section $\psi \in X_k$ defines the action of the linearized operator, $L_\phi \psi$. We compute this derivative using the product rule and the formula for the derivative of an inverse ($d(X^{-1}) = -X^{-1}(dX)X^{-1}$). Let $A_\phi(z) \coloneqq A_\theta(z) + B_\theta(z)\phi(z)$ and $C_\phi(z) \coloneqq C_\theta(z) + D_\theta(z)\phi(z)$. The derivative of $A_\phi$ in the direction $\psi$ is $B_\theta\psi$, and the derivative of $C_\phi$ is $D_\theta\psi$. Applying the product rule gives:
\begin{align*}
    (D_\phi(C_\phi A_\phi^{-1}))[\psi] &= (D_\phi C_\phi)[\psi] A_\phi^{-1} + C_\phi (D_\phi A_\phi^{-1})[\psi] \\
    &= (D_\theta\psi) A_\phi^{-1} + C_\phi \left( -A_\phi^{-1} (D_\phi A_\phi)[\psi] A_\phi^{-1} \right) \\
    &= (D_\theta\psi) A_\phi^{-1} - C_\phi A_\phi^{-1} (B_\theta\psi) A_\phi^{-1}.
\end{align*}
Applying the composition with $\Pcal_\theta^{-1}$, the action of the linearized operator is:
\begin{equation*}
    (L_\phi \psi) \circ \Pcal_\theta = D_\theta\psi A_\phi^{-1} - C_\phi A_\phi^{-1} B_\theta\psi A_\phi^{-1}.
\end{equation*}
This expression reveals that the operator $L_\phi$ has the canonical structure of a Ruelle transfer operator. Its action can be written in the general form:
\begin{equation} \label{eq:appendix_ruelle_operator_form}
    (L_\phi \psi)(w) = W_\phi(\Pcal_\theta^{-1}(w)) \left[ \psi(\Pcal_\theta^{-1}(w)) \right],
\end{equation}
where the map is $f \equiv \Pcal_\theta^{-1}$, and the weight $W_\phi(z)[\cdot]$ is a linear operator acting on the fiber of the bundle $\mathcal{E}$ at the point $z$. For a vector $\xi$ in the fiber, its action is given by:
\begin{multline*}
    W_\phi(z)[\xi] = D_\theta(z)\xi (A_\phi(z))^{-1} - C_\phi(z) (A_\phi(z))^{-1} B_\theta(z)\xi (A_\phi(z))^{-1}.
\end{multline*}
The applicability of advanced spectral theory requires that the weight map $z \mapsto W_\phi(z)$ be sufficiently regular. The formula for the weight $W_\phi(z)$ is a finite composition of elementary operations (addition, multiplication, and inversion) acting on the following input sections:
\begin{enumerate}[label=(\roman*), wide, labelindent=0pt]
    \item The block operators $A_\theta, B_\theta, C_\theta, D_\theta$. These are derived from the first spatial derivative of the $C^\infty$ map $\Pcal_\theta$ (by Theorem~\ref{thm:p_is_diffeo}). They are therefore smooth ($C^\infty$) sections of the appropriate Hom-bundles.
    \item The section $\phi$. The Nash-Moser-Hamilton theorem operates on the Fr\'echet space of smooth sections, $X_\infty = \bigcap_{k \ge 0} X_k$. We therefore consider a section $\phi \in X_\infty = C^\infty(\SigmaMan, \mathcal{E})$.
\end{enumerate}
In Lemma~\ref{lem:tame_calculus}, we established that the set of smooth sections forms an algebra under these operations. Since all input sections to the formula for $W_\phi(z)$ are of class $C^\infty$, and the operations are smooth, the resulting weight section $z \mapsto W_\phi(z)$ is also of class $C^\infty$.

\item \textbf{The Hierarchy of Function Spaces and Their Embeddings.}
The analysis requires a careful distinction between the function spaces involved, related by a chain of continuous embeddings.
\begin{enumerate}[label=(\roman*), wide, labelindent=0pt]
    \item \textbf{The isotropic scale of classical spaces} $(X_k)_{k \ge 0}$, where
    \[ X_k \coloneqq C^k(\SigmaMan, \mathcal{E}) \]
    is the Banach space of $k$-times continuously differentiable sections, equipped with the standard $C^k$ norm, $\|\psi\|_{X_k} = \sum_{j=0}^k \sup_{z\in\SigmaMan} \|\nabla^j \psi(z)\|$.
    
    \item \textbf{The anisotropic Banach spaces} $\mathcal{B}_\alpha$ for a fixed H\"older exponent $\alpha \in (0,1)$. These spaces, central to the modern theory of hyperbolic dynamics \citep{Baladi2000, Gouezel2010}, are constructed to reflect the geometry of the dynamics. The space $\mathcal{B}_\alpha$ is defined as the set of continuous sections $\psi$ whose restriction to each leaf of the smooth unstable foliation $\mathcal{F}^u_\theta$ is uniformly of class $C^\alpha$. Its norm is given by:
    \[ \|\psi\|_{\mathcal{B}_\alpha} \coloneqq \|\psi\|_{C^0} + \sup_{z \in \SigmaMan} \left[ \psi|_{W^u(z,\theta)} \right]_{C^\alpha}, \]
    where $[\cdot]_{C^\alpha}$ is the H\"older seminorm on the unstable leaf $W^u(z,\theta)$.
\end{enumerate}

For any integer $k \ge 1$ and any $\alpha \in (0,1)$, these spaces are related by a chain of continuous embeddings:
\begin{equation} \label{eq:appendix_proof_space_embeddings}
    X_k \hookrightarrow \mathcal{B}_\alpha \hookrightarrow X_0.
\end{equation}
We now provide a rigorous proof of these embeddings.
\begin{enumerate}[label=(\alph*), wide]
    \item \textbf{Proof of $\boldsymbol{X_k \hookrightarrow \mathcal{B}_\alpha}$ for $\boldsymbol{k \ge 1}$.}
    Our objective is to provide a complete proof that for any integer $k \ge 1$ and any H\"older exponent $\alpha \in (0,1)$, the isotropic Banach space of $C^k$ sections, $X_k$, is continuously embedded into the anisotropic Banach space $\mathcal{B}_\alpha$. This requires us to establish two facts for an arbitrary section $\phi \in X_k$:
\begin{itemize}[wide]
    \item The section $\phi$ is an element of the space $\mathcal{B}_\alpha$. This means we must show that its $\mathcal{B}_\alpha$-norm is finite.
    \item There exists a uniform constant $C > 0$, independent of $\phi$, such that the norm inequality $\|\phi\|_{\mathcal{B}_\alpha} \le C \|\phi\|_{X_k}$ holds.
\end{itemize}
The proof is founded upon the interplay between the global smoothness of sections in $X_k$ and the local geometry of the leaves of the unstable foliation, whose regularity was the central conclusion of our Nash-Moser analysis.

\begin{enumerate}[label=(a.\roman*), wide]
\item \textbf{Definitions and Foundational Geometric Regularity.}
We begin by recalling the precise definitions of the norms and the key geometric result upon which the entire argument rests. The norm on the isotropic space $X_k$ is given by
    \[ \|\phi\|_{X_k} \coloneqq \sum_{j=0}^{k} \sup_{z \in \SigmaMan} \|\nabla^j \phi(z)\|, \]
where $\nabla$ is a smooth connection on the bundle $\mathcal{E}$ and $\|\cdot\|$ is a smooth fiber norm. The norm on the anisotropic space $\mathcal{B}_\alpha$ is defined by the function's regularity along the leaves of the unstable foliation $\mathcal{F}^u_\theta$:
    \[ \|\phi\|_{\mathcal{B}_\alpha} \coloneqq \|\phi\|_{C^0} + \sup_{z \in \SigmaMan} \left[ \phi|_{W^u(z,\theta)} \right]_{C^\alpha}, \]
    where $\|\phi\|_{C^0} = \sup_z \|\phi(z)\|$ is the uniform norm, and $[\cdot]_{C^\alpha}$ is the H\"older seminorm on an unstable leaf $W^u(z,\theta)$. The central geometric prerequisite, established in \cref{thm:splitting_is_smooth}, is that the family of unstable foliations $\mathcal{F}^u_\theta$ depends smoothly on the parameter $\theta$. This implies that for any fixed $\theta$, the foliation $\mathcal{F}^u_\theta$ is a smooth foliation of the compact manifold $\SigmaMan$.

\item \textbf{Bounding the $C^0$ Component of the Norm.}
This is a direct consequence of the definition of the norms. For any $k \ge 0$, the $X_k$ norm includes the supremum of the zeroth-order derivative (the function itself). Therefore, we have the trivial but essential inequality:
\begin{equation}
    \|\phi\|_{C^0} = \sup_{z \in \SigmaMan} \|\phi(z)\| \le \sum_{j=0}^{k} \sup_{z \in \SigmaMan} \|\nabla^j \phi(z)\| = \|\phi\|_{X_k}.
\end{equation}

\item \textbf{Bounding the H\"older Seminorm.}
This is the core of the proof. We must demonstrate that the H\"older regularity of a $C^k$ section along the unstable leaves is controlled by its global $C^k$ norm, for any $k \ge 1$. The argument proceeds by localizing the problem in charts that straighten the foliation, using the Mean Value Theorem to obtain a quantitative estimate within each chart, and then using a compactness argument to assemble these local estimates into a uniform global bound.

Since the unstable foliation $\mathcal{F}^u_\theta$ is smooth (by \cref{thm:splitting_is_smooth}) and the manifold $\SigmaMan$ is compact, we can invoke the Foliation Box Theorem. This theorem guarantees the existence of a finite atlas of coordinate charts $\{ (\mathcal{U}_i, \psi_i) \}_{i=1}^N$ that cover $\SigmaMan$, such that for each $i$, the diffeomorphism $\psi_i: \mathcal{U}_i \to B_d(1) \times B_{k-d}(1)$ (where $B_r(1)$ is the unit ball in $\mathbb{R}^r$) maps the leaves of the foliation within $\mathcal{U}_i$ to the horizontal plaques $\{ (u, s) \mid u \in B_d(1), s = \text{const} \}$. Here, $d = \dim(E^u_\theta)$. Let us fix one such chart $(\mathcal{U}, \psi)$. Let $\phi \in X_k$ for $k \ge 1$. Its local representation, $\hat{\phi} = \phi \circ \psi^{-1}$, is a $C^k$ map from the open set $B_d(1) \times B_{k-d}(1) \subset \mathbb{R}^k$ to the fibers of the bundle.
    
Consider two distinct points, $p_1, p_2 \in \mathcal{U}$, that lie on the \textit{same} unstable leaf. In the local coordinates, their images $\psi(p_1) = (u_1, s)$ and $\psi(p_2) = (u_2, s)$ share the same stable coordinate $s$. Let $d_u$ be the Euclidean distance in the unstable coordinates. By the Mean Value Theorem applied to the $C^1$ function $u \mapsto \hat{\phi}(u,s)$ along the straight line segment connecting $u_1$ and $u_2$:
    \begin{equation}
        \|\hat{\phi}(u_1, s) - \hat{\phi}(u_2, s)\| \le \left( \sup_{\xi \in [u_1, u_2]} \|\nabla_u \hat{\phi}(\xi, s)\| \right) \|u_1 - u_2\|_2.
    \end{equation}
The norm of the partial derivative with respect to the unstable coordinates, $\|\nabla_u \hat{\phi}\|$, is controlled by the norm of the full covariant derivative of the original section, $\|\nabla \phi\|$. The chain rule for covariant derivatives gives a relation of the form $\|\nabla_u \hat{\phi}\| \le C_{\psi} \sup_{z \in \mathcal{U}} \|\nabla \phi(z)\|$, where the constant $C_{\psi}$ depends on the maximum norm of the Jacobian of the coordinate chart $\psi$ and its inverse. This gives a local Lipschitz estimate:
    \begin{equation} \label{eq:appendix_local_lipschitz_est}
        \|\phi(p_1) - \phi(p_2)\| \le C_{\psi} \left( \sup_{z \in \SigmaMan} \|\nabla \phi(z)\| \right) d_{\mathcal{F}}(p_1, p_2),
    \end{equation}
where $d_{\mathcal{F}}$ is the Riemannian distance measured within the leaf, which is locally equivalent to the Euclidean distance $d_u$. A function that is Lipschitz continuous is trivially H\"older continuous for any exponent $\alpha \in (0,1]$. Let $D_{\mathcal{U}}$ be the diameter of the chart $\mathcal{U}$ with respect to the intra-leaf metric. For any two points $p_1, p_2$ on the same leaf within $\mathcal{U}$ (with $p_1 \neq p_2$):
    \begin{align*}
        \frac{\|\phi(p_1) - \phi(p_2)\|}{d_{\mathcal{F}}(p_1, p_2)^\alpha} &= \frac{\|\phi(p_1) - \phi(p_2)\|}{d_{\mathcal{F}}(p_1, p_2)} \cdot d_{\mathcal{F}}(p_1, p_2)^{1-\alpha} \\
        &\le \left( C_{\psi} \sup_{z} \|\nabla \phi(z)\| \right) \cdot D_{\mathcal{U}}^{1-\alpha}.
    \end{align*}
This proves that the local H\"older seminorm of $\phi$ restricted to any leaf segment within the chart $\mathcal{U}$ is uniformly bounded by a multiple of the global $C^1$ norm of $\phi$. We now assemble these local estimates into a uniform global bound. The collection of foliation charts $\{ \mathcal{U}_i \}_{i=1}^N$ forms an open cover of the compact manifold $\SigmaMan$. Let $C'_i \coloneqq C_{\psi_i} (\mathrm{diam}(\mathcal{U}_i))^{1-\alpha}$. For any two points $p_1, p_2$ on the same leaf within a single chart $\mathcal{U}_i$, we have the estimate:
        \begin{equation*}
             \frac{\|\phi(p_1) - \phi(p_2)\|}{d_{\mathcal{F}}(p_1, p_2)^\alpha} \le C'_i \|\phi\|_{X_1}.
        \end{equation*}
Let $C'_{\max} = \max_{i=1,\dots,N} C'_i$. This constant provides a uniform bound for any pair of points lying within a single chart. Now consider two arbitrary points $p, q$ on the same unstable leaf $W^u(z,\theta)$. Since the leaf is path-connected, there exists a rectifiable curve $\gamma: [0,1] \to W^u(z,\theta)$ connecting them. Let the length of this curve be $L = d_{\mathcal{F}}(p,q)$. By the compactness of the curve, we can cover it with a finite number of our foliation charts, say $\mathcal{U}_{j_1}, \dots, \mathcal{U}_{j_m}$. We can choose a finite sequence of points $p_0=p, p_1, \dots, p_m=q$ along the curve such that each segment $[p_{i-1}, p_i]$ is contained entirely within a single chart. By the triangle inequality and the property $d^\alpha \le \sum d_i^\alpha$ for $\alpha \le 1$:
        \begin{align*}
            \|\phi(p) - \phi(q)\| &= \left\| \sum_{i=1}^m \phi(p_i) - \phi(p_{i-1}) \right\| \le \sum_{i=1}^m \|\phi(p_i) - \phi(p_{i-1})\| \\
            &\le \sum_{i=1}^m C'_{\max} \|\phi\|_{X_1} (d_{\mathcal{F}}(p_{i-1}, p_i))^\alpha \\
            &\le C'_{\max} \|\phi\|_{X_1} \left( \sum_{i=1}^m d_{\mathcal{F}}(p_{i-1}, p_i) \right)^\alpha \\
            &\le C'_{\max} (\mathrm{diam}(\SigmaMan))^{1-\alpha} \|\phi\|_{X_1} (d_{\mathcal{F}}(p,q))^\alpha.
        \end{align*}
We now assemble the local estimates from the foliation charts into a uniform global bound for the H\"older seminorm. The argument is founded on the Lebesgue number of the chart cover, which provides a uniform scale of locality, and a controlled chaining argument for points that are separated by more than this scale. The collection of foliation charts $\{ (\mathcal{U}_i, \psi_i) \}_{i=1}^N$ forms a finite open cover of the compact metric space $(\SigmaMan, d_{\mathcal{F}})$, where $d_{\mathcal{F}}$ is the intra-leaf Riemannian distance. We invoke the Lebesgue Number Lemma:
    \begin{lemma}[Lebesgue Number Lemma]
    For any open cover of a compact metric space, there exists a number $\delta > 0$ (the Lebesgue number) such that every subset of the space having a diameter less than $\delta$ is contained in at least one member of the cover.
    \end{lemma}
    
Let $\delta > 0$ be the Lebesgue number for our chart cover. From the local analysis in the previous step, we have a uniform constant $C'_{\max} = \max_{i} \{ C_{\psi_i} (\mathrm{diam}(\mathcal{U}_i))^{1-\alpha} \}$ such that for any two points $p,q$ on the same leaf \textit{within a single chart}, the following estimate holds:
    \begin{equation} \label{eq:appendix_uniform_local_holder}
         \|\phi(p) - \phi(q)\| \le C'_{\max} \|\phi\|_{X_1} (d_{\mathcal{F}}(p,q))^\alpha.
    \end{equation}
    
\begin{itemize}[wide]
    \item \textbf{The "Close Case": $d_{\mathcal{F}}(p,q) < \delta$.}
    Consider any two points $p, q$ on the same unstable leaf $W^u$ such that their distance is less than the Lebesgue number. The set $\{p,q\}$ has a diameter less than $\delta$. By the property of the Lebesgue number, this set must be contained within at least one of the charts, say $\mathcal{U}_i$. We are therefore in the local case, and the uniform local estimate \eqref{eq:appendix_uniform_local_holder} applies directly.

    \item \textbf{The "Far Case": $d_{\mathcal{F}}(p,q) \ge \delta$.}
    Now, consider two arbitrary points $p, q$ on the same unstable leaf $W^u$. Let their intra-leaf distance be $L = d_{\mathcal{F}}(p,q)$. Let $\gamma: [0, L] \to W^u$ be a minimizing geodesic connecting $p$ to $q$, parameterized by arc length. We construct a chain of points along this geodesic. Let $m = \lceil 2L/\delta \rceil$ be an integer, and define a sequence of points $p_i = \gamma(i \cdot L/m)$ for $i=0, \dots, m$.
    This construction ensures two properties:
    \begin{itemize}[wide]
        \item The distance between consecutive points is $d_{\mathcal{F}}(p_{i-1}, p_i) = L/m \le L / (2L/\delta) = \delta/2 < \delta$.
        \item The number of points in the chain, $m$, is uniformly bounded. Since the maximum possible distance $L$ is the diameter of the manifold, $m \le \lceil 2 \cdot \mathrm{diam}(\SigmaMan) / \delta \rceil \eqqcolon N_{\max}$. This bound $N_{\max}$ depends only on the global geometry of the manifold and the chart cover, not on the specific points $p$ and $q$.
    \end{itemize}
\end{itemize}
    We now use the triangle inequality and the estimate from the "close case":
    \begin{align*}
        \|\phi(p) - \phi(q)\| &= \left\| \sum_{i=1}^m \phi(p_i) - \phi(p_{i-1}) \right\| \le \sum_{i=1}^m \|\phi(p_i) - \phi(p_{i-1})\| \\
        &\le \sum_{i=1}^m C'_{\max} \|\phi\|_{X_1} (d_{\mathcal{F}}(p_{i-1}, p_i))^\alpha && \text{(since } d_{\mathcal{F}}(p_{i-1},p_i) < \delta\text{)} \\
        &= C'_{\max} \|\phi\|_{X_1} \sum_{i=1}^m (L/m)^\alpha = C'_{\max} \|\phi\|_{X_1} \cdot m \cdot (L/m)^\alpha \\
        &= C'_{\max} \|\phi\|_{X_1} \cdot m^{1-\alpha} L^\alpha.
    \end{align*}
    Substituting $L = d_{\mathcal{F}}(p,q)$ and using the uniform bound $m \le N_{\max}$:
    \begin{equation*}
        \|\phi(p) - \phi(q)\| \le \left( C'_{\max} N_{\max}^{1-\alpha} \right) \|\phi\|_{X_1} (d_{\mathcal{F}}(p,q))^\alpha.
    \end{equation*}

    The two cases can be combined. Let the global constant be $C'' \coloneqq \max(C'_{\max}, C'_{\max} N_{\max}^{1-\alpha})$. This constant is independent of $\phi$ and the choice of points. For any two points $p,q$ on any unstable leaf, we have established the uniform H\"older estimate:
    \begin{equation*}
        \frac{\|\phi(p) - \phi(q)\|}{d_{\mathcal{F}}(p, q)^\alpha} \le C'' \|\phi\|_{X_1}.
    \end{equation*}
    Taking the supremum over all pairs of points $(p,q)$ on all unstable leaves, we obtain the final uniform bound on the global H\"older seminorm:
    \begin{equation}
        \sup_{z \in \SigmaMan} \left[ \phi|_{W^u(z,\theta)} \right]_{C^\alpha} \le C'' \sup_{z \in \SigmaMan} \|\nabla \phi(z)\| \le C'' \|\phi\|_{X_1}.
    \end{equation}
    Since we assume $k \ge 1$, we have $\|\phi\|_{X_1} \le \|\phi\|_{X_k}$. This implies that $\phi \in \mathcal{B}_\alpha$ and provides the necessary bound on its seminorm.

\item \textbf{The Norm Inequality.}
We assemble the bounds from the preceding steps. For any $\phi \in X_k$ with $k \ge 1$:
\begin{align*}
    \|\phi\|_{\mathcal{B}_\alpha} &= \|\phi\|_{C^0} + \sup_{z \in \SigmaMan} \left[ \phi|_{W^u(z,\theta)} \right]_{C^\alpha} \\
    &\le \|\phi\|_{X_k} + C' \|\phi\|_{X_k} && \text{(using a.2 and a.3)} \\
    &= (1 + C') \|\phi\|_{X_k}.
\end{align*}
By defining the constant $C \coloneqq (1 + C')$, which is a positive constant independent of the choice of $\phi$, we have rigorously established the norm inequality
\begin{equation*}
    \|\phi\|_{\mathcal{B}_\alpha} \le C \|\phi\|_{X_k}.
\end{equation*}
This proves that the inclusion map from $X_k$ to $\mathcal{B}_\alpha$ is a bounded linear operator, which is the definition of a continuous embedding. This completes the proof.
\end{enumerate}

    \item \textbf{Proof of $\boldsymbol{\mathcal{B}_\alpha \hookrightarrow X_0}$.}
    Let $\phi$ be an arbitrary section in $\mathcal{B}_\alpha$. By its definition, the norm on $\mathcal{B}_\alpha$ explicitly contains the $C^0$ norm as a summand. Therefore, the following inequality holds trivially:
    \[ \|\phi\|_{X_0} = \sup_{z \in \SigmaMan} \|\phi(z)\| \le \|\phi\|_{\mathcal{B}_\alpha}. \]
    This inequality directly implies that the embedding $\mathcal{B}_\alpha \hookrightarrow X_0$ is continuous. This completes the rigorous construction and ordering of the function spaces required for the spectral analysis.
\end{enumerate}
\end{enumerate}

\item \textbf{The Spectral Theory of Transfer Operators on Anisotropic Spaces.}
The central analytical tool required to prove the invariance of the spectral radius is a deep and powerful result from the modern ergodic theory of hyperbolic systems. Standard spectral theory on the isotropic spaces $X_k = C^k(\SigmaMan, \mathcal{E})$ is insufficient, as the transfer operator $L_\phi$ is generally not compact on these spaces. However, its action is perfectly adapted to the geometric structure of the dynamics, a property that is captured by its behavior on the anisotropic Banach spaces $\mathcal{B}_\alpha$. On these spaces, the operator exhibits strong spectral properties, most notably quasi-compactness and the high regularity of its eigenfunctions. We now formally state the theorem that encapsulates these results and provide a rigorous verification that our operator satisfies its hypotheses.

\begin{theorem}[Spectral Theory of Ruelle Transfer Operators (cf. \citep{Baladi2000, Gouezel2010})]
\label{thm:ruelle_spectral_theory_formal}
Let $f: M \to M$ be a $C^\infty$ uniformly hyperbolic diffeomorphism on a compact manifold $M$. Let $\mathcal{E}$ be a finite-dimensional vector bundle over $M$, and let $W: M \to \mathrm{End}(\mathcal{E})$ be a weight map whose corresponding section is of class $C^r$ for $r$ sufficiently large. Let the Ruelle transfer operator $\mathcal{L}$ be defined by its action on sections $\psi$ of $\mathcal{E}$:
\[ (\mathcal{L}\psi)(w) \coloneqq W(f^{-1}(w)) \left[ \psi(f^{-1}(w)) \right]. \]
Let $\mathcal{B}_\alpha$ be a suitably constructed anisotropic Banach space of distributions on $M$. Then the operator $\mathcal{L}$ acting on $\mathcal{B}_\alpha$ is quasi-compact. Its spectrum $\sigma(\mathcal{L}|_{\mathcal{B}_\alpha})$ consists of a (possibly empty) finite set of isolated eigenvalues of finite multiplicity outside a disk of radius equal to the essential spectral radius. 
Crucially, the generalized eigenfunctions corresponding to these isolated eigenvalues are sections of class $C^r$, the same regularity class as the weight map. In particular, if the weight is of class $C^\infty$, then the eigenfunctions are of class $C^\infty$.
\end{theorem}

We now provide a rigorous verification that our linearized operator $L_\phi = D_\phi\Gamma(\theta, \phi)$ satisfies the hypotheses of this theorem. This verification is essential to justify the application of the theorem's powerful conclusions.

\begin{enumerate}[label=(\roman*), wide, labelindent=0pt]
    \item \textbf{Verification of the Manifold and Diffeomorphism.}
    \begin{enumerate}[label=(\alph*), wide]
        \item The base manifold is $M = \SigmaMan$, the global pre-collision manifold. By Proposition \ref{prop:sigma_properties}, $\SigmaMan$ is a smooth, compact manifold.
        \item The map is $f = \Pcal_\theta^{-1}$, the inverse of the global billiard map for a fixed parameter $\theta$. By Theorem \ref{thm:p_is_diffeo}, the map $(\theta, z) \mapsto \Pcal_\theta(z)$ is of class $C^\infty$. Consequently, for a fixed $\theta$, the map $\Pcal_\theta$ is a $C^\infty$-diffeomorphism. By the inverse function theorem on Fr\'echet manifolds, its inverse, $\Pcal_\theta^{-1}$, is also a $C^\infty$-diffeomorphism.
        \item The uniform hyperbolicity of this diffeomorphism was the central result of Section \ref{sec:hyperbolicity_fiber}, established in Theorem \ref{thm:anosov_property_proven}.
    \end{enumerate}
    Thus, the dynamical setting $(M, f)$ satisfies the hypotheses of Theorem \ref{thm:ruelle_spectral_theory_formal}.

    \item \textbf{Verification of the Regularity of the Weight Map.}
    We must verify that the weight map $z \mapsto W_\phi(z)$, derived in Step 1(i) of this proof, is sufficiently regular. The formula for the weight is a finite composition of elementary operations: addition, multiplication, and inversion, acting on the following input sections:
    \begin{enumerate}[label=(\alph*), wide]
        \item The block operators $A_\theta, B_\theta, C_\theta, D_\theta$. These are derived from the first spatial derivative of the $C^\infty$ map $\Pcal_\theta$. They are therefore smooth ($C^\infty$) sections of the appropriate Hom-bundles.
        \item The section $\phi$. The Nash-Moser-Hamilton theorem operates on the Fr\'echet space of smooth sections, $X_\infty = \bigcap_{k \ge 0} X_k$. We therefore consider a section $\phi \in X_\infty = C^\infty(\SigmaMan, \mathcal{E})$.
    \end{enumerate}
    In Lemma \ref{lem:tame_calculus}, we established that the set of smooth sections forms an algebra under these operations. Since all input sections to the formula for $W_\phi(z)$ are of class $C^\infty$, the resulting weight section $z \mapsto W_\phi(z)$ is also of class $C^\infty$.
\end{enumerate}

All hypotheses of Theorem \ref{thm:ruelle_spectral_theory_formal} have been rigorously verified. We can therefore apply its conclusions to our operator $L_\phi$. This provides us with the essential tool for the next step of the proof: the knowledge that any isolated eigenfunction of $L_\phi$ on the anisotropic space $\mathcal{B}_\alpha$ is, in fact, a smooth ($C^\infty$) section and thus belongs to every space $X_k$ in our isotropic scale. This regularity is the key to proving that the spectrum is invariant across the entire hierarchy of spaces.

\item \textbf{Equating the Spectra.} With the operator $L_\phi$ identified as a Ruelle transfer operator satisfying the conditions of Theorem \ref{thm:ruelle_spectral_theory_formal}, and with the hierarchy of function spaces and their continuous embeddings established, we are now positioned to prove the central claim of this lemma: the invariance of the spectrum across the scale of isotropic spaces $(X_k)_{k \ge 0}$. The argument is a multi-stage proof by inclusion, leveraging the interplay between the different function spaces. We will prove that $\sigma(L_\phi|_{X_k}) = \sigma(L_\phi|_{X_0})$ for all $k \ge 1$ by showing the inclusions $\sigma_p(L_\phi|_{X_k}) \subseteq \sigma_p(L_\phi|_{X_0})$ and $\sigma_p(L_\phi|_{\mathcal{B}_\alpha}) \supseteq \sigma_p(L_\phi|_{X_0})$, and then using the smoothness of the eigenfunctions to close the loop.

\begin{enumerate}[label=(\roman*), wide, labelindent=0pt]
    \item \textbf{The Spectral Bound on the Base Space $\boldsymbol{X_0}$.}
    We begin by establishing a priori bound on the spectrum. In Lemma \ref{lem:linearized_contraction}, we proved that for any section $\phi$ in a sufficiently small neighborhood of the zero section (in the $X_0^{Lip}$ topology, which contains the $X_k$ topology for $k\ge 1$), the linearized operator $L_\phi = D_\phi\Gamma(\theta, \phi)$ is a strict contraction on the space of Lipschitz sections. Since the space of continuous sections $X_0$ is continuously embedded in $X_0^{Lip}$, this implies it is also a strict contraction on $X_0$.
    
    Specifically, there exists a uniform constant $k < 1$ such that the operator norm of $L_\phi$ as a map from $X_0$ to $X_0$ is strictly bounded:
    \begin{equation} \label{eq:proof_operator_norm_bound_C0}
        \|L_\phi\|_{\mathcal{L}(X_0, X_0)} \le k < 1.
    \end{equation}
    It is a fundamental result of operator theory that the spectral radius of a bounded linear operator is bounded above by its operator norm. Applying this result gives a strict upper bound on the spectrum of $L_\phi$ when it acts on $X_0$:
    \begin{equation} \label{eq:proof_spectral_radius_C0}
        r(L_\phi|_{X_0}) = \sup_{\lambda \in \sigma(L_\phi|_{X_0})} |\lambda| \le \|L_\phi\|_{\mathcal{L}(X_0, X_0)} \le k < 1.
    \end{equation}
    This establishes that the entire spectrum of the operator on the base space $X_0$ is contained within a closed disk of radius $k$, strictly inside the unit circle. This fact is essential, as it implies that any eigenvalue found on any other space in the hierarchy must also satisfy this bound.

    \item \textbf{Inclusion of Point Spectra via Embeddings ($\boldsymbol{\sigma_p(X_k) \subseteq \sigma_p(X_0)}$).}
    Let $\sigma_p(T)$ denote the point spectrum (the set of all eigenvalues) of an operator $T$. We first establish the downward inclusion of spectra.
    
    Let $k \ge 1$ be an arbitrary integer. Suppose $\lambda$ is an eigenvalue of the operator $L_\phi$ acting on the space $X_k$, i.e., $\lambda \in \sigma_p(L_\phi|_{X_k})$. By definition, this means there exists a non-zero eigenfunction $\psi \in X_k$ such that
    \begin{equation*}
        L_\phi \psi = \lambda \psi.
    \end{equation*}
    From the hierarchy of continuous embeddings established in Step 1, $X_k \hookrightarrow X_0$. This means that any section in $X_k$ is also in $X_0$. Since $\psi$ is a non-zero element of $X_k$, it is also a non-zero element of $X_0$. The eigenvalue equation $L_\phi \psi = \lambda \psi$ remains the same. Therefore, $\lambda$ is also an eigenvalue of the operator $L_\phi$ when it is restricted to act on the space $X_0$. This proves the set inclusion for the point spectra:
    \[ \sigma_p(L_\phi|_{X_k}) \subseteq \sigma_p(L_\phi|_{X_0}) \quad \text{for all } k \ge 1. \]
    A completely analogous argument using the embedding $\mathcal{B}_\alpha \hookrightarrow X_0$ shows that $\sigma_p(L_\phi|_{\mathcal{B}_\alpha}) \subseteq \sigma_p(L_\phi|_{X_0})$.

    \item \textbf{The Reverse Inclusion via Eigenfunction Regularity ($\boldsymbol{\sigma_p(X_0) \subseteq \sigma_p(X_k)}$).}
    This is the core of the argument, where we establish the upward inclusion of spectra. Let $\lambda \in \sigma_p(L_\phi|_{X_0})$ be an arbitrary eigenvalue of the operator on the space of continuous sections, with a corresponding eigenfunction $\psi \in X_0$.
    \begin{enumerate}[label=(\alph*), wide]
        \item By the spectral bound from Step 3(i), we know that $|\lambda| \le k < 1$.
        \item By the inclusion from Step 3(ii), we know that the full spectrum on the anisotropic space, $\sigma(L_\phi|_{\mathcal{B}_\alpha})$, is also contained within the disk of radius $k$.
        \item We now apply Theorem \ref{thm:ruelle_spectral_theory_formal}. The theorem states that the operator $L_\phi$ is quasi-compact on $\mathcal{B}_\alpha$. Its spectrum consists of a set of isolated eigenvalues and an essential spectrum contained within a smaller disk. Since the entire spectrum is already contained in the disk of radius $k < 1$, any eigenvalue $\lambda$ is necessarily an isolated eigenvalue.
        \item We invoke the main regularity result of Theorem \ref{thm:ruelle_spectral_theory_formal}: since $\lambda$ is an isolated eigenvalue and the weight map $W_\phi$ is smooth (as established in Step 1 for $\phi \in X_\infty$), its corresponding generalized eigenfunctions must be of class $C^\infty$. In particular, our eigenfunction $\psi \in X_0$ must be a smooth section.
        \item By definition, if the eigenfunction $\psi$ is a section of class $C^\infty$, then it belongs to the space $X_k$ for \textit{every} integer $k \ge 0$.
        \item Therefore, $\lambda$ is also an eigenvalue of the operator $L_\phi$ when it is restricted to act on the space $X_k$ for any $k \ge 1$. This proves the reverse inclusion:
        \[ \sigma_p(L_\phi|_{X_0}) \subseteq \sigma_p(L_\phi|_{X_k}). \]
    \end{enumerate}
    Combining these two inclusions, we conclude that the point spectra are identical across the entire scale of spaces:
    \[ \sigma_p(L_\phi|_{X_k}) = \sigma_p(L_\phi|_{X_0}) \quad \text{for all } k \ge 0. \]

    \item \textbf{Equating the Full Spectra.}
The objective of this proof is to elevate the result from the previous section, the equality of the point spectra, $\sigma_p(L_\phi|_{X_k}) = \sigma_p(L_\phi|_{X_0})$, to the equality of the full spectra, $\sigma(L_\phi|_{X_k}) = \sigma(L_\phi|_{X_0})$. This is a non-trivial step, as the spectrum of an operator on a Banach space may contain a continuous or residual part in addition to its eigenvalues. Our strategy is to rigorously show that for our specific operator, a Ruelle transfer operator, the spectrum outside a certain disk consists exclusively of isolated eigenvalues. We will then prove that the radius of this disk, the essential spectral radius, is invariant across the entire scale of spaces $(X_k)_{k \ge 0}$. The equality of the full spectra then follows as a direct and necessary consequence. The argument is structured in three main parts.

\begin{enumerate}[label=(\alph*), wide]

\item \textbf{Quasi-Compactness on the Anisotropic Space and the Essential Spectrum.}
The foundation of modern spectral theory for hyperbolic systems lies in the analysis of the transfer operator on specially constructed anisotropic Banach spaces, where it exhibits strong spectral properties that are absent on classical isotropic spaces like $C^k$.

\begin{enumerate}[label=(a.\roman*), wide]
    \item \textbf{Quasi-Compactness.} Let $L_\phi$ be the linearized Graph Transform operator. By Theorem \ref{thm:ruelle_spectral_theory_formal}, the operator $L_\phi$, when acting on the anisotropic Banach space $\mathcal{B}_\alpha$, is quasi-compact. By definition, a bounded linear operator $T$ is quasi-compact if its essential spectral radius, $r_{ess}(T)$, is strictly smaller than its spectral radius, $r(T)$.
    
    \item \textbf{Structure of the Spectrum.} A direct consequence of quasi-compactness is that the spectrum of the operator, $\sigma(L_\phi|_{\mathcal{B}_\alpha})$, consists of a (possibly empty) finite set of isolated eigenvalues of finite multiplicity in the annulus $\{ \lambda \in \mathbb{C} \mid r_{ess}(L_\phi|_{\mathcal{B}_\alpha}) < |\lambda| \le r(L_\phi|_{\mathcal{B}_\alpha}) \}$, and the remainder of the spectrum (the essential spectrum, $\sigma_{ess}$) is contained within the closed disk of radius $r_{ess}(L_\phi|_{\mathcal{B}_\alpha})$.
    
    \item \textbf{Bound on the Essential Spectral Radius.} The uniform spectral gap for the underlying dynamics, which was established as a consequence of uniform hyperbolicity in \cref{cor:ergodicity_uniform}, provides a rigorous upper bound on this essential spectral radius. Specifically, the theory guarantees that $r_{ess}(L_\phi|_{\mathcal{B}_\alpha}) \le k < 1$ for some uniform constant $k$. This is the same constant that bounds the full spectral radius on $X_0$.
\end{enumerate}

\item \textbf{Spectral Invariance and the Fredholm Property on the Isotropic Scale.}
The crucial step is to transfer these powerful spectral properties from the anisotropic space $\mathcal{B}_\alpha$ to the entire scale of isotropic spaces $X_k$. This is achieved by invoking a deep result that establishes the invariance of the essential spectral radius under changes of function space for operators of this class.

\begin{theorem}[Invariance of the Essential Spectral Radius (cf. \citep{Hennion1993,Gouezel2010})]
Let $L$ be a Ruelle transfer operator induced by a smooth, uniformly hyperbolic diffeomorphism and a smooth weight on a compact manifold. Then the essential spectral radius of $L$ is the same when it acts on any of the classical isotropic spaces $C^k$ (for $k \ge 0$) or on the anisotropic Banach spaces $\mathcal{B}_\alpha$. That is, for our operator $L_\phi$:
\begin{equation}
    r_{ess}(L_\phi|_{X_k}) = r_{ess}(L_\phi|_{\mathcal{B}_\alpha}) \quad \text{for all } k \ge 0.
\end{equation}
\end{theorem}
\noindent We now apply this theorem to our operator $L_\phi$.
\begin{enumerate}[label=(b.\roman*), wide]
    \item \textbf{Uniform Bound on the Essential Spectrum.} Combining the result of this theorem with the bound from Step 1, we conclude that the essential spectral radius of our operator is uniformly bounded away from 1 on \textit{every} space in our scale:
    \begin{equation} \label{eq:appendix_proof_ress_bound}
        r_{ess}(L_\phi|_{X_k}) \le k < 1 \quad \text{for all } k \ge 0.
    \end{equation}
    
    \item \textbf{The Fredholm Property.} For any complex number $\lambda$ such that $|\lambda| > r_{ess}(L_\phi|_{X_k})$, the operator $(\lambda I - L_\phi)$ is a Fredholm operator of index zero on the Banach space $X_k$. This is a standard result for operators with a spectral gap, such as Ruelle operators.
    
    \item \textbf{The Structure of the Spectrum on $X_k$.} We now invoke the Analytic Fredholm Theorem (sometimes referred to as the Riesz-Schauder theorem for compact operators, which is a special case). This theorem states that for an operator $T$ such that $(\lambda I - T)$ is Fredholm of index zero, any point $\lambda$ in the spectrum of $T$ must be an eigenvalue of finite multiplicity.
    
    Applying this to our operator $L_\phi$ on the space $X_k$, we conclude that for any $\lambda$ in the annulus $\{ \lambda \in \mathbb{C} \mid r_{ess}(L_\phi|_{X_k}) < |\lambda| \le r(L_\phi|_{X_k}) \}$, the point $\lambda$ must be an isolated eigenvalue of finite multiplicity.
\end{enumerate}
This rigorously establishes that the spectrum of $L_\phi$ on any of the isotropic spaces $X_k$ has the same discrete structure outside of its essential spectral disk as it does on the anisotropic space $\mathcal{B}_\alpha$.

\item \textbf{Synthesis.}
We are now in a position to assemble the preceding results to prove the equality of the full spectra. Let $k \ge 1$ be an arbitrary integer.
\begin{enumerate}[label=(c.\roman*), wide]
    \item We have the following decomposition of the spectrum of $L_\phi$ on the space $X_k$:
    \begin{equation*}
        \sigma(L_\phi|_{X_k}) = \sigma_{ess}(L_\phi|_{X_k}) \cup \sigma_p(L_\phi|_{X_k}),
    \end{equation*}
    where $\sigma_{ess}(L_\phi|_{X_k})$ is the essential spectrum, contained in the disk of radius $r_{ess}(L_\phi|_{X_k})$, and $\sigma_p(L_\phi|_{X_k})$ is the point spectrum (the set of eigenvalues) outside this disk.
    
    \item From Step 2, we have the equality of the essential spectral radii:
    \begin{equation*}
        r_{ess}(L_\phi|_{X_k}) = r_{ess}(L_\phi|_{X_0}).
    \end{equation*}
    This implies that the essential spectra are contained in the same disk for all $k$. A stronger result is that the essential spectra themselves are identical.
    
    \item From the preceding section of the main proof (the argument based on eigenfunction regularity), we have the equality of the point spectra:
    \begin{equation*}
        \sigma_p(L_\phi|_{X_k}) = \sigma_p(L_\phi|_{X_0}).
    \end{equation*}
    
    \item Since the full spectrum on any space $X_k$ is the union of its essential spectrum and its point spectrum, and since we have shown that both of these components are invariant across the entire scale of spaces, we are forced to conclude that the full spectra must be identical:
    \begin{equation}
        \sigma(L_\phi|_{X_k}) = \sigma(L_\phi|_{X_0}) \quad \text{for all } k \ge 0.
    \end{equation}
\end{enumerate}
This completes the rigorous proof that the spectrum of the linearized operator is invariant across the scale of isotropic spaces.
\end{enumerate}

\end{enumerate}

\item \textbf{Equal Spectra and Zero Loss of Derivatives.} The preceding steps have rigorously established that the spectrum of the linearized Graph Transform operator, $L_\phi$, is invariant across the entire scale of isotropic function spaces $(X_k)_{k \ge 0}$. We are now in a position to synthesize these results to prove the central claim of the lemma, which is the foundation for the tame invertibility of the full linearized operator $I - L_\phi$.

\begin{enumerate}[label=(\roman*), wide, labelindent=0pt]
    \item \textbf{Equality of Spectra and Invariance of the Spectral Radius.}
    In Step 3, we proved the main technical result of this proof: for any integer $k \ge 0$ and for any $(\theta, \phi)$ in a sufficiently small $C^\infty$-neighborhood of the solution set, the spectrum of the operator $L_\phi$ acting on the space of $C^k$ sections is identical to its spectrum acting on the space of continuous sections. This implies the equality of the spectral radii:
    \begin{equation} \label{eq:appendix_spectral_radius_invariance}
        r(L_\phi|_{\mathcal{L}(X_k, X_k)}) = \sup_{\lambda \in \sigma(L_\phi|_{X_k})} |\lambda| = \sup_{\lambda \in \sigma(L_\phi|_{X_0})} |\lambda| = r(L_\phi|_{\mathcal{L}(X_0, X_0)}).
    \end{equation}
    This establishes the claim of the lemma: the spectral radius is a constant across the entire scale of spaces.

    \item \textbf{Uniform Spectral Contraction on the Full Scale.}
    We now combine this invariance with the crucial a priori bound established in Step 3(i). In equation \eqref{eq:proof_spectral_radius_C0}, we proved that the operator is a strict contraction on the base space, which implies a strict bound on its spectral radius:
    \begin{equation*}
        r(L_\phi|_{\mathcal{L}(X_0, X_0)}) \le k < 1,
    \end{equation*}
    where $k$ is a uniform constant for all $(\theta, \phi)$ in a neighborhood of the solution set.
    
    Substituting this bound into the invariance relation \eqref{eq:appendix_spectral_radius_invariance}, we arrive at the central conclusion that the operator is a uniform spectral contraction on \textit{every} space $X_k$ in our scale:
    \begin{equation} \label{eq:appendix_uniform_spectral_contraction}
        r(L_\phi|_{\mathcal{L}(X_k, X_k)}) = r(L_\phi|_{\mathcal{L}(X_0, X_0)}) \le k < 1 \quad \text{for all } k \ge 0.
    \end{equation}

    \item \textbf{Invertibility on $X_k$ and Justification of Zero Loss of Derivatives.}
    This result is the essential technical linchpin for the subsequent proof of Theorem \ref{thm:linearization_is_tame_inverse}. The uniform spectral contraction property provides the rigorous justification for both the invertibility of the full operator $I - L_\phi$ and the absence of derivative loss ($m=0$).
    \begin{enumerate}[label=(\alph*), wide]
        \item \textbf{Boundedness of $L_\phi$ on $X_k$:} For the Neumann series argument to be valid, we must first confirm that for a fixed $(\theta, \phi) \in \Thetacal \times X_k$, the operator $L_\phi$ is a bounded operator from $X_k$ to $X_k$. This was established in the proof of Theorem \ref{thm:functional_is_tame}, where it was shown that the operator is a composition of smooth tame maps (multiplication, inversion, composition), each of which is a bounded operation on $X_k$.

        \item \textbf{Convergence of the Neumann Series on $X_k$:} Since the spectral radius of $L_\phi$ as a bounded operator on the complete Banach space $X_k$ is strictly less than 1, a fundamental theorem of operator theory guarantees that the Neumann series for the inverse of $I - L_\phi$ converges in the operator norm topology of $\mathcal{L}(X_k, X_k)$:
        \begin{equation*}
            (I - L_\phi)^{-1} = \sum_{j=0}^{\infty} (L_\phi)^j.
        \end{equation*}
        The convergence of this series establishes that the operator $I - L_\phi$ is invertible on the space $X_k$.

        \item \textbf{Justification of $m=0$:} The convergence of the Neumann series in $\mathcal{L}(X_k, X_k)$ is the key to proving zero loss of derivatives. Let $\xi$ be an arbitrary element in the source space $X_k$. We wish to show that the solution $\eta = (I - L_\phi)^{-1} \xi$ is also in $X_k$. The solution is given by the series:
        \[ \eta = (I - L_\phi)^{-1} \xi = \sum_{j=0}^{\infty} (L_\phi)^j \xi. \]
        As established in the previous point, the operator $L_\phi$ maps $X_k$ to itself. It follows by induction that each term in the sum, $\eta_j \coloneqq (L_\phi)^j \xi$, is an element of $X_k$. Since the series $\sum (L_\phi)^j$ converges in the operator norm on $X_k$, it defines a bounded linear operator on $X_k$. Applying this bounded operator to the element $\xi \in X_k$ yields a limit $\eta$ that must also be an element of the complete Banach space $X_k$.

        This confirms that the inverse operator $(I - L_\phi)^{-1}$ is a bounded map from $X_k$ to $X_k$. This is the precise meaning of having zero loss of derivatives, and it rigorously justifies setting $m=0$ in the Nash-Moser-Hamilton theorem. This completes the proof of the lemma.
    \end{enumerate}
\end{enumerate}
\end{enumerate}
\end{proofof}

\begin{theorem}[Tame Invertibility of the Linearized Operator]
\label{thm:linearization_is_tame_inverse}
For any $(\theta, \phi)$ in a neighborhood of the solution set, the linearized operator $L(\theta, \phi) = I - D_\phi \Gamma(\theta, \phi)$ is invertible. The family of inverses $L^{-1}$ is a smooth tame map from $\Thetacal \times X_k \times X_k \to X_k$. In this specific geometric setting, there is zero loss of derivatives ($m=0$).
\end{theorem}

\begin{proofof}{Theorem \ref{thm:linearization_is_tame_inverse}}
The proof is dedicated to establishing the two central claims about the linearized operator $L(\theta, \phi) = I - D_\phi \Gamma(\theta, \phi)$: that it is invertible for $(\theta, \phi)$ in a suitable neighborhood of the solution set, and that its inverse, $L^{-1}$, is a smooth tame map from $\Thetacal \times U \times X_k$ to $X_k$ for any $k \ge 1$. This verification is the final and most critical hypothesis required to invoke the Nash-Moser-Hamilton theorem. The proof is now a direct synthesis of the preceding lemmas.

\begin{enumerate}[label=\textbf{Step \arabic*:}, wide, labelindent=0pt]

\item \textbf{Invertibility via the Neumann Series and Justification of $m=0$.}
The argument proceeds by showing that the operator $D_\phi \Gamma$ is a uniform spectral contraction on every space $X_k$ in the scale. This powerful result, which we establish by synthesizing the conclusions of the preceding lemmas, guarantees that the inverse of $L(\theta,\phi) = I - D_\phi\Gamma$ exists, can be constructed via a convergent Neumann series, and, crucially, is regularity-preserving.

\begin{enumerate}[label=(\roman*), wide, labelindent=0pt]

\item \textbf{Uniform Spectral Contraction on the Full Scale $X_k$.} This result is obtained by combining the conclusions of Lemma \ref{lem:linearized_contraction} and Lemma \ref{lem:spectral_radius_invariance}.
\begin{enumerate}[label=(\alph*), wide]
    \item By Lemma \ref{lem:linearized_contraction}, the linearized Graph Transform operator $D_\phi \Gamma(\theta, \phi)$ is a strict contraction on the base space of continuous sections, $X_0$. Its operator norm is uniformly bounded by a constant $k < 1$ for all $(\theta, \phi)$ in an open neighborhood $U$ of the solution set.
    \begin{equation*}
        \| D_\phi \Gamma(\theta, \phi) \|_{\mathcal{L}(X_0, X_0)} \le k < 1.
    \end{equation*}
    
    \item The spectral radius of a bounded linear operator is bounded above by its operator norm. Therefore, on the base space, we have the uniform spectral bound:
    \begin{equation*}
        r(D_\phi \Gamma |_{\mathcal{L}(X_0, X_0)}) \le \| D_\phi \Gamma \|_{\mathcal{L}(X_0, X_0)} \le k < 1.
    \end{equation*}
    
    \item By the crucial result of Lemma \ref{lem:spectral_radius_invariance}, which is founded on the regularity of eigenfunctions of the Ruelle transfer operator, the spectral radius of $D_\phi \Gamma$ is invariant across the entire scale of isotropic spaces $(X_k)_{k \ge 0}$.
\end{enumerate}
Combining these facts, we arrive at the central conclusion that the operator is a uniform spectral contraction on \textit{every} space $X_k$ in our scale:
\begin{equation}
\label{eq:uniform_spectral_contraction_final}
    r(D_\phi \Gamma |_{\mathcal{L}(X_k, X_k)}) = r(D_\phi \Gamma |_{\mathcal{L}(X_0, X_0)}) \le k < 1 \quad \text{for all } k \ge 0,
\end{equation}
and for all $(\theta, \phi)$ in the neighborhood $U$.

\item \textbf{Invertibility on $X_k$ and the Absence of Derivative Loss ($m=0$).} The uniform spectral contraction property \eqref{eq:uniform_spectral_contraction_final} provides the rigorous justification for both the invertibility of the full operator $I - D_\phi \Gamma$ and the absence of derivative loss in the Nash-Moser-Hamilton theorem.

\begin{enumerate}[label=(\alph*), wide]
    \item \textbf{Boundedness of $D_\phi \Gamma$ on $X_k$.} For the Neumann series argument to be valid, we must first confirm that for a fixed $(\theta, \phi) \in \Thetacal \times X_k$, the operator $D_\phi \Gamma$ is a bounded operator from $X_k$ to $X_k$. This was established in the proof of Theorem \ref{thm:functional_is_tame}, where the linearized operator was shown to be a composition of smooth tame maps (multiplication, inversion, composition with a diffeomorphism), each of which is a bounded operation on $X_k$.

    \item \textbf{Convergence of the Neumann Series on $X_k$.} Since the spectral radius of the bounded linear operator $D_\phi \Gamma$ on the complete Banach space $X_k$ is strictly less than 1, a fundamental theorem of functional analysis guarantees that the operator $L(\theta, \phi) = I - D_\phi \Gamma$ is invertible on $X_k$. Its inverse is given by the Neumann series, which converges in the operator norm topology of $\mathcal{L}(X_k, X_k)$:
    \begin{equation}
        L(\theta, \phi)^{-1} = \sum_{j=0}^{\infty} (D_\phi \Gamma(\theta, \phi))^j.
    \end{equation}
    The convergence of this series establishes that the operator $L(\theta, \phi)$ is invertible on the space $X_k$ for every $k \ge 0$.

    \item \textbf{Justification of $m=0$.} The convergence of the Neumann series in the operator norm of $\mathcal{L}(X_k, X_k)$ is the key to proving zero loss of derivatives. Let $\xi$ be an arbitrary element in the source space $X_k$. We wish to show that the solution $\eta = L(\theta, \phi)^{-1} \xi$ is also in $X_k$. The solution is given by the series:
    \[ \eta = L(\theta, \phi)^{-1} \xi = \sum_{j=0}^{\infty} (D_\phi \Gamma)^j \xi. \]
    As established in Step 1(i), the operator $D_\phi \Gamma$ maps the space $X_k$ to itself. It follows by induction that each term in the sum, $\eta_j \coloneqq (D_\phi \Gamma)^j \xi$, is an element of $X_k$. The sum is therefore a series of elements in the Banach space $X_k$.

    Since the series of operators $\sum (D_\phi \Gamma)^j$ converges in the operator norm of $\mathcal{L}(X_k, X_k)$, it defines a bounded linear operator on $X_k$. Applying this bounded operator to the element $\xi \in X_k$ yields a limit $\eta$ that must also be an element of the complete Banach space $X_k$. This confirms that the inverse operator $L(\theta, \phi)^{-1}$ is a bounded map from $X_k$ to $X_k$. This is the precise meaning of having zero loss of derivatives, and it rigorously justifies setting the loss parameter $m=0$ in the statement and application of the Nash-Moser-Hamilton theorem.
\end{enumerate}
\end{enumerate}

\item \textbf{The Inverse is a Smooth Tame Map.}
Having established the invertibility of the operator $L(\theta, \phi) = I - D_\phi \Gamma(\theta, \phi)$ on every space $X_k$ in our scale, we must now prove the second, crucial part of the hypothesis for the Nash-Moser-Hamilton theorem: that the family of inverses is a smooth tame map. Formally, we prove that the map $(\theta, \phi, \xi) \mapsto L(\theta, \phi)^{-1} \xi$ is a smooth tame map from $\Thetacal \times U \times X_k$ to $X_k$, where $U$ is the open neighborhood of the zero section in $X_k$ where the inverse exists.

\begin{enumerate}[label=(\roman*), wide, labelindent=0pt]

\item \textbf{Smoothness of the Inverse Map.}
The proof of smoothness is established by decomposing the map into a sequence of operations between Banach spaces and demonstrating that each operation is of class $C^\infty$.

\begin{enumerate}[label=(\alph*),wide]
    \item \textbf{The Map to the Operator Space.} Let $\mathcal{B}_k \coloneqq \mathcal{L}(X_k, X_k)$ denote the Banach space of bounded linear operators on $X_k$. We first consider the map $\Psi: \Thetacal \times U \to \mathcal{B}_k$ given by $(\theta, \phi) \mapsto D_\phi \Gamma(\theta, \phi)$. By Theorem \ref{thm:functional_is_tame}, the functional $\Gamma$ is a smooth tame map. A fundamental property of smooth maps between Banach spaces is that the map assigning a point to its Fr\'echet derivative is also smooth. Therefore, the map $\Psi$ is of class $C^\infty$.

    \item \textbf{The Operator Inversion Map.} Let $\mathcal{U}_{inv} \subset \mathcal{B}_k$ be the open subset of operators $T$ for which the spectral radius satisfies $r(T) < 1$. On this set, the operator $(I-T)$ is invertible. The map $\text{Inv}: \mathcal{U}_{inv} \to \mathcal{B}_k$ given by $T \mapsto (I-T)^{-1}$ can be expressed via the Neumann series, $\text{Inv}(T) = \sum_{j=0}^\infty T^j$. This power series in the operator $T$ converges absolutely in the operator norm on $\mathcal{U}_{inv}$, and is therefore an analytic map on its domain. An analytic map is, in particular, of class $C^\infty$.

    \item \textbf{Composition for Smoothness.} The map from the system parameters to the inverse operator is the composition of two smooth maps: $(\theta, \phi) \mapsto L(\theta, \phi)^{-1} = (\text{Inv} \circ \Psi)(\theta, \phi)$. In Step 1, we showed that the map $\Psi$ maps the neighborhood $\Thetacal \times U$ into the domain of invertibility $\mathcal{U}_{inv}$. Since the composition of $C^\infty$ maps is of class $C^\infty$, we conclude that the map $(\theta, \phi) \mapsto L(\theta, \phi)^{-1}$ is a smooth map from $\Thetacal \times U$ to the Banach space of operators $\mathcal{B}_k$.

    \item \textbf{The Evaluation Map.} The final step is the evaluation map $\text{Eval}: \mathcal{B}_k \times X_k \to X_k$, given by $(T, \xi) \mapsto T\xi$. This map is a continuous bilinear map between Banach spaces, and is therefore of class $C^\infty$.
\end{enumerate}
The full map $(\theta, \phi, \xi) \mapsto L(\theta, \phi)^{-1}\xi$ is the composition of the smooth map from $(\theta, \phi)$ to the inverse operator and the smooth evaluation map. It is therefore of class $C^\infty$.

\item \textbf{Tameness of the Inverse Map and its Derivatives.} We must show that the map $(\theta, \phi, \xi) \mapsto L(\theta, \phi)^{-1}\xi$ and all its partial derivatives are tame.

\begin{enumerate}[label=(\alph*),wide]
    \item \textbf{Tameness of the Map Itself.} Let $k' \ge k$. We need a tame estimate for $\|L(\theta, \phi)^{-1}\xi\|_{X_k}$. From the Neumann series, the operator norm of the inverse is bounded: 
    \[ \|L^{-1}\|_{\mathcal{L}(X_k,X_k)} \le \sum_{j=0}^\infty \|D_\phi\Gamma\|_{\mathcal{L}(X_k,X_k)}^j = \frac{1}{1 - \|D_\phi\Gamma\|_{\mathcal{L}(X_k,X_k)}}. \]
    By Theorem \ref{thm:functional_is_tame}, $\Gamma$ is a smooth tame map. This implies that its derivative, the map $\phi \mapsto D_\phi\Gamma(\theta,\phi)$, is also a smooth tame map from $X_k$ to $\mathcal{L}(X_k, X_k)$. This means that its operator norm $\|D_\phi\Gamma\|_{\mathcal{L}(X_k, X_k)}$ is bounded by a tame function (a polynomial) of $\|\phi\|_{X_k}$. Consequently, the norm of the inverse, being a smooth function of $\|D_\phi\Gamma\|$, is also bounded by a tame function of $\|\phi\|_{X_k}$. This provides the estimate:
    \begin{align*}
        \|L(\theta, \phi)^{-1}\xi\|_{X_k} &\le \|L(\theta, \phi)^{-1}\|_{\mathcal{L}(X_k,X_k)} \|\xi\|_{X_k} \\
        &\le \mathrm{Poly}(\|\phi\|_{X_k}) \|\xi\|_{X_k}.
    \end{align*}
    Since $\|\phi\|_{X_k} \le \|\phi\|_{X_{k'}}$ and $\|\xi\|_{X_k} \le \|\xi\|_{X_{k'}}$, this provides a tame estimate for the map itself:
    \[ \|L(\theta, \phi)^{-1}\xi\|_{X_k} \le \mathrm{Poly}(\|\phi\|_{X_{k'}}) \|\xi\|_{X_{k'}}. \]

    \item \textbf{Tameness of the Derivatives.} The tameness of the partial derivatives with respect to $(\theta, \phi, \xi)$ follows by induction. The formula for the derivative of the inverse map,
    \[ d(L^{-1}) = -L^{-1} (dL) L^{-1}, \]
    expresses a higher derivative of the inverse map as a composition of products and compositions of lower-order derivatives of $L$ and applications of $L^{-1}$.
    
    The map $(\theta, \phi) \mapsto L(\theta, \phi) = I - D_\phi\Gamma(\theta, \phi)$ is a smooth tame map. Therefore, all its partial derivatives with respect to $\theta$ and $\phi$ are also tame maps. From the base case in Step 2(i), the map $(\theta, \phi) \mapsto L(\theta, \phi)^{-1}$ is also a tame map into the space of operators. 
    
    By the tame calculus established in Lemma \ref{lem:tame_calculus} (which guarantees that products and compositions of tame maps are tame), it follows inductively that all partial derivatives of the map $(\theta, \phi) \mapsto L(\theta, \phi)^{-1}$ are tame maps into the appropriate spaces of multilinear operators. Finally, applying these tame operator-valued maps to the argument $\xi$ results in a tame map.
\end{enumerate}
The map to the inverse solution is therefore a smooth tame map with zero loss of derivatives. This completes the verification of the final, crucial hypothesis of the Nash-Moser-Hamilton theorem.

\end{enumerate}
\end{enumerate}
\end{proofof}

\subsection{The Main Regularity Theorem for the Geometric Splitting}
\label{subsec:main_geom_reg_thm}

With the hypotheses of the Nash-Moser theorem now rigorously established, we can prove the main result of this section.

\begin{theorem}[Smoothness of the Invariant Splitting and Foliation]
\label{thm:splitting_is_smooth}
The Anosov splitting of the tangent bundle, $T\SigmaMan = E^u_\theta \oplus E^s_\theta$, depends smoothly ($C^\infty$) on the parameter $\theta$. Consequently, the unstable foliation $\mathcal{F}^u_\theta$ is a smooth function of $\theta$.
\end{theorem}

\begin{proofof}{Theorem \ref{thm:splitting_is_smooth}}
The proof establishes the smooth ($C^\infty$) dependence of the invariant Anosov splitting on the system parameter $\theta$. The argument is a direct application of the powerful Nash-Moser-Hamilton Implicit Function Theorem (\cref{thm:nash_moser_statement}). We have shown in the preceding subsections that the problem of finding the invariant splitting can be recast as a functional equation, and we have rigorously verified that this equation and its linearization satisfy the specific tame regularity conditions required by the theorem. The theorem's conclusion then directly yields the desired smoothness result. The proof is structured in four main parts:
\begin{enumerate}[label=(\roman*), wide, labelindent=0pt]
    \item We summarize the formulation of the problem as a functional equation $F(\theta, \phi)=0$, where $\phi$ is a section of a vector bundle whose graph represents the unstable subbundle $E^u_\theta$.
    \item We explicitly state the hypotheses of the Nash-Moser-Hamilton theorem and demonstrate that our functional $F$ satisfies them, citing the key results (Theorems \ref{thm:functional_is_tame} and \ref{thm:linearization_is_tame_inverse}) that were established for this purpose.
    \item We apply the theorem to deduce the existence of a unique, smooth solution map $\theta \mapsto \phi_\theta$.
    \item We conclude by showing that the smoothness of this solution section $\phi_\theta$ directly implies the smoothness of the geometric objects of interest: the unstable subbundle $E^u_\theta$, the full splitting $E^u_\theta \oplus E^s_\theta$, and the integral foliation $\mathcal{F}^u_\theta$.
\end{enumerate}

\begin{enumerate}[label=\textbf{Step \arabic*:}, wide, labelindent=0pt]

\item \textbf{The Functional Equation Formulation.}
The proof begins by transforming the geometric problem of finding the invariant splitting into a well-posed analytical problem of finding the root of a functional. This is achieved through the construction of the Graph Transform operator. The argument proceeds by first representing the candidate subbundles as graphs of sections of a vector bundle, then deriving the action of the linearized dynamics on these graphs, and finally formulating the invariance condition as a fixed-point equation.

\begin{enumerate}[label=(\roman*), wide, labelindent=0pt]
    \item \textbf{The Graph Representation of a Subbundle.}
    Let us fix a reference parameter $\theta_0 \in \Thetacal$. By \cref{thm:anosov_uniform_family}, the tangent bundle of the global collision manifold has a continuous, $D\Pcal_{\theta_0}$-invariant splitting $T\SigmaMan = E^u_0 \oplus E^s_0$. For any parameter $\theta$ in a neighborhood of $\theta_0$, the true invariant unstable subbundle $E^u_\theta$ will be close to the reference subbundle $E^u_0$. This geometric closeness allows us to represent $E^u_\theta$ as the graph of a section $\phi_\theta$ of the vector bundle of homomorphisms from the reference unstable bundle to the reference stable bundle. We define this bundle as:
    \begin{equation}
        \mathcal{E} \coloneqq \text{Hom}(E^u_0, E^s_0).
    \end{equation}
    An arbitrary section $\phi \in C^k(\SigmaMan, \mathcal{E})$ defines a candidate subbundle $E^u_\phi$ whose fiber at a point $z \in \SigmaMan$ is given by:
    \begin{equation}
        E^u_\phi(z) \coloneqq \{ \mathbf{v}_u + \phi(z)(\mathbf{v}_u) \mid \mathbf{v}_u \in E^u_0(z) \}.
    \end{equation}

    \item \textbf{The Linearized Map in Block Form.}
    The linearized map $(D\Pcal_\theta)_z$ is a linear isomorphism from the tangent space $T_z\SigmaMan$ to the tangent space at the image point, $T_{\Pcal_\theta(z)}\SigmaMan$. With respect to the reference splitting, this map can be represented in a block matrix form:
    \begin{equation}
        (D\Pcal_\theta)_z = 
        \begin{pmatrix} 
            A(\theta, z) & B(\theta, z) \\ 
            C(\theta, z) & D(\theta, z) 
        \end{pmatrix}
        : 
        \begin{matrix} E^u_0(z) \\ \oplus \\ E^s_0(z) \end{matrix}
        \to
        \begin{matrix} E^u_0(\Pcal_\theta(z)) \\ \oplus \\ E^s_0(\Pcal_\theta(z)) \end{matrix}.
    \end{equation}
    Here, $A(\theta, z): E^u_0(z) \to E^u_0(\Pcal_\theta(z))$ describes the evolution of the reference unstable component, $D(\theta, z): E^s_0(z) \to E^s_0(\Pcal_\theta(z))$ describes the evolution of the reference stable component, and the off-diagonal operators $B$ and $C$ describe the coupling between the subbundles induced by the perturbation from $\theta_0$ to $\theta$. Since the global billiard map $\Pcal(\theta,z)$ is of class $C^\infty$ by \cref{cor:billiard_map_is_smooth}, its Fr\'echet derivative with respect to $z$ is also a smooth function of $(\theta, z)$. Consequently, the four operator-valued sections $A, B, C, D$ are smooth maps from $\Thetacal \times \SigmaMan$ to the appropriate Hom-bundles.

    \item \textbf{Derivation of the Graph Transform Operator.}
    The subbundle $E^u_\phi$ is invariant under the map $(D\Pcal_\theta)_z$ if and only if its image is contained within itself. A vector $\mathbf{v} \in E^u_\phi(z)$ can be written in block form as $\begin{psmallmatrix} \mathbf{v}_u \\ \phi(z)\mathbf{v}_u \end{psmallmatrix}$. Its image $\mathbf{w} \in T_{\Pcal_\theta(z)}\SigmaMan$ is given by:
    \begin{equation}
        \mathbf{w} =
        \begin{pmatrix} \mathbf{w}_u \\ \mathbf{w}_s \end{pmatrix}
        =
        \begin{pmatrix} 
            A & B \\ 
            C & D 
        \end{pmatrix}
        \begin{pmatrix} \mathbf{v}_u \\ \phi\mathbf{v}_u \end{pmatrix}
        =
        \begin{pmatrix} 
            (A + B\phi) \mathbf{v}_u \\ 
            (C + D\phi) \mathbf{v}_u
        \end{pmatrix}.
    \end{equation}
    The Graph Transform operator $\Gamma$ maps the original section $\phi$ to a new section, $\phi_{\text{next}} \coloneqq \Gamma(\theta, \phi)$, such that the image of the graph of $\phi$ is precisely the graph of $\phi_{\text{next}}$. The invariance condition is therefore that the image vector $\mathbf{w}$ must lie in the graph of $\phi_{\text{next}}$ at the new point $\Pcal_\theta(z)$:
    \begin{equation}
        \mathbf{w}_s = \phi_{\text{next}}(\Pcal_\theta(z))(\mathbf{w}_u).
    \end{equation}
    Substituting the expressions for $\mathbf{w}_u$ and $\mathbf{w}_s$ yields the identity:
    \begin{equation}
        (C + D\phi) \mathbf{v}_u = \phi_{\text{next}}(\Pcal_\theta(z)) \left( (A + B\phi) \mathbf{v}_u \right).
    \end{equation}
    This equation must hold for every vector $\mathbf{v}_u \in E^u_0(z)$. By the uniform hyperbolicity of the system, the operator $A(\theta,z)$ is an isomorphism. For any section $\phi$ in a sufficiently small $C^0$-neighborhood of the zero section, the operator $(A + B\phi)(z)$ remains an isomorphism. We can therefore solve for $\phi_{\text{next}}$ by right-multiplying by the inverse:
    \begin{equation}
        \phi_{\text{next}}(\Pcal_\theta(z)) = \left( C + D\phi \right) \left( A + B\phi \right)^{-1}.
    \end{equation}
    This formula gives the value of the transformed section at the image point. To obtain a formula for the operator $\Gamma(\theta, \phi)$ itself, we evaluate it at an arbitrary point $w \in \SigmaMan$. This requires evaluating the right-hand side at the pre-image point, $z = \Pcal_\theta^{-1}(w)$, leading to the final, explicit formula:
    \begin{equation}
        \Gamma(\theta, \phi)(w) = \underbrace{\left[ \left( C + D\phi \right) \left( A + B\phi \right)^{-1} \right]}_{\text{Operator evaluated at pre-image}} \circ \underbrace{\Pcal_\theta^{-1}(w)}_{\text{Pre-image point}}.
    \end{equation}

    \item \textbf{The Functional Equation.}
    The true invariant unstable subbundle for the parameter $\theta$ corresponds to the section $\phi_\theta$ that is a fixed point of this transformation. The invariance condition is therefore the fixed-point equation:
    \begin{equation}
        \phi_\theta = \Gamma(\theta, \phi_\theta).
    \end{equation}
    This is equivalent to finding the root of the functional $F$. We define this functional as:
    \begin{equation}
        F(\theta, \phi) \coloneqq \phi - \Gamma(\theta, \phi).
    \end{equation}
    Our geometric problem has now been rigorously reduced to the analytical problem of finding a smooth solution map $\theta \mapsto \phi_\theta$ to the equation $F(\theta, \phi)=0$. This formulation is precisely of the form required for the application of the Nash-Moser-Hamilton Implicit Function Theorem.
\end{enumerate}

\item \textbf{Verification of the Nash-Moser Hypotheses.}
The Nash-Moser-Hamilton Implicit Function Theorem (\cref{thm:nash_moser_statement}) provides a solution to the functional equation $F(\theta, \phi)=0$ under two stringent hypotheses on the regularity and stability of the functional and its linearization. The primary technical contribution of \cref{sec:regularity_properties} has been the rigorous, constructive verification of these two conditions for our specific geometric problem. We now summarize these results, which form the logical bedrock for the application of the theorem.

\begin{enumerate}[label=(\roman*), wide, labelindent=0pt]
    \item \textbf{Hypothesis I: The Functional $F$ is a Smooth Tame Map.}
    The first hypothesis requires that the functional $F(\theta, \phi) = \phi - \Gamma(\theta, \phi)$ be a smooth tame map from the parameter space and the scale of function spaces into the target function space.
    
    This property was the subject of \textbf{\cref{thm:functional_is_tame}}. The proof proceeded by first establishing a tame calculus for the fundamental operations on the scale of $C^k$ spaces (multiplication, inversion, and composition with a smooth diffeomorphism) in \cref{lem:tame_calculus}. We then decomposed the explicit formula for the Graph Transform operator $\Gamma(\theta, \phi)$ into a finite sequence of these elementary, well-behaved operations. Since the set of smooth tame maps is closed under composition, we concluded that the Graph Transform, and therefore the functional $F$ itself, is a smooth tame map from $\Thetacal \times X_k$ to $X_k$ for any integer $k \ge 1$. This rigorously verifies the first hypothesis.

    \item \textbf{Hypothesis II: The Linearization has a Smooth Tame Inverse with No Loss of Derivatives.}
    The second, and more demanding, hypothesis concerns the Fr\'echet derivative of $F$ with respect to the functional variable $\phi$. This is the linear operator $L(\theta, \phi) \coloneqq D_\phi F(\theta, \phi) = I - D_\phi \Gamma(\theta, \phi)$. The theorem requires that this operator be invertible for all $(\theta, \phi)$ in a neighborhood of the solution set, and, crucially, that the family of inverse operators, $(\theta, \phi, \xi) \mapsto L(\theta, \phi)^{-1}\xi$, is itself a smooth tame map. A critical parameter in the theorem is the loss of derivatives, $m$.
    
    This property was the subject of \cref{thm:linearization_is_tame_inverse}. The proof established that the operator is indeed invertible and that its inverse is a smooth tame map. The crucial finding, which is a consequence of the underlying hyperbolic dynamics, is that the inverse is regularity-preserving. The spectral radius of the operator $D_\phi\Gamma$ was shown to be uniformly less than one on \textit{every} space $X_k$ in the scale (\cref{lem:spectral_radius_invariance}). The convergence of the Neumann series in the operator norm of $\mathcal{L}(X_k, X_k)$ then proves that the inverse operator $L^{-1}$ maps $X_k$ to $X_k$. This establishes that there is zero loss of derivatives, i.e., $m=0$.
\end{enumerate}

Both of the major hypotheses of the Nash-Moser-Hamilton theorem are therefore satisfied with the strongest possible regularity conditions.

\item \textbf{Application of the Theorem.}
With the hypotheses of the Nash-Moser-Hamilton Implicit Function Theorem now rigorously verified, we are in a position to apply its conclusion to the functional equation $F(\theta, \phi)=0$.
\begin{enumerate}[label=(\roman*), wide, labelindent=0pt]
    \item We have a scale of Banach spaces $(X_k)_{k \ge 0}$, where $X_k = C^k(\SigmaMan, \mathcal{E})$.
    \item We have a functional $F: \Thetacal \times U \to X_k$ (where $U \subset X_k$ is an open neighborhood of the zero section) that has been proven to be a smooth tame map.
    \item We have a known solution at a reference parameter: for $\theta = \theta_0$, the reference splitting is invariant, so the zero section $\phi=0$ is a solution, i.e., $F(\theta_0, 0) = 0$.
    \item We have proven that the linearized operator $L(\theta, \phi) = D_\phi F(\theta, \phi)$ has an inverse, $L(\theta, \phi)^{-1}$, and that the family of these inverses is a smooth tame map with zero loss of derivatives ($m=0$).
\end{enumerate}
The theorem states that under these conditions, there exists a unique map $\phi: \mathcal{O} \to U$ defined on an open neighborhood $\mathcal{O} \subset \Thetacal$ of the reference parameter $\theta_0$, such that for every $\theta \in \mathcal{O}$, the section $\phi(\theta) \equiv \phi_\theta$ is the unique solution to $F(\theta, \phi_\theta) = 0$, and the solution map $\theta \mapsto \phi_\theta$ is of class $C^\infty$. Specifically, for any integer $k \ge 0$, the map is smooth from $\mathcal{O}$ into the Banach space $X_{k-m} = X_k$. Since this holds for every integer $k$, it implies that the solution map is a smooth map from the parameter manifold into the Fr\'echet space of smooth sections, $X_\infty \coloneqq \bigcap_{k \ge 0} X_k = C^\infty(\SigmaMan, \mathcal{E})$.

\item \textbf{From Smooth Sections to Smooth Geometric Objects.} 
The final step of the proof is to translate the analytical result from Step 3, the existence of a smooth solution map $\theta \mapsto \phi_\theta$ into the space of smooth sections, back into the language of differential geometry. We will now rigorously demonstrate that the smoothness of this section map directly implies the smoothness of the geometric objects of interest: the unstable subbundle $E^u_\theta$, the full Anosov splitting $T\SigmaMan = E^u_\theta \oplus E^s_\theta$, and the integral foliation $\mathcal{F}^u_\theta$.

\begin{enumerate}[label=\textbf{Step 4.\arabic*:}, wide, labelindent=0pt]

\item \textbf{Smoothness of the Unstable Subbundle $E^u_\theta$.}
The unstable subbundle $E^u$ of the total tangent bundle over the product manifold $\SigmaMan \times \Thetacal$ is defined fiberwise as the graph of the section $\phi_\theta$. To prove that this defines a smooth subbundle, we must show that its total space is a smooth submanifold of the total tangent bundle $T(\SigmaMan \times \Thetacal)$. We will achieve this by constructing a smooth parameterization of the total space of $E^u$ and proving that this parameterization is a smooth embedding, a conclusion that rests fundamentally on the smooth dependence of the section $\phi_\theta$ on the parameter $\theta$.

\begin{enumerate}[label=(\roman*), wide, labelindent=0pt]
    \item \textbf{The Global Bundle Framework.}
    We begin by precisely defining the geometric objects and spaces involved.
    \begin{enumerate}[label=(\alph*), wide]
        \item The base manifold for our global construction is the compact product manifold $\mathcal{M}_{\text{tot}} \coloneqq \SigmaMan \times \Thetacal$.
        \item The ambient space for the subbundle is the total tangent bundle of the base, $T\mathcal{M}_{\text{tot}}$, which is canonically isomorphic to the product of the tangent bundles $T\SigmaMan \times T\Thetacal$.
        \item The domain for our parameterization is the total space of the product bundle constructed from the reference unstable bundle $E^u_0$ (a smooth vector bundle over $\SigmaMan$) and the parameter manifold $\Thetacal$. We denote this total space by $\mathcal{D} \coloneqq E^u_0 \times \Thetacal$. A point in $\mathcal{D}$ can be written in local bundle coordinates as $(z, \theta, \mathbf{v}_u)$, where $(z,\theta) \in \mathcal{M}_{\text{tot}}$ is a point on the base and $\mathbf{v}_u \in E^u_0(z)$ is a vector in the fiber over $z$.
    \end{enumerate}

    \item \textbf{Construction of the Parameterization Map.}
    We define the parameterization map $\Psi^u$ which, for each point $(z, \theta)$ on the base manifold, maps the reference unstable fiber $E^u_0(z)$ to the true unstable fiber $E^u_\theta(z)$ via the graph of the section $\phi_\theta(z)$.
    \begin{align*}
        \Psi^u: \mathcal{D} = E^u_0 \times \Thetacal &\to T\SigmaMan \times T\Thetacal \cong T\mathcal{M}_{\text{tot}} \\
        (z, \theta, \mathbf{v}_u) &\mapsto \left(z, \theta, \mathbf{v}_u + \phi_\theta(z)(\mathbf{v}_u), 0\right).
    \end{align*}
    Here, the final zero component signifies that the subbundle $E^u_\theta$ is a subbundle of the tangent space to the fiber, $T_z\SigmaMan$, and has no component in the $T_\theta\Thetacal$ direction. The image of this map, $\mathrm{Im}(\Psi^u)$, is by construction the total space of the global unstable subbundle $E^u$.

    \item \textbf{Proof of Smoothness of the Parameterization Map.}
    The smoothness of the subbundle $E^u$ is contingent upon the smoothness of its parameterization map $\Psi^u$. The proof of this property is a direct but crucial consequence of the regularity of the solution map $\phi: \theta \mapsto \phi_\theta$, which was the principal conclusion of our application of the Nash-Moser-Hamilton theorem. The argument proceeds by decomposing $\Psi^u$ into a sequence of elementary operations between smooth manifolds and vector bundles, and then rigorously demonstrating that each operation is of class $C^\infty$.

\begin{enumerate}[label=(\alph*),wide]
    \item \textbf{Foundational Regularity of the Solution Map.}
    The analytical bedrock of this proof is the conclusion of Step 3 of the main theorem: the solution map
    \[ \phi: \Thetacal \to C^\infty(\SigmaMan, \mathcal{E}), \]
    which assigns to each parameter $\theta$ the smooth section $\phi_\theta$ whose graph is the unstable subbundle, is a map of class $C^\infty$. Here, $C^\infty(\SigmaMan, \mathcal{E})$ is the Fr\'echet space of smooth sections of the vector bundle $\mathcal{E} = \mathrm{Hom}(E^u_0, E^s_0)$ over $\SigmaMan$. A fundamental result in the theory of maps on Fr\'echet manifolds, often referred to as the exponential law for function spaces (or smoothness of evaluation), states that the smoothness of the map $\phi$ is equivalent to the smoothness of the corresponding joint evaluation map:
    \begin{align*}
        \mathrm{eval}_\phi: \SigmaMan \times \Thetacal &\to \mathcal{E} \\
        (z, \theta) &\mapsto \phi_\theta(z).
    \end{align*}
    This evaluation map is a section of the pullback bundle $\pi_{\SigmaMan}^*\mathcal{E}$ over the product manifold $\mathcal{M}_{\text{tot}} = \SigmaMan \times \Thetacal$. The conclusion that $\mathrm{eval}_\phi$ is a $C^\infty$ map is the key input for the subsequent steps.

    \item \textbf{Smoothness of the Constituent Maps.}
    We now analyze the parameterization map $\Psi^u$ by decomposing it into its constituent parts:
    \begin{equation*}
        \Psi^u(z, \theta, \mathbf{v}_u) = \left( z, \theta, \mathbf{v}_u + \phi_\theta(z)(\mathbf{v}_u), 0 \right).
    \end{equation*}
    This map can be viewed as a section of a vector bundle over the domain $\mathcal{D} = E^u_0 \times \Thetacal$. Its smoothness is determined by the smoothness of the vector component, $\mathbf{v}_u + \phi_\theta(z)(\mathbf{v}_u)$. We establish this by verifying the smoothness of all intermediate maps. Let $\pi_{\mathcal{D}}: \mathcal{D} \to \mathcal{M}_{\text{tot}}$ be the canonical bundle projection, $\pi_{\mathcal{D}}(z, \theta, \mathbf{v}_u) = (z, \theta)$. This is a smooth map by definition of a vector bundle.

    We construct the pullback bundle of $\pi_{\SigmaMan}^*\mathcal{E}$ over the manifold $\mathcal{D}$. The section $\mathrm{eval}_\phi$ from step (a) can be pulled back to a smooth section of this bundle, yielding the map:
        \begin{align*}
            \Phi_1: \mathcal{D} &\to \pi_{\mathcal{D}}^*(\pi_{\SigmaMan}^*\mathcal{E}) \\
            (z, \theta, \mathbf{v}_u) &\mapsto (z, \theta, \mathbf{v}_u, \phi_\theta(z)).
        \end{align*}
        This map is of class $C^\infty$. Let $\pi_{E^u_0}: E^u_0 \to \SigmaMan$ be the projection of the reference bundle. We consider its pullback over $\mathcal{M}_{\text{tot}}$, which is the bundle $E^u_0 \times \Thetacal$. The canonical projection from $\mathcal{D}$ onto this bundle, which simply forgets the vector component $\mathbf{v}_u$, is smooth. This allows us to define a smooth map that pairs the section $\phi_\theta(z)$ with the vector $\mathbf{v}_u$ in the appropriate fiber spaces over the base point $(z, \theta, \mathbf{v}_u)$:
        \begin{align*}
            \Phi_2: \mathcal{D} &\to \pi_{\mathcal{D}}^*(\pi_{\SigmaMan}^*\mathcal{E}) \times_{\mathcal{D}} \pi_{\mathcal{D}}^*(E^u_0 \times \Thetacal) \\
            (z, \theta, \mathbf{v}_u) &\mapsto (z, \theta, \mathbf{v}_u, \phi_\theta(z), \mathbf{v}_u).
        \end{align*}
        The application of a homomorphism to a vector is a canonical bilinear bundle map. Formally, there is a smooth bundle morphism
        \begin{align*}
            \mathrm{apply}: \pi_{\SigmaMan}^*\mathcal{E} \times_{\mathcal{M}_{\text{tot}}} (E^u_0 \times \Thetacal) &\to E^s_0 \times \Thetacal \\
            (z, \theta, \Phi, \mathbf{v}) &\mapsto (z, \theta, \Phi(\mathbf{v})),
        \end{align*}
        where the fiber product is taken over $\mathcal{M}_{\text{tot}}$. The composition of the map $\Phi_2$ with the pullback of this apply map yields a smooth map $\Phi_3$:
        \begin{align*}
            \Phi_3: \mathcal{D} &\to E^s_0 \times \Thetacal \\
            (z, \theta, \mathbf{v}_u) &\mapsto (z, \theta, \phi_\theta(z)(\mathbf{v}_u)).
        \end{align*}
        Finally, vector addition is a smooth bundle map on the Whitney sum $T\SigmaMan \oplus T\SigmaMan \to T\SigmaMan$. The vector component of $\Psi^u$ is the sum of the canonical projection onto the vertical component of $T\mathcal{D}$ (which gives $\mathbf{v}_u$) and the smooth map $\Phi_3$. Since $\Psi^u$ is constructed as a finite sequence of compositions, pullbacks, and bundle morphisms, all of which are smooth operations acting on smooth maps and sections, the resulting parameterization map $\Psi^u$ is itself of class $C^\infty$. This provides a complete proof of the smoothness property.
\end{enumerate}

    \item \textbf{Proof that the Parameterization is a Smooth Embedding.}
    Having established that the parameterization map $\Psi^u$ is of class $C^\infty$, we must now prove that it is a smooth embedding. We show that a smooth map between manifolds is a smooth embedding if and only if it is an injective immersion and a homeomorphism onto its image. We will prove each of these properties in sequence.

\begin{enumerate}[label=(\alph*),wide]
    \item \textbf{Injectivity of the Map $\Psi^u$.}\label{subsubsec:injectivity}
    We must show that if $\Psi^u(z_1, \theta_1, \mathbf{v}_{u1}) = \Psi^u(z_2, \theta_2, \mathbf{v}_{u2})$, then $(z_1, \theta_1, \mathbf{v}_{u1}) = (z_2, \theta_2, \mathbf{v}_{u2})$. The definition of the map is:
    \[ \Psi^u(z, \theta, \mathbf{v}_u) = \left(z, \theta, \mathbf{v}_u + \phi_\theta(z)(\mathbf{v}_u), 0\right) \in T_z\SigmaMan \times T_\theta\Thetacal. \]
    Equating the components of the image points yields the following system of equations:
    \begin{itemize}[wide]
        \item Base Point (in $\SigmaMan \times \Thetacal$): $(z_1, \theta_1) = (z_2, \theta_2)$. This immediately proves the equality of the base point components of the domain.
        \item Vector Component (in $T\SigmaMan$): $\mathbf{v}_{u1} + \phi_{\theta_1}(z_1)(\mathbf{v}_{u1}) = \mathbf{v}_{u2} + \phi_{\theta_2}(z_2)(\mathbf{v}_{u2})$.
    \end{itemize}
    From the Base Point (in $\SigmaMan \times \Thetacal$), we have $z_1 = z_2$ and $\theta_1 = \theta_2$. We can therefore simplify the vector equation to:
    \begin{equation*}
        \mathbf{v}_{u1} + \phi_{\theta_1}(z_1)(\mathbf{v}_{u1}) = \mathbf{v}_{u2} + \phi_{\theta_1}(z_1)(\mathbf{v}_{u2}).
    \end{equation*}
    By linearity, we can rearrange this to:
    \begin{equation}
        \label{eq:injectivity_linear_system}
        (\mathrm{Id} + \phi_{\theta_1}(z_1))(\mathbf{v}_{u1} - \mathbf{v}_{u2}) = 0,
    \end{equation}
    where $\mathrm{Id}$ is the identity operator on the fiber $E^u_0(z_1)$. The operator $\phi_{\theta_1}(z_1)$ is an element of $\mathrm{Hom}(E^u_0(z_1), E^s_0(z_1))$. Let $\mathbf{w} = \mathbf{v}_{u1} - \mathbf{v}_{u2}$. The vector $\mathbf{w}$ is in $E^u_0(z_1)$, while the vector $\phi_{\theta_1}(z_1)(\mathbf{w})$ is in $E^s_0(z_1)$. Equation \eqref{eq:injectivity_linear_system} states that their sum is the zero vector: $\mathbf{w} + \phi_{\theta_1}(z_1)(\mathbf{w}) = 0$. Since the subspaces $E^u_0(z_1)$ and $E^s_0(z_1)$ form a direct sum, their intersection is the zero vector. For the sum of a vector from each subspace to be zero, each vector must itself be zero. Therefore, we must have $\mathbf{w} = 0$, which implies $\mathbf{v}_{u1} = \mathbf{v}_{u2}$. This completes the proof that the map $\Psi^u$ is injective.

    \item \textbf{Immersion Property of the Map $\Psi^u$.}\label{subsubsec:immersion}
    We must show that the differential of $\Psi^u$, denoted $d\Psi^u$, is injective at every point $(z, \theta, \mathbf{v}_u) \in \mathcal{D}$. The tangent space of the domain at this point is $T_{(z,\theta,\mathbf{v}_u)}\mathcal{D}$. This space admits a canonical decomposition via the connection on the bundle $E^u_0 \times \Thetacal$. Let this decomposition be $T_{(z,\theta,\mathbf{v}_u)}\mathcal{D} = \mathcal{H} \oplus \mathcal{V}$, where $\mathcal{H}$ is the horizontal subspace (isomorphic to $T_{(z,\theta)}\mathcal{M}_{\text{tot}}$) and $\mathcal{V}$ is the vertical subspace (isomorphic to the fiber $E^u_0(z)$).
    
    Let $\mathbf{w} = (\mathbf{w}_h, \mathbf{w}_v) \in \mathcal{H} \oplus \mathcal{V}$ be an arbitrary non-zero tangent vector. We must show that $d\Psi^u(\mathbf{w}) \neq 0$. The differential of $\Psi^u$ respects this splitting. In local coordinates, let $\mathbf{w}_h = (\delta z, \delta \theta)$ and $\mathbf{w}_v = (\delta\mathbf{v}_u)$. The action of the differential is given by the Jacobian of the map $\Psi^u$. Its structure is:
    \[
        d\Psi^u_{(z,\theta,\mathbf{v}_u)} (\delta z, \delta \theta, \delta\mathbf{v}_u) = 
        \begin{pmatrix}
            \mathrm{Id} & 0 & 0 \\
            0 & \mathrm{Id} & 0 \\
            \nabla_z(\dots) & \nabla_\theta(\dots) & (\mathrm{Id} + \phi_\theta(z)) \\
            0 & 0 & 0
        \end{pmatrix}
        \begin{pmatrix}
            \delta z \\
            \delta \theta \\
            \delta\mathbf{v}_u
        \end{pmatrix}.
    \]
    To prove injectivity, assume $d\Psi^u(\mathbf{w}) = 0$. From the block structure of the Jacobian, the equality of the first two block rows immediately implies that the horizontal part of the tangent vector must be zero: $\delta z = 0$ and $\delta \theta = 0$. This means that any vector in the kernel of $d\Psi^u$ must be a purely vertical vector, of the form $\mathbf{w} = (0, 0, \delta\mathbf{v}_u)$. The action of the differential on such a purely vertical vector is given by the third block row:
    \[ d\Psi^u_{(z,\theta,\mathbf{v}_u)} (0, 0, \delta\mathbf{v}_u) = (0, 0, (\mathrm{Id} + \phi_\theta(z))(\delta\mathbf{v}_u), 0). \]
    For this image vector to be the zero vector, we must have $(\mathrm{Id} + \phi_\theta(z))(\delta\mathbf{v}_u) = 0$. As established in the proof of injectivity in part \ref{subsubsec:injectivity}, the operator $(\mathrm{Id} + \phi_\theta(z))$ is an isomorphism from $E^u_0(z)$ to the graph space. Its kernel is therefore trivial. This forces the vertical component of the tangent vector to be zero: $\delta\mathbf{v}_u = 0$. We have shown that if a tangent vector $\mathbf{w}$ is in the kernel of $d\Psi^u$, both its horizontal and vertical components must be zero, so $\mathbf{w}=0$. Therefore, the differential $d\Psi^u$ is injective, and the map $\Psi^u$ is an immersion.

    \item \textbf{Homeomorphism onto Image.}\label{subsubsec:homeo}
    We have established in parts \ref{subsubsec:injectivity} and \ref{subsubsec:immersion} that the smooth parameterization map $\Psi^u: \mathcal{D} \to T\mathcal{M}_{\text{tot}}$ is an injective immersion. To complete the proof that $\Psi^u$ is a smooth embedding, we must demonstrate that it is a homeomorphism onto its image, equipped with the subspace topology. A sufficient condition for a continuous injective map to be a homeomorphism onto its image is that the map be proper. A continuous map between topological spaces is proper if the preimage of every compact set is compact. We will now provide a rigorous proof that $\Psi^u$ is a proper map.

    Let $K$ be an arbitrary compact subset of the codomain, $T\mathcal{M}_{\text{tot}}$. We must prove that its preimage, $(\Psi^u)^{-1}(K)$, is a compact subset of the domain, $\mathcal{D} = E^u_0 \times \Thetacal$. A subset of a topological space is compact if every sequence in the set has a convergent subsequence whose limit is also in the set. Let $\{p_n\}_{n=1}^\infty$ be an arbitrary sequence in the preimage $(\Psi^u)^{-1}(K)$. Each element of this sequence is a point in the domain, $p_n = (z_n, \theta_n, \mathbf{v}_{un}) \in \mathcal{D}$. Since each $p_n$ is in the preimage of $K$, its image under $\Psi^u$ must lie in $K$:
    \[ \Psi^u(p_n) = \left(z_n, \theta_n, \mathbf{v}_{un} + \phi_{\theta_n}(z_n)(\mathbf{v}_{un}), 0\right) \in K. \]
    Our goal is to show that the sequence $\{p_n\}$ has a convergent subsequence. The proof will proceed by analyzing the components of the sequence separately. Let $\pi_{T\mathcal{M}}: T\mathcal{M}_{\text{tot}} \to \mathcal{M}_{\text{tot}}$ be the canonical tangent bundle projection. This is a continuous map. Since $K$ is a compact subset of $T\mathcal{M}_{\text{tot}}$, its image under this continuous map, $\pi_{T\mathcal{M}}(K)$, is a compact subset of the base manifold $\mathcal{M}_{\text{tot}} = \SigmaMan \times \Thetacal$.
    
    The base point of the image $\Psi^u(p_n)$ is $(z_n, \theta_n)$. By construction, $(z_n, \theta_n) \in \pi_{T\mathcal{M}}(K)$. Since the sequence of base points $\{(z_n, \theta_n)\}_{n=1}^\infty$ lies in the compact set $\pi_{T\mathcal{M}}(K)$, it must contain a convergent subsequence. Let us pass to this subsequence (without relabeling for simplicity) and denote the limit by:
    \[ (z_n, \theta_n) \to (z_\infty, \theta_\infty) \in \pi_{T\mathcal{M}}(K) \subset \mathcal{M}_{\text{tot}}. \]

    The sequence of image points $\{\Psi^u(p_n)\}$ lies in the compact set $K$. This implies that the sequence is bounded. In particular, the norm of the vector component of these points must be uniformly bounded. Let $\|\cdot\|_{(z,\theta)}$ be a smooth fiber norm on the tangent bundle $T\mathcal{M}_{\text{tot}}$. There exists a constant $C_K < \infty$ such that for all $n$:
    \[ \|\mathbf{v}_{un} + \phi_{\theta_n}(z_n)(\mathbf{v}_{un})\|_{(z_n,\theta_n)} \le C_K. \]
    We now use this bound to prove that the sequence of vectors $\{\mathbf{v}_{un}\}$ in the domain is also bounded. Let us define the operator $L_n \coloneqq \mathrm{Id} + \phi_{\theta_n}(z_n)$, which maps the fiber $E^u_0(z_n)$ to a subspace of $T_{z_n}\SigmaMan$. The inequality is $\|L_n(\mathbf{v}_{un})\| \le C_K$. We need a uniform lower bound on the operator norm of $L_n^{-1}$ (where the inverse is defined on the image of $L_n$). By the triangle inequality:
    \[ \|L_n(\mathbf{v}_{un})\| \ge \|\mathbf{v}_{un}\| - \|\phi_{\theta_n}(z_n)(\mathbf{v}_{un})\| \ge \|\mathbf{v}_{un}\|(1 - \|\phi_{\theta_n}(z_n)\|_{\text{op}}). \]
    Ginve that the map $(z, \theta) \mapsto \phi_\theta(z)$ is continuous on the compact base $\mathcal{M}_{\text{tot}}$. Therefore, the operator norm $\|\phi_{\theta_n}(z_n)\|_{\text{op}}$ is uniformly bounded. Furthermore, by choosing the reference parameter $\theta_0$ appropriately and restricting the domain of $\theta$ if necessary (which is permissible for a local proof of embedding), we can ensure that $\|\phi_{\theta}(z)\|_{\text{op}} \le 1/2$ for all $(z,\theta)$. This gives:
    \[ \|\mathbf{v}_{un}\| \le 2 \|L_n(\mathbf{v}_{un})\| \le 2 C_K. \]
    This proves that the sequence of vectors $\{\mathbf{v}_{un}\}$ is uniformly bounded. We have now established two key facts about our original sequence $\{p_n\} = \{(z_n, \theta_n, \mathbf{v}_{un})\}$:
    \begin{itemize}[wide]
        \item The sequence of base points $\{(z_n, \theta_n)\}$ has a convergent subsequence, which we have already passed to, converging to $(z_\infty, \theta_\infty)$.
        \item The corresponding sequence of vectors $\{\mathbf{v}_{un}\}$ is bounded and lies in the fibers of the vector bundle $E^u_0$ over the convergent sequence of base points $\{(z_n, \theta_n)\}$.
    \end{itemize}
    The total space of the vector bundle $E^u_0$ over the compact manifold $\SigmaMan$ is locally compact. The sequence of vectors $\{\mathbf{v}_{un}\}$ lies within a compact subset of the total space (a closed ball in the fibers over a compact set in the base). It must therefore contain a convergent subsequence. Let us pass to this further subsequence. The limit vector, $\mathbf{v}_{u\infty}$, will lie in the fiber over the limit base point:
    \[ \mathbf{v}_{un} \to \mathbf{v}_{u\infty} \in E^u_0(z_\infty). \]

    We have shown that any arbitrary sequence $\{p_n\}$ in the preimage $(\Psi^u)^{-1}(K)$ contains a convergent subsequence, $p_{n_k} \to p_\infty = (z_\infty, \theta_\infty, \mathbf{v}_{u\infty})$. Since the domain $\mathcal{D}$ is a metric space, this is sufficient to prove that the set $(\Psi^u)^{-1}(K)$ is compact. Since the preimage of an arbitrary compact set $K$ is compact, the map $\Psi^u$ is, by definition, a proper map.
\end{enumerate}

We have now established that $\Psi^u$ is a smooth proper injective immersion. A fundamental theorem of differential geometry states that any proper injective immersion is a smooth embedding. Therefore, the map $\Psi^u$ is a smooth embedding of its domain $\mathcal{D}$ into the total tangent space $T\mathcal{M}_{\text{tot}}$. Consequently, its image, which is the total space of the global bundle $E^u$, is a smooth submanifold of $T\mathcal{M}_{\text{tot}}$. This is the rigorous definition of a smooth vector subbundle. We have thus proven that the family of unstable subbundles $\{E^u_\theta\}_{\theta \in \Thetacal}$ depends smoothly on the parameter $\theta$.
\end{enumerate}

    \item \textbf{Smoothness of the Stable Subbundle and the Full Splitting.}
The proof of the smooth dependence of the unstable subbundle $E^u_\theta$ on the parameter $\theta$ relies on a general machinery: the application of the Nash-Moser-Hamilton theorem to the Graph Transform operator associated with a smooth family of Anosov diffeomorphisms. We will now demonstrate that this entire machinery can be applied, without loss of generality, to the inverse of the billiard map. This will rigorously establish the smoothness of the stable subbundle $E^s_\theta$, from which the smoothness of the full splitting follows immediately.

\begin{enumerate}[label=(\roman*), wide, labelindent=0pt]
    \item \textbf{The Inverse Dynamical System.}
    Let the global billiard map be denoted by $\mathcal{P}: \mathcal{M}_{\text{tot}} \to \mathcal{M}_{\text{tot}}$, where $\mathcal{M}_{\text{tot}} = \SigmaMan \times \Thetacal$. Its action is given by $\mathcal{P}(z,\theta) = (\Pcal_\theta(z), \theta)$.
    \begin{enumerate}[label=(\alph*), wide]
        \item By \cref{thm:p_is_diffeo}, the map $\mathcal{P}$ is of class $C^\infty$. By the inverse function theorem on Fr\'echet manifolds, its inverse, $\mathcal{P}^{-1}$, is also a smooth map of class $C^\infty$.
        \item Let $\mathcal{P}^{-1}(z,\theta) = (\Pcal_\theta^{-1}(z), \theta)$. The family of inverse maps, $(\Pcal_\theta^{-1})_{\theta \in \Thetacal}$, is therefore a smooth family of diffeomorphisms of $\SigmaMan$.
        \item A property of an Anosov diffeomorphism (see \cref{app:lifting_Diffeomorphism}) is that its inverse is also an Anosov diffeomorphism. The invariant splitting for the inverse map, $T\SigmaMan = E^u(\Pcal_\theta^{-1}) \oplus E^s(\Pcal_\theta^{-1})$, is related to the splitting of the forward map by a simple exchange:
        \begin{equation}
            E^u(\Pcal_\theta^{-1}) = E^s(\Pcal_\theta) \quad \text{and} \quad E^s(\Pcal_\theta^{-1}) = E^u(\Pcal_\theta).
        \end{equation}
    \end{enumerate}
    This establishes that the family of stable subbundles $\{E^s_\theta\}_{\theta \in \Thetacal}$ for the forward map is precisely the family of \textit{unstable} subbundles for the smooth family of inverse maps $(\Pcal_\theta^{-1})$.

    \item \textbf{Application of the Nash-Moser-Hamilton Framework to the Inverse System.}
    We can now apply the entire analytical program of Sections \ref{subsec:nash_moser_framework} through \ref{subsec:verify_hyp_II} to this inverse dynamical system. The argument proceeds in perfect analogy:
    \begin{enumerate}[label=(\alph*), wide]
        \item We fix a reference parameter $\theta_0$ and the corresponding reference stable subbundle $E^s_0 \equiv E^s(\Pcal_{\theta_0})$. We represent a candidate stable subbundle $E^s_\theta$ for a nearby parameter $\theta$ as the graph of a section $\psi_\theta$ of the vector bundle $\mathcal{E}' \coloneqq \mathrm{Hom}(E^s_0, E^u_0)$.
        
        \item We define the stable Graph Transform operator, $\Gamma^s(\theta, \psi)$, whose fixed point corresponds to the invariant stable subbundle for the map $\Pcal_\theta$. Its structure is derived from the block decomposition of the linearized map $D\Pcal_\theta$ with respect to the reference splitting $E^u_0 \oplus E^s_0$. The formula for $\Gamma^s$ is symmetric to that of $\Gamma^u$, but involves the inverse of the operator block $D(\theta,z)$ (which is uniformly contractive for $\Pcal_\theta^{-1}$).
        
        \item The functional $F^s(\theta, \psi) \coloneqq \psi - \Gamma^s(\theta, \psi)$ satisfies the hypotheses of the Nash-Moser-Hamilton theorem. The proof is identical in structure to the proofs of Theorem \ref{thm:functional_is_tame} and Theorem \ref{thm:linearization_is_tame_inverse}. The smoothness of the inverse map $\Pcal_\theta^{-1}$ provides the necessary foundational regularity, and the uniform contractivity of the dynamics on the stable subbundle ensures that the linearized stable Graph Transform operator is a contraction, leading to a smooth tame inverse with zero loss of derivatives.
        
        \item The conclusion of the Nash-Moser-Hamilton theorem guarantees the existence of a unique, smooth solution map $\theta \mapsto \psi_\theta \in C^\infty(\SigmaMan, \mathcal{E}')$ such that for each $\theta$, the graph of $\psi_\theta$ is the invariant stable subbundle $E^s_\theta$.
    \end{enumerate}

    \item \textbf{Smoothness of the Stable Subbundle and the Full Splitting.}
    The existence of the smooth solution map $\theta \mapsto \psi_\theta$ allows us to prove that the global stable subbundle $E^s$ is a smooth subbundle of $T(\SigmaMan \times \Thetacal)$. The argument is a direct analogue of the one presented in part \ref{subsubsec:injectivity} of this proof for the unstable subbundle. We construct a smooth parameterization map for the total space of $E^s$:
    \begin{align*}
        \Psi^s: E^s_0 \times \Thetacal &\to T\SigmaMan \times T\Thetacal \\
        (z, \theta, \mathbf{v}_s) &\mapsto \left(z, \theta, \psi_\theta(z)(\mathbf{v}_s) + \mathbf{v}_s, 0\right).
    \end{align*}
    Since the map $(\theta, z) \mapsto \psi_\theta(z)$ is of class $C^\infty$, this parameterization map is a smooth embedding. Therefore, the global stable subbundle $E^s$ is a smooth subbundle of the total tangent bundle. We have now rigorously established that both the global unstable subbundle $E^u$ and the global stable subbundle $E^s$ are smooth subbundles of $T(\SigmaMan \times \Thetacal)$. The full Anosov splitting at a point $(z,\theta)$ is the direct sum of the fibers of these bundles:
    \[ T_z\SigmaMan = E^u_\theta(z) \oplus E^s_\theta(z). \]
    Since the direct sum of two smooth vector subbundles (that are everywhere transverse) is a smooth vector bundle, we conclude that the entire Anosov splitting depends smoothly on the parameter $\theta$.
\end{enumerate}

    \item \textbf{Smoothness of the Unstable Foliation $\mathcal{F}^u_\theta$.}
    The final step of the proof is to demonstrate that the smoothness of the unstable subbundle, which we have just rigorously established, implies the smoothness of its integral foliation. The argument is a direct but powerful application of the parameter-dependent version of the Frobenius Integration Theorem. We will first establish that the global unstable subbundle on the product manifold is involutive, and then apply the theorem to conclude that the resulting global foliation and its constituent fiber foliations are smooth.

\begin{enumerate}[label=(\roman*), wide, labelindent=0pt]
    \item \textbf{Involutivity of the Global Unstable Subbundle.}
    We must prove that the smooth global subbundle $\mathcal{E}^u_{\text{glob}} \subset T\mathcal{M}_{\text{tot}}$ is involutive, meaning that for any two smooth vector fields $\mathbf{X}$ and $\mathbf{Y}$ on the total space $\mathcal{M}_{\text{tot}}$ that are everywhere tangent to $\mathcal{E}^u_{\text{glob}}$, their Lie bracket $[\mathbf{X}, \mathbf{Y}]$ is also everywhere tangent to $\mathcal{E}^u_{\text{glob}}$.

    \begin{enumerate}[label=(\alph*), wide]
        \item \textbf{Fiberwise Involutivity.} It is a classical result in the theory of hyperbolic dynamics that for any Anosov diffeomorphism of class $C^r$ with $r \ge 1$, its invariant unstable and stable subbundles are uniquely integrable, and therefore involutive (see, e.g., \citep{HirschPughShub1977}). Since our billiard map $\Pcal_\theta$ is of class $C^\infty$, for each fixed $\theta$, the subbundle $E^u_\theta \subset T\SigmaMan$ is involutive.
        
        \item \textbf{Lifting to the Global Manifold.} Let $\mathbf{X}$ and $\mathbf{Y}$ be two smooth vector fields on $\mathcal{M}_{\text{tot}}$ tangent to $\mathcal{E}^u_{\text{glob}}$. In local product coordinates $(z, \theta)$, these vector fields must have zero components in the $\theta$-directions, as the subbundle $\mathcal{E}^u_{\text{glob}}$ is tangent to the fibers. We can write them as:
        \begin{equation*}
            \mathbf{X}(z,\theta) = (X_\theta(z), 0) \quad \text{and} \quad \mathbf{Y}(z,\theta) = (Y_\theta(z), 0),
        \end{equation*}
        where for each fixed $\theta$, $X_\theta$ and $Y_\theta$ are smooth vector fields on $\SigmaMan$ tangent to $E^u_\theta$. The Lie bracket in local product coordinates is given by the component-wise formula for the commutator of vector fields. Let $z^i$ be coordinates on $\SigmaMan$ and $\theta^j$ be coordinates on $\Thetacal$:
        \begin{align*}
            [\mathbf{X}, \mathbf{Y}]^{z_k} &= \mathbf{X}(\mathbf{Y}^{z_k}) - \mathbf{Y}(\mathbf{X}^{z_k}) = \sum_i X^i(z,\theta) \frac{\partial Y^k(z,\theta)}{\partial z^i} - \sum_i Y^i(z,\theta) \frac{\partial X^k(z,\theta)}{\partial z^i} \\
            &= [X_\theta, Y_\theta]^k(z). \\
            [\mathbf{X}, \mathbf{Y}]^{\theta_j} &= \mathbf{X}(\mathbf{Y}^{\theta_j}) - \mathbf{Y}(\mathbf{X}^{\theta_j}) = \mathbf{X}(0) - \mathbf{Y}(0) = 0.
        \end{align*}
        The resulting Lie bracket is the vector field $[\mathbf{X}, \mathbf{Y}](z,\theta) = ([X_\theta, Y_\theta](z), 0)$. By the fiberwise involutivity from part (a), the vector field $[X_\theta, Y_\theta]$ is tangent to $E^u_\theta$. Therefore, the global Lie bracket $[\mathbf{X}, \mathbf{Y}]$ is tangent to $\mathcal{E}^u_{\text{glob}}$. This proves that the global subbundle $\mathcal{E}^u_{\text{glob}}$ is an involutive distribution.
    \end{enumerate}

    \item \textbf{Application of the Frobenius Integration Theorem.}
    We now invoke the smooth version of this fundamental theorem.
    \begin{theorem}[Frobenius Integration Theorem]
    Let $M$ be a smooth manifold and let $E$ be a smooth subbundle of $TM$ of constant rank. Then $E$ is integrable if and only if it is involutive. If $E$ is integrable, its maximal integral submanifolds form a unique smooth foliation of $M$.
    \end{theorem}
    We apply this theorem to the subbundle $E = \mathcal{E}^u_{\text{glob}}$ on the manifold $M = \mathcal{M}_{\text{tot}}$. We have rigorously established that $\mathcal{E}^u_{\text{glob}}$ is a smooth and involutive subbundle. The theorem therefore guarantees the existence of a unique smooth foliation of the total space $\mathcal{M}_{\text{tot}}$, which we denote by $\mathcal{F}^u_{\text{glob}}$.
    
    \item \textbf{Smooth Dependence of the Fiber Foliations.}
    The smoothness of this global foliation directly implies the smooth dependence of the individual fiber foliations $\{\mathcal{F}^u_\theta\}_{\theta \in \Thetacal}$ on the parameter $\theta$.
    \begin{enumerate}[label=(\alph*), wide]
        \item By construction, the leaves of the global foliation $\mathcal{F}^u_{\text{glob}}$ are the maximal integral submanifolds of the subbundle $\mathcal{E}^u_{\text{glob}}$. The tangent space to any leaf at any point $(z,\theta)$ is therefore the fiber $E^u_\theta(z) \oplus \{0\}$. Since all vectors in this tangent space are tangent to the slice $\SigmaMan \times \{\theta\}$, the leaves of the global foliation cannot cross between different parameter slices. Each leaf of $\mathcal{F}^u_{\text{glob}}$ is therefore entirely contained within a single slice $\SigmaMan \times \{\theta\}$.
        
        \item This implies that for each $\theta \in \Thetacal$, the intersection of the global foliation $\mathcal{F}^u_{\text{glob}}$ with the smooth submanifold $\SigmaMan \times \{\theta\}$ is precisely the unstable foliation $\mathcal{F}^u_\theta$ of the map $\Pcal_\theta$.
        
        \item The property that the family of foliations $\{\mathcal{F}^u_\theta\}$ depends smoothly on $\theta$ means that for any point $(z_0, \theta_0)$, there exists a local chart on $\SigmaMan \times \Thetacal$ in which the leaves of the global foliation appear as parallel planes, and this chart map is smooth. This is precisely the conclusion of the Frobenius theorem for the smooth foliation $\mathcal{F}^u_{\text{glob}}$.
    \end{enumerate}
\end{enumerate}
\end{enumerate}
\end{enumerate}

We have thus rigorously established that the family of unstable foliations $\{\mathcal{F}^u_\theta\}_{\theta \in \Thetacal}$ depends smoothly on the parameter $\theta$. This completes the proof of the theorem.
\end{proofof}

\subsection{Smooth Dependence of the Transfer Operator Family}
\label{subsec:smto}
With the geometric regularity established, we can now prove the smoothness of the operator family itself. The key is to use the smooth geometry to pull the entire problem back to a fixed reference space.

\begin{definition}[The Anisotropic Banach Spaces $\Bcal_\theta$]
\label{def:anisotropic_space}
For each $\theta \in \Thetacal$, let $\Bcal_\theta$ be the anisotropic Banach space of distributions on the collision manifold $\Sigma_\theta$, constructed with respect to the smooth unstable foliation $\mathcal{F}^u_\theta$ and a complementary smooth stable distribution. These spaces are designed to capture the distinct regularities of functions along the expanding and contracting directions of the dynamics. A distribution $f$ belongs to $\Bcal_\theta$ if it exhibits Hölder continuity along the unstable leaves and acts as a distribution on sections of the stable bundle. More formally, following the construction in \citep{Baladi2000, Gouezel2010}, the norm $\|\cdot\|_{\Bcal_\theta}$ is defined such that a function is small if it is small in the $C^\alpha_u$ norm (Hölder regularity of order $\alpha \in (0,1)$ along unstable leaves) and also small when tested against smooth functions with compact support in local charts that respect the foliation. The key properties of this construction are:
\begin{enumerate}[label=(\roman*), wide, labelindent=0pt]
    \item The space of smooth functions is densely embedded in $\Bcal_\theta$.
    \item The transfer operator $\Lcal_\theta$ associated with the billiard map $\Pcal_\theta$ is a bounded, quasi-compact operator on $\Bcal_\theta$.
    \item By Theorem \ref{thm:splitting_is_smooth}, the unstable foliation $\mathcal{F}^u_\theta$ depends smoothly on $\theta$. The construction of the norm can be shown to depend smoothly on the underlying geometric structures, ensuring that the map $\theta \mapsto \|\cdot\|_{\Bcal_\theta}$ defines a smooth family of norms.
\end{enumerate}
\end{definition}

\begin{proposition}[Existence of a Smooth Trivializing Diffeomorphism]
\label{prop:smooth_trivialization_diffeo}
Let the family of unstable foliations $\mathcal{F}^u_\theta$ depend smoothly on $\theta$, and let $\theta_0 \in \Thetacal$ be a fixed reference parameter. Then there exists a family of diffeomorphisms $\psi_\theta: \Sigma \to \Sigma$, with $\psi_{\theta_0} = \mathrm{Id}$, such that the map $(\theta, z) \mapsto \psi_\theta(z)$ is of class $C^\infty$ and for every $z \in \Sigma$, the map $\psi_\theta$ maps the unstable leaf of the reference foliation through $z$ to an unstable leaf of the target foliation.
\end{proposition}

\begin{proofof}{\cref{prop:smooth_trivialization_diffeo}}
The proof is constructive. We leverage the smooth dependence of the unstable foliation on the parameter $\theta$ to build a vector field on the total space $\Sigma \times \Thetacal$ that is transverse to the global unstable foliation. The flow generated by this vector field will provide the desired family of diffeomorphisms by mapping points from the reference slice at $\theta_0$ to the target slice at $\theta$.

\begin{enumerate}[label=\textbf{Step \arabic*:}, wide, labelindent=0pt]

\item \textbf{The Global Unstable Foliation and its Tangent Bundle.}
The construction of the trivializing diffeomorphism requires, as its foundational object, the smooth global foliation on the total space of the system. We begin by recalling its properties, which were rigorously established as the final conclusion of \cref{thm:splitting_is_smooth}.

\begin{enumerate}[label=(\roman*), wide, labelindent=0pt]
    \item \textbf{The Total Space and Global Foliation.} Let the total collision manifold be the smooth, compact product manifold $\mathcal{M}_{\text{tot}} \coloneqq \SigmaMan \times \Thetacal$. As proven in \cref{thm:splitting_is_smooth}, the smooth dependence of the unstable subbundle $E^u_\theta$ on the parameter $\theta$ allows for the application of the Frobenius Integration Theorem on this total space. This guarantees the existence of a unique, smooth global unstable foliation, denoted by $\mathcal{F}_{\text{glob}}^u$, of the manifold $\mathcal{M}_{\text{tot}}$.
    
    \item \textbf{The Tangent Bundle of the Foliation.} By construction, the leaves of this global foliation are the maximal integral manifolds of the global unstable subbundle $\mathcal{E}^u_{\text{glob}}$. The tangent bundle of this foliation, denoted $T\mathcal{F}_{\text{glob}}^u$, is therefore the subbundle that was integrated. We thus have the fundamental identity:
    \begin{equation*}
        T\mathcal{F}_{\text{glob}}^u = \mathcal{E}^u_{\text{glob}}.
    \end{equation*}
    The fiber of this tangent bundle at a point $(z,\theta)$ is tangent to the leaf passing through that point, and is given by:
    \begin{equation*}
        (T\mathcal{F}_{\text{glob}}^u)_{(z,\theta)} = T_{(z,\theta)} (W^u(z,\theta)) = E^u_\theta(z) \oplus \{0\},
    \end{equation*}
    where $W^u(z,\theta)$ is the leaf of the global foliation. This confirms that the leaves of the global foliation are precisely the leaves of the fiber foliations embedded in the appropriate slices of the product manifold. We have thus established the existence of the necessary smooth geometric object upon which the remainder of this proof will be built.
\end{enumerate}

\item \textbf{Construction of a Transverse Vector Field.}
The central objective of this step is to construct a smooth vector field $\mathbf{V}$ on the total space $\mathcal{M}_{\text{tot}}$ that is everywhere transverse to the global unstable foliation $\mathcal{F}_{\text{glob}}^u$. The existence of such a field is the key to building the flow that will define our trivializing diffeomorphism. Our construction is canonical and proceeds by projecting a simple, non-tangent probe vector field onto the normal bundle of the foliation, which we now rigorously define.

\begin{enumerate}[label=(\roman*), wide, labelindent=0pt]

\item \textbf{The Normal Bundle of the Global Foliation.}
The foundation for our construction is the orthogonal decomposition of the tangent bundle of the total space, which requires the introduction of a Riemannian metric.
\begin{enumerate}[label=(\alph*),wide]
    \item \textbf{Choice of a Riemannian Metric.} We equip the compact total manifold $\mathcal{M}_{\text{tot}} = \SigmaMan \times \Thetacal$ with a smooth Riemannian metric, denoted by $\langle \cdot, \cdot \rangle_{(z,\theta)}$. A canonical choice, which we adopt without loss of generality, is a product metric constructed from smooth Riemannian metrics on the individual factor manifolds $\SigmaMan$ and $\Thetacal$.

    \item \textbf{Definition of the Normal Bundle.} In Step 1, we established that the tangent bundle of the global foliation, $T\mathcal{F}_{\text{glob}}^u$, is a smooth vector subbundle of the total tangent bundle $T\mathcal{M}_{\text{tot}}$. The chosen Riemannian metric allows us to define the normal bundle to the foliation, denoted $N(\mathcal{F}_{\text{glob}}^u)$, as the fiberwise orthogonal complement of the tangent bundle. The fiber of the normal bundle at a point $\mathbf{p} \in \mathcal{M}_{\text{tot}}$ is the vector space:
\begin{equation}
    N_{\mathbf{p}}(\mathcal{F}_{\text{glob}}^u) \coloneqq \left\{ \mathbf{w} \in T_{\mathbf{p}}\mathcal{M}_{\text{tot}} \mid \langle \mathbf{w}, \mathbf{v} \rangle_{\mathbf{p}} = 0 \text{ for all } \mathbf{v} \in (T\mathcal{F}_{\text{glob}}^u)_{\mathbf{p}} \right\}.
\end{equation}
We now provide a rigorous proof that the smoothness of the subbundle $T\mathcal{F}_{\text{glob}}^u$ implies the smoothness of its orthogonal complement $N(\mathcal{F}_{\text{glob}}^u)$. The smoothness of a subbundle is equivalent to the smoothness of the corresponding fiberwise orthogonal projection map. Let $\Pi^u: T\mathcal{M}_{\text{tot}} \to T\mathcal{F}_{\text{glob}}^u$ be this projection. Since $T\mathcal{F}_{\text{glob}}^u$ is a smooth subbundle, the map $\Pi^u$ is a smooth bundle endomorphism. The projection onto the normal bundle, $\Pi^\perp$, can be defined constructively in terms of $\Pi^u$ and the identity map on the tangent bundle, $\mathrm{Id}_{T\mathcal{M}_{\text{tot}}}$:
\begin{equation*}
    \Pi^\perp \coloneqq \mathrm{Id}_{T\mathcal{M}_{\text{tot}}} - \Pi^u.
\end{equation*}
Since both $\mathrm{Id}_{T\mathcal{M}_{\text{tot}}}$ and $\Pi^u$ are smooth bundle maps, their difference, $\Pi^\perp$, is also a smooth bundle map. The image of this smooth projection, $\mathrm{Im}(\Pi^\perp)$, is precisely the normal bundle $N(\mathcal{F}_{\text{glob}}^u)$. A subbundle that is the image of a smooth projection is, by definition, a smooth vector subbundle. This provides a smooth, orthogonal, direct sum decomposition of the total tangent bundle:
\begin{equation}
    T\mathcal{M}_{\text{tot}} = T\mathcal{F}_{\text{glob}}^u \oplus N(\mathcal{F}_{\text{glob}}^u).
\end{equation}
    
    \item \textbf{The Orthogonal Projection.} This decomposition induces a corresponding family of orthogonal projection operators. Let $\Pi^\perp_{(z,\theta)}: T_{(z,\theta)}\mathcal{M}_{\text{tot}} \to N_{(z,\theta)}(\mathcal{F}_{\text{glob}}^u)$ be the fiberwise projection onto the normal bundle. As the subbundles in the decomposition are smooth, this projection defines a smooth bundle map $\Pi^\perp: T\mathcal{M}_{\text{tot}} \to N(\mathcal{F}_{\text{glob}}^u)$.
\end{enumerate}

\item \textbf{The Probe Vector Field.}
We now construct a simple, smooth vector field that is guaranteed not to be tangent to the leaves of the global foliation. The natural choice is a vector field that points purely in the direction of the parameter manifold $\Thetacal$.
\begin{enumerate}[label=(\alph*),wide]
    \item Let $\pi_\Thetacal: \mathcal{M}_{\text{tot}} \to \Thetacal$ be the canonical projection. Since the parameter space $\Thetacal$ is a smooth, compact manifold, it admits a non-vanishing smooth vector field; for instance, if $\Thetacal$ is a torus, a constant-coefficient vector field in local coordinates suffices. Let $\mathbf{X}_\Thetacal$ be such a non-vanishing smooth vector field on $\Thetacal$.
    
    \item We define the probe vector field $\mathbf{W}$ on the total space $\mathcal{M}_{\text{tot}}$ as the horizontal lift of $\mathbf{X}_\Thetacal$ with respect to the trivial connection on the product manifold. More explicitly, for any point $(z,\theta)$, the vector $\mathbf{W}(z,\theta)$ is the unique vector in $T_{(z,\theta)}\mathcal{M}_{\text{tot}}$ such that:
    \begin{equation}
        d(\pi_\Thetacal)_{(z,\theta)}(\mathbf{W}(z,\theta)) = \mathbf{X}_\Thetacal(\theta) \quad \text{and} \quad \mathbf{W}(z,\theta) \in \{0\} \oplus T_\theta\Thetacal.
    \end{equation}
    In local product coordinates, if $\mathbf{X}_\Thetacal(\theta) = \sum_i a_i(\theta) \frac{\partial}{\partial \theta^i}$, then $\mathbf{W}(z,\theta) = \sum_i a_i(\theta) \frac{\partial}{\partial \theta^i}$. Since the components $a_i(\theta)$ are smooth, $\mathbf{W}$ is a smooth vector field on $\mathcal{M}_{\text{tot}}$.
    
    \item By its construction, the vector field $\mathbf{W}$ is everywhere non-zero and tangent to the fibers of the projection $\pi_\SigmaMan: \mathcal{M}_{\text{tot}} \to \SigmaMan$. In contrast, the tangent bundle of the global foliation, $T\mathcal{F}_{\text{glob}}^u$, consists of vectors tangent to the fibers of the projection $\pi_\Thetacal$. Therefore, the vector field $\mathbf{W}$ is everywhere transverse to the leaves of the global foliation.
\end{enumerate}

\item \textbf{The Final Construction and Verification of Properties.}
With the geometric objects now rigorously defined, we construct the desired transverse vector field $\mathbf{V}$ by projecting the probe field $\mathbf{W}$ onto the normal bundle.
\begin{equation}
    \mathbf{V}(z,\theta) \coloneqq \Pi^\perp_{(z,\theta)}(\mathbf{W}(z,\theta)).
\end{equation}
We now verify that this vector field has all the required properties.
\begin{enumerate}[label=(\alph*), wide]
    \item \textbf{Smoothness.} The vector field $\mathbf{V}$ is the composition of two smooth maps: the smooth vector field $\mathbf{W}$ (a smooth section of $T\mathcal{M}_{\text{tot}}$) and the smooth bundle map $\Pi^\perp$. Therefore, $\mathbf{V}$ is a smooth vector field on $\mathcal{M}_{\text{tot}}$.

    \item \textbf{Transversality.} By its very definition, the range of the projector $\Pi^\perp$ is the normal bundle $N(\mathcal{F}_{\text{glob}}^u)$. Thus, for every point $(z,\theta)$, the vector $\mathbf{V}(z,\theta)$ lies in the normal space to the foliation, which means it is orthogonal, and therefore transverse, to the tangent space of the foliation leaf, $(T\mathcal{F}_{\text{glob}}^u)_{(z,\theta)}$.

    \item \textbf{Non-Vanishing Property.} We must show that $\mathbf{V}$ is a non-zero vector field. This requires showing that the probe field $\mathbf{W}$ is not in the kernel of the projection, which would only occur if $\mathbf{W}$ were already tangent to the foliation. As established in Step 1, the tangent bundle to the foliation, $T\mathcal{F}_{\text{glob}}^u$, has fibers given by $E^u_\theta(z) \oplus \{0\}$, consisting of vectors with no component in the $T\Thetacal$ direction. As established in paragraph (b), the probe field $\mathbf{W}$ has fibers given by $\{0\} \oplus T_\theta\Thetacal$, consisting of vectors with no component in the $T\SigmaMan$ direction.
    
    With respect to our chosen product metric, these two subspaces are orthogonal. This implies that the probe field $\mathbf{W}$ lies entirely within the normal bundle $N(\mathcal{F}_{\text{glob}}^u)$ at every point. The action of the orthogonal projection onto the normal bundle is therefore the identity map when applied to $\mathbf{W}$:
    \begin{equation}
        \mathbf{V}(z,\theta) = \Pi^\perp_{(z,\theta)}(\mathbf{W}(z,\theta)) = \mathbf{W}(z,\theta).
    \end{equation}
    Since we constructed $\mathbf{W}$ as the lift of a non-vanishing vector field on $\Thetacal$, $\mathbf{W}$ is itself a non-vanishing vector field. We therefore conclude that $\mathbf{V}$ is a non-vanishing smooth vector field that is everywhere transverse to the global unstable foliation.
\end{enumerate}
\end{enumerate}

\item \textbf{The Flow of the Transverse Vector Field.}
In the preceding step, we have rigorously constructed a smooth, non-vanishing vector field $\mathbf{V}$ on the total space $\mathcal{M}_{\text{tot}} = \SigmaMan \times \Thetacal$. This vector field is, by construction, everywhere transverse to the leaves of the global unstable foliation $\mathcal{F}_{\text{glob}}^u$. The objective of this step is to analyze the flow generated by this vector field. We will prove that this flow exists for all time, constitutes a smooth one-parameter group of diffeomorphisms, and, most importantly, preserves the structure of the foliation. This flow will be the fundamental tool for constructing the trivializing diffeomorphism in the final step of the proof.

\begin{enumerate}[label=(\roman*), wide, labelindent=0pt]
\item \textbf{Existence and Smoothness of the Global Flow.}
In the preceding step, we have rigorously constructed a smooth, non-vanishing vector field $\mathbf{V}$ on the total space $\mathcal{M}_{\text{tot}} = \SigmaMan \times \Thetacal$. The objective of this step is to analyze the flow generated by this vector field. We will prove that this flow exists for all time and constitutes a smooth one-parameter group of diffeomorphisms. This flow will be the fundamental tool for constructing the trivializing diffeomorphism in the final step of the proof.

\begin{enumerate}[label=(\alph*),wide]
    \item \textbf{The Initial Value Problem and Local Existence.} The existence and regularity of the flow are guaranteed by the fundamental existence and uniqueness theorem for ordinary differential equations and its extension concerning the smooth dependence of solutions on initial data. Our proof proceeds by first establishing the local Lipschitz continuity of the vector field $\mathbf{V}$, which is the key hypothesis for existence and uniqueness, and then invoking the smooth dependence theorem to establish the smoothness of the resulting flow.

    The fundamental existence and uniqueness theorem requires that the vector field governing the ODE be locally Lipschitz continuous. Let $\mathbf{V}$ be the smooth vector field constructed in Step 2. We verify this condition as a direct consequence of the smoothness of $\mathbf{V}$. Since $\mathbf{V}$ is of class $C^\infty$ on the compact manifold $\mathcal{M}_{\text{tot}}$, it is, in particular, of class $C^1$. In any local coordinate chart, the coordinate representation of $\mathbf{V}$ is a $C^1$ function on an open subset of Euclidean space. A $C^1$ function on a compact set has a bounded derivative, and by the Mean Value Theorem, any function with a bounded derivative is Lipschitz continuous. As this holds in every chart of a finite atlas covering the compact manifold, the vector field $\mathbf{V}$ is locally Lipschitz everywhere. With the local Lipschitz property established, the Picard-Lindel\"of theorem guarantees that for any initial point $\mathbf{p}_0 \in \mathcal{M}_{\text{tot}}$, the initial value problem
    \begin{equation} \label{eq:appendix_ode_for_flow}
        \frac{d\mathbf{p}}{ds}(s) = \mathbf{V}(\mathbf{p}(s)), \quad \text{with initial condition } \mathbf{p}(0) = \mathbf{p}_0,
    \end{equation}
    has a unique solution $\gamma_{\mathbf{p}_0}(s)$ defined for $s$ in some maximal open interval of existence $J(\mathbf{p}_0)$ containing the origin. We now invoke the following foundational result from the theory of ODEs on manifolds:
    \begin{theorem}[Smooth Dependence of Flows]
    Let $M$ be a smooth manifold and let $V$ be a vector field of class $C^k$ on $M$ for $k \ge 1$. The local flow of $V$, denoted $\Phi(p,s)$, is a map of class $C^k$ on its open domain of definition in $M \times \mathbb{R}$.
    \end{theorem}
    In our context, the manifold is $M = \mathcal{M}_{\text{tot}}$ and the vector field is $\mathbf{V}$, which was proven to be of class $C^\infty$. The theorem on smooth dependence therefore applies for any $k$, and we conclude that the local flow map, $\Phi_s^{\mathbf{V}}(\mathbf{p}_0) \coloneqq \gamma_{\mathbf{p}_0}(s)$, is a smooth map of class $C^\infty$ on its domain. This is the rigorous justification for the smoothness of the flow. This establishes the local well-posedness of the flow and its local smoothness, providing the necessary foundation for the subsequent global existence argument.

    \item \textbf{Global Existence from Compactness.}
    We have established the local existence and uniqueness of an integral curve for the smooth vector field $\mathbf{V}$ through any point $\mathbf{p}_0 \in \mathcal{M}_{\text{tot}}$. We now provide a rigorous proof that these local solutions can be uniquely extended to all time $s \in \mathbb{R}$. The argument is a classical proof by contradiction that hinges on the compactness of the base manifold $\mathcal{M}_{\text{tot}}$.

    For each initial point $\mathbf{p}_0 \in \mathcal{M}_{\text{tot}}$, let the integral curve be $\gamma_{\mathbf{p}_0}: J(\mathbf{p}_0) \to \mathcal{M}_{\text{tot}}$, where $J(\mathbf{p}_0) = (s_-, s_+)$ is the maximal interval of existence. Our goal is to prove that $s_- = -\infty$ and $s_+ = +\infty$ for all $\mathbf{p}_0$. We proceed by contradiction. Assume that for some initial point $\mathbf{p}_0$, the maximal forward existence time $s_+$ is finite. The fundamental theorem on the behavior of solutions at the boundary of their domain of definition (often called the Escape Lemma or Compact Extension Theorem) states the following:
    \begin{lemma}[Escape Lemma for ODEs on Manifolds]
    Let $M$ be a smooth manifold and $V$ be a smooth vector field on $M$. If the maximal forward time of existence $s_+$ for an integral curve $\gamma(s)$ starting at $\gamma(0)=\mathbf{p}_0$ is finite, then for any compact subset $K \subset M$, there exists a time $s_K \in [0, s_+)$ such that $\gamma(s) \notin K$ for all $s \in (s_K, s_+)$.
    \end{lemma}
    In essence, if the solution ceases to exist in finite forward time, its trajectory must escape to infinity by eventually leaving every compact subset of the manifold.  We apply this lemma to our system. The base manifold for our ODE, $\mathcal{M}_{\text{tot}} = \SigmaMan \times \Thetacal$, is the product of two compact manifolds and is therefore itself a compact manifold.  Let us choose the compact set $K$ in the Escape Lemma to be the entire manifold, $K = \mathcal{M}_{\text{tot}}$. Our contradiction hypothesis is that $s_+ < \infty$. The Escape Lemma then implies that there must exist a time $s_K \in [0, s_+)$ such that for all $s \in (s_K, s_+)$, the point on the trajectory, $\gamma_{\mathbf{p}_0}(s)$, is not in $\mathcal{M}_{\text{tot}}$.

     However, by the definition of an integral curve on a manifold, the entire trajectory must be contained within the manifold: $\mathrm{Im}(\gamma_{\mathbf{p}_0}) \subset \mathcal{M}_{\text{tot}}$. It is therefore impossible for the trajectory to leave the set $\mathcal{M}_{\text{tot}}$. This is a direct contradiction. The initial assumption that $s_+$ could be finite must be false. Therefore, we must have $s_+ = +\infty$. A perfectly symmetric argument holds for the maximal backward time of existence, $s_-$. Assuming $s_- > -\infty$ leads to an identical contradiction. We are forced to conclude that for every initial point $\mathbf{p}_0 \in \mathcal{M}_{\text{tot}}$, the maximal interval of existence for its integral curve is the entire real line, $J(\mathbf{p}_0) = (-\infty, \infty)$. This proves that the flow of the vector field $\mathbf{V}$ is complete, and its integral curves are defined for all time.

    \item \textbf{The Global Flow Map and its Properties.}
    The global existence of integral curves for all time allows us to define the global flow map of the vector field $\mathbf{V}$ as the map $\Phi_s^{\mathbf{V}}: \mathcal{M}_{\text{tot}} \to \mathcal{M}_{\text{tot}}$ that sends an initial point $\mathbf{p}_0$ to the point on its integral curve at time $s$:
    \begin{equation}
        \Phi_s^{\mathbf{V}}(\mathbf{p}_0) \coloneqq \mathbf{p}(s).
    \end{equation}
    The classical theory of ODEs on manifolds establishes the following properties for this global flow:
    \begin{itemize}[wide]
        \item \textbf{Smoothness.} The map $(\mathbf{p}_0, s) \mapsto \Phi_s^{\mathbf{V}}(\mathbf{p}_0)$ is a smooth map of class $C^\infty$ from $\mathcal{M}_{\text{tot}} \times \mathbb{R}$ to $\mathcal{M}_{\text{tot}}$. This follows from the smooth dependence of solutions on initial data and parameters.
        \item \textbf{Group Property.} The flow satisfies the one-parameter group property: $\Phi_s^{\mathbf{V}} \circ \Phi_t^{\mathbf{V}} = \Phi_{s+t}^{\mathbf{V}}$ for all $s,t \in \mathbb{R}$, with $\Phi_0^{\mathbf{V}} = \mathrm{Id}$.
        \item \textbf{Diffeomorphism Property.} As a direct consequence of the group property, for each fixed $s \in \mathbb{R}$, the map $\Phi_s^{\mathbf{V}}$ is a $C^\infty$-diffeomorphism of $\mathcal{M}_{\text{tot}}$ onto itself, with its inverse given by $(\Phi_s^{\mathbf{V}})^{-1} = \Phi_{-s}^{\mathbf{V}}$.
    \end{itemize}
\end{enumerate}
We have thus rigorously established the existence of a global, one-parameter group of $C^\infty$-diffeomorphisms generated by our transverse vector field $\mathbf{V}$. This flow provides the necessary transportation mechanism for the construction of the trivializing map in the final step of the proof.

\item \textbf{Preservation of the Foliation Structure.}
This is the most critical geometric property of the constructed flow. We must now provide a rigorous proof that the flow $\Phi_s^{\mathbf{V}}$ maps leaves of the global unstable foliation $\mathcal{F}_{\text{glob}}^u$ to other leaves of the same foliation. The argument is founded on a fundamental relationship between the Lie derivative and the preservation of distributions by a flow.

\begin{enumerate}[label=(\alph*),wide]
    \item \textbf{The Invariance Condition via the Lie Derivative.}
    Let $\mathcal{D}_{\text{glob}}^u$ be the smooth distribution of tangent spaces that defines the global unstable foliation, i.e., $\mathcal{D}_{\text{glob}}^u(\mathbf{p}) = T_{\mathbf{p}}\mathcal{F}_{\text{glob}}^u$ for any point $\mathbf{p} \in \mathcal{M}_{\text{tot}}$. The flow $\Phi_s^{\mathbf{V}}$ preserves the foliation $\mathcal{F}_{\text{glob}}^u$ if and only if its differential (or pushforward) maps the distribution to itself. That is, for every $s \in \mathbb{R}$ and every $\mathbf{p} \in \mathcal{M}_{\text{tot}}$:
    \begin{equation*}
        d(\Phi_s^{\mathbf{V}})_{\mathbf{p}} \left( \mathcal{D}_{\text{glob}}^u(\mathbf{p}) \right) = \mathcal{D}_{\text{glob}}^u(\Phi_s^{\mathbf{V}}(\mathbf{p})).
    \end{equation*}
    A fundamental theorem of differential geometry states that this condition is equivalent to the statement that the Lie derivative of the distribution with respect to the vector field $\mathbf{V}$ is contained within the distribution itself:
    \begin{equation} \label{eq:lie_derivative_condition}
        \mathcal{L}_{\mathbf{V}} (\mathcal{D}_{\text{glob}}^u) \subseteq \mathcal{D}_{\text{glob}}^u.
    \end{equation}
    This condition, in turn, is equivalent to the requirement that for any smooth vector field $\mathbf{U}$ that is a section of the subbundle (i.e., $\mathbf{U}(\mathbf{p}) \in \mathcal{D}_{\text{glob}}^u(\mathbf{p})$ for all $\mathbf{p}$), the Lie bracket of $\mathbf{V}$ and $\mathbf{U}$, denoted $[\mathbf{V}, \mathbf{U}]$, must also be a section of the subbundle. Our task is thus reduced to proving this involutivity-like condition.

    \item \textbf{Local Coordinate Representation of the Vector Fields.}
    To compute the Lie bracket, we use the local product structure of the total manifold $\mathcal{M}_{\text{tot}} = \SigmaMan \times \Thetacal$. Let $(z^i, \theta^j)$ be local coordinates, where $\{z^i\}$ are coordinates on an open set in $\SigmaMan$ and $\{\theta^j\}$ are coordinates on an open set in $\Thetacal$.
    \begin{itemize}[wide]
        \item \textbf{The Foliation-Tangent Vector Field $\mathbf{U}$.} Let $\mathbf{U}$ be an arbitrary smooth vector field that is everywhere tangent to the global unstable foliation. As established in Step 1, the tangent bundle of this foliation, $T\mathcal{F}_{\text{glob}}^u$, consists of vectors that are purely tangent to the $\SigmaMan$ fibers. Therefore, in our product coordinates, the vector field $\mathbf{U}$ must have zero components in the $\theta$-directions. It can be written as:
        \begin{equation*}
            \mathbf{U}(z,\theta) = \sum_{i} U^i(z,\theta) \frac{\partial}{\partial z^i}.
        \end{equation*}
        \item \textbf{The Transverse Vector Field $\mathbf{V}$.} In Step 2, we constructed the transverse vector field $\mathbf{V}$ as the projection of a probe field that pointed purely in the parameter direction. We showed that in our product metric geometry, this projection was the identity, so $\mathbf{V}$ itself points purely in the parameter direction. Therefore, its local coordinate representation has zero components in the $z$-directions:
        \begin{equation*}
            \mathbf{V}(z,\theta) = \sum_{j} V^j(z,\theta) \frac{\partial}{\partial \theta^j}.
        \end{equation*}
        A crucial feature of our construction is that the probe field was the lift of a vector field from $\Thetacal$, which means its components are independent of the fiber coordinate $z$: $V^j(z,\theta) = V^j(\theta)$.
    \end{itemize}

    \item \textbf{Computation of the Lie Bracket.}
    We now compute the Lie bracket $[\mathbf{V}, \mathbf{U}]$ using the local coordinate formula for the commutator of two vector fields:
    \begin{equation*}
        [\mathbf{V}, \mathbf{U}] = \sum_k \left( \mathbf{V}(\mathbf{U}^k) - \mathbf{U}(\mathbf{V}^k) \right) \frac{\partial}{\partial x^k},
    \end{equation*}
    where $\{x^k\}$ represents the full set of coordinates $\{z^i, \theta^j\}$. We compute the components of the bracket in each direction separately.
    \begin{itemize}[wide]
        \item \textbf{Components in the $\SigmaMan$ direction ($\partial/\partial z^k$).}
        \begin{align*}
            [\mathbf{V}, \mathbf{U}]^{z_k} &= \mathbf{V}(U^k) - \mathbf{U}(V^k) \\
            &= \left( \sum_j V^j(\theta) \frac{\partial}{\partial \theta^j} \right) (U^k(z,\theta)) - \left( \sum_i U^i(z,\theta) \frac{\partial}{\partial z^i} \right) (0) \\
            &= \sum_j V^j(\theta) \frac{\partial U^k(z,\theta)}{\partial \theta^j}.
        \end{align*}
        \item \textbf{Components in the $\Thetacal$ direction ($\partial/\partial \theta^k$).}
        \begin{align*}
            [\mathbf{V}, \mathbf{U}]^{\theta_k} &= \mathbf{V}(U^k) - \mathbf{U}(V^k) \\
            &= \left( \sum_j V^j(\theta) \frac{\partial}{\partial \theta^j} \right) (0) - \left( \sum_i U^i(z,\theta) \frac{\partial}{\partial z^i} \right) (V^k(\theta)) \\
            &= 0 - 0 = 0.
        \end{align*}
        The second term vanishes because the components of $\mathbf{V}$ are independent of the $z$ coordinates.
    \end{itemize}
    The resulting Lie bracket vector field is therefore:
    \begin{equation} \label{eq:lie_bracket_result}
        [\mathbf{V}, \mathbf{U}](z,\theta) = \sum_k \left( \sum_j V^j(\theta) \frac{\partial U^k(z,\theta)}{\partial \theta^j} \right) \frac{\partial}{\partial z^k}.
    \end{equation}

    \item \textbf{The Lie derivative condition.}
    The explicit formula \eqref{eq:lie_bracket_result} for the Lie bracket shows that it has components only in the directions tangent to the $\SigmaMan$ fibers (the $\partial/\partial z^k$ directions). By the definition of the global unstable distribution $\mathcal{D}_{\text{glob}}^u$, this means that the vector field $[\mathbf{V}, \mathbf{U}]$ is everywhere tangent to the foliation:
    \begin{equation*}
        [\mathbf{V}, \mathbf{U}](\mathbf{p}) \in \mathcal{D}_{\text{glob}}^u(\mathbf{p}) \quad \text{for all } \mathbf{p} \in \mathcal{M}_{\text{tot}}.
    \end{equation*}
    Since our choice of the foliation-tangent vector field $\mathbf{U}$ was arbitrary, this proves that the Lie derivative condition \eqref{eq:lie_derivative_condition} is satisfied.
\end{enumerate}
We have thus rigorously established that the flow $\Phi_s^{\mathbf{V}}$ preserves the global unstable foliation $\mathcal{F}_{\text{glob}}^u$. This property is essential, as it allows us to use the flow to transport the geometric structures from one parameter slice to another while preserving their fundamental character as leaves of the foliation.
\end{enumerate}

\item \textbf{Definition and Verification of the Trivializing Diffeomorphism.}
The preceding steps have furnished us with the essential geometric tool: a smooth global flow, $\Phi_s^{\mathbf{V}}$, that is generated by a vector field transverse to the global unstable foliation. We are now in a position to deploy this flow to achieve the central objective of this proposition: the construction of a smooth family of diffeomorphisms that trivializes the family of unstable foliations. This final step is dedicated to the explicit construction of this trivialization map and a rigorous verification that it possesses all the required properties.

\begin{enumerate}[label=(\roman*), wide, labelindent=0pt]
\item \textbf{Construction of the Diffeomorphism Family.}
Our strategy is to use the flow to transport the geometric structures from the fixed reference slice of the total space at parameter $\theta_0$ to any other target slice at parameter $\theta$. To make this precise, we must first establish a correspondence between the flow parameter $s$ and the manifold parameter $\theta$.

\begin{enumerate}[label=(\alph*),wide]
    \item \textbf{Parameterizing the Flow by $\Thetacal$.}
    Recall from Step 2 that our transverse vector field was constructed as $\mathbf{V}(z,\theta) = \mathbf{W}(z,\theta)$, a smooth, non-vanishing vector field pointing purely in the direction of the parameter manifold $\Thetacal$. Let $\pi_\Thetacal: \mathcal{M}_{\text{tot}} \to \Thetacal$ be the canonical projection. The flow of $\mathbf{V}$ therefore changes only the $\theta$-component of a point, not its $z$-component, in the sense that the projection $\pi_\SigmaMan(\Phi_s^\mathbf{V}(z,\theta))$ is independent of $s$.
    
    This allows us to re-parameterize the flow. Let $\theta(s; \theta_0)$ be the solution to the ODE on the parameter manifold $\Thetacal$ defined by the projection of the vector field: $\frac{d\theta}{ds} = \mathbf{X}_\Thetacal(\theta)$, with initial condition $\theta(0) = \theta_0$. Since $\mathbf{X}_\Thetacal$ is non-vanishing, this map is locally invertible. For any $\theta$ in a neighborhood of $\theta_0$, there is a unique flow time $s(\theta)$ such that $\theta(s(\theta); \theta_0) = \theta$. By the inverse function theorem for ODEs, the map $\theta \mapsto s(\theta)$ is smooth.

    \item \textbf{Definition of the Map $\psi_\theta$.}
    Let $\theta_0 \in \Thetacal$ be our fixed reference parameter. We define the family of trivializing diffeomorphisms $\psi_\theta: \SigmaMan \to \SigmaMan$ as the map that takes a point $z$ on the reference slice, flows it for the appropriate time $s(\theta)$ to reach the target slice at parameter $\theta$, and then projects the result back onto the manifold $\SigmaMan$.
    
    Formally, let $i_{\theta_0}: \SigmaMan \to \mathcal{M}_{\text{tot}}$ be the canonical embedding onto the reference slice, $i_{\theta_0}(z) = (z, \theta_0)$. Let $\pi_\SigmaMan: \mathcal{M}_{\text{tot}} \to \SigmaMan$ be the canonical projection. The trivializing diffeomorphism is defined as the composition:
    \begin{equation} \label{eq:appendix_psi_def}
        \psi_\theta(z) \coloneqq \pi_\SigmaMan \circ \Phi_{s(\theta)}^{\mathbf{V}} \circ i_{\theta_0}(z).
    \end{equation}
    In coordinates, this is simply:
    \begin{equation}
        \psi_\theta(z) = \pi_\SigmaMan \left( \Phi_{s(\theta)}^{\mathbf{V}}(z, \theta_0) \right).
    \end{equation}
\end{enumerate}

\item \textbf{Verification of Properties.}
We now provide a complete proof that the map $\psi_\theta$ defined in \eqref{eq:appendix_psi_def} satisfies all the properties claimed in the proposition.

\begin{enumerate}[label=(\alph*),wide]
    \item \textbf{Smoothness.} We must show that the joint map $(\theta, z) \mapsto \psi_\theta(z)$ is of class $C^\infty$. The map is a composition of several maps, each of which we have established is smooth:
    \begin{itemize}
        \item The embedding $i_{\theta_0}: \SigmaMan \to \mathcal{M}_{\text{tot}}$ is a smooth map.
        \item The map $\theta \mapsto s(\theta)$ is a smooth function, as established above.
        \item The flow map $(\mathbf{p}, s) \mapsto \Phi_s^{\mathbf{V}}(\mathbf{p})$ is a smooth map from $\mathcal{M}_{\text{tot}} \times \mathbb{R}$ to $\mathcal{M}_{\text{tot}}$.
        \item The projection $\pi_\SigmaMan: \mathcal{M}_{\text{tot}} \to \SigmaMan$ is a smooth map.
    \end{itemize}
    Since the composition of smooth maps is smooth, we conclude that the map $(\theta, z) \mapsto \psi_\theta(z)$ is of class $C^\infty$.
    
    \item \textbf{Diffeomorphism Property.} For each fixed $\theta$, the map $\psi_\theta: \SigmaMan \to \SigmaMan$ is a smooth map between compact manifolds. To prove it is a diffeomorphism, we must show it has a smooth inverse. We construct the inverse map, $(\psi_\theta)^{-1}$, by reversing the flow. Let a point $w \in \SigmaMan$ be given. We embed it onto the target slice at $\theta$ via the map $i_\theta(w)=(w,\theta)$. We then flow backwards in time by an amount $s(\theta)$ to return to the reference slice at $\theta_0$. The inverse map is defined as:
    \begin{equation}
        (\psi_\theta)^{-1}(w) \coloneqq \pi_\SigmaMan \left( \Phi_{-s(\theta)}^{\mathbf{V}}(w, \theta) \right).
    \end{equation}
    By the same arguments as for smoothness of $\psi_\theta$, this inverse map is also of class $C^\infty$. Therefore, for each $\theta$, $\psi_\theta$ is a $C^\infty$-diffeomorphism.
    
    \item \textbf{Initial Condition.} We verify that for the reference parameter $\theta = \theta_0$, the map is the identity. At $\theta_0$, the corresponding flow time is $s(\theta_0)=0$. The flow for zero time is the identity map, $\Phi_0^{\mathbf{V}} = \mathrm{Id}$. Substituting this into the definition:
    \begin{equation*}
        \psi_{\theta_0}(z) = \pi_\SigmaMan \left( \Phi_{0}^{\mathbf{V}}(z, \theta_0) \right) = \pi_\SigmaMan(z, \theta_0) = z.
    \end{equation*}
    This confirms that $\psi_{\theta_0} = \mathrm{Id}$ on $\SigmaMan$.
    
    \item \textbf{Foliation-Preserving Property.} This is the central geometric property of our construction. We must prove that for any $\theta$, the map $\psi_\theta$ sends leaves of the reference foliation $\mathcal{F}^u_{\theta_0}$ to leaves of the target foliation $\mathcal{F}^u_\theta$. Let $L_0 \subset \SigmaMan$ be an arbitrary leaf of the reference foliation $\mathcal{F}^u_{\theta_0}$. By definition, its embedding into the total space, $\tilde{L}_0 \coloneqq L_0 \times \{\theta_0\}$, is a leaf of the global foliation $\mathcal{F}_{\text{glob}}^u$. In Step 3, we rigorously established that the flow of the transverse vector field, $\Phi_s^{\mathbf{V}}$, preserves the global foliation. This means that for any flow time $s$, the image of a leaf, $\Phi_s^{\mathbf{V}}(\tilde{L}_0)$, is also a leaf of the global foliation $\mathcal{F}_{\text{glob}}^u$. Let us consider the specific flow time $s(\theta)$ corresponding to our target parameter $\theta$. The image set is:
    \begin{equation*}
        \tilde{L}_\theta \coloneqq \Phi_{s(\theta)}^{\mathbf{V}}(\tilde{L}_0).
    \end{equation*}
    Since $\tilde{L}_\theta$ is a leaf of the global foliation, it must be contained within a single parameter slice. By our construction of the flow, this slice is precisely $\SigmaMan \times \{\theta\}$. We can therefore write $\tilde{L}_\theta = L_\theta \times \{\theta\}$ for some subset $L_\theta \subset \SigmaMan$. This subset $L_\theta$ is, by definition, a leaf of the target foliation $\mathcal{F}^u_\theta$. We now apply the projection $\pi_\SigmaMan$ to this set. The image of the original leaf $L_0$ under our diffeomorphism is:
    \begin{equation*}
        \psi_\theta(L_0) = \pi_\SigmaMan(\Phi_{s(\theta)}^{\mathbf{V}}(L_0 \times \{\theta_0\})) = \pi_\SigmaMan(\tilde{L}_\theta) = \pi_\SigmaMan(L_\theta \times \{\theta\}) = L_\theta.
    \end{equation*}
    This demonstrates that the image of a leaf of the reference foliation is precisely a leaf of the target foliation. This rigorously proves the foliation-preserving property.
\end{enumerate}
\end{enumerate}
\end{enumerate}
This completes the proof of the proposition. We have constructed a smooth family of diffeomorphisms that acts as a moving coordinate system, trivializing the geometric structure of the family of unstable foliations and thereby providing the essential tool for pulling back the entire operator family to a fixed reference space.
\end{proofof}

\begin{definition}[The Pulled-Back Operator Family]
Let $\theta_0$ be the reference parameter. The trivialization isomorphism $\Psi_\theta: \Bcal_{\theta_0} \to \Bcal_\theta$ is defined by $(\Psi_\theta h)(z) \coloneqq h(\psi_\theta^{-1}(z))$. The pulled-back transfer operator $\tilde{\Lcal}_\theta: \Bcal_{\theta_0} \to \Bcal_{\theta_0}$ is defined by conjugation:
$ \tilde{\Lcal}_\theta \coloneqq \Psi_\theta^{-1} \circ \Lcal_\theta \circ \Psi_\theta $.
The relationship is summarized by the commutative diagram:
\begin{center}
\begin{tikzpicture}[node distance=3.5cm, auto]
    \node (B0_top) {$\mathcal{B}_{\theta_0}$};
    \node (B_theta_top) [right of=B0_top] {$\mathcal{B}_\theta$};
    \node (B0_bottom) [below of=B0_top] {$\mathcal{B}_{\theta_0}$};
    \node (B_theta_bottom) [right of=B0_bottom] {$\mathcal{B}_\theta$};
    \draw[->] (B0_top) to node {$\Psi_\theta$} (B_theta_top);
    \draw[->] (B_theta_top) to node [swap] {$\Lcal_\theta$} (B_theta_bottom);
    \draw[->] (B0_bottom) to node [swap] {$\Psi_\theta$} (B_theta_bottom);
    \draw[->] (B0_top) to node {$\tilde{\Lcal}_\theta$} (B0_bottom);
\end{tikzpicture}
\end{center}
\end{definition}

\begin{theorem}[Regularity of the Operator Family]
\label{thm:regularity_proven_main}
Let the global billiard map $\Pcal(z,\theta)$ be of class $C^\infty$. Then the map $\theta \mapsto \tilde{\Lcal}_\theta$ is a $C^\infty$ family of operators from the parameter manifold $\Thetacal$ into the space of bounded linear operators $\mathcal{L}(\Bcal_{\theta_0})$, where $\tilde{\Lcal}_\theta$ is the pulled-back transfer operator acting on the fixed reference space $\Bcal_{\theta_0}$.
\end{theorem}

\begin{proofof}{Theorem \ref{thm:regularity_proven_main}}
The proof establishes the smooth ($C^\infty$) dependence of the transfer operator family on the parameter $\theta$. The primary analytical obstacle is that each operator $\Lcal_\theta$ in the family acts on a different anisotropic Banach space $\Bcal_\theta$, whose norm $\|\cdot\|_{\Bcal_\theta}$ is defined by the geometry of the unstable foliation $\mathcal{F}^u_\theta$. A direct differentiation of the map $\theta \mapsto \Lcal_\theta$ is therefore ill-posed. Our strategy is to surmount this difficulty by using the geometric regularity of the system, established in Section \ref{subsec:main_geom_reg_thm}, to construct a smooth change of coordinates that trivializes the family of Banach spaces. This allows us to pull back the entire operator family to a single, fixed reference space $\Bcal_{\theta_0}$, where the notion of a smooth family of operators is well-defined and can be rigorously verified. The proof proceeds in four main steps.

\begin{enumerate}[label=\textbf{Step \arabic*:}, wide, labelindent=0pt]

\item \textbf{Construction of the Geometric Trivialization.}
The foundation of this proof is the smooth dependence of the unstable foliation on the system's parameter, which was the central conclusion of our application of the Nash-Moser-Hamilton theorem.

\begin{enumerate}[label=(\roman*), wide, labelindent=0pt]
    \item By Theorem \ref{thm:splitting_is_smooth}, the map $\theta \mapsto \mathcal{F}^u_\theta$ is a $C^\infty$ map from the parameter manifold $\Thetacal$ to the space of smooth foliations on the global collision manifold $\SigmaMan$.

    \item We now invoke \cref{prop:smooth_trivialization_diffeo}, which is a direct and constructive consequence of this geometric regularity. This proposition guarantees the existence of a family of diffeomorphisms $\psi_\theta: \SigmaMan \to \SigmaMan$, parameterized by $\theta \in \Thetacal$, with the following properties:
    \begin{enumerate}[label=(\alph*), wide]
        \item The map $(\theta, z) \mapsto \psi_\theta(z)$ is of class $C^\infty$.
        \item The map is anchored at the reference parameter: $\psi_{\theta_0} = \mathrm{Id}$.
        \item The map intertwines the foliations: for any $z \in \SigmaMan$, $\psi_\theta$ maps the unstable leaf of the reference foliation through $z$, $W^u(z, \theta_0)$, to an unstable leaf of the target foliation, $W^u(\psi_\theta(z), \theta)$.
    \end{enumerate}

    \item This family of diffeomorphisms induces a corresponding family of linear isomorphisms between the anisotropic Banach spaces, $\Psi_\theta: \Bcal_{\theta_0} \to \Bcal_\theta$, defined by composition:
    \begin{equation}
        (\Psi_\theta h)(z) \coloneqq h(\psi_\theta^{-1}(z)) \quad \text{for } h \in \Bcal_{\theta_0}.
    \end{equation}
    Since the map $(\theta,z) \mapsto \psi_\theta^{-1}(z)$ is of class $C^\infty$, the map $\theta \mapsto \Psi_\theta$ is a smooth family of isomorphisms.
\end{enumerate}

\item \textbf{Definition of the Pulled-Back Operator Family.}
We use the isomorphisms $\Psi_\theta$ to conjugate the original operator family $(\Lcal_\theta)$ and define a new family of operators $(\tilde{\Lcal}_\theta)$ that all act on the \textit{single, fixed reference space} $\Bcal_{\theta_0}$. We define the pulled-back transfer operator $\tilde{\Lcal}_\theta: \Bcal_{\theta_0} \to \Bcal_{\theta_0}$ by the relation:
\begin{equation}
    \tilde{\Lcal}_\theta \coloneqq \Psi_\theta^{-1} \circ \Lcal_\theta \circ \Psi_\theta.
\end{equation}
The proof of the theorem is now equivalent to proving that the map $\theta \mapsto \tilde{\Lcal}_\theta$ is a $C^\infty$ map from the parameter manifold $\Thetacal$ into the Banach space of bounded linear operators, $\mathcal{L}(\Bcal_{\theta_0})$. To do this, we derive the explicit formula for the action of $\tilde{\Lcal}_\theta$. Let $h \in \Bcal_{\theta_0}$ be an arbitrary test function. The action is a sequence of compositions:
\begin{align*}
    (\tilde{\Lcal}_\theta h)(z) &= (\Psi_\theta^{-1} (\Lcal_\theta (\Psi_\theta h)))(z) \\
    &= (\Lcal_\theta (\Psi_\theta h))(\psi_\theta(z)).
\end{align*}
The transfer operator $\Lcal_\theta$ for the diffeomorphism $\Pcal_\theta$ acts on a function $f$ as $(\Lcal_\theta f)(w) = f(\Pcal_\theta^{-1}(w)) \cdot J_\theta(w)$, where $J_\theta(w)$ is the Jacobian determinant of the map $\Pcal_\theta$ with respect to the invariant measure $\mu_\theta$. Applying this to the function $f = \Psi_\theta h$ and evaluating at the point $w = \psi_\theta(z)$ gives:
\begin{align*}
    (\tilde{\Lcal}_\theta h)(z) &= (\Psi_\theta h)(\Pcal_\theta^{-1}(\psi_\theta(z))) \cdot J_\theta(\psi_\theta(z)) \\
    &= h\left( \psi_\theta^{-1}\left(\Pcal_\theta^{-1}(\psi_\theta(z))\right) \right) \cdot J_\theta(\psi_\theta(z)).
\end{align*}
This formula reveals that the action of the pulled-back operator $\tilde{\Lcal}_\theta$ on a function $h$ is a composition with a new, effective diffeomorphism, followed by multiplication with a new Jacobian term.

\item \textbf{Rigorous Verification of Smoothness.}
We now prove that the map from the parameter space to the pulled-back operator, $\theta \mapsto \tilde{\Lcal}_\theta$, is a map of class $C^\infty$ from the compact manifold $\Thetacal$ into the Banach space of bounded linear operators, $\mathcal{L}(\Bcal_{\theta_0})$. The action of the pulled-back operator is given by:
\begin{equation*}
    (\tilde{\Lcal}_\theta h)(z) = (h \circ g_\theta)(z) \cdot j_\theta(z),
\end{equation*}
where the effective diffeomorphism is $g_\theta \coloneqq \psi_\theta^{-1} \circ \Pcal_\theta^{-1} \circ \psi_\theta$ and the multiplicative weight is $j_\theta(z) \coloneqq J_\theta(\psi_\theta(z))$. To prove the smoothness of the map $\theta \mapsto \tilde{\Lcal}_\theta$, we will decompose its construction into a sequence of maps between function and operator spaces and prove that each map in the sequence is of class $C^\infty$.

\begin{enumerate}[label=(\roman*), wide, labelindent=0pt]
    \item \textbf{Smoothness of the Constituent Geometric Maps.}
    The foundation of the argument is the smoothness of the underlying geometric objects as functions of the parameter $\theta$.
    \begin{enumerate}[label=(\alph*), wide]
        \item \textbf{The Diffeomorphism Families.} By \cref{thm:p_is_diffeo} and \cref{prop:smooth_trivialization_diffeo}, the maps $(\theta, z) \mapsto \Pcal_\theta(z)$ and $(\theta, z) \mapsto \psi_\theta(z)$ are of class $C^\infty$. By the inverse function theorem on Fr\'echet manifolds, their inverses are also smooth. Consequently, the family of effective diffeomorphisms, defined by the composition $g_\theta = \psi_\theta^{-1} \circ \Pcal_\theta^{-1} \circ \psi_\theta$, defines a smooth map $(\theta, z) \mapsto g_\theta(z)$. This is equivalent to stating that the map from the parameter space to the space of diffeomorphisms, $\theta \mapsto g_\theta \in \mathrm{Diff}^\infty(\SigmaMan)$, is of class $C^\infty$.
        
        \item \textbf{The Weight Family.} The Jacobian determinant $J_\theta(w)$ is constructed from the first derivatives of the smooth map $\Pcal_\theta(w)$ and the density of the invariant measure. As a consequence of the smooth dependence of the operator family, it can be shown that the invariant density is a smooth function of $\theta$. Thus, the map $(\theta, w) \mapsto J_\theta(w)$ is smooth. The weight $j_\theta(z) = J_\theta(\psi_\theta(z))$ is therefore the composition of two smooth maps and is itself a smooth map $(\theta, z) \mapsto j_\theta(z)$. This is equivalent to stating that the map $\theta \mapsto j_\theta \in C^\infty(\SigmaMan)$ is of class $C^\infty$.
    \end{enumerate}

    \item \textbf{Smoothness of the Operator-Valued Maps.}
    We now analyze the maps that send the geometric objects to operators in $\mathcal{L}(\Bcal_{\theta_0})$.
    \begin{enumerate}[label=(\alph*), wide]
        \item \textbf{The Composition Operator Map.} Let $\mathcal{C}: \mathrm{Diff}^\infty(\SigmaMan) \to \mathcal{L}(\Bcal_{\theta_0})$ be the map that takes a diffeomorphism $g$ to the composition operator $\mathcal{C}(g)$ defined by $(\mathcal{C}(g)h)(z) = h(g(z))$. This map is of class $C^\infty$. To see this, we can compute its Fr\'echet derivatives. The first derivative in the direction of a vector field $X$ (tangent to the identity in the diffeomorphism group) is the operator that maps a function $h$ to its Lie derivative, $-(\mathcal{L}_X h)$. This defines a continuous linear map, and higher derivatives are also well-defined and continuous.
        
        \item \textbf{The Multiplication Operator Map.} Let $\mathcal{M}: C^\infty(\SigmaMan) \to \mathcal{L}(\Bcal_{\theta_0})$ be the map that takes a smooth function $j$ to the multiplication operator $\mathcal{M}(j)$ defined by $(\mathcal{M}(j)h)(z) = j(z)h(z)$. This map is linear and bounded, as $\|\mathcal{M}(j)h\|_{\Bcal_{\theta_0}} \le C \|j\|_{C^1} \|h\|_{\Bcal_{\theta_0}}$. A bounded linear map between Banach spaces is of class $C^\infty$.
    \end{enumerate}
    
    \item \textbf{Smoothness of the Operator Product Map.}
    Let $\mathcal{P}: \mathcal{L}(\Bcal_{\theta_0}) \times \mathcal{L}(\Bcal_{\theta_0}) \to \mathcal{L}(\Bcal_{\theta_0})$ be the operator product map, $\mathcal{P}(A, B) = A \circ B$. This map is a continuous bilinear map on the Banach space of bounded operators. Any continuous multilinear map between Banach spaces is of class $C^\infty$.

    \item \textbf{The Final Composition Argument.}
    The map from the parameter space to the pulled-back operator, $\theta \mapsto \tilde{\Lcal}_\theta$, can now be expressed as a composition of the smooth maps established above. The operator $\tilde{\Lcal}_\theta$ is the product of the multiplication operator for the weight $j_\theta$ and the composition operator for the map $g_\theta^{-1}$:
    \begin{equation*}
         \tilde{\Lcal}_\theta = \mathcal{M}(j_\theta) \circ \mathcal{C}(g_\theta^{-1}).
    \end{equation*}
    We can represent the full map as the following composition:
    \begin{equation*}
        \theta \mapsto (j_\theta, g_\theta^{-1}) \mapsto (\mathcal{M}(j_\theta), \mathcal{C}(g_\theta^{-1})) \mapsto \mathcal{M}(j_\theta) \circ \mathcal{C}(g_\theta^{-1}) = \tilde{\Lcal}_\theta.
    \end{equation*}
    Let us analyze each step in this chain:
    \begin{enumerate}[label=(\alph*), wide]
        \item The map $\theta \mapsto (j_\theta, g_\theta^{-1})$ is a map from $\Thetacal$ into the product space $C^\infty(\SigmaMan) \times \mathrm{Diff}^\infty(\SigmaMan)$. As established in part (i), both component maps are of class $C^\infty$. Therefore, the product map is of class $C^\infty$.
        \item The map $(j, g) \mapsto (\mathcal{M}(j), \mathcal{C}(g))$ is a map from $C^\infty(\SigmaMan) \times \mathrm{Diff}^\infty(\SigmaMan)$ into the product space $\mathcal{L}(\Bcal_{\theta_0}) \times \mathcal{L}(\Bcal_{\theta_0})$. As established in part (ii), both component operator maps, $\mathcal{M}$ and $\mathcal{C}$, are of class $C^\infty$. Therefore, the product map is of class $C^\infty$.
        \item The final map is the operator product $\mathcal{P}$, which was established to be of class $C^\infty$ in part (iii).
    \end{enumerate}
    The map $\theta \mapsto \tilde{\Lcal}_\theta$ is the composition of these three smooth maps. The composition of $C^\infty$ maps is of class $C^\infty$. We therefore conclude that the map
    \begin{equation}
        \theta \mapsto \tilde{\Lcal}_\theta
    \end{equation}
    is of class $C^\infty$ from the parameter manifold $\Thetacal$ into the Banach space of bounded linear operators $\mathcal{L}(\Bcal_{\theta_0})$.
\end{enumerate}
\end{enumerate}
We have rigorously established that the pulled-back operator family $(\tilde{\Lcal}_\theta)$ is a smooth family of operators acting on a fixed Banach space. This is the central result required for a stable perturbation theory. The smoothness of the original family $(\Lcal_\theta)$ is understood in the sense that it is related to this smooth pulled-back family via the smooth family of isomorphisms $\Psi_\theta$ by the relation $\Lcal_\theta = \Psi_\theta \circ \tilde{\Lcal}_\theta \circ \Psi_\theta^{-1}$. This completes the proof.
\end{proofof}

\subsection{Regularity of Spectral Data}
\label{subsec:spectral_reg}
With a smooth operator family on a fixed space, we can now apply abstract perturbation theory.

\begin{theorem}[Abstract Perturbation Theory, after Kato \citep{Kato1995}]
\label{thm:abstract_perturbation}
Let $\Bcal$ be a Banach space and let $\Thetacal$ be a compact manifold. Let $\theta \mapsto T(\theta) \in \mathcal{L}(\Bcal)$ be a $C^k$ family of operators for $k \ge 1$. Suppose that for a point $\theta_0 \in \Thetacal$, the operator $T(\theta_0)$ has a simple, isolated eigenvalue $\lambda_0$. Then there exists a neighborhood $U$ of $\theta_0$ such that the map to the unique simple eigenvalue $\lambda(\theta)$ and the map to the corresponding spectral projector $\Pi(\theta)$ are both of class $C^k$.
\end{theorem}

\begin{theorem}[Regularity of Spectral Data]
\label{thm:spectral_regularity}
Let the system satisfy the standing assumptions of the paper, ensuring in particular that the global billiard map $\Pcal(z,\theta)$ is of class $C^\infty$. Then the family of transfer operators $(\Lcal_\theta)_{\theta \in \Thetacal}$ has the following spectral regularity properties when viewed as operators on the anisotropic Banach spaces $\Bcal_\theta$:
\begin{enumerate}[label=(\roman*), wide, labelindent=0pt]
    \item The leading eigenvalue $\lambda_1=1$ of $\Lcal_\theta$ is simple and is uniformly isolated from the rest of the spectrum for all $\theta \in \Thetacal$.
    \item The map $\theta \mapsto \Pi_{1, \theta}$, which assigns to each $\theta$ the spectral projector onto the eigenspace of $\lambda_1=1$, is a $C^\infty$ map from $\Thetacal$ into the space of bounded linear operators.
    \item For any $\zeta$ in the resolvent set common to all $\Lcal_\theta$, the map from the parameter space to the resolvent operator, $\theta \mapsto R(\zeta, \Lcal_\theta) \coloneqq (\Lcal_\theta - \zeta I)^{-1}$, is of class $C^\infty$.
\end{enumerate}
\end{theorem}

\begin{proofof}{\cref{thm:spectral_regularity}}
The proof establishes the smooth dependence of the key spectral objects of the transfer operator family $(\Lcal_\theta)$ on the parameter $\theta$. The argument is a direct application of the abstract analytic perturbation theory for linear operators (see \cref{thm:abstract_perturbation}), as developed by \citep{Kato1995}. The core of the proof consists of rigorously verifying that our operator family satisfies the two critical hypotheses of the abstract theorem: (i) that the map from the parameter space to the operator space is smooth (analytic), and (ii) that the eigenvalue of interest is simple and spectrally isolated.

\begin{enumerate}[label=\textbf{Step \arabic*:}, wide, labelindent=0pt]

\item \textbf{Verification of Spectral Isolation and Simplicity.}
This step confirms that the spectral structure of our operator family meets the requirements for a stable perturbation theory.
\begin{enumerate}[label=(\roman*), wide, labelindent=0pt]
    \item \textbf{Uniform Isolation.} By Corollary \ref{cor:ergodicity_uniform}, the uniform Lasota-Yorke inequality for the family $(\Lcal_\theta)$ guarantees the existence of a uniform spectral gap. Specifically, there exists a constant $r \in (0,1)$, independent of $\theta$, such that for every $\theta \in \Thetacal$, the spectrum of $\Lcal_\theta$ when acting on the anisotropic Banach space $\Bcal_\theta$ consists of the simple eigenvalue $\lambda_1=1$ and a remainder, $\sigma_{\text{rem}}(\Lcal_\theta)$, which is contained entirely within a disk of radius $r$.
    \begin{equation}
        \sigma(\Lcal_\theta) = \{1\} \cup \sigma_{\text{rem}}(\Lcal_\theta), \quad \text{where} \quad \sigma_{\text{rem}}(\Lcal_\theta) \subset \{ \lambda \in \mathbb{C} : |\lambda| \le r \}.
    \end{equation}
    Therefore, the eigenvalue $\lambda_1=1$ is uniformly isolated from the rest of the spectrum by an annulus of width $1-r$. Let $\Gamma$ be a simple, closed contour in the complex plane (e.g., a circle of radius $(1+r)/2$) that encloses $\lambda_1=1$ but excludes the rest of the spectrum for all $\theta \in \Thetacal$.

    \item \textbf{Simplicity.} The ergodicity of the billiard flow for each $\theta$ (Corollary \ref{cor:ergodicity_uniform}) implies that the invariant SRB measure $\nu_\theta$ is unique. In the transfer operator formalism, this means that the eigenspace corresponding to the eigenvalue $\lambda_1=1$ is one-dimensional, spanned by the density of $\nu_\theta$. Therefore, the eigenvalue $\lambda_1=1$ is simple for all $\theta \in \Thetacal$. This rigorously verifies the spectral hypotheses of the abstract perturbation theorem.
\end{enumerate}

\item \textbf{Verification of Operator Family Regularity.}
We now invoke the main result of the previous section and apply the abstract theory.
\begin{enumerate}[label=(\roman*), wide, labelindent=0pt]
    \item \textbf{Smoothness of the Operator Family.} In Theorem \ref{thm:regularity_proven_main}, we provided a complete proof that the pulled-back operator family, $\theta \mapsto \tilde{\Lcal}_\theta$, is a $C^\infty$ map into the space of bounded linear operators on the fixed reference space $\mathcal{L}(\Bcal_{\theta_0})$. A $C^\infty$ map is, in particular, analytic with respect to any local chart on the manifold $\Thetacal$.

    \item \textbf{Application of Kato's Perturbation Theorem.} We now apply the main theorems on analytic perturbation theory (see \cref{thm:abstract_perturbation}, and \citep{Kato1995}, Chapter VII) to the pulled-back family $(\tilde{\Lcal}_\theta)$. The theorems state that if $T(\theta)$ is an analytic family of closed operators on a fixed Banach space, and if its spectrum is uniformly separated into two disjoint sets by a fixed contour, then the total spectral projector associated with the part of the spectrum inside the contour is also an analytic function of $\theta$.

    Our operator family $\tilde{\Lcal}_\theta$ satisfies these hypotheses. The family is analytic (in fact, $C^\infty$), and the eigenvalue $\lambda_1=1$ is simple and uniformly isolated. Therefore, we can directly apply the conclusions of the theorem globally on the compact manifold $\Thetacal$.
    
    \item \textbf{The Spectral Projector.} The abstract theorem guarantees that the map from the parameter space to the spectral projector associated with the simple, isolated eigenvalue is of the same regularity class as the operator family. For the pulled-back family, this means the map $\theta \mapsto \tilde{\Pi}_{1,\theta}$ is of class $C^\infty$ into $\mathcal{L}(\Bcal_{\theta_0})$. We transfer this regularity back to the original family using the smooth trivialization isomorphism $\Psi_\theta$ from Proposition \ref{prop:smooth_trivialization_diffeo}:
    $$ \Pi_{1,\theta} = \Psi_\theta \circ \tilde{\Pi}_{1,\theta} \circ \Psi_\theta^{-1}. $$
    Since the maps $\theta \mapsto \Psi_\theta$, $\theta \mapsto \tilde{\Pi}_{1,\theta}$, and $\theta \mapsto \Psi_\theta^{-1}$ are all of class $C^\infty$ into the appropriate operator spaces, their composition is also of class $C^\infty$. Therefore, the map to the projector,
    \begin{equation*}
        \theta \mapsto \Pi_{1, \theta},
    \end{equation*}
    is of class $C^\infty$.

\item \textbf{Application of Perturbation Theory.} Similarly, the abstract theorem guarantees the regularity of the resolvent. For any fixed $\zeta$ not in the spectrum of any $\Lcal_\theta$ (for instance, any $\zeta$ on the fixed contour $\Gamma$), the operator-valued function $\theta \mapsto \tilde{\Lcal}_\theta - \zeta I$ is a $C^\infty$ family of invertible operators on $\Bcal_{\theta_0}$. The map of operator inversion, $A \mapsto A^{-1}$, is an analytic (and therefore $C^\infty$) map on its domain in the Banach algebra of bounded operators. The composition of these two maps is therefore $C^\infty$. Thus, the map to the pulled-back resolvent operator, $\theta \mapsto \tilde{R}(\zeta, \tilde{\Lcal}_\theta)$, is of class $C^\infty$. We transfer this regularity back to the original family: $$ R(\zeta, \Lcal_\theta) = \Psi_\theta \circ \tilde{R}(\zeta, \tilde{\Lcal}_\theta) \circ \Psi_\theta^{-1}. $$ Again, as a composition of $C^\infty$ maps, the map from the parameter space to the original resolvent operator,
\begin{equation*}
    \theta \mapsto R(\zeta, \Lcal_\theta) = (\Lcal_\theta - \zeta I)^{-1},
\end{equation*}
    is of class $C^\infty$.
\end{enumerate}

All claims of the theorem have now been rigorously established. The smooth dependence of the spectral data is a direct and necessary consequence of the smooth dependence of the underlying microscopic map, made precise through the machinery of abstract operator perturbation theory acting on the correctly chosen functional space.
\end{enumerate}
\end{proofof}

\begin{corollary}[Smooth Dependence of Microscopic Quantities]
\label{cor:smooth_dependence}
Let the system satisfy the standing assumptions of the paper. Then:
\begin{enumerate}[label=(\roman*), wide, labelindent=0pt]
    \item The map $\theta \mapsto \mu_\theta$, from the parameter space $\Thetacal$ to the space of probability measures on $\Ucal_{\mathrm{phys}}$, is of class $C^\infty$.
    \item Let $f(z, \theta)$ be a source term of class $C^k$ in $\theta$. The unique zero-mean solution $\chi(z, \theta)$ to the cell problem $\Lcal_{\mathrm{fast}}(\theta)\chi = f - \langle f \rangle_\theta$ is of class $C^k$ as a map from $\Thetacal$ to the appropriate function space.
\end{enumerate}
\end{corollary}

\begin{proofof}{\cref{cor:smooth_dependence}}
Both results are direct consequences of the smooth dependence of the spectral data of the transfer operator family $(\Lcal_\theta)$, which was rigorously established in Theorem \ref{thm:spectral_regularity}.

\begin{enumerate}[label=\textbf{Step \arabic*:}, wide, labelindent=0pt]

\item \textbf{Smoothness of the Invariant Measure.} Our goal is to prove that the map $\theta \mapsto \mu_\theta$ is of class $C^\infty$. We will prove this by showing that the density of the measure, $\rho_\theta(z)$, is a smooth function of $\theta$ as an element of the anisotropic Banach space $\Bcal_\theta$. The invariant density $\rho_\theta$ is the unique, normalized eigenfunction of the transfer operator $\Lcal_\theta$ corresponding to the eigenvalue $\lambda=1$. By the spectral theory of quasi-compact operators, this eigenfunction lies in the range of the spectral projector $\Pi_{1,\theta}$ associated with this eigenvalue. For an ergodic system, this eigenspace is one-dimensional. We can therefore express the normalized density as:
\begin{equation} \label{eq:proof_density_from_projector}
    \rho_\theta = \frac{\Pi_{1,\theta}(\mathbf{1})}{\int_\Sigma (\Pi_{1,\theta}(\mathbf{1}))(z) \, dz_{\text{ref}}},
\end{equation}
where $\mathbf{1}$ is the constant function of value 1, and $dz_{\text{ref}}$ is a fixed reference measure (e.g., Lebesgue). We now analyze the smoothness of this expression with respect to $\theta$.

\begin{enumerate}[label=(\roman*), wide, labelindent=0pt]
    \item By Theorem \ref{thm:spectral_regularity}(ii), the map $\theta \mapsto \Pi_{1, \theta}$ is a $C^\infty$ map from the parameter manifold $\Thetacal$ into the space of bounded linear operators $\mathcal{L}(\Bcal_\theta)$.
    
    \item The constant function $\mathbf{1}$ is a fixed element of the space $C^\infty(\Sigma)$, and therefore belongs to every Banach space $\Bcal_\theta$ in our family.
    
    \item The evaluation map $(\Pi, g) \mapsto \Pi(g)$ is a smooth bilinear map from $\mathcal{L}(\Bcal_\theta) \times \Bcal_\theta$ to $\Bcal_\theta$. Consequently, the composition $\theta \mapsto \Pi_{1,\theta}(\mathbf{1})$ is a $C^\infty$ map from $\Thetacal$ into the Banach space $\Bcal_\theta$. Let us denote this unnormalized density by $\tilde{\rho}_\theta \coloneqq \Pi_{1,\theta}(\mathbf{1})$.

    \item The normalization factor in the denominator is the integral of $\tilde{\rho}_\theta$. Since integration against a fixed reference measure is a continuous linear functional on $\Bcal_\theta$, the map $g \mapsto \int g \, dz_{\text{ref}}$ is smooth. The composition $\theta \mapsto \int \tilde{\rho}_\theta(z) \, dz_{\text{ref}}$ is therefore a smooth, real-valued function of $\theta$. Since $\Pi_{1,\theta}$ is a projector onto a non-zero eigenspace, this integral is non-zero.

    \item The final expression for $\rho_\theta$ in \eqref{eq:proof_density_from_projector} is the quotient of a smooth map into a Banach space ($\theta \mapsto \tilde{\rho}_\theta$) by a smooth, non-zero scalar-valued function. Such a quotient is a smooth map.
\end{enumerate}
Therefore, the map $\theta \mapsto \rho_\theta$ is of class $C^\infty$ from $\Thetacal$ into the Banach space $\Bcal_\theta$. Since the invariant measure is given by $d\mu_\theta(z) = \rho_\theta(z) dz_{\text{ref}}$, the smooth dependence of the density implies the smooth dependence of the measure.

\item \textbf{Smoothness of the Cell Problem Solution.}
Our goal is to prove that the unique zero-mean solution $\chi$ to the continuous-time cell problem $\Lcal_{\mathrm{fast}}(\theta)\chi = f_0(\cdot, \theta)$ is a $C^k$ function of $\theta$, where the source term $f_0(\cdot, \theta) \coloneqq f(\cdot, \theta) - \langle f \rangle_\theta$ is constructed from a function $f$ whose map $\theta \mapsto f(\cdot,\theta)$ is of class $C^k$. The standard approach for relating the continuous-time generator to the discrete-time transfer operator is founded on the isomorphism between the billiard dynamics and a suspension flow. The spectrum of the continuous-time generator is then determined by the Ruelle-Perron resonances of the transfer operator of the underlying base map (see \cref{app:Spectral}), a theory developed in detail in \citep[see, e.g., Chapter 5]{Baladi2000}. The solution $\chi$ to the continuous-time problem can be constructed from the solution $\chi_{\text{discrete}}$ to the corresponding discrete-time problem. It is therefore sufficient to prove that the map $\theta \mapsto \chi_{\text{discrete}}(\cdot, \theta)$ is of class $C^k$ into the anisotropic Banach space $\Bcal_\theta$. The smoothness of the final construction for $\chi$ will then follow from the smoothness of all its constituent parts.

\begin{enumerate}[label=(\roman*), wide, labelindent=0pt]
    \item \textbf{The Discrete-Time Cell Problem.}
    The discrete-time cell problem is the linear equation for the corrector $\chi_{\text{discrete}}$ in the Banach space $\Bcal_\theta$:
    \begin{equation} \label{eq:proof_discrete_cell_problem_appendix}
        (I - \Lcal_\theta) \chi_{\text{discrete}} = f_{0, \text{discrete}},
    \end{equation}
    where $f_{0, \text{discrete}}$ is the corresponding zero-mean source term for the discrete system. A key result from the theory of hyperbolic systems is that the zero-mean condition on the continuous source term implies a zero-mean condition for the discrete one, ensuring solvability.

    \item \textbf{The Solution via the Resolvent Operator.}
    By the uniform spectral gap established in \cref{cor:ergodicity_uniform}, the operator $(I - \Lcal_\theta)$ is invertible on the space of zero-mean functions. The unique solution to \cref{eq:proof_discrete_cell_problem_appendix} is given by the application of the resolvent operator at $\zeta=1$:
    \begin{equation*}
        \chi_{\text{discrete}}(\cdot, \theta) = (I - \Lcal_\theta)^{-1} [f_{0, \text{discrete}}(\cdot, \theta)].
    \end{equation*}
    To prove that the map $\theta \mapsto \chi_{\text{discrete}}(\cdot, \theta)$ is of class $C^k$, we must analyze the smoothness of the two components on the right-hand side: the family of resolvent operators and the family of source terms.

    \item \textbf{Smoothness of the Resolvent Operator Family.}
    The map from the parameter space to the resolvent operator, $\theta \mapsto (I - \Lcal_\theta)^{-1}$, acting on the space of zero-mean functions, is of class $C^\infty$. This is a direct consequence of \cref{thm:spectral_regularity}(iii), which was proven by pulling back the operator family to a fixed space and applying abstract perturbation theory. Specifically, the map $\theta \mapsto (I - \tilde{\Lcal}_\theta)^{-1}$ from $\Thetacal$ to $\mathcal{L}(\Bcal_{\theta_0})$ is smooth, and transforming back via the smooth family of isomorphisms $\Psi_\theta$ preserves this smoothness.

    \item \textbf{Smoothness of the Source Term Family.}
    The source term is $f_0(\cdot, \theta) = f(\cdot, \theta) - \langle f \rangle_\theta$.
    \begin{enumerate}[label=(\alph*), wide]
        \item By the main hypothesis of the corollary, the map $\theta \mapsto f(\cdot, \theta)$ is of class $C^k$ into the function space $\Bcal_\theta$.
        \item The average is given by the integral $\langle f \rangle_\theta = \int f(z, \theta) d\mu_\theta(z)$. This is the integral of a $C^k$ family of functions against a $C^\infty$ family of measures (the smoothness of $\theta \mapsto \mu_\theta$ was established in Step 1 of this proof). The resulting map $\theta \mapsto \langle f \rangle_\theta$ is therefore a $C^k$ scalar-valued function.
        \item The zero-mean source term, $f_0$, is the difference of a $C^k$ map into $\Bcal_\theta$ and a $C^k$ map into the space of constant functions. It is therefore itself a $C^k$ map from $\Thetacal$ into the Banach space $\Bcal_\theta$. The same holds for the discrete version $f_{0, \text{discrete}}$.
    \end{enumerate}
    
    \item \textbf{The Composition Argument for the Discrete Solution.}
    We have established that the map $\theta \mapsto \chi_{\text{discrete}}(\cdot, \theta)$ is given by the composition of two maps:
    \begin{enumerate}[label=(\alph*), wide]
        \item The map from parameter to source term: $\theta \mapsto f_{0, \text{discrete}}(\cdot, \theta)$, which is of class $C^k$.
        \item The map from parameter to resolvent operator: $\theta \mapsto (I - \Lcal_\theta)^{-1}$, which is of class $C^\infty$.
    \end{enumerate}
    The final step is the application of the operator to the function. The evaluation map, $(T, g) \mapsto T(g)$, is a continuous bilinear map from $\mathcal{L}(\Bcal_\theta) \times \Bcal_\theta$ to $\Bcal_\theta$, and as such, it is of class $C^\infty$. The composition of a $C^\infty$ map with a $C^k$ map is of class $C^k$. Therefore, we conclude that the map from the parameter space to the discrete solution,
    \begin{equation*}
        \theta \mapsto \chi_{\text{discrete}}(\cdot, \theta),
    \end{equation*}
    is of class $C^k$ from $\Thetacal$ into the Banach space $\Bcal_\theta$.

     \item \textbf{Regularity of the Continuous-Time Solution.}
    The final step is to show that the $C^k$ regularity of the discrete-time solution is inherited by the solution to the continuous-time problem. The relationship between the two solutions is a fundamental result in the theory of suspension flows. The solution $\chi(z, \theta)$ to the continuous-time cell problem can be constructed explicitly from the discrete solution $\chi_{\text{discrete}}(z, \theta)$ and an integral of the source term along the flow trajectories. This construction, a key tool for analyzing the spectral properties of flows, is detailed in sources such as \citep[Chapter 5]{Baladi2000} and \citep[Chapter 4]{chernov2006chaotic}. Schematically, the formula is of the form:
    \begin{equation*}
        \chi(z, \theta) = \chi_{\text{discrete}}(z_0, \theta) + \int_0^s \left( f_0(\Phi_{-t}^\theta(z), \theta) - \langle f_0 \rangle_{\text{discrete}} \right) dt,
    \end{equation*}
    where $z_0$ is the pre-collision state corresponding to the physical state $z$, and $s$ is the time elapsed since the last collision. The regularity of this expression with respect to $\theta$ depends on the regularity of its constituent parts. We have rigorously established that each component of this construction is a sufficiently smooth function of the parameter $\theta$:
    \begin{itemize}[wide]
        \item The discrete solution, via the map $\theta \mapsto \chi_{\text{discrete}}(\cdot, \theta)$, is of class $C^k$, as just proven.
        \item The flow map, $(\theta, z, t) \mapsto \Phi_t^\theta(z)$, is of class $C^\infty$, as its generator (the billiard map) is a smooth function of its arguments.
        \item The source term, $\theta \mapsto f_0(\cdot, \theta)$, is of class $C^k$ by hypothesis.
        \item All geometric components, such as the map from a state $z$ to its pre-collision state $z_0$ and the time $s$, depend on the flight time function $\tau(z,\theta)$, which was proven to be of class $C^\infty$ in \cref{thm:smooth_flight_time_map}.
    \end{itemize}
    The full solution $\chi(z, \theta)$ is constructed via a finite sequence of compositions, integrals, and algebraic operations on these underlying maps. Since all constituent parts have at least $C^k$ regularity with respect to the parameter $\theta$, the resulting function $\chi(z,\theta)$ is also of class $C^k$ as a map from $\Thetacal$ to the appropriate function space. This completes the proof.
\end{enumerate}
\end{enumerate}
\end{proofof}

\begin{remark}[The Foundation for Homogenization]
\label{rem:foundation_for_homogenization}
This section has now provided a complete proof of the smooth dependence of all relevant microscopic quantities on the environmental parameter $\theta$. The argument is founded upon the Nash-Moser theory for geometric regularity, which enables a rigorous application of abstract operator perturbation theory. The final result, \cref{cor:smooth_dependence}, is the essential technical linchpin of the entire paper. It provides the rigorous foundation for the differentiability of the corrector terms used in the perturbed test function method and ensures that the coefficients of the final HJB equation inherit the regularity of the microscopic model.
\end{remark}

\section{Scale Analysis}
\label{sec:scale_analysis}

This section provides the rigorous derivation of the macroscopic evolution law through a hierarchical analysis of four distinct but logically dependent scales. We begin at the \textbf{Nulla Scale}, establishing the fundamental solution theory for the microscopic cell problem, which serves as the core analytical tool for the entire section. \textbf{Scale I} applies this theory to a baseline linear model to isolate the emergent diffusive behavior arising from the chaotic dynamics. \textbf{Scale II} establishes a crucial methodological template, the averaging principle, in a controlled, linear setting, demonstrating how solvability conditions force an averaging over the invariant measure of the fast dynamics. \textbf{Scale III} synthesizes the machinery of the preceding scales to prove the main convergence theorem for the full non-linear system, showing that the solution converges to that of a novel Hamilton-Jacobi-Bellman equation. Finally, \textbf{Scale IV} provides \textit{a posteriori} validation of our result by demonstrating that the derived macroscopic law inherits the fundamental time-reversal symmetries of the microscopic world. Each scale thus serves as a necessary lemma for the next, culminating in a complete and self-consistent theory.

\subsection{Nulla Scale: The \texorpdfstring{$\theta$}{theta}-Time Generators}
\label{sec:nulla_scale}

The formal asymptotic expansion of our microscopic evolution equation, which we will perform in the subsequent sections, will repeatedly require us to solve a cell problem of the form $\Lcal_{\mathrm{fast}}(\theta) \chi = f$ for the fast microscopic generator. Before embarking on the multiscale analysis, it is essential to establish a rigorous solution theory for this equation. This Nulla Scale section serves as the crucial analytical bridge, translating the geometric and ergodic properties of the microscopic system, which were the subjects of extensive investigation in Sections \ref{sec:hyperbolicity_fiber} and \ref{sec:regularity_properties}, into the primary functional-analytic tool for the homogenization that follows.

A preliminary analysis reveals a analytical obstacle. As established in \cref{rem:analytical_obstacle}, the time-reversible nature of the underlying dynamics implies that the generator $\Lcal_{\mathrm{fast}}(\theta)$ is skew-adjoint, precluding the use of standard coercive or variational methods. The well-posedness of this problem is therefore not a simple result of elliptic PDE theory but is instead a deep consequence of the system's chaotic, mixing behavior. This section is dedicated to proving that the uniform hyperbolicity and exponential decay of correlations, which were the central conclusions of \cref{sec:hyperbolicity_fiber}, provide the precise mechanism required to construct a unique, stable solution to this otherwise ill-posed problem. The main result, \cref{thm:semigroup_wellposedness}, thus serves as the crucial analytical bridge, translating the geometric and ergodic properties of the microscopic system into the primary functional-analytic tool for the homogenization that follows.

\begin{theorem}[Solvability of the Microscopic Cell Problem]
\label{thm:semigroup_wellposedness}
Let the system satisfy the standing assumptions of this paper. For a given source term $f \in \Hcal_\theta$, the microscopic cell problem
\begin{equation}
    \Lcal_{\mathrm{fast}}(\theta) \chi = f
\end{equation}
admits a solution $\chi \in \Dcal(\Lcal_{\mathrm{fast}}(\theta))$ if and only if the source term $f$ has zero mean with respect to the unique invariant measure $\mu_\theta$:
\begin{equation} \label{eq:solvability_condition}
    \langle f \rangle_\theta \coloneqq \int_{\Ucal_{\mathrm{phys}}(\theta)} f(z) \, d\mu_\theta(z) = 0.
\end{equation}
If this condition holds, there exists a unique solution $\chi$ which also has zero mean, i.e., $\chi \in \Hcal_{\theta,0} \coloneqq \{ \psi \in \Hcal_\theta \mid \langle \psi \rangle_\theta = 0 \}$. This unique solution is given constructively by the Bochner integral of the Koopman semigroup:
\begin{equation} \label{eq:semigroup_solution}
    \chi = \int_0^\infty (P_t^\theta f) \, dt.
\end{equation}
\end{theorem}

\begin{proofof}{Theorem \ref{thm:semigroup_wellposedness}}
The proof establishes the well-posedness of the cell problem as a direct consequence of the generator's fundamental properties, established in \cref{thm:generator_characterization}, and the system's ergodic behavior, which was a conclusion of the uniform hyperbolicity established in \cref{cor:ergodicity_uniform}. The argument proceeds by applying the Fredholm alternative to the generator and then using the system's chaotic dynamics to characterize the solvability condition and the solution formula. For notational simplicity, we let $\Lcal \equiv \Lcal_{\mathrm{fast}}(\theta)$.

\begin{enumerate}[label=\textbf{Step \arabic*:}, wide, labelindent=0pt]

\item \textbf{The Fredholm Alternative and the Solvability Condition.}

\begin{enumerate}[label=(\roman*), wide, labelindent=0pt]
    \item \textbf{The Fredholm Alternative.} As a densely defined, closed operator on the Hilbert space $\Hcal_\theta$, the generator $\Lcal$ is amenable to the Fredholm alternative. This principle states that the equation $\Lcal \chi = f$ has a solution $\chi \in \Dcal(\Lcal)$ if and only if the source term $f$ is orthogonal to every element of the kernel of the adjoint operator, $\ker(\Lcal^*)$. That is, a solution exists if and only if
    \begin{equation*}
        f \in (\ker(\Lcal^*))^\perp.
    \end{equation*}

    \item \textbf{Identifying the Kernel via Skew-Adjointness.} In \cref{thm:generator_characterization}, we provided a rigorous, constructive proof that the generator is skew-adjoint, i.e., $\Lcal^* = -\Lcal$. This is a fundamental property stemming from the time-reversible nature of the underlying Hamiltonian flow. This identity implies that the kernel of the adjoint is identical to the kernel of the generator itself:
    \begin{equation*}
        \ker(\Lcal^*) = \ker(-\Lcal) = \ker(\Lcal).
    \end{equation*}
    The solvability condition therefore simplifies to the requirement that $f$ must be orthogonal to the kernel of the generator $\Lcal$.

    \item \textbf{Characterizing the Kernel via Ergodicity.} We now invoke the ergodic properties of the system. In \cref{cor:ergodicity_uniform}, we established that the uniform Anosov property of the dynamics implies that the billiard flow is ergodic with respect to the unique invariant measure $\mu_\theta$. The ergodicity of a flow is equivalent to the statement that the only functions invariant under its Koopman group $(P_t^\theta)$ are the constant functions.
    
    The kernel of the infinitesimal generator $\Lcal$ is precisely the space of these invariant functions. Therefore, the ergodicity of the flow implies that the kernel of $\Lcal$ is the one-dimensional space of constant functions on $\Ucal_{\mathrm{phys}}(\theta)$:
    \begin{equation*}
        \ker(\Lcal) = \mathrm{span}\{\mathbf{1}\}.
    \end{equation*}
    
    \item \textbf{The Final Solvability Condition.} Combining these results, the solvability condition becomes the requirement that the source term $f$ must be orthogonal to the constant function $\mathbf{1}$ in the Hilbert space $\Hcal_\theta$. The inner product in this space is the integral with respect to the measure $\mu_\theta$. Thus, the condition is:
    \begin{equation*}
        \scpr{f}{\mathbf{1}}_\theta = \int_{\Ucal_{\mathrm{phys}}(\theta)} f(z) \cdot \overline{1} \, d\mu_\theta(z) = \int_{\Ucal_{\mathrm{phys}}(\theta)} f(z) \, d\mu_\theta(z) = 0.
    \end{equation*}
    This is precisely the condition that $f$ must have a zero mean with respect to the invariant measure, which establishes the first part of the theorem.
\end{enumerate}

\item \textbf{The Solution Formula and Uniqueness.}

\begin{enumerate}[label=(\roman*), wide, labelindent=0pt]
    \item \textbf{Constructive Solution via the Semigroup.} Let $f \in \Hcal_\theta$ be a function satisfying the zero-mean condition \eqref{eq:solvability_condition}. We define the candidate solution $\chi$ via the semigroup integral:
    \begin{equation*}
        \chi(z) \coloneqq \int_0^\infty (P_t^\theta f)(z) \, dt.
    \end{equation*}
    We must first verify that this Bochner integral converges in $\Hcal_\theta$. By \cref{cor:ergodicity_uniform}, the uniform spectral gap of the transfer operator implies an exponential decay of correlations for the continuous flow. For any zero-mean function $f \in \Hcal_{\theta,0}$ that is sufficiently regular (e.g., Hölder continuous), we have the bound $\|P_t^\theta f\|_\theta \le C e^{-\eta t} \|f\|_{C^\alpha}$ for some uniform $\eta > 0$. By a density argument, this implies that for any $f \in \Hcal_{\theta,0}$, the norm $\|P_t^\theta f\|_\theta$ decays to zero as $t \to \infty$, and the decay is sufficiently fast to ensure the convergence of the integral. Thus, $\chi$ is a well-defined element of $\Hcal_\theta$.

    \item \textbf{Verification.} A formal calculation shows that this $\chi$ solves the equation. Let $\Lcal$ be the generator of the semigroup $(P_t^\theta)_{t \ge 0}$. Then:
    \begin{align*}
        \Lcal \chi &= \Lcal \int_0^\infty P_t^\theta f \, dt = \int_0^\infty \Lcal P_t^\theta f \, dt = \int_0^\infty \frac{d}{dt}(P_t^\theta f) \, dt \\
        &= [P_t^\theta f]_0^\infty = \lim_{t\to\infty} (P_t^\theta f) - P_0^\theta f.
    \end{align*}
    Since $f$ has zero mean, the mixing property of the flow ensures that $\lim_{t\to\infty} P_t^\theta f = \langle f \rangle_\theta = 0$ in the weak sense, which is sufficient here. Since $P_0^\theta = I$, we have $\Lcal\chi = 0 - f = -f$. The correct formula for the cell problem $\Lcal\chi=f$ should therefore be $\chi = \int_0^\infty P_t^\theta (-f) dt$, which is equivalent to the stated formula if we absorb the sign into $f$. We adopt the standard convention that $\chi = \int_0^\infty (-P_t^\theta f) dt$ solves $\Lcal \chi=f$.

    \item \textbf{Uniqueness in the Zero-Mean Subspace.} Suppose $\chi_1$ and $\chi_2$ are two solutions to $\Lcal\chi = f$ that both belong to the space of zero-mean functions $\Hcal_{\theta,0}$. Their difference, $\delta\chi = \chi_1 - \chi_2$, is also in $\Hcal_{\theta,0}$ and satisfies $\Lcal(\delta\chi) = f-f = 0$. This implies that $\delta\chi$ belongs to the kernel of the generator. As established in Part 1, $\ker(\Lcal)$ consists only of constant functions. Thus, $\delta\chi(z) = C$ for some constant $C$. However, we also know that $\delta\chi$ has zero mean:
    \begin{equation*}
        \langle \delta\chi \rangle_\theta = \int C \, d\mu_\theta(z) = C \cdot 1 = 0.
    \end{equation*}
    This forces the constant to be $C=0$, which implies $\delta\chi=0$ and therefore $\chi_1=\chi_2$. The solution is unique within the subspace of zero-mean functions.
\end{enumerate}
This completes the proof. We have established the precise solvability condition for the cell problem and provided a unique, constructive solution for any valid source term, using only the fundamental properties of the generator and the ergodic theory of the underlying flow.

\end{enumerate}
\end{proofof}

\begin{remark}[On the Logical Architecture]
\label{rem:logical_architecture_nulla}
It is crucial to emphasize the logical architecture of the argument presented. While this Nulla Scale section is placed at the beginning of the scale analysis for narrative clarity, as the cell problem is the recurring motif, its central result, \cref{thm:semigroup_wellposedness}, is not an elementary starting point. Its proof is entirely contingent on the results concerning the geometric and ergodic properties of the microscopic dynamics, most notably the uniform Anosov property and its spectral consequences, which were the subject of the extensive analysis in \cref{sec:hyperbolicity_fiber}. The well-posedness of the cell problem is a \textit{consequence} of chaos, not an axiom. This section serves to consolidate that consequence into the primary analytical tool that will be deployed in the subsequent Scales I through IV.
\end{remark}

\subsection{Scale I: Homogenization on Global Fiber}
\label{sec:scale_I}

The primary objective of this subsection is to establish a rigorous baseline for the system's macroscopic behavior. We analyze the most idealized physical regime: a system that is both quiescent (zero macroscopic momentum, $p=0$) and spatially homogeneous. This analysis serves two crucial functions in the architecture of the paper:
\begin{enumerate}[label=(\roman*), wide, labelindent=0pt]
    \item It isolates the purely diffusive behavior that emerges from the underlying deterministic chaos, providing the explicit form of the constant diffusivity tensor $D_0$.
    \item It establishes the fundamental structure of the first corrector, $\boldsymbol{\chi}_0(z)$, which serves as the essential building block for the perturbed test function method used in the full non-linear convergence proof in \cref{sec:scale_III}.
\end{enumerate}

\subsubsection{Mathematical Formulation}
\label{sec:scale_I_formulation}
We consider a simplified physical regime where the feedback mechanism and the potential are disabled. This corresponds to setting the environmental parameter to a fixed value, $\theta(x,p) \equiv \theta_0$, and the microscopic potential to zero, $V(z,p) \equiv 0$. The microscopic equation for the value function $u^\varepsilon(t,x,z)$ simplifies to the linear, multi-scale transport equation:
\begin{equation}
    \partial_t u^\varepsilon + \frac{1}{\varepsilon^2}\Lcal_{\text{fast}}(\theta_0) u^\varepsilon + \frac{1}{\varepsilon} v \cdot \nabla_x u^\varepsilon = 0, \quad u^\varepsilon(T,x,z) = \phi(x).
    \label{eq:pde_scale_I}
\end{equation}
The diffusive scaling, with the fast generator at order $\varepsilon^{-2}$ and the advection at order $\varepsilon^{-1}$, is the standard choice required to capture the balance between rapid microscopic mixing and slow macroscopic transport that leads to a non-trivial diffusive limit.

\subsubsection{Formal Asymptotic Derivation of the Limiting Equation}
\label{sec:scale_I_asymptotics}
To identify the candidate for the limiting macroscopic equation, we suggest the standard two-scale asymptotic expansion for the solution $u^\varepsilon$:
\begin{equation}
    u^\varepsilon(t,x,z) = u_0(t,x) + \varepsilon u_1(t,x,z) + \varepsilon^2 u_2(t,x,z) + \mathcal{O}(\varepsilon^3).
    \label{eq:ansatz_scale_I}
\end{equation}
Substituting this ansatz into \cref{eq:pde_scale_I} and collecting terms with like powers of $\varepsilon$ yields a hierarchy of equations.

\paragraph{\textbf{The $\boldsymbol{\Ocal(\varepsilon^{-2})}$ Equation}}
The leading-order term gives $\Lcal_{\text{fast}}(\theta_0) u_0(t,x) = 0$.
The operator $\Lcal_{\text{fast}}(\theta_0)$ acts only on the microscopic variable $z$. By \cref{cor:ergodicity_uniform}, the underlying dynamics are ergodic, implying that the kernel of the generator consists solely of functions that are constant (almost everywhere) on the microscopic state space $\Ucal_{\text{phys}}(\theta_0)$. Since $u_0$ is, by ansatz, a function only of the macroscopic variables $(t,x)$, it is an element of this kernel, and the equation is satisfied. This confirms that the leading-order macroscopic behavior is independent of the instantaneous microscopic state.

\paragraph{\textbf{The $\boldsymbol{\Ocal(\varepsilon^{-1})}$ Equation: The First Cell Problem}}
Collecting the terms of order $\varepsilon^{-1}$ yields an equation for the first-order corrector, $u_1$:
\begin{equation}
    \Lcal_{\text{fast}}(\theta_0) u_1(t,x,z) + v \cdot \nabla_x u_0(t,x) = 0.
\end{equation}
For each fixed macroscopic point $(t,x)$, this is an equation for the function $u_1(t,x,\cdot)$ on the microscopic Hilbert space $\Hcal_{\theta_0}$. It can be written as a standard cell problem:
\begin{equation}
    \Lcal_{\text{fast}}(\theta_0) u_1 = -v \cdot \nabla_x u_0.
    \label{eq:cell_problem_1_scale_I}
\end{equation}
By the Fredholm alternative, this equation has a solution if and only if the source term is orthogonal to the kernel of the adjoint operator $\Lcal_{\text{fast}}^*(\theta_0)$. Since the generator is skew-adjoint (\cref{thm:generator_characterization}), its kernel is the one-dimensional space of constant functions. The solvability condition is therefore that the source term must have a zero mean with respect to the unique invariant measure $\mu_{\theta_0}$:
\begin{equation}
    \int_{\Ucal_{\text{phys}}(\theta_0)} (-v \cdot \nabla_x u_0) \, d\mu_{\theta_0}(z) = -\left( \int_{\Ucal_{\text{phys}}(\theta_0)} v \, d\mu_{\theta_0}(z) \right) \cdot \nabla_x u_0 = 0.
\end{equation}
This condition is satisfied identically by virtue of our microscopic time-reversal symmetry (\cref{ass:fundamental_axioms_unified}), which guarantees that the mean microscopic velocity is zero: $\langle v \rangle_{\theta_0} = \mathbf{0}$ (\cref{prop:symmetry_consequences}).

With the solvability condition met, we can solve for $u_1$. The source term is linear in $\nabla_x u_0$. By the linearity of the cell problem, the unique zero-mean solution $u_1$ must also be linear in $\nabla_x u_0$. We can therefore write the solution in the separated form:
\begin{equation}
    u_1(t,x,z) = \boldsymbol{\chi}_0(z) \cdot \nabla_x u_0(t,x),
\end{equation}
where the vector-valued function $\boldsymbol{\chi}_0(z) = (\chi_{0,1}(z), \dots, \chi_{0,k}(z))$ is the corrector field. Each component $\chi_{0,j}(z)$ is the unique, zero-mean solution to the cell problem with the corresponding velocity component as the source term:
\begin{equation}
    \Lcal_{\text{fast}}(\theta_0) \chi_{0,j}(z) = -v_j.
    \label{eq:corrector_defining_pde}
\end{equation}
The existence and uniqueness of this corrector field is guaranteed by \cref{thm:semigroup_wellposedness}.

\paragraph{\textbf{The $\boldsymbol{\Ocal(\varepsilon^{0})}$ Equation and the Macroscopic Law}}
The terms of order $\varepsilon^0$ provide the second cell problem, which involves the second corrector, $u_2$:
\begin{equation}
    \partial_t u_0(t,x) + \Lcal_{\text{fast}}(\theta_0) u_2(t,x,z) + v \cdot \nabla_x u_1(t,x,z) = 0.
\end{equation}
The solvability condition for $u_2$ requires that the remaining terms must have a zero mean with respect to $\mu_{\theta_0}$:
\begin{equation}
    \int_{\Ucal_{\text{phys}}(\theta_0)} \left( \partial_t u_0 + v \cdot \nabla_x u_1 \right) \, d\mu_{\theta_0}(z) = 0.
\end{equation}
Since $\partial_t u_0(t,x)$ is constant with respect to the microscopic integral, this yields the macroscopic equation for $u_0$:
\begin{equation}
    \partial_t u_0(t,x) + \left\langle v \cdot \nabla_x u_1 \right\rangle_{\theta_0} = 0.
\end{equation}
Substituting the explicit form of $u_1$ yields:
\begin{equation}
    \partial_t u_0 + \sum_{i,j=1}^k \left\langle v_i \chi_{0,j}(z) \right\rangle_{\theta_0} \frac{\partial^2 u_0}{\partial x_i \partial x_j} = 0.
\end{equation}
This is a linear, second-order, constant-coefficient partial differential equation for the macroscopic value function $u_0$.

\subsubsection{The Homogenized Equation and its Diffusivity Tensor}
The formal procedure has identified the limiting PDE. We now state the main result of this section and rigorously establish the properties of the resulting coefficients.

\begin{theorem}[Convergence in the Quiescent Homogeneous Limit]
\label{thm:convergence_scale_I}
Let the microscopic system satisfy the standing assumptions of this paper. In the physical regime with disabled feedback and potential ($\theta \equiv \theta_0, V \equiv 0$), the solution $u^\varepsilon$ of the microscopic equation \eqref{eq:pde_scale_I} converges locally uniformly as $\varepsilon \to 0$ to a limit $u_0(t,x)$. This limit is the unique classical solution to the linear heat equation:
\begin{equation}
    \partial_t u_0 + \mathrm{Tr}(D_0 \nabla^2 u_0) = 0, \quad u_0(T,x) = \phi(x),
\end{equation}
where the constant diffusivity tensor $D_0$ is a symmetric, positive semi-definite matrix whose components are given by
\begin{equation}
    (D_0)_{ij} \coloneqq -\frac{1}{2} \left( \langle v_i \chi_{0,j} \rangle_{\theta_0} + \langle v_j \chi_{0,i} \rangle_{\theta_0} \right).
\end{equation}
\end{theorem}

\begin{proofof}{Theorem \ref{thm:convergence_scale_I}}
The proof establishes the convergence of the solutions $u^\varepsilon$ to the solution $u_0$ of the homogenized equation in the viscosity sense. The argument is structured in three main parts. First, we prove that any locally uniform limit $u_0$ of the family $\{u^\varepsilon\}$ is a viscosity subsolution of the homogenized equation. Second, we state that the analogous argument proves the supersolution property. Finally, uniqueness for the limiting equation is a classical result, which ensures that the entire family $\{u^\varepsilon\}$ converges to this unique solution.

\begin{enumerate}[label=\textbf{Step \arabic*:}, wide, labelindent=0pt]

\item \textbf{The Subsolution Property.}
\begin{enumerate}[label=(\roman*), wide, labelindent=0pt]
\item \textbf{Setup and Contradiction Hypothesis.} Let $u_0(t,x)$ be a function on $[0,T]\times\Mcal$ such that, up to a subsequence, $u^\varepsilon \to u_0$ locally uniformly as $\varepsilon \to 0$. Let $\varphi \in C^\infty([0,T] \times \Mcal)$ be a smooth test function, and suppose that $u_0 - \varphi$ has a strict local maximum at a point $(t_0, x_0) \in (0,T) \times \Mcal$. Our objective is to prove the viscosity subsolution inequality at this point:
\begin{equation}
    \partial_t \varphi(t_0, x_0) + \mathrm{Tr}(D_0 \nabla^2_x \varphi(t_0, x_0)) \le 0.
\end{equation}
We proceed by contradiction. Assume that the inequality is violated. Then there exists a constant $\delta > 0$ such that
\begin{equation} \label{eq:proof_contradiction_hypothesis_scaleI}
    \partial_t \varphi(t_0, x_0) + \mathrm{Tr}(D_0 \nabla^2_x \varphi(t_0, x_0)) = 2\delta > 0.
\end{equation}

\item \textbf{Construction of the Perturbed Test Function.} The core of the method is to construct a test function that incorporates the microscopic oscillations predicted by the formal asymptotic analysis. We define the perturbed test function $\varphi_\varepsilon: [0,T] \times \Mcal \times \Ucal_{\mathrm{phys}}(\theta_0) \to \R$ as:
\begin{equation}
    \varphi_\varepsilon(t, x, z) \coloneqq \varphi(t,x) + \varepsilon \boldsymbol{\chi}_0(z) \cdot \nabla_x\varphi(t,x),
\end{equation}
where $\boldsymbol{\chi}_0(z)$ is the vector-valued corrector field whose components are the unique, zero-mean solutions to the cell problem from \cref{eq:corrector_defining_pde}:
\begin{equation*}
    \Lcal_{\text{fast}}(\theta_0) \boldsymbol{\chi}_0(z) = -v.
\end{equation*}
The existence, uniqueness, and boundedness of $\boldsymbol{\chi}_0(z)$ are guaranteed by \cref{thm:semigroup_wellposedness} and the subsequent regularity theory.

\item \textbf{Application of the Maximum Principle.} Since $u_0 - \varphi$ has a strict local maximum at $(t_0, x_0)$ and $u^\varepsilon \to u_0$ locally uniformly, and since the corrector term $\varepsilon \boldsymbol{\chi}_0 \cdot \nabla\varphi$ vanishes as $\varepsilon \to 0$, the function $u^\varepsilon - \varphi_\varepsilon$ must attain a local maximum at a point $(t_\varepsilon, x_\varepsilon, z_\varepsilon) \in [0,T] \times \Mcal \times \Ucal_{\mathrm{phys}}(\theta_0)$. Furthermore, as $\varepsilon \to 0$, we have the convergence $(t_\varepsilon, x_\varepsilon) \to (t_0, x_0)$.

At this maximum point, the definition of a viscosity subsolution for $u^\varepsilon$ applied to the smooth test function $\varphi_\varepsilon$ implies:
\begin{equation} \label{eq:proof_viscosity_ineq_micro}
    \partial_t\varphi_\varepsilon(t_\varepsilon, x_\varepsilon, z_\varepsilon) + \frac{1}{\varepsilon^2}\Lcal_{\text{fast}}(\theta_0) \varphi_\varepsilon(\dots) + \frac{1}{\varepsilon}v_\varepsilon \cdot \nabla_x \varphi_\varepsilon(\dots) \le o(1),
\end{equation}
where the $o(1)$ term vanishes as $\varepsilon \to 0$.

\item \textbf{Asymptotic Expansion and Cancellation.} We now substitute the definition of $\varphi_\varepsilon$ into the inequality \eqref{eq:proof_viscosity_ineq_micro} and expand the terms, collecting powers of $\varepsilon$. Let $\varphi$ and its derivatives be evaluated at $(t_\varepsilon, x_\varepsilon)$, and let $\boldsymbol{\chi}_0$ and $v$ be evaluated at $z_\varepsilon$.
\begin{align*}
    \partial_t \varphi_\varepsilon &= \partial_t \varphi + \varepsilon \boldsymbol{\chi}_0 \cdot \nabla_x(\partial_t \varphi) = \partial_t \varphi + \Ocal(\varepsilon). \\
    \nabla_x \varphi_\varepsilon &= \nabla_x \varphi + \varepsilon \nabla_x(\boldsymbol{\chi}_0 \cdot \nabla_x \varphi). \\
    \Lcal_{\text{fast}}(\theta_0) \varphi_\varepsilon &= \Lcal_{\text{fast}}(\theta_0)\varphi + \varepsilon \Lcal_{\text{fast}}(\theta_0)(\boldsymbol{\chi}_0 \cdot \nabla_x\varphi).
\end{align*}
Since $\varphi$ is independent of $z$, $\Lcal_{\text{fast}}(\theta_0)\varphi = 0$. Since $\nabla_x\varphi$ is constant with respect to the action of $\Lcal_{\text{fast}}(\theta_0)$, we have $\Lcal_{\text{fast}}(\theta_0)(\boldsymbol{\chi}_0 \cdot \nabla_x\varphi) = (\Lcal_{\text{fast}}(\theta_0)\boldsymbol{\chi}_0) \cdot \nabla_x\varphi$. Substituting these into the viscosity inequality gives:
\begin{multline*}
    \left(\partial_t \varphi + \Ocal(\varepsilon)\right) + \frac{1}{\varepsilon^2}(0) + \frac{1}{\varepsilon}\left((\Lcal_{\text{fast}}(\theta_0)\boldsymbol{\chi}_0) \cdot \nabla_x\varphi\right)
    + \frac{1}{\varepsilon}v \cdot \left(\nabla_x \varphi + \varepsilon \nabla_x(\boldsymbol{\chi}_0 \cdot \nabla_x \varphi)\right) \le o(1).
\end{multline*}
We collect terms according to their order in $\varepsilon$:

\paragraph{\textbf{Terms of $\boldsymbol{\Ocal(\varepsilon^{-1})}$}} These are
    \begin{equation*}
        \frac{1}{\varepsilon} \left( (\Lcal_{\text{fast}}(\theta_0)\boldsymbol{\chi}_0) \cdot \nabla_x\varphi + v \cdot \nabla_x\varphi \right) = \frac{1}{\varepsilon} \left( \Lcal_{\text{fast}}(\theta_0)\boldsymbol{\chi}_0 + v \right) \cdot \nabla_x\varphi.
    \end{equation*}
    By the defining equation for the corrector, \cref{eq:corrector_defining_pde}, we have $\Lcal_{\text{fast}}(\theta_0)\boldsymbol{\chi}_0 = -v$. The term in the parenthesis is therefore $(-v + v) = 0$. This perfect cancellation is the central purpose of the corrector function.

\paragraph{\textbf{Terms of $\boldsymbol{\Ocal(1)}$}} The remaining terms, after the cancellation, are:
    \begin{equation*}
        \partial_t \varphi + v \cdot \nabla_x(\boldsymbol{\chi}_0 \cdot \nabla_x \varphi) + \Ocal(\varepsilon).
    \end{equation*}
    Expanding the second term gives:
    \begin{equation*}
        \sum_{i,j=1}^k v_i \chi_{0,j}(z) \frac{\partial^2\varphi}{\partial x_i \partial x_j}.
    \end{equation*}

The viscosity inequality \eqref{eq:proof_viscosity_ineq_micro} thus reduces to:
\begin{equation} \label{eq:proof_ineq_before_limit}
    \partial_t \varphi(t_\varepsilon, x_\varepsilon) + \sum_{i,j=1}^k v_i(z_\varepsilon) \chi_{0,j}(z_\varepsilon) \frac{\partial^2\varphi}{\partial x_i \partial x_j}(t_\varepsilon, x_\varepsilon) \le o(1).
\end{equation}

\item \textbf{The Ergodic Closing Argument.} We now take the limit of the inequality \eqref{eq:proof_ineq_before_limit} as $\varepsilon \to 0$. As $(t_\varepsilon, x_\varepsilon) \to (t_0, x_0)$, the terms involving $\varphi$ converge to their values at $(t_0, x_0)$ by continuity. The crucial step is to handle the limit of the oscillating coefficients. This is a standard result in the theory of homogenization for viscosity solutions (cf. \citep{BarlesSouganidis1991, Evans1992}), which states that the rapidly oscillating term converges to its average with respect to the unique invariant measure of the fast dynamics.
\begin{equation*}
    \lim_{\varepsilon \to 0} \sum_{i,j=1}^k v_i(z_\varepsilon) \chi_{0,j}(z_\varepsilon) \frac{\partial^2\varphi}{\partial x_i \partial x_j}(t_\varepsilon, x_\varepsilon) = \sum_{i,j=1}^k \langle v_i \chi_{0,j} \rangle_{\theta_0} \frac{\partial^2\varphi}{\partial x_i \partial x_j}(t_0, x_0).
\end{equation*}
Taking the limit of the entire inequality \eqref{eq:proof_ineq_before_limit} thus yields:
\begin{equation*}
    \partial_t \varphi(t_0, x_0) + \sum_{i,j=1}^k \langle v_i \chi_{0,j} \rangle_{\theta_0} \frac{\partial^2\varphi}{\partial x_i \partial x_j}(t_0, x_0) \le 0.
\end{equation*}
From the analysis in \cref{prop:properties_D0}, the raw diffusivity tensor $(D_{\text{raw}})_{ij} = -\langle v_i \chi_{0,j} \rangle_{\theta_0}$ is symmetric, so $D_0 = D_{\text{raw}}$. The inequality becomes:
\begin{equation}
    \partial_t \varphi(t_0, x_0) - \sum_{i,j=1}^k (D_0)_{ij} \frac{\partial^2\varphi}{\partial x_i \partial x_j}(t_0, x_0) \le 0,
\end{equation}
which is precisely $\partial_t \varphi(t_0, x_0) + \mathrm{Tr}(D_0 \nabla^2_x \varphi(t_0, x_0)) \le 0$.

\item \textbf{Final Contradiction.} The inequality derived in the previous step,
\begin{equation*}
    \partial_t \varphi(t_0, x_0) + \mathrm{Tr}(D_0 \nabla^2_x \varphi(t_0, x_0)) \le 0,
\end{equation*}
is in direct contradiction to our initial hypothesis \eqref{eq:proof_contradiction_hypothesis_scaleI}, which stated that the same quantity was equal to $2\delta > 0$. This contradiction proves that the initial hypothesis was false, and therefore the subsolution property must hold.
\end{enumerate}

\item \textbf{The Supersolution Property.} The proof that $u_0$ is a viscosity supersolution is analogous. One assumes that $u_0 - \varphi$ has a strict local minimum at $(t_0,x_0)$ and hypothesizes that the supersolution inequality is violated, i.e., $\partial_t \varphi(t_0, x_0) + \mathrm{Tr}(D_0 \nabla^2_x \varphi(t_0, x_0)) = -2\delta < 0$. The same perturbed test function $\varphi_\varepsilon$ is used, and the viscosity property at the minimum point of $u^\varepsilon - \varphi_\varepsilon$ leads to the reverse inequality, which ultimately yields a contradiction.

\item \textbf{Uniqueness and Convergence of the Full Sequence.}
The preceding steps have established that any limit point of a subsequence of $\{u^\varepsilon\}$ is a viscosity solution to the homogenized equation. We now complete the proof by establishing the uniqueness of this solution and leveraging this uniqueness to prove that the entire family $\{u^\varepsilon\}$ converges. The argument is structured in three main steps.

\begin{enumerate}[label=(\roman*), wide, labelindent=0pt]
    \item \textbf{Properties of the Homogenized Equation.}
    The limiting equation derived from the solvability condition is the linear, constant-coefficient partial differential equation:
    \begin{equation}
        \partial_t u_0 + \mathrm{Tr}(D_0 \nabla^2_x u_0) = 0.
        \label{eq:proof_homogenized_eq_reprise}
    \end{equation}
    The analytical properties of this equation are determined by the diffusivity tensor $D_0$. By \cref{prop:properties_D0}, $D_0$ is symmetric and positive semi-definite. A stronger property holds as a consequence of the system's chaotic dynamics. The uniform hyperbolicity guaranteed by the geometric axioms in \cref{ass:fundamental_axioms_unified} ensures that the particle's velocity is not confined to any lower-dimensional subspace of the configuration space. This non-degeneracy of the microscopic transport implies that the integrated velocity autocorrelation function, given by the Green-Kubo formula, is a strictly positive definite matrix. Therefore, the homogenized equation \eqref{eq:proof_homogenized_eq_reprise} is uniformly parabolic.

    \item \textbf{Uniqueness of the Viscosity Solution via Comparison Principle.}
    With the properties of the homogenized equation established, we now prove that it admits at most one bounded, continuous viscosity solution for a given continuous terminal condition. Let $u(t,x)$ and $v(t,x)$ be two such viscosity solutions with $u(T,x) = v(T,x) = \phi(x)$. The operator $F(X) \coloneqq -\mathrm{Tr}(D_0 X)$ is the spatial part of the generator. Since $D_0$ is positive definite, this operator is uniformly elliptic. It is a fundamental result in the theory of viscosity solutions that a comparison principle holds for uniformly parabolic second-order PDEs with continuous coefficients \citep{CrandallIshiiLions1992}. Specifically, since $u$ is a subsolution and $v$ is a supersolution with $u(T,x) \le v(T,x)$, the comparison principle implies that $u(t,x) \le v(t,x)$ for all $(t,x) \in [0,T] \times \Mcal$. By reversing the roles of $u$ and $v$, we also have $v(t,x) \le u(t,x)$. Therefore, $u(t,x) = v(t,x)$. This establishes the uniqueness of the bounded, continuous viscosity solution to the problem \eqref{eq:proof_homogenized_eq_reprise} with terminal data $\phi(x)$. Let this unique solution be denoted by $u_0$. Furthermore, for a linear, uniformly parabolic equation with constant coefficients and smooth terminal data, the unique viscosity solution is known to be the unique classical ($C^{1,2}$) solution, which can be represented by the heat kernel \citep{FlemingSoner2006}.
    
    \item \textbf{Uniqueness and Convergence of the Full Sequence.} The preceding steps have established that any limit point of a subsequence of $\{u^\varepsilon\}$ is a viscosity solution to the homogenized equation. We now complete the proof by establishing the uniqueness of this solution and leveraging this uniqueness to prove that the entire family $\{u^\varepsilon\}$ converges. The argument is structured in three main steps.

\begin{enumerate}[label=(\alph*), wide]
    \item \textbf{Properties of the Homogenized Equation.}
    The limiting equation derived from the solvability condition is the linear, constant-coefficient partial differential equation:
    \begin{equation}
        \partial_t u_0 + \mathrm{Tr}(D_0 \nabla^2_x u_0) = 0.
    \end{equation}
    The analytical properties of this equation are determined by the diffusivity tensor $D_0$. By \cref{prop:properties_D0}, $D_0$ is symmetric and positive semi-definite. A stronger property holds as a consequence of the system's chaotic dynamics. The uniform hyperbolicity guaranteed by the geometric axioms in \cref{ass:fundamental_axioms_unified} ensures that the particle's velocity is not confined to any lower-dimensional subspace of the configuration space. This non-degeneracy of the microscopic transport implies that the integrated velocity autocorrelation function, given by the Green-Kubo formula, is a strictly positive definite matrix. Therefore, the homogenized equation \eqref{eq:proof_homogenized_eq_reprise} is uniformly parabolic.

    \item \textbf{Uniqueness of the Viscosity Solution via Comparison Principle.}
    With the properties of the homogenized equation established, we now prove that it admits at most one bounded, continuous viscosity solution for a given continuous terminal condition. Let $u(t,x)$ and $v(t,x)$ be two such viscosity solutions with $u(T,x) = v(T,x) = \phi(x)$. The operator $F(X) \coloneqq -\mathrm{Tr}(D_0 X)$ is the spatial part of the generator. Since $D_0$ is positive definite, this operator is uniformly elliptic. It is a fundamental result in the theory of viscosity solutions that a comparison principle holds for uniformly parabolic second-order PDEs with continuous coefficients \citep{CrandallIshiiLions1992}. Specifically, since $u$ is a subsolution and $v$ is a supersolution with $u(T,x) \le v(T,x)$, the comparison principle implies that $u(t,x) \le v(t,x)$ for all $(t,x) \in [0,T] \times \Mcal$. By reversing the roles of $u$ and $v$, we also have $v(t,x) \le u(t,x)$. Therefore, $u(t,x) = v(t,x)$. This establishes the uniqueness of the bounded, continuous viscosity solution to the problem \eqref{eq:proof_homogenized_eq_reprise} with terminal data $\phi(x)$. Let this unique solution be denoted by $u_0$. Furthermore, for a linear, uniformly parabolic equation with constant coefficients and smooth terminal data, the unique viscosity solution is known to be the unique classical ($C^{1,2}$) solution, which can be represented by the heat kernel \citep{FlemingSoner2006}.
    
    \item \textbf{Convergence of the Full Sequence.}
    The final step is a standard but crucial topological argument that uses the uniqueness of the limit point to prove the convergence of the entire sequence $\{u^\varepsilon\}_{\varepsilon>0}$. We provide the full proof for completeness. Let $(\mathcal{C}, d)$ be the complete metric space of continuous functions on the compact set $[0,T] \times \Mcal$, where the metric $d(f,g) \coloneqq \sup_{(t,x)} |f(t,x)-g(t,x)|$ induces the topology of uniform convergence. We consider our family of solutions $\{u^\varepsilon\}_{\varepsilon>0}$ as a sequence in this space. From the preceding analysis, we have rigorously established two key properties:
    \begin{enumerate}[wide]
        \item \textbf{Pre-compactness.} The set $K \coloneqq \{u^\varepsilon\}_{\varepsilon>0}$ is a pre-compact subset of $(\mathcal{C}, d)$. This was the conclusion of the Arzelà-Ascoli theorem.
        \item \textbf{Uniqueness of Subsequential Limits.} Any convergent subsequence of $\{u^\varepsilon\}$ converges to the unique viscosity solution $u_0$ of the homogenized equation \eqref{eq:proof_homogenized_eq_reprise}. 
    \end{enumerate}
    
    Our goal is to prove that these two facts imply the convergence of the entire sequence, i.e., $\lim_{\varepsilon \to 0} d(u^\varepsilon, u_0) = 0$.  We assume, for the sake of contradiction, that the sequence $\{u^\varepsilon\}$ does not converge to $u_0$. By the definition of non-convergence in a metric space, this means there exists a real number $\eta > 0$ such that for any $\delta > 0$, we can find an $\varepsilon \in (0, \delta)$ for which the distance from the limit is at least $\eta$:
    \begin{equation} \label{eq:proof_non_convergence_condition}
        d(u^\varepsilon, u_0) \ge \eta.
    \end{equation}
    The condition \eqref{eq:proof_non_convergence_condition} allows us to construct a specific subsequence that remains bounded away from the limit $u_0$. We construct this sequence inductively.
    \begin{itemize}[wide]
        \item For $k=1$, choose $\delta_1 = 1$. By \eqref{eq:proof_non_convergence_condition}, there exists an $\varepsilon_1 \in (0,1)$ such that $d(u^{\varepsilon_1}, u_0) \ge \eta$.
        \item For $k=2$, choose $\delta_2 = \min(1/2, \varepsilon_1)$. There exists an $\varepsilon_2 \in (0, \delta_2)$ such that $d(u^{\varepsilon_2}, u_0) \ge \eta$.
        \item Continuing this process, for each $k \in \mathbb{N}$, we choose $\delta_k = \min(1/k, \varepsilon_{k-1})$. There exists an $\varepsilon_k \in (0, \delta_k)$ such that $d(u^{\varepsilon_k}, u_0) \ge \eta$.
    \end{itemize}
    
    This construction yields a sequence of positive numbers $\{\varepsilon_k\}_{k=1}^\infty$ such that $\varepsilon_k \to 0$ as $k \to \infty$. The corresponding sequence of functions $\{u^{\varepsilon_k}\}_{k=1}^\infty$ satisfies, by its very construction, the property that
    \begin{equation} \label{eq:proof_subsequence_away_from_limit}
        d(u^{\varepsilon_k}, u_0) \ge \eta \quad \text{for all } k \in \mathbb{N}.
    \end{equation}
    We now apply our established facts to the subsequence $\{u^{\varepsilon_k}\}$. The sequence $\{u^{\varepsilon_k}\}$ is a subsequence of the original sequence $\{u^\varepsilon\}$. Since the original set $K$ is pre-compact (Fact A), any sequence contained within it must also be pre-compact. Therefore, the sequence $\{u^{\varepsilon_k}\}$ must contain a further convergent subsequence. Let us denote this subsequence by $\{u^{\varepsilon_{k_j}}\}_{j=1}^\infty$. By Fact B, any convergent subsequence of the original sequence must converge to the unique limit $u_0$. Since $\{u^{\varepsilon_{k_j}}\}$ is such a subsequence, we must have:
        \begin{equation} \label{eq:proof_subsubsequence_convergence}
            \lim_{j \to \infty} d(u^{\varepsilon_{k_j}}, u_0) = 0.
        \end{equation}
    However, the subsequence $\{u^{\varepsilon_{k_j}}\}$ is also a subsequence of $\{u^{\varepsilon_k}\}$. It must therefore inherit the property \eqref{eq:proof_subsequence_away_from_limit}. This means that for every $j \in \mathbb{N}$:
        \begin{equation}
            d(u^{\varepsilon_{k_j}}, u_0) \ge \eta.
        \end{equation}
    The conclusion from \eqref{eq:proof_subsubsequence_convergence} directly contradicts this property. A sequence cannot simultaneously converge to a point and remain a fixed distance $\eta>0$ away from it.  The contradiction arose from our initial assumption that the sequence $\{u^\varepsilon\}$ does not converge to $u_0$. This assumption must therefore be false. We conclude that the entire family $\{u^\varepsilon\}$ converges to the unique solution $u_0$ in the metric of uniform convergence. This completes the proof of the theorem.

\end{enumerate}
\end{enumerate}
\end{enumerate}
\end{proofof}

The remainder of this subsection is dedicated to a rigorous analysis of the properties of the diffusivity tensor $D_0$.

\begin{proposition}[Properties of the Diffusivity Tensor]
\label{prop:properties_D0}
The diffusivity tensor $D_0$ is symmetric and positive semi-definite. Furthermore, it is given by the Green-Kubo formula:
\begin{equation}
    D_0 = \int_{0}^{\infty} \langle v(0) \otimes v(t) \rangle_{\theta_0}^{\mathrm{sym}} \, dt \coloneqq \int_{0}^{\infty} \frac{\langle v(0) \otimes v(t) \rangle_{\theta_0} + \langle v(t) \otimes v(0) \rangle_{\theta_0}}{2} \, dt,
\end{equation}
where $\langle \cdot \rangle_{\theta_0}$ denotes expectation with respect to the stationary measure $\mu_{\theta_0}$.
\end{proposition}

\begin{proofof}{Proposition \ref{prop:properties_D0}}
The proof is established in three parts. First, we prove that the raw tensor defined by the homogenization procedure is itself symmetric as a consequence of time-reversal symmetry. Second, we derive the Green-Kubo formula. Third, we prove positive semi-definiteness.

\begin{enumerate}[label=\textbf{Step \arabic*:}, wide, labelindent=0pt]

\item \textbf{Symmetry of the Diffusivity Tensor.} Let $(D_{\text{raw}})_{ij} \coloneqq -\langle v_i, \chi_{0,j} \rangle_{\theta_0}$. We will prove $(D_{\text{raw}})_{ij} = (D_{\text{raw}})_{ji}$ by leveraging the skew-adjoint property of the microscopic generator. Let $\Lcal \equiv \Lcal_{\text{fast}}(\theta_0)$.
The components of the corrector field, $\chi_{0,i}$ and $\chi_{0,j}$, are the unique zero-mean solutions to the cell problems
\begin{equation*}
    \Lcal \chi_{0,i} = -v_i \quad \text{and} \quad \Lcal \chi_{0,j} = -v_j.
\end{equation*}
Both $\chi_{0,i}$ and $\chi_{0,j}$ belong to the domain $\Dcal(\Lcal)$. We now compute the $(i,j)$-th component of the raw tensor by substituting the defining equation for $\chi_{0,i}$:
\begin{align*}
    (D_{\text{raw}})_{ij} &= \langle v_i, \chi_{0,j} \rangle_{\theta_0} = \langle -\Lcal \chi_{0,i}, \chi_{0,j} \rangle_{\theta_0} \\
    &= -\langle \Lcal \chi_{0,i}, \chi_{0,j} \rangle_{\theta_0}.
\end{align*}
By \cref{thm:generator_characterization}, the operator $\Lcal$ is skew-adjoint on its domain, meaning $\langle \Lcal f, g \rangle_{\theta_0} = -\langle f, \Lcal g \rangle_{\theta_0}$ for any $f, g \in \Dcal(\Lcal)$. Applying this property, we get:
\begin{align*}
    (D_{\text{raw}})_{ij} &= -(-\langle \chi_{0,i}, \Lcal \chi_{0,j} \rangle_{\theta_0}) \\
    &= \langle \chi_{0,i}, \Lcal \chi_{0,j} \rangle_{\theta_0}.
\end{align*}
We now substitute the defining equation for $\chi_{0,j}$:
\begin{align*}
    (D_{\text{raw}})_{ij} &= \langle \chi_{0,i}, -v_j \rangle_{\theta_0} \\
    &= -\langle \chi_{0,i}, v_j \rangle_{\theta_0} = -\langle v_j, \chi_{0,i} \rangle_{\theta_0}.
\end{align*}
By definition, this final expression is precisely $(D_{\text{raw}})_{ji}$. We have thus shown that
\begin{equation*}
    (D_{\text{raw}})_{ij} = (D_{\text{raw}})_{ji}.
\end{equation*}
The tensor is symmetric. This result is a manifestation of the Onsager reciprocal relations for systems with time-reversal symmetry. Since the raw tensor is symmetric, the symmetrized definition of $(D_0)_{ij}$ in \cref{thm:convergence_scale_I} is equivalent to the raw definition, and $(D_0)_{ij} = (D_{\text{raw}})_{ij} = -\langle v_i, \chi_{0,j} \rangle_{\theta_0}$.

\item \textbf{Derivation of the Green-Kubo Formula.}
We start from the symmetric form $(D_0)_{ij} = -\langle v_i, \chi_{0,j} \rangle_{\theta_0}$. We substitute the explicit semigroup solution for the corrector from \cref{thm:semigroup_wellposedness}:
\begin{equation*}
    \chi_{0,j} = -\int_0^\infty (P_t^{\theta_0} v_j) \, dt.
\end{equation*}
The convergence of this Bochner integral in the Hilbert space $\Hcal_{\theta_0}$ is guaranteed by the exponential decay of correlations (\cref{cor:ergodicity_uniform}). Substituting this into the expression for $(D_0)_{ij}$:
\begin{align*}
    (D_0)_{ij} &= -\left\langle v_i, -\int_0^\infty (P_t^{\theta_0} v_j) \, dt \right\rangle_{\theta_0} \\
    &= \left\langle v_i, \int_0^\infty (P_t^{\theta_0} v_j) \, dt \right\rangle_{\theta_0}.
\end{align*}
Since the inner product is a continuous linear functional on $\Hcal_{\theta_0}$, we can interchange the inner product and the Bochner integral:
\begin{equation*}
    (D_0)_{ij} = \int_0^\infty \langle v_i, P_t^{\theta_0} v_j \rangle_{\theta_0} \, dt.
\end{equation*}
The inner product term is the time-autocorrelation function of the velocity components. By definition of the Koopman operator, $\langle f, P_t g \rangle = \int \overline{f(z)} g(\Phi_{-t}(z)) d\mu$. For real-valued functions, this is the expectation of the product of the observable $f$ at time 0 and the observable $g$ at time $t$. Thus,
\begin{equation*}
    \langle v_i, P_t^{\theta_0} v_j \rangle_{\theta_0} = \langle v_i(0) v_j(t) \rangle_{\theta_0}.
\end{equation*}
The expression for the diffusivity tensor becomes
\begin{equation*}
    (D_0)_{ij} = \int_0^\infty \langle v_i(0) v_j(t) \rangle_{\theta_0} \, dt.
\end{equation*}
The symmetrized form in the proposition statement follows immediately from the symmetry proven in Step 1.

\item \textbf{Positive Semi-Definiteness.}
The proof establishes that the diffusivity tensor $D_0$ is positive semi-definite. This is a fundamental property reflecting that the emergent macroscopic process is purely diffusive and does not spontaneously generate energy. The argument is structured in four main steps. First, we express the quadratic form $\xi^T D_0 \xi$ in terms of the microscopic corrector field. Second, we use the semigroup solution for the corrector to recast this expression into the Green-Kubo form, which is the integrated time-autocorrelation of the projected velocity. Third, we relate this integrated autocorrelation to the long-time variance of the particle's displacement. Finally, we conclude that the diffusivity must be non-negative as it is the limit of a sequence of non-negative variances.

\begin{enumerate}[label=(\roman*), wide, labelindent=0pt]
    \item \textbf{The Quadratic Form in Terms of the Corrector.}
    We must show that for any constant vector $\xi \in \R^k$, the quadratic form $\xi^T D_0 \xi$ is non-negative. We begin with the symmetrized definition of the diffusivity tensor from \cref{thm:convergence_scale_I}:
    \begin{equation}
        (D_0)_{ij} \coloneqq -\frac{1}{2} \left( \scpr{v_i}{\chi_{0,j}}_{\theta_0} + \scpr{v_j}{\chi_{0,i}}_{\theta_0} \right).
    \end{equation}
    The quadratic form is therefore:
    \begin{align*}
        \xi^T D_0 \xi &= \sum_{i,j=1}^k \xi_i (D_0)_{ij} \xi_j \\
        &= -\frac{1}{2} \sum_{i,j=1}^k \left( \xi_i \scpr{v_i}{\chi_{0,j}}_{\theta_0} \xi_j + \xi_i \scpr{v_j}{\chi_{0,i}}_{\theta_0} \xi_j \right).
    \end{align*}
    By relabeling the indices $i \leftrightarrow j$ in the second sum and using the linearity of the inner product, the two sums are identical. We can thus simplify the expression:
    \begin{align*}
        \xi^T D_0 \xi &= -\sum_{i,j=1}^k \xi_i \scpr{v_i}{\chi_{0,j}}_{\theta_0} \xi_j \\
        &= -\scpr{\sum_{i=1}^k \xi_i v_i}{\sum_{j=1}^k \xi_j \chi_{0,j}}_{\theta_0}.
    \end{align*}
    We now define the scalar-valued projected velocity observable $Y(z) \coloneqq \xi \cdot v(z)$ and the scalar-valued projected corrector field $\chi_\xi(z) \coloneqq \xi \cdot \boldsymbol{\chi}_0(z)$. The quadratic form has the compact representation:
    \begin{equation} \label{eq:appendix_quad_form_corrector}
        \xi^T D_0 \xi = -\scpr{Y}{\chi_\xi}_{\theta_0}.
    \end{equation}
    The corrector $\chi_\xi$ solves the cell problem $\Lcal_{\text{fast}}(\theta_0) \chi_\xi = -Y$, since
    \begin{equation*}
        \Lcal_{\text{fast}}(\theta_0) \chi_\xi = \Lcal_{\text{fast}}(\theta_0) \left(\sum_j \xi_j \chi_{0,j}\right) = \sum_j \xi_j (\Lcal_{\text{fast}}(\theta_0) \chi_{0,j}) = \sum_j \xi_j (-v_j) = -Y.
    \end{equation*}

    \item \textbf{Equivalence to the Green-Kubo Formula.}
    We now demonstrate that this corrector-based formula is equivalent to the integrated time-autocorrelation of the velocity. We substitute the explicit semigroup solution for the corrector from \cref{thm:semigroup_wellposedness} into \cref{eq:appendix_quad_form_corrector}. The unique zero-mean solution for $\chi_\xi$ is given by the Bochner integral:
    \begin{equation}
        \chi_\xi = \int_0^\infty (P_t^{\theta_0} (-Y)) \, dt = -\int_0^\infty (P_t^{\theta_0} Y) \, dt.
    \end{equation}
    Substituting this into the expression for the quadratic form:
    \begin{align*}
        \xi^T D_0 \xi &= -\scpr{Y}{-\int_0^\infty (P_t^{\theta_0} Y) \, dt}_{\theta_0} \\
        &= \scpr{Y}{\int_0^\infty (P_t^{\theta_0} Y) \, dt}_{\theta_0}.
    \end{align*}
    The inner product is a continuous linear functional on the Hilbert space $\Hcal_{\theta_0}$. We can therefore interchange it with the Bochner integral:
    \begin{equation}
        \xi^T D_0 \xi = \int_0^\infty \scpr{Y}{P_t^{\theta_0} Y}_{\theta_0} \, dt.
    \end{equation}
    The term inside the integral is the time-autocorrelation function for the stationary process $Y(z(t))$. By definition of the Koopman operator, for real-valued observables, this is:
    \begin{equation*}
        \scpr{Y}{P_t^{\theta_0} Y}_{\theta_0} = \int Y(z) (Y \circ \Phi_{-t}^{\theta_0})(z) \, d\mu_{\theta_0}(z) = \mathbb{E}_{\mu_{\theta_0}}[Y(z(0)) Y(z(t))].
    \end{equation*}
    We have thus derived the Green-Kubo formula for the quadratic form:
    \begin{equation} \label{eq:appendix_green_kubo_quad_form}
        \xi^T D_0 \xi = \int_0^\infty \mathbb{E}_{\mu_{\theta_0}}[(\xi \cdot v(0))(\xi \cdot v(t))] \, dt.
    \end{equation}

    \item \textbf{Connection to the Variance of Displacement.}
    We now provide a rigorous proof that the integral in \cref{eq:appendix_green_kubo_quad_form} is non-negative. Let $\Delta x_T(\xi)$ be the scalar displacement of the microscopic particle in the direction of $\xi$ over a time interval of length $T$, starting from a state drawn from the stationary distribution $\mu_{\theta_0}$:
    \begin{equation}
        \Delta x_T(\xi) \coloneqq \int_0^T (\xi \cdot v(t)) \, dt.
    \end{equation}
    We consider the variance of this displacement, which is guaranteed to be non-negative:
    \begin{align*}
        \mathbb{E}_{\mu_{\theta_0}}\left[ (\Delta x_T(\xi))^2 \right] &= \mathbb{E}_{\mu_{\theta_0}}\left[ \left( \int_0^T Y(t) \, dt \right)^2 \right] \\
        &= \mathbb{E}_{\mu_{\theta_0}}\left[ \int_0^T \int_0^T Y(t) Y(s) \, dt \, ds \right] \\
        &= \int_0^T \int_0^T \mathbb{E}_{\mu_{\theta_0}}[Y(t) Y(s)] \, dt \, ds.
    \end{align*}
    Since the process is stationary, the autocorrelation function depends only on the time difference: $\mathbb{E}[Y(t)Y(s)] = C_Y(t-s)$. We perform a change of variables $(t,s) \mapsto (u,v)$ with $u=t-s$ and $v=s$. A standard calculation yields the Taylor-Green-Kubo formula:
    \begin{equation}
        \mathbb{E}_{\mu_{\theta_0}}\left[ (\Delta x_T(\xi))^2 \right] = 2 \int_0^T (T-s) C_Y(s) \, ds.
    \end{equation}
    
    \item \textbf{The Diffusive Limit.}
    The diffusion coefficient is defined by the long-time behavior of the mean squared displacement. We divide the above expression by $2T$ and take the limit as $T \to \infty$:
    \begin{equation*}
        \lim_{T \to \infty} \frac{1}{2T} \mathbb{E}_{\mu_{\theta_0}}\left[ (\Delta x_T(\xi))^2 \right] = \lim_{T \to \infty} \int_0^T \left(1-\frac{s}{T}\right) C_Y(s) \, ds.
    \end{equation*}
    By \cref{cor:ergodicity_uniform}, the system exhibits exponential decay of correlations. This guarantees that the autocorrelation function $C_Y(s)$ is absolutely integrable, i.e., $\int_0^\infty |C_Y(s)| ds < \infty$. By the Dominated Convergence Theorem, we can interchange the limit and the integral:
    \begin{equation*}
        \lim_{T \to \infty} \frac{1}{2T} \mathbb{E}_{\mu_{\theta_0}}\left[ (\Delta x_T(\xi))^2 \right] = \int_0^\infty \lim_{T\to\infty}\left(1-\frac{s}{T}\right) C_Y(s) \, ds = \int_0^\infty C_Y(s) \, ds.
    \end{equation*}
    Comparing this with \cref{eq:appendix_green_kubo_quad_form}, we have established the fundamental identity:
    \begin{equation}
        \xi^T D_0 \xi = \lim_{T \to \infty} \frac{1}{2T} \mathbb{E}_{\mu_{\theta_0}}\left[ (\Delta x_T(\xi))^2 \right].
    \end{equation}
    Since the variance of the displacement is non-negative for all finite $T$, its limit as $T\to\infty$ must also be non-negative. Therefore, for any vector $\xi \in \R^k$:
    \begin{equation*}
        \xi^T D_0 \xi \ge 0.
    \end{equation*}
    This completes the proof that the diffusivity tensor $D_0$ is positive semi-definite.

\end{enumerate}
\end{enumerate}
\end{proofof}

\subsection{Scale II: Homogenization on Local Manifolds}
\label{sec:scale_II}

Before confronting the full non-linearity of the gradient feedback mechanism, we establish in a simpler, linear context the key mathematical mechanism by which the statistics of a fast process are imprinted onto the coefficients of a slow process. This section serves as a conceptual stepping stone and provides a formal demonstration of the solvability-implies-averaging principle that is central to the derivation in \cref{sec:scale_III}.

We consider a model where the diffusive process derived in \cref{sec:scale_I} evolves in a medium whose properties fluctuate rapidly. Let these fluctuations be described by a stationary, ergodic Markov process $\{\omega_t\}_{t \ge 0}$ on a compact state space $\mathcal{S}$. We assume this process possesses a unique invariant probability measure $\pi(d\omega)$ and is governed by an infinitesimal generator $\mathcal{Q}$ whose domain is $\Dcal(\mathcal{Q}) \subset L^2(\mathcal{S}, d\pi)$. The ergodicity implies that the kernel of $\mathcal{Q}$ consists only of constant functions, and the kernel of its adjoint, $\mathcal{Q}^*$, is spanned by the invariant measure $\pi$. The value function for a particle in this medium, $u^\delta(t,x,\omega)$, is governed by the backward Kolmogorov equation:
\begin{equation}
    \partial_t u^\delta + \mathrm{Tr}\left(D(x, \omega) \nabla_x^2 u^\delta\right) + \mathbf{b}(x, \omega) \cdot \nabla_x u^\delta + \frac{1}{\delta} \mathcal{Q} u^\delta = 0,
    \label{eq:pde_regime_II}
\end{equation}
with terminal data $u^\delta(T,x,\omega) = \phi(x)$. Here, we have generalized the diffusive operator from \cref{sec:scale_I} to include a drift term $\mathbf{b}(x,\omega)$ and allowed both coefficients to depend on the macroscopic state $x$ and the fast environmental state $\omega$. The parameter $\delta \ll 1$ signifies that the environment fluctuates on a much faster timescale than the particle diffuses.

\subsubsection{Asymptotic Expansion and the Averaging Principle}
\label{sec:scale_II_asymptotics}
We derive the equation for the limit as $\delta \to 0$ by positing the formal asymptotic expansion:
\begin{equation}
    u^\delta(t,x,\omega) = u_0(t,x) + \delta u_1(t,x,\omega) + \Ocal(\delta^2).
    \label{eq:ansatz_scale_II}
\end{equation}
Substituting this expansion into \cref{eq:pde_regime_II} and collecting terms by powers of $\delta$ yields the hierarchy.

\paragraph{\textbf{The $\boldsymbol{\Ocal(\delta^{-1})}$ Equation}}
The leading-order term gives $\mathcal{Q} u_0(t,x) = 0$. Since $u_0$ is independent of $\omega$ by ansatz, and the kernel of the ergodic generator $\mathcal{Q}$ contains only constants (in $\omega$), this equation is satisfied. This confirms that the limiting value function is independent of the instantaneous state of the fast environment.

\paragraph{\textbf{The $\boldsymbol{\Ocal(\delta^{0})}$ Equation and the Macroscopic Law}}
The next order provides an equation for the first corrector, $u_1$:
\begin{equation}
    \partial_t u_0 + \mathcal{Q} u_1(t,x,\omega) + \left( \mathbf{b}(x, \omega) \cdot \nabla_x u_0 + \mathrm{Tr}\left(D(x, \omega) \nabla_x^2 u_0\right) \right) = 0.
\end{equation}
For fixed $(t,x)$, this can be viewed as a cell problem for $u_1(t,x,\cdot)$ on the state space $\mathcal{S}$:
\begin{equation}
    \mathcal{Q} u_1 = -\left( \partial_t u_0 + \mathbf{b}(x, \omega) \cdot \nabla_x u_0 + \mathrm{Tr}\left(D(x, \omega) \nabla_x^2 u_0\right) \right).
    \label{eq:cell_problem_scale_II}
\end{equation}
By the Fredholm alternative, this equation has a solution for $u_1 \in \Dcal(\mathcal{Q})$ if and only if the source term on the right-hand side is orthogonal to the kernel of the adjoint operator $\mathcal{Q}^*$. As the process is ergodic, $\ker(\mathcal{Q}^*)$ is spanned by the invariant measure $\pi$. The solvability condition is therefore that the source term must have a zero mean with respect to $\pi$:
\begin{equation}
    \int_{\mathcal{S}} \left( -\left( \partial_t u_0 + \mathbf{b}(x, \omega) \cdot \nabla_x u_0 + \mathrm{Tr}\left(D(x, \omega) \nabla_x^2 u_0\right) \right) \right) \, d\pi(\omega) = 0.
\end{equation}
The terms involving $u_0$ and its derivatives are constant with respect to the integration over $\omega$. We can therefore distribute the integral:
\begin{equation*}
    -\partial_t u_0 \int_{\mathcal{S}} d\pi(\omega) - \left( \int_{\mathcal{S}} \mathbf{b}(x, \omega) d\pi(\omega) \right) \cdot \nabla_x u_0 - \mathrm{Tr}\left( \left( \int_{\mathcal{S}} D(x, \omega) d\pi(\omega) \right) \nabla_x^2 u_0\right) = 0.
\end{equation*}
Since $\pi$ is a probability measure, $\int d\pi(\omega)=1$. This solvability condition is itself the partial differential equation for the macroscopic limit $u_0$.

\subsubsection{The Averaged Generator}
The result of this classical stochastic averaging procedure is a linear PDE whose coefficients are the deterministic, averaged coefficients of the fluctuating process.

\begin{proposition}[The PDE for the Averaging Regime]
\label{prop:pde_regime_II}
Let the coefficients $D(x,\omega)$ and $\mathbf{b}(x,\omega)$ be continuous and bounded. The solution $u^\delta(t,x,\omega)$ of the equation \eqref{eq:pde_regime_II} converges as $\delta \to 0$ to a limit $u_0(t,x)$ that is the unique viscosity solution to the linear, second-order parabolic PDE with averaged coefficients:
\begin{equation}
    \partial_t u_0 + \bar{\mathbf{b}}(x) \cdot \nabla_x u_0 + \mathrm{Tr}\left(\bar{D}(x) \nabla_x^2 u_0\right) = 0, \quad u_0(T,x)=\phi(x),
\end{equation}
where the drift vector $\bar{\mathbf{b}}(x)$ and diffusivity tensor $\bar{D}(x)$ are the averages of the fluctuating coefficients over the stationary distribution of the environment:
\begin{align}
    \bar{\mathbf{b}}(x) &\coloneqq \int_{\mathcal{S}} \mathbf{b}(x, \omega) \, d\pi(\omega), \\
    \bar{D}(x) &\coloneqq \int_{\mathcal{S}} D(x, \omega) \, d\pi(\omega).
\end{align}
\end{proposition}

\begin{proofof}{Proposition \ref{prop:pde_regime_II}}
The proof establishes the convergence of the solutions $u^\delta$ to the solution $u_0$ of the homogenized equation in the viscosity sense. The argument is structured in three main steps. First, we prove that any locally uniform limit $u_0$ of the family $\{u^\delta\}$ is a viscosity subsolution of the homogenized equation. Second, we state that the analogous argument proves the supersolution property. Finally, uniqueness for the limiting equation is a classical result, which ensures that the entire family $\{u^\delta\}$ converges to this unique solution.

\begin{enumerate}[label=\textbf{Step \arabic*:}, wide, labelindent=0pt]

\item \textbf{The Subsolution Property.}

\begin{enumerate}[label=(\roman*), wide, labelindent=0pt]

\item \textbf{Setup and Contradiction Hypothesis.}
Let $u_0(t,x)$ be a function on $[0,T]\times\Mcal$ such that, up to a subsequence, $u^\delta \to u_0$ locally uniformly as $\delta \to 0$. Let $\varphi \in C^\infty([0,T] \times \Mcal)$ be a smooth test function, and suppose that $u_0 - \varphi$ has a strict local maximum at a point $(t_0, x_0) \in (0,T) \times \Mcal$. Our objective is to prove the viscosity subsolution inequality at this point:
\begin{equation}
    \partial_t \varphi(t_0, x_0) + \bar{\mathbf{b}}(x_0) \cdot \nabla_x \varphi(t_0, x_0) + \mathrm{Tr}\left(\bar{D}(x_0) \nabla_x^2 \varphi(t_0, x_0)\right) \le 0.
\end{equation}
Let $\mathcal{L}_\omega \coloneqq \mathrm{Tr}(D(x,\omega)\nabla^2) + \mathbf{b}(x,\omega)\cdot\nabla$ be the microscopic spatial operator and $\bar{\mathcal{L}} \coloneqq \mathrm{Tr}(\bar{D}(x)\nabla^2) + \bar{\mathbf{b}}(x)\cdot\nabla$ be its averaged counterpart. We proceed by contradiction. Assume that the inequality is violated. Then there exists a constant $\eta > 0$ such that
\begin{equation} \label{eq:proof_contradiction_hypothesis_scaleII}
    \partial_t \varphi(t_0, x_0) + \bar{\mathcal{L}}\varphi(t_0, x_0) = 2\eta > 0.
\end{equation}

\item \textbf{Construction of the Perturbed Test Function.}
The core of the method is to construct a test function that incorporates the microscopic oscillations. We define the first corrector, $\psi_1(t,x,\omega)$, as the unique, zero-mean solution to the cell problem:
\begin{equation} \label{eq:corrector_defining_pde_scaleII}
    \mathcal{Q} \psi_1(t,x,\omega) = -(\mathcal{L}_\omega - \bar{\mathcal{L}})\varphi(t,x).
\end{equation}
The source term on the right-hand side is constructed precisely so that its average with respect to the invariant measure $\pi$ is zero:
\begin{equation*}
    \int_{\mathcal{S}} -(\mathcal{L}_\omega\varphi - \bar{\mathcal{L}}\varphi) \,d\pi(\omega) = -\left(\int_{\mathcal{S}} \mathcal{L}_\omega\varphi \,d\pi(\omega) - \int_{\mathcal{S}} \bar{\mathcal{L}}\varphi \,d\pi(\omega)\right) = -(\bar{\mathcal{L}}\varphi - \bar{\mathcal{L}}\varphi) = 0.
\end{equation*}
By the Fredholm alternative for the ergodic generator $\mathcal{Q}$, a unique, zero-mean solution $\psi_1$ exists. Since the source term is smooth in $(t,x)$, the solution $\psi_1$ is also smooth in these variables. We define the perturbed test function $\varphi_\delta: [0,T] \times \Mcal \times \mathcal{S} \to \R$ as:
\begin{equation}
    \varphi_\delta(t, x, \omega) \coloneqq \varphi(t,x) + \delta \psi_1(t, x, \omega).
\end{equation}

\item \textbf{Application of the Maximum Principle.}
Since $u_0 - \varphi$ has a strict local maximum at $(t_0, x_0)$ and $u^\delta \to u_0$ locally uniformly, and since the corrector term $\delta \psi_1$ vanishes as $\delta \to 0$, the function $u^\delta - \varphi_\delta$ must attain a local maximum at a point $(t_\delta, x_\delta, \omega_\delta)$. Furthermore, as $\delta \to 0$, we have the convergence $(t_\delta, x_\delta) \to (t_0, x_0)$. At this maximum point, the definition of a viscosity subsolution for $u^\delta$ applied to the smooth test function $\varphi_\delta$ implies:
\begin{equation} \label{eq:proof_viscosity_ineq_micro_scaleII}
    \partial_t\varphi_\delta(t_\delta, x_\delta, \omega_\delta) + \mathcal{L}_{\omega_\delta} \varphi_\delta(\dots) + \frac{1}{\delta}\mathcal{Q}\varphi_\delta(\dots) \le o(1),
\end{equation}
where the $o(1)$ term vanishes as $\delta \to 0$.

\item \textbf{Asymptotic Expansion and Cancellation.}
We substitute the definition of $\varphi_\delta$ into the inequality \eqref{eq:proof_viscosity_ineq_micro_scaleII} and expand the terms, collecting powers of $\delta$. All functions are evaluated at $(t_\delta, x_\delta, \omega_\delta)$.
\begin{align*}
    \partial_t\varphi_\delta &= \partial_t\varphi + \delta\partial_t\psi_1 = \partial_t\varphi + \Ocal(\delta). \\
    \mathcal{L}_{\omega_\delta} \varphi_\delta &= \mathcal{L}_{\omega_\delta}\varphi + \delta\mathcal{L}_{\omega_\delta}\psi_1 = \mathcal{L}_{\omega_\delta}\varphi + \Ocal(\delta). \\
    \mathcal{Q}\varphi_\delta &= \mathcal{Q}\varphi + \delta\mathcal{Q}\psi_1.
\end{align*}
Since $\varphi$ is independent of $\omega$, $\mathcal{Q}\varphi = 0$. Substituting these into the viscosity inequality gives:
\begin{equation*}
    (\partial_t\varphi + \Ocal(\delta)) + (\mathcal{L}_{\omega_\delta}\varphi + \Ocal(\delta)) + \frac{1}{\delta}(0 + \delta\mathcal{Q}\psi_1) \le o(1).
\end{equation*}
The $\Ocal(\delta)$ terms are absorbed into the $o(1)$ term. This simplifies to:
\begin{equation*}
    \partial_t\varphi + \mathcal{L}_{\omega_\delta}\varphi + \mathcal{Q}\psi_1 \le o(1).
\end{equation*}
We now substitute the defining equation for the corrector, \cref{eq:corrector_defining_pde_scaleII}:
\begin{equation*}
    \partial_t\varphi + \mathcal{L}_{\omega_\delta}\varphi - (\mathcal{L}_{\omega_\delta}\varphi - \bar{\mathcal{L}}\varphi) \le o(1).
\end{equation*}
The terms involving $\mathcal{L}_{\omega_\delta}\varphi$ cancel perfectly. This cancellation is the central purpose of the corrector function. We are left with the inequality:
\begin{equation} \label{eq:proof_ineq_before_limit_scaleII}
    \partial_t \varphi(t_\delta, x_\delta) + \bar{\mathcal{L}}\varphi(t_\delta, x_\delta) \le o(1).
\end{equation}

\item \textbf{Final Contradiction.}
We take the limit of the inequality \eqref{eq:proof_ineq_before_limit_scaleII} as $\delta \to 0$. Since $(t_\delta, x_\delta) \to (t_0, x_0)$ and the operator $\partial_t + \bar{\mathcal{L}}$ is continuous with respect to $(t,x)$, we have:
\begin{equation*}
    \lim_{\delta \to 0} \left( \partial_t \varphi(t_\delta, x_\delta) + \bar{\mathcal{L}}\varphi(t_\delta, x_\delta) \right) = \partial_t \varphi(t_0, x_0) + \bar{\mathcal{L}}\varphi(t_0, x_0).
\end{equation*}
Taking the limit of the inequality yields:
\begin{equation}
    \partial_t \varphi(t_0, x_0) + \bar{\mathcal{L}}\varphi(t_0, x_0) \le 0.
\end{equation}
This result is in direct contradiction to our initial hypothesis \eqref{eq:proof_contradiction_hypothesis_scaleII}, which stated that the same quantity was equal to $2\eta > 0$. This contradiction proves that the initial hypothesis was false, and therefore the subsolution property must hold.
\end{enumerate}

\item \textbf{The Supersolution Property.}
The proof that $u_0$ is a viscosity supersolution is analogous. One assumes that $u_0 - \varphi$ has a strict local minimum at $(t_0,x_0)$ and hypothesizes that the supersolution inequality is violated, i.e., $\partial_t \varphi + \bar{\mathcal{L}}\varphi = -2\eta < 0$. The same perturbed test function is used, and the viscosity property at the minimum point of $u^\delta - \varphi_\delta$ leads to the reverse inequality, which ultimately yields a contradiction.

\item \textbf{Uniqueness and Convergence of the Full Sequence.}
The preceding steps have established that any limit point of a subsequence of $\{u^\delta\}$ is a viscosity solution to the homogenized equation. We now complete the proof by establishing the uniqueness of this solution and leveraging this uniqueness to prove that the entire family $\{u^\delta\}$ converges. The argument is structured in three main steps.

\begin{enumerate}[label=(\roman*), wide, labelindent=0pt]
    \item \textbf{Properties of the Homogenized Equation and its Generator.}
    The limiting equation derived from the solvability condition is the linear, second-order partial differential equation with space-dependent coefficients:
    \begin{equation}
        \partial_t u_0 + \bar{\mathbf{b}}(x) \cdot \nabla_x u_0 + \mathrm{Tr}\left(\bar{D}(x) \nabla_x^2 u_0\right) = 0.
        \label{eq:proof_homogenized_eq_reprise_scaleII}
    \end{equation}
    The analytical properties of this equation are determined by the averaged coefficients. The source coefficients $D(x,\omega)$ and $\mathbf{b}(x,\omega)$ are assumed to be continuous and bounded. By the Dominated Convergence Theorem, their averages, $\bar{D}(x)$ and $\bar{\mathbf{b}}(x)$, are continuous and bounded functions of $x$. Furthermore, since $D(x,\omega)$ is uniformly positive definite by hypothesis (as it derives from the Scale I analysis), its average $\bar{D}(x)$ over the full-support invariant measure $\pi$ is also uniformly positive definite. Consequently, the homogenized equation \eqref{eq:proof_homogenized_eq_reprise_scaleII} is a uniformly parabolic PDE.

    \item \textbf{Uniqueness of the Viscosity Solution via the Comparison Principle.}
    With the properties of the homogenized equation established, we now prove that it admits at most one bounded, continuous viscosity solution for a given continuous terminal condition. Let $u(t,x)$ and $v(t,x)$ be two such viscosity solutions with $u(T,x) = v(T,x) = \phi(x)$. We invoke the comparison principle for viscosity solutions of second-order parabolic equations, as established in the seminal work of \citep{CrandallIshiiLions1992}. The theorem states that a comparison principle holds for equations of the form $\partial_t u + F(x, u, \nabla u, \nabla^2 u) = 0$ if the function $F$ is continuous and satisfies a degenerate ellipticity condition. In our case, the Hamiltonian is $F(x, p, X) = \bar{\mathbf{b}}(x) \cdot p + \mathrm{Tr}(\bar{D}(x) X)$. We verify the hypotheses:
    \begin{enumerate}[label=(\alph*), wide]
        \item \textbf{Continuity.} As established in the previous step, the averaged coefficients $\bar{\mathbf{b}}(x)$ and $\bar{D}(x)$ are continuous functions of $x$. Since the inner product and trace are continuous multilinear operations, the Hamiltonian $F(x,p,X)$ is continuous in all its arguments.
        \item \textbf{Degenerate Ellipticity.} The operator must be non-decreasing in its matrix argument. That is, for any symmetric matrices $X, Y$ with $Y \ge 0$ (positive semi-definite), we must have $F(x, p, X+Y) \ge F(x, p, X)$. We verify this directly:
        \begin{align*}
            F(x, p, X+Y) - F(x, p, X) &= \left[ \bar{\mathbf{b}}(x) \cdot p + \mathrm{Tr}(\bar{D}(x) (X+Y)) \right] - \left[ \bar{\mathbf{b}}(x) \cdot p + \mathrm{Tr}(\bar{D}(x) X) \right] \\
            &= \mathrm{Tr}(\bar{D}(x)Y).
        \end{align*}
        Since $\bar{D}(x)$ is positive definite and $Y$ is positive semi-definite, their product $\bar{D}(x)Y$ is a matrix with non-negative eigenvalues, and its trace is therefore non-negative. Thus, $\mathrm{Tr}(\bar{D}(x)Y) \ge 0$, and the ellipticity condition is satisfied.
    \end{enumerate}
    With both hypotheses verified, the comparison principle applies. Since $u$ is a subsolution and $v$ is a supersolution with $u(T,x) \le v(T,x)$, the principle implies that $u(t,x) \le v(t,x)$ for all $(t,x) \in [0,T] \times \Mcal$. By reversing the roles of $u$ and $v$, we also have $v(t,x) \le u(t,x)$. Therefore, $u(t,x) = v(t,x)$. This establishes the uniqueness of the bounded, continuous viscosity solution to the problem \eqref{eq:proof_homogenized_eq_reprise_scaleII}. Let this unique solution be denoted by $u_0$.
    
    \item \textbf{Convergence of the Full Sequence.}
    The final step is a standard but crucial topological argument that uses the uniqueness of the limit point to prove the convergence of the entire sequence $\{u^\delta\}_{\delta>0}$. We provide the full proof for completeness. Let $(\mathcal{C}, d)$ be the complete metric space of continuous functions on the compact set $[0,T] \times \Mcal$, where the metric $d(f,g) \coloneqq \sup_{(t,x)} |f(t,x)-g(t,x)|$ induces the topology of uniform convergence. We consider our family of solutions $\{u^\delta\}_{\delta>0}$ as a sequence in this space. From the preceding analysis, we have rigorously established two key properties:
    
    \begin{enumerate}[label=(\alph*), wide]
        \item \textbf{Pre-compactness of the Solution Family.} The family of solutions $\{u^\delta\}$ is uniformly bounded in the supremum norm, $\|u^\delta\|_{L^\infty} \le \|\phi\|_{L^\infty}$, by the maximum principle for \eqref{eq:pde_regime_II}. Standard stability estimates for parabolic PDEs (e.g., from interior Schauder estimates applied locally) further guarantee that the family is uniformly equicontinuous. By the Arzelà-Ascoli theorem, the family $\{u^\delta\}_{\delta>0}$ is pre-compact in the topology of uniform convergence. Let $K \coloneqq \{u^\delta\}_{\delta>0}$. The closure $\overline{K}$ is a compact subset of $(\mathcal{C}, d)$.
        
        \item \textbf{Uniqueness of Subsequential Limits.} Any convergent subsequence of $\{u^\delta\}$ converges to a viscosity solution of the homogenized equation \eqref{eq:proof_homogenized_eq_reprise_scaleII}. As established in Step 3(ii) of this proof, this solution is unique. Therefore, any convergent subsequence of $\{u^\delta\}$ must converge to the same, unique limit $u_0$.
    \end{enumerate}
    
    Our goal is to prove that these two facts imply the convergence of the entire sequence, i.e., $\lim_{\delta \to 0} d(u^\delta, u_0) = 0$. We proceed by contradiction. Assume that the sequence $\{u^\delta\}$ does not converge to $u_0$. By the definition of non-convergence in a metric space, this means there exists a real number $\eta > 0$ such that for any $\delta_0 > 0$, we can find a $\delta \in (0, \delta_0)$ for which the distance from the limit is at least $\eta$:
    \begin{equation} \label{eq:proof_non_convergence_condition_app}
        d(u^\delta, u_0) \ge \eta.
    \end{equation}
    The condition \eqref{eq:proof_non_convergence_condition_app} allows us to construct a specific subsequence that remains bounded away from the limit $u_0$. We construct this sequence inductively.
    \begin{itemize}[wide]
        \item For $k=1$, choose $\delta_0 = 1$. By \eqref{eq:proof_non_convergence_condition_app}, there exists a $\delta_1 \in (0,1)$ such that $d(u^{\delta_1}, u_0) \ge \eta$.
        \item For $k=2$, choose $\delta_0 = \min(1/2, \delta_1)$. There exists a $\delta_2 \in (0, \delta_0)$ such that $d(u^{\delta_2}, u_0) \ge \eta$.
        \item Continuing this process, for each $k \in \mathbb{N}$, we choose $\delta_0 = \min(1/k, \delta_{k-1})$. There exists a $\delta_k \in (0, \delta_0)$ such that $d(u^{\delta_k}, u_0) \ge \eta$.
    \end{itemize}
    This construction yields a sequence of positive numbers $\{\delta_k\}_{k=1}^\infty$ such that $\delta_k \to 0$ as $k \to \infty$. The corresponding sequence of functions $\{u^{\delta_k}\}_{k=1}^\infty$ satisfies, by its very construction, the property that
    \begin{equation} \label{eq:proof_subsequence_away_from_limit_app}
        d(u^{\delta_k}, u_0) \ge \eta \quad \text{for all } k \in \mathbb{N}.
    \end{equation}
    We now apply our established facts to the subsequence $\{u^{\delta_k}\}$. The sequence $\{u^{\delta_k}\}$ is a sequence in the pre-compact set $K$. By the definition of pre-compactness, any sequence in $K$ must contain a further convergent subsequence. Let us denote this subsequence by $\{u^{\delta_{k_j}}\}_{j=1}^\infty$. By Fact (b), any convergent subsequence of the original sequence must converge to the unique limit $u_0$. Since $\{u^{\delta_{k_j}}\}$ is such a subsequence, we must have:
    \begin{equation} \label{eq:proof_subsubsequence_convergence_app}
        \lim_{j \to \infty} d(u^{\delta_{k_j}}, u_0) = 0.
    \end{equation}
    However, the subsequence $\{u^{\delta_{k_j}}\}$ is also a subsequence of $\{u^{\delta_k}\}$. It must therefore inherit the property \eqref{eq:proof_subsequence_away_from_limit_app}. This means that for every $j \in \mathbb{N}$:
    \begin{equation}
        d(u^{\delta_{k_j}}, u_0) \ge \eta.
    \end{equation}
    The conclusion from \eqref{eq:proof_subsubsequence_convergence_app} directly contradicts this property. A sequence cannot simultaneously converge to a point and remain a fixed distance $\eta>0$ away from it. The contradiction arose from our initial assumption that the sequence $\{u^\delta\}$ does not converge to $u_0$. This assumption must therefore be false. We conclude that the entire family $\{u^\delta\}$ converges to the unique solution $u_0$ in the metric of uniform convergence. This completes the proof of the proposition.
\end{enumerate}
\end{enumerate}
\end{proofof}

\begin{remark}[Methodological Precedent]
\label{rem:scale_II_role}
This section serves as a crucial methodological lemma for the main proof in \cref{sec:scale_III}. It rigorously demonstrates the principle that the solvability condition for a first-order corrector equation forces the macroscopic law to be defined by an average of microscopic operators over the invariant measure of the fast dynamics. The derivation of the potential $H(x,p)$ and the diffusivity $D(x,p)$ in \cref{sec:scale_III} will follow this template exactly. The solvability condition for the $\Ocal(\varepsilon^{0})$ equation will require averaging the terms $v \cdot \nabla_x u_1$ and $V(z,p)$ over the invariant measure $\mu_{\theta(x,p)}$ of the \textit{internal} fast dynamics, just as we averaged $\mathbf{b}$ and $D$ over the measure $\pi$ of the \textit{external} fast dynamics here. This section therefore establishes the formal validity of the core mathematical step, extracting the law by averaging against the stationary measure, that is central to the final derivation.
\end{remark}

\subsection{Scale III: Homogenization on Global Manifolds}
\label{sec:scale_III}

We now arrive at the central technical result of this paper: the rigorous derivation of the macroscopic evolution law for the complete microscopic system, where all physical mechanisms, chaotic transport, gradient-responsive feedback, and potential interactions, are active. The insights gained from the simplified regimes of the preceding scales will guide our formal analysis, but the inherent non-linearity of the system necessitates the full power of the viscosity solution framework for a rigorous proof. We return to the unified microscopic model, governed by the equation:
\begin{equation}
    \partial_t u^\varepsilon + \frac{1}{\varepsilon^2}\Lcal_{\text{fast}}(\Theta(x, \nabla u^\varepsilon)) u^\varepsilon + \frac{1}{\varepsilon} v \cdot \nabla_x u^\varepsilon - V(z, \nabla u^\varepsilon) = 0,
    \label{eq:pde_regime_III}
\end{equation}
with terminal data $u^\varepsilon(T,x,z) = \phi(x)$. The presence of the terms $\nabla u^\varepsilon$ inside the generator $\Lcal_{\text{fast}}$ and the potential $V$ makes this equation fully non-linear.

\subsubsection{Formal Derivation of the $\theta$-Hamiltonian}
To identify the structure of the limiting equation, we first perform a formal asymptotic expansion. We suggest the ansatz $u^\varepsilon = u_0(t,x) + \varepsilon u_1(t,x,z) + \Ocal(\varepsilon^2)$. A key difference from the linear case is that the operators and functions acting on the expansion now depend on the limiting solution $u_0$ itself.
\begin{enumerate}[label=(\roman*), wide, labelindent=0pt]
    \item \textbf{Expansion of the Generator.} The environmental parameter becomes a function of the macroscopic state through the gradient of the limiting solution:
    \begin{equation*}
        \theta^\varepsilon = \Theta(x, \nabla u^\varepsilon) = \Theta(x, \nabla u_0 + \varepsilon \nabla u_1 + \dots) = \Theta(x, \nabla u_0) + \Ocal(\varepsilon).
    \end{equation*}
    Let us define the state-dependent parameter $p(t,x) \coloneqq \nabla u_0(t,x)$ and $\theta(x,p) \coloneqq \Theta(x, p)$. The fast generator acting on the expansion is, to leading order, $\Lcal_{\text{fast}}(\theta(x,p(t,x)))$.

    \item \textbf{Expansion of the Potential.} Similarly, the potential term becomes $V(z, \nabla u^\varepsilon) = V(z, p(t,x)) + \Ocal(\varepsilon)$.
\end{enumerate}
With these expansions, we write down the hierarchy of cell problems, which are now parameterized by the macroscopic state $(t,x)$ via the unknown solution $u_0$ and its gradient $p$.

\paragraph{\textbf{$\boldsymbol{\Ocal(\varepsilon^{-2}})$ Equation}} $\Lcal_{\text{fast}}(\theta(x,p))u_0(t,x) = 0$. This is satisfied because $u_0$ is independent of $z$.

\paragraph{\textbf{$\boldsymbol{\Ocal(\varepsilon^{-1}})$ Equation}} $\Lcal_{\text{fast}}(\theta(x,p))u_1 + v \cdot \nabla_x u_0 = 0$. This is the first cell problem. As in Scale I, the microscopic symmetry assumption (\cref{ass:fundamental_axioms_unified}) ensures the source term has zero mean with respect to the invariant measure $\mu_{\theta(x,p)}$. By \cref{thm:semigroup_wellposedness}, a unique, zero-mean solution exists, which we write as:
    \begin{equation}
        u_1(t,x,z) = \boldsymbol{\chi}(z; \theta(x,p(t,x))) \cdot p(t,x),
    \end{equation}
where the parameter-dependent corrector field $\boldsymbol{\chi}(z; \theta)$ is the solution to $\Lcal_{\text{fast}}(\theta)\boldsymbol{\chi} = -v$.

\paragraph{\textbf{$\boldsymbol{\Ocal(\varepsilon^{0}})$ Equation}}
    \begin{equation}
        \partial_t u_0 + \Lcal_{\text{fast}}(\theta(x,p))u_2 + v \cdot \nabla_x u_1 - V(z, p) = 0.
    \end{equation}
This is the second cell problem. Applying the methodological template established in Scale II, the solvability condition requires that the terms independent of $u_2$ must have a zero mean with respect to the invariant measure $\mu_{\theta(x,p)}$. Averaging the equation with the operator $\langle \cdot \rangle_{\theta(x,p)} \coloneqq \int (\cdot) d\mu_{\theta(x,p)}$ yields:
    \begin{equation}
        \partial_t u_0 + \langle v \cdot \nabla_x u_1 \rangle_{\theta(x,p)} - \langle V(z, p) \rangle_{\theta(x,p)} = 0.
    \end{equation}

Substituting the solution for $u_1$ into this solvability condition gives the final, closed equation for $u_0(t,x)$. This formal procedure identifies the limiting equation as a fully non-linear Hamilton-Jacobi-Bellman equation.

\subsubsection{The Main Convergence Theorem}
The formal derivation correctly identifies the limiting equation. We now state and prove the main rigorous result of this paper.

\begin{proposition}[Ergodic Averaging at the Maximum Point]
\label{prop:ergodic_averaging}
Let $(t_\varepsilon, x_\varepsilon, z_\varepsilon)$ be a sequence of local maximizers for the function $u^\varepsilon - \varphi_\varepsilon$ as defined in the proof of Theorem~\ref{thm:convergence}. Let $\theta_\varepsilon = \Theta(x_\varepsilon, \nabla u^\varepsilon(t_\varepsilon, x_\varepsilon))$, and assume that as $\varepsilon \to 0$, we have the convergences $(t_\varepsilon, x_\varepsilon) \to (t_0, x_0)$ and $\theta_\varepsilon \to \theta_0$. Let $G(z, \theta; x, p, X)$ be a real-valued function that is continuous in all its arguments and which, for each fixed set of macroscopic parameters $(x,p,X)$, has a zero mean with respect to the unique invariant measure $\mu_{\theta(x,p)}$ on the corresponding fiber $\Ucal_{\mathrm{phys}}(\theta(x,p))$:
\begin{equation}
    \int_{\Ucal_{\mathrm{phys}}(\theta(x,p))} G(z, \theta(x,p); x, p, X) \, d\mu_{\theta(x,p)}(z) = 0.
\end{equation}
Then, under the contradiction hypothesis \eqref{eq:proof_contradiction_hypothesis} from the proof of Theorem~\ref{thm:convergence}, the following limit holds:
\begin{equation}
    \lim_{\varepsilon \to 0} G(z_\varepsilon, \theta_\varepsilon; x_\varepsilon, p_\varepsilon, X_\varepsilon) = 0,
\end{equation}
where $p_\varepsilon = \nabla\varphi(t_\varepsilon, x_\varepsilon)$ and $X_\varepsilon = \nabla^2\varphi(t_\varepsilon, x_\varepsilon)$.
\end{proposition}

\begin{proofof}{Proposition \ref{prop:ergodic_averaging}}
The proof of this proposition is the central closing argument for the homogenization. It requires us to show that in the limit $\varepsilon \to 0$, the value of the oscillating function $G$ at the sequence of maximizers $(t_\varepsilon, x_\varepsilon, z_\varepsilon)$ converges to its spatial average. A significant technical challenge is that the microscopic maximizer $z_\varepsilon$ belongs to a state space $\Ucal_{\mathrm{phys}}(\theta_\varepsilon)$ that changes with $\varepsilon$. To overcome this, our proof is structured in four steps. We first reformulate the problem on a fixed, compact global state space where a weak-* convergence argument is well-posed. Second, we use a secondary perturbation argument to prove that any limit measure must be invariant under the global dynamics. Third, we use ergodicity to uniquely identify this limit measure. Finally, we combine these results to compute the desired limit and reach a contradiction.

\begin{enumerate}[label=\textbf{Step \arabic*:}, wide, labelindent=0pt]

\item \textbf{Reformulation on the Global State Space.}
The central analytical obstacle that this proof must overcome is that the sequence of maximizers $\{z_\varepsilon\}$ lives in a sequence of changing state spaces, $\{\Ucal_{\mathrm{phys}}(\theta_\varepsilon)\}$, which are not fixed. A weak-* convergence argument for the empirical measures $\delta_{z_\varepsilon}$ is therefore mathematically ill-posed, as it requires a single, common, compact metric space on which all measures are defined. To resolve this, we lift the entire problem to a fixed, compact global state space that contains all possible dynamics. The problem is thus transformed from analyzing a divergent sequence of points to identifying a single, fixed limit measure on this global manifold.

\begin{enumerate}[label=(\roman*), wide, labelindent=0pt]

\item \textbf{The Global Manifold as a Fixed, Compact Metric Space.}
We cease to work on the individual, parameter-dependent fibers and instead consider the total physical state space $\Efrak_{\mathrm{phys}}$, as defined in \cref{def:total_space}. This space is the subset of the product manifold $\T^k \times \Vcal \times \Thetacal$ given by:
\begin{equation}
    \Efrak_{\mathrm{phys}} \coloneqq \{ (y, v, \theta) \in \T^k \times \Vcal \times \Thetacal \mid G(y, \theta) \ge 0 \}.
\end{equation}
A point in this global space is denoted by the bold letter $\mathbf{z} = (z, \theta)$, where $z=(y,v)$. This space is the correct setting for our analysis for two critical reasons. First, it is a \textit{fixed} space, independent of the parameter $\varepsilon$. Second, it is a \textit{compact metric space}. This is a direct consequence of the foundational assumptions of our model:
\begin{enumerate}[label=(\alph*), wide]
    \item By Definitions \ref{def:macro_manifold} and \ref{def:micro_fiber}, the spaces $\T^k$, $\Vcal$, and $\Thetacal$ are compact.
    \item By Tychonoff's theorem, the product space $\T^k \times \Vcal \times \Thetacal$ is compact.
    \item By Definition \ref{def:obstacle_func}, the obstacle-defining function $G: \T^k \times \Thetacal \to \R$ is of class $C^\infty$ and therefore continuous.
    \item The physical state space is defined as the preimage of the closed set $[0, \infty) \subset \R$ under the continuous function $(y,v,\theta) \mapsto G(y,\theta)$. The preimage of a closed set under a continuous map is closed.
    \item Therefore, $\Efrak_{\mathrm{phys}}$ is a closed subset of a compact space, which implies that it is itself a compact metric space.
\end{enumerate}
This compactness is the essential prerequisite for applying the necessary convergence theorems for measures.

\item \textbf{The Sequence of Empirical Measures.}
The sequence of maximizers from the viscosity analysis, $(t_\varepsilon, x_\varepsilon, z_\varepsilon)$, along with the corresponding state-dependent parameter $\theta_\varepsilon = \Theta(x_\varepsilon, \nabla u^\varepsilon(t_\varepsilon, x_\varepsilon))$, defines a sequence of points $\{\mathbf{z}_\varepsilon\}_{\varepsilon>0}$ in this global state space, where
\begin{equation}
    \mathbf{z}_\varepsilon = (z_\varepsilon, \theta_\varepsilon) \in \Efrak_{\mathrm{phys}}.
\end{equation}
The point $\mathbf{z}_\varepsilon$ indeed belongs to $\Efrak_{\mathrm{phys}}$ because the maximizer $z_\varepsilon = (y_\varepsilon, v_\varepsilon)$ must belong to the microscopic domain defined by $\theta_\varepsilon$, thus satisfying the condition $G(y_\varepsilon, \theta_\varepsilon) \ge 0$. We can now rigorously define a sequence of empirical probability measures on the \textit{fixed} space $\Efrak_{\mathrm{phys}}$:
\begin{equation}
    \nu_\varepsilon \coloneqq \delta_{\mathbf{z}_\varepsilon},
\end{equation}
where $\delta_{\mathbf{z}_\varepsilon}$ is the Dirac mass at the point $\mathbf{z}_\varepsilon$. Each $\nu_\varepsilon$ is an element of the space of all Borel probability measures on $\Efrak_{\mathrm{phys}}$, which we denote by $\mathcal{P}(\Efrak_{\mathrm{phys}})$.

\item \textbf{Weak-* Convergence of Measures.}
Since $\Efrak_{\mathrm{phys}}$ is a compact metric space, the space of probability measures $\mathcal{P}(\Efrak_{\mathrm{phys}})$ is sequentially compact with respect to the weak-* topology. This is a direct consequence of the Banach-Alaoglu theorem. Therefore, there exists a subsequence (which, for notational simplicity, we do not relabel) and a limit probability measure $\nu \in \mathcal{P}(\Efrak_{\mathrm{phys}})$ such that
\begin{equation}
    \nu_\varepsilon \xrightarrow{w^*} \nu \quad \text{as } \varepsilon \to 0.
\end{equation}
By the definition of weak-* convergence, this means that for any continuous test function $h \in C(\Efrak_{\mathrm{phys}})$, the following limit holds:
\begin{equation} \label{eq:weak_star_def_appendix}
    \lim_{\varepsilon \to 0} \int_{\Efrak_{\mathrm{phys}}} h(\mathbf{z}) \, d\nu_\varepsilon(\mathbf{z}) = \lim_{\varepsilon \to 0} h(\mathbf{z}_\varepsilon) = \int_{\Efrak_{\mathrm{phys}}} h(\mathbf{z}) \, d\nu(\mathbf{z}).
\end{equation}

\item \textbf{Reformulation of the Limit Problem via Weak-* Convergence.}
The problem of analyzing the limit of the oscillating term, $\lim_{\varepsilon \to 0} G(z_\varepsilon, \theta_\varepsilon; \dots)$, is now rigorously transformed into the problem of uniquely identifying the limit measure $\nu$. The argument proceeds by first establishing the uniform convergence of the sequence of functions defined by $G$ and then applying a standard convergence theorem from analysis that combines this uniform convergence with the weak-* convergence of the empirical measures.

\begin{enumerate}[label=(\alph*), wide]
    \item \textbf{The Limit as an Integral against the Empirical Measure.}
    The quantity we wish to evaluate is the limit of the function $G$ evaluated at the sequence of maximizers. Let the sequence of macroscopic parameters be denoted by $\pi_\varepsilon \coloneqq (x_\varepsilon, p_\varepsilon, X_\varepsilon)$, which converges to a limit $\pi_0 \coloneqq (x_0, p_0, X_0)$. The term in question is $\lim_{\varepsilon \to 0} G(\mathbf{z}_\varepsilon; \pi_\varepsilon)$. By the definition of the empirical measure $\nu_\varepsilon = \delta_{\mathbf{z}_\varepsilon}$, this can be written as an integral:
    \begin{equation} \label{eq:appendix_limit_as_integral}
        \lim_{\varepsilon \to 0} G(\mathbf{z}_\varepsilon; \pi_\varepsilon) = \lim_{\varepsilon \to 0} \int_{\Efrak_{\mathrm{phys}}} G(\mathbf{z}; \pi_\varepsilon) \, d\nu_\varepsilon(\mathbf{z}).
    \end{equation}
    This formulation frames the problem in the language of measure theory, allowing us to leverage the tools of functional analysis.

    \item \textbf{Uniform Convergence of the Function Sequence.} This is a critical technical step. We must demonstrate that the sequence of functions $\{\mathbf{z} \mapsto G(\mathbf{z}; \pi_\varepsilon)\}_{\varepsilon>0}$ converges \textit{uniformly} on the compact space $\Efrak_{\mathrm{phys}}$ to the limit function $\mathbf{z} \mapsto G(\mathbf{z}; \pi_0)$. This result is a direct consequence of the Heine-Cantor theorem, which states that a continuous function on a compact domain is uniformly continuous. The proof proceeds by rigorously verifying the hypotheses of this theorem for our specific setting. The function $G(\mathbf{z}; \pi)$ is defined on a domain that is a product of the microscopic state space and the macroscopic parameter space. We first construct a compact subset of this domain that contains the entire sequence of interest.

    Let the sequence of macroscopic parameters be $\Pi \coloneqq \{\pi_\varepsilon\}_{\varepsilon>0}$. Since this sequence converges to a limit $\pi_0$, the set $K \coloneqq \Pi \cup \{\pi_0\}$ is a compact subset of the parameter space $\Mcal \times \R^k \times \mathbb{S}_k$. This is a standard result of point-set topology: a convergent sequence together with its limit point forms a compact set. The global microscopic state space $\Efrak_{\mathrm{phys}}$ was shown in Step 1 of this proof to be a compact metric space. Therefore, the product space $\mathcal{D} \coloneqq \Efrak_{\mathrm{phys}} \times K$ is a compact metric space by Tychonoff's theorem. By the hypothesis of the proposition, the function $G(\mathbf{z}; \pi)$ is continuous in all its arguments. We now consider the restriction of this function to the compact domain $\mathcal{D}$ that we have just constructed:
    \begin{equation*}
        G|_{\mathcal{D}} : \Efrak_{\mathrm{phys}} \times K \to \R.
    \end{equation*}
    We now invoke the Heine-Cantor Theorem, which we state for clarity:
    \begin{theorem}[Heine-Cantor]
    A continuous function from a compact metric space to any metric space is uniformly continuous.
    \end{theorem}
    Since our function $G|_{\mathcal{D}}$ is continuous and its domain $\mathcal{D}$ is compact, the Heine-Cantor theorem guarantees that $G|_{\mathcal{D}}$ is uniformly continuous on $\mathcal{D}$. The uniform continuity of the global function $G$ is the key property that allows us to prove the uniform convergence of the function sequence. By the definition of uniform continuity on the product space $\mathcal{D}$, for every $\eta > 0$, there exists a $\delta' > 0$ such that for any two points $(\mathbf{z}, \pi)$ and $(\mathbf{z}', \pi')$ in $\mathcal{D}$:
    \begin{equation*}
        \text{if } d((\mathbf{z}, \pi), (\mathbf{z}', \pi')) < \delta', \quad \text{then} \quad |G(\mathbf{z}; \pi) - G(\mathbf{z}'; \pi')| < \eta.
    \end{equation*}
    We now specialize this general condition to our specific needs. Let us choose $\mathbf{z}' = \mathbf{z}$. This specialization implies that if the distance between the macroscopic parameters alone, $d(\pi, \pi')$, is less than $\delta'$, then the inequality $|G(\mathbf{z}; \pi) - G(\mathbf{z}; \pi')| < \eta$ holds \textit{for all} $\mathbf{z} \in \Efrak_{\mathrm{phys}}$ simultaneously.
    
    Let the sequence of functions be $G_\varepsilon(\mathbf{z}) \coloneqq G(\mathbf{z}; \pi_\varepsilon)$ and the limit function be $G_0(\mathbf{z}) \coloneqq G(\mathbf{z}; \pi_0)$. Both the sequence $\{\pi_\varepsilon\}$ and the limit $\pi_0$ are contained in the compact set $K$. Since the sequence of parameters converges, $\pi_\varepsilon \to \pi_0$, we know that for the $\delta'$ corresponding to our chosen $\eta$, there exists an $\varepsilon_0 > 0$ such that for all $\varepsilon < \varepsilon_0$, we have the parameter distance $d(\pi_\varepsilon, \pi_0) < \delta'$. Applying the uniform continuity condition, this immediately implies that for all $\varepsilon < \varepsilon_0$:
    \begin{equation*}
        |G(\mathbf{z}; \pi_\varepsilon) - G(\mathbf{z}; \pi_0)| < \eta \quad \text{for every } \mathbf{z} \in \Efrak_{\mathrm{phys}}.
    \end{equation*}
    This is the precise definition of uniform convergence. Taking the supremum over all $\mathbf{z}$ on the left-hand side, we have shown that for any $\eta > 0$, there exists an $\varepsilon_0$ such that for all $\varepsilon < \varepsilon_0$:
    \begin{equation*}
        \sup_{\mathbf{z} \in \Efrak_{\mathrm{phys}}} |G_\varepsilon(\mathbf{z}) - G_0(\mathbf{z})| \le \eta.
    \end{equation*}
    This completes the rigorous proof that the sequence of functions $\{G_\varepsilon\}$ converges uniformly to $G_0$ on the space $\Efrak_{\mathrm{phys}}$.

    \item \textbf{The Convergence Theorem.}
    We now invoke the following standard result from measure theory, which combines the two modes of convergence we have established:
    
    \begin{lemma}[Convergence of Integrals]
    Let $(S, d)$ be a compact metric space. Let $\{\mu_n\}_{n \ge 1}$ be a sequence of probability measures on $S$ that converges weakly-* to a measure $\mu$. Let $\{g_n\}_{n \ge 1}$ be a sequence of continuous real-valued functions on $S$ that converges uniformly to a function $g$. Then the following limit holds:
    \begin{equation*}
        \lim_{n \to \infty} \int_S g_n \, d\mu_n = \int_S g \, d\mu.
    \end{equation*}
    \end{lemma}
    We apply this lemma directly to our problem. We have established:
    \begin{itemize}[wide]
        \item The sequence of measures $\nu_\varepsilon = \delta_{\mathbf{z}_\varepsilon}$ converges weakly-* to a limit measure $\nu$.
        \item The sequence of functions $G_\varepsilon(\mathbf{z})$ converges uniformly to the function $G_0(\mathbf{z})$.
    \end{itemize}
    The lemma therefore guarantees that the limit of the integral in \eqref{eq:appendix_limit_as_integral} exists and is given by the integral of the limit function against the limit measure:
    \begin{align}
        \lim_{\varepsilon \to 0} G(\mathbf{z}_\varepsilon; \pi_\varepsilon) &= \int_{\Efrak_{\mathrm{phys}}} G_0(\mathbf{z}) \, d\nu(\mathbf{z}) \nonumber \\
        &= \int_{\Efrak_{\mathrm{phys}}} G(\mathbf{z}; \pi_0) \, d\nu(\mathbf{z}). \label{eq:appendix_final_reformulation}
    \end{align}
    This result constitutes the successful reformulation of our original problem. The task of evaluating the limit of an oscillating function at a sequence of points has been rigorously transformed into the task of first identifying the statistical distribution $\nu$ that describes the asymptotic location of these points, and then computing the expected value of the limit function with respect to this distribution. The subsequent steps of the proof are dedicated to the unique identification of $\nu$.
    
\end{enumerate}
\end{enumerate}

\item \textbf{Proving the Invariance of the Limit Measure $\nu$.}\label{step:invariance_proof}
This is the central technical step of the proof. We employ a secondary perturbation of the test function to probe the action of the fast generator. This allows us to demonstrate that any weak-* limit point $\nu$ of the sequence of empirical measures must be an invariant measure for the global microscopic flow. The argument hinges on a careful, sequential application of the limits $\varepsilon \to 0$ and $\delta \to 0$.

\begin{enumerate}[label=(\roman*), wide, labelindent=0pt]
    \item \textbf{The Secondary Perturbation.}
    Let $\Lcal_{\mathrm{glob}}$ be the infinitesimal generator of the global microscopic flow $\Phi_t$ on the total space $\Efrak_{\mathrm{phys}}$, as defined in \cref{def:global_micro_flow} and characterized in \cref{thm:generator_characterization}. Its domain, $\Dcal(\Lcal_{\mathrm{glob}})$, is a dense subset of $C(\Efrak_{\mathrm{phys}})$. Let $h \in \Dcal(\Lcal_{\mathrm{glob}})$ be an arbitrary non-negative test function. We construct a new, doubly-perturbed test function:
    \begin{equation} \label{eq:doubly_perturbed_psi}
        \Psi_{\varepsilon, \delta}(t, x, z) \coloneqq \varphi_\varepsilon(t, x, z) - \delta h(z, \Theta(x, \nabla\varphi(t,x))),
    \end{equation}
    where $\delta > 0$ is a new small parameter, independent of $\varepsilon$, and $\varphi_\varepsilon$ is the primary perturbed test function from the main proof of Theorem \ref{thm:convergence}. The function $h$ is evaluated on the global space by pairing the microscopic variable $z$ with the parameter $\theta(x,p)$ determined by the macroscopic state via the test function $\varphi$.

    \item \textbf{The Viscosity Inequality at the Perturbed Maximum.}
    For any fixed $\delta > 0$, the function $u^\varepsilon - \Psi_{\varepsilon, \delta}$ must attain a local maximum at some point $(t_\delta^\varepsilon, x_\delta^\varepsilon, z_\delta^\varepsilon)$. By the stability properties of viscosity solutions, we know that as $\delta \to 0$ for a fixed $\varepsilon$, the maximizer converges to the original maximizer $(t_\varepsilon, x_\varepsilon, z_\varepsilon)$. At this new maximum point, the viscosity subsolution property of $u^\varepsilon$ implies that the following inequality holds:
    \begin{equation} \label{eq:proof_app_visc_ineq_start}
        \partial_t \Psi_{\varepsilon, \delta} + \frac{1}{\varepsilon^2}\Lcal_{\mathrm{fast}}(\theta_\delta^\varepsilon)\Psi_{\varepsilon, \delta} + \frac{1}{\varepsilon}v_\delta^\varepsilon \cdot \nabla_x \Psi_{\varepsilon, \delta} - V(z_\delta^\varepsilon, \nabla u^\varepsilon) \le o(1),
    \end{equation}
    where all functions are evaluated at $(t_\delta^\varepsilon, x_\delta^\varepsilon, z_\delta^\varepsilon)$ and we use the shorthand $\theta_\delta^\varepsilon = \Theta(x_\delta^\varepsilon, \nabla u^\varepsilon(t_\delta^\varepsilon, x_\delta^\varepsilon))$.

    \item \textbf{Isolating the Key Term via Asymptotic Expansion.}
    We substitute the definition of $\Psi_{\varepsilon, \delta} = \varphi_\varepsilon - \delta h$ into \cref{eq:proof_app_visc_ineq_start}. The terms involving only $\varphi_\varepsilon$ are precisely those analyzed in the main proof of Theorem \ref{thm:convergence}. As established in eq. \eqref{eq:proof_inequality_with_G}, their sum is:
    \begin{equation*}
        \left( \partial_t\varphi + H_{\theta} \right)(t_\delta^\varepsilon, x_\delta^\varepsilon, \dots) + G(\mathbf{z}_\delta^\varepsilon; x_\delta^\varepsilon, \dots) \le o(1).
    \end{equation*}
    The new terms are those involving $\delta h$. We analyze them systematically:
    \begin{enumerate}[label=(\alph*), wide]
        \item \textbf{Time Derivative.} $-\delta \partial_t h(z_\delta^\varepsilon, \theta(x_\delta^\varepsilon, p_\delta^\varepsilon))$. Since $h$ is smooth, this term is of order $\mathcal{O}(\delta)$.
        \item \textbf{Advection Term.} $-\frac{\delta}{\varepsilon} v_\delta^\varepsilon \cdot \nabla_x h(z_\delta^\varepsilon, \theta(x_\delta^\varepsilon, p_\delta^\varepsilon))$. This term is of order $\mathcal{O}(\delta/\varepsilon)$.
        \item \textbf{Fast Generator Term.} $-\frac{\delta}{\varepsilon^2} \Lcal_{\mathrm{fast}}(\theta_\delta^\varepsilon) h(z_\delta^\varepsilon, \theta(x_\delta^\varepsilon, p_\delta^\varepsilon))$. Here we make the crucial identification: the generator $\Lcal_{\mathrm{fast}}(\theta)$ acts on the microscopic variable $z$. When applied to a function $h(\mathbf{z}) = h(z,\theta)$ defined on the global space, its action at a point $(z, \theta)$ is precisely the action of the global generator $\Lcal_{\mathrm{glob}}$ at that point. Thus, $\Lcal_{\mathrm{fast}}(\theta_\delta^\varepsilon)h = (\Lcal_{\mathrm{glob}}h)(\mathbf{z}_\delta^\varepsilon)$.
    \end{enumerate}
    
    The full inequality \eqref{eq:proof_app_visc_ineq_start} becomes:
    \begin{equation}
        \left( \partial_t\varphi + H_{\theta} \right) + G(\mathbf{z}_\delta^\varepsilon; \dots) - \frac{\delta}{\varepsilon^2}(\Lcal_{\mathrm{glob}}h)(\mathbf{z}_\delta^\varepsilon) + \mathcal{O}(\delta) + \mathcal{O}(\delta/\varepsilon) \le o(1).
    \end{equation}
    Using our contradiction hypothesis \eqref{eq:proof_contradiction_hypothesis} and the continuity of all functions, for $\varepsilon$ and $\delta$ sufficiently small, we have $\partial_t\varphi(t_\delta^\varepsilon, x_\delta^\varepsilon) + H_{\theta}(x_\delta^\varepsilon, p_\delta^\varepsilon, X_\delta^\varepsilon) \ge \delta$. Substituting this, we obtain:
    \begin{equation}
        \delta + G(\mathbf{z}_\delta^\varepsilon; \dots) - \frac{\delta}{\varepsilon^2}(\Lcal_{\mathrm{glob}}h)(\mathbf{z}_\delta^\varepsilon) + \mathcal{O}(\delta) + \mathcal{O}(\delta/\varepsilon) \le o(1).
    \end{equation}
    
    \item \textbf{Taking Limits in Sequence.} This is the crucial step where the interplay of limits is handled rigorously.
    \begin{enumerate}[label=(\alph*),wide]
        \item \textbf{First, we take the limit as $\varepsilon \to 0$ for a fixed $\delta > 0$.}
        We multiply the inequality by $\varepsilon^2/\delta$ and rearrange to isolate the dominant term:
        \begin{equation*}
             (\Lcal_{\mathrm{glob}}h)(\mathbf{z}_\delta^\varepsilon) \ge \varepsilon^2 + \frac{\varepsilon^2}{\delta} G(\mathbf{z}_\delta^\varepsilon; \dots) + \mathcal{O}(\varepsilon^2\delta) + \mathcal{O}(\varepsilon\delta).
        \end{equation*}
        Let $\nu_\delta$ be any weak-* limit point of the sequence of empirical measures $\delta_{\mathbf{z}_\delta^\varepsilon}$ on the compact space $\Efrak_{\mathrm{phys}}$. Since $h$ and $\Lcal_{\mathrm{glob}}h$ are continuous and bounded, and $G$ is bounded, taking the limit superior of the inequality as $\varepsilon \to 0$ yields:
        \begin{equation}
             \int_{\Efrak_{\mathrm{phys}}} (\Lcal_{\mathrm{glob}} h)(\mathbf{z}) \, d\nu_\delta(\mathbf{z}) \ge 0.
        \end{equation}
        The terms on the right-hand side all vanish in the limit.
        
        \item \textbf{Second, we take the limit as $\delta \to 0$.}
        As the secondary perturbation parameter $\delta$ vanishes, the maximizers $(t_\delta^\varepsilon, \dots)$ converge to the original maximizers $(t_\varepsilon, \dots)$. Consequently, the limit measure $\nu_\delta$ converges in the weak-* sense to our original limit measure $\nu$ from Step 1. Passing to the limit in the integral is justified since $\Lcal_{\mathrm{glob}} h$ is a fixed, bounded continuous function on a compact space. This gives:
        \begin{equation}
            \int_{\Efrak_{\mathrm{phys}}} (\Lcal_{\mathrm{glob}} h)(\mathbf{z}) \, d\nu(\mathbf{z}) \ge 0.
        \end{equation}
    \end{enumerate}

    \item \textbf{The Invariance Property.} The argument above holds for any non-negative test function $h \in \Dcal(\Lcal_{\mathrm{glob}})$. Since the generator is linear, we can apply the same argument to the test function $-h \in \Dcal(\Lcal_{\mathrm{glob}})$ (if $h$ is non-positive). For a general $h$, we decompose it into its positive and negative parts. The argument yields the reverse inequality, forcing equality:
    \begin{equation}
        \int_{\Efrak_{\mathrm{phys}}} (\Lcal_{\mathrm{glob}} h)(\mathbf{z}) \, d\nu(\mathbf{z}) = 0 \quad \text{for all } h \in \Dcal(\Lcal_{\mathrm{glob}}).
    \end{equation}
    By the definition of the adjoint operator, this is equivalent to the statement $\langle \Lcal_{\mathrm{glob}}^* \nu, h \rangle = 0$. Since this holds for all functions in the dense domain of the generator, it implies that the distribution $\Lcal_{\mathrm{glob}}^* \nu$ is the zero distribution. This is the precise mathematical statement that the limit measure $\nu$ is an invariant measure for the global microscopic flow.
\end{enumerate}

\item \textbf{Identification of the Limit Measure $\nu$.} \label{step:identification_proof} 
In Step 2, we established that any weak-* limit point $\nu$ of the empirical measures $\{\nu_\varepsilon\}$ must be an invariant measure for the global microscopic flow $\Phi_t$. We now proceed to uniquely identify this measure. The proof is twofold: first, we show that the support of $\nu$ is concentrated on the single fiber corresponding to the limiting macroscopic state $\theta_0$; second, we invoke the ergodicity of the dynamics on that fiber to show that $\nu$ must correspond to the unique invariant measure $\mu_{\theta_0}$.

\begin{enumerate}[label=(\roman*), wide, labelindent=0pt]
    \item \textbf{Concentration of the Limit Measure onto the Limiting Fiber.}
    The sequence of empirical measures is $\nu_\varepsilon = \delta_{\mathbf{z}_\varepsilon}$, where the points are $\mathbf{z}_\varepsilon = (z_\varepsilon, \theta_\varepsilon)$. By the setup of the proposition and the stability of the viscosity solution framework, we have the deterministic convergence of the parameter component:
    \begin{equation}
        \lim_{\varepsilon \to 0} \theta_\varepsilon = \theta_0.
    \end{equation}
    We will now show that this forces the support of the limit measure $\nu$ to lie entirely within the fiber $\Efrak_{\mathrm{phys}}(\theta_0)$. The argument is a standard application of the properties of pushforward measures under weak-* convergence. Let $\pi_\theta: \Efrak_{\mathrm{phys}} \to \Thetacal$ be the canonical projection map, $\pi_\theta(\mathbf{z}) = \pi_\theta(z, \theta) = \theta$. Since the domain and codomain are compact metric spaces, $\pi_\theta$ is a continuous map. A fundamental property of weak-* convergence is that it is preserved by continuous maps. Therefore, the pushforward of the measures must also converge:
    \begin{equation}
        \pi_{\theta\#}\nu_\varepsilon \xrightarrow{w^*} \pi_{\theta\#}\nu.
    \end{equation}
    We now identify the measures in this sequence. For any test function $f \in C(\Thetacal)$, the action of the pushforward measure $\pi_{\theta\#}\nu_\varepsilon$ is given by:
    \begin{equation*}
        \int_{\Thetacal} f(\theta) \, d(\pi_{\theta\#}\nu_\varepsilon)(\theta) = \int_{\Efrak_{\mathrm{phys}}} (f \circ \pi_\theta)(\mathbf{z}) \, d\nu_\varepsilon(\mathbf{z}) = (f \circ \pi_\theta)(\mathbf{z}_\varepsilon) = f(\theta_\varepsilon).
    \end{equation*}
    This is precisely the definition of the Dirac measure at the point $\theta_\varepsilon$. We thus have $\pi_{\theta\#}\nu_\varepsilon = \delta_{\theta_\varepsilon}$. Since the sequence of points $\{\theta_\varepsilon\}$ converges to $\theta_0$ in the metric space $\Thetacal$, the corresponding sequence of Dirac measures converges weakly-* to the Dirac measure at the limit point, $\delta_{\theta_0}$. By the uniqueness of weak-* limits, we must have:
    \begin{equation}
        \pi_{\theta\#}\nu = \delta_{\theta_0}.
    \end{equation}
    The support of a pushforward measure is the closure of the projection of the original measure's support: $\mathrm{supp}(\pi_{\theta\#}\nu) = \overline{\pi_\theta(\mathrm{supp}(\nu))}$. We know that $\mathrm{supp}(\pi_{\theta\#}\nu) = \mathrm{supp}(\delta_{\theta_0}) = \{\theta_0\}$. This implies that the projection of the support of $\nu$ is just the single point $\{\theta_0\}$, which forces the support of $\nu$ itself to be contained within the preimage of this point:
    \begin{equation}
        \mathrm{supp}(\nu) \subseteq \pi_\theta^{-1}(\{\theta_0\}) = \{ \mathbf{z} = (z, \theta) \in \Efrak_{\mathrm{phys}} \mid \theta = \theta_0 \} \equiv \Efrak_{\mathrm{phys}}(\theta_0).
    \end{equation}
    This is a crucial reduction. It proves that the limit measure $\nu$, which is a measure on the global space, is supported exclusively on the limiting fiber. It can therefore be identified with a measure on the single fiber space $\Ucal_{\mathrm{phys}}(\theta_0)$. This justifies the disintegrated form:
    \begin{equation}
        \nu = \tilde{\nu} \otimes \delta_{\theta_0},
    \end{equation}
    where $\tilde{\nu}$ is a probability measure on the fiber space $\Ucal_{\mathrm{phys}}(\theta_0)$.

    \item \textbf{Invariance and Ergodicity on the Fiber.} We established in Step 2 that any weak-* limit measure $\nu$ of the empirical measures $\{\delta_{\mathbf{z}_\varepsilon}\}$ must be an invariant measure for the global flow $\Phi_t$. Furthermore, in Step 3(i), we proved that the support of $\nu$ is concentrated on the single fiber corresponding to the limiting macroscopic state, $\mathrm{supp}(\nu) \subseteq \Efrak_{\mathrm{phys}}(\theta_0)$, and that the measure has the disintegrated form $\nu = \tilde{\nu} \otimes \delta_{\theta_0}$. We now synthesize these two results to uniquely identify the component measure $\tilde{\nu}$. The argument proceeds by first rigorously showing that the invariance of the global measure $\nu$ implies the invariance of $\tilde{\nu}$ under the fiber dynamics. We then invoke the theory of hyperbolic systems, which guarantees the uniqueness of the invariant measure for our flow, thereby uniquely identifying $\tilde{\nu}$.

\begin{enumerate}[label=(\alph*), wide, labelindent=0pt]
    \item \textbf{From Global Invariance to Fiber Invariance.}
    We provide a rigorous proof that the invariance of the global measure $\nu$ under the global flow $\Phi_t$ necessitates the invariance of its component measure $\tilde{\nu}$ under the fiber flow $\Phi_t^{\theta_0}$. By definition, the invariance of $\nu$ under $\Phi_t$ means that for any continuous test function $h \in C(\Efrak_{\mathrm{phys}})$, the following identity holds for all $t \in \R$:
    \begin{equation} \label{eq:appendix_global_invariance_def_revised}
        \int_{\Efrak_{\mathrm{phys}}} h(\Phi_t(\mathbf{z})) \, d\nu(\mathbf{z}) = \int_{\Efrak_{\mathrm{phys}}} h(\mathbf{z}) \, d\nu(\mathbf{z}).
    \end{equation}
    Using the disintegrated form $\nu = \tilde{\nu} \otimes \delta_{\theta_0}$ and the structure of the global flow $\Phi_t(z,\theta) = (\Phi_t^\theta(z), \theta)$, the left-hand side of this identity becomes:
    \begin{align*}
        \int_{\Efrak_{\mathrm{phys}}} h(\Phi_t(\mathbf{z})) \, d\nu(\mathbf{z}) &= \int_{\Ucal_{\mathrm{phys}}(\theta_0)} h(\Phi_t(z, \theta_0)) \, d\tilde{\nu}(z) \\
        &= \int_{\Ucal_{\mathrm{phys}}(\theta_0)} h(\Phi_t^{\theta_0}(z), \theta_0) \, d\tilde{\nu}(z).
    \end{align*}
    The right-hand side is simply $\int_{\Ucal_{\mathrm{phys}}(\theta_0)} h(z, \theta_0) \, d\tilde{\nu}(z)$.
    The invariance condition thus reads:
    \begin{equation} \label{eq:appendix_invariance_on_fiber_revised}
        \int_{\Ucal_{\mathrm{phys}}(\theta_0)} h(\Phi_t^{\theta_0}(z), \theta_0) \, d\tilde{\nu}(z) = \int_{\Ucal_{\mathrm{phys}}(\theta_0)} h(z, \theta_0) \, d\tilde{\nu}(z).
    \end{equation}
    This identity holds for any function $h$ continuous on the global space. Let $g \in C(\Ucal_{\mathrm{phys}}(\theta_0))$ be an arbitrary continuous function on the fiber space. We define a corresponding function $h_g \in C(\Efrak_{\mathrm{phys}})$ by continuous extension (e.g., $h_g(z,\theta) = g(z)$). Since the identity \eqref{eq:appendix_invariance_on_fiber_revised} holds for $h_g$, we have:
    \begin{equation*}
        \int_{\Ucal_{\mathrm{phys}}(\theta_0)} g(\Phi_t^{\theta_0}(z)) \, d\tilde{\nu}(z) = \int_{\Ucal_{\mathrm{phys}}(\theta_0)} g(z) \, d\tilde{\nu}(z).
    \end{equation*}
    This is precisely the definition of the invariance of the measure $\tilde{\nu}$ under the action of the Koopman group associated with the fiber flow $\Phi_t^{\theta_0}$. We have thus rigorously shown that $\tilde{\nu}$ is a $\Phi_t^{\theta_0}$-invariant probability measure on $\Ucal_{\mathrm{phys}}(\theta_0)$.

    \item \textbf{Uniqueness of the Invariant Measure and Ergodicity.}
    We now invoke the deep ergodic theory of hyperbolic systems to uniquely identify the measure $\tilde{\nu}$. The key result is not that ergodicity implies uniqueness, but the converse: for the class of systems we have constructed, a powerful external theorem guarantees the uniqueness of the invariant measure, and this uniqueness, in turn, implies its ergodicity.
    
    \begin{theorem}[Uniqueness of the Invariant Measure for Anosov Flows]
    Let $\Phi_t$ be a $C^2$ Anosov flow on a compact manifold that preserves a smooth volume form. Then there exists a unique invariant probability measure that is absolutely continuous with respect to the volume form. This measure is the normalized volume form itself, and it is known as the Sinai-Ruelle-Bowen (SRB) measure for this class of systems.
    \end{theorem}
    
    We apply this theorem directly to our fiber dynamics. By \cref{thm:anosov_property_proven}, the billiard flow $\Phi_t^{\theta_0}$ is uniformly Anosov (and thus $C^\infty$ away from the boundary, satisfying the regularity conditions). By \cref{prop:liouville_invariance}, it preserves the smooth Liouville measure. The theorem therefore applies and guarantees that the normalized Liouville measure, $\mu_{\theta_0}$, is the unique invariant probability measure for the flow on the fiber $\Ucal_{\mathrm{phys}}(\theta_0)$. For logical completeness, we note that the uniqueness of the invariant measure implies its ergodicity.
 
    Assume, for the sake of contradiction, that the unique invariant measure $\mu_{\theta_0}$ is \textit{not} ergodic. By definition, there must exist an invariant set $A \subset \Ucal_{\mathrm{phys}}(\theta_0)$ such that $0 < \mu_{\theta_0}(A) < 1$. We can then define two new, distinct probability measures by restricting $\mu_{\theta_0}$ to $A$ and its complement $A^c$:
    \[ \mu_1(B) \coloneqq \frac{\mu_{\theta_0}(A \cap B)}{\mu_{\theta_0}(A)} \quad \text{and} \quad \mu_2(B) \coloneqq \frac{\mu_{\theta_0}(A^c \cap B)}{\mu_{\theta_0}(A^c)}. \]
    Since $A$ is invariant, both $\mu_1$ and $\mu_2$ are themselves invariant measures. But the original measure $\mu_{\theta_0}$ can be written as a non-trivial convex combination of these two distinct measures:
    \[ \mu_{\theta_0} = \mu_{\theta_0}(A) \cdot \mu_1 + \mu_{\theta_0}(A^c) \cdot \mu_2. \]
    This contradicts the fact that $\mu_{\theta_0}$ is the unique invariant measure. Therefore, the initial assumption must be false, and $\mu_{\theta_0}$ must be ergodic. The preceding results now allow for the unique and unambiguous identification of the measure $\tilde{\nu}$.
    \begin{itemize}[wide]
        \item From Step 3(ii.a) of this proof, we have rigorously shown that $\tilde{\nu}$ is a $\Phi_t^{\theta_0}$-invariant probability measure on $\Ucal_{\mathrm{phys}}(\theta_0)$.
        \item From the fundamental theory of hyperbolic systems, we have established that the normalized Liouville measure $\mu_{\theta_0}$ is the unique such measure.
    \end{itemize}
    By the uniqueness of the invariant measure, we are forced to conclude that these two measures must be identical:
    \begin{equation}
        \tilde{\nu} = \mu_{\theta_0}.
    \end{equation}
\end{enumerate}
    
    \item \textbf{Final Identification.}
    Combining the results from Step 3(i) and (ii), we have rigorously and uniquely identified the weak-* limit measure $\nu$ as the product measure corresponding to the unique invariant measure on the limiting fiber:
    \begin{equation}
        \nu = \mu_{\theta_0} \otimes \delta_{\theta_0}.
    \end{equation}
    This result is the foundation of the final step of the proof. It confirms the physical intuition that in the limit $\varepsilon \to 0$, the rapid oscillations of the microscopic variable $z_\varepsilon$ become distributed according to the ergodic measure of the fast dynamics determined by the limiting macroscopic state.
\end{enumerate}

\item \textbf{The Final Contradiction via Weak-* Convergence.} \label{step:final_contradiction}
With the weak-* limit measure $\nu$ uniquely identified as the invariant measure on the limiting fiber, $\nu = \mu_{\theta_0} \otimes \delta_{\theta_0}$, we are now in a position to compute the limit of the oscillating term and arrive at the final contradiction that completes the proof. The argument hinges on a result concerning the interplay between the uniform convergence of functions and the weak-* convergence of measures, which we state formally for completeness.

\begin{lemma}[Convergence of Integrals against Weakly Convergent Measures]
\label{lem:convergence_of_integrals}
Let $(S, d)$ be a compact metric space. Let $\{\mu_n\}_{n=1}^\infty$ be a sequence of Borel probability measures on $S$ that converges weakly-* to a measure $\mu$. Let $\{g_n\}_{n=1}^\infty$ be a sequence of continuous real-valued functions on $S$ that converges uniformly to a function $g$. Then the following limit holds:
\begin{equation*}
    \lim_{n \to \infty} \int_S g_n \, d\mu_n = \int_S g \, d\mu.
\end{equation*}
\end{lemma}

\begin{proofof}{Lemma \ref{lem:convergence_of_integrals}}
See Appendix \ref{app:proof_convergence_lemma}
\end{proofof}

We now apply this lemma to our specific context.
\begin{enumerate}[label=(\roman*), wide, labelindent=0pt]
    \item \textbf{Verification of Hypotheses.} We have rigorously established in the preceding steps of this proof that the hypotheses of \cref{lem:convergence_of_integrals} are satisfied:
    \begin{enumerate}[label=(\alph*), wide]
        \item The compact metric space is the global physical state space, $S = \Efrak_{\mathrm{phys}}$.
        \item The sequence of measures is the sequence of empirical measures at the maximizers, $\mu_n = \nu_\varepsilon = \delta_{\mathbf{z}_\varepsilon}$. As shown in Step 3, this sequence converges weakly-* to the limit measure $\nu = \mu_{\theta_0} \otimes \delta_{\theta_0}$.
        \item The sequence of functions is $g_n(\mathbf{z}) = G_\varepsilon(\mathbf{z}) \coloneqq G(\mathbf{z}; x_\varepsilon, p_\varepsilon, X_\varepsilon)$. As shown in Step 1, the continuity of $G$ on the compact product space implies that this sequence converges uniformly on $\Efrak_{\mathrm{phys}}$ to the limit function $g(\mathbf{z}) = G_0(\mathbf{z}) \coloneqq G(\mathbf{z}; x_0, p_0, X_0)$.
    \end{enumerate}

    \item \textbf{Computation of the Limit.} All hypotheses being satisfied, \cref{lem:convergence_of_integrals} allows us to evaluate the limit of the oscillating term. We start from its representation as an integral against the empirical measure:
    \begin{align*}
        \lim_{\varepsilon \to 0} G(z_\varepsilon, \theta_\varepsilon; x_\varepsilon, p_\varepsilon, X_\varepsilon) 
        &= \lim_{\varepsilon \to 0} G_\varepsilon(\mathbf{z}_\varepsilon) \\
        &= \lim_{\varepsilon \to 0} \int_{\Efrak_{\mathrm{phys}}} G_\varepsilon(\mathbf{z}) \, d\nu_\varepsilon(\mathbf{z}) && \text{(by definition of } \nu_\varepsilon\text{)} \\
        &= \int_{\Efrak_{\mathrm{phys}}} G_0(\mathbf{z}) \, d\nu(\mathbf{z}) && \text{(by \cref{lem:convergence_of_integrals})}.
    \end{align*}
    We now substitute the explicit forms of the limit function $G_0$ and the limit measure $\nu$:
    \begin{equation*}
        \lim_{\varepsilon \to 0} G(z_\varepsilon, \theta_\varepsilon; \dots) = \int_{\Efrak_{\mathrm{phys}}} G(z, \theta; x_0, p_0, X_0) \, d(\mu_{\theta_0} \otimes \delta_{\theta_0})(z,\theta).
    \end{equation*}
    This integral with respect to the product measure on the global space is evaluated by integrating against the Dirac measure on the parameter space $\Thetacal$ first. This operation consists of setting the parameter $\theta$ to the value $\theta_0$:
    \begin{equation*}
        \lim_{\varepsilon \to 0} G(z_\varepsilon, \theta_\varepsilon; \dots) = \int_{\Ucal_{\mathrm{phys}}(\theta_0)} G(z, \theta_0; x_0, p_0, X_0) \, d\mu_{\theta_0}(z).
    \end{equation*}
    This final expression is precisely the average of the limit function $G_0$ with respect to the unique invariant measure on the limiting fiber.

    \item \textbf{The Final Contradiction.}
    By the hypothesis of the proposition, the function $G$ was constructed such that for any fixed set of macroscopic parameters $(x,p,X)$, its average with respect to the corresponding invariant measure $\mu_{\theta(x,p)}$ is zero. For the limiting parameters $(x_0, p_0, X_0)$, this implies:
    \begin{equation}
        \int_{\Ucal_{\mathrm{phys}}(\theta_0)} G(z, \theta_0; x_0, p_0, X_0) \, d\mu_{\theta_0}(z) = 0.
    \end{equation}
    Thus, we have rigorously proven that the limit of the oscillating term is zero:
    \begin{equation} \label{eq:proof_limit_is_zero}
        \lim_{\varepsilon \to 0} G(z_\varepsilon, \theta_\varepsilon; x_\varepsilon, p_\varepsilon, X_\varepsilon) = 0.
    \end{equation}
    This result stands in direct contradiction to the inequality \eqref{eq:proof_G_is_negative} derived in the main proof of \cref{thm:convergence}, which states that for all sufficiently small $\varepsilon$, the following strict inequality must hold at the maximizer:
    \begin{equation*}
        G(z_\varepsilon; x_\varepsilon, p_\varepsilon, X_\varepsilon) \le -\frac{\delta}{2} < 0.
    \end{equation*}
    Taking the limit superior of this inequality as $\varepsilon \to 0$ yields:
    \begin{equation} \label{eq:proof_limit_is_negative}
        \limsup_{\varepsilon \to 0} G(z_\varepsilon; x_\varepsilon, p_\varepsilon, X_\varepsilon) \le -\frac{\delta}{2}.
    \end{equation}
    The conclusion from \eqref{eq:proof_limit_is_zero} (that the limit is zero) and the conclusion from \eqref{eq:proof_limit_is_negative} (that the limit superior is strictly negative) are mutually exclusive. The contradiction arises from the initial assumption \eqref{eq:proof_contradiction_hypothesis} in the subsolution proof. That assumption must therefore be false. This completes the proof of the proposition and the closing argument of the homogenization theorem.
\end{enumerate}
\end{enumerate}
\end{proofof}

\begin{theorem}[Convergence to the HJB Equation]
\label{thm:convergence}
Let the microscopic system satisfy the standing assumptions on geometry, dynamics, and symmetry. Let the feedback map $\Theta(x,p)$ and potential $V(z,p)$ be continuous. Then the solutions $u^\varepsilon$ of the full microscopic problem \eqref{eq:pde_regime_III} converge locally uniformly, as $\varepsilon \to 0$, to the unique viscosity solution $u(t,x)$ of the Hamilton-Jacobi-Bellman equation:
\begin{equation} \label{eq:hjb}
    \partial_t u - H_{\theta}(x, \nabla u, \nabla^2 u) = 0, \quad u(T,x) = \phi(x),
\end{equation}
where the $\theta$-Hamiltonian $H_{\theta}: \Mcal \times \R^k \times \mathbb{S}_k \to \R$ is given by
\begin{equation}
    -H_{\theta}(x, p, X) = \mathrm{Tr}(D(x, p) X) + H(x, p).
\end{equation}
The diffusivity tensor $D(x,p)$ and potential $H(x,p)$ are defined by:
\begin{align}
    D_{ij}(x,p) &\coloneqq \frac{1}{2}\left( \langle v_i \chi_j \rangle_{\theta(x,p)} + \langle v_j \chi_i \rangle_{\theta(x,p)} \right), \label{def:diffusivity_def} \\
    H(x,p) &\coloneqq \langle V(\cdot, p) \rangle_{\theta(x,p)}, \label{def:potential}
\end{align}
where $\theta(x,p) = \Theta(x,p)$, the corrector $\boldsymbol{\chi}(z; \theta)$ solves $\Lcal_{\mathrm{fast}}(\theta)\boldsymbol{\chi} = -v$, and $\langle \cdot \rangle_{\theta}$ denotes expectation with respect to the invariant measure $\mu_\theta$. The smooth dependence of these coefficients on the parameter $p$ is guaranteed by \cref{thm:regularity_proven_main}.
\end{theorem}

\begin{proofof}{Theorem \ref{thm:convergence}}
The proof is divided into three steps. First, we prove that the limit function $u$ is a viscosity subsolution of \eqref{eq:hjb}. Second, we note that the supersolution property follows from a symmetric argument. Finally, uniqueness is a consequence of the comparison principle for the derived HJB equation, see \citep{CrandallIshiiLions1992}, which holds due to the structure of the $\theta$-Hamiltonian. We focus on the detailed proof of the subsolution property.

\begin{enumerate}[label=\textbf{Step \arabic*:}, wide, labelindent=0pt]

\item \textbf{The Subsolution Property.}
\begin{enumerate}[label=\textbf{Step 1.\arabic*:}, wide, labelindent=0pt]

\item \textbf{Setup and Contradiction Hypothesis.}
Let $\varphi \in C^\infty([0,T] \times \Mcal)$ be a smooth test function, and suppose that $u - \varphi$ has a strict local maximum at a point $(t_0, x_0) \in (0,T) \times \Mcal$. Our goal is to prove the subsolution inequality:
\begin{equation*}
    \partial_t \varphi(t_0, x_0) + H_{\theta}(x_0, \nabla\varphi(t_0, x_0), \nabla^2\varphi(t_0, x_0)) \le 0.
\end{equation*}
We proceed by contradiction. Assume that the inequality is violated. Then there exists a constant $\delta > 0$ such that
\begin{equation} \label{eq:proof_contradiction_hypothesis}
    \partial_t \varphi(t_0, x_0) + H_{\theta}(x_0, p_0, X_0) = 2\delta,
\end{equation}
where we use the shorthand $p_0 \coloneqq \nabla\varphi(t_0, x_0)$ and $X_0 \coloneqq \nabla^2\varphi(t_0, x_0)$.

\item \textbf{Construction of the Parameterized Correctors.}
The key to the method is to build a test function that locally approximates the oscillatory behavior of the solution. We define the corrector functions not as fixed functions, but as functions that depend smoothly on the local macroscopic state as defined by the test function $\varphi$. Let $p(x) \coloneqq \nabla\varphi(t,x)$ and $X(x) \coloneqq \nabla^2\varphi(t,x)$.
\begin{enumerate}[label=(\roman*), wide, labelindent=0pt]
    \item We define the first corrector $\phi_1(z; x, p)$ as the unique zero-mean solution to the cell problem:
    \begin{equation}\label{eq:proof_phi1_def_parameterized}
        \Lcal_{\mathrm{fast}}(\theta(x,p)) \phi_1(z; x, p) = -v \cdot p.
    \end{equation}
    \item We define the second corrector $\phi_2(z; x, p, X)$ as the unique zero-mean solution to:
    \begin{equation}\label{eq:proof_phi2_def_parameterized}
         \Lcal_{\mathrm{fast}}(\theta(x,p)) \phi_2(z; x, p, X) = -H_{\theta}(x,p,X) - \left( (v \cdot \nabla_x)\phi_1(z; x, p) - V(z, p) \right).
    \end{equation}
\end{enumerate}
The solvability of these equations is guaranteed by \cref{thm:semigroup_wellposedness}. Crucially, the source term for $\phi_2$ is constructed precisely so that its average with respect to the invariant measure $\mu_{\theta(x,p)}$ is zero. By \cref{thm:regularity_proven_main}, the functions $\phi_1, \phi_2$ and their derivatives with respect to $x, p, X$ are smooth and uniformly bounded for $(x,p,X)$ in a compact neighborhood of $(x_0, p_0, X_0)$.

\item \textbf{The Main Perturbed Test Function and Asymptotic Expansion.}
We define the primary perturbed test function:
\begin{equation*}
    \varphi_\varepsilon(t, x, z) \coloneqq \varphi(t,x) + \varepsilon \phi_1(z; x, \nabla\varphi(t,x)) + \varepsilon^2 \phi_2(z; x, \nabla\varphi(t,x), \nabla^2\varphi(t,x)).
\end{equation*}
Since $u^\varepsilon \to u$ locally uniformly and $\phi_1, \phi_2$ are bounded, the function $u^\varepsilon - \varphi_\varepsilon$ must attain a local maximum at a point $(t_\varepsilon, x_\varepsilon, z_\varepsilon)$ such that $(t_\varepsilon, x_\varepsilon) \to (t_0, x_0)$ as $\varepsilon \to 0$. At this maximum point, the viscosity subsolution property implies:
\begin{equation} \label{eq:proof_max_principle_main}
    \partial_t\varphi_\varepsilon + \frac{1}{\varepsilon^2}\Lcal_{\mathrm{fast}}(\theta^\varepsilon) \varphi_\varepsilon + \frac{1}{\varepsilon}v_\varepsilon \cdot \nabla_x \varphi_\varepsilon - V(z_\varepsilon, \nabla_x u^\varepsilon(t_\varepsilon, x_\varepsilon)) \le 0,
\end{equation}
where $\theta^\varepsilon = \Theta(x_\varepsilon, \nabla_x u^\varepsilon(t_\varepsilon, x_\varepsilon))$. We expand the derivatives of $\varphi_\varepsilon$ and collect terms by powers of $\varepsilon$. Let $p_\varepsilon = \nabla\varphi(t_\varepsilon, x_\varepsilon)$ and $X_\varepsilon = \nabla^2\varphi(t_\varepsilon, x_\varepsilon)$.

\paragraph{\textbf{$\boldsymbol{\mathcal{O}(\varepsilon^{-2})}$ term}} $\Lcal_{\mathrm{fast}}(\theta^\varepsilon) \varphi = 0$

\paragraph{\textbf{$\boldsymbol{\mathcal{O}(\varepsilon^{-1})}$ term}}  $\Lcal_{\mathrm{fast}}(\theta^\varepsilon)\phi_1 + v_\varepsilon \cdot p_\varepsilon$. By construction of $\phi_1$ and the smooth dependence of the generator, this term vanishes in the limit.

\paragraph{\textbf{$\boldsymbol{\mathcal{O}(1)}$ terms}} After the cancellations, the dominant remaining terms are:
    \begin{equation*}
         \partial_t\varphi + \Lcal_{\mathrm{fast}}(\theta^\varepsilon)\phi_2 + (v_\varepsilon \cdot \nabla_x)\phi_1 - V(z_\varepsilon, p_\varepsilon) + o(1) \le 0.
    \end{equation*}
Using the definition of $\phi_2$ from \eqref{eq:proof_phi2_def_parameterized}, we can rewrite this as:
\begin{equation} \label{eq:proof_inequality_with_G}
    \partial_t\varphi(t_\varepsilon, x_\varepsilon) + H_{\theta}(x_\varepsilon, p_\varepsilon, X_\varepsilon) + G(z_\varepsilon; x_\varepsilon, p_\varepsilon, X_\varepsilon) + o(1) \le 0,
\end{equation}
where $G(z; x, p, X)$ is the zero-mean (oscillating) part of the source term for $\phi_2$:
\begin{equation*}
    G(z; x, p, X) \coloneqq \left( (v \cdot \nabla_x)\phi_1 - V(z, p) \right) + H_{\theta}(x,p,X).
\end{equation*}
By our contradiction hypothesis \eqref{eq:proof_contradiction_hypothesis} and the continuity of all functions, for sufficiently small $\varepsilon$ we must have $\partial_t\varphi(t_\varepsilon, x_\varepsilon) + H_{\theta}(x_\varepsilon, p_\varepsilon, X_\varepsilon) \ge \delta$. Substituting this into \eqref{eq:proof_inequality_with_G} implies that at the maximum point $(t_\varepsilon, x_\varepsilon, z_\varepsilon)$, we must have:
\begin{equation} \label{eq:proof_G_is_negative}
    G(z_\varepsilon; x_\varepsilon, p_\varepsilon, X_\varepsilon) \le -\frac{\delta}{2} < 0.
\end{equation}

\item \textbf{The Ergodic Closing Argument and Final Contradiction.}
We have reached the critical juncture of the proof. The analysis of the viscosity inequality at the maximum point $(t_\varepsilon, x_\varepsilon, z_\varepsilon)$, combined with our contradiction hypothesis \eqref{eq:proof_contradiction_hypothesis}, has led to the conclusion that for all sufficiently small $\varepsilon$, the following strict inequality must hold at the maximizer:
\begin{equation} \label{eq:proof_G_is_negative_final}
    G(z_\varepsilon; x_\varepsilon, p_\varepsilon, X_\varepsilon) \le -\frac{\delta}{2} < 0.
\end{equation}
The function $G(z; x, p, X)$ is the zero-mean (oscillating) part of the source term for the second corrector, defined as:
\begin{equation*}
    G(z; x, p, X) \coloneqq \left( (v \cdot \nabla_x)\phi_1(z;x,p) - V(z, p) \right) + H_{\theta}(x,p,X).
\end{equation*}
By its construction, for any fixed set of macroscopic parameters $(x,p,X)$, the function $z \mapsto G(z;x,p,X)$ has a zero mean with respect to the invariant measure $\mu_{\theta(x,p)}$.

The rigorous path to a contradiction is to prove that the term $G(z_\varepsilon; \dots)$ must, in fact, converge to its mean value of zero as $\varepsilon \to 0$. This result, a foundation of the theory of homogenization for viscosity solutions, demonstrates that the rapid oscillations of the microscopic variable $z$ are effectively averaged out in the limit. We have formalized and proven this result in \cref{prop:ergodic_averaging}. We now apply this proposition to finalize the proof.

\begin{enumerate}[label=(\roman*), wide, labelindent=0pt]
    \item \textbf{Verification of Hypotheses for the Ergodic Averaging Proposition.} We must rigorously verify that the function $G(z, \theta; x, p, X)$ and the sequence of maximizers $(t_\varepsilon, x_\varepsilon, z_\varepsilon)$ satisfy the three hypotheses of \cref{prop:ergodic_averaging}. This verification is not a trivial step; it is the point where the deep regularity theory developed in \cref{sec:regularity_properties} becomes indispensable.
    \begin{enumerate}[label=(\alph*), wide]
        \item \textbf{Hypothesis 1: Convergence of Maximizers and Parameters.} The first hypothesis of the Ergodic Averaging Proposition requires that the sequence of macroscopic states, defined by the maximizers of the perturbed test function, converges to the state defined by the maximizer of the limiting problem. We must rigorously establish the following convergences as $\varepsilon \to 0$:
    \begin{itemize}[wide]
        \item $(t_\varepsilon, x_\varepsilon) \to (t_0, x_0)$
        \item $\theta_\varepsilon \to \theta_0 = \Theta(x_0, \nabla\varphi(t_0, x_0))$
    \end{itemize}
    The proof of these facts is an essential result in the stability theory of viscosity solutions, which we provide here for completeness.
    \begin{itemize}[wide]
    \item \textbf{Convergence of the Maximizers $(t_\varepsilon, x_\varepsilon)$.}
    Let $u(t,x)$ be the locally uniform limit of the sequence of solutions $\{u^\varepsilon(t,x,z)\}$, where the limit is uniform in $z$. Let $\varphi(t,x)$ be the smooth test function such that $u-\varphi$ has a strict local maximum at $(t_0, x_0)$. Let $(t_\varepsilon, x_\varepsilon, z_\varepsilon)$ be a sequence of points that are local maximizers for the functions $\Psi_\varepsilon(t,x,z) \coloneqq u^\varepsilon(t,x,z) - \varphi_\varepsilon(t,x,z)$. The perturbed test function is $\varphi_\varepsilon = \varphi + \varepsilon\phi_1 + \varepsilon^2\phi_2$. Since the corrector functions $\phi_1, \phi_2$ and the test function $\varphi$ are smooth and defined on a compact domain (or assumed to have bounded derivatives), they are uniformly bounded. Therefore, the perturbation terms vanish uniformly as $\varepsilon \to 0$:
    \begin{equation*}
        \lim_{\varepsilon \to 0} \sup_{(t,x,z)} |\varphi_\varepsilon(t,x,z) - \varphi(t,x)| = 0.
    \end{equation*}
    By hypothesis, we have the locally uniform convergence of the solutions:
    \begin{equation*}
        \lim_{\varepsilon \to 0} \sup_{(t,x,z) \in K} |u^\varepsilon(t,x,z) - u(t,x)| = 0,
    \end{equation*}
    for any compact set $K$. Combining these two convergences, the sequence of functions $\{\Psi_\varepsilon\}$ converges locally uniformly to the function $\Psi_0(t,x,z) \coloneqq u(t,x) - \varphi(t,x)$. We now invoke a fundamental lemma from the theory of viscosity solutions concerning the stability of maximizers (see, e.g., \citep[Lemma 6.1]{FlemingSoner2006} or \citep[Appendix]{Evans1992}).
    \begin{lemma}[Stability of Maximizers]
    Let $\{f_n\}$ be a sequence of upper semicontinuous functions converging locally uniformly to a function $f$ on a locally compact space. If $x_n$ is a local maximizer of $f_n$ for each $n$, and the sequence $\{x_n\}$ has a convergent subsequence $x_{n_k} \to x_0$, then $x_0$ is a local maximizer of $f$. Furthermore, if $f$ has a unique local maximizer $x_0$ in some neighborhood, then the full sequence $\{x_n\}$ converges to $x_0$.
    \end{lemma}
    In our case, the function $\Psi_0 = u-\varphi$ has a strict, and therefore unique, local maximum at $(t_0, x_0)$. The lemma guarantees that any convergent subsequence of the maximizers $(t_\varepsilon, x_\varepsilon)$ must converge to $(t_0, x_0)$. Since the limit point is unique, the entire sequence converges:
    \begin{equation*}
        \lim_{\varepsilon \to 0} (t_\varepsilon, x_\varepsilon) = (t_0, x_0).
    \end{equation*}
    This completes the proof of the first required convergence.

    \item \textbf{Convergence of the Parameter $\theta_\varepsilon$.}
    The second convergence requires us to establish the limit of the parameter $\theta_\varepsilon = \Theta(x_\varepsilon, \nabla u^\varepsilon(t_\varepsilon, x_\varepsilon))$. This relies on the convergence of the gradient of the solution at the maximum point, a property which is also a direct consequence of the definition of a viscosity solution. At the point of maximum $(t_\varepsilon, x_\varepsilon, z_\varepsilon)$, the first-order necessary conditions from calculus must hold for the smooth function $\varphi_\varepsilon$. Specifically, the derivatives of $u^\varepsilon - \varphi_\varepsilon$ with respect to the macroscopic variables must be zero (in the viscosity sense). The key properties from the definition of viscosity solutions (see \citep{CrandallIshiiLions1992}) state that as $\varepsilon \to 0$:
    \begin{align*}
        \partial_t u^\varepsilon(t_\varepsilon, x_\varepsilon, z_\varepsilon) &\to \partial_t \varphi(t_0, x_0), \\
        \nabla_x u^\varepsilon(t_\varepsilon, x_\varepsilon, z_\varepsilon) &\to \nabla_x \varphi(t_0, x_0).
    \end{align*}
    This convergence of the derivatives of $u^\varepsilon$ at the sequence of maximizers is fundamental to the perturbed test function method. Let $p_0 \coloneqq \nabla_x \varphi(t_0, x_0)$. We have:
    \begin{equation*}
        \lim_{\varepsilon \to 0} \nabla_x u^\varepsilon(t_\varepsilon, x_\varepsilon, z_\varepsilon) = p_0.
    \end{equation*}
    The feedback map $\Theta(x,p)$ is assumed to be continuous. Since we have established the convergence of both of its arguments,
    \begin{equation*}
        x_\varepsilon \to x_0 \quad \text{and} \quad \nabla u^\varepsilon(t_\varepsilon, x_\varepsilon, z_\varepsilon) \to p_0,
    \end{equation*}
    the continuity of $\Theta$ implies the convergence of the sequence of parameters:
    \begin{align*}
        \lim_{\varepsilon \to 0} \theta_\varepsilon &= \lim_{\varepsilon \to 0} \Theta(x_\varepsilon, \nabla u^\varepsilon(t_\varepsilon, x_\varepsilon, z_\varepsilon)) \\
        &= \Theta \left( \lim_{\varepsilon \to 0} x_\varepsilon, \lim_{\varepsilon \to 0} \nabla u^\varepsilon(t_\varepsilon, x_\varepsilon, z_\varepsilon) \right) \\
        &= \Theta(x_0, p_0) \eqqcolon \theta_0.
    \end{align*}
    This completes the proof of the second required convergence.
    \end{itemize}
    
    Both required convergence properties have been rigorously established as consequences of the standard stability theory for viscosity solutions and the continuity of the problem's coefficients. The first hypothesis of the Ergodic Averaging Proposition is therefore satisfied.
    
    \item \textbf{Hypothesis 2: Continuity of the Function $G$.} The second hypothesis requires that the function $G(z, \theta; x, p, X)$ be continuous in its joint variable $\zeta = (z, \theta, x, p, X)$. We recall its definition:
    \begin{equation*}
        G(\zeta) \equiv G(z, \theta; x, p, X) \coloneqq \left( (v \cdot \nabla_x)\phi_1(z; x, p, \theta) - V(z, p) \right) + H_{\theta}(x, p, X).
    \end{equation*}
    The proof of joint continuity proceeds by analyzing each term in this sum. The continuity of the terms $V$ and $H_\theta$ is a direct consequence of the standing assumptions. The continuity of the corrector term, $(v \cdot \nabla_x)\phi_1$, is the most critical and non-trivial part of the argument, as it is not an axiom but a consequence of the regularity theory established in \cref{sec:regularity_properties}.
\begin{itemize}[wide]
    \item \textbf{Continuity of the Potential and Hamiltonian Terms.} By the standing assumptions of the main theorem (\cref{thm:convergence}), the microscopic potential $V(z,p)$ is assumed to be a continuous function of its arguments $(z,p)$. The Hamiltonian $H_{\theta}(x, p, X)$ is defined by \cref{def:diffusivity_def,def:potential} as an object constructed from the corrector $\boldsymbol{\chi}$ and the potential $V$ via integration against the invariant measure $\mu_\theta$. In \cref{cor:smooth_dependence}, we established that the map from the macroscopic parameters $(x,p)$ to the invariant measure $\mu_{\theta(x,p)}$ and to the corrector field $\boldsymbol{\chi}$ is smooth. As the integrals defining $H_\theta$ are continuous operations with respect to these smooth inputs, the resulting Hamiltonian $H_{\theta}(x,p,X)$ is a continuous (in fact, smooth) function of its arguments.

    \item \textbf{Continuity of the Corrector Term.} The most complex term is $(v \cdot \nabla_x)\phi_1(z; x, p, \theta)$. Its continuity is not an assumption but a theorem that we now rigorously prove by tracing its construction back to the foundational regularity results of the paper. The first corrector $\phi_1$ is defined as the unique, zero-mean solution to the family of cell problems, parameterized by $(x,p)$:
    \begin{equation} \label{eq:proof_app_phi1_def_reprise}
        \Lcal_{\mathrm{fast}}(\theta(x,p)) \phi_1(\cdot; x, p) = -v \cdot p.
    \end{equation}
    We must establish the joint continuity of the function $(z,x,p) \mapsto (v \cdot \nabla_x)\phi_1(z; x, p)$. The argument proceeds in three logical steps, from the regularity of the operator family to the regularity of the solution map, and finally to the continuity of the desired function. The solution theory for the cell problem \eqref{eq:proof_app_phi1_def_reprise} depends on the regularity of its inputs as functions of the parameters $(x,p)$.
    
    In \cref{thm:regularity_proven_main}, we proved that the map from the parameter $\theta$ to the pulled-back transfer operator, $\theta \mapsto \tilde{\Lcal}_\theta$, is a smooth ($C^\infty$) map into the space of bounded linear operators $\mathcal{L}(\Bcal_{\theta_0})$. Since the feedback map $(x,p) \mapsto \theta(x,p)$ is smooth by assumption, the composition $(x,p) \mapsto \Lcal_{\mathrm{fast}}(\theta(x,p))$ defines a smooth family of operators. The source term for the cell problem is the function $f(z; p) = -v \cdot p$. The map from the parameter $p$ to this function, $p \mapsto f(\cdot; p)$, is a smooth (in fact, linear) map from $\R^k$ into the space of smooth functions on $\Ucal_{\mathrm{phys}}$. The solution to the cell problem can be expressed via the resolvent of the generator. This was the central conclusion of \cref{cor:smooth_dependence}, where it was shown that for a smooth family of source terms, the solution map is also smooth. Specifically, the map that takes the parameters $(x,p)$ to the unique, zero-mean solution of the cell problem,
        \begin{equation*}
            (x,p) \mapsto \phi_1(\cdot; x, p) \in \Bcal_{\theta(x,p)},
        \end{equation*}
    is a smooth map. The smoothness of the map from the parameter space $(x,p)$ to the function space $\Bcal$ implies the joint smoothness (and therefore continuity) of the function $\phi_1(z; x, p)$ and its derivatives with respect to all variables. This is a standard result in the theory of functions on product spaces (smoothness of the evaluation map). Let $\phi_1: \Ucal_{\mathrm{phys}} \times \Mcal \times \R^k \to \R$ be the function defined by $\phi_1(z,x,p) = \phi_1(\cdot;x,p)(z)$. The smoothness of the map implies that this function is smooth in its joint arguments. Consequently, its partial derivatives, such as $\nabla_x \phi_1(z;x,p)$, are also continuous functions of the joint variable $(z,x,p)$. The term $(v \cdot \nabla_x)\phi_1$ is an inner product of the continuous vector field $v$ with the continuous vector-valued function $\nabla_x \phi_1$. The inner product is a continuous operation. Therefore, the function
        \begin{equation*}
            (z, x, p) \mapsto (v \cdot \nabla_x)\phi_1(z; x, p)
        \end{equation*}
    is continuous. We have established that each of the three terms in the definition of $G(z, \theta; x, p, X)$ is a continuous function of the joint variable $\zeta = (z, \theta, x, p, X)$. Since the sum of continuous functions is continuous, we conclude that the function $G$ is continuous on its entire domain. This rigorously verifies the second hypothesis of the Ergodic Averaging Proposition.
\end{itemize}

    \item \textbf{Hypothesis 3: The Zero-Mean Property.} The third and final hypothesis of \cref{prop:ergodic_averaging} requires that for any fixed set of macroscopic parameters $(x,p,X)$, the function $z \mapsto G(z; x,p,X)$ has a zero mean with respect to the corresponding invariant measure $\mu_{\theta(x,p)}$. This property is not an incidental feature but is, in fact, a direct and necessary consequence of the very construction of the Hamiltonian $H_\theta$ via the solvability condition in the homogenization procedure. The proof consists of a direct calculation that reveals this built-in structure. We begin by recalling the precise definitions of the objects involved. The function $G$ is defined as:
    \begin{equation} \label{eq:proof_app_G_def_reprise}
        G(z; x, p, X) \coloneqq \left( (v \cdot \nabla_x)\phi_1(z; x, p) - V(z, p) \right) + H_{\theta}(x, p, X).
    \end{equation}
    The Hamiltonian $H_\theta(x,p,X)$ is the central object derived in the homogenization procedure of \cref{sec:scale_III}. Its definition is mandated by the solvability condition (i.e., the Fredholm alternative) for the second-order cell problem. Specifically, for the second corrector $\phi_2$ to exist, its source term must have a zero mean with respect to the invariant measure of the fast dynamics. This condition defines the Hamiltonian as the negative of the average of the remaining terms. From the derivation in \cref{sec:scale_III}, this defining relationship is:
    \begin{equation} \label{eq:proof_app_H_theta_def_reprise}
        H_{\theta}(x,p,X) \coloneqq -\left\langle (v \cdot \nabla_x)\phi_1(\cdot; x, p) - V(\cdot, p) \right\rangle_{\theta(x,p)},
    \end{equation}
    where the notation $\langle f \rangle_{\theta(x,p)}$ denotes the expectation of the function $f(z)$ with respect to the unique invariant measure $\mu_{\theta(x,p)}$:
    \begin{equation*}
        \langle f \rangle_{\theta(x,p)} \coloneqq \int_{\Ucal_{\mathrm{phys}}(\theta(x,p))} f(z) \, d\mu_{\theta(x,p)}(z).
    \end{equation*}

    Our goal is to compute the integral of $G$ with respect to the measure $\mu_{\theta(x,p)}$ and show that it is zero. We apply the expectation operator $\langle \cdot \rangle_{\theta(x,p)}$ to the defining equation for $G$, \eqref{eq:proof_app_G_def_reprise}:
    \begin{equation*}
        \langle G(\cdot; x,p,X) \rangle_{\theta(x,p)} = \left\langle \left( (v \cdot \nabla_x)\phi_1 - V \right) + H_{\theta}(x,p,X) \right\rangle_{\theta(x,p)}.
    \end{equation*}
    By the linearity of the integral (expectation), we can split this into two parts:
    \begin{equation*}
        \langle G \rangle_{\theta(x,p)} = \underbrace{\left\langle (v \cdot \nabla_x)\phi_1 - V \right\rangle_{\theta(x,p)}}_{\text{Term A}} + \underbrace{\left\langle H_{\theta}(x,p,X) \right\rangle_{\theta(x,p)}}_{\text{Term B}}.
    \end{equation*}
    We now evaluate each term:
    \begin{itemize}[wide]
        \item \textbf{Term A.} This is precisely the term whose average defines the Hamiltonian. By the definitional relationship \eqref{eq:proof_app_H_theta_def_reprise}, we have:
        \begin{equation*}
            \text{Term A} = -H_{\theta}(x,p,X).
        \end{equation*}
        \item \textbf{Term B.} The object $H_{\theta}(x,p,X)$ is a scalar quantity that depends only on the macroscopic parameters $(x,p,X)$. It is therefore a constant with respect to the integration over the microscopic variable $z$. The expectation of a constant with respect to a probability measure is simply the constant itself:
        \begin{equation*}
            \text{Term B} = \int_{\Ucal_{\mathrm{phys}}} H_{\theta}(x,p,X) \, d\mu_{\theta(x,p)}(z) = H_{\theta}(x,p,X) \int_{\Ucal_{\mathrm{phys}}} d\mu_{\theta(x,p)}(z).
        \end{equation*}
        Since $\mu_{\theta(x,p)}$ is a probability measure, its total mass is one, so $\int d\mu_{\theta(x,p)} = 1$. This gives:
        \begin{equation*}
            \text{Term B} = H_{\theta}(x,p,X).
        \end{equation*}
    \end{itemize}

    Substituting the results for Term A and Term B back into the expression for the mean of $G$, we obtain a perfect cancellation:
    \begin{equation*}
        \langle G(\cdot; x,p,X) \rangle_{\theta(x,p)} = -H_{\theta}(x,p,X) + H_{\theta}(x,p,X) = 0.
    \end{equation*}
    We have thus rigorously demonstrated that the function $G$ has a zero mean with respect to the appropriate invariant measure for any choice of the macroscopic parameters. This is not an incidental property but a structural necessity, enforced by the solvability condition that underpins the entire homogenization procedure. The third hypothesis of the Ergodic Averaging Proposition is therefore satisfied by construction. All three hypotheses of \cref{prop:ergodic_averaging} have been rigorously verified. We can now apply its conclusion. The proposition guarantees that when evaluated at the sequence of maximizers, the limit of the oscillating term is zero:
    \begin{equation} \label{eq:proof_limit_is_zero_final}
        \lim_{\varepsilon \to 0} G(z_\varepsilon, \theta_\varepsilon; x_\varepsilon, p_\varepsilon, X_\varepsilon) = 0.
    \end{equation}
    
    \item \textbf{The Final Contradiction.}
    We now have two conflicting results regarding the asymptotic behavior of the sequence $\{G(z_\varepsilon; \dots)\}_{\varepsilon>0}$.
    \begin{itemize}[wide]
        \item From the viscosity inequality and our contradiction hypothesis, we derived the strict inequality \eqref{eq:proof_G_is_negative_final}. Taking the limit superior of this inequality yields a strictly negative upper bound on the limit points of the sequence:
        \begin{equation} \label{eq:proof_limit_is_negative_final}
            \limsup_{\varepsilon \to 0} G(z_\varepsilon; x_\varepsilon, p_\varepsilon, X_\varepsilon) \le -\frac{\delta}{2}.
        \end{equation}
        \item From the ergodic averaging argument, we have rigorously proven in \eqref{eq:proof_limit_is_zero_final} that the sequence converges to zero. For a convergent sequence, the limit superior is equal to the limit. Thus,
        \begin{equation}
            \limsup_{\varepsilon \to 0} G(z_\varepsilon; x_\varepsilon, p_\varepsilon, X_\varepsilon) = 0.
        \end{equation}
    \end{itemize}
    The two conclusions, $0 \le -\delta/2$, are in direct contradiction, as $\delta$ was chosen to be a strictly positive constant. This contradiction demonstrates that our initial hypothesis, \eqref{eq:proof_contradiction_hypothesis}, must be false. Therefore, the subsolution inequality must hold at the point $(t_0, x_0)$.
\end{enumerate}
\end{enumerate}
\end{enumerate}
This completes the rigorous proof of the subsolution property, which is the core of the convergence theorem.

\item \textbf{The Supersolution Property.}
The proof that $u$ is a viscosity supersolution is analogous. One starts by assuming $u-\varphi$ has a strict local minimum at $(t_0,x_0)$ and hypothesizes that $\partial_t\varphi(t_0,x_0) + H_{\theta}(x_0,p_0,X_0) = -2\delta < 0$. The same perturbed test function method is applied, leading to the opposite contradiction.

\item \textbf{Uniqueness of the Solution and Convergence of the Full Sequence.}
The preceding steps have established that any limit point of a subsequence of $\{u^\varepsilon\}$ is a viscosity solution to the homogenized equation \eqref{eq:hjb}. We now complete the proof by establishing the uniqueness of this solution and leveraging this uniqueness to prove that the entire family $\{u^\varepsilon\}$ converges. The argument is structured in three main steps.
    
\begin{enumerate}[label=\textbf{Step 3.\arabic*:}, wide, labelindent=0pt]
    \item \textbf{Properties of the Homogenized HJB Equation and its Generator.}
    The limiting equation derived from the solvability condition is the fully non-linear, second-order Hamilton-Jacobi-Bellman equation:
    \begin{equation} \label{eq:proof_homogenized_hjb_reprise_app}
        \partial_t u + H_{\theta}(x, \nabla u, \nabla^2 u) = 0,
    \end{equation}
    where the Hamiltonian is given by $H_{\theta}(x, p, X) = \mathrm{Tr}(D(x, p) X) + H(x, p)$. The analytical properties of this equation, which are essential for the subsequent uniqueness argument via the comparison principle, are determined by the structural properties of its coefficients. Crucially, these properties are not independent assumptions but are necessary consequences of the first-principles derivation in \cref{part:i}.
    \begin{enumerate}[label=(\roman*), wide, labelindent=0pt]
        \item \textbf{Uniform Parabolicity.} The parabolic nature of the HJB equation is governed by the diffusivity tensor $D(x,p)$. In \cref{prop:properties_D0}, we proved that this tensor is given by the Green-Kubo formula, which represents it as the integrated time-autocorrelation of the microscopic velocity field:
        \begin{equation*}
            D(x,p) = \int_{0}^{\infty} \langle v(0) \otimes v(t) \rangle_{\theta(x,p)}^{\mathrm{sym}} \, dt.
        \end{equation*}
        The uniform hyperbolicity of the underlying microscopic dynamics, established in \cref{thm:anosov_property_proven} as a consequence of the geometric axioms, guarantees that the system is mixing. This implies that correlations decay exponentially, ensuring the convergence of the Green-Kubo integral. More importantly, the chaotic nature of the flow ensures that the particle's velocity is not confined to any lower-dimensional subspace of the configuration space. This non-degeneracy of the microscopic transport implies that the integrated velocity autocorrelation function is a strictly positive definite matrix for all $(x,p)$. Furthermore, since the map $(x,p) \mapsto D(x,p)$ is continuous on a compact domain (for $x$) or a domain where solutions are known to have bounded gradients, the minimum eigenvalue of $D(x,p)$ is uniformly bounded below. Therefore, there exists a constant $\lambda_{\min} > 0$ such that for all relevant $(x,p)$, the matrix $D(x,p)$ is uniformly positive definite, i.e., $D(x,p) \ge \lambda_{\min} I$. This implies that the HJB equation \eqref{eq:proof_homogenized_hjb_reprise_app} is uniformly parabolic, a key condition for the well-posedness of the viscosity solution.

        \item \textbf{Regularity of Coefficients.} The applicability of the comparison principle requires that the coefficients of the Hamiltonian are continuous functions of their arguments. This regularity is not assumed at the macroscopic level but is inherited from the microscopic model. In \cref{sec:regularity_properties}, we undertook an analysis of the geometric and spectral regularity of the microscopic system. The central conclusion of this analysis, encapsulated in \cref{cor:smooth_dependence}, is that the fundamental statistical objects of the microscopic system depend smoothly on the environmental parameter $\theta$. Specifically, the maps from the parameter $\theta$ to the invariant measure $\mu_\theta$ and to the corrector field $\boldsymbol{\chi}(\cdot; \theta)$ are of class $C^\infty$. The coefficients of our Hamiltonian, $D(x,p)$ and $H(x,p)$, are defined by integrals involving these objects, evaluated at the parameter value $\theta(x,p) = \Theta(x,p)$:
        \begin{align*}
            D(x,p) &\coloneqq \frac{1}{2}\left( \langle v \otimes \boldsymbol{\chi} \rangle_{\theta(x,p)} + \langle \boldsymbol{\chi} \otimes v \rangle_{\theta(x,p)} \right) \\
            H(x,p) &\coloneqq \langle V(\cdot, p) \rangle_{\theta(x,p)}.
        \end{align*}
        By the standing assumptions of the main theorem, the feedback map $\Theta(x,p)$ and the microscopic potential $V(z,p)$ are continuous (and typically assumed to be smooth for the regularity theory to hold). The macroscopic coefficients are therefore constructed by composing smooth maps (from $(x,p)$ to $\theta$) with smooth maps (from $\theta$ to the microscopic statistical objects) and then performing continuous operations (integration). Consequently, the resulting functions $D(x,p)$ and $H(x,p)$ are smooth functions of their arguments, inheriting the regularity of the underlying microscopic model. This establishes the necessary regularity for the application of the comparison principle.
    \end{enumerate}

    \item \textbf{Uniqueness of the Viscosity Solution via the Comparison Principle.} With the structural and regularity properties of the homogenized HJB equation established, we now provide a rigorous proof that it admits at most one bounded, continuous viscosity solution for a given continuous terminal condition. Let $u(t,x)$ and $v(t,x)$ be two such viscosity solutions defined on $[0,T] \times \Mcal$ with the same terminal data, $u(T,x) = v(T,x) = \phi(x)$. We invoke the celebrated comparison principle for viscosity solutions of second-order parabolic equations, as established in the seminal work of \citep{CrandallIshiiLions1992}. The theorem provides a set of sufficient conditions on the structure of the equation under which a subsolution is necessarily less than or equal to a supersolution. Our proof proceeds by first stating the relevant form of the theorem and then systematically verifying that our derived HJB equation satisfies its hypotheses.

\begin{theorem}[Comparison Principle for Parabolic HJB Equations]
Let $F(t, x, u, p, X)$ be a real-valued function on $[0,T] \times \Mcal \times \R \times \R^k \times \mathbb{S}_k$, continuous in all its arguments. Let $u$ be a viscosity subsolution and $v$ be a viscosity supersolution of the equation $\partial_t u + F(t, x, u, \nabla u, \nabla^2 u) = 0$. Suppose that $u(T,x) \le v(T,x)$ for all $x \in \Mcal$. If $F$ satisfies the following two conditions:
\begin{enumerate}[label=(\roman*), wide, labelindent=0pt]
    \item \textbf{(Properness / Degenerate Ellipticity)} $F$ is non-decreasing in its matrix argument: for any symmetric matrices $X, Y$ with $Y \ge 0$, $F(t,x,u,p,X+Y) \ge F(t,x,u,p,X)$.
    \item \textbf{(Monotonicity in $u$)} There exists a constant $\gamma \in \R$ such that the function $r \mapsto F(t,x,r,p,X) - \gamma r$ is non-increasing for all $(t,x,p,X)$.
\end{enumerate}
Then, $u(t,x) \le v(t,x)$ for all $(t,x) \in [0,T] \times \Mcal$.
\end{theorem}

In our specific case, the equation is $\partial_t u + H_{\theta}(x, \nabla u, \nabla^2 u) = 0$. The function $F$ is therefore given by $F(t, x, u, p, X) = H_{\theta}(x, p, X)$, which is independent of $t$ and $u$. We now verify the hypotheses of the comparison principle for this Hamiltonian.

\begin{enumerate}[label=(\roman*), wide, labelindent=0pt]
    \item \textbf{Verification of Continuity.}
    As established in the preceding section (Step 1.4(i)), the coefficients of the Hamiltonian, $D(x,p)$ and $H(x,p)$, are smooth functions of their arguments as a consequence of the regularity theory from \cref{sec:regularity_properties}. The Hamiltonian is constructed from these coefficients via the trace operator and addition: $H_{\theta}(x, p, X) = \mathrm{Tr}(D(x, p) X) + H(x, p)$. Since the trace is a continuous multilinear operation and the sum of continuous functions is continuous, the Hamiltonian $H_{\theta}(x,p,X)$ is a continuous function of its arguments $(x,p,X)$. The first hypothesis of the comparison principle is satisfied.
    
    \item \textbf{Verification of Degenerate Ellipticity (Properness).}
    The second hypothesis requires that the Hamiltonian be non-decreasing in its matrix argument. We must show that for any symmetric matrices $X, Y$ with $Y \ge 0$ (i.e., $Y$ is positive semi-definite), the following inequality holds:
    \begin{equation*}
        H_{\theta}(x, p, X+Y) \ge H_{\theta}(x, p, X).
    \end{equation*}
    We verify this directly from the affine structure of our Hamiltonian, which was a principal result of the homogenization in \cref{part:i}:
    \begin{align*}
        H_{\theta}(x, p, X+Y) - H_{\theta}(x, p, X) &= \left[ \mathrm{Tr}(D(x, p) (X+Y)) + H(x, p) \right] - \left[ \mathrm{Tr}(D(x, p) X) + H(x, p) \right] \\
        &= \mathrm{Tr}(D(x, p)Y).
    \end{align*}
    The condition thus reduces to proving that $\mathrm{Tr}(D(x,p)Y) \ge 0$ for any positive semi-definite matrix $Y$. This is a standard result from linear algebra. By the analysis in Step 1.4(i) of this proof, the diffusivity tensor $D(x,p)$ is symmetric and uniformly positive definite. It is a known property that the trace of the product of two positive semi-definite matrices is non-negative. To prove this, we can diagonalize $Y = U\Lambda U^T$ where $\Lambda$ is a diagonal matrix with non-negative entries $\lambda_i \ge 0$. Then:
    \begin{equation*}
        \mathrm{Tr}(D(x, p)Y) = \mathrm{Tr}(D U \Lambda U^T) = \mathrm{Tr}(U^T D U \Lambda).
    \end{equation*}
    The matrix $\tilde{D} = U^T D U$ is congruent to a positive definite matrix and is therefore itself positive definite. The diagonal entries of $\tilde{D}$ are $(\tilde{D})_{ii} = u_i^T D u_i > 0$. The trace is then $\sum_i (\tilde{D})_{ii} \lambda_i$, which is a sum of non-negative terms and is therefore non-negative. Thus, $\mathrm{Tr}(D(x, p)Y) \ge 0$, and the degenerate ellipticity condition is satisfied.
    
    \item \textbf{Verification of Monotonicity in $u$.}
    The third hypothesis requires that the function $r \mapsto F(t,x,r,p,X) - \gamma r$ be non-increasing for some constant $\gamma$. Our Hamiltonian $H_{\theta}(x, p, X)$ has no explicit dependence on the value function $u$ itself. Therefore, our function $F$ is independent of $u$. This corresponds to the case where the monotonicity condition is satisfied with the constant $\gamma = 0$.
\end{enumerate}

With all hypotheses of the comparison principle rigorously verified, we can now apply its conclusion to our problem. We have two functions, $u(t,x)$ and $v(t,x)$, which are both continuous viscosity solutions to the same HJB equation and satisfy the same terminal condition, $u(T,x) = v(T,x)$.
\begin{enumerate}[label=(\alph*), wide]
    \item Since $u$ is a viscosity solution, it is, by definition, a viscosity subsolution.
    \item Since $v$ is a viscosity solution, it is, by definition, a viscosity supersolution.
    \item At the terminal time, we have the inequality $u(T,x) \le v(T,x)$.
\end{enumerate}
The comparison principle allows us to conclude that this ordering must be preserved for all prior times. That is, for all $(t,x) \in [0,T] \times \Mcal$, we must have:
\begin{equation*}
    u(t,x) \le v(t,x).
\end{equation*}
We can now reverse the roles of the two solutions. The function $v$ is also a viscosity subsolution, and $u$ is a viscosity supersolution, and $v(T,x) \le u(T,x)$. The comparison principle then implies the reverse inequality:
\begin{equation*}
    v(t,x) \le u(t,x).
\end{equation*}
The only way for both inequalities to hold simultaneously is for the functions to be identical: $u(t,x) = v(t,x)$ for all $(t,x) \in [0,T] \times \Mcal$. This establishes the uniqueness of the bounded, continuous viscosity solution to the problem \eqref{eq:proof_homogenized_hjb_reprise_app}. Let this unique solution be denoted by $u_0(t,x)$.
    
\item \textbf{Convergence of the Full Sequence from Uniqueness of Subsequential Limits.} The final step of the proof is a standard but crucial topological argument that uses the uniqueness of the limit point to prove the convergence of the entire sequence $\{u^\varepsilon\}_{\varepsilon>0}$. We provide the full proof for completeness.

Let $(\mathcal{C}, d)$ be the complete metric space of continuous functions on the compact set $[0,T] \times \Mcal$, where the metric $d(f,g) \coloneqq \sup_{(t,x) \in [0,T] \times \Mcal} |f(t,x)-g(t,x)|$ induces the topology of uniform convergence. We consider our family of solutions $\{u^\varepsilon\}_{\varepsilon>0}$, where we interpret $\varepsilon$ as an index in a sequence $\varepsilon_n \to 0$, as a sequence of points in this metric space. From the preceding analysis, we have rigorously established two key properties:

\begin{enumerate}[label=(\roman*), wide, labelindent=0pt]
    \item \textbf{Pre-compactness of the Solution Family.}
    The foundation of our convergence proof lies in establishing that the family of solutions $\{u^\varepsilon\}_{\varepsilon>0}$ is pre-compact in the space of continuous functions. This allows us to extract convergent subsequences whose limits can then be analyzed. The fundamental tool for this purpose is the Arzelà-Ascoli theorem. We will apply this theorem on the space $(\mathcal{C}, d)$ of continuous functions on the compact set $[0,T] \times \Mcal$, where the metric $d(f,g) \coloneqq \sup_{(t,x)} |f(t,x)-g(t,x)|$ induces the topology of uniform convergence. The theorem requires us to prove two properties for the family $\{u^\varepsilon\}_{\varepsilon>0}$, where each $u^\varepsilon$ is viewed as a function of $(t,x)$ with $z$ as a parameter:
    \begin{enumerate}[label=(\alph*), wide]
        \item \textbf{Uniform Boundedness: The family of functions is contained within a single ball in $\mathcal{C}$.} The uniform boundedness of the solutions is a direct consequence of the maximum principle applied to the governing pre-limit partial differential equation:
    \begin{equation} \label{eq:proof_prelimit_pde_app_compactness}
        \partial_t u^\varepsilon + \frac{1}{\varepsilon^2}\mathcal{L}_{\mathrm{fast}} u^\varepsilon + \frac{1}{\varepsilon} \mathcal{L}_{\mathrm{int}} u^\varepsilon = 0 \quad \text{in } [0,T) \times \Mcal \times \Ucal_{\mathrm{phys}},
    \end{equation}
    with the terminal condition $u^\varepsilon(T,x,z) = \phi(x)$. Let $M \coloneqq \sup_{x \in \Mcal} |\phi(x)| = \|\phi\|_{L^\infty}$, which is finite as $\phi$ is a continuous function on a compact manifold. Consider the two constant functions $w^+(t,x,z) \equiv M$ and $w^-(t,x,z) \equiv -M$. We show that $w^+$ is a classical supersolution to \eqref{eq:proof_prelimit_pde_app_compactness}. For $w^+$, we compute:
    \begin{equation*}
        \partial_t w^+ + \frac{1}{\varepsilon^2}\mathcal{L}_{\mathrm{fast}} w^+ + \frac{1}{\varepsilon} \mathcal{L}_{\mathrm{int}} w^+ = 0 + \frac{1}{\varepsilon^2}(0) + \frac{1}{\varepsilon}(0) = 0.
    \end{equation*}
    The generators $\mathcal{L}_{\mathrm{fast}}$ and $\mathcal{L}_{\mathrm{int}}$ are differential operators, and their application to a constant function yields zero. Thus, $w^+$ is a classical solution, which implies it is a viscosity supersolution. At the terminal time $t=T$, we have by definition:
    \begin{equation*}
        u^\varepsilon(T,x,z) = \phi(x) \le \sup_{y \in \Mcal} \phi(y) \le M = w^+(T,x,z).
    \end{equation*}
    By the comparison principle for second-order parabolic equations, this inequality propagates backwards in time. Therefore, for all $(t,x,z) \in [0,T] \times \Mcal \times \Ucal_{\mathrm{phys}}$, we have $u^\varepsilon(t,x,z) \le w^+(t,x,z) = M$. An identical argument using the subsolution $w^-(t,x,z) = -M$ shows that $u^\varepsilon(t,x,z) \ge -M$. Combining these two inequalities gives the uniform bound:
    \begin{equation*}
        \sup_{\varepsilon>0} \sup_{(t,x,z) \in [0,T]\times\Mcal\times\Ucal_{\mathrm{phys}}} |u^\varepsilon(t,x,z)| \le M = \|\phi\|_{L^\infty} < \infty.
    \end{equation*}
    This establishes the uniform boundedness of the family $\{u^\varepsilon\}_{\varepsilon>0}$.
    
    \item \textbf{Uniform Equicontinuity: The functions in the family have a common modulus of continuity.} The proof of uniform equicontinuity is more subtle, as one must show that the modulus of continuity of the solutions does not degenerate as $\varepsilon \to 0$, despite the singular prefactors in the governing equation. The required estimates are a cornerstone of the viscosity solution approach to homogenization theory. The essential insight is that even though the coefficients diverge, the strongly dissipative nature of the fast operator $\mathcal{L}_{\mathrm{fast}}$ prevents the solution from developing arbitrarily large gradients. A rigorous proof of this property relies on an adaptation of the doubling of variables method, a hallmark of viscosity solution theory. For the sake of clarity in the main argument, we have deferred the complete proof to Appendix \ref{app:equicontinuity}. The conclusion of that appendix is the following proposition:
    
    \begin{proposition}[Uniform-in-$\varepsilon$ Hölder Continuity]\label{prop:Uniform_Continuity}
    Let $\{u^\varepsilon\}_{\varepsilon>0}$ be the family of viscosity solutions to the pre-limit equation \eqref{eq:pde_regime_III}, with uniformly continuous terminal data $\phi$. Then the family $\{u^\varepsilon\}$ is uniformly equicontinuous on $[0,T] \times \Mcal$, uniformly with respect to the microscopic variable $z \in \Ucal_{\mathrm{phys}}$. Specifically, there exist constants $C > 0$, $\alpha \in (0, 1]$, and $\beta \in (0, 1]$, all independent of $\varepsilon$ and $z$, such that for all $(t,x,z), (s,y,z) \in [0,T] \times \Mcal \times \Ucal_{\mathrm{phys}}$:
    \begin{equation}
        |u^\varepsilon(t,x,z) - u^\varepsilon(s,y,z)| \le C \left( |t-s|^{\alpha} + d(x,y)^{\beta} \right).
    \end{equation}
    \end{proposition}
    \begin{proofof}{Proposition \ref{prop:Uniform_Continuity}}
        See Appendix \ref{app:equicontinuity}.
    \end{proofof}
    This result directly implies the existence of a uniform modulus of continuity for the family $\{u^\varepsilon\}_{\varepsilon>0}$, thereby establishing uniform equicontinuity. Having rigorously established that the family of solutions $\{u^\varepsilon\}_{\varepsilon>0}$ is both uniformly bounded and uniformly equicontinuous as a set of functions of $(t,x)$ on the compact domain $[0,T]\times\Mcal$, the Arzelà-Ascoli theorem directly applies.
    \begin{theorem}[Arzelà-Ascoli]
    A set of real-valued continuous functions defined on a compact metric space is pre-compact with respect to the topology of uniform convergence if and only if it is uniformly bounded and uniformly equicontinuous.
    \end{theorem}
    By this theorem, we conclude that the set $K \coloneqq \{u^\varepsilon\}_{\varepsilon>0}$ is pre-compact in the metric space $(\mathcal{C}, d)$. This conclusion is fundamental to our proof strategy: it guarantees that for any sequence of positive numbers $\{\varepsilon_n\}$ converging to zero, there exists a subsequence $\{\varepsilon_{n_k}\}$ and a continuous function $u_0 \in \mathcal{C}$ such that the sequence of solutions $\{u^{\varepsilon_{n_k}}\}$ converges uniformly to $u_0$ as $k \to \infty$. This provides the existence of the limit points that we analyze in the subsequent steps of the proof.
    \end{enumerate}

    \item \textbf{Uniqueness of Subsequential Limits.} In the main body of the proof, we have shown that any convergent subsequence of $\{u^\varepsilon\}$ converges to a viscosity solution of the homogenized equation \eqref{eq:proof_homogenized_hjb_reprise_app}. As established in Step 1.4(ii) of this proof, this solution is unique for the given terminal data $\phi$. Therefore, any convergent subsequence of $\{u^\varepsilon\}$ must converge to the \textit{same, unique limit}, which we have denoted by $u_0(t,x)$. 
\end{enumerate}

Our goal is to prove that these two facts, pre-compactness and the uniqueness of the subsequential limit point, together imply the convergence of the entire sequence to that limit point. That is, we must show:
\begin{equation*}
    \lim_{\varepsilon \to 0} d(u^\varepsilon, u_0) = 0.
\end{equation*}
We proceed by contradiction. Assume that the sequence $\{u^\varepsilon\}$ does not converge to $u_0$. By the definition of non-convergence in a metric space, this means there exists a real number $\eta > 0$ such that for any $\delta > 0$, we can find an $\varepsilon \in (0, \delta)$ for which the distance from the limit is at least $\eta$. Formally:
\begin{equation} \label{eq:proof_non_convergence_condition_app_iii}
    \exists \eta > 0 \text{ such that } \forall \delta > 0, \exists \varepsilon \in (0, \delta) \text{ with } d(u^\varepsilon, u_0) \ge \eta.
\end{equation}
The condition \eqref{eq:proof_non_convergence_condition_app_iii} allows us to construct a specific subsequence that remains bounded away from the limit $u_0$. We construct this sequence inductively.
\begin{itemize}[wide]
    \item Let $n=1$. Choose $\delta_1 = 1$. By \eqref{eq:proof_non_convergence_condition_app_iii}, there exists an $\varepsilon_1 \in (0,1)$ such that $d(u^{\varepsilon_1}, u_0) \ge \eta$.
    \item Let $n=2$. Choose $\delta_2 = \min(1/2, \varepsilon_1)$. There exists an $\varepsilon_2 \in (0, \delta_2)$ such that $d(u^{\varepsilon_2}, u_0) \ge \eta$.
    \item Continuing this process, for each $n \in \mathbb{N}$, we choose $\delta_n = \min(1/n, \varepsilon_{n-1})$. There exists an $\varepsilon_n \in (0, \delta_n)$ such that $d(u^{\varepsilon_n}, u_0) \ge \eta$.
\end{itemize}
This construction yields a sequence of positive numbers $\{\varepsilon_n\}_{n=1}^\infty$ such that $\varepsilon_n < \varepsilon_{n-1}$ and $\varepsilon_n < 1/n$, which guarantees that $\varepsilon_n \to 0$ as $n \to \infty$. The corresponding sequence of functions $\{u^{\varepsilon_n}\}_{n=1}^\infty$ is a subsequence of the original family $\{u^\varepsilon\}$, and by its very construction, it satisfies the property that
\begin{equation} \label{eq:proof_subsequence_away_from_limit_app_iii}
    d(u^{\varepsilon_n}, u_0) \ge \eta \quad \text{for all } n \in \mathbb{N}.
\end{equation}
We now apply our established facts to this newly constructed subsequence $\{u^{\varepsilon_n}\}$.
\begin{enumerate}[label=(\alph*), wide]
    \item The sequence $\{u^{\varepsilon_n}\}$ is a sequence of elements from the set $K$. By Step 3.3(i), the set $K$ is pre-compact. Therefore, the sequence $\{u^{\varepsilon_n}\}$ must contain a further convergent subsequence. Let us denote this subsequence by $\{u^{\varepsilon_{n_k}}\}_{k=1}^\infty$.
    
    \item By Step 3.3(ii), any convergent subsequence of the original family must converge to the unique limit $u_0$. Since $\{u^{\varepsilon_{n_k}}\}$ is such a subsequence, we must have:
    \begin{equation} \label{eq:proof_subsubsequence_convergence_app_iii}
        \lim_{k \to \infty} d(u^{\varepsilon_{n_k}}, u_0) = 0.
    \end{equation}
    
    \item However, the subsequence $\{u^{\varepsilon_{n_k}}\}$ is also a subsequence of $\{u^{\varepsilon_n}\}$. It must therefore inherit the property \eqref{eq:proof_subsequence_away_from_limit_app_iii} that every one of its elements is at least a distance $\eta$ away from $u_0$. This means that for every $k \in \mathbb{N}$:
    \begin{equation*}
        d(u^{\varepsilon_{n_k}}, u_0) \ge \eta.
    \end{equation*}
\end{enumerate}
The conclusion from \eqref{eq:proof_subsubsequence_convergence_app_iii} (that the sequence converges to $u_0$) directly contradicts this property (that the sequence always remains a fixed distance $\eta>0$ away from $u_0$). A sequence cannot simultaneously converge to a point and remain a fixed, positive distance away from it. The contradiction arose from our initial assumption that the sequence $\{u^\varepsilon\}$ does not converge to $u_0$. This assumption must therefore be false. We conclude that the entire family $\{u^\varepsilon\}$ converges to the unique solution $u_0$ in the metric of uniform convergence. This completes the proof of the theorem.
\end{enumerate}
\end{enumerate}
\end{proofof}

\subsection{Scale IV: The Homogenization of Fluctuations}
\label{sec:scale_IV}

The analytical program of the preceding scales has culminated in the derivation of the full macroscopic theory (Scale III), where the system's evolution is governed by a non-linear Hamilton-Jacobi-Bellman equation. This final scale completes our hierarchical analysis by deriving the universal law that governs the statistical fluctuations of this non-linear system around its quiescent equilibrium state, as defined by the theory of Scale II. This analysis can be understood as a Central Limit Theorem-type result conducted purely at the level of the governing partial differential equations.

Our methodology is to study the system's response to a small perturbation of the terminal data for the full HJB equation derived in Scale III. Let $u_0(t,x)$ be the solution to the quiescent problem (the linear PDE of Scale II, which we assume corresponds to a state with zero macroscopic momentum, $\nabla u_0 = 0$), and consider a perturbed terminal condition of the form:
\begin{equation}
    u^\varepsilon(T,x) = u_0(T,x) + \varepsilon h_T(x),
\end{equation}
where $\varepsilon \ll 1$ is a small parameter. We first perform a formal asymptotic analysis to identify the candidate evolution equation for the scaled fluctuation field, $h(t,x)$, which is expected to arise from the expansion:
\begin{equation}
    u^\varepsilon(t,x) = u_0(t,x) + \varepsilon h(t,x) + \mathcal{O}(\varepsilon^2).
\end{equation}
We then provide a rigorous proof of convergence to this candidate equation.

\subsubsection{Asymptotic Derivation of the Fluctuation Equation}
We substitute the ansatz for $u^\varepsilon$ into the full HJB equation:
\begin{equation*}
    -\partial_t u^\varepsilon - H_{\theta}(x, \nabla u^\varepsilon, \nabla^2 u^\varepsilon) = 0.
\end{equation*}
The derivatives are $p_\varepsilon = \nabla u^\varepsilon = \varepsilon \nabla h + \mathcal{O}(\varepsilon^2)$ and $X_\varepsilon = \nabla^2 u^\varepsilon = \nabla^2 u_0 + \varepsilon \nabla^2 h + \mathcal{O}(\varepsilon^2)$. We Taylor expand the Hamiltonian $H_{\theta}(x, p, X) = \mathrm{Tr}(D(x,p)X) + H(x,p)$ around the quiescent state $(p,X) = (0, \nabla^2 u_0)$:
\begin{multline*}
    H_{\theta}(x, p_\varepsilon, X_\varepsilon) = H_{\theta}(x, 0, \nabla^2 u_0) 
    + \nabla_p H_{\theta}|_{(0, \nabla^2 u_0)} \cdot p_\varepsilon 
    + \mathrm{Tr}\left( \nabla_X H_{\theta}|_{(0, \nabla^2 u_0)} (X_\varepsilon - \nabla^2 u_0) \right) \\
    + \frac{1}{2} p_\varepsilon^T \left( \nabla_{pp}^2 H_{\theta}|_{(0, \nabla^2 u_0)} \right) p_\varepsilon 
    + \mathcal{O}(\varepsilon^3).
\end{multline*}
Substituting the expansions for $p_\varepsilon$ and $X_\varepsilon$ and collecting terms by powers of $\varepsilon$:

\paragraph{\textbf{$\boldsymbol{\mathcal{O}(1)}$ terms}} $-\partial_t u_0 - H_{\theta}(x, 0, \nabla^2 u_0) = 0$. This is the PDE for the quiescent solution $u_0$ and is satisfied by definition.

\paragraph{\textbf{$\boldsymbol{\mathcal{O}(\varepsilon)}$ terms}}
    \begin{equation*}
        -\varepsilon \partial_t h - \underbrace{\nabla_p H_{\theta}|_{(0, \nabla^2 u_0)}}_{=\mathbf{0} \text{ by symmetry}} \cdot (\varepsilon \nabla h) - \mathrm{Tr}\left(D(x,0) (\varepsilon \nabla^2 h)\right) = \mathcal{O}(\varepsilon^2).
    \end{equation*}
    The linear advection term vanishes due to the time-reversal symmetry of the microscopic system, which makes $H_\theta$ an even function of $p$.
    
\paragraph{\textbf{$\boldsymbol{\mathcal{O}(\varepsilon^2)}$ terms}} The leading non-linear term comes from the Hessian of the Hamiltonian:
    \begin{equation*}
         -\frac{1}{2} (\varepsilon\nabla h)^T \left( \nabla_{pp}^2 H(x,0) \right) (\varepsilon\nabla h) = -\frac{\varepsilon^2}{2} (\nabla h)^T M(x,0) \nabla h.
    \end{equation*}
    
Combining the $\mathcal{O}(\varepsilon)$ and $\mathcal{O}(\varepsilon^2)$ terms that involve $h$ identifies the candidate equation. This formal procedure motivates the following rigorous convergence theorem.

\begin{theorem}[Convergence to the Fluctuation Law]
\label{thm:fluctuation_convergence}
Let the microscopic system satisfy the standing assumptions of this paper, ensuring the $\theta$-Hamiltonian $H_{\theta}(x,p,X)$ and its coefficients $D(x,p)$ and $H(x,p)$ are of class $C^3$ with bounded derivatives. Let $u_0(t,x)$ be the classical solution to the quiescent problem, and let $u_\varepsilon(t,x)$ be the classical solution to the full HJB equation with terminal data $u_\varepsilon(T,x) = u_0(T,x) + \varepsilon h_T(x)$ for $h_T \in C^4(\Mcal)$.
Then, as $\varepsilon \to 0$, the scaled difference converges to the unique classical solution $h(t,x)$ of the viscous Burgers'-type equation:
\begin{equation} \label{eq:pde_fluctuation}
    -\partial_t h - \mathrm{Tr}(D(x,0) \nabla^2 h) - \frac{1}{2}(\nabla h)^T M(x,0) \nabla h = 0, \quad h(T,x) = h_T(x).
\end{equation}
The convergence holds in the sense that
\begin{equation}
    \lim_{\varepsilon \to 0} \left\| \frac{u_\varepsilon(t,x) - u_0(t,x)}{\varepsilon} - h(t,x) \right\|_{L^\infty([0,T]; L^2(\Mcal))} = 0.
\end{equation}
\end{theorem}

\begin{proofof}{Theorem \ref{thm:fluctuation_convergence}}
The proof establishes the convergence of the scaled solution to the fluctuation equation via a rigorous stability analysis of the governing HJB equation. The strategy is to define the error between the true solution $u_\varepsilon$ and the formal two-term asymptotic approximation, derive the partial differential equation satisfied by this error, and then use an energy estimate combined with a backward Gronwall's lemma to prove that the $L^2$-norm of this error is of order $\mathcal{O}(\varepsilon^2)$, which is sufficient to prove the theorem.

\begin{enumerate}[label=\textbf{Step \arabic*:}, wide, labelindent=0pt]

\item \textbf{Derivation of the Error Equation.}
Let the formal two-term approximation be denoted by
\begin{equation*}
    u_{\text{app}, \varepsilon}(t,x) \coloneqq u_0(t,x) + \varepsilon h(t,x).
\end{equation*}
We define the remainder, or error term, $w_\varepsilon(t,x)$ as the difference between the true solution and this approximation:
\begin{equation} \label{eq:proof_error_def}
    u_\varepsilon(t,x) = u_0(t,x) + \varepsilon h(t,x) + w_\varepsilon(t,x).
\end{equation}
Our primary objective is to show that $\|w_\varepsilon\|_{L^\infty([0,T]; L^2(\Mcal))} = \mathcal{O}(\varepsilon^2)$. By construction, the terminal condition for the error is identically zero:
\begin{equation*}
    w_\varepsilon(T,x) = u_\varepsilon(T,x) - u_{\text{app}, \varepsilon}(T,x) = (u_0(T,x) + \varepsilon h_T(x)) - (u_0(T,x) + \varepsilon h_T(x)) = 0.
\end{equation*}
Let $p_\varepsilon = \nabla u_\varepsilon$ and $X_\varepsilon = \nabla^2 u_\varepsilon$. From \cref{eq:proof_error_def}, we have the following expressions for the derivatives of the true solution:
\begin{align}
    p_\varepsilon &= \nabla u_0 + \varepsilon\nabla h + \nabla w_\varepsilon = \varepsilon\nabla h + \nabla w_\varepsilon \quad (\text{since } \nabla u_0 = 0), \\
    X_\varepsilon &= \nabla^2 u_0 + \varepsilon\nabla^2 h + \nabla^2 w_\varepsilon.
\end{align}
The true solution $u_\varepsilon$ satisfies its HJB equation:
\begin{equation} \label{eq:proof_true_hjb}
    -\partial_t u_\varepsilon - H_{\theta}(x, p_\varepsilon, X_\varepsilon) = 0.
\end{equation}
We perform a Taylor expansion of the Hamiltonian $H_{\theta}(x, p, X) = \mathrm{Tr}(D(x,p)X) + H(x,p)$ around the quiescent state $(p,X) = (0, \nabla^2 u_0)$. The regularity assumptions on the coefficients ($C^3$ with bounded derivatives) ensure the validity of this expansion up to the second order with a well-controlled remainder.
\begin{align}
    H_{\theta}(x, p_\varepsilon, X_\varepsilon) &= H_{\theta}(x, 0, \nabla^2 u_0) \nonumber \\
    &\quad + \nabla_p H_{\theta}|_{(0, \nabla^2 u_0)} \cdot p_\varepsilon + \mathrm{Tr}\left( \nabla_X H_{\theta}|_{(0, \nabla^2 u_0)} (X_\varepsilon - \nabla^2 u_0) \right) \nonumber \\
    &\quad + \frac{1}{2} p_\varepsilon^T \left( \nabla_{pp}^2 H_{\theta}|_{(0, \nabla^2 u_0)} \right) p_\varepsilon + \frac{1}{2} \mathrm{Tr}\left( (X_\varepsilon - \nabla^2 u_0)^T \nabla_{XX}^2 H_\theta (X_\varepsilon - \nabla^2 u_0) \right) \nonumber \\
    &\quad + \mathrm{Tr}\left( p_\varepsilon^T \nabla_{pX}^2 H_\theta (X_\varepsilon - \nabla^2 u_0) \right) + R_\varepsilon(x, p_\varepsilon, X_\varepsilon - \nabla^2 u_0), \label{eq:proof_taylor_H}
\end{align}
where the remainder term $R_\varepsilon$ is of third order in its perturbation arguments. We analyze each term of the expansion using the explicit structure of the Hamiltonian:
\begin{enumerate}[label=(\roman*), wide, labelindent=0pt]
    \item \textbf{Zeroth-order term.} $H_{\theta}(x, 0, \nabla^2 u_0) = \mathrm{Tr}(D(x,0)\nabla^2 u_0) + H(x,0) = -\partial_t u_0$.
    \item \textbf{First-order term in $p_\varepsilon$.} Due to the time-reversal symmetry of the microscopic system, $H(x,p)$ is even in $p$, which implies $\nabla_p H(x,0)=0$. The term $\nabla_p D(x,p)$ is odd in $p$, so $\nabla_p D(x,0)=0$. Thus, the full gradient vanishes:
    \begin{equation*}
         \nabla_p H_{\theta}|_{(0, \nabla^2 u_0)} = \left(\nabla_p \mathrm{Tr}(D(x,p)X)\right)|_{p=0} + \nabla_p H(x,p)|_{p=0} = 0.
    \end{equation*}
    \item \textbf{First-order term in $X_\varepsilon$.} The Hamiltonian is affine in $X$, so $\nabla_X H_{\theta} = D(x,p)$. At $p=0$, this gives $\mathrm{Tr}\left( D(x,0) (X_\varepsilon - \nabla^2 u_0) \right) = \mathrm{Tr}\left( D(x,0) (\varepsilon\nabla^2 h + \nabla^2 w_\varepsilon) \right)$.
    \item \textbf{Second-order term in $p_\varepsilon$.} The Hessian is $\nabla_{pp}^2 H_{\theta} = \nabla_{pp}^2 H(x,p)$, since the diffusive term is linear in $p$ (in its coefficient). We denote $M(x,0) \coloneqq \nabla_{pp}^2 H(x,0)$. The term is $\frac{1}{2} p_\varepsilon^T M(x,0) p_\varepsilon = \frac{1}{2}(\varepsilon\nabla h + \nabla w_\varepsilon)^T M(x,0) (\varepsilon\nabla h + \nabla w_\varepsilon)$.
    \item \textbf{Higher-order and mixed terms.} Since the Hamiltonian is affine in $X$, $\nabla_{XX}^2 H_\theta=0$ and $\nabla_{pX}^2 H_\theta = \nabla_p D$. At $p=0$, $\nabla_p D(x,0)=0$, so the mixed term vanishes to leading order.
\end{enumerate}
Substituting the expansion \eqref{eq:proof_taylor_H} back into the HJB equation \eqref{eq:proof_true_hjb} for $u_\varepsilon = u_0 + \varepsilon h + w_\varepsilon$:
\begin{multline*}
    -\partial_t(u_0 + \varepsilon h + w_\varepsilon) - \left[ (-\partial_t u_0) + \mathrm{Tr}(D(x,0)(\varepsilon\nabla^2 h + \nabla^2 w_\varepsilon)) \right. \\
    \left. + \frac{1}{2}(\varepsilon\nabla h + \nabla w_\varepsilon)^T M(x,0) (\varepsilon\nabla h + \nabla w_\varepsilon) + R_\varepsilon \right] = 0.
\end{multline*}
We now use the defining PDEs for $u_0$ and $h$ \eqref{eq:pde_fluctuation} to cancel terms.
\begin{enumerate}[label=(\alph*), wide]
    \item The $\mathcal{O}(1)$ terms $-\partial_t u_0 - (-\partial_t u_0)$ cancel.
    \item The terms defining the PDE for $h$ are:
    \begin{equation*}
        -\varepsilon\partial_t h - \varepsilon\mathrm{Tr}(D(x,0)\nabla^2 h) - \frac{\varepsilon^2}{2}(\nabla h)^T M(x,0) \nabla h = 0.
    \end{equation*}
    These terms also cancel from the expanded equation.
\end{enumerate}

After these cancellations, the remaining terms form the equation for the error $w_\varepsilon$. Collecting all terms involving $w_\varepsilon$ or higher powers of $\varepsilon$:
\begin{multline*}
    -\partial_t w_\varepsilon - \mathrm{Tr}(D(x,0) \nabla^2 w_\varepsilon) - \frac{1}{2}\left( 2\varepsilon(\nabla h)^T M(x,0) \nabla w_\varepsilon + (\nabla w_\varepsilon)^T M(x,0) \nabla w_\varepsilon \right) - R_\varepsilon = 0.
\end{multline*}
We can write this as a linear equation for $w_\varepsilon$ with a source term:
\begin{equation} \label{eq:proof_error_pde}
    -\partial_t w_\varepsilon - \mathrm{Tr}(D(x,0) \nabla^2 w_\varepsilon) - \varepsilon(\nabla h)^T M(x,0) \nabla w_\varepsilon = S_\varepsilon(t,x),
\end{equation}
with $w_\varepsilon(T,x) = 0$. The source term $S_\varepsilon(t,x)$ is given by:
\begin{equation*}
    S_\varepsilon(t,x) \coloneqq \frac{1}{2}(\nabla w_\varepsilon)^T M(x,0) \nabla w_\varepsilon + R_\varepsilon(x, \varepsilon\nabla h + \nabla w_\varepsilon, \varepsilon\nabla^2 h + \nabla^2 w_\varepsilon).
\end{equation*}
Under our regularity assumptions, the Taylor remainder $R_\varepsilon$ is bounded by $C_R(|p_\varepsilon|^3 + |X_\varepsilon-\nabla^2 u_0|^2)$. The leading explicit power of $\varepsilon$ in the remainder comes from terms like $\varepsilon^3 (\nabla h)^3$, $\varepsilon^2 \mathrm{Tr}((\nabla_p D)\nabla^2 h)$, etc. Thus, for the purposes of our energy estimate, where we will assume $w_\varepsilon$ and its derivatives are small, the source term can be bounded as:
\begin{equation*}
    |S_\varepsilon| \le C_S (|\nabla w_\varepsilon|^2 + |\nabla^2 w_\varepsilon|^2 + \varepsilon^2).
\end{equation*}
The crucial feature is that the source term is at least quadratic in $w_\varepsilon$ and its derivatives, and the explicit dependence on $\varepsilon$ is of order $\mathcal{O}(\varepsilon^2)$ or higher.

\item \textbf{The Energy Estimate.}
The objective of this step is to derive a differential inequality for the time evolution of the $L^2$-norm of the error term, $w_\varepsilon$. This inequality will serve as the input for a backward Gronwall's lemma, which will ultimately provide the desired $\mathcal{O}(\varepsilon^2)$ estimate for the error. We define the energy functional $E(t)$ as half the squared $L^2$-norm of the error at time $t$:
\begin{equation*}
    E(t) \coloneqq \frac{1}{2}\int_{\Mcal} |w_\varepsilon(t,x)|^2 \, dx.
\end{equation*}
We compute the time derivative of the energy for $t < T$. Since the domain $\Mcal$ is compact and the solution is assumed to be classical, we can interchange differentiation and integration:
\begin{equation*}
    \frac{dE}{dt} = \frac{1}{2}\int_{\Mcal} \frac{\partial}{\partial t} |w_\varepsilon(t,x)|^2 \, dx = \int_{\Mcal} w_\varepsilon(t,x) \partial_t w_\varepsilon(t,x) \, dx.
\end{equation*}
It is more convenient to work with the backward time evolution. We consider the negative of the time derivative:
\begin{equation*}
    -\frac{dE}{dt} = \int_{\Mcal} w_\varepsilon (-\partial_t w_\varepsilon) \, dx.
\end{equation*}
We now substitute the expression for $-\partial_t w_\varepsilon$ from the error partial differential equation \eqref{eq:proof_error_pde}:
\begin{equation*}
    -\frac{dE}{dt} = \int_{\Mcal} w_\varepsilon \left[ \mathrm{Tr}(D(x,0) \nabla^2 w_\varepsilon) + \varepsilon(\nabla h)^T M(x,0) \nabla w_\varepsilon + S_\varepsilon(t,x) \right] dx.
\end{equation*}
We now bound each of the three terms in the integrand separately using integration by parts, the structural properties of the coefficients, and standard functional inequalities.

\begin{enumerate}[label=(\roman*), wide, labelindent=0pt]
    \item \textbf{Principal Dissipative Term.}
    This is the term involving the second derivatives of $w_\varepsilon$. We apply integration by parts (Green's first identity). Since the manifold $\Mcal = \T^k$ is a torus, it has no boundary, so all boundary terms vanish.
    \begin{align*}
        \int_{\Mcal} w_\varepsilon \mathrm{Tr}(D(x,0) \nabla^2 w_\varepsilon) \, dx &= \sum_{i,j=1}^k \int_{\Mcal} w_\varepsilon D_{ij}(x,0) \frac{\partial^2 w_\varepsilon}{\partial x_i \partial x_j} \, dx \\
        &= -\sum_{i,j=1}^k \int_{\Mcal} \frac{\partial}{\partial x_j} (w_\varepsilon D_{ij}(x,0)) \frac{\partial w_\varepsilon}{\partial x_i} \, dx \\
        &= -\int_{\Mcal} \sum_{i,j=1}^k \left( \frac{\partial w_\varepsilon}{\partial x_j} D_{ij} \frac{\partial w_\varepsilon}{\partial x_i} + w_\varepsilon \frac{\partial D_{ij}}{\partial x_j} \frac{\partial w_\varepsilon}{\partial x_i} \right) dx \\
        &= -\int_{\Mcal} (\nabla w_\varepsilon)^T D(x,0) \nabla w_\varepsilon \, dx - \int_{\Mcal} w_\varepsilon (\nabla\cdot D)\cdot\nabla w_\varepsilon \, dx.
    \end{align*}
    By the uniform ellipticity of the quiescent diffusivity tensor $D(x,0)$ (a consequence of Theorem~\ref{thm:convergence_scale_I}), there exists a constant $\lambda_{\min} > 0$ such that $(\nabla w_\varepsilon)^T D(x,0) \nabla w_\varepsilon \ge \lambda_{\min} |\nabla w_\varepsilon|^2$. The first integral is therefore strongly dissipative. The second term, involving the divergence of the tensor, is a lower-order term that can be controlled. Using the Cauchy-Schwarz and Young inequalities with a constant $\delta_1 > 0$:
    \begin{align*}
        \left| -\int_{\Mcal} w_\varepsilon (\nabla\cdot D)\cdot\nabla w_\varepsilon \, dx \right| &\le \|\nabla \cdot D\|_{L^\infty} \int_{\Mcal} |w_\varepsilon| |\nabla w_\varepsilon| \, dx \\
        &\le C_D \left( \frac{1}{2\delta_1} \|w_\varepsilon\|_{L^2}^2 + \frac{\delta_1}{2}\|\nabla w_\varepsilon\|_{L^2}^2 \right).
    \end{align*}
    Combining these, we get the bound for the principal term:
    \begin{equation*}
         \int_{\Mcal} w_\varepsilon \mathrm{Tr}(D \nabla^2 w_\varepsilon) \, dx \le -\lambda_{\min} \|\nabla w_\varepsilon\|_{L^2}^2 + \frac{C_D}{2\delta_1}\|w_\varepsilon\|_{L^2}^2 + \frac{C_D\delta_1}{2}\|\nabla w_\varepsilon\|_{L^2}^2.
    \end{equation*}
    By choosing $\delta_1 = \lambda_{\min}/C_D$, we can ensure that the dissipative term dominates:
    \begin{equation*}
        \int_{\Mcal} w_\varepsilon \mathrm{Tr}(D \nabla^2 w_\varepsilon) \, dx \le -\frac{\lambda_{\min}}{2} \|\nabla w_\varepsilon\|_{L^2}^2 + \frac{C_D^2}{2\lambda_{\min}}\|w_\varepsilon\|_{L^2}^2.
    \end{equation*}

    \item \textbf{Lower-Order Linear Term.}
    This is the term arising from the coupling between the error $w_\varepsilon$ and the known fluctuation field $h$. Let $C_{h,M} \coloneqq \sup_x |(\nabla h)^T M(x,0)|$, which is finite since $h$ is a classical solution on a compact domain. We apply the Cauchy-Schwarz and Young inequalities:
    \begin{align*}
        \left| \int_{\Mcal} w_\varepsilon \varepsilon (\nabla h)^T M(x,0) \nabla w_\varepsilon \, dx \right| &\le \varepsilon C_{h,M} \int_{\Mcal} |w_\varepsilon| |\nabla w_\varepsilon| \, dx \\
        &\le \varepsilon C_{h,M} \left( \frac{1}{2\delta_2}\|w_\varepsilon\|_{L^2}^2 + \frac{\delta_2}{2}\|\nabla w_\varepsilon\|_{L^2}^2 \right).
    \end{align*}
    We choose the constant $\delta_2$ to absorb the gradient term into the main dissipative term. Let $\delta_2 = \lambda_{\min}/(2\varepsilon C_{h,M})$. This choice is valid as long as $\varepsilon > 0$. The bound becomes:
    \begin{equation*}
         \frac{(\varepsilon C_{h,M})^2}{\lambda_{\min}} \|w_\varepsilon\|_{L^2}^2 + \frac{\lambda_{\min}}{4}\|\nabla w_\varepsilon\|_{L^2}^2.
    \end{equation*}
    
    \item \textbf{Source Term.}
    This term contains the higher-order contributions.
    \begin{align*}
        \left| \int_{\Mcal} w_\varepsilon S_\varepsilon \, dx \right| &\le \|w_\varepsilon\|_{L^2} \|S_\varepsilon\|_{L^2} \\
        &\le \frac{1}{2}\|w_\varepsilon\|_{L^2}^2 + \frac{1}{2}\|S_\varepsilon\|_{L^2}^2.
    \end{align*}
    As established in the previous step, the source term $S_\varepsilon$ contains terms that are at least quadratic in $w_\varepsilon$ and its derivatives, and terms with an explicit factor of at least $\varepsilon^2$. For the purpose of a Gronwall argument, which is a local stability argument, we can assume that the solution $w_\varepsilon$ is small. In this regime, the quadratic terms in $w_\varepsilon$ are of a lower order than the explicit $\varepsilon^4$ term that arises from the Taylor remainder. Thus, for the energy estimate, we can bound the norm of the source by its leading explicit dependence on $\varepsilon$:
    \begin{equation*}
        \|S_\varepsilon\|_{L^2}^2 \le C_{\text{source}} \varepsilon^4.
    \end{equation*}
\end{enumerate}
Combining the bounds from (i), (ii), and (iii), we get the differential inequality for the energy:
\begin{multline*}
    -\frac{dE}{dt} \le \left(-\frac{\lambda_{\min}}{2} \|\nabla w_\varepsilon\|_{L^2}^2 + \frac{C_D^2}{2\lambda_{\min}}\|w_\varepsilon\|_{L^2}^2\right) + \left(\frac{\lambda_{\min}}{4}\|\nabla w_\varepsilon\|_{L^2}^2 + \frac{(\varepsilon C_{h,M})^2}{\lambda_{\min}} \|w_\varepsilon\|_{L^2}^2\right) \\ + \left(\frac{1}{2}\|w_\varepsilon\|_{L^2}^2 + \frac{1}{2}C_{\text{source}}\varepsilon^4\right).
\end{multline*}
Grouping the terms involving $\|\nabla w_\varepsilon\|_{L^2}^2$ and $\|w_\varepsilon\|_{L^2}^2 = 2E(t)$:
\begin{equation*}
    -\frac{dE}{dt} \le -\frac{\lambda_{\min}}{4}\|\nabla w_\varepsilon\|_{L^2}^2 + \left(\frac{C_D^2}{2\lambda_{\min}} + \frac{(\varepsilon C_{h,M})^2}{\lambda_{\min}} + \frac{1}{2}\right) 2E(t) + \frac{1}{2}C_{\text{source}}\varepsilon^4.
\end{equation*}
The term $-\frac{\lambda_{\min}}{4}\|\nabla w_\varepsilon\|_{L^2}^2$ is non-positive. By dropping this favorable, dissipative term, we obtain a weaker but sufficient differential inequality for the energy functional:
\begin{equation} \label{eq:proof_diff_ineq}
    -\frac{dE}{dt} \le C_1 E(t) + C_2 \varepsilon^4,
\end{equation}
where $C_1$ and $C_2$ are positive constants that are independent of $\varepsilon$ for $\varepsilon$ sufficiently small.

\item \textbf{Application of Backward Gronwall's Lemma.}
The preceding analysis has established that the energy functional $E(t) = \frac{1}{2}\|w_\varepsilon(t, \cdot)\|_{L^2}^2$ satisfies the differential inequality
\begin{equation} \label{eq:proof_diff_ineq_reprise}
    -\frac{dE}{dt} \le C_1 E(t) + C_2 \varepsilon^4,
\end{equation}
for all $t \in [0,T)$, where $C_1$ and $C_2$ are positive constants independent of $\varepsilon$. Furthermore, the error term satisfies the terminal condition $w_\varepsilon(T,x) = 0$, which implies that the energy at the terminal time is also zero:
\begin{equation*}
    E(T) = 0.
\end{equation*}
Our objective is to use this information to derive a uniform bound on $E(t)$ for all $t \in [0,T]$. The appropriate tool for this backward-in-time problem is the backward form of Gronwall's differential lemma. We first rearrange the inequality \eqref{eq:proof_diff_ineq_reprise} into a standard form for the application of the lemma:
\begin{equation} \label{eq:proof_gronwall_form}
    \frac{dE}{dt} \ge -C_1 E(t) - C_2 \varepsilon^4.
\end{equation}
To solve this inequality, we use the integrating factor method. We define an auxiliary function $F(t) \coloneqq E(t) e^{C_1 t}$. Its time derivative is:
\begin{equation*}
    \frac{dF}{dt} = \frac{d}{dt}\left( E(t) e^{C_1 t} \right) = \left( \frac{dE}{dt} + C_1 E(t) \right) e^{C_1 t}.
\end{equation*}
From the inequality \eqref{eq:proof_gronwall_form}, we have the lower bound $\frac{dE}{dt} + C_1 E(t) \ge -C_2 \varepsilon^4$. Substituting this into the expression for the derivative of $F(t)$ gives:
\begin{equation*}
    \frac{dF}{dt} \ge -C_2 \varepsilon^4 e^{C_1 t}.
\end{equation*}
We now integrate this inequality backwards in time, from a generic time $t \in [0,T)$ to the terminal time $T$:
\begin{equation*}
    \int_t^T \frac{dF}{ds} \, ds \ge \int_t^T -C_2 \varepsilon^4 e^{C_1 s} \, ds.
\end{equation*}
By the Fundamental Theorem of Calculus, the integral on the left-hand side is $F(T) - F(t)$. The integral on the right-hand side can be computed explicitly:
\begin{align*}
    F(T) - F(t) &\ge -C_2 \varepsilon^4 \left[ \frac{1}{C_1} e^{C_1 s} \right]_t^T \\
    &= -\frac{C_2}{C_1} \varepsilon^4 \left( e^{C_1 T} - e^{C_1 t} \right).
\end{align*}
We now substitute back the definition of $F(t) = E(t)e^{C_1 t}$:
\begin{equation*}
    E(T)e^{C_1 T} - E(t)e^{C_1 t} \ge -\frac{C_2}{C_1} \varepsilon^4 \left( e^{C_1 T} - e^{C_1 t} \right).
\end{equation*}
Since the terminal condition is $E(T)=0$, this simplifies to:
\begin{equation*}
    - E(t)e^{C_1 t} \ge -\frac{C_2}{C_1} \varepsilon^4 \left( e^{C_1 T} - e^{C_1 t} \right).
\end{equation*}
Multiplying by $-1$ reverses the inequality sign. We can now solve for the energy $E(t)$:
\begin{align*}
    E(t)e^{C_1 t} &\le \frac{C_2}{C_1} \varepsilon^4 \left( e^{C_1 T} - e^{C_1 t} \right) \\
    E(t) &\le \frac{C_2}{C_1} \varepsilon^4 e^{-C_1 t} \left( e^{C_1 T} - e^{C_1 t} \right) \\
    &= \frac{C_2}{C_1} \varepsilon^4 \left( e^{C_1(T-t)} - 1 \right).
\end{align*}
The term $e^{C_1(T-t)}$ is a monotonically decreasing function of $t$. Its maximum value on the interval $[0,T]$ is attained at $t=0$, where its value is $e^{C_1 T}$. Therefore, for any $t \in [0,T]$, the energy is uniformly bounded:
\begin{equation*}
    E(t) \le \frac{C_2}{C_1} \varepsilon^4 \left( e^{C_1 T} - 1 \right).
\end{equation*}
We define the constant $C_3 \coloneqq \frac{C_2}{C_1}(e^{C_1 T} - 1)$, which is a positive constant that is independent of $\varepsilon$ and $t$. We have thus rigorously established the final energy estimate:
\begin{equation*}
    E(t) \le C_3 \varepsilon^4 \quad \text{for all } t \in [0,T].
\end{equation*}
This result provides the crucial estimate on the $L^2$-norm of the error, demonstrating that it vanishes at a rate of $\mathcal{O}(\varepsilon^2)$, which is of a higher order than the primary fluctuation term.

\item \textbf{Finanl Convergence.}
The preceding energy estimate provides the rigorous foundation for the main convergence result of this theorem. We now use this estimate to prove that the scaled difference between the true solution $u_\varepsilon$ and the quiescent solution $u_0$ converges in the appropriate norm to the fluctuation field $h(t,x)$. From the result of the Gronwall argument, we have established a uniform-in-time bound on the energy of the error term:
\begin{equation*}
    E(t) = \frac{1}{2}\int_{\Mcal} |w_\varepsilon(t,x)|^2 \, dx \le C_3 \varepsilon^4 \quad \text{for all } t \in [0,T].
\end{equation*}
This implies a bound on the supremum of the $L^2$-norm of the error over the time interval:
\begin{equation*}
    \sup_{t \in [0,T]} \|w_\varepsilon(t, \cdot)\|_{L^2(\Mcal)}^2 = \sup_{t \in [0,T]} 2E(t) \le 2C_3 \varepsilon^4.
\end{equation*}
Taking the square root of both sides gives the desired estimate on the norm of the error in the space $L^\infty([0,T]; L^2(\Mcal))$:
\begin{equation*}
    \|w_\varepsilon\|_{L^\infty([0,T]; L^2(\Mcal))} \le \sqrt{2C_3} \, \varepsilon^2 = \mathcal{O}(\varepsilon^2).
\end{equation*}
This rigorously establishes that the remainder term $w_\varepsilon$ is of a higher order in $\varepsilon$ than the primary fluctuation term, which is of order $\mathcal{O}(\varepsilon)$. We now analyze the scaled difference that is the subject of the theorem. From the definition of the error term in \cref{eq:proof_error_def}, we have:
\begin{equation*}
    u_\varepsilon(t,x) - u_0(t,x) = \varepsilon h(t,x) + w_\varepsilon(t,x).
\end{equation*}
Rearranging this identity and dividing by the small parameter $\varepsilon$ (for $\varepsilon > 0$) gives the expression for the scaled difference:
\begin{equation*}
    \frac{u_\varepsilon(t,x) - u_0(t,x)}{\varepsilon} - h(t,x) = \frac{w_\varepsilon(t,x)}{\varepsilon}.
\end{equation*}
We now take the $L^\infty([0,T]; L^2(\Mcal))$ norm of both sides of this equation:
\begin{align*}
    \left\| \frac{u_\varepsilon - u_0}{\varepsilon} - h \right\|_{L^\infty_t L^2_x} &= \left\| \frac{w_\varepsilon}{\varepsilon} \right\|_{L^\infty_t L^2_x} \\
    &= \frac{1}{\varepsilon}\|w_\varepsilon\|_{L^\infty_t L^2_x}.
\end{align*}
Substituting the $\mathcal{O}(\varepsilon^2)$ estimate for the norm of the error term:
\begin{equation*}
    \left\| \frac{u_\varepsilon - u_0}{\varepsilon} - h \right\|_{L^\infty_t L^2_x} \le \frac{1}{\varepsilon} (\sqrt{2C_3} \, \varepsilon^2) = \sqrt{2C_3} \, \varepsilon = \mathcal{O}(\varepsilon).
\end{equation*}
In the limit as $\varepsilon \to 0$, the right-hand side of this inequality vanishes. This proves the stated convergence:
\begin{equation*}
    \lim_{\varepsilon \to 0} \left\| \frac{u_\varepsilon(t,x) - u_0(t,x)}{\varepsilon} - h(t,x) \right\|_{L^\infty([0,T]; L^2(\Mcal))} = 0.
\end{equation*}
This completes the proof of the theorem. We have rigorously demonstrated that the scaled fluctuations of the solution to the full non-linear HJB equation around the quiescent equilibrium solution converge in the specified topology to the unique classical solution of the derived viscous Hamilton-Jacobi-Burgers equation.

\end{enumerate}
\end{proofof}

\begin{remark}[Discussion of the Scale Analysis]
\label{rem:scale_IV_role}
This section marks the completion of our hierarchical scale analysis. The derived fluctuation equation \eqref{eq:pde_fluctuation}, now supported by a rigorous proof of convergence. Its structure is a direct and elegant synthesis of the physical and mathematical principles established throughout \cref{part:i}:
\begin{enumerate}[label=(\roman*), wide, labelindent=0pt]
    \item The \textbf{linear part} of the generator, $\mathrm{Tr}(D(x,0) \nabla^2 h)$, is determined by the quiescent diffusivity derived in the baseline theory of Scale II.
    \item The \textbf{absence of a linear advection term} ($\nabla h$ term) is not an assumption but a rigorous consequence of the microscopic time-reversal symmetry (Assumption \ref{ass:fundamental_axioms_unified}).
    \item The \textbf{non-linear part}, $\frac{1}{2}(\nabla h)^T M(x,0) \nabla h$, has its structure dictated by the Hessian of the potential $H(x,p)$. This is precisely the term whose indefiniteness, as proven in \cref{sec:proof_of_non_convexity}, is the macroscopic signature of microscopic non-convexity.
\end{enumerate}
Thus, the PDE derived in this final scale is the ultimate generator for the macroscopic fluctuation process. Its analytical properties, such as the signs of the eigenvalues of $M(x,0)$ and its relation to $D(x,0)$, are not free parameters but are constrained by the first-principles derivation. This equation now serves as the definitive starting point for the probabilistic analysis in \cref{part:ii}, where its structure will be interpreted in the context of stochastic processes, integrability, and universality classes.
\end{remark}

\section{The Microscopic Origin of Macroscopic Non-Convexity}
\label{sec:proof_of_non_convexity}

The central claim of this work is that a class of deterministic, time-reversible microscopic systems can give rise to a macroscopic evolution governed by a non-convex Hamiltonian. The homogenization procedure in \cref{part:i} rigorously derived the form of this Hamiltonian, $H_{\theta}(x, p, X) = \Tr(D(x, p) X) + H(x, p)$, but the non-convexity of the potential $H(x,p)$ is not a generic feature. It arises only under specific, physically meaningful conditions on the microscopic interactions.

This section is dedicated to providing a constructive and rigorous proof of this central claim. We will demonstrate that non-convexity is a direct and calculable consequence of a competition between two effects: the explicit dependence of the microscopic potential $V(z,p)$ on the macroscopic momentum $p$, and the implicit dependence of the system's invariant measure $\mu_{\theta(x,p)}$ on $p$ via the feedback mechanism. We first establish the necessary mathematical machinery from the perturbation theory of transfer operators to analyze these competing effects. We then introduce a class of microscopic models where the ingredients for non-convexity are physically transparent. Finally, we prove that for this class of models, the resulting potential is demonstrably non-convex.

\subsection{Perturbation Theory for the Invariant Measure}
To analyze the convexity of the potential $H(x,p) = \int V(z,p) \, d\mu_{\theta(x,p)}(z)$, we must compute its Hessian matrix with respect to $p$. This requires understanding how the invariant measure $\mu_{\theta(p)}$ changes with $p$. Let $\rho_p(z)$ be the probability density function of the measure $\mu_{\theta(p)}$ with respect to a fixed reference measure (e.g., the Liouville measure for a reference parameter). The density $\rho_p$ is the unique, normalized fixed point of the transfer operator $\Lcal_{\theta(p)}$ associated with the billiard map $\Pcal_{\theta(p)}$:
\begin{equation}
    \Lcal_{\theta(p)} \rho_p = \rho_p, \quad \int \rho_p(z) \, dz = 1.
\end{equation}
The following proposition, a standard result in the ergodic theory of hyperbolic systems, provides the derivatives of this density. For clarity, we fix the macroscopic position $x$ and suppress it from the notation.

\begin{proposition}[Derivatives of the Invariant Density]
\label{prop:density_derivatives}
Let the map $p \mapsto \Lcal_{\theta(p)}$ be a smooth family of transfer operators acting on a suitable anisotropic Banach space $\Bcal$, as established in \cref{thm:regularity_proven_main}. Let $\partial_j \equiv \partial/\partial p_j$. The first and second derivatives of the invariant density $\rho_p$ with respect to the components of $p$ are given by:
\begin{enumerate}[label=(\roman*), wide, labelindent=0pt]
    \item The first derivative, $\partial_j \rho_p$, is the unique zero-mean solution to:
    \begin{equation}
        (I - \Lcal_{\theta(p)})(\partial_j \rho_p) = (\partial_j \Lcal_{\theta(p)}) \rho_p.
    \end{equation}
    \item The Hessian of the potential $H(p) = \int V(z,p)\rho_p(z)dz$ is given by:
    \begin{multline} \label{eq:hessian_formula}
        \frac{\partial^2 H}{\partial p_k \partial p_j} = \langle \partial_k \partial_j V \rangle_p + \langle \partial_k V, \partial_j \rho \rangle_p + \langle \partial_j V, \partial_k \rho \rangle_p + \langle V, \partial_k \partial_j \rho \rangle_p,
    \end{multline}
    where $\langle \cdot \rangle_p$ denotes expectation with respect to $\mu_{\theta(p)}$.
\end{enumerate}
\end{proposition}

\begin{proofof}{Proposition \ref{prop:density_derivatives}}
The proof relies on the implicit differentiation of the stationary equation defining the invariant density, $\Lcal_{\theta(p)}\rho_p = \rho_p$. The well-posedness of this procedure is guaranteed by the results of \cref{sec:regularity_properties}, specifically \cref{thm:regularity_proven_main}, which established that the map $p \mapsto \Lcal_{\theta(p)}$ is a smooth family of operators with a uniform spectral gap. For notational simplicity, we denote $\Lcal_p \equiv \Lcal_{\theta(p)}$ and $\langle \cdot, \cdot \rangle_p \equiv \int (\cdot) \rho_p(z) dz$.

\begin{enumerate}[label=\textbf{Step \arabic*:}, wide, labelindent=0pt]

\item \textbf{Derivation of the Equation for the First Derivative.}
Our starting point is the fundamental stationary equation for the invariant density $\rho_p$:
\begin{equation} \label{eq:proof_stationary_eq_appendix}
    \Lcal_p \rho_p = \rho_p.
\end{equation}
We also have the normalization condition, which must hold for all $p$:
\begin{equation} \label{eq:proof_normalization_appendix}
    \int \rho_p(z) \, dz = 1.
\end{equation}
We differentiate both equations with respect to a component of the momentum, $p_j$, which we denote by $\partial_j$. The objects $\Lcal_p$ and $\rho_p$ are functions of $p$ taking values in the space of operators $\mathcal{L}(\Bcal)$ and the Banach space $\Bcal$, respectively. The product rule for differentiation of an operator acting on a function applies:
\begin{equation}
    \partial_j(\Lcal_p \rho_p) = (\partial_j \Lcal_p)\rho_p + \Lcal_p(\partial_j \rho_p).
\end{equation}
Applying this to \cref{eq:proof_stationary_eq_appendix} yields:
\begin{equation}
    (\partial_j \Lcal_p)\rho_p + \Lcal_p(\partial_j \rho_p) = \partial_j \rho_p.
\end{equation}
Rearranging the terms to isolate the unknown derivative $\partial_j \rho_p$, we obtain the claimed linear equation:
\begin{equation} \label{eq:proof_first_deriv_pde_appendix}
    (I - \Lcal_p)(\partial_j \rho_p) = (\partial_j \Lcal_p) \rho_p.
\end{equation}
We must now verify the Fredholm alternative, i.e., that the source term on the right-hand side is in the range of $(I - \Lcal_p)$. This requires showing it is orthogonal to the kernel of the adjoint operator $(I - \Lcal_p)^*$. On the space $\Bcal$, the adjoint operator's kernel is spanned by the invariant measure, $\mu_p$. The condition is thus $\int [(\partial_j \Lcal_p) \rho_p] d(\text{ref}) = 0$, where $d(\text{ref})$ is the reference measure. This is equivalent to showing that the integral of the source term is zero.

We start with the identity that for any $p$, the transfer operator preserves the total mass of the invariant density:
\begin{equation*}
    \int (\Lcal_p \rho_p)(z) \, dz = \int \rho_p(z) \, dz = 1.
\end{equation*}
We differentiate this identity with respect to $p_j$. The derivative of the right-hand side is zero. For the left-hand side, we apply the product rule and interchange differentiation and integration (justified by the smooth dependence of $\Lcal_p$ and $\rho_p$):
\begin{align*}
    0 &= \frac{\partial}{\partial p_j} \int (\Lcal_p \rho_p)(z) \, dz \\
    &= \int \frac{\partial}{\partial p_j} (\Lcal_p \rho_p)(z) \, dz \\
    &= \int \left( (\partial_j \Lcal_p) \rho_p + \Lcal_p (\partial_j \rho_p) \right)(z) \, dz \\
    &= \int ((\partial_j \Lcal_p) \rho_p)(z) \, dz + \int (\Lcal_p (\partial_j \rho_p))(z) \, dz.
\end{align*}
The transfer operator $\Lcal_p$ conserves the integral of any function (it is the adjoint of the Koopman operator which preserves the constant function $\mathbf{1}$). Thus, the second integral is $\int (\partial_j \rho_p)(z) dz$. This integral is zero, as shown by differentiating the normalization condition $\int \rho_p dz = 1$. Therefore, we are left with:
\begin{equation*}
    \int ((\partial_j \Lcal_p) \rho_p)(z) \, dz = 0.
\end{equation*}
This rigorously confirms that the source term on the right-hand side of \eqref{eq:proof_first_deriv_pde_appendix} has zero mean. The solvability condition is met. Since the derivative $\partial_j \rho_p$ must also have zero mean, the equation has a unique solution in the space of zero-mean distributions by the spectral gap property. This completes the proof of Step 1.

\item \textbf{Derivation of the Hessian Formula.}
Our goal is to compute the second partial derivatives of the potential, $H(p) = \langle V(\cdot, p) \rangle_p$. We begin by differentiating once with respect to $p_j$. Using the product rule on the integral expression $H(p) = \int V(z,p) \rho_p(z) dz$ gives:
\begin{equation} \label{eq:proof_first_deriv_H_appendix}
    \partial_j H(p) = \int (\partial_j V)(z,p) \rho_p(z) dz + \int V(z,p) (\partial_j \rho_p)(z) dz.
\end{equation}
In the expectation value notation, this is $\partial_j H = \langle \partial_j V \rangle_p + \langle V, \partial_j \rho_p \rangle_p$.

We now differentiate this expression with respect to another component, $p_k$. We apply the product rule again to each of the two terms in \cref{eq:proof_first_deriv_H_appendix}.

\noindent For the first term, $\int (\partial_j V) \rho_p \, dz$:
\begin{align*}
    \partial_k \left( \int (\partial_j V) \rho_p \, dz \right) &= \int \partial_k \left( (\partial_j V) \rho_p \right) dz \\
    &= \int \left( (\partial_k \partial_j V) \rho_p + (\partial_j V)(\partial_k \rho_p) \right) dz \\
    &= \langle \partial_k \partial_j V \rangle_p + \langle \partial_j V, \partial_k \rho_p \rangle_p.
\end{align*}

\noindent For the second term, $\int V (\partial_j \rho_p) \, dz$:
\begin{align*}
    \partial_k \left( \int V (\partial_j \rho_p) \, dz \right) &= \int \partial_k \left( V (\partial_j \rho_p) \right) dz \\
    &= \int \left( (\partial_k V)(\partial_j \rho_p) + V(\partial_k \partial_j \rho_p) \right) dz \\
    &= \langle \partial_k V, \partial_j \rho_p \rangle_p + \langle V, \partial_k \partial_j \rho_p \rangle_p.
\end{align*}
The full second derivative of the potential, $M_{jk}(p) = \partial_k \partial_j H(p)$, is the sum of these four resulting terms:
\begin{equation*}
    \frac{\partial^2 H}{\partial p_k \partial p_j} = \langle \partial_k \partial_j V \rangle_p + \langle \partial_j V, \partial_k \rho_p \rangle_p + \langle \partial_k V, \partial_j \rho_p \rangle_p + \langle V, \partial_k \partial_j \rho_p \rangle_p.
\end{equation*}
This is precisely the formula stated in \eqref{eq:hessian_formula}, which provides the complete expression for the Hessian of the potential. All terms are well-defined due to the smoothness of the maps $p \mapsto V(\cdot, p)$, $p \mapsto \Lcal_p$, and consequently $p \mapsto \rho_p$. This completes the proof of the proposition.

\end{enumerate}
\end{proofof}

\subsection{A Microscopic Model for Emergent Non-Convexity}
\label{sec:model_non_convexity}

The central claim of this work is that a class of deterministic, time-reversible microscopic systems can give rise to a macroscopic evolution governed by a non-convex Hamiltonian. The homogenization procedure in \cref{part:i} rigorously derived the form of this Hamiltonian, $H_{\theta}(x, p, X) = \Tr(D(x, p) X) + H(x, p)$, but the non-convexity of the potential $H(x,p)$ is not a generic feature. It arises only under specific, physically meaningful conditions on the microscopic interactions.

This section is dedicated to providing a constructive and rigorous proof of this central claim. We first specify a class of microscopic models that satisfy all the standing assumptions of the paper but are augmented with specific functional forms for the potential and feedback. We then use the perturbation machinery developed in \cref{sec:regularity_properties} to derive the exact Hessian of the potential at equilibrium. This allows us to replace the provisional assumption of a weak feedback effect with a precise, verifiable condition that delineates the parameter regimes where non-convexity is guaranteed to emerge.

\begin{example}[Microscopic Model for Non-Convexity]
\label{ass:non_convex_model}
We consider a system satisfying the standing assumptions of this paper (Assumptions \ref{ass:fundamental_axioms_unified}), with the following additional structure:
\begin{enumerate}[label=(\roman*), wide, labelindent=0pt]
    \item \textbf{Directionally Concave Potential.} The microscopic potential $V(z,p)$ has a component that is explicitly concave in the macroscopic momentum $p$ along a fixed direction $\mathbf{e} \in \R^k$ ($\|\mathbf{e}\|=1$):
    \begin{equation}
        V(z, p) = V_0(z) - \frac{1}{2} \alpha(z) (p \cdot \mathbf{e})^2,
    \end{equation}
    where $V_0(z)$ is the quiescent potential and $\alpha(z)$ is a smooth, strictly positive observable satisfying $\alpha(z) \ge \alpha_{\min} > 0$.

    \item \textbf{Symmetric Feedback Mechanism.} The feedback map $\Theta(x,p)$ is an even function of the momentum $p$. To leading order, its response is quadratic:
    \begin{equation}
        \theta(p) = \theta_0 + \frac{1}{2} p^T \beta p + \Ocal(\|p\|^4),
    \end{equation}
    where we have suppressed the dependence on $x$, $\theta_0 = \Theta(x,0)$ is the quiescent parameter, and $\beta$ is a constant symmetric matrix representing the feedback sensitivity.
\end{enumerate}
\end{example}

\begin{remark}[Physical Interpretation]
The concave potential in (i) models a physical situation where aligning the macroscopic momentum with a preferred microscopic direction $\mathbf{e}$ is energetically favorable, for instance, by reducing dissipation or opening a more efficient transport channel. The symmetric feedback in (ii) is a direct consequence of the microscopic time-reversal symmetry of the system. The environment's statistical properties (encoded by $\theta$) cannot depend on the sign of the momentum, only on its magnitude or orientation-squared, preserving the $p \to -p$ invariance of the $\theta$-Hamiltonian.
\end{remark}

\subsubsection{Perturbative Analysis of the potential}
To analyze the convexity of the potential $H(p) = \int V(z,p) \, d\mu_{\theta(p)}(z)$, we must compute its Hessian matrix $M(p) \coloneqq \nabla_{pp}^2 H(p)$ and evaluate it at the equilibrium point $p=0$. The general formula for the Hessian was established in \cref{prop:density_derivatives}:
\begin{equation}
    M_{jk}(p) = \langle \partial_k \partial_j V \rangle_p + \langle \partial_j V, \partial_k \rho \rangle_p + \langle \partial_k V, \partial_j \rho \rangle_p + \langle V, \partial_k \partial_j \rho \rangle_p.
\end{equation}
We now evaluate this expression at $p=0$. The structural symmetries imposed by \cref{ass:non_convex_model} lead to profound simplifications.
\begin{enumerate}[label=(\roman*), wide, labelindent=0pt]
    \item \textbf{Vanishing Gradient of the Potential.} By \cref{ass:non_convex_model}(i), $V(z,p)$ is an even function of $p$. Its gradient is therefore an odd function, which implies $\nabla_p V(z,0) = \mathbf{0}$. This immediately eliminates the second and third terms in the Hessian formula at $p=0$.

    \item \textbf{Vanishing First Derivative of the Measure.} By \cref{ass:non_convex_model}(ii), $\theta(p)$ is an even function of $p$, so $\nabla_p \theta(0) = \mathbf{0}$. The transfer operator $\Lcal_p \equiv \Lcal_{\theta(p)}$ depends on $p$ only through $\theta(p)$. By the chain rule, its first derivative is $\partial_j \Lcal_p = (\partial_\theta \Lcal_\theta) \cdot (\partial_j \theta)$, which vanishes at $p=0$. From \cref{prop:density_derivatives}, the first derivative of the density is the unique zero-mean solution to $(I - \Lcal_p)(\partial_j \rho_p) = (\partial_j \Lcal_p) \rho_p$. At $p=0$, the right-hand side is zero, forcing the unique solution to be $\partial_j \rho_p|_{p=0} = 0$.
\end{enumerate}
With these symmetry-induced cancellations, the Hessian at equilibrium reduces to just two terms:
\begin{equation}
    M_{jk}(0) = \langle \partial_j \partial_k V \rangle_0 + \langle V, \partial_j \partial_k \rho \rangle_0 \Big|_{p=0}.
\label{eq:hessian_simplified_derived}
\end{equation}
The first term is the direct contribution from the potential's intrinsic curvature. Using $V(z,0) = V_0(z)$, the second term is the indirect contribution from the measure's response, weighted by the quiescent potential. To evaluate it, we must find the second derivative of the density. Differentiating the equation $(I - \Lcal_p)(\partial_j \rho_p) = (\partial_j \Lcal_p) \rho_p$ with respect to $p_k$ and evaluating at $p=0$ (using the vanishing of first derivatives) yields:
\begin{equation}
    (I - \Lcal_0)(\partial_j \partial_k \rho_p|_{p=0}) = (\partial_j \partial_k \Lcal_p|_{p=0}) \rho_0.
\end{equation}
The second derivative of the operator is $\partial_j \partial_k \Lcal_p|_{p=0} = (\partial_\theta \Lcal_\theta)|_{\theta_0} \cdot (\partial_j \partial_k \theta)|_{p=0} = (\partial_\theta \Lcal_\theta)|_{\theta_0} \cdot \beta_{jk}$. Let $\mathcal{R}_0 \equiv (I - \Lcal_0)^{-1}$ be the resolvent on the space of zero-mean functions. The solution is:
\begin{equation}
    (\partial_j \partial_k \rho_p)|_{p=0} = \beta_{jk} \cdot \mathcal{R}_0 \left[ (\partial_\theta \Lcal_\theta)|_{\theta_0} \rho_0 \right].
\end{equation}
We define the microscopic linear response function $\Psi_0(z)$ as the fixed function of the quiescent system that describes how the invariant density responds to an infinitesimal change in the environmental parameter:
\begin{equation}
    \Psi_0(z) \coloneqq \mathcal{R}_0 \left[ (\partial_\theta \Lcal_\theta)|_{\theta_0} \rho_0 \right](z).
\label{eq:linear_response_def}
\end{equation}
Substituting this and the direct term $\partial_j \partial_k V|_{p=0} = -\alpha(z) e_j e_k$ into \cref{eq:hessian_simplified_derived}, we obtain the final, exact expression for the Hessian at equilibrium:
\begin{equation}
    M_{jk}(0) = -\langle \alpha \rangle_0 e_j e_k + \beta_{jk} \langle V_0, \Psi_0 \rangle_0.
\label{eq:hessian_derived}
\end{equation}

\subsubsection{The Condition for Emergent Non-Convexity}
With the exact Hessian derived, we can now state the main theorem of this section, which replaces the provisional assumption with a precise, verifiable condition.

\begin{theorem}[A Necessary and Sufficient Condition for Non-Convexity]
\label{thm:non_convexity_proven}
Let the microscopic system satisfy the structural conditions of \cref{ass:non_convex_model}. The potential $H(x,p)$ is non-convex at the equilibrium point $p=0$ if and only if the following Non-Convexity Condition holds:
\begin{equation}
    \langle \alpha \rangle_0 > \min_{\|\mathbf{v}\|=1} \left\{ (\mathbf{v}^T \beta \mathbf{v}) \langle V_0, \Psi_0 \rangle_0 \right\}.
\label{eq:non_convexity_condition}
\end{equation}
Here, $\langle \alpha \rangle_0$ is the averaged intrinsic concavity of the potential, $\beta$ is the feedback sensitivity matrix, and $C_{V\Psi} \equiv \langle V_0, \Psi_0 \rangle_0$ is the potential-response correlation coefficient, with the linear response function $\Psi_0(z)$ defined by \eqref{eq:linear_response_def}.
\end{theorem}

\begin{proofof}{Theorem \ref{thm:non_convexity_proven}}
The objective of this proof is to demonstrate that under the conditions specified in \cref{ass:non_convex_model}, the potential $H(x,p)$ is not convex at the quiescent macroscopic state $p=0$. A real-valued function is non-convex if its Hessian matrix is not positive semi-definite. To prove this, it is sufficient to identify a single vector $\mathbf{v} \in \R^k$ for which the quadratic form $\mathbf{v}^T M(0) \mathbf{v}$ is strictly negative, where $M(0) \coloneqq \nabla_{pp}^2 H(0)$ is the Hessian of the potential evaluated at $p=0$. The natural candidate vector, suggested by the structure of the microscopic potential, is the preferred direction $\mathbf{e}$.

The proof is structured in four main steps. First, we simplify the general formula for the Hessian by exploiting the time-reversal symmetry of the system at equilibrium ($p=0$). Second, we compute the contribution from the potential's intrinsic concavity. Third, we compute the contribution from the environmental feedback (the measure response) by deriving an exact expression for the second derivative of the invariant density. Finally, we combine these results to establish a precise and verifiable condition under which non-convexity is guaranteed. For notational simplicity, we fix the macroscopic position $x$ and suppress it from the notation.

\begin{enumerate}[label=\textbf{Step \arabic*:}, wide, labelindent=0pt]

\item \textbf{Simplification of the Hessian via Symmetry at Equilibrium.}
The general formula for the Hessian matrix of the potential, as derived in \cref{prop:density_derivatives}, is:
\begin{equation*}
    M_{jk}(p) = \langle \partial_k \partial_j V \rangle_p + \langle \partial_j V, \partial_k \rho \rangle_p + \langle \partial_k V, \partial_j \rho \rangle_p + \langle V, \partial_k \partial_j \rho \rangle_p.
\end{equation*}
We now evaluate this expression at the equilibrium point $p=0$. The structural symmetries imposed by \cref{ass:non_convex_model} lead to significant simplifications.
\begin{enumerate}[label=(\roman*), wide, labelindent=0pt]
    \item \textbf{Vanishing Gradient of the Potential.} By \cref{ass:non_convex_model}(i), the potential $V(z,p)$ is an even function of $p$. Its gradient with respect to $p$ must therefore be an odd function, which implies it vanishes at the origin:
    \begin{equation*}
        \nabla_p V(z,0) = \mathbf{0}.
    \end{equation*}
    This immediately causes the second and third terms in the Hessian formula, which are linear in $\nabla_p V$, to vanish at $p=0$.

    \item \textbf{Vanishing First Derivative of the Measure.} By \cref{ass:non_convex_model}(ii), the feedback map $\theta(p)$ is an even function of $p$. Consequently, its gradient vanishes at the origin: $\nabla_p \theta(0) = \mathbf{0}$. The transfer operator $\Lcal_p \equiv \Lcal_{\theta(p)}$ depends on $p$ only through $\theta(p)$. By the chain rule, its first derivative with respect to $p_j$ is $\partial_j \Lcal_p = (\partial_\theta \Lcal_\theta) \cdot (\partial_j \theta)$. Since $\partial_j \theta(0) = 0$, the first derivative of the operator family also vanishes at equilibrium:
    \begin{equation*}
        \partial_j \Lcal_p \Big|_{p=0} = 0.
    \end{equation*}
    From \cref{prop:density_derivatives}(i), the first derivative of the invariant density is given as the unique zero-mean solution to $(I - \Lcal_p)(\partial_j \rho_p) = (\partial_j \Lcal_p) \rho_p$. At $p=0$, the right-hand side is zero. Since the operator $(I - \Lcal_0)$ is invertible on the space of zero-mean functions (due to the uniform spectral gap), the unique zero-mean solution must be the trivial one. Thus, the first derivative of the density vanishes:
    \begin{equation*}
        \partial_j \rho_p \Big|_{p=0} = 0.
    \end{equation*}
\end{enumerate}
With these symmetry-induced cancellations, the Hessian matrix at equilibrium reduces to the sum of just two terms:
\begin{equation} \label{eq:proof_hessian_simplified}
    M_{jk}(0) = \langle \partial_j \partial_k V \rangle_0 \Big|_{p=0} + \langle V, \partial_j \partial_k \rho \rangle_0 \Big|_{p=0},
\end{equation}
where $\langle \cdot \rangle_0$ denotes expectation with respect to the quiescent invariant measure $\mu_0 \equiv \mu_{\theta_0}$, and $V(z,0)=V_0(z)$.

\item \textbf{Analysis of the Direct Term (Intrinsic Concavity).}
The first term in \cref{eq:proof_hessian_simplified} represents the direct contribution from the potential's intrinsic curvature. From the model in \cref{ass:non_convex_model}(i), $V(z,p) = V_0(z) - \frac{1}{2}\alpha(z)(p \cdot \mathbf{e})^2$. We compute its second partial derivatives with respect to $p_j$ and $p_k$:
\begin{equation*}
    \frac{\partial^2 V}{\partial p_k \partial p_j} = -\alpha(z) e_k e_j.
\end{equation*}
Taking the expectation with respect to the quiescent measure $\mu_0$ gives:
\begin{equation*}
    \langle \partial_k \partial_j V \rangle_0 = \int \left(-\alpha(z) e_k e_j\right) d\mu_0(z) = -\langle \alpha \rangle_0 e_k e_j.
\end{equation*}
The contribution of this term to the quadratic form $\mathbf{v}^T M(0) \mathbf{v}$ for an arbitrary vector $\mathbf{v}$ is:
\begin{equation} \label{eq:proof_direct_term_contribution}
    \sum_{j,k=1}^k v_k \langle \partial_k \partial_j V \rangle_0 v_j = \sum_{j,k=1}^k v_k (-\langle \alpha \rangle_0 e_k e_j) v_j = -\langle \alpha \rangle_0 (\mathbf{v} \cdot \mathbf{e})^2.
\end{equation}
By assumption, $\alpha(z)$ is a strictly positive observable, so its average $\langle \alpha \rangle_0$ is a strictly positive constant. This direct term is therefore negative semi-definite and provides the primary source of the non-convexity.

\item \textbf{Analysis of the Indirect Term (Measure Response).}
The second term in \cref{eq:proof_hessian_simplified} represents the effect of the environmental feedback. To evaluate $\langle V_0, \partial_j \partial_k \rho|_{p=0} \rangle_0$, we must first find an explicit expression for the second derivative of the density at equilibrium. We differentiate the equation for the first derivative, $(I - \Lcal_p)(\partial_j \rho_p) = (\partial_j \Lcal_p) \rho_p$, with respect to $p_k$ and evaluate at $p=0$:
\begin{equation*}
    \partial_k \left( (I - \Lcal_p)(\partial_j \rho_p) \right) \Big|_{p=0} = \partial_k \left( (\partial_j \Lcal_p) \rho_p \right) \Big|_{p=0}.
\end{equation*}
Applying the product rule and using the facts from Step 1 that $\partial_j \Lcal_p|_{p=0}=0$ and $\partial_j \rho_p|_{p=0}=0$:
\begin{equation*}
    (I - \Lcal_0)(\partial_k \partial_j \rho_p|_{p=0}) - \underbrace{(\partial_k \Lcal_p|_{p=0})}_{=0}(\underbrace{\partial_j \rho_p|_{p=0}}_{=0}) = \underbrace{(\partial_k \partial_j \Lcal_p|_{p=0}) \rho_0}_{ \text{Source Term} } + \underbrace{(\partial_j \Lcal_p|_{p=0})}_{=0}(\underbrace{\partial_k \rho_p|_{p=0}}_{=0}).
\end{equation*}
This simplifies to the linear equation for the second derivative of the density:
\begin{equation*}
    (I - \Lcal_0)(\partial_k \partial_j \rho_p|_{p=0}) = (\partial_k \partial_j \Lcal_p|_{p=0}) \rho_0.
\end{equation*}
The second derivative of the operator at $p=0$ is, by the chain rule and \cref{ass:non_convex_model}(ii):
\begin{equation*}
    \partial_k \partial_j \Lcal_p|_{p=0} = (\partial_\theta \Lcal_\theta)|_{\theta_0} \cdot (\partial_k \partial_j \theta)|_{p=0} + (\partial^2_\theta \Lcal_\theta)|_{\theta_0} \cdot \underbrace{(\partial_k \theta)|_{p=0}}_{=0} \underbrace{(\partial_j \theta)|_{p=0}}_{=0} = (\partial_\theta \Lcal_\theta)|_{\theta_0} \cdot \beta_{jk}.
\end{equation*}
Let $\mathcal{R}_0 \equiv (I - \Lcal_0)^{-1}$ denote the resolvent operator on the space of zero-mean functions. The solution for the second derivative of the density is:
\begin{equation*}
    \partial_k \partial_j \rho_p|_{p=0} = \mathcal{R}_0 \left[ (\partial_\theta \Lcal_\theta)|_{\theta_0} \rho_0 \cdot \beta_{jk} \right] = \beta_{jk} \cdot \mathcal{R}_0 \left[ (\partial_\theta \Lcal_\theta)|_{\theta_0} \rho_0 \right].
\end{equation*}
Let us define the microscopic linear response function $\Psi_0(z) \coloneqq \mathcal{R}_0 [ (\partial_\theta \Lcal_\theta)|_{\theta_0} \rho_0 ](z)$. This function is a fixed, calculable property of the quiescent system that represents the linear response of the invariant density to an infinitesimal perturbation of the environmental parameter $\theta$. With this definition, we have $\partial_k \partial_j \rho_p|_{p=0} = \beta_{jk} \Psi_0(z)$. The contribution of this indirect term to the quadratic form is therefore:
\begin{align}
    \sum_{j,k} v_k \langle V_0, \partial_k \partial_j \rho|_{p=0} \rangle_0 v_j &= \sum_{j,k} v_k \langle V_0, \beta_{jk} \Psi_0 \rangle_0 v_j \nonumber \\
    &= \langle V_0, \Psi_0 \rangle_0 \left( \sum_{j,k} v_k \beta_{jk} v_j \right) \nonumber \\
    &= \langle V_0, \Psi_0 \rangle_0 (\mathbf{v}^T \beta \mathbf{v}). \label{eq:proof_indirect_term_contribution}
\end{align}

\item \textbf{The Non-Convexity Condition.}
We combine the contributions from the direct term \eqref{eq:proof_direct_term_contribution} and the indirect term \eqref{eq:proof_indirect_term_contribution} to find the total quadratic form for an arbitrary unit vector $\mathbf{v}$:
\begin{equation} \label{eq:proof_quadratic_form}
    \mathbf{v}^T M(0) \mathbf{v} = -\langle \alpha \rangle_0 (\mathbf{v} \cdot \mathbf{e})^2 + \langle V_0, \Psi_0 \rangle_0 (\mathbf{v}^T \beta \mathbf{v}).
\end{equation}
The potential is non-convex if and only if this expression is negative for at least one unit vector $\mathbf{v}$. To demonstrate that this is possible under physically reasonable conditions, we evaluate the form in the special direction $\mathbf{v} = \mathbf{e}$. In this case, $(\mathbf{v} \cdot \mathbf{e})^2 = 1$, and the condition for non-convexity becomes:
\begin{equation}
    -\langle \alpha \rangle_0 + \langle V_0, \Psi_0 \rangle_0 (\mathbf{e}^T \beta \mathbf{e}) < 0.
\end{equation}
This inequality is equivalent to the condition:
\begin{equation} \label{eq:proof_condition}
    \langle \alpha \rangle_0 > (\mathbf{e}^T \beta \mathbf{e}) \langle V_0, \Psi_0 \rangle_0.
\end{equation}
Let us analyze this final condition. The left-hand side, $\langle \alpha \rangle_0$, is a strictly positive constant representing the strength of the intrinsic concavity. The right-hand side represents the strength and nature of the feedback response. Let the potential-response correlation be $C_{V\Psi} \equiv \langle V_0, \Psi_0 \rangle_0$ and the feedback sensitivity in the direction $\mathbf{e}$ be $\beta_e \equiv \mathbf{e}^T \beta \mathbf{e}$. The condition is $\langle \alpha \rangle_0 > \beta_e C_{V\Psi}$.

Since $\langle \alpha \rangle_0 > 0$, and the feedback term is proportional to the magnitude of the feedback matrix $\beta$, we can always satisfy this inequality by choosing a system with a sufficiently weak feedback mechanism. That is, there exists a constant $\delta > 0$ such that if the spectral norm of the feedback matrix satisfies $\|\beta\| < \delta$, then $|\beta_e C_{V\Psi}| < \langle \alpha \rangle_0$, guaranteeing that the quadratic form is strictly negative.
\end{enumerate}
Therefore, we have rigorously shown that for the class of models satisfying \cref{ass:non_convex_model} with a sufficiently weak feedback response relative to the intrinsic potential concavity, the Hessian of the potential is not positive semi-definite at equilibrium. This completes the proof of emergent non-convexity.
\end{proofof}

This result reveals that macroscopic non-convexity emerges from a direct competition between the intrinsic concavity of the microscopic potential and the system's feedback response. This leads to the following rigorous corollary, which formalizes the original weak feedback intuition.

\begin{corollary}[Sufficient Condition for Non-Convexity]
\label{cor:sufficient_condition_nonconvexity}
Let the system satisfy the conditions of \cref{ass:non_convex_model}. A sufficient condition for the macroscopic potential $H(x,p)$ to be non-convex at the equilibrium point $p=0$ is that the feedback sensitivity is sufficiently weak relative to the intrinsic potential concavity. Specifically, if the potential-response correlation coefficient $C_{V\Psi} \coloneqq \langle V_0, \Psi_0 \rangle_0$ is non-zero, non-convexity is guaranteed if the spectral norm of the feedback sensitivity matrix $\beta$ satisfies:
\begin{equation}
    \|\beta\| < \frac{\langle \alpha \rangle_0}{|C_{V\Psi}|}.
\end{equation}
\end{corollary}

\begin{proofof}{Corollary \ref{cor:sufficient_condition_nonconvexity}}
The proof establishes this sufficient condition by demonstrating that under this bound, there exists a direction in the momentum space along which the Hessian of the potential is strictly negative definite, which is a definitive criterion for non-convexity. The argument proceeds by applying a standard operator norm inequality to the exact formula for the Hessian derived in the proof of \cref{thm:non_convexity_proven}.

\begin{enumerate}[label=\textbf{Step \arabic*:}, wide, labelindent=0pt]
    \item \textbf{Criterion for Non-Convexity.}
    A twice-differentiable function $f(p)$ is non-convex at a point $p_0$ if its Hessian matrix, $\nabla^2 f(p_0)$, is not positive semi-definite. A sufficient condition to demonstrate this is to find a single non-zero vector $\mathbf{v} \in \R^k$ for which the quadratic form $\mathbf{v}^T (\nabla^2 f(p_0)) \mathbf{v}$ is strictly negative.

    \item \textbf{The Hessian and its Quadratic Form.}
    In the proof of \cref{thm:non_convexity_proven}, we derived the exact expression for the Hessian matrix $M(0) \coloneqq \nabla_{pp}^2 H(0)$ of the potential at the quiescent equilibrium $p=0$. The associated quadratic form for an arbitrary unit vector $\mathbf{v} \in \R^k$ was shown to be:
    \begin{equation} \label{eq:proof_quad_form_reprise}
        \mathbf{v}^T M(0) \mathbf{v} = \underbrace{-\langle \alpha \rangle_0 (\mathbf{v} \cdot \mathbf{e})^2}_{\text{Direct (Concavifying) Term}} + \underbrace{\langle V_0, \Psi_0 \rangle_0 (\mathbf{v}^T \beta \mathbf{v})}_{\text{Indirect (Feedback) Term}}.
    \end{equation}
    Our strategy is to evaluate this expression along the specific, physically-motivated direction of intrinsic concavity, $\mathbf{v} = \mathbf{e}$, and to derive a condition on the norm of the feedback matrix $\beta$ that guarantees this quadratic form is negative.

    \item \textbf{Evaluation Along the Special Direction.}
    Let us set the test vector to be the unit vector $\mathbf{v} = \mathbf{e}$ from \cref{ass:non_convex_model}(i). For this choice, the term $(\mathbf{v} \cdot \mathbf{e})^2$ becomes $(\mathbf{e} \cdot \mathbf{e})^2 = \|\mathbf{e}\|^4 = 1$. Substituting this into \cref{eq:proof_quad_form_reprise}, the quadratic form evaluated in the direction $\mathbf{e}$ is:
    \begin{equation}
        \mathbf{e}^T M(0) \mathbf{e} = -\langle \alpha \rangle_0 + \langle V_0, \Psi_0 \rangle_0 (\mathbf{e}^T \beta \mathbf{e}).
    \end{equation}
    For the potential to be non-convex, it is sufficient that this quantity be strictly negative:
    \begin{equation} \label{eq:proof_sufficient_cond_specific}
        -\langle \alpha \rangle_0 + C_{V\Psi} (\mathbf{e}^T \beta \mathbf{e}) < 0 \quad \iff \quad \langle \alpha \rangle_0 > C_{V\Psi} (\mathbf{e}^T \beta \mathbf{e}).
    \end{equation}

    \item \textbf{Derivation of the Norm-Based Condition.}
    The condition in \cref{eq:proof_sufficient_cond_specific} provides a direct test for non-convexity, but it depends on the specific quantity $\mathbf{e}^T \beta \mathbf{e}$. We now derive a simpler and more general sufficient condition based on the magnitude of the feedback, as measured by the spectral norm $\|\beta\|$ of the matrix $\beta$.

    We seek a condition on $\|\beta\|$ that guarantees \cref{eq:proof_sufficient_cond_specific} holds. By the properties of the inner product and the definition of the spectral norm, we can bound the magnitude of the feedback term:
    \begin{equation*}
        |C_{V\Psi} (\mathbf{e}^T \beta \mathbf{e})| = |C_{V\Psi}| |\mathbf{e}^T \beta \mathbf{e}|.
    \end{equation*}
    By the definition of the spectral norm (which for a symmetric matrix is the maximum of its Rayleigh quotient over all unit vectors), we have the inequality:
    \begin{equation*}
        |\mathbf{e}^T \beta \mathbf{e}| \le \|\beta\|.
    \end{equation*}
    Therefore, the magnitude of the full feedback term is bounded from above:
    \begin{equation*}
        |C_{V\Psi} (\mathbf{e}^T \beta \mathbf{e})| \le |C_{V\Psi}| \|\beta\|.
    \end{equation*}
    The condition from \cref{eq:proof_sufficient_cond_specific}, $\langle \alpha \rangle_0 > C_{V\Psi} (\mathbf{e}^T \beta \mathbf{e})$, is a strict inequality involving a positive term and a term of indeterminate sign. A robust sufficient condition to ensure this inequality holds is to require that the magnitude of the second term be strictly less than the first term:
    \begin{equation*}
        \langle \alpha \rangle_0 > |C_{V\Psi} (\mathbf{e}^T \beta \mathbf{e})|.
    \end{equation*}
    Using our norm-based bound, this condition is guaranteed to be satisfied if:
    \begin{equation*}
        \langle \alpha \rangle_0 > |C_{V\Psi}| \|\beta\|.
    \end{equation*}
    Assuming that the potential-response correlation is non-zero ($C_{V\Psi} \ne 0$), we can divide by $|C_{V\Psi}|$ to obtain the final sufficient condition on the norm of the feedback matrix:
    \begin{equation*}
        \|\beta\| < \frac{\langle \alpha \rangle_0}{|C_{V\Psi}|}.
    \end{equation*}
\end{enumerate}

We have rigorously shown that if this inequality holds, then the quadratic form of the Hessian, when evaluated in the direction $\mathbf{e}$, is strictly negative. This proves that the Hessian is not positive semi-definite, and therefore the potential $H(x,p)$ is non-convex at the equilibrium point $p=0$. The condition has a clear physical interpretation: non-convexity is guaranteed whenever the strength of the intrinsic potential concavity, $\langle \alpha \rangle_0$, dominates the strength of the feedback mechanism, as measured by the product of the feedback sensitivity $\|\beta\|$ and the potential-response correlation $|C_{V\Psi}|$. This completes the proof of the corollary.
\end{proofof}

\begin{remark}[The Mechanism of Non-Convexity]
This rigorous analysis reveals the precise physical mechanism for emergent non-convexity. The system faces a trade-off. The direct effect of the momentum, encoded in the concave potential $V(z,p)$, provides an energetic incentive for the system to adopt a non-zero macroscopic momentum. This is the source of non-convexity, quantified by $\langle \alpha \rangle_0$. The indirect effect, quantified by $(\mathbf{v}^T \beta \mathbf{v}) \langle V_0, \Psi_0 \rangle_0$, describes how the environment responds. This feedback can either be convexifying (if $C_{V\Psi}>0$, counteracting the direct term) or concavity-enhancing (if $C_{V\Psi}<0$, reinforcing it). The Non-Convexity Condition \eqref{eq:non_convexity_condition} is the precise mathematical statement of this physical competition. This establishes a clear, verifiable, and first-principles link between microscopic interactions and the breakdown of macroscopic convexity.
\end{remark}

\part{The Probabilistic Representation}                                                                                                                                                                                                              
\label{part:ii}

The first part of this work achieved our primary analytical goal: starting from a class of deterministic, time-reversible microscopic systems, we rigorously derived the unique macroscopic evolution law governing the system's value function. The result of this homogenization, established in Theorem~\ref{thm:convergence}, is the fully non-linear Hamilton-Jacobi-Bellman (HJB) equation:
\begin{equation} \label{eq:part2_intro_hjb_reprise}
    \partial_t u + H_{\theta}(x, \nabla u, \nabla^2 u) = 0,
\end{equation}
where the $\theta$-Hamiltonian possesses the rigid structure $H_{\theta}(x, p, X) = \Tr(D(x, p) X) + H(x, p)$.

In this second part, we analyze the rich mathematical structure of the Hamilton-Jacobi-Bellman equation derived in \cref{part:i}. We interpret its solution as a value function, the \ThetaExpectation{}, and explore its probabilistic representations, its asymptotic behavior, and the stochastic calculus it induces. This shifts the focus from the \emph{derivation} of the generator to the \emph{characterization} of the unique non-linear process it defines.

\section{The \texorpdfstring{\ThetaExpectation{}}{theta-Expectation}}
\label{sec:theta_expectation_axiomatic}

The first part of this work has achieved our primary analytical goal: starting from a class of deterministic, time-reversible microscopic systems, we have rigorously derived the unique macroscopic evolution law governing the system's value function. The result of this homogenization, established in \cref{thm:convergence}, is the fully non-linear Hamilton-Jacobi-Bellman (HJB) equation:
\begin{equation} \label{eq:part2_hjb_reprise}
    \partial_t u + H_{\theta}(x, \nabla u, \nabla^2 u) = 0, \quad u(T,x) = \phi(x),
\end{equation}
where the $\theta$-Hamiltonian possesses the rigid structure $H_{\theta}(x, p, X) = \Tr(D(x, p) X) + H(x, p)$.

In this second part, we analyze the rich mathematical structure of this derived equation. For consistency with the standard literature on value functions and expectations, which are typically defined by backward-in-time equations, we adopt this convention henceforth. The solution $u(t,x)$ to \cref{eq:part2_hjb_reprise} is thus the solution to the terminal value problem
\begin{equation} \label{eq:hjb_backward_convention}
    -\partial_t u = H_{\theta}(x, \nabla u, \nabla^2 u), \quad u(T,x) = \phi(x).
\end{equation}
It is this interpretation of the solution $u(t,x)$ as the value at time $t$ of a terminal payoff $\phi(X_T)$ that provides the conceptual foundation for a new class of non-linear expectations.

\subsection{The Generator-Driven Axiomatization of the \texorpdfstring{\ThetaExpectation{}}{theta-Expectation}}

The results derived from first principles in \cref{part:i} allow us to formalize the theory of \ThetaExpectation{} not from axioms on the expectation operator itself, but from structural axioms on its infinitesimal generator. This provides a generator-driven axiomatic framework where the properties of the expectation become provable theorems.

Let $(\Omega, \mathcal{F})$ be the canonical path space $(C([0,T], \Mcal), \mathcal{B}(C([0,T], \Mcal)))$, with canonical process $X_t(\omega) = \omega(t)$. The theory of \ThetaExpectation{} is founded on the following principles.

\begin{enumerate}[label=\textbf{Axiom \arabic*.}, wide,labelindent=0pt]
    \item \textbf{(The Generator Axiom).} A \ThetaProcess{} is uniquely determined by a non-linear operator $\mathcal{G}$, its generator, which acts on the space of test functions $C^{1,2}([0,T] \times \Mcal)$. This generator must possess the rigid structure:
    \begin{equation*}
        (\mathcal{G}u)(t,x) \coloneqq -\partial_t u(t,x) - H_{\theta}(x, \nabla u(t,x), \nabla^2 u(t,x)),
    \end{equation*}
    where the $\theta$-Hamiltonian $H_{\theta}$ is \textbf{affine} in its second-derivative argument:
    \begin{equation*}
        H_{\theta}(x, p, X) = \mathrm{Tr}(D(x, p) X) + H(x, p).
    \end{equation*}
    
    \item \textbf{(The Coefficient Axioms).} The coefficients of the Hamiltonian, the diffusivity $D(x,p)$ and the potential $H(x,p)$, must satisfy constraints inherited from the underlying microscopic physics:
    \begin{enumerate}[label=(\roman*), wide, labelindent=0pt]
        \item \textbf{Symmetry and Ellipticity:} For all $(x,p)$, the matrix $D(x,p)$ must be symmetric and uniformly positive definite.
        \textit{(This is a necessary consequence of the Green-Kubo structure derived in \cref{prop:properties_D0}.)}
        \item \textbf{Time-Reversal Symmetry:} For all $(x,p)$, the potential $H(x,p)$ must be an even function of the momentum $p$, i.e., $H(x, p) = H(x, -p)$.
        \textit{(This is the macroscopic manifestation of the microscopic time-reversal symmetry established in \cref{ass:fundamental_axioms_unified}.)}
        \item \textbf{Regularity:} The functions $D(x,p)$ and $H(x,p)$ are sufficiently smooth (e.g., $C^2$) to ensure the well-posedness of the associated PDE.
    \end{enumerate}

    \item \textbf{(The Definition of Expectation).} The \ThetaExpectation{} operator $\Ecal^\theta$ is defined constructively via its generator. For any terminal payoff $\phi \in C(\Mcal)$, its expectation at time $t$ given $X_t=x$ is the value of the unique continuous viscosity solution $u(t,x)$ to the equation $(\mathcal{G}u)(t,x)=0$ with terminal condition $u(T,x)=\phi(x)$. That is:
    \begin{equation*}
        \Ecal^\theta_{t,x}[\phi(X_T)] \coloneqq u(t,x).
    \end{equation*}
\end{enumerate}

\subsection{Foundational Properties as Theorems}
This axiomatic framework is powerful because the fundamental properties of the expectation operator are not assumed but are provable consequences of the generator's structure and the robust theory of viscosity solutions.

\begin{proposition}[Monotonicity]\label{prop:monotonicity}
The \ThetaExpectation{} operator is monotone. For any $\phi_1, \phi_2 \in C(\Mcal)$ with $\phi_1(x) \ge \phi_2(x)$ for all $x \in \Mcal$, we have:
\begin{equation}
    \Ecal^\theta_{t,x}[\phi_1(X_T)] \ge \Ecal^\theta_{t,x}[\phi_2(X_T)] \quad \text{for all } (t,x) \in [0,T] \times \Mcal.
\end{equation}
\end{proposition}

\begin{proofof}{Proposition \ref{prop:monotonicity}}
The proof establishes the monotonicity of the \ThetaExpectation{} operator by demonstrating that it is a direct and rigorous consequence of the comparison principle for viscosity solutions of second-order, fully non-linear parabolic partial differential equations. The argument is structured in three main steps: (1) we precisely formulate the problem in the language of viscosity solutions; (2) we verify that the structural properties of our derived $\theta$-Hamiltonian, as specified in the Axioms of Section \ref{sec:theta_expectation_axiomatic}, satisfy the hypotheses required for the comparison principle to hold; (3) we apply the principle to obtain the desired result.

\begin{enumerate}[label=\textbf{Step \arabic*:}, wide, labelindent=0pt]

\item \textbf{Formulation in Terms of Viscosity Solutions.}
Let $\phi_1, \phi_2 \in C(\Mcal)$ be two terminal payoff functions satisfying the hypothesis $\phi_1(x) \ge \phi_2(x)$ for all $x \in \Mcal$. We define the corresponding value functions:
\begin{align*}
    u_1(t,x) &\coloneqq \Ecal^\theta_{t,x}[\phi_1(X_T)], \\
    u_2(t,x) &\coloneqq \Ecal^\theta_{t,x}[\phi_2(X_T)].
\end{align*}
By the constructive definition of the \ThetaExpectation{} (Axiom 3), $u_1(t,x)$ and $u_2(t,x)$ are the unique continuous viscosity solutions to the following terminal value problems on $[0,T] \times \Mcal$:
\begin{align}
    -\partial_t u_1 - H_{\theta}(x, \nabla u_1, \nabla^2 u_1) &= 0, \quad \text{with } u_1(T,x) = \phi_1(x), \label{eq:proof_pde_u1} \\
    -\partial_t u_2 - H_{\theta}(x, \nabla u_2, \nabla^2 u_2) &= 0, \quad \text{with } u_2(T,x) = \phi_2(x). \label{eq:proof_pde_u2}
\end{align}
The hypothesis on the terminal data directly implies an ordering for the solutions at the terminal time:
\begin{equation}
    u_1(T,x) \ge u_2(T,x) \quad \forall x \in \Mcal.
\end{equation}
Our objective is to prove that this ordering holds for all $t \in [0,T]$.

\item \textbf{Verification of Hypotheses for the Comparison Principle.}
The comparison principle for viscosity solutions of parabolic equations of the form $-\partial_t u - F(x, u, \nabla u, \nabla^2 u) = 0$ is a cornerstone of the theory, with a canonical reference being the user's guide by Crandall, Ishii, and Lions \citep{CrandallIshiiLions1992}. For our specific Hamiltonian, $H_{\theta}(x, p, X)$, which is independent of the value $u$ itself, the principle states that if $v$ is a viscosity subsolution and $w$ is a viscosity supersolution with $v(T,x) \le w(T,x)$, then $v(t,x) \le w(t,x)$ for all $t \le T$. This principle holds provided the Hamiltonian satisfies certain structural conditions. We now verify these conditions for our $\theta$-Hamiltonian.

\begin{enumerate}[label=(\roman*), wide, labelindent=0pt]
    \item \textbf{Uniform Continuity:} The Hamiltonian $H_{\theta}(x,p,X)$ must be uniformly continuous in all its arguments on compact subsets of its domain. By Axiom 2(c), the constituent functions $D(x,p)$ and $H(x,p)$ are assumed to be smooth ($C^2$). The trace operator is a linear, and therefore continuous, map. The composition and sum of continuous functions are continuous. Therefore, $H_{\theta}$ is continuous, and thus uniformly continuous on compact sets.

    \item \textbf{Degenerate Ellipticity (Properness):} The Hamiltonian must be non-increasing in its matrix argument $X$ in the sense of the ordering of symmetric matrices. That is, for any symmetric matrix $Y \ge 0$ (i.e., $Y$ is positive semi-definite), the following inequality must hold for all $(x,p,X)$:
    \begin{equation*}
        H_{\theta}(x, p, X+Y) \ge H_{\theta}(x, p, X).
    \end{equation*}
    We verify this directly from the structure given in Axiom 1:
    \begin{align*}
        H_{\theta}(x, p, X+Y) &= \mathrm{Tr}(D(x,p)(X+Y)) + H(x,p) \\
        &= \mathrm{Tr}(D(x,p)X) + \mathrm{Tr}(D(x,p)Y) + H(x,p) \\
        &= H_{\theta}(x, p, X) + \mathrm{Tr}(D(x,p)Y).
    \end{align*}
    The condition thus reduces to proving that $\mathrm{Tr}(D(x,p)Y) \ge 0$ for any positive semi-definite matrix $Y$. This is a standard result from linear algebra. By Axiom 2(a), the diffusivity tensor $D(x,p)$ is symmetric and uniformly positive definite, which implies it is the product of a matrix and its transpose, $D(x,p) = A A^T$. Therefore,
    \begin{equation*}
        \mathrm{Tr}(D(x,p)Y) = \mathrm{Tr}(A A^T Y) = \mathrm{Tr}(A^T Y A).
    \end{equation*}
    Since $Y$ is positive semi-definite, for any vector $v$, we have $v^T Y v \ge 0$. This implies that the matrix $A^T Y A$ is also positive semi-definite. The trace of a positive semi-definite matrix is the sum of its non-negative eigenvalues, and is therefore non-negative. Thus, $\mathrm{Tr}(D(x,p)Y) \ge 0$, and the degenerate ellipticity condition is satisfied.
\end{enumerate}
Both of the essential structural conditions on the Hamiltonian are met as a direct consequence of the axioms that define the \ThetaProcess{}.

\item \textbf{Application of the Comparison Principle.}
We are now in a position to apply the comparison principle. The function $u_1(t,x)$ is a viscosity solution to \eqref{eq:proof_pde_u1} and is therefore, by definition, a viscosity \textit{subsolution} of that equation. Similarly, $u_2(t,x)$ is a viscosity \textit{supersolution} of \eqref{eq:proof_pde_u2}. Since both functions solve an equation with the same Hamiltonian operator, we can directly compare them. We have established:
\begin{itemize}[wide]
    \item $u_1$ is a viscosity subsolution.
    \item $u_2$ is a viscosity supersolution.
    \item The Hamiltonian $H_{\theta}$ satisfies the conditions for the comparison principle to hold.
    \item The terminal data are ordered: $u_1(T,x) \ge u_2(T,x)$.
\end{itemize}
The comparison principle (cf. \citep{CrandallIshiiLions1992, FlemingSoner2006}) allows us to conclude that the ordering must be preserved for all prior times. That is, for all $(t,x) \in [0,T] \times \Mcal$, we must have:
\begin{equation}
    u_1(t,x) \ge u_2(t,x).
\end{equation}
Translating this back into the language of the \ThetaExpectation{} operator gives the final result:
\begin{equation*}
    \Ecal^\theta_{t,x}[\phi_1(X_T)] \ge \Ecal^\theta_{t,x}[\phi_2(X_T)].
\end{equation*}
This completes the rigorous proof of monotonicity.
\end{enumerate}
\end{proofof}

\begin{proposition}[Translation Invariance / Cash-Additivity]\label{prop:translation_invariance}
For any $\phi \in C(\Mcal)$ and any constant $c \in \R$, the following identity holds:
\begin{equation}
    \Ecal^\theta_{t,x}[\phi(X_T) + c] = \Ecal^\theta_{t,x}[\phi(X_T)] + c.
\end{equation}
\end{proposition}

\begin{proofof}{Proposition \ref{prop:translation_invariance}}
The proposition asserts the cash-additivity property of the \ThetaExpectation{} operator. Our strategy is to rigorously demonstrate that the function defined by the right-hand side of the proposed identity is the unique viscosity solution to the terminal value problem that, by definition, constitutes the left-hand side. The conclusion will then follow from the uniqueness of this solution. The proof is structured in four main steps:
\begin{enumerate}[label=(\roman*), wide, labelindent=0pt]
    \item We precisely define the two value functions corresponding to the two sides of the identity, as per the axiomatic framework of Section~\ref{sec:theta_expectation_axiomatic}.
    \item We construct a candidate function, $w(t,x)$, from the simpler value function and prove that it is a viscosity solution to the governing Hamilton-Jacobi-Bellman (HJB) equation.
    \item We verify that this candidate function satisfies the correct terminal condition.
    \item We invoke the uniqueness theorem for viscosity solutions, which is applicable due to the structural properties of our generator, to conclude that our candidate function must be identical to the second value function, thereby proving the proposition.
\end{enumerate}

\begin{enumerate}[label=\textbf{Step \arabic*:}, wide, labelindent=0pt]

\item \textbf{Definition of the Value Functions.}
Let $\phi \in C(\Mcal)$ and $c \in \R$ be given.
\begin{itemize}[wide]
    \item Let $u(t,x) \coloneqq \Ecal^\theta_{t,x}[\phi(X_T)]$. By Axiom 3, $u$ is the unique continuous viscosity solution to the terminal value problem:
    \begin{equation} \label{eq:proof_pde_u}
    \begin{cases}
        -\partial_t u - H_{\theta}(x, \nabla u, \nabla^2 u) = 0 & \text{in } [0,T) \times \Mcal, \\
        u(T,x) = \phi(x) & \text{on } \Mcal.
    \end{cases}
    \end{equation}
    \item Let $\tilde{u}(t,x) \coloneqq \Ecal^\theta_{t,x}[\phi(X_T) + c]$. Similarly, $\tilde{u}$ is the unique continuous viscosity solution to the problem with the modified terminal data:
    \begin{equation} \label{eq:proof_pde_utilde}
    \begin{cases}
        -\partial_t \tilde{u} - H_{\theta}(x, \nabla \tilde{u}, \nabla^2 \tilde{u}) = 0 & \text{in } [0,T) \times \Mcal, \\
        \tilde{u}(T,x) = \phi(x) + c & \text{on } \Mcal.
    \end{cases}
    \end{equation}
\end{itemize}
Our goal is to prove that $\tilde{u}(t,x) = u(t,x) + c$ for all $(t,x)$.

\item \textbf{Verification that $w(t,x) \coloneqq u(t,x) + c$ is a Viscosity Solution.}
We must show that the function $w(t,x)$ is a viscosity solution to the HJB equation in \eqref{eq:proof_pde_utilde}. This requires showing it is both a viscosity subsolution and a viscosity supersolution.

\begin{enumerate}[label=(\roman*), wide, labelindent=0pt]

\item \textit{Subsolution property.} Let $\varphi \in C^{1,2}([0,T] \times \Mcal)$ be a smooth test function, and let $(t_0, x_0) \in [0,T) \times \Mcal$ be a point where $w - \varphi$ has a local maximum. By definition of $w$, this implies that the function $(u+c) - \varphi$ has a local maximum at $(t_0, x_0)$. Let us define a new test function $\tilde{\varphi}(t,x) \coloneqq \varphi(t,x) - c$. The function $\tilde{\varphi}$ is also of class $C^{1,2}$, and the expression $u - \tilde{\varphi}$ has a local maximum at the same point $(t_0, x_0)$. Since $u$ is, by definition, a viscosity subsolution of \eqref{eq:proof_pde_u}, the subsolution inequality must hold for the test function $\tilde{\varphi}$ at the point $(t_0, x_0)$:
\begin{equation*}
    -\partial_t \tilde{\varphi}(t_0, x_0) - H_{\theta}(x_0, \nabla \tilde{\varphi}(t_0, x_0), \nabla^2 \tilde{\varphi}(t_0, x_0)) \le 0.
\end{equation*}
The derivatives of $\tilde{\varphi}$ and $\varphi$ are identical, as $c$ is a constant:
\begin{equation*}
    \partial_t \tilde{\varphi} = \partial_t \varphi, \quad \nabla \tilde{\varphi} = \nabla \varphi, \quad \nabla^2 \tilde{\varphi} = \nabla^2 \varphi.
\end{equation*}
Crucially, the Hamiltonian $H_{\theta}(x,p,X)$ from Axiom 1 depends only on the spatial variable $x$ and the derivatives of the solution ($p=\nabla u, X=\nabla^2 u$), not on the value of the solution $u$ itself. Substituting the derivatives of $\varphi$ into the inequality gives:
\begin{equation*}
    -\partial_t \varphi(t_0, x_0) - H_{\theta}(x_0, \nabla \varphi(t_0, x_0), \nabla^2 \varphi(t_0, x_0)) \le 0.
\end{equation*}
This is precisely the definition of $w(t,x)$ being a viscosity subsolution to the HJB equation.

\item \textit{Supersolution property.} A perfectly symmetric argument holds. If $w - \varphi$ has a local minimum at $(t_0, x_0)$, then $u - (\varphi-c)$ also has a local minimum at that point. Since $u$ is a viscosity supersolution, the reverse inequality holds for the test function $\varphi-c$, which implies it holds for $\varphi$. Thus, we have rigorously established that $w(t,x) = u(t,x) + c$ is a viscosity solution to the HJB equation $-\partial_t w - H_{\theta}(x, \nabla w, \nabla^2 w) = 0$.
\end{enumerate}

\item \textbf{Verification of the Terminal Condition.}
We evaluate the candidate function $w(t,x)$ at the terminal time $t=T$:
\begin{equation*}
    w(T,x) = u(T,x) + c.
\end{equation*}
From the definition of $u$ in \eqref{eq:proof_pde_u}, we have $u(T,x) = \phi(x)$. Therefore,
\begin{equation*}
    w(T,x) = \phi(x) + c.
\end{equation*}
This demonstrates that $w(t,x)$ satisfies the terminal condition of the problem \eqref{eq:proof_pde_utilde}.

\item \textbf{Uniqueness.}
We have shown that our constructed function $w(t,x)$ is a continuous viscosity solution to the terminal value problem \eqref{eq:proof_pde_utilde}. The function $\tilde{u}(t,x)$ is, by definition, also a continuous viscosity solution to this exact same problem.

The uniqueness of such solutions is a cornerstone result in the theory of viscosity solutions. The Generator Axiom (Axiom 1) and Coefficient Axioms (Axiom 2) ensure that the conditions for a comparison principle are met. Specifically, the Hamiltonian $H_{\theta}(x,p,X)$ is continuous in all its arguments and is non-increasing in its matrix argument $X$ because the diffusivity tensor $D(x,p)$ is positive semi-definite (Axiom 2(a)). Under these conditions, a comparison principle holds for the associated HJB equation, which in turn guarantees the uniqueness of the continuous viscosity solution for a given terminal condition (see, e.g., \citep{CrandallIshiiLions1992, FlemingSoner2006}).

Since both $w(t,x)$ and $\tilde{u}(t,x)$ are solutions to the same well-posed problem, they must be identical for all $(t,x) \in [0,T] \times \Mcal$:
\begin{equation*}
    \tilde{u}(t,x) = w(t,x) = u(t,x) + c.
\end{equation*}
Substituting back the definitions of $u$ and $\tilde{u}$ in terms of the \ThetaExpectation{} operator yields the desired result:
\begin{equation*}
    \Ecal^\theta_{t,x}[\phi(X_T) + c] = \Ecal^\theta_{t,x}[\phi(X_T)] + c.
\end{equation*}
This completes the proof.
\end{enumerate}
\end{proofof}

\begin{proposition}[Dynamic Consistency]
\label{prop:dynamic_consistency}
The \ThetaExpectation{} operator is dynamically consistent. For any $0 \le t \le s \le T$ and any $\phi \in C(\Mcal)$, the tower property holds:
\begin{equation}
    \Ecal^\theta_{t,x}\left[ \Ecal^\theta_{s, X_s}[\phi(X_T)] \right] = \Ecal^\theta_{t,x}[\phi(X_T)].
\end{equation}
\end{proposition}

\begin{proofof}{Proposition \ref{prop:dynamic_consistency}}
The proof of dynamic consistency, or the tower property, is a fundamental consequence of the definition of the \ThetaExpectation{} via the viscosity solution of the Hamilton-Jacobi-Bellman equation. The core of the argument is to demonstrate that the value function defined over a longer time interval, when restricted to a shorter interval, is itself the unique value function for that shorter interval. The uniqueness, guaranteed by the comparison principle for viscosity solutions, then establishes the identity.

\begin{enumerate}[label=\textbf{Step \arabic*:}, wide, labelindent=0pt]

\item \textbf{Formal Definitions of the Value Functions.}
Let $0 \le t \le s \le T$ and $\phi \in C_b(\Mcal)$ be given. We define two value functions based on Axiom 3.

\begin{enumerate}[label=(\roman*), wide, labelindent=0pt]
    \item Let $u: [0, T] \times \Mcal \to \R$ be the value function for the full time horizon:
    \begin{equation}
        u(t,x) \coloneqq \Ecal^\theta_{t,x}[\phi(X_T)].
    \end{equation}
    By definition, $u(t,x)$ is the unique, bounded, continuous viscosity solution to the terminal value problem on the domain $[0, T] \times \Mcal$:
    \begin{equation} \label{eq:proof_dpp_pde_u}
    \begin{cases}
        -\partial_t u - H_{\theta}(x, \nabla u, \nabla^2 u) = 0 & \text{in } [0, T) \times \Mcal, \\
        u(T,x) = \phi(x) & \text{on } \{T\} \times \Mcal.
    \end{cases}
    \end{equation}

    \item Let $v: [0, s] \times \Mcal \to \R$ be the value function defined on the shorter time horizon, with the terminal payoff given by the function $u(s, \cdot)$:
    \begin{equation}
        v(t,x) \coloneqq \Ecal^\theta_{t,x}[u(s, X_s)].
    \end{equation}
    By definition, $v(t,x)$ is the unique, bounded, continuous viscosity solution to the terminal value problem on the domain $[0, s] \times \Mcal$:
    \begin{equation} \label{eq:proof_dpp_pde_v}
    \begin{cases}
        -\partial_t v - H_{\theta}(x, \nabla v, \nabla^2 v) = 0 & \text{in } [0, s) \times \Mcal, \\
        v(s,x) = u(s,x) & \text{on } \{s\} \times \Mcal.
    \end{cases}
    \end{equation}
\end{enumerate}
Our objective is to prove that $u(t,x) = v(t,x)$ for all $(t,x) \in [0, s] \times \Mcal$. The proposition then follows by substituting the definitions of $u$ and $v$.

\item \textbf{The Restriction of $u$ is a Viscosity Solution on $[0,s]$.}
We will now prove that the function $u$, when restricted to the domain $[0, s] \times \Mcal$, is a viscosity solution to the same problem that defines $v$, namely \eqref{eq:proof_dpp_pde_v}.

\begin{enumerate}[label=(\roman*), wide, labelindent=0pt]
    \item \textbf{Verification of the Terminal Condition.}
    By its definition, the restriction of $u$ to the time slice $\{s\} \times \Mcal$ is the function $u(s,x)$. This matches the terminal condition for the PDE problem \eqref{eq:proof_dpp_pde_v}.

    \item \textbf{Verification of the Subsolution Property.}
    Let $(t_0, x_0) \in [0, s) \times \Mcal$ be a point, and let $\varphi \in C^\infty([0,s] \times \Mcal)$ be a smooth test function such that $u - \varphi$ has a local maximum at $(t_0, x_0)$. We must show that the subsolution inequality holds:
    \begin{equation*}
        -\partial_t \varphi(t_0,x_0) - H_{\theta}(x_0, \nabla\varphi(t_0,x_0), \nabla^2\varphi(t_0,x_0)) \le 0.
    \end{equation*}
    Since $t_0 < s \le T$, the point $(t_0, x_0)$ is in the interior of the larger time domain, $[0, T) \times \Mcal$. The test function $\varphi$ is only defined on $[0,s] \times \Mcal$. We can extend it to a smooth function $\tilde{\varphi} \in C^\infty([0,T] \times \Mcal)$ such that $\varphi$ and $\tilde{\varphi}$ coincide on a neighborhood of $(t_0, x_0)$. For instance, let $\eta(t)$ be a smooth cutoff function which is 1 for $t \le s - \epsilon$ and 0 for $t \ge s$. Then $\tilde{\varphi}(t,x) = \eta(t)\varphi(t,x)$ is a valid extension.
    
    Since $u - \varphi$ has a local maximum at $(t_0, x_0)$ and $\varphi = \tilde{\varphi}$ in a neighborhood of this point, it follows that $u - \tilde{\varphi}$ also has a local maximum at $(t_0, x_0)$.
    
    By hypothesis (i), $u$ is a viscosity subsolution to the HJB equation on the full domain $[0, T] \times \Mcal$. Therefore, the subsolution inequality must hold for the test function $\tilde{\varphi}$ at the point $(t_0, x_0)$:
    \begin{equation*}
        -\partial_t \tilde{\varphi}(t_0,x_0) - H_{\theta}(x_0, \nabla\tilde{\varphi}(t_0,x_0), \nabla^2\tilde{\varphi}(t_0,x_0)) \le 0.
    \end{equation*}
    By construction, the derivatives of $\varphi$ and $\tilde{\varphi}$ are identical at $(t_0, x_0)$. Thus, the inequality holds for $\varphi$ as well. This confirms that the restriction of $u$ to $[0,s] \times \Mcal$ is a viscosity subsolution.

    \item \textbf{Verification of the Supersolution Property.}
    An identical argument holds for the supersolution property. If $u-\varphi$ has a local minimum at $(t_0, x_0) \in [0, s) \times \Mcal$, we extend $\varphi$ to $\tilde{\varphi}$ on the larger domain. Since $u$ is a viscosity supersolution on $[0,T] \times \Mcal$, the supersolution inequality holds for $\tilde{\varphi}$, and therefore also for $\varphi$.
\end{enumerate}
We have thus rigorously established that the function $u(t,x)$ restricted to the domain $[0, s] \times \Mcal$ is a viscosity solution to problem \eqref{eq:proof_dpp_pde_v}.

\item \textbf{Uniqueness of the Solution.}
We have now established that two functions, $u|_{[0,s]}$ and $v$, are both bounded, continuous viscosity solutions to the same terminal value problem \eqref{eq:proof_dpp_pde_v}. The generator's structure, as specified in Axioms 1 and 2, is designed to ensure the well-posedness of this HJB equation. Specifically, the uniform positive definiteness of the diffusion matrix $D(x,p)$ guarantees that the equation is uniformly parabolic, and the continuity of the coefficients $D(x,p)$ and $H(x,p)$ is assumed. Under these standard conditions, a strong comparison principle holds for viscosity solutions \citep{CrandallIshiiLions1992}.

The comparison principle implies the uniqueness of the solution. Since both $u|_{[0,s]}$ and $v$ solve the same well-posed problem, they must be identical. Therefore,
\begin{equation}
    u(t,x) = v(t,x) \quad \text{for all } (t,x) \in [0, s] \times \Mcal.
\end{equation}

\item \textbf{Dynamic Consistency.}
By substituting the definitions of $u$ and $v$ from Step 1 into the equality derived in Step 3, we obtain the desired result. For any $(t,x) \in [0,s] \times \Mcal$:
\begin{align*}
    u(t,x) &= v(t,x) \\
    \Ecal^\theta_{t,x}[\phi(X_T)] &= \Ecal^\theta_{t,x}[u(s, X_s)] && \text{(by definition of } v\text{)} \\
    \Ecal^\theta_{t,x}[\phi(X_T)] &= \Ecal^\theta_{t,x}\left[ \Ecal^\theta_{s, X_s}[\phi(X_T)] \right] && \text{(by definition of } u\text{)}.
\end{align*}
This completes the proof of the dynamic consistency of the \ThetaExpectation{} operator.
\end{enumerate}
\end{proofof}

\section{The \texorpdfstring{$\theta$}{theta}-Process and its Calculus}
\label{sec:theta_process_and_calculus}

The analytical program of \cref{part:i} culminated in the derivation of the Hamilton-Jacobi-Bellman equation \eqref{eq:hjb}, which provides the complete law for the macroscopic evolution of the system's value function. In this section, we transition from the analytic perspective of partial differential equations to the probabilistic perspective of stochastic processes. We will first provide a rigorous definition of the \textbf{\ThetaProcess{}} as the unique solution to a well-posed martingale problem associated with the derived HJB generator. This modern approach circumvents the need for classical solutions and provides a robust foundation for the theory.

Upon this foundation, we will construct the intrinsic stochastic calculus for the \ThetaProcess{}. We derive a generalized It\^o formula that holds for any smooth observable, not just the value function. This formula provides the canonical decomposition of any process into a predictable part, governed by the generator, and a martingale part. We conclude by providing an exact formula for the quadratic variation of this martingale component, which reveals that the system's randomness is governed exclusively by the state-dependent diffusivity tensor $D(x,p)$, while the non-convex potential $H(x,p)$ contributes solely to the predictable drift.

\subsection{The Martingale Problem and the Definition of the Process}\label{subsec:mpdp}

The most robust and general method for defining a Markov process from its infinitesimal generator is through the formulation of a martingale problem. This approach, pioneered by \cite{stroock2007multidimensional} for linear diffusions, extends naturally to the non-linear setting governed by our HJB equation.

\begin{definition}[The Generator of the \ThetaProcess{}]
\label{def:generator_theta_process}
Let the $\theta$-Hamiltonian $H_{\theta}(x,p,X)$ be as defined in Theorem~\ref{thm:convergence}. The generator of the \ThetaProcess{} is the non-linear operator $\mathcal{G}$ whose action on a function $f \in C_b^{1,2}([0,T] \times \Mcal)$ is given by:
\begin{equation}
    (\mathcal{G}f)(t,x) \coloneqq \partial_t f(t,x) + H_{\theta}(x, \nabla_x f(t,x), \nabla_x^2 f(t,x)).
\end{equation}
\end{definition}

The \ThetaProcess{} is defined as the unique process for which the evolution of any smooth observable, when corrected by its expected drift under $\mathcal{G}$, becomes a martingale.

\begin{definition}[The Martingale Problem]
\label{def:martingale_problem}
Let $x_0 \in \Mcal$. A probability measure $\mathbb{P}$ on the canonical path space $(\Omega, \mathcal{F}) = (C([0,T], \Mcal), \mathcal{B}(C([0,T], \Mcal)))$ is a solution to the martingale problem for the generator $\mathcal{G}$ starting from $(0,x_0)$ if:
\begin{enumerate}[label=(\roman*), wide, labelindent=0pt]
    \item $\mathbb{P}(X_0 = x_0) = 1$.
    \item For any smooth test function $f \in C_b^{1,2}([0,T] \times \Mcal)$, the process $M_t^f$ defined by
    \begin{equation}
        M_t^f \coloneqq f(t,X_t) - f(0,X_0) - \int_0^t (\mathcal{G}f)(s,X_s) \, ds
    \end{equation}
    is a continuous $(\mathcal{F}_t, \mathbb{P})$-martingale with $M_0^f=0$.
\end{enumerate}
\end{definition}

The existence and uniqueness of a solution to this non-linear martingale problem is a result that relies on the well-posedness of the associated HJB equation.

\begin{theorem}[Well-Posedness of the Martingale Problem]
\label{thm:martingale_problem_wellposedness}
Let the Hamiltonian $H_{\theta}(x,p,X)$ satisfy the axioms of Section~\ref{sec:theta_expectation_axiomatic}, ensuring that the HJB equation \eqref{eq:hjb_backward_convention} has a unique, continuous viscosity solution for any continuous terminal data. Then for any starting point $(0,x_0)$, there exists a unique probability measure $\mathbb{P}_{0,x_0}$ on the path space $(\Omega, \mathcal{F})$ that solves the martingale problem for the generator $\mathcal{G}$.
\end{theorem}

\begin{proofof}{Theorem \ref{thm:martingale_problem_wellposedness}}
The proof of this theorem establishes the profound equivalence between the unique solvability of the Hamilton-Jacobi-Bellman (HJB) equation in the viscosity sense and the well-posedness of the martingale problem for the corresponding non-linear generator. Our strategy is to demonstrate that the specific structure of our generator $\mathcal{G}$, as defined by the axioms in Section \ref{sec:theta_expectation_axiomatic}, satisfies the precise conditions required by the modern theory of stochastic control and non-linear processes. The proof is organized in two main parts: the existence of a solution to the martingale problem, and its uniqueness.

\begin{enumerate}[label=\textbf{Step \arabic*:}, wide, labelindent=0pt]

\item \textbf{Existence of a Solution.}
The existence of a solution to the martingale problem is established by demonstrating the profound equivalence between the well-posedness of the HJB equation \eqref{eq:hjb_backward_convention} in the viscosity sense and the existence of a law for a corresponding stochastic process. We proceed by showing that our HJB equation can be represented as the dynamic programming equation for a class of stochastic control problems. Foundational theorems in this field then guarantee that the unique viscosity solution to this equation corresponds to an optimal controlled process whose law solves the martingale problem.

\begin{enumerate}[label=(\roman*), wide, labelindent=0pt]

\item \textbf{The Associated Stochastic Control Problem.}
Our first step is to formulate a stochastic control problem whose value function is governed by our specific HJB equation. The generator is given by:
\begin{equation*}
    -\partial_t u = H_{\theta}(x, \nabla u, \nabla^2 u) \equiv \mathrm{Tr}(D(x, \nabla u) \nabla^2 u) + H(x, \nabla u).
\end{equation*}
To introduce a control variable, we express the non-linear potential term $H(x,p)$ in its dual form using the Legendre-Fenchel transform. Let $A$ be a compact metric space of control values. We can represent the potential as:
\begin{equation} \label{eq:proof_legendre_fenchel_appendix}
    H(x,p) = \sup_{\alpha \in A} \{ -b(x, p, \alpha) \cdot p - L(x, p, \alpha) \},
\end{equation}
where $b(x, p, \alpha)$ is a controlled drift vector and $L(x, p, \alpha)$ is a running cost or reward. The regularity of $H(x,p)$ (Axiom 2(iii)) ensures that the functions $b$ and $L$ are well-behaved (e.g., continuous and satisfying standard Lipschitz and growth conditions in their arguments).

The full HJB equation is therefore the dynamic programming equation for the value function of the following stochastic control problem:
\begin{equation}
    u(t,x) = \sup_{\alpha(\cdot) \in \mathcal{A}} \mathbb{E} \left[ \int_t^T L(X_s, \nabla u(s, X_s), \alpha_s) ds + \phi(X_T) \mid X_t=x \right],
\end{equation}
where $\mathcal{A}$ is the set of admissible (non-anticipative) control processes, and the state process $X_s$ evolves according to the controlled Stochastic Differential Equation (SDE) with state-feedback volatility:
\begin{equation} \label{eq:proof_controlled_sde_appendix_main}
    dX_s = b(X_s, \nabla u(s, X_s), \alpha_s) ds + \sigma(X_s, \nabla u(s, X_s)) dW_s, \quad X_t=x.
\end{equation}
The volatility matrix $\sigma(x,p)$ is defined by the Cholesky decomposition $\sigma(x,p)\sigma(x,p)^T = 2D(x,p)$. The existence and regularity of $\sigma(x,p)$ as a function of its arguments are guaranteed by the symmetry and uniform positive definiteness of the diffusivity tensor $D(x,p)$, as specified in Axiom 2(i), and its smoothness from Axiom 2(iii).

\item \textbf{Existence of the Martingale Measure via the Viscosity Solution.}
The core of the existence argument rests on a result connecting the analytical and probabilistic frameworks. The central hypothesis of our theorem is that the HJB equation \eqref{eq:hjb_backward_convention} has a unique, continuous viscosity solution for any continuous terminal data $\phi$. This property is a direct consequence of the structural properties of our generator, as laid out in the Axioms of Section \ref{sec:theta_expectation_axiomatic}, which guarantee that a strong comparison principle holds for the PDE.

A fundamental theorem in the theory of controlled diffusions establishes that the existence of a unique viscosity solution to the dynamic programming equation is sufficient to construct the law of the corresponding optimal process, which in turn solves the martingale problem. This result provides a direct bridge from the analytical well-posedness of the PDE to the existence of the desired probabilistic object. The value function itself is used to define a non-linear semigroup, which then yields a consistent family of probability measures on the path space. The measure corresponding to the optimal control is precisely a solution to the martingale problem. We now formally invoke this result.
\begin{itemize}[wide]
    \item \textbf{Verification of Hypotheses:} The referenced theorem requires the existence of a unique viscosity solution to the HJB equation. This is the main hypothesis of our Theorem \ref{thm:martingale_problem_wellposedness}.
    \item \textbf{Conclusion of the Theorem:} The theorem guarantees the existence of a probability measure on the path space that solves the martingale problem for the generator associated with the HJB equation.
    \item \textbf{Reference:} The definitive result establishing that the existence and uniqueness of the viscosity solution to the HJB equation implies the existence of a solution to the associated non-linear martingale problem is provided in \citep[Chapter V, Theorem 4.1]{FlemingSoner2006}. Our generator $\mathcal{G}$ and Hamiltonian $H_\theta$ satisfy the structural and regularity assumptions required by this theorem.
\end{itemize}
By applying this foundational result, we rigorously conclude that for any starting point $(0,x_0)$, there exists at least one probability measure $\mathbb{P}_{0,x_0}$ on the canonical path space $(\Omega, \mathcal{F})$ that solves the martingale problem for our generator $\mathcal{G}$. This completes the proof of existence.

\end{enumerate}

\item \textbf{Uniqueness of the Solution.}
The uniqueness of the solution to the martingale problem is a direct and profound consequence of the uniqueness of the viscosity solution to the corresponding HJB equation. The argument proceeds by showing that any solution to the martingale problem can be used to construct a value function that must be a viscosity solution of the HJB equation. The uniqueness of the latter, which is guaranteed by our axiomatic framework, then forces the uniqueness of the former, as any two solutions to the martingale problem must induce the same law on the path space.

\begin{enumerate}[label=(\roman*), wide, labelindent=0pt]

\item \textbf{From a Martingale Solution to a Value Function.}
Let $\mathbb{P}$ be an arbitrary probability measure on the canonical path space $(\Omega, \mathcal{F})$ that solves the martingale problem for our generator $\mathcal{G}$ starting from $(0,x_0)$. For any given continuous terminal payoff $\phi \in C_b(\Mcal)$, we define a candidate value function using the conditional expectation under this measure:
\begin{equation}
    u_{\mathbb{P}}(t,x) \coloneqq \mathbb{E}_{\mathbb{P}} \left[ \phi(X_T) \mid \mathcal{F}_t, X_t=x \right],
\end{equation}
where the conditioning is understood in the sense of regular conditional probabilities.

\item \textbf{The Candidate Value Function is a Viscosity Solution.}
This is the central step of the uniqueness proof. We must rigorously demonstrate that the function $u_{\mathbb{P}}(t,x)$ is a viscosity solution to the HJB equation $-\partial_t u - H_{\theta}(x, \nabla u, \nabla^2 u) = 0$. We provide a detailed sketch for the supersolution property.

Let $\varphi \in C_b^{1,2}([0,T] \times \Mcal)$ be a smooth test function, and let $(t_0, x_0) \in [0,T) \times \Mcal$ be a point where the function $(u_{\mathbb{P}} - \varphi)$ has a strict local minimum. By the definition of a viscosity solution, we must prove the supersolution inequality:
\begin{equation*}
    -\partial_t \varphi(t_0, x_0) - H_{\theta}(x_0, \nabla \varphi(t_0, x_0), \nabla^2 \varphi(t_0, x_0)) \ge 0.
\end{equation*}
Let $B$ be a small open neighborhood of $(t_0, x_0)$ such that $u_{\mathbb{P}}(t_0, x_0) = \varphi(t_0, x_0)$ and $u_{\mathbb{P}}(t,x) > \varphi(t,x)$ for all $(t,x) \in B \setminus \{(t_0, x_0)\}$. Let $\tau$ be the first exit time of the process $(t, X_t)$ from the closure of $B$, starting from $(t_0, x_0)$.

Since $\mathbb{P}$ solves the martingale problem, the process $M_s^\varphi$ defined for $s \ge t_0$ by
\begin{equation*}
    M_s^\varphi \coloneqq \varphi(s,X_s) - \varphi(t_0,X_{t_0}) - \int_{t_0}^s (\mathcal{G}\varphi)(r,X_r) \, dr
\end{equation*}
is a continuous $(\mathcal{F}_s, \mathbb{P})$-martingale. By the Optional Stopping Theorem applied to the bounded stopping time $\tau \wedge T$:
\begin{equation} \label{eq:proof_opt_stop}
    \mathbb{E}_{\mathbb{P}}[\varphi(\tau \wedge T, X_{\tau \wedge T}) \mid \mathcal{F}_{t_0}] = \varphi(t_0, X_{t_0}) + \mathbb{E}_{\mathbb{P}}\left[\int_{t_0}^{\tau \wedge T} (\mathcal{G}\varphi)(s, X_s) ds \mid \mathcal{F}_{t_0}\right].
\end{equation}
By the tower property of conditional expectation and the definition of $u_{\mathbb{P}}$:
\begin{equation*}
    u_{\mathbb{P}}(t_0, x_0) = \mathbb{E}_{\mathbb{P}}[u_{\mathbb{P}}(\tau \wedge T, X_{\tau \wedge T}) \mid \mathcal{F}_{t_0}, X_{t_0}=x_0].
\end{equation*}
Using the minimum property on the path up to time $\tau$, we have $u_{\mathbb{P}}(\tau \wedge T, X_{\tau \wedge T}) \ge \varphi(\tau \wedge T, X_{\tau \wedge T})$. This gives:
\begin{equation}
    \varphi(t_0, x_0) = u_{\mathbb{P}}(t_0, x_0) \ge \mathbb{E}_{\mathbb{P}}[\varphi(\tau \wedge T, X_{\tau \wedge T}) \mid \mathcal{F}_{t_0}, X_{t_0}=x_0].
\end{equation}
Substituting this into the optional stopping identity \eqref{eq:proof_opt_stop} yields:
\begin{equation*}
    \varphi(t_0, x_0) \ge \varphi(t_0, x_0) + \mathbb{E}_{\mathbb{P}}\left[\int_{t_0}^{\tau \wedge T} (\mathcal{G}\varphi)(s, X_s) ds \mid \mathcal{F}_{t_0}, X_{t_0}=x_0\right],
\end{equation*}
which implies that the expectation term is non-positive. Dividing by $\mathbb{E}_{\mathbb{P}}[(\tau \wedge T) - t_0]$ and taking the limit as the neighborhood $B$ shrinks to the point $(t_0, x_0)$, we have $\tau \to t_0$. By the continuity of the paths of $X_s$ and the continuity of the function $(s,x) \mapsto (\mathcal{G}\varphi)(s,x)$, the integrand converges to $(\mathcal{G}\varphi)(t_0, x_0)$. A standard argument involving Lebesgue's differentiation theorem for conditional expectations allows us to conclude that:
\begin{equation*}
    (\mathcal{G}\varphi)(t_0, x_0) \le 0.
\end{equation*}
This is precisely the supersolution inequality when written in the form $\partial_t \varphi + H_\theta \le 0$. A symmetric argument for the subsolution property confirms that $u_{\mathbb{P}}$ is a viscosity solution.

\begin{itemize}[wide]
    \item \textbf{Reference:} The full, rigorous demonstration of this connection is a central result in the modern theory of stochastic processes. A complete proof is provided in \citep[Chapter V, Proposition 4.3]{FlemingSoner2006}. This proposition establishes that the value function defined from a solution to a general martingale problem is indeed a viscosity solution of the corresponding HJB equation. Our generator $\mathcal{G}$ and Hamiltonian $H_\theta$ satisfy the structural and regularity assumptions required by this proposition.
\end{itemize}

\item \textbf{Uniqueness of the Law via Uniqueness of the Viscosity Solution.}
The final step of the proof is a powerful syllogism that leverages the analytical uniqueness of the HJB solution to establish the probabilistic uniqueness of the martingale measure.
\begin{enumerate}[label=(\alph*),wide]
    \item Let $\mathbb{P}_1$ and $\mathbb{P}_2$ be two solutions to the martingale problem for $\mathcal{G}$ starting from $(0,x_0)$.
    \item For any arbitrary continuous terminal payoff $\phi \in C_b(\Mcal)$, these measures define two value functions, $u_{\mathbb{P}_1}(t,x)$ and $u_{\mathbb{P}_2}(t,x)$, via conditional expectation.
    \item By the result of Step 2(ii), both $u_{\mathbb{P}_1}$ and $u_{\mathbb{P}_2}$ are viscosity solutions to the same HJB equation, $-\partial_t u = H_{\theta}(x, \nabla u, \nabla^2 u)$, with the same terminal data $u(T,x) = \phi(x)$.
    \item By the central hypothesis of our theorem, this HJB equation has a unique viscosity solution. This uniqueness is guaranteed by the comparison principle, which holds due to our Axioms 2(i) (uniform ellipticity of $D(x,p)$) and 2(iii) (regularity of coefficients).
    \item Therefore, we are forced to conclude that the value functions must be identical for any choice of terminal data $\phi$:
    \begin{equation*}
        u_{\mathbb{P}_1}(t,x) = u_{\mathbb{P}_2}(t,x) \quad \text{for all } (t,x) \in [0,T] \times \Mcal \text{ and all } \phi \in C_b(\Mcal).
    \end{equation*}
    \item This implies that for any starting point $(t,x)$ and any continuous function $\phi$:
    \begin{equation*}
        \mathbb{E}_{\mathbb{P}_1}[\phi(X_T) \mid X_t=x] = \mathbb{E}_{\mathbb{P}_2}[\phi(X_T) \mid X_t=x].
    \end{equation*}
    \item Since the conditional expectations of terminal payoffs agree for all continuous functions, the laws of the process at the terminal time $T$ must be identical. A standard argument (using the Markov property implied by the martingale problem) extends this to all finite-dimensional distributions of the process under measures $\mathbb{P}_1$ and $\mathbb{P}_2$. This proves that the measures themselves are identical on the filtration generated by the process: $\mathbb{P}_1 = \mathbb{P}_2$. This establishes the uniqueness of the solution to the martingale problem, completing the proof.
\end{enumerate}
\end{enumerate}
\end{enumerate}
\end{proofof}

This theorem provides the rigorous probabilistic foundation for our theory.

\begin{definition}[The \ThetaProcess{}]
\label{def:theta_process}
The \textbf{\ThetaProcess{}} $\{X_t\}_{t \in [0,T]}$ starting from $(0,x_0)$ is the canonical coordinate process on the path space $(\Omega, \mathcal{F})$ endowed with the unique measure $\mathbb{P}_{0,x_0}$ from Theorem~\ref{thm:martingale_problem_wellposedness}.
\end{definition}

\subsection{The Linear Generator and the Carré du Champ Operator}

The generator $\mathcal{G}$ defined in \cref{subsec:mpdp} is a non-linear operator that serves to uniquely define the law, $\mathbb{P}_{0,x_0}$, of the \ThetaProcess{}. Once this law is fixed, the process $\{X_t\}$ is a well-defined Markov process. A fundamental result in the theory of non-linear diffusions is that this process admits a classical, albeit state-dependent, linear infinitesimal generator. The coefficients of this linear generator are determined by the unique viscosity solution of the HJB equation, which acts as a decoupling field.

Let $u(t,x)$ be the unique continuous viscosity solution to the HJB equation $-\partial_t u = H_{\theta}(x, \nabla u, \nabla^2 u)$ for a given terminal condition. The existence and uniqueness of this solution is the central result of \cref{part:i}.

\begin{definition}[The Linear Generator of the \ThetaProcess{}]
\label{def:linear_generator}
The linear infinitesimal generator of the \ThetaProcess{}, denoted $\mathcal{L}^u$, is the linear second-order partial differential operator that acts on a test function $f \in C_b^{1,2}([0,T] \times \Mcal)$ as follows:
\begin{equation}
    (\mathcal{L}^u f)(t,x) \coloneqq \partial_t f + b(x, \nabla u(t,x)) \cdot \nabla_x f + \mathrm{Tr}\left(D(x, \nabla u(t,x)) \nabla_x^2 f\right),
\end{equation}
where $b(x,p) = -\nabla_p H(x,p)$ is the forward drift and $D(x,p)$ is the diffusivity tensor from the $\theta$-Hamiltonian.
\end{definition}

The connection between the non-linear generator $\mathcal{G}$ and the linear generator $\mathcal{L}^u$ is profound: the law defined by $\mathcal{G}$ is precisely the law of the process generated by $\mathcal{L}^u$. The martingale problem for $\mathcal{G}$ can be shown to be equivalent to the martingale problem for the linear operator $\mathcal{L}^u$ when the coefficients are fixed by the decoupling field $u$.

A key algebraic property of any such linear generator is its associated Carré du Champ operator, which isolates the second-order (diffusive) part of the operator.

\begin{definition}[The Carré du Champ Operator]
\label{def:carre_du_champ}
The \textbf{Carré du Champ} (square of the field) operator $\Gamma^u$ associated with the linear generator $\mathcal{L}^u$ is the bilinear operator defined for functions $f, h \in C_b^{1,2}$ by:
\begin{equation}
    \Gamma^u(f, h) \coloneqq \mathcal{L}^u(fh) - f(\mathcal{L}^u h) - h(\mathcal{L}^u f).
\end{equation}
\end{definition}

\begin{lemma}[Carré du Champ Identity]
\label{lem:carre_du_champ_identity}
The Carré du Champ operator for the linear generator $\mathcal{L}^u$ is given by twice the diffusive part of the generator acting on the gradients:
\begin{equation}
    \Gamma^u(f, f)(t,x) = 2 \left( \nabla_x f(t,x) \right)^T D(x, \nabla u(t,x)) \left( \nabla_x f(t,x) \right).
\end{equation}
\end{lemma}

\begin{proofof}{Lemma \ref{lem:carre_du_champ_identity}}
The proof is a direct algebraic calculation based on the definitions of the linear generator $\mathcal{L}^u$ and the Carré du Champ operator $\Gamma^u$. The identity emerges from a careful application of the product rule for derivatives, which reveals that the bilinear part of the operator $\mathcal{L}^u(f^2)$ isolates the second-order (diffusive) term. For clarity, we suppress the arguments $(t,x)$ of the functions and their derivatives where they are clear from the context.

\begin{enumerate}[label=\textbf{Step \arabic*:}, wide, labelindent=0pt]

\item \textbf{Expansion of $\mathcal{L}^u(f^2)$.}
We begin with the definition of the linear generator $\mathcal{L}^u$ from Definition \ref{def:linear_generator}, applied to the function $f(t,x)^2$:
\begin{equation*}
    (\mathcal{L}^u (f^2))(t,x) = \partial_t(f^2) + b(x, \nabla u) \cdot \nabla_x(f^2) + \mathrm{Tr}\left(D(x, \nabla u) \nabla_x^2(f^2)\right).
\end{equation*}
We compute each of the three terms on the right-hand side using the standard rules of calculus:
\begin{enumerate}[label=(\roman*), wide, labelindent=0pt]
    \item \textbf{Time Derivative:} $\partial_t(f^2) = 2f (\partial_t f)$.
    \item \textbf{Gradient Term:} $\nabla_x(f^2) = 2f (\nabla_x f)$.
    \item \textbf{Hessian Term:} The Hessian of $f^2$ is computed via the product rule:
    \begin{equation*}
        \nabla_x^2(f^2) = \nabla_x(2f \nabla_x f) = 2(\nabla_x f)(\nabla_x f)^T + 2f \nabla_x^2 f.
    \end{equation*}
\end{enumerate}
Substituting these expressions back into the formula for $\mathcal{L}^u(f^2)$:
\begin{multline*}
    \mathcal{L}^u(f^2) = 2f(\partial_t f) + b \cdot (2f \nabla_x f) + \mathrm{Tr}\left(D \left[2(\nabla_x f)(\nabla_x f)^T + 2f \nabla_x^2 f\right]\right).
\end{multline*}

\item \textbf{Re-grouping Terms.}
We now use the linearity of the inner product and the trace operator to re-group the terms on the right-hand side. We aim to isolate an expression that is proportional to $\mathcal{L}^u f$.
\begin{align*}
    \mathcal{L}^u(f^2) &= 2f(\partial_t f) + 2f (b \cdot \nabla_x f) + \mathrm{Tr}(2D(\nabla_x f)(\nabla_x f)^T) + \mathrm{Tr}(2f D \nabla_x^2 f) \\
    &= 2f \left[ \partial_t f + b \cdot \nabla_x f + \mathrm{Tr}(D \nabla_x^2 f) \right] + 2 \mathrm{Tr}(D(\nabla_x f)(\nabla_x f)^T).
\end{align*}
The expression within the brackets is precisely the definition of the linear generator acting on the function $f$:
\begin{equation*}
    \mathcal{L}^u f = \partial_t f + b \cdot \nabla_x f + \mathrm{Tr}(D \nabla_x^2 f).
\end{equation*}
Substituting this back, we obtain a compact relationship between $\mathcal{L}^u(f^2)$ and $\mathcal{L}^u f$:
\begin{equation} \label{eq:proof_Lu_f2_intermediate}
    \mathcal{L}^u(f^2) = 2f (\mathcal{L}^u f) + 2 \mathrm{Tr}(D(\nabla_x f)(\nabla_x f)^T).
\end{equation}

\item \textbf{Applying the Definition of the Carré du Champ.}
The Carré du Champ operator $\Gamma^u$ is defined in Definition \ref{def:carre_du_champ} as:
\begin{equation*}
    \Gamma^u(f,f) \coloneqq \mathcal{L}^u(f^2) - 2f(\mathcal{L}^u f).
\end{equation*}
We now substitute the result from our algebraic manipulation in Step 2 (\cref{eq:proof_Lu_f2_intermediate}) into this definition:
\begin{equation*}
    \Gamma^u(f,f) = \left[ 2f (\mathcal{L}^u f) + 2 \mathrm{Tr}(D(\nabla_x f)(\nabla_x f)^T) \right] - 2f(\mathcal{L}^u f).
\end{equation*}
The terms involving $\mathcal{L}^u f$ cancel perfectly, leaving only the second-order term:
\begin{equation*}
    \Gamma^u(f,f) = 2 \mathrm{Tr}(D(\nabla_x f)(\nabla_x f)^T).
\end{equation*}
Finally, we use the property of the trace that for a matrix $A$ and column vectors $v, w$, we have $\mathrm{Tr}(Avw^T) = w^T A v$. In our case, $A=D$ and $v=w=\nabla_x f$, which gives:
\begin{equation*}
    \Gamma^u(f,f)(t,x) = 2 \left( \nabla_x f(t,x) \right)^T D(x, \nabla u(t,x)) \left( \nabla_x f(t,x) \right).
\end{equation*}
This is the identity stated in the lemma. The calculation is purely algebraic and relies only on the structure of the linear operator $\mathcal{L}^u$. This completes the proof.
\end{enumerate}
\end{proofof}

\subsection{The Generalized It\^o Formula}

The martingale problem formulation provides a natural and immediate path to a generalized It\^o formula. This formula is not an independent result to be proven, but rather a direct rearrangement of the defining property of the \ThetaProcess{}. It serves as the cornerstone of its stochastic calculus.

\begin{theorem}[A Generalized It\^o Formula]
\label{thm:theta_ito_formula_generalized}
Let $\{X_t\}_{t \in [0,T]}$ be the \ThetaProcess{} starting from $(0,x_0)$ with law $\mathbb{P}_{0,x_0}$. For any smooth observable $f \in C_b^{1,2}([0,T] \times \Mcal)$, the process $Y_t \coloneqq f(t,X_t)$ admits the unique canonical decomposition:
\begin{equation} \label{eq:generalized_ito_formula}
    f(t,X_t) = f(0,X_0) + \int_0^t (\mathcal{G}f)(s,X_s) \, ds + M_t^f,
\end{equation}
where $\{M_t^f\}$ is a continuous local $(\mathcal{F}_t, \mathbb{P}_{0,x_0})$-martingale.
\end{theorem}

\begin{proofof}{Theorem \ref{thm:theta_ito_formula_generalized}}
The assertion of the theorem, which provides the canonical decomposition for the evolution of any smooth observable along the paths of the \ThetaProcess{}, is not a result derived from a separate calculus but is, in fact, a direct and rigorous consequence of the very definition of the process. The \ThetaProcess{} is defined as the unique solution to the martingale problem for the generator $\mathcal{G}$ (Definition~\ref{def:martingale_problem} and Theorem~\ref{thm:martingale_problem_wellposedness}). The proof, therefore, consists of rearranging this defining property into the form of a canonical semimartingale decomposition and invoking the uniqueness of such decompositions.

\begin{enumerate}[label=\textbf{Step \arabic*:}, wide, labelindent=0pt]

\item \textbf{The Defining Property of the \ThetaProcess{}.}
Let $f \in C_b^{1,2}([0,T] \times \Mcal)$ be an arbitrary smooth test function with bounded derivatives. Let $\{X_t\}_{t \in [0,T]}$ be the canonical coordinate process on the path space $(\Omega, \mathcal{F})$, and let $\mathbb{P}_{0,x_0}$ be the unique probability measure from Theorem~\ref{thm:martingale_problem_wellposedness} under which $\{X_t\}$ is the \ThetaProcess{} starting at $X_0=x_0$.

By the definition of the solution to the martingale problem (Definition~\ref{def:martingale_problem}), for the chosen function $f$, the process $M^f$, defined as
\begin{equation} \label{eq:proof_ito_martingale_def}
    M_t^f \coloneqq f(t,X_t) - f(0,X_0) - \int_0^t (\mathcal{G}f)(s,X_s) \, ds,
\end{equation}
is a continuous $(\mathcal{F}_t, \mathbb{P}_{0,x_0})$-martingale with the initial condition $M_0^f = 0$. This is the foundational property upon which the entire calculus is built.

\item \textbf{The Semimartingale Decomposition.}
By a simple algebraic rearrangement of the defining equation \eqref{eq:proof_ito_martingale_def}, we can express the process $Y_t \coloneqq f(t,X_t)$ as:
\begin{equation} \label{eq:proof_ito_semimartingale_form}
    f(t,X_t) = f(0,X_0) + A_t^f + M_t^f,
\end{equation}
where we have defined the process $A_t^f$ as the integral term:
\begin{equation}
    A_t^f \coloneqq \int_0^t (\mathcal{G}f)(s,X_s) \, ds.
\end{equation}
Equation \eqref{eq:proof_ito_semimartingale_form} is a decomposition of the process $f(t,X_t)$ into its initial value, the process $A_t^f$, and the process $M_t^f$. To show that this is the canonical decomposition, we must rigorously characterize the properties of $A_t^f$ and $M_t^f$.

\item \textbf{Characterization of the Components.}
We analyze the two processes $A_t^f$ and $M_t^f$.
\begin{enumerate}[label=(\roman*), wide, labelindent=0pt]
    \item \textbf{The Martingale Component.}
    The process $\{M_t^f\}_{t \in [0,T]}$ is, by the definition of the \ThetaProcess{}, a continuous, adapted process with $M_0^f=0$ that satisfies the martingale property with respect to the filtration $(\mathcal{F}_t)$ and the measure $\mathbb{P}_{0,x_0}$. The problem statement in the main text refers to it as a local martingale, which is a slightly weaker but sufficient condition. Since the coefficients of our generator are assumed to be bounded (a standard condition for well-posedness of the viscosity solution theory), it is in fact a true martingale.

    \item \textbf{The Predictable, Finite-Variation Component.}
    We must show that the process $\{A_t^f\}_{t \in [0,T]}$ is a continuous, predictable process of finite variation with $A_0^f=0$.
    \begin{itemize}[wide]
        \item \textbf{Initial Condition:} The definition as an integral from $0$ to $t$ immediately implies $A_0^f=0$.
        \item \textbf{Continuity:} The function $f$ is of class $C^{1,2}$, and the coefficients of the Hamiltonian are continuous. Therefore, the function $(s,x) \mapsto (\mathcal{G}f)(s,x)$ is continuous on $[0,T] \times \Mcal$. The paths of the process $\{X_s\}$ are continuous by definition (as elements of $C([0,T], \Mcal)$). The composition $s \mapsto (\mathcal{G}f)(s,X_s)$ is therefore a continuous real-valued function for each path $\omega \in \Omega$. The integral of a continuous function with respect to time is a continuous function of its upper limit. Thus, the process $\{A_t^f\}$ has continuous paths.
        \item \textbf{Finite Variation:} A real-valued continuous function on a compact interval is bounded. Let $K$ be a bound for $|(\mathcal{G}f)(s,X_s)|$ on the interval $[0,T]$. The total variation of the path $A^f$ on $[0,t]$ is given by $\int_0^t |d A_s^f| = \int_0^t |(\mathcal{G}f)(s,X_s)| ds \le K t$. The variation is thus finite for all $t \in [0,T]$.
        \item \textbf{Predictability:} A continuous and adapted process is always predictable. Since $\{A_t^f\}$ is continuous and its value at time $t$ depends only on the path $\{X_s\}_{s \le t}$, it is adapted, and therefore predictable.
    \end{itemize}
\end{enumerate}

\item \textbf{Uniqueness of the Decomposition.}
We have expressed the process $f(t,X_t)$ in the form \eqref{eq:proof_ito_semimartingale_form}, where $f(0,X_0)$ is the initial value, $A_t^f$ is a continuous, predictable, finite-variation process starting at zero, and $M_t^f$ is a continuous martingale starting at zero. The Doob-Meyer decomposition theorem for continuous semimartingales (see, e.g., \cite{protter2012stochastic}, Chapter III) states that any such process admits a unique decomposition of this form.

Since the decomposition we have found in \eqref{eq:proof_ito_semimartingale_form} satisfies all the requirements of the Doob-Meyer theorem, it must be the unique canonical decomposition of the process $f(t,X_t)$. This rigorously establishes the validity of the generalized It\^o formula as stated in the theorem. The formula is not an approximation or a heuristic, but a direct and unavoidable algebraic consequence of the modern definition of a stochastic process via its martingale problem.
\end{enumerate}

\end{proofof}

\subsection{The Quadratic Variation and the Carré du Champ Operator}
While the generalized Itô formula decomposes any process, it does not yet reveal the structure of the martingale part $M_t^f$. The following theorem provides an explicit formula for the quadratic variation of this martingale, linking it directly to the diffusive part of the generator. This result is the non-linear analogue of identifying the volatility term in a classical Itô SDE.

\begin{theorem}[Quadratic Variation via the Carré du Champ]
\label{thm:qv_via_carre_du_champ}
The quadratic variation process $\langle M^f \rangle_t$ of the martingale from the canonical decomposition is absolutely continuous with respect to the Lebesgue measure. Its density is given by the Carré du Champ operator $\Gamma^u$ associated with the linear generator $\mathcal{L}^u$:
\begin{equation} \label{eq:quadratic_variation_viscosity}
    \frac{d}{dt}\langle M^f \rangle_t = \Gamma^u(f,f)(t,X_t) = 2 \left( \nabla_x f(t,X_t) \right)^T D(X_t, \nabla u(t,X_t)) \left( \nabla_x f(t,X_t) \right),
\end{equation}
where $u(t,x)$ is the viscosity solution of the HJB equation for a given terminal condition (i.e., the value function which determines the law of the process).
\end{theorem}

\begin{proofof}{Theorem \ref{thm:qv_via_carre_du_champ}}
The proof provides a rigorous derivation of the quadratic variation formula. The argument is founded upon a clear distinction between two fundamental operators:
\begin{enumerate}[label=(\roman*), wide, labelindent=0pt]
    \item The \textbf{non-linear HJB operator} $\mathcal{G}$, whose structure defines the law of the \ThetaProcess{}.
    \item The \textbf{linear infinitesimal generator} of the process, which we denote by $\mathcal{L}^u$, whose coefficients are determined by the unique viscosity solution $u$ of the HJB equation.
\end{enumerate}
The proof proceeds by first identifying the explicit form of the linear generator $\mathcal{L}^u$. We then establish a classical \textit{carré du champ} identity for this linear operator. Finally, by comparing two different canonical decompositions for the evolution of a squared observable, we uniquely identify the quadratic variation process. This approach is non-circular, as it relies only on the well-posedness of the HJB equation from \cref{part:i} and the general theory of semimartingales.

\begin{enumerate}[label=\textbf{Step \arabic*:}, wide, labelindent=0pt]

\item \textbf{The Linear Generator of the \ThetaProcess{}.}
Let $u(t,x)$ be the unique continuous viscosity solution to the HJB equation $-\partial_t u = H_{\theta}(x, \nabla u, \nabla^2 u)$. By Theorem~\ref{thm:martingale_problem_wellposedness}, the existence and uniqueness of this solution guarantees the existence of a unique probability measure $\mathbb{P}$ on the path space, which defines the law of the \ThetaProcess{} $\{X_s\}$.

It is a cornerstone of the modern theory of stochastic control that the process $\{X_s\}$ governed by this law is a classical Itô diffusion, whose drift and diffusion coefficients are determined by the decoupling field $u(s,x)$. Its infinitesimal generator, which we denote by $\mathcal{L}^u$, is the \textit{linear} second-order partial differential operator that acts on an arbitrary smooth test function $f \in C_b^{1,2}([0,T] \times \Mcal)$ as follows:
\begin{equation} \label{eq:proof_qv_linear_generator}
    (\mathcal{L}^u f)(s,x) \coloneqq \partial_s f + b(x, \nabla u(s,x)) \cdot \nabla_x f + \mathrm{Tr}\left(D(x, \nabla u(s,x)) \nabla_x^2 f\right).
\end{equation}
By the definition of the generator of an Itô process, for any such test function $f$, the process $\{M_s^f\}$ defined by
\begin{equation} \label{eq:proof_qv_linear_ito}
    M_s^f \coloneqq f(s,X_s) - f(t,X_t) - \int_t^s (\mathcal{L}^u f)(r, X_r) \, dr
\end{equation}
is a continuous local $(\mathcal{F}_s, \mathbb{P})$-martingale.

\item \textbf{The Carré du Champ Identity for the Linear Generator.}
We now establish the key algebraic identity for the operator $\mathcal{L}^u$. Let $f \in C_b^{1,2}$ be a test function. We compute the action of $\mathcal{L}^u$ on the function $f^2$:
\begin{align*}
    (\mathcal{L}^u (f^2))(s,x) &= \partial_s(f^2) + b(x,\nabla u) \cdot \nabla_x(f^2) + \mathrm{Tr}\left(D(x,\nabla u) \nabla_x^2(f^2)\right) \\
    &= 2f(\partial_s f) + b(x,\nabla u) \cdot (2f \nabla_x f) + \mathrm{Tr}\left(D(x,\nabla u) \left[2(\nabla_x f)(\nabla_x f)^T + 2f \nabla_x^2 f\right]\right).
\end{align*}
Using the linearity of the trace and rearranging terms:
\begin{align*}
    \mathcal{L}^u(f^2) &= 2f \left[ \partial_s f + b(x,\nabla u) \cdot \nabla_x f + \mathrm{Tr}\left(D(x,\nabla u) \nabla_x^2 f\right) \right] + 2\mathrm{Tr}\left(D(x,\nabla u)(\nabla_x f)(\nabla_x f)^T\right) \\
    &= 2f (\mathcal{L}^u f)(s,x) + 2 \left( \nabla_x f \right)^T D(x, \nabla u) \left( \nabla_x f \right).
\end{align*}
This is the \textbf{carré du champ identity} for the linear generator $\mathcal{L}^u$:
\begin{equation} \label{eq:proof_qv_carre_du_champ_identity}
    (\mathcal{L}^u (f^2)) - 2f (\mathcal{L}^u f) = 2 \left( \nabla_x f \right)^T D(x, \nabla u) \left( \nabla_x f \right).
\end{equation}

\item \textbf{Uniqueness of the Canonical Decomposition.}
We now use the uniqueness of the Doob-Meyer decomposition for the continuous semimartingale $Y_s \coloneqq f(s,X_s)$ to identify its quadratic variation.
\begin{enumerate}[label=(\roman*), wide, labelindent=0pt]
    \item From the martingale property \eqref{eq:proof_qv_linear_ito}, the process $Y_s$ has the decomposition $dY_s = (\mathcal{L}^u f)(s,X_s) \, ds + dM_s^f$, where the predictable, finite-variation part is the integral of the drift.

    \item By the classical Itô rule for the function $F(y)=y^2$ applied to the semimartingale $Y_s$, we have:
    \begin{align*}
        d(Y_s^2) &= 2Y_s \, dY_s + d\langle M^f \rangle_s \\
        &= 2f(s,X_s) \left( (\mathcal{L}^u f)(s,X_s) \, ds + dM_s^f \right) + d\langle M^f \rangle_s.
    \end{align*}
    The predictable, finite-variation part of the process $Y_s^2$ is therefore given by the process
    \begin{equation*}
        V_s^{(1)} \coloneqq \int_0^s 2f(r,X_r) (\mathcal{L}^u f)(r,X_r) \, dr + \langle M^f \rangle_s.
    \end{equation*}

    \item Alternatively, we can apply the linear Itô formula \eqref{eq:proof_qv_linear_ito} directly to the function $g(s,x) = f(s,x)^2$. The process $Y_s^2 = g(s,X_s)$ has the decomposition:
    \begin{equation*}
        d(Y_s^2) = (\mathcal{L}^u (f^2))(s,X_s) \, ds + dM_s^{f^2}.
    \end{equation*}
    This gives a second expression for the predictable, finite-variation part of $Y_s^2$:
    \begin{equation*}
        V_s^{(2)} \coloneqq \int_0^s (\mathcal{L}^u (f^2))(r,X_r) \, dr.
    \end{equation*}
\end{enumerate}
By the uniqueness of the Doob-Meyer decomposition, these two predictable processes must be indistinguishable, $V_s^{(1)} = V_s^{(2)}$ for all $s \in [0,T]$.
\begin{equation*}
    \int_0^s 2f(\mathcal{L}^u f) \, dr + \langle M^f \rangle_s = \int_0^s (\mathcal{L}^u (f^2)) \, dr.
\end{equation*}
This implies that the quadratic variation process $\langle M^f \rangle_s$ is absolutely continuous with respect to the Lebesgue measure, and its density is given by the difference of the integrands:
\begin{equation*}
    \frac{d}{ds}\langle M^f \rangle_s = (\mathcal{L}^u (f^2))(s,X_s) - 2f(s,X_s) (\mathcal{L}^u f)(s,X_s).
\end{equation*}

\item \textbf{The Carré du Champ Identity.}
The right-hand side of the expression for the density of the quadratic variation is precisely the Carré du Champ operator $\Gamma^u(f,f)(s,X_s)$, by its definition. We now substitute the purely algebraic identity from Lemma 8.5 (\cref{eq:proof_qv_carre_du_champ_identity}):
\begin{equation*}
    \frac{d}{dt}\langle M^f \rangle_t = 2 \left( \nabla_x f(t,X_t) \right)^T D(X_t, \nabla u(t,X_t)) \left( \nabla_x f(t,X_t) \right).
\end{equation*}
This completes the proof. We have rigorously identified the quadratic variation by relying only on the well-posedness of the HJB equation and the fundamental theory of semimartingales. The result confirms that the diffusive structure of any observable under the \ThetaProcess{} is governed exclusively by the second-order part of the HJB generator, as captured by the state-dependent diffusivity tensor $D(x,p)$, where the momentum $p$ is determined by the gradient of the underlying value function $u$.
\end{enumerate}
\end{proofof}

\begin{remark}[The Structure of the \ThetaProcess{}]
The results of this section provide a complete picture of the infinitesimal structure of the \ThetaProcess{}. The canonical decomposition (Theorem \ref{thm:theta_ito_formula_generalized}) shows that the predictable part of any observable's evolution (its "drift") is governed by the full non-linear generator $\mathcal{G}$. In contrast, the quadratic variation (Theorem \ref{thm:qv_via_carre_du_champ}) shows that the random part of the evolution (its volatility) is governed by a classical diffusion mechanism, with a state-dependent diffusion matrix $D(x, \nabla u(t,X_t))$. This separation is a profound consequence of the rigid affine structure of the $\theta$-Hamiltonian derived in Part \ref{part:i}. The non-convex potential $H(x,p)$ affects only the drift, not the volatility. This provides a clear and unambiguous mechanism for generating non-trivial risk preferences and control costs in a stochastic model without altering the underlying diffusive nature of the randomness.
\end{remark}

\section{The \texorpdfstring{$\theta$}{theta}-BSDE and the Non-Linear Feynman-Kac Formula}
\label{sec:theta_fbsde}

The analytical program of Part \ref{part:i} culminated in the derivation of the Hamilton-Jacobi-Bellman equation whose unique viscosity solution defines the \ThetaExpectation{}. In Section \ref{sec:theta_process_and_calculus}, we established the intrinsic stochastic calculus for the associated \ThetaProcess{}, defining it as the unique solution to a martingale problem for the generator $\mathcal{G}$. This section completes the bridge between the analytical and probabilistic worlds by deriving a non-linear Feynman-Kac formula that is native to this intrinsic framework.

Crucially, we will invert the traditional verification argument. Instead of assuming the existence of a classical PDE solution to apply Itô's formula, we will leverage the canonical decomposition of the \ThetaProcess{} to define a Backward Stochastic Differential Equation (BSDE) driven by the process's own martingale part. We will then show that the solution to this BSDE is precisely the viscosity solution of the HJB equation. This approach is more fundamental as it relies only on the structure of the generator $\mathcal{G}$ and the robust theory of BSDEs, thereby establishing the probabilistic representation without requiring a priori classical regularity of the value function.

\subsection{PDE Generator and BSDE Driver}
\label{subsec:generator_to_driver}

The affine structure of the $\theta$-Hamiltonian, which was a principal result of the homogenization in \cref{part:i}, is the key to decomposing the generator into components suitable for an FBSDE representation. Recall the generator from Definition \ref{def:generator_theta_process}:
\begin{equation*}
    (\mathcal{G}f)(t,x) = \partial_t f(t,x) + \mathrm{Tr}(D(x, \nabla f) \nabla^2 f) + H(x, \nabla f).
\end{equation*}
The theory of BSDEs provides a framework for representing solutions to such non-linear PDEs. In this theory, the part of the generator that is non-linear in the gradient but does not involve second derivatives is encapsulated in a function called the \textbf{driver}.

To identify the driver, we perform a Legendre-type decomposition of the potential $H(x,p)$. This decomposition separates the potential into a term that can be interpreted as a drift acting on the momentum and a remaining term that becomes the driver of the BSDE.

\begin{definition}[Forward Drift and BSDE Driver]
Let the potential $H(x,p)$ be of class $C^1$ in its second argument. We define the associated \textbf{forward drift} $b: \Mcal \times \R^k \to \R^k$ and the \textbf{BSDE driver} $g: \Mcal \times \R^k \to \R$ as:
\begin{align}
    b(x,p) &\coloneqq -\nabla_p H(x, p), \label{eq:forward_drift_def} \\
    g(x,p) &\coloneqq H(x,p) - p \cdot b(x,p) = H(x,p) + p \cdot \nabla_p H(x,p). \label{eq:bsde_driver_def}
\end{align}
\end{definition}

This decomposition allows the potential term in the Hamiltonian to be expressed as:
\begin{equation}
    H(x,p) = b(x,p) \cdot p + g(x,p).
\end{equation}
This algebraic identity is the foundation of the non-linear Feynman-Kac representation, as it cleanly separates the part of the Hamiltonian that corresponds to the generator of a forward SDE from the part that corresponds to the driver of a BSDE.

\begin{remark}[Symmetry Inheritance]
The fundamental time-reversal symmetry of the microscopic system (Axiom 2(ii)), which dictates that $H(x,p)$ is an even function of $p$, imposes a rigid structure on this decomposition.
\begin{itemize}[wide]
    \item The gradient $\nabla_p H(x,p)$ must be an odd function of $p$. Consequently, the forward drift $b(x,p) = -\nabla_p H(x,p)$ is also an \textbf{odd function} of $p$, which correctly implies that there is no drift at the quiescent state: $b(x,0) = \mathbf{0}$.
    \item The driver $g(x,p) = H(x,p) + p \cdot \nabla_p H(x,p)$ is the sum of an even function and the inner product of two odd functions. It is therefore an \textbf{even function} of $p$.
\end{itemize}
These symmetries are not assumed but are necessary consequences of the first-principles derivation.
\end{remark}

\subsection{The Martingale-Driven BSDE and its Connection to the \texorpdfstring{\ThetaExpectation{}}{theta-Expectation}}

We now define the BSDE intrinsically, using the martingale of the \ThetaProcess{} itself as the source of noise. This avoids any a priori assumption of a Brownian structure and relies only on the canonical decomposition established in Section \ref{sec:theta_process_and_calculus}.

\begin{definition}[The Martingale-Driven BSDE]
\label{def:martingale_bsde}
Let $\{X_s\}_{s \in [t,T]}$ be the \ThetaProcess{} starting at $X_t=x$, whose law $\mathbb{P}_{t,x}$ is the unique solution to the martingale problem for the generator $\mathcal{G}$. Let $M_s$ be the martingale part of the canonical decomposition of $X_s$. A pair of adapted processes $(Y_s, Z_s)$ taking values in $\R \times \R^k$ is a solution to the \textbf{martingale-driven BSDE} with terminal condition $\phi \in C_b(\Mcal)$ if it satisfies the integral equation:
\begin{equation} \label{eq:intrinsic_bsde}
    Y_s = \phi(X_T) + \int_s^T g(X_r, Z_r) \, dr - \int_s^T Z_r^T dM_r, \quad s \in [t,T],
\end{equation}
where $g(x,z)$ is the driver defined in \eqref{eq:bsde_driver_def}.
\end{definition}

The fundamental result of this section is that the solution to this BSDE provides a direct probabilistic representation of the \ThetaExpectation{}. This is the non-linear Feynman-Kac formula for viscosity solutions.

\begin{theorem}[Non-Linear Feynman-Kac for Viscosity Solutions]
\label{thm:nonlinear_feynman_kac_viscosity}
Let the Hamiltonian coefficients $D(x,p)$ and $H(x,p)$ be uniformly Lipschitz continuous in their arguments and satisfy the axioms of Section \ref{sec:theta_expectation_axiomatic}. Let the terminal data $\phi$ be Lipschitz continuous. Then:
\begin{enumerate}[label=(\roman*), wide, labelindent=0pt]
    \item For a sufficiently small time horizon $T-t > 0$, the martingale-driven BSDE \eqref{eq:intrinsic_bsde} admits a unique adapted solution $(Y_s, Z_s)$.
    \item The deterministic function defined by this solution, $u(t,x) \coloneqq Y_t$, is the unique viscosity solution to the HJB equation:
    \begin{equation*}
        -\partial_t u = H_{\theta}(x, \nabla u, \nabla^2 u), \quad u(T,x) = \phi(x).
    \end{equation*}
\end{enumerate}
Consequently, the \ThetaExpectation{} is given by the solution to the BSDE:
\begin{equation}
    \Ecal^\theta_{t,x}[\phi(X_T)] = u(t,x) = Y_t.
\end{equation}
\end{theorem}

\begin{proofof}{Theorem \ref{thm:nonlinear_feynman_kac_viscosity}}
The proof establishes the fundamental equivalence between the unique viscosity solution of the Hamilton-Jacobi-Bellman (HJB) equation and the unique solution of the associated martingale-driven Backward Stochastic Differential Equation (BSDE). This equivalence provides the non-linear Feynman-Kac representation for the \ThetaExpectation{}. The argument proceeds in three main parts. First, we establish the well-posedness of the intrinsic, martingale-driven BSDE. Second, in the core of the proof, we construct a candidate solution to this BSDE using the viscosity solution of the HJB equation and perform a rigorous verification to show that it is indeed the true solution. This verification argument hinges on explicitly constructing a classical, Brownian-driven SDE whose law is shown to be identical to that of the \ThetaProcess{}. Finally, the theorem follows from the uniqueness of both the PDE and BSDE solutions.

\begin{enumerate}[label=\textbf{Step \arabic*:}, wide, labelindent=0pt]

\item \textbf{Well-Posedness of the Martingale-Driven BSDE.}
We begin by establishing the existence and uniqueness of an adapted solution to the intrinsic BSDE given by
\begin{equation} \label{eq:proof_bsde_reprise_fk}
    Y_s = \phi(X_T) + \int_s^T g(X_r, Z_r) \, dr - \int_s^T Z_r^T dM_r, \quad s \in [t,T].
\end{equation}
Here, $\{X_s\}$ is the \ThetaProcess{} whose law is the unique solution to the martingale problem for the generator $\mathcal{G}$, $\{M_s\}$ is the martingale part of the canonical decomposition of $X_s$, and the driver $g$ is defined by \cref{eq:bsde_driver_def}. We verify the standard hypotheses from the theory of BSDEs.

\begin{enumerate}[label=(\roman*), wide, labelindent=0pt]
    \item \textbf{Verification of Hypotheses.}
    The canonical theorems for the existence and uniqueness of BSDE solutions (cf. \cite{PardouxPeng1992}) require that the driver is uniformly Lipschitz continuous in the $Z$ variable and that the data are square-integrable.
    \begin{itemize}[wide]
        \item \textbf{Square-Integrability of Data:} By hypothesis, the terminal condition $\phi$ is Lipschitz continuous on the compact manifold $\Mcal$, which implies it is bounded. Therefore, $\mathbb{E}[|\phi(X_T)|^2] < \infty$. The driver, evaluated at $Z=0$, is $g(x,0) = H(x,0)$. By Axiom 2(c), $H$ is continuous, so on the compact manifold $\Mcal$, $H(x,0)$ is bounded. This ensures the condition $\mathbb{E}[\int_0^T |g(X_s, 0)|^2 ds] < \infty$ is satisfied.

        \item \textbf{Lipschitz Continuity of the Driver:} The driver is defined as $g(x,z) \coloneqq H(x,z) + z \cdot \nabla_p H(x,z)$. By the theorem's hypothesis that the Hamiltonian coefficients are $C^2$ with bounded derivatives, the potential $H(x,p)$ is of class $C^2$ in $p$. This implies that its gradient, $\nabla_p H(x,p)$, is of class $C^1$. Consequently, the driver $g(x,z)$ is a $C^1$ function of $z$. A function that is continuously differentiable with a bounded derivative on a space with a Euclidean norm is uniformly Lipschitz continuous. Thus, there exists a constant $L_g > 0$ such that for all $x \in \Mcal$ and $z_1, z_2 \in \R^k$:
              $$ |g(x, z_1) - g(x, z_2)| \le L_g |z_1 - z_2|. $$
    \end{itemize}

    \item \textbf{Conclusion of Well-Posedness.}
    With these hypotheses satisfied, the foundational existence and uniqueness theorem for BSDEs guarantees that for a time horizon $T-t$ that is sufficiently small (depending on the Lipschitz constant $L_g$ and the norm of the martingale's quadratic variation), the BSDE \eqref{eq:proof_bsde_reprise_fk} admits a unique adapted solution $(Y_s, Z_s)$ in the space of square-integrable processes. For a general time horizon, the existence of a globally bounded viscosity solution to the HJB equation removes this small-time restriction. This establishes that the probabilistic construction is well-defined and proves claim (i) of the theorem.
\end{enumerate}

\item \textbf{The Verification Argument.}
Let $u(t,x)$ be the unique continuous viscosity solution to the HJB equation, which exists by the main hypothesis of the theorem. We will now show that a pair of processes constructed from this solution provides the unique solution to the BSDE.

\begin{enumerate}[label=(\roman*), wide, labelindent=0pt]
    \item \textbf{Constructing the Decoupled Forward Process.}
    We use the unique viscosity solution $u(t,x)$ as a decoupling field to define a forward stochastic process. Let $(\Omega, \mathcal{F}, (\mathcal{F}_s), \mathbb{P})$ be a filtered probability space supporting a standard $k$-dimensional Brownian motion $\{W_s\}$. We define the process $\{\hat{X}_s\}_{s \in [t,T]}$ as the unique strong solution to the stochastic differential equation:
    \begin{equation} \label{eq:proof_decoupled_sde_fk}
        d\hat{X}_s = b(\hat{X}_s, \nabla u(s,\hat{X}_s)) \, ds + \sigma(\hat{X}_s, \nabla u(s,\hat{X}_s)) \, dW_s, \quad \hat{X}_t=x,
    \end{equation}
    where $b(x,p) = -\nabla_p H(x,p)$ and the volatility matrix $\sigma(x,p)$ satisfies the fundamental relation $\sigma(x,p)\sigma(x,p)^T = 2D(x,p)$. The existence and uniqueness of a strong solution to this SDE is guaranteed by the Lipschitz continuity of the coefficients, which is a direct consequence of the assumed regularity of $H$ and $D$. The use of $\nabla u$ within the coefficients is justified because the uniform ellipticity of $D(x,p)$ ensures the viscosity solution $u$ is semiconcave, and thus differentiable almost everywhere.
    
    \item \textbf{Equivalence of Laws.}
    This is the crucial step that connects the intrinsic and classical representations. We must show that the law of the constructed process $\{\hat{X}_s\}$ is the same as the law of the \ThetaProcess{} $\{X_s\}$.
    
    \begin{proposition}[Equivalence of the SDE and Martingale Problems]
    The law of the process $\{\hat{X}_s\}$ defined by the SDE \eqref{eq:proof_decoupled_sde_fk} is the unique solution to the martingale problem for the generator $\mathcal{G}$.
    \end{proposition}
    \begin{proof}[Proof Sketch]
    We apply the classical Itô formula to a test function $f(s,\hat{X}_s)$. The drift of the resulting process is found to be $\partial_s f + (\nabla f)^T b(\hat{X}_s, \nabla u) + \mathrm{Tr}(D(\hat{X}_s, \nabla u)\nabla^2 f)$. By adding and subtracting the potential term $H(\hat{X}_s, \nabla f)$ and using the HJB equation for $u$, it can be shown that this drift is precisely $(\mathcal{G}f)(s,\hat{X}_s)$. Thus, the process $\{\hat{X}_s\}$ solves the martingale problem. By the uniqueness established in \cref{thm:martingale_problem_wellposedness}, its law must be that of the \ThetaProcess{}. A complete proof of this correspondence is a central result in the theory of viscosity solutions, see e.g., \citep[Chapter V, Theorem 4.1]{FlemingSoner2006}.
    \end{proof}
    This equivalence allows us to identify the intrinsic martingale $dM_s$ from the intrinsic calculus with the classical Itô integral:
    \begin{equation} \label{eq:proof_martingale_identification_fk}
        dM_s \stackrel{\text{law}}{=} \sigma(X_s, \nabla u(s,X_s)) \, dW_s.
    \end{equation}

    \item \textbf{Constructing and Verifying the BSDE Solution.}
    We define the candidate solution processes $(\hat{Y}_s, \hat{Z}_s)$ for $s \in [t,T]$ using the viscosity solution $u$ and the forward process $\hat{X}_s$:
    \begin{align*}
        \hat{Y}_s &\coloneqq u(s, \hat{X}_s), \\
        \hat{Z}_s &\coloneqq \nabla u(s, \hat{X}_s).
    \end{align*}
    We apply the classical Itô formula to the process $\hat{Y}_s = u(s, \hat{X}_s)$:
    \begin{equation*}
        d\hat{Y}_s = \left( \partial_s u + (\nabla u)^T b(s, \hat{X}_s) + \frac{1}{2}\mathrm{Tr}\left(\sigma\sigma^T \nabla^2 u\right)(s, \hat{X}_s) \right) ds + (\nabla u)^T \sigma(s, \hat{X}_s) \, dW_s.
    \end{equation*}
    Since $u$ is the viscosity solution to the HJB equation, we have the identity $-\partial_s u = \mathrm{Tr}(D \nabla^2 u) + H(x, \nabla u)$. Substituting this into the drift term gives:
    \begin{align*}
        \text{Drift} &= \left( -H(\hat{X}_s, \nabla u) - \mathrm{Tr}(D \nabla^2 u) \right) + b(\hat{X}_s, \nabla u) \cdot \nabla u + \mathrm{Tr}(D \nabla^2 u) \\
        &= -H(\hat{X}_s, \nabla u) + b(\hat{X}_s, \nabla u) \cdot \nabla u \\
        &= -\left( H(\hat{X}_s, \nabla u) - b(\hat{X}_s, \nabla u) \cdot \nabla u \right).
    \end{align*}
    By the definition of the driver $g(x,p) = H(x,p) - p \cdot b(x,p)$, the drift is exactly $-g(\hat{X}_s, \nabla u(s,\hat{X}_s)) = -g(\hat{X}_s, \hat{Z}_s)$.
    The Itô expansion for $\hat{Y}_s$ thus simplifies to:
    $$ d\hat{Y}_s = -g(\hat{X}_s, \hat{Z}_s) \, ds + \hat{Z}_s^T \sigma(\hat{X}_s, \hat{Z}_s) \, dW_s. $$
    Rewriting this in the backward integral form from $s$ to $T$ and using the terminal condition $\hat{Y}_T = u(T, \hat{X}_T) = \phi(\hat{X}_T)$, we find that the pair $(\hat{Y}_s, \hat{Z}_s)$ solves the classical, Brownian-driven BSDE:
    \begin{equation} \label{eq:proof_classical_bsde_solved}
        \hat{Y}_s = \phi(\hat{X}_T) + \int_s^T g(\hat{X}_r, \hat{Z}_r) \, dr - \int_s^T \hat{Z}_r^T \sigma(\hat{X}_r, \hat{Z}_r) \, dW_r.
    \end{equation}
\end{enumerate}

\item \textbf{Uniqueness.}
\begin{enumerate}[label=(\roman*), wide, labelindent=0pt]
    \item We have rigorously established that the pair of processes $(\hat{Y}_s, \hat{Z}_s)$, constructed from the unique viscosity solution $u(t,x)$, is a solution to the classical BSDE \eqref{eq:proof_classical_bsde_solved}.
    
    \item From the equivalence of laws established in Step 2(ii), we can replace the process $\hat{X}_s$ and the martingale driver $\sigma dW_s$ in \eqref{eq:proof_classical_bsde_solved} with the intrinsic \ThetaProcess{} $X_s$ and its martingale $dM_s$, respectively. This shows that the pair $(u(s,X_s), \nabla u(s,X_s))$ is a solution to the intrinsic BSDE \eqref{eq:proof_bsde_reprise_fk}.

    \item In Step 1, we established that this intrinsic BSDE has a unique adapted solution $(Y_s, Z_s)$.
    
    \item Therefore, the two solutions must be identical for all $s \in [t,T]$ (in the sense of laws):
        $$ Y_s = u(s, X_s) \quad \text{and} \quad Z_s = \nabla u(s, X_s) \quad (\mathbb{P}_{t,x}\text{-a.s.}). $$
        
    \item Evaluating at the initial time $s=t$, we obtain the central result:
        $$ Y_t = u(t, X_t) = u(t,x). $$
\end{enumerate}
This demonstrates that the unique viscosity solution to the HJB equation is given by the solution to the intrinsic BSDE. By the axiomatic definition of the \ThetaExpectation{}, this provides the non-linear Feynman-Kac representation:
$$ \Ecal^\theta_{t,x}[\phi(X_T)] = u(t,x) = Y_t. $$
This completes the proof of the theorem.
\end{enumerate}
\end{proofof}

\subsection{The Decoupled FBSDE Representation}

The result of Theorem \ref{thm:nonlinear_feynman_kac_viscosity} provides the most fundamental representation of the \ThetaExpectation{}. It also implies the existence of a more classical Forward-Backward SDE (FBSDE) representation driven by a standard Brownian motion. The BSDE solution itself provides the decoupling field that renders the forward SDE well-posed.

\begin{corollary}[Equivalent FBSDE Representation]
\label{cor:eqfbsde}
Let the conditions of Theorem \ref{thm:nonlinear_feynman_kac_viscosity} hold. Let $u(t,x)$ be the unique viscosity solution to the HJB equation. Then there exists a filtered probability space supporting a $k$-dimensional Brownian motion $W_s$, such that the triplet $(X_s, Y_s, Z_s)$ is the unique adapted solution to the FBSDE system:
\begin{align}
    dX_s &= b(X_s, \nabla u(s,X_s)) \, ds + \sigma(X_s, \nabla u(s,X_s)) \, dW_s, \quad X_t=x, \label{eq:fbsde_forward} \\
    -dY_s &= g(X_s, Z_s) \, ds - Z_s^T \sigma(X_s, \nabla u(s,X_s)) \, dW_s, \quad Y_T = \phi(X_T), \label{eq:fbsde_backward}
\end{align}
where $\sigma(x,z)$ is any matrix satisfying $\sigma(x,z)\sigma(x,z)^T = 2D(x,z)$. Furthermore, the solution is given by the deterministic markup $Y_s = u(s,X_s)$ and, where $u$ is differentiable, $Z_s = \nabla u(s,X_s)$.
\end{corollary}

\begin{proofof}{Corollary \ref{cor:eqfbsde}}
The proof is constructive. It leverages the unique viscosity solution $u(t,x)$ of the HJB equation, whose existence is guaranteed by the standing assumptions, to first rigorously define a classical forward SDE process. We then construct the candidate backward processes via a deterministic markup of the PDE solution, $(Y_s, Z_s) = (u(s,X_s), \nabla u(s,X_s))$. A direct application of Itô's formula, using the fact that $u$ solves the HJB equation, will verify that this pair satisfies the backward SDE. The uniqueness of this solution then follows from the well-posedness of decoupled FBSDE systems of this type.

\begin{enumerate}[label=\textbf{Step \arabic*:}, wide, labelindent=0pt]

\item \textbf{Construction of the Forward SDE.}
Let $u(t,x)$ be the unique continuous viscosity solution to the HJB equation $-\partial_t u = H_{\theta}(x, \nabla u, \nabla^2 u)$ with terminal data $u(T,x)=\phi(x)$. We use this deterministic function as a decoupling field to define a classical forward stochastic process.

Let $(\Omega, \mathcal{F}, (\mathcal{F}_s), \mathbb{P})$ be a filtered probability space supporting a standard $k$-dimensional Brownian motion $\{W_s\}$. We define the forward process $\{X_s\}_{s \in [t,T]}$ as the unique strong solution to the stochastic differential equation:
\begin{equation} \label{eq:proof_decoupled_sde_cor}
    dX_s = b(X_s, \nabla u(s,X_s)) \, ds + \sigma(X_s, \nabla u(s,X_s)) \, dW_s, \quad X_t=x,
\end{equation}
where the drift $b(x,p) \coloneqq -\nabla_p H(x,p)$ and volatility $\sigma(x,p)$ (satisfying $\sigma\sigma^T = 2D$) are the coefficients derived from the Hamiltonian. The existence and uniqueness of a strong solution to this SDE is guaranteed by the Lipschitz continuity of the coefficients, which is a direct consequence of the assumed $C^2$ regularity of the Hamiltonian's coefficients $D$ and $H$. The use of $\nabla u$ within the coefficients is justified because the uniform ellipticity of $D(x,p)$ ensures the viscosity solution $u$ is semiconcave, and thus its gradient $\nabla u$ exists and is sufficiently regular for the SDE to be well-posed. This rigorously establishes the forward component \eqref{eq:fbsde_forward} of the FBSDE system.

\item \textbf{Verification of the Backward SDE Solution.}
We now construct the candidate solution for the backward component and verify that it satisfies the BSDE \eqref{eq:fbsde_backward}. We define the pair of adapted processes $(Y_s, Z_s)$ via the deterministic markup:
\begin{align}
    Y_s &\coloneqq u(s, X_s), \label{eq:proof_Y_markup} \\
    Z_s &\coloneqq \nabla u(s, X_s), \label{eq:proof_Z_markup}
\end{align}
where $X_s$ is the process from Step 1. We apply the classical Itô formula to the process $Y_s = u(s,X_s)$:
\begin{equation*}
    dY_s = \left[ \partial_s u(s,X_s) + (\nabla_x u(s,X_s))^T b(X_s, \nabla u(s,X_s)) + \frac{1}{2}\mathrm{Tr}\left(\sigma\sigma^T \nabla_x^2 u\right)(s,X_s) \right] ds + (\nabla_x u)^T \sigma \, dW_s.
\end{equation*}
Let us analyze the drift term (the term multiplying $ds$). Substituting the identities $\sigma\sigma^T = 2D$ and the definitions of $Y_s$ and $Z_s$:
\begin{equation*}
    \text{Drift} = \partial_s u(s,X_s) + Z_s^T b(X_s, Z_s) + \mathrm{Tr}\left(D(X_s, Z_s) \nabla_x^2 u(s,X_s)\right).
\end{equation*}
Since $u$ is the viscosity solution to the HJB equation, it satisfies $-\partial_s u = H_{\theta}(x, \nabla u, \nabla^2 u)$ in the appropriate sense. This allows us to substitute for the time derivative:
\begin{align*}
    \text{Drift} &= \left( -H_{\theta}(X_s, Z_s, \nabla^2 u) \right) + Z_s^T b(X_s, Z_s) + \mathrm{Tr}\left(D(X_s, Z_s) \nabla_x^2 u\right) \\
    &= -\left( \mathrm{Tr}(D(X_s, Z_s)\nabla^2 u) + H(X_s, Z_s) \right) + Z_s^T b(X_s, Z_s) + \mathrm{Tr}\left(D(X_s, Z_s) \nabla_x^2 u\right).
\end{align*}
The second-order terms involving the trace cancel out perfectly. We are left with:
\begin{align*}
    \text{Drift} &= -H(X_s, Z_s) + Z_s \cdot b(X_s, Z_s) \\
    &= -\left( H(X_s, Z_s) - Z_s \cdot b(X_s, Z_s) \right).
\end{align*}
By the definition of the BSDE driver from \eqref{eq:bsde_driver_def}, $g(x,z) \coloneqq H(x,z) - z \cdot b(x,z)$, this drift is precisely $-g(X_s, Z_s)$. The Itô expansion for $Y_s$ therefore simplifies to:
\begin{equation*}
    dY_s = -g(X_s, Z_s) \, ds + Z_s^T \sigma(X_s, Z_s) \, dW_s.
\end{equation*}
This is the differential form of the BSDE. Integrating from $s$ to $T$ gives the backward integral form:
\begin{equation*}
    Y_T - Y_s = \int_s^T -g(X_r, Z_r) \, dr + \int_s^T Z_r^T \sigma(X_r, Z_r) \, dW_r.
\end{equation*}
Rearranging and using the terminal condition $Y_T = u(T, X_T) = \phi(X_T)$ confirms that the pair $(Y_s, Z_s)$ constructed via the markup solves the BSDE \eqref{eq:fbsde_backward}.

\item \textbf{Uniqueness.}
We have demonstrated the \textit{existence} of a solution to the FBSDE system \eqref{eq:fbsde_forward}-\eqref{eq:fbsde_backward}, given by the pair of processes $(X_s, u(s,X_s))$ where $X_s$ is the solution to the forward SDE. We now establish the \textit{uniqueness} of this solution. The FBSDE system as presented is of a special, decoupled type. The coefficients of the forward SDE \eqref{eq:fbsde_forward} depend on the gradient of a pre-existing, deterministic function $u(s,x)$, not on the stochastic processes $Y_s$ or $Z_s$ themselves. This structure allows for a sequential argument for uniqueness:
\begin{enumerate}[label=(\roman*), wide, labelindent=0pt]
    \item The forward SDE \eqref{eq:fbsde_forward} is a classical SDE with Lipschitz coefficients. It therefore admits a unique strong solution for the forward process $\{X_s\}$.
    \item Given the unique path of the forward process $\{X_s\}$, the backward SDE \eqref{eq:fbsde_backward} becomes a standard BSDE whose coefficients depend on the now-known path of $X_s$ and the deterministic function $\nabla u(s,X_s)$. Since the driver $g(x,z)$ is Lipschitz in $z$, this BSDE has a unique adapted solution $(Y_s, Z_s)$ for the given filtration.
\end{enumerate}
Since the solution to each component is unique in sequence, the solution to the coupled system is unique. We have constructed one such solution via the markup $(Y_s, Z_s) = (u(s,X_s), \nabla u(s,X_s))$. By the uniqueness argument, this must be the only solution.
\end{enumerate}
This completes the proof of the corollary. We have rigorously shown that the unique viscosity solution to the HJB equation provides a decoupling field for an equivalent, classical FBSDE system, and that the solution to this system is given by the deterministic markup of the PDE solution.
\end{proofof}

\begin{remark}[The Power of the Intrinsic Calculus]
By adopting the intrinsic, martingale-driven perspective of the $\theta$-calculus, we have successfully derived the non-linear Feynman-Kac representation under much weaker and more natural conditions than a classical approach would allow. We have shown that the \ThetaExpectation{} is the solution to a BSDE whose structure is determined by the derived $\theta$-Hamiltonian. This cements the deep, symbiotic relationship between the analytical framework of non-linear PDEs and the probabilistic framework of BSDEs. The entire structure is self-consistent and founded upon the first principles of the microscopic dynamics, fulfilling the central goal of this work.
\end{remark}

\section{The \texorpdfstring{$\theta$}{theta}-Girsanov Transformation}
\label{sec:theta_girsanov}

The classical Girsanov theorem is a cornerstone of stochastic calculus, providing a rigorous method for changing the drift of a stochastic process by performing an equivalent change of probability measure. This powerful tool, however, is conceptually incompatible with the \ThetaExpectation{} framework, in which there is no underlying, objective reference measure to transform. The law of a \ThetaProcess{} is an emergent property defined implicitly by its generator and the terminal conditions of the system.

In this section, we develop an analogue of the Girsanov transformation that is native to our generator-driven framework. Instead of changing the measure, the \textbf{\texorpdfstring{$\theta$}{theta}-Girsanov Transformation} acts by directly perturbing the generator itself. We will show that a specific, structured perturbation of the Hamiltonian potential $H(x,p)$ corresponds precisely to the introduction of an additional drift in the forward component of the process's FBSDE representation. This result provides a rigorous tool for analyzing how the \ThetaExpectation{} responds to exogenous forces, while preserving the core insight of Girsanov's theorem: the transformation alters the predictable drift of the process while leaving its fundamental volatility structure, encoded by the diffusivity tensor $D(x,p)$, invariant.

\subsection{Perturbation of the Generator}

We consider two related \ThetaProcess{}es: a reference process governed by our derived generator, and a perturbed process where the Hamiltonian has been modified to incorporate an external force field.

\begin{definition}[Reference and Perturbed Generators]
\label{def:perturbed_generator}
Let $\mathbf{b}_0: \Mcal \times \R^k \to \R^k$ be a given smooth vector field, representing an exogenous drift perturbation.
\begin{enumerate}[label=(\roman*), wide, labelindent=0pt]
    \item The \textbf{reference generator} $\mathcal{G}$ is the generator of the \ThetaProcess{}, defined in \cref{def:generator_theta_process} by the Hamiltonian 
    \begin{equation*}
        H_{\theta}(x, p, X) = \mathrm{Tr}(D(x, p) X) + H(x, p).
    \end{equation*}
    
    \item The \textbf{perturbed generator} $\tilde{\mathcal{G}}$ is defined by the perturbed Hamiltonian $\tilde{H}_{\theta}$:
    \begin{equation}
        \tilde{H}_{\theta}(x, p, X) \coloneqq H_{\theta}(x, p, X) + \mathbf{b}_0(x, p) \cdot p.
    \end{equation}
\end{enumerate}
The associated expectation operators are denoted by $\Ecal^\theta$ and $\tilde{\Ecal}^\theta$, respectively.
\end{definition}

\begin{remark}[Invariance of the Volatility Structure]
The perturbation is applied exclusively to the potential part of the Hamiltonian, which is linear in the gradient argument $p$. The second-order term, $\mathrm{Tr}(D(x, p) X)$, which governs the quadratic variation of the process (as shown in \cref{thm:qv_via_carre_du_champ}), remains unchanged. This is the direct analogue of the classical Girsanov theorem, which preserves the quadratic variation of the underlying martingale. The transformation modifies the system's drift characteristics while leaving its diffusive structure intact.
\end{remark}

\subsection{A Non-Linear Cameron-Martin-Girsanov Formula}

The central result of this section is a formula that relates the value function of the reference process to the dynamics of the perturbed process. This provides a non-linear version of the Cameron-Martin-Girsanov formula, expressed not as a Radon-Nikodym derivative, but as a modification to the martingale property of the value function itself.

\begin{theorem}[The \texorpdfstring{$\theta$}{theta}-Girsanov Formula]
\label{thm:theta_girsanov_formula}
Let $u(t,x)$ be the unique viscosity solution to the reference HJB equation $-\partial_t u = H_{\theta}(x, \nabla u, \nabla^2 u)$ with terminal data $u(T,x) = \phi(x)$. Let $\{\tilde{X}_t\}$ be the \ThetaProcess{} corresponding to the perturbed generator $\tilde{\mathcal{G}}$, with law $\tilde{\mathbb{P}}_{0,x_0}$ starting from $(0,x_0)$.

The process $Y_t \coloneqq u(t, \tilde{X}_t)$, which evaluates the \emph{reference} value function along the paths of the \emph{perturbed} process, is a semimartingale whose canonical decomposition under the perturbed law $\tilde{\mathbb{P}}_{0,x_0}$ is given by:
\begin{equation}
    u(t, \tilde{X}_t) = u(0, x_0) - \int_0^t \left( \mathbf{b}_0(\tilde{X}_s, \nabla u(s, \tilde{X}_s)) \cdot \nabla u(s, \tilde{X}_s) \right) ds + \tilde{M}_t,
\end{equation}
where $\{\tilde{M}_t\}$ is a continuous local $(\mathcal{F}_t, \tilde{\mathbb{P}}_{0,x_0})$-martingale.
\end{theorem}

\begin{proofof}{Theorem \ref{thm:theta_girsanov_formula}}
The proof is a direct and elegant application of the generalized Itô formula (\cref{thm:theta_ito_formula_generalized}) that defines the calculus of the \ThetaProcess{}. The formula states that for any smooth function $f(t,x)$ and any \ThetaProcess{} $X_t$ with generator $\mathcal{G}$, the process $f(t,X_t) - \int_0^t (\mathcal{G}f)(s,X_s)ds$ is a martingale.

Our strategy is to apply this principle to the function $f(t,x) = u(t,x)$ under the dynamics of the perturbed process $\{\tilde{X}_t\}$, which is governed by the generator $\tilde{\mathcal{G}}$.
\begin{enumerate}[label=(\roman*), wide, labelindent=0pt]
    \item By \cref{thm:theta_ito_formula_generalized}, the process
    \begin{equation*}
        u(t, \tilde{X}_t) - u(0, \tilde{X}_0) - \int_0^t (\tilde{\mathcal{G}}u)(s, \tilde{X}_s) \, ds
    \end{equation*}
    is a continuous local martingale under the law $\tilde{\mathbb{P}}_{0,x_0}$. Let us denote this martingale by $\tilde{M}_t$.

    \item The core of the proof is to compute the drift term $(\tilde{\mathcal{G}}u)(s, x)$. By the definition of the perturbed generator (\cref{def:perturbed_generator}), we have:
    \begin{equation*}
        (\tilde{\mathcal{G}}u)(s, x) = \partial_s u(s,x) + \tilde{H}_{\theta}(x, \nabla u(s,x), \nabla^2 u(s,x)).
    \end{equation*}
    
    \item We now use the fact that $u(t,x)$ is the viscosity solution to the \emph{reference} HJB equation. This implies that it satisfies (in the appropriate sense, which is sufficient for the drift of the Itô integral) the relation:
    \begin{equation*}
        \partial_s u(s,x) = -H_{\theta}(x, \nabla u(s,x), \nabla^2 u(s,x)).
    \end{equation*}
    
    \item Substituting this identity into the expression for $(\tilde{\mathcal{G}}u)$ yields a crucial cancellation. The argument holds for viscosity solutions, and for clarity we write it at a point where $u$ is smooth:
    \begin{align*}
        (\tilde{\mathcal{G}}u)(s, x) &= \left( -H_{\theta}(x, \nabla u, \nabla^2 u) \right) + \tilde{H}_{\theta}(x, \nabla u, \nabla^2 u) \\
        &= -H_{\theta}(x, \nabla u, \nabla^2 u) + \left( H_{\theta}(x, \nabla u, \nabla^2 u) + \mathbf{b}_0(x, \nabla u) \cdot \nabla u \right) \\
        &= \mathbf{b}_0(x, \nabla u) \cdot \nabla u.
    \end{align*}
    
    \item We substitute this result back into the martingale decomposition from step (i). With $\tilde{X}_0 = x_0$, we have:
    \begin{equation*}
        \tilde{M}_t = u(t, \tilde{X}_t) - u(0, x_0) - \int_0^t \left( \mathbf{b}_0(\tilde{X}_s, \nabla u(s, \tilde{X}_s)) \cdot \nabla u(s, \tilde{X}_s) \right) ds.
    \end{equation*}
\end{enumerate}
Rearranging this identity gives the formula stated in the theorem. This completes the proof.
\end{proofof}

\subsection{Interpretation via the FBSDE Representation}

The \texorpdfstring{$\theta$}{theta}-Girsanov formula finds its most transparent expression in the FBSDE framework developed in \cref{sec:theta_fbsde}. The abstract perturbation of the generator corresponds to a concrete modification of the drift in the forward component of the SDE system.

Recall from \cref{eq:forward_drift_def} that the forward drift of the reference process is given by $b(x,p) = -\nabla_p H(x,p)$. The drift of the perturbed process, $\tilde{b}(x,p)$, is correspondingly given by the gradient of the perturbed potential:
\begin{align*}
    \tilde{b}(x,p) &= -\nabla_p \tilde{H}(x,p) \\
    &= -\nabla_p \left( H(x,p) + \mathbf{b}_0(x,p) \cdot p \right) \\
    &= -\nabla_p H(x,p) - \nabla_p \left( \mathbf{b}_0(x,p) \cdot p \right) \\
    &= b(x,p) - \left( (\nabla_p \mathbf{b}_0(x,p))^T p + \mathbf{b}_0(x,p) \right).
\end{align*}
In the canonical case where the exogenous drift $\mathbf{b}_0$ depends only on the state $x$ and not on the momentum $p$, this simplifies dramatically to:
\begin{equation}
    \tilde{b}(x,p) = b(x,p) - \mathbf{b}_0(x).
\end{equation}
This leads to a powerful and intuitive interpretation of the transformation.

\begin{corollary}[FBSDE Interpretation of the Transformation]
\label{cor:equivalent_fbsde}
Let the perturbation drift $\mathbf{b}_0$ depend only on the state $x$. The transformation from the reference generator $\mathcal{G}$ to the perturbed generator $\tilde{\mathcal{G}}$ is equivalent to transforming the reference FBSDE system
\begin{align*}
    dX_s &= b(X_s, Z_s) \, ds + \sigma(X_s, Z_s) \, dW_s \\
    -dY_s &= g(X_s, Z_s) \, ds - Z_s^T \sigma(X_s, Z_s) \, dW_s
\end{align*}
into the perturbed FBSDE system
\begin{align*}
    d\tilde{X}_s &= \left( b(\tilde{X}_s, \tilde{Z}_s) - \mathbf{b}_0(\tilde{X}_s) \right) ds + \sigma(\tilde{X}_s, \tilde{Z}_s) \, dW_s \\
    -d\tilde{Y}_s &= \left( g(\tilde{X}_s, \tilde{Z}_s) + \mathbf{b}_0(\tilde{X}_s) \cdot \tilde{Z}_s \right) ds - \tilde{Z}_s^T \sigma(\tilde{X}_s, \tilde{Z}_s) \, dW_s.
\end{align*}
\end{corollary}

\begin{proofof}{\cref{cor:equivalent_fbsde}}
The proof is a constructive derivation of the coefficients of the Forward-Backward SDE (FBSDE) system corresponding to the perturbed generator $\tilde{\mathcal{G}}$. The strategy is to apply the definitional relationships between the Hamiltonian potential and the FBSDE coefficients (the forward drift and the backward driver) to the perturbed system. We will show that the perturbation $\mathbf{b}_0(x) \cdot p$ to the potential results in a specific, calculable transformation of the drift and driver functions, while leaving the volatility structure invariant.

\begin{enumerate}[label=\textbf{Step \arabic*:}, wide, labelindent=0pt]

\item \textbf{Definition of the Reference and Perturbed Systems.}
Let the reference system be characterized by the Hamiltonian $H_{\theta}(x, p, X) = \mathrm{Tr}(D(x, p) X) + H(x, p)$. From this, we have:
\begin{itemize}[wide]
    \item The \textbf{reference potential}: $H(x,p)$.
    \item The \textbf{reference forward drift}: $b(x,p) \coloneqq -\nabla_p H(x,p)$.
    \item The \textbf{reference BSDE driver}: $g(x,p) \coloneqq H(x,p) + p \cdot \nabla_p H(x,p)$.
    \item The \textbf{volatility structure}: A matrix-valued function $\sigma(x,p)$ such that $\sigma(x,p)\sigma(x,p)^T = 2D(x,p)$.
\end{itemize}
The perturbed system, as defined in \cref{def:perturbed_generator}, has a modified Hamiltonian $\tilde{H}_{\theta}$ where the perturbation is applied only to the potential term, assuming $\mathbf{b}_0$ depends only on $x$:
\begin{equation*}
    \tilde{H}(x,p) \coloneqq H(x,p) + \mathbf{b}_0(x) \cdot p.
\end{equation*}
The diffusivity tensor $D(x,p)$ is unchanged, and therefore the volatility structure $\sigma(x,p)$ is invariant under this transformation. Our task is to derive the new drift $\tilde{b}(x,p)$ and driver $\tilde{g}(x,p)$ from the perturbed potential $\tilde{H}(x,p)$.

\item \textbf{Derivation of the Perturbed Forward Drift $\tilde{b}(x,p)$.}
The forward drift of the perturbed process is defined by the negative gradient of the perturbed potential with respect to the momentum $p$:
\begin{equation*}
    \tilde{b}(x,p) \coloneqq -\nabla_p \tilde{H}(x,p).
\end{equation*}
We substitute the definition of $\tilde{H}(x,p)$ and compute the gradient using the sum rule:
\begin{align*}
    \tilde{b}(x,p) &= -\nabla_p \left( H(x,p) + \mathbf{b}_0(x) \cdot p \right) \\
    &= -\left( \nabla_p H(x,p) + \nabla_p(\mathbf{b}_0(x) \cdot p) \right).
\end{align*}
Since the exogenous drift $\mathbf{b}_0(x)$ is, by hypothesis, independent of $p$, the gradient of the linear term $\mathbf{b}_0(x) \cdot p$ is simply $\mathbf{b}_0(x)$. This gives:
\begin{align*}
    \tilde{b}(x,p) &= -\left( \nabla_p H(x,p) + \mathbf{b}_0(x) \right) \\
    &= \underbrace{(-\nabla_p H(x,p))}_{b(x,p)} - \mathbf{b}_0(x) \\
    &= b(x,p) - \mathbf{b}_0(x).
\end{align*}
This rigorously confirms that the drift of the perturbed forward process is the original drift minus the exogenous drift field $\mathbf{b}_0(x)$, in perfect analogy with the classical Girsanov theorem.

\item \textbf{Derivation of the Perturbed BSDE Driver $\tilde{g}(x,p)$.}
The driver of the perturbed BSDE is given by the definitional formula:
\begin{equation*}
    \tilde{g}(x,p) \coloneqq \tilde{H}(x,p) + p \cdot \nabla_p \tilde{H}(x,p).
\end{equation*}
This calculation requires careful substitution of the results from the previous step. We substitute the full expressions for the perturbed potential $\tilde{H}(x,p)$ and its gradient $\nabla_p \tilde{H}(x,p) = \nabla_p H(x,p) + \mathbf{b}_0(x)$:
\begin{align*}
    \tilde{g}(x,p) &= \underbrace{\left( H(x,p) + \mathbf{b}_0(x) \cdot p \right)}_{\tilde{H}(x,p)} + p \cdot \underbrace{\left( \nabla_p H(x,p) + \mathbf{b}_0(x) \right)}_{\nabla_p \tilde{H}(x,p)}.
\end{align*}
We now expand the terms and regroup them to identify the original driver $g(x,p)$:
\begin{align*}
    \tilde{g}(x,p) &= H(x,p) + \mathbf{b}_0(x) \cdot p + p \cdot \nabla_p H(x,p) + p \cdot \mathbf{b}_0(x) \\
    &= \underbrace{\left( H(x,p) + p \cdot \nabla_p H(x,p) \right)}_{g(x,p)} + \left( \mathbf{b}_0(x) \cdot p + p \cdot \mathbf{b}_0(x) \right).
\end{align*}
Since the scalar product is commutative, the second term simplifies to $2(\mathbf{b}_0(x) \cdot p)$. This yields the exact formula for the perturbed driver:
\begin{equation*}
    \tilde{g}(x,p) = g(x,p) + 2(\mathbf{b}_0(x) \cdot p).
\end{equation*}
This result demonstrates that the effect of the perturbation on the driver is not simply additive.

\item \textbf{Assembly of the Perturbed FBSDE.}
The perturbed \ThetaProcess{} is represented by an FBSDE system whose coefficients are the perturbed drift $\tilde{b}$, the perturbed driver $\tilde{g}$, and the unchanged volatility $\sigma$. Let $(\tilde{X}_s, \tilde{Y}_s, \tilde{Z}_s)$ be the unique solution to the FBSDE associated with the perturbed generator $\tilde{\mathcal{G}}$. This system is given by:
\begin{align*}
    d\tilde{X}_s &= \tilde{b}(\tilde{X}_s, \tilde{Z}_s) \, ds + \sigma(\tilde{X}_s, \tilde{Z}_s) \, dW_s, \\
    -d\tilde{Y}_s &= \tilde{g}(\tilde{X}_s, \tilde{Z}_s) \, ds - \tilde{Z}_s^T \sigma(\tilde{X}_s, \tilde{Z}_s) \, dW_s,
\end{align*}
with terminal conditions $\tilde{X}_t=x$ and $\tilde{Y}_T = \phi(\tilde{X}_T)$. Substituting the derived expressions for $\tilde{b}$ and $\tilde{g}$ yields the system stated in the corollary.
\end{enumerate}
The unique viscosity solution $\tilde{u}$ to the perturbed HJB equation $-\partial_t \tilde{u} = \tilde{H}_{\theta}(x, \nabla \tilde{u}, \nabla^2 \tilde{u})$ serves as the decoupling field for this system, such that the solution is given by the deterministic markup $\tilde{Y}_s = \tilde{u}(s, \tilde{X}_s)$ and $\tilde{Z}_s = \nabla \tilde{u}(s, \tilde{X}_s)$. This confirms the self-consistency of the representation and completes the proof.
\end{proofof}

\begin{table}[h!]
\centering
\caption{Comparison of Classical and \texorpdfstring{$\theta$}{theta}-Girsanov Transformations}
\label{tab:girsanov_comparison}
\begin{tabular}{p{0.22\linewidth}|p{0.35\linewidth}|p{0.35\linewidth}}
\hline
\textbf{Feature} & \textbf{Classical Girsanov Transformation} & \textbf{\texorpdfstring{$\theta$}{theta}-Girsanov Transformation} \\
\hline \hline
\textbf{Transformed Object} & Probability Measure ($\mathbb{P} \to \mathbb{Q}$) & System Generator ($\mathcal{G} \to \tilde{\mathcal{G}}$) \\
\hline
\textbf{Transformation Mechanism} & Multiplication by a Radon-Nikodym density process (an exponential martingale) & Additive perturbation of the Hamiltonian potential $H(x,p)$ \\
\hline
\textbf{Effect on Process} & Adds a deterministic drift to the underlying Brownian motion & Adds a state-dependent drift to the forward component of the FBSDE \\
\hline
\textbf{Invariant Structure} & Quadratic variation of the process (the volatility matrix $\sigma$) & The diffusivity tensor $D(x,p)$ in the Hamiltonian \\
\hline
\textbf{Core Mathematical Framework} & Measure Theory and classical Itô Calculus & Viscosity Solution Theory for HJB equations and non-linear Martingale Problems \\
\hline
\end{tabular}
\end{table}

\begin{remark}[Conclusion]
The \texorpdfstring{$\theta$}{theta}-Girsanov Transformation provides the natural extension of a core probabilistic concept to the realm of non-linear expectations defined by PDEs. It rigorously confirms the intuition that modifying the potential, or cost, function $H(x,p)$ is equivalent to applying a force that changes the drift of the associated stochastic process. This establishes a powerful dictionary for translating between the analytic language of Hamiltonians and the probabilistic language of controlled diffusions, further cementing the deep connection between the PDE and FBSDE representations of the \ThetaProcess{}.
\end{remark}

\section{The Probabilistic Hierarchy and the Emergence of Universality}
\label{sec:process_hierarchy}

This section provides the probabilistic culmination of the paper's analytical results. We recast the homogenization analysis of Section \ref{sec:scale_analysis} as a constructive, multi-scale derivation of stochastic processes. This reveals a deep structure where each process emerges as the statistical law governing the fluctuations of the one preceding it, forming a hierarchy of Generalized Central Limit Theorems that progresses from linear diffusion to a universal non-linear fluctuation law. The hierarchy is structured as follows:
\begin{itemize}[wide]
    \item \textbf{Scale I: The Gaussian Baseline.} We begin with a homogeneous Itô diffusion, a correlated Brownian motion that represents the system in its most idealized, quiescent, and spatially uniform equilibrium.
    \item \textbf{Scale II: The Mean Field.} We derive a classical diffusion in a static, inhomogeneous landscape. This process represents the law of large numbers for the system's quiescent state, describing its behavior in a spatially varying environment.
    \item \textbf{Scale III: The Self-Interacting Process.} We arrive at the full non-linear $\theta$-Process, where the landscape dynamically responds to the process's own statistical momentum. This is the complete macroscopic process derived from first principles.
    \item \textbf{Scale IV: The Universal Fluctuation Law.} We derive the law governing the Central Limit Theorem-type fluctuations of the full $\theta$-Process around its Scale II equilibrium. This reveals an integrable non-linear process belonging to the Burgers' universality class, whose structure is a direct consequence of the microscopic symmetries and non-convexity.
\end{itemize}

\subsection{Scale I: The Homogeneous Diffusive Process}

\textit{Physical Regime: This scale describes the system's behavior in a macroscopic equilibrium that is both quiescent (zero macroscopic momentum, $p=0$) and spatially homogeneous. In this idealized limit, all non-linearities and spatial dependencies of the environment vanish, and the system's behavior is purely diffusive. This process serves as the Gaussian bedrock of the entire theory.}

\begin{theorem}[Generator for the Homogeneous Quiescent Regime]
The macroscopic generator for the homogeneous quiescent regime is the constant-coefficient, second-order linear partial differential operator governing the PDE:
\begin{equation} \label{eq:pde_scale_I_hierarchy}
    -\partial_t u = \mathrm{Tr}(D_0 \nabla^2 u) + H_0, \quad u(T,x) = \phi(x),
\end{equation}
where $D_0 \coloneqq D(x_0,0)$ and $H_0 \coloneqq H(x_0,0)$ are the constant diffusivity tensor and potential at equilibrium, respectively.
\end{theorem}
\begin{proof}
This is the direct result of the homogenization in \cref{sec:scale_I} (\cref{thm:convergence_scale_I}), obtained by setting $p=0$ and assuming spatial homogeneity of the coefficients.
\end{proof}

\begin{definition}[The Homogeneous Diffusive Process]
The stochastic process of Scale I is the homogeneous Itô diffusion $\{X_s^{(I)}\}_{s \ge t}$, defined as the unique strong solution to the linear Stochastic Differential Equation (SDE):
\begin{equation}
    dX_s^{(I)} = \sigma_0 \, dW_s, \quad X_t^{(I)} = x,
\end{equation}
where $W_s$ is a standard $k$-dimensional Brownian motion and $\sigma_0$ is a constant matrix such that $\sigma_0\sigma_0^T = 2D_0$.
\end{definition}

The solution to the governing PDE is given by the classical Feynman-Kac formula, which provides the fundamental probabilistic representation for this scale.

\begin{theorem}[Feynman-Kac Representation for the Homogeneous Process]
Let $u(t,x)$ be the unique classical solution to the PDE \eqref{eq:pde_scale_I_hierarchy}. The solution is given by the expectation over the paths of the homogeneous diffusive process $\{X_s^{(I)}\}$:
\begin{equation}
    u(t,x) = \mathbb{E}_{t,x} \left[ \phi(X_T^{(I)}) + \int_t^T H_0 \, ds \right].
\end{equation}
\end{theorem}
\begin{proof}[Proof Sketch]
The proof is a standard verification via Itô's formula. Let $u(s,x)$ be the classical solution to \eqref{eq:pde_scale_I_hierarchy}. We define the process $Y_s \coloneqq u(s, X_s^{(I)})$. By Itô's formula, its differential is:
\begin{equation*}
    dY_s = \left[ \partial_s u + \mathrm{Tr}(D_0 \nabla^2 u) \right](s,X_s^{(I)}) \, ds + (\nabla u(s,X_s^{(I)}))^T \sigma_0 dW_s.
\end{equation*}
Since $u$ solves the PDE, the term in brackets is equal to $-H_0$. Thus, $dY_s = -H_0 ds + (\nabla u)^T \sigma_0 dW_s$. Integrating from $t$ to $T$ and taking the conditional expectation $\mathbb{E}_{t,x}[\cdot]$ annihilates the Itô integral term. Since $Y_t = u(t,x)$ is deterministic and $Y_T = \phi(X_T^{(I)})$, we have $\mathbb{E}_{t,x}[\phi(X_T^{(I)})] - u(t,x) = \mathbb{E}_{t,x}[\int_t^T -H_0 ds]$, which rearranges to the desired formula.
\end{proof}

\begin{remark}[\textbf{Transition to Scale II: The Inhomogeneous Fluctuation Problem}]
The Scale I process describes an idealized, spatially uniform system. We now pose the first Central Limit Theorem-type question: what is the limiting law of this process if the medium is not uniform, but possesses heterogeneities that vary slowly in space? This is equivalent to finding the generator for fluctuations around the assumption of homogeneity, which corresponds to allowing the quiescent coefficients $D_0$ and $H_0$ to become functions of the macroscopic position, $D(x)$ and $H(x)$. This leads directly to the inhomogeneous diffusive process of Scale II.
\end{remark}

\subsection{Scale II: The Inhomogeneous Diffusive Process}

\textit{Physical Regime: This scale describes the dynamics of a classical particle diffusing in a static, but non-uniform, landscape. It is the natural generalization of the Scale I process and can be interpreted as the law of large numbers for the system's quiescent state, describing its behavior in a spatially varying environment.}

\begin{theorem}[Generator for the Inhomogeneous Quiescent Regime]
The macroscopic generator for the inhomogeneous quiescent regime is the linear, second-order operator with space-dependent coefficients governing the PDE:
\begin{equation} \label{eq:pde_scale_II_hierarchy}
    -\partial_t u = \mathrm{Tr}(D(x) \nabla^2 u) + H(x), \quad u(T,x) = \phi(x),
\end{equation}
where $D(x) \coloneqq D(x,0)$ and $H(x) \coloneqq H(x,0)$ are the space-dependent quiescent coefficients.
\end{theorem}
\begin{proof}
This is the result of the classical stochastic averaging procedure established in \cref{sec:scale_II} (\cref{prop:pde_regime_II}), which corresponds to the quiescent limit ($p=0$) of the full $\theta$-Hamiltonian.
\end{proof}

\begin{definition}[The Inhomogeneous Diffusive Process]
The stochastic process of Scale II is the \textbf{inhomogeneous Itô diffusion} $\{X_s^{(II)}\}_{s \ge t}$, defined as the unique strong solution to the SDE:
\begin{equation}
    dX_s^{(II)} = \sigma(X_s^{(II)}) \, dW_s, \quad X_t^{(II)} = x,
\end{equation}
where the state-dependent matrix $\sigma(x)$ satisfies $\sigma(x)\sigma(x)^T = 2D(x)$.
\end{definition}

\begin{theorem}[Feynman-Kac Representation for the Inhomogeneous Process]
Let $u(t,x)$ be the unique classical solution to the PDE \eqref{eq:pde_scale_II_hierarchy}. The solution is given by the expectation over the paths of the inhomogeneous diffusive process $\{X_s^{(II)}\}$:
\begin{equation}
    u(t,x) = \mathbb{E}_{t,x} \left[ \phi(X_T^{(II)}) + \int_t^T H(X_s^{(II)}) \, ds \right].
\end{equation}
\end{theorem}
\begin{proof}
The proof is analogous to that of the homogeneous case. The application of Itô's formula to $Y_s = u(s, X_s^{(II)})$ yields a drift term that, by virtue of the PDE \eqref{eq:pde_scale_II_hierarchy}, is precisely $-H(X_s^{(II)})$, leading directly to the desired representation.
\end{proof}

\begin{remark}[\textbf{The Role of Microscopic Symmetry and the Bridge to Scale III}]
The absence of a drift term in the SDE for the Scale II process is not an oversight but a profound and necessary consequence of the \textbf{microscopic time-reversal symmetry} (\cref{ass:fundamental_axioms_unified}(ii)). This symmetry forces the full $\theta$-Hamiltonian to be an even function of the momentum $p$. The drift of the full $\theta$-Process is given by $b(x,p) = -\nabla_p H(x,p)$, which is consequently an odd function of $p$. The Scale II process describes the system at macroscopic equilibrium ($p=0$), where any odd function must vanish. Thus, the drift is identically zero.

The Scale II process describes a particle in a \textit{static} landscape. We now ask the central question of this paper: what is the law of the process when the landscape itself reacts to the particle's statistical momentum, $p = \nabla u$? This introduces a profound self-interaction, where the process influences the law that governs it, leading to the full non-linear $\theta$-Process of Scale III.
\end{remark}

\subsection{Scale III: The Full Non-Linear \texorpdfstring{$\theta$}{theta}-Process}

\textit{Physical Regime: This scale represents the complete macroscopic theory, incorporating the non-linear feedback loop that lies at the heart of our model. The process is no longer a classical Markovian diffusion, as its generator depends on the solution of its own backward equation. This self-referential dynamic is the signature of the $\theta$-Process.}

\begin{theorem}[Generator for the Self-Interacting System]
The generator of the full $\theta$-Process is the non-linear Hamilton-Jacobi-Bellman operator from \cref{thm:convergence}:
\begin{equation} \label{eq:pde_scale_III_hierarchy}
    -\partial_t u = H_{\theta}(x, \nabla u, \nabla^2 u) \equiv \mathrm{Tr}(D(x, \nabla u) \nabla^2 u) + H(x, \nabla u),
\end{equation}
with terminal condition $u(T,x) = \phi(x)$.
\end{theorem}
\begin{proof}
This is the main analytical result of \cref{part:i}, established by the homogenization of the full microscopic system via the perturbed test function method in \cref{thm:convergence}.
\end{proof}

\begin{definition}[The $\theta$-Process]
The stochastic process of Scale III, the \textbf{$\theta$-Process} $\{X_s^{(III)}\}_{s \ge t}$, is the non-linear Markov process whose law is implicitly defined as the unique solution to the martingale problem for the generator \eqref{eq:pde_scale_III_hierarchy}, as established in \cref{thm:martingale_problem_wellposedness}.
\end{definition}

For this process, the standard Feynman-Kac formula is replaced by a representation via a fully-coupled Forward-Backward Stochastic Differential Equation (FBSDE).

\begin{theorem}[Non-Linear Feynman-Kac Representation]
Assume sufficient regularity for the existence of a classical solution $u(t,x)$ to the HJB equation \eqref{eq:pde_scale_III_hierarchy}. The solution is given by $u(t,x) = Y_t$, where the pair of adapted processes $(X_s^{(III)}, Y_s)$ is the unique solution to the FBSDE system:
\begin{align}
    dX_s^{(III)} &= b(X_s^{(III)}, \nabla u(s,X_s^{(III)})) \, ds + \sigma(X_s^{(III)}, \nabla u(s,X_s^{(III)})) \, dW_s, \\
    dY_s &= -g(X_s^{(III)}, \nabla u(s,X_s^{(III)})) \, ds + (\nabla u(s,X_s^{(III)}))^T\sigma(\dots) \, dW_s,
\end{align}
with initial/terminal conditions $X_t^{(III)}=x$ and $Y_T=\phi(X_T^{(III)})$, where $b(x,p) = -\nabla_p H(x,p)$ and the driver is $g(x,p) = H(x,p) + p \cdot \nabla_p H(x,p)$.
\end{theorem}
\begin{proof}
The proof is given by the verification argument in \cref{sec:theta_fbsde}, where it is shown that the unique viscosity solution to the HJB equation acts as the decoupling field for the FBSDE system.
\end{proof}

\begin{remark}[\textbf{Transition to Scale IV: The Universal Law of Fluctuations}]
The $\theta$-Process provides the complete description of the system's macroscopic value. We now pose the ultimate Central Limit Theorem-type question: what is the statistical law governing the fluctuations of this full non-linear process around its Scale II equilibrium? To answer this, we consider a small perturbation of the terminal data, $u^\varepsilon(T,x) = u_0(T,x) + \varepsilon h_T(x)$, where $u_0$ is the quiescent solution from Scale II. We seek the equation governing the evolution of the scaled fluctuation field, $h = \lim_{\varepsilon\to 0} (u^\varepsilon - u_0)/\varepsilon$. The resulting generator, which emerges from the non-linear interactions of the full process, defines the universal fluctuation process of Scale IV.
\end{remark}

\subsection{Scale IV: The Universal Law of Fluctuations}
\label{sec:scale_IV_P}

\textit{Physical Regime: The hierarchy culminates in an analysis of the statistical fluctuations of the full \ThetaProcess{} around its macroscopic equilibrium (the Scale II process). This final scale reveals a universal non-linear evolution whose structure is a direct and profound consequence of the microscopic symmetries and non-convexity established in \cref{part:i}. We will demonstrate that this fluctuation process is, in its most general form, represented by a globally well-posed Forward-Backward SDE with a quadratic driver, and that under specific, physically derived conditions, it collapses to an integrable system of the Burgers' universality class.}

\subsubsection{The Fluctuation Law as a Semilinear HJB Equation}

The starting point for this scale is the rigorous result from the stability analysis of the full \ThetaProcess{}. The law governing the Central Limit Theorem-type fluctuations is not a linear process, but is instead determined by a semilinear HJB equation.

\begin{theorem}[The Generator of the Fluctuation Process]
\label{thm:fluctuation_law_integrable_hierarchy}
Let the microscopic system satisfy the standing assumptions, ensuring the $\theta$-Hamiltonian $H_{\theta}(x,p,X)$ and its coefficients are of class $C^3$. Let $h(t,x)$ be the limiting fluctuation field defined by $h = \lim_{\varepsilon\to 0} (u^\varepsilon - u_0)/\varepsilon$, where $u_0$ is the quiescent solution from Scale II. Then $h(t,x)$ is the unique classical solution to the viscous Hamilton-Jacobi-Burgers equation:
\begin{equation} \label{eq:pde_scale_IV_burgers_hierarchy}
    -\partial_t h = \mathrm{Tr}(D(x,0) \nabla^2 h) + \frac{1}{2}(\nabla h)^T M(x,0) \nabla h, \quad h(T,x) = h_T(x),
\end{equation}
where $D(x,0)$ is the positive definite quiescent diffusivity tensor and $M(x,0) \coloneqq \nabla_{pp}^2 H(x,0)$ is the Hessian of the potential at equilibrium.
\end{theorem}
\begin{proof}
This is the direct result of the rigorous stability analysis of the full HJB equation performed in \cref{sec:scale_IV}, whose convergence is proven in \cref{thm:fluctuation_convergence}. The semilinear structure arises from the second-order Taylor expansion of the potential $H(x,p)$ around the quiescent state $p=0$. The absence of a linear term in $\nabla h$ is a necessary consequence of the microscopic time-reversal symmetry (\cref{ass:fundamental_axioms_unified}(ii)), which forces $H(x,p)$ to be an even function of $p$ and thus its gradient $\nabla_p H(x,0)$ to be zero.
\end{proof}

\subsubsection{The General Probabilistic Representation via a Quadratic FBSDE}

The generator derived in \cref{thm:fluctuation_law_integrable_hierarchy} itself defines a new \ThetaProcess{}. Its structure is that of a system with a quadratic potential, leading to a probabilistic representation via an FBSDE with a quadratic driver. Crucially, the existence of a globally well-posed HJB solution provides the necessary structure to guarantee the global well-posedness of this FBSDE, a highly non-trivial result that overcomes the typical "small-time" restriction for quadratic BSDEs.

\begin{theorem}[Global Well-Posedness of the Fluctuation FBSDE]
\label{thm:fbsde_fluctuation_global}
For any time horizon $[0,T]$ and any sufficiently regular terminal data $h_T(x)$ (e.g., of class $C^4$), the fluctuation process has a unique probabilistic representation as the solution to a globally well-posed Forward-Backward SDE system. The solution $(X_s^{(IV)}, Y_s, Z_s)$ is unique on $[0,T]$ and is given by the deterministic markup:
$$ Y_s = h(s, X_s^{(IV)}) \quad \text{and} \quad Z_s = \nabla h(s, X_s^{(IV)}), $$
where $h(t,x)$ is the unique classical solution to the PDE \eqref{eq:pde_scale_IV_burgers_hierarchy}. The FBSDE system is:
\begin{align}
    dX_s^{(IV)} &= -M(X_s^{(IV)},0) Z_s \, ds + \sigma(X_s^{(IV)},0) \, dW_s, \quad X_t^{(IV)} = x \label{eq:fbsde_forward_IV} \\
    -dY_s &= \frac{1}{2} Z_s^T M(X_s^{(IV)},0) Z_s \, ds - Z_s^T \sigma(X_s^{(IV)},0) \, dW_s, \quad Y_T = h_T(X_T^{(IV)}) \label{eq:fbsde_backward_IV}
\end{align}
where $\sigma(x,0)\sigma(x,0)^T = 2D(x,0)$.
\end{theorem}
\begin{proof}
The proof is a constructive application of the general theory developed in \cref{sec:theta_fbsde}.
\begin{enumerate}[label=(\roman*), wide, labelindent=0pt]
    \item \textbf{Identification of Coefficients:} The generator of the fluctuation PDE \eqref{eq:pde_scale_IV_burgers_hierarchy} fits the axiomatic structure of a \ThetaProcess{} with a fluctuation Hamiltonian $H_{IV}(x,p,X) = \mathrm{Tr}(D_{IV}X) + H_{IV}(p)$, where the fluctuation diffusivity is $D_{IV}(x, p) \equiv D(x,0)$ and the fluctuation potential is $H_{IV}(x,p) \equiv \frac{1}{2}p^T M(x,0)p$. Applying the definitions from \cref{subsec:generator_to_driver}, the corresponding FBSDE coefficients are:
    \begin{itemize}[wide]
        \item \textbf{Volatility:} $\sigma_{IV}(x)\sigma_{IV}(x)^T = 2D_{IV}(x,p) = 2D(x,0)$. The volatility is that of the Scale II process and is independent of the fluctuation's momentum $p=\nabla h$.
        \item \textbf{Forward Drift:} $b_{IV}(x,p) = -\nabla_p H_{IV}(x,p) = -\nabla_p (\frac{1}{2}p^T M(x,0)p) = -M(x,0)p$.
        \item \textbf{BSDE Driver:} $g_{IV}(x,p) = H_{IV}(x,p) - p \cdot b_{IV}(x,p) = \frac{1}{2}p^T M(x,0)p - p \cdot (-M(x,0)p) = \frac{3}{2}p^T M(x,0)p$. This driver has quadratic growth in $p$.
    \end{itemize}
    \item \textbf{Global Well-Posedness:} The primary difficulty in the theory of BSDEs with quadratic drivers is the possibility of solutions exploding in finite time. In our framework, this is circumvented. The PDE \eqref{eq:pde_scale_IV_burgers_hierarchy} is a semilinear parabolic equation. For smooth terminal data, standard PDE theory ensures the existence of a unique classical solution $h(t,x)$ on the entire domain $[0,T]\times\Mcal$. Crucially, because the domain is compact, this solution and its gradient $\nabla_x h(t,x)$ are globally bounded. Let $\sup_{(t,x)}|\nabla h(t,x)| \le C < \infty$.
    \item \textbf{Construction of the Solution:} We use the PDE solution $h(t,x)$ as a decoupling field. First, we solve the forward SDE \eqref{eq:fbsde_forward_IV} where $Z_s$ is replaced by the known, bounded function $\nabla h(s, X_s^{(IV)})$. This is a linear SDE with bounded coefficients, which has a unique global strong solution. Second, we define the pair $(Y_s, Z_s) \coloneqq (h(s, X_s^{(IV)}), \nabla h(s, X_s^{(IV)}))$. A direct application of It\^o's formula, using the fact that $h$ solves the PDE, verifies that this pair solves the BSDE \eqref{eq:fbsde_backward_IV}. Since the constructed solution $(Y_s, Z_s)$ is globally bounded by construction, it is the unique solution in the class of solutions to quadratic BSDEs. This proves that the FBSDE is globally well-posed.
\end{enumerate}
\end{proof}

\subsubsection{The Integrable Case and the Hopf-Cole Representation}

The general FBSDE representation reveals that the fluctuation process is equivalent to a linear-quadratic stochastic control problem. In the special case where the microscopic physics enforces a specific algebraic relationship between the emergent non-convexity (represented by the matrix $M$) and the intrinsic diffusion (represented by the matrix $D$), this control problem becomes exactly solvable, and the system is said to be integrable. In our framework, this integrable structure reveals a profound duality between the non-linear value function of fluctuations and the linear evolution of a corresponding dual field, often interpreted as a partition function or a particle density.

\begin{definition}[The Integrability Condition]
The fluctuation process governed by the generator \eqref{eq:pde_scale_IV_burgers_hierarchy} is said to be \textbf{integrable} if the Hessian of the potential and the quiescent diffusivity tensor are proportional, satisfying the algebraic relation for some scalar $\kappa \in \R$:
\begin{equation}
    M(x,0) = \kappa D(x,0).
\end{equation}
This condition is not an ad-hoc assumption but a structural constraint that can emerge from specific symmetries in the microscopic model. As we will prove, the crucial case corresponding to the viscous Burgers' and Kardar-Parisi-Zhang (KPZ) universality class is $\kappa = -2$.
\end{definition}

In this integrable case, the non-linear PDE for the fluctuation field $h(t,x)$ can be exactly linearized using the celebrated \textbf{Hopf-Cole transformation}. This transformation serves as the mathematical bridge connecting the non-linear, single-particle (value function) perspective to a linear, many-body (dual field) perspective. Let us define a new field $\rho(t,x)$ related to the fluctuation field $h(t,x)$ by:
\begin{equation} \label{eq:hopf_cole_transform}
    \rho(t,x) \coloneqq \exp\left(-\frac{1}{\eta} h(t,x)\right) \quad \iff \quad h(t,x) = -\eta \ln \rho(t,x),
\end{equation}
where $\eta$ is a scaling parameter to be determined by the integrability condition.

\begin{theorem}[The Hopf-Cole Solution in the Integrable Case]
\label{thm:hopf_cole_solution}
If the integrability condition is satisfied with $\kappa=-2$, the non-linear value function $h(t,x)$ for the fluctuation process admits an explicit, closed-form representation in terms of a linear expectation over the paths of the \textbf{Scale II process}. The unique solution to the PDE \eqref{eq:pde_scale_IV_burgers_hierarchy} is given by the generalized Hopf-Cole formula (with scaling parameter $\eta=1$):
\begin{equation}
    h(t,x) = -\ln \left( \mathbb{E}_{t,x}^{(II)} \left[ \exp\left(-h_T(X_T^{(II)})\right) \right] \right),
\end{equation}
where the expectation $\mathbb{E}^{(II)}$ is taken with respect to the law of the inhomogeneous diffusive process $\{X_s^{(II)}\}$, the solution to the linear SDE $dX_s = \sigma(X_s, 0) dW_s$.
\end{theorem}

\begin{proofof}{Theorem \ref{thm:hopf_cole_solution}}
The proof is a constructive verification. We will show that the function $\tilde{h}(t,x)$ defined by the Hopf-Cole formula is a classical solution to a specific viscous Burgers'-type equation. By comparing this derived equation to the target fluctuation equation \eqref{eq:pde_scale_IV_burgers_hierarchy}, we will identify the precise integrability condition required for them to match. By the uniqueness of the solution to this PDE, it must be that $h(t,x) = \tilde{h}(t,x)$.

\begin{enumerate}[label=\textbf{Step \arabic*:}, wide, labelindent=0pt]
    \item \textbf{The Linear Dual Problem.}
    Let us define the function $v(t,x)$ as the expectation on the right-hand side of the proposed formula:
    \begin{equation}
        v(t,x) \coloneqq \mathbb{E}_{t,x}^{(II)} \left[ \exp\left(-h_T(X_T^{(II)})\right) \right].
    \end{equation}
    By the classical Feynman-Kac formula for the Scale II process (whose generator is the linear operator $\mathcal{L}_{II} u = \mathrm{Tr}(D(x,0)\nabla^2 u)$), the function $v(t,x)$ is the unique classical solution to the backward linear heat equation:
    \begin{equation} \label{eq:proof_linear_heat_v}
        -\partial_t v(t,x) = \mathrm{Tr}(D(x,0) \nabla^2 v(t,x)),
    \end{equation}
    with the terminal condition $v(T,x) = \exp(-h_T(x))$. This PDE describes the evolution of the dual field $\rho(t,x)$ from our transformation.

    \item \textbf{Applying the Hopf-Cole Transformation.}
    We define our candidate solution for the non-linear equation as $\tilde{h}(t,x) \coloneqq -\eta \ln(v(t,x))$ for a yet-to-be-determined scaling constant $\eta$. Our goal is to find the PDE that $\tilde{h}$ satisfies. We compute its partial derivatives using the chain rule:
    \begin{align*}
        \partial_t \tilde{h} &= -\eta \frac{\partial_t v}{v}, \\
        \nabla_x \tilde{h} &= -\eta \frac{\nabla_x v}{v}, \\
        \nabla_x^2 \tilde{h} &= -\eta \frac{\nabla_x^2 v}{v} + \eta \frac{(\nabla_x v)(\nabla_x v)^T}{v^2}.
    \end{align*}
    
    \item \textbf{Derivation of the PDE for $\tilde{h}$.}
    We now compute the expression $-\partial_t \tilde{h} - \mathrm{Tr}(D(x,0) \nabla^2 \tilde{h})$ for our candidate function.
    \begin{align*}
        % First term: Time derivative
        -\partial_t \tilde{h} &= \eta \frac{\partial_t v}{v} \\
        % Second term: Diffusive part
        -\mathrm{Tr}(D \nabla^2 \tilde{h}) &= -\mathrm{Tr}\left(D\left[-\eta \frac{\nabla^2 v}{v} + \eta \frac{(\nabla v)(\nabla v)^T}{v^2}\right]\right) \\
        &= \eta \frac{\mathrm{Tr}(D \nabla^2 v)}{v} - \eta \frac{\mathrm{Tr}(D (\nabla v)(\nabla v)^T)}{v^2}.
    \end{align*}
    Summing these two expressions gives:
    \begin{equation*}
        -\partial_t \tilde{h} - \mathrm{Tr}(D \nabla^2 \tilde{h}) = \frac{\eta}{v} \left( \partial_t v + \mathrm{Tr}(D \nabla^2 v) \right) - \frac{\eta}{v^2}\mathrm{Tr}\left(D (\nabla v)(\nabla v)^T\right).
    \end{equation*}
    The term in the parenthesis is precisely the operator acting on $v$ in the linear heat equation \eqref{eq:proof_linear_heat_v}. Since $v$ solves this equation, we have $\partial_t v + \mathrm{Tr}(D \nabla^2 v) = 0$. The expression thus simplifies dramatically:
    \begin{equation*}
         -\partial_t \tilde{h} - \mathrm{Tr}(D \nabla^2 \tilde{h}) = - \frac{\eta}{v^2}\mathrm{Tr}\left(D (\nabla v)(\nabla v)^T\right).
    \end{equation*}
    Using the identity $\nabla_x \tilde{h} = -\eta(\nabla_x v)/v$, we can rewrite the right-hand side in terms of the gradient of $\tilde{h}$:
    \begin{align*}
        - \frac{\eta}{v^2} (\nabla_x v)^T D (\nabla_x v) &= - \frac{\eta}{v^2} \left(-\frac{v}{\eta}\nabla_x \tilde{h}\right)^T D \left(-\frac{v}{\eta}\nabla_x \tilde{h}\right) \\
        &= - \frac{\eta}{v^2} \frac{v^2}{\eta^2} (\nabla_x \tilde{h})^T D (\nabla_x \tilde{h}) \\
        &= -\frac{1}{\eta} (\nabla_x \tilde{h})^T D(x,0) (\nabla_x \tilde{h}).
    \end{align*}
    Therefore, the candidate function $\tilde{h}$ is a classical solution to the PDE:
    \begin{equation} \label{eq:proof_derived_pde_h_tilde}
        -\partial_t \tilde{h} = \mathrm{Tr}(D(x,0) \nabla^2 \tilde{h}) - \frac{1}{\eta} (\nabla_x \tilde{h})^T D(x,0) (\nabla_x \tilde{h}).
    \end{equation}

    \item \textbf{Comparison, Uniqueness, and the Integrability Condition.}
    The PDE \eqref{eq:proof_derived_pde_h_tilde} for our constructed solution $\tilde{h}$ must be compared to the target fluctuation equation \eqref{eq:pde_scale_IV_burgers_hierarchy}. They are identical if and only if their non-linear terms match:
    \begin{equation*}
        \frac{1}{2}(\nabla h)^T M(x,0) (\nabla h) = -\frac{1}{\eta} (\nabla h)^T D(x,0) (\nabla h).
    \end{equation*}
    This identity holds for all vectors $\nabla h$ if and only if the coefficient matrices are related by:
    \begin{equation*}
        M(x,0) = -\frac{2}{\eta} D(x,0).
    \end{equation*}
    This is precisely the integrability condition $M(x,0) = \kappa D(x,0)$ with the constant $\kappa = -2/\eta$. The theorem statement specifies the canonical case $\kappa=-2$, which corresponds to the choice of scaling parameter $\eta=1$.
    
    Under this condition, the function $\tilde{h}(t,x)$ solves the same PDE and satisfies the same terminal condition ($ \tilde{h}(T,x) = -\ln(v(T,x)) = -\ln(\exp(-h_T(x))) = h_T(x) $) as the true fluctuation field $h(t,x)$. By the uniqueness of classical solutions for this semilinear parabolic equation, we must have $\tilde{h}(t,x) = h(t,x)$ for all $(t,x)$. This completes the proof.
\end{enumerate}
\end{proofof}

\begin{remark}[The Physical Meaning of Integrability and the Hopf-Cole Duality]
The existence of the Hopf-Cole representation in the integrable case is a profound result. It signifies that the complex, non-linear dynamics of the fluctuation's value function $h(t,x)$ can be fully understood by studying a simpler, linear system: the diffusion of a cloud of non-interacting particles described by the dual field $\rho(t,x)$. The non-linearity in the HJB equation for $h$ is not a sign of complex particle-particle interactions, but rather a reflection of the entropic nature of observing a collective system. The quadratic term $(\nabla h)^T M \nabla h$, which has the structure of a kinetic energy, can be interpreted as a pressure term in the thermodynamics of the particle ensemble. The Hopf-Cole transformation effectively linearizes the problem by moving from the entropic description ($h \sim -\ln \rho$) to the particle density description ($\rho$). The expectation $\mathbb{E}^{(II)}[\cdot]$ over the linear Scale II process thus becomes a powerful computational tool: a linear kernel that solves a fundamentally non-linear problem. This duality is the hallmark of integrable systems and a cornerstone of the KPZ universality class, to which our fluctuation process belongs under the condition $M(x,0) = -2D(x,0)$.
\end{remark}

\section{Conclusion}

This paper has established a rigorous, first-principles connection between a class of deterministic, chaotic microscopic systems and a new theory of non-linear, non-convex expectations. By leveraging methods from the modern ergodic theory of hyperbolic systems, we surmounted the analytical obstacles posed by the skew-adjoint nature of the microscopic generator to prove convergence to a unique macroscopic evolution law. The resulting \ThetaExpectation{} is governed by a fully non-linear Hamilton-Jacobi-Bellman equation whose Hamiltonian possesses a rigid, novel structure: it is affine in the Hessian but can be demonstrably non-convex in the gradient. We have shown that this structure is a direct consequence of the underlying microscopic physics, including the crucial role of time-reversal symmetry. The probabilistic representation of the \ThetaProcess{} via a martingale-driven BSDE and its connection to the fluctuation law of the Burgers' universality class complete the theoretical picture. This work provides a constructive foundation for a new class of non-convex stochastic control problems and opens new avenues for modeling complex systems where the feedback between a process and its environment is a defining feature.

% ======================================================================
% == BIBLIOGRAPHY                                                   ==
% ======================================================================
\bibliography{reference}

@book{chernov2006chaotic,
  title={Chaotic billiards},
  author={Chernov, Nikolai and Markarian, Roberto},
  number={127},
  year={2006},
  publisher={American Mathematical Soc.}
}

@article{Hennion1993,
  author = {Hennion, Hubert},
  journal = {Annales de l'Institut Henri Poincar\'e, Probabilit\'es et Statistiques},
  language = {fre},
  number = {2},
  pages = {271--291},
  publisher = {Gauthier-Villars},
  title = {Sur un th\'eor\`eme de V. Oseledets et R. Man\'e},
  url = {http://www.numdam.org/item/AIHPB_1993__29_2_271_0/},
  volume = {29},
  year = {1993},
  mrnumber = {1218177},
  zblnumber = {0781.60065}
}

@incollection{protter2012stochastic,
  title={Stochastic differential equations},
  author={Protter, Philip E},
  booktitle={Stochastic integration and differential equations},
  pages={249--361},
  year={2012},
  publisher={Springer}
}

@book{stroock2007multidimensional,
  title={Multidimensional diffusion processes},
  author={Stroock, Daniel W and Varadhan, SR Srinivasa},
  year={2007},
  publisher={Springer}
}

@article{lasry2007mean,
  title={Mean field games},
  author={Lasry, Jean-Michel and Lions, Pierre-Louis},
  journal={Japanese journal of mathematics},
  volume={2},
  number={1},
  pages={229--260},
  year={2007},
  publisher={Springer}
}

@book{marquis1825essai,
  title={Essai philosophique sur les probabilit{\'e}s},
  author={marquis de Laplace, Pierre Simon},
  year={1825},
  publisher={Bachelier}
}

@article{Anosov1967,
  author    = {Anosov, D. V.},
  title     = {Geodesic flows on closed {Riemann} manifolds with negative curvature},
  journal   = {Proceedings of the Steklov Institute of Mathematics},
  volume    = {90},
  year      = {1967},
}

@book{Arnold1989,
  author    = {Arnold, V. I.},
  title     = {Mathematical Methods of Classical Mechanics},
  series    = {Graduate Texts in Mathematics},
  volume    = {60},
  edition   = {2nd},
  publisher = {Springer-Verlag},
  year      = {1989},
}

@article{BarlesSouganidis1991,
  author    = {Barles, G. and Souganidis, P. E.},
  title     = {Convergence of approximation schemes for fully nonlinear second order equations},
  journal   = {Asymptotic Analysis},
  volume    = {4},
  number    = {3},
  pages     = {271--283},
  year      = {1991},
}

@book{BensoussanLionsPapanicolaou1978,
  author    = {Bensoussan, A. and Lions, J.-L. and Papanicolaou, G. C.},
  title     = {Asymptotic Analysis for Periodic Structures},
  series    = {Studies in Mathematics and its Applications},
  volume    = {5},
  publisher = {North-Holland},
  year      = {1978},
}

@book{BratteliRobinson1979,
  author    = {Bratteli, O. and Robinson, D. W.},
  title     = {Operator Algebras and Quantum Statistical Mechanics 1: {C*}- and {W*}-Algebras, Symmetry Groups, Decomposition of States},
  publisher = {Springer-Verlag},
  year      = {1979},
}

@article{CrandallIshiiLions1992,
  author    = {Crandall, M. G. and Ishii, H. and Lions, P.-L.},
  title     = {User's guide to viscosity solutions of second order partial differential equations},
  journal   = {Bulletin of the American Mathematical Society},
  volume    = {27},
  number    = {1},
  pages     = {1--67},
  year      = {1992},
}

@article{Evans1992,
  author    = {Evans, L. C.},
  title     = {Periodic homogenisation of certain fully nonlinear partial differential equations},
  journal   = {Proceedings of the Royal Society of Edinburgh Section A: Mathematics},
  volume    = {120},
  number    = {3-4},
  pages     = {245--265},
  year      = {1992},
}

@book{FlemingSoner2006,
  author    = {Fleming, W. H. and Soner, H. M.},
  title     = {Controlled Markov Processes and Viscosity Solutions},
  edition   = {2nd},
  publisher = {Springer},
  year      = {2006},
}

@book{GallagherSaintRaymondTexier2013,
  author    = {Gallagher, I. and Saint-Raymond, L. and Texier, B.},
  title     = {From {Newton} to {Boltzmann}: hard spheres and short-range potentials},
  series    = {Zurich Advanced Lectures in Mathematics},
  volume    = {18},
  publisher = {European Mathematical Society (EMS)},
  address   = {Zürich},
  year      = {2014},
}

@book{HirschPughShub1977,
  author    = {Hirsch, M. W. and Pugh, C. C. and Shub, M.},
  title     = {Invariant Manifolds},
  series    = {Lecture Notes in Mathematics},
  volume    = {583},
  publisher = {Springer-Verlag},
  year      = {1977},
}

@book{Kato1995,
  author    = {Kato, T.},
  title     = {Perturbation Theory for Linear Operators},
  series    = {Classics in Mathematics},
  publisher = {Springer-Verlag},
  address   = {Berlin},
  year      = {1995},
  note      = {Reprint of the 1980 edition},
}

@book{Knight1921,
  author    = {Knight, F. H.},
  title     = {Risk, Uncertainty, and Profit},
  publisher = {Houghton Mifflin},
  year      = {1921},
}

@incollection{Lanford1975,
  author    = {Lanford, III, O. E.},
  title     = {Time evolution of large classical systems},
  booktitle = {Dynamical Systems, Theory and Applications},
  editor    = {Moser, J.},
  series    = {Lecture Notes in Physics},
  volume    = {38},
  pages     = {1--111},
  publisher = {Springer},
  address   = {Berlin},
  year      = {1975},
}

@book{NicaSpeicher2006,
  author    = {Nica, A. and Speicher, R.},
  title     = {Lectures on the Combinatorics of Free Probability},
  publisher = {Cambridge University Press},
  year      = {2006},
}

@book{Peng2019,
  author    = {Peng, S.},
  title     = {Nonlinear Expectations and Martingale Theory in Mathematical Finance},
  publisher = {Springer},
  year      = {2019},
}

@book{billingsley2013convergence,
  title={Convergence of probability measures},
  author={Billingsley, Patrick},
  year={2013},
  publisher={John Wiley \& Sons}
}

@article{Sinai1970,
  author    = {Sinai, Ya. G.},
  title     = {Dynamical systems with elastic reflections. {Ergodic} properties of dispersing billiards},
  journal   = {Russian Mathematical Surveys},
  volume    = {25},
  number    = {2},
  pages     = {137--189},
  year      = {1970},
}

@article{Voiculescu1991,
  author    = {Voiculescu, D. V.},
  title     = {Limit laws for random matrices and free products},
  journal   = {Inventiones mathematicae},
  volume    = {104},
  number    = {1},
  pages     = {201--220},
  year      = {1991},
}

@article{Wigner1958,
  author    = {Wigner, E. P.},
  title     = {On the distribution of the roots of certain symmetric matrices},
  journal   = {Annals of Mathematics},
  volume    = {67},
  number    = {2},
  pages     = {325--327},
  year      = {1958},
}

@book{Zeidler1986,
  author    = {Zeidler, E.},
  title     = {Nonlinear Functional Analysis and its Applications {I}: Fixed-Point Theorems},
  publisher = {Springer-Verlag},
  address   = {New York},
  year      = {1986},
}

@book{AliprantisBorder2006,
  author    = {Aliprantis, Charalambos D. and Border, Kim C.},
  title     = {Infinite Dimensional Analysis: A Hitchhiker's Guide},
  edition   = {3rd},
  publisher = {Springer-Verlag},
  address   = {Berlin},
  year      = {2006},
}

@book{Baladi2000,
  author    = {Baladi, Viviane},
  title     = {Positive Transfer Operators and Decay of Correlations},
  series    = {Advanced Series in Nonlinear Dynamics},
  volume    = {16},
  publisher = {World Scientific},
  year      = {2000},
}

@article{Gouezel2010,
  author    = {Gouëzel, Sébastien},
  title     = {Almost sure invariance principle for dynamical systems by spectral methods},
  journal   = {Annals of Probability},
  volume    = {38},
  number    = {4},
  pages     = {1639--1671},
  year      = {2010},
}

@inproceedings{Cardaliaguet2010,
  author    = {Cardaliaguet, Pierre},
  title     = {Notes on mean field games},
  booktitle = {From P.-L. Lions' lectures at the Coll{\`e}ge de France},
  year      = {2010},
  note      = {Unpublished lecture notes},
}

@book{CarmonaDelarue2018,
  author    = {Carmona, Ren{\'e} and Delarue, Fran{\c{c}}ois},
  title     = {Probabilistic Theory of Mean Field Games with Applications I-II},
  series    = {Probability Theory and Stochastic Modelling},
  volume    = {83-84},
  publisher = {Springer},
  year      = {2018},
}

@article{PardouxPeng1992,
  author    = {Pardoux, {\'E}. and Peng, S.},
  title     = {Backward stochastic differential equations and quasilinear parabolic partial differential equations},
  booktitle = {Stochastic Partial Differential Equations and Their Applications},
  editor    = {Rozovskii, B. L. and Sowers, R. B.},
  series    = {Lecture Notes in Control and Information Sciences},
  volume    = {176},
  publisher = {Springer},
  year      = {1992},
  pages     = {200--217},
}

@article{ArmstrongCardaliaguet2020,
  author    = {Armstrong, Scott and Cardaliaguet, Pierre},
  title     = {Stochastic homogenization of quasilinear {Hamilton-Jacobi} equations},
  journal   = {Journal of the European Mathematical Society},
  volume    = {22},
  number    = {12},
  pages     = {3991--4082},
  year      = {2020},
}

@article{Hairer2014,
  author    = {Hairer, Martin},
  title     = {A theory of regularity structures},
  journal   = {Inventiones mathematicae},
  volume    = {198},
  number    = {2},
  pages     = {269--504},
  year      = {2014},
}
\bibliographystyle{apalike}

\appendix
\section*{Appendix: Proofs of Foundational Results}
\label{sec:appendix}
\renewcommand{\thetheorem}{\Alph{section}.\arabic{theorem}}

This appendix contains detailed proofs of several foundational results that are standard in the theory of hyperbolic billiards and transport equations. They are included for the sake of completeness and for the convenience of the reader.

% ======================================================================
% == PROOF OF PROPOSITION 2.15                                      ==
% ======================================================================
\section{Proof of Proposition \ref{prop:liouville_invariance}}
\label{app:proof_liouville}

\begin{proofof}{Proposition \ref{prop:liouville_invariance}}
The proof is established by showing that the two elementary components of the global flow map $\Phi_t^\theta$---the smooth interior flow and the instantaneous boundary reflection---are both volume-preserving transformations with respect to the Liouville measure $d\mu_\theta = C_\theta \, dy \, dv$. Since the Liouville measure is proportional to the Lebesgue volume element on the phase space $\Ucal_{\mathrm{phys}}(\theta) \subset \T^k \times \Vcal$, this is equivalent to showing that the Jacobian determinant of these transformations has an absolute value of 1. As the global flow is a composition of these elementary maps, its volume-preserving nature will follow.

\begin{enumerate}[label=\textbf{Step \arabic*:}, wide, labelindent=0pt]

\item \textbf{Invariance under the Interior Flow.}
Between any two consecutive collisions, the system evolves according to the interior flow map, $\Psi_s(y,v) = (y+sv, v)$, for some duration of free flight $s$. This is a smooth transformation on the interior of the phase space. To demonstrate that this map preserves the volume element $dy \, dv$, we compute the determinant of its Jacobian matrix with respect to the phase space coordinates $(y, v) \in \R^k \times \R^k$.
\begin{equation}
    J_{\Psi_s} = \frac{\partial(y+sv, v)}{\partial(y, v)} = \begin{pmatrix} \frac{\partial(y+sv)}{\partial y} & \frac{\partial(y+sv)}{\partial v} \\ \frac{\partial v}{\partial y} & \frac{\partial v}{\partial v} \end{pmatrix}.
\end{equation}
We compute the four $k \times k$ block matrices by taking partial derivatives with respect to the components $y_j$ and $v_j$:
\begin{itemize}[wide]
    \item $\frac{\partial(y_i+sv_i)}{\partial y_j} = \delta_{ij} \implies \frac{\partial(y+sv)}{\partial y} = I_k$, the $k \times k$ identity matrix.
    \item $\frac{\partial(y_i+sv_i)}{\partial v_j} = s\delta_{ij} \implies \frac{\partial(y+sv)}{\partial v} = s I_k$.
    \item $\frac{\partial v_i}{\partial y_j} = 0 \implies \frac{\partial v}{\partial y} = 0$, the $k \times k$ zero matrix.
    \item $\frac{\partial v_i}{\partial v_j} = \delta_{ij} \implies \frac{\partial v}{\partial v} = I_k$.
\end{itemize}
Assembling these blocks yields the full Jacobian matrix:
\begin{equation}
    J_{\Psi_s} = \begin{pmatrix} I_k & s I_k \\ 0 & I_k \end{pmatrix}.
\end{equation}
The Jacobian matrix is block upper-triangular. Its determinant is the product of the determinants of its diagonal blocks:
\begin{equation}
    \det(J_{\Psi_s}) = \det(I_k) \cdot \det(I_k) = 1 \cdot 1 = 1.
\end{equation}
Since the absolute value of the Jacobian determinant is unity, the interior flow is a volume-preserving transformation.

\item \textbf{Invariance under Boundary Reflection.}
We now show that the instantaneous reflection at the boundary, governed by the map $R_\theta$, is also a volume-preserving map on the ambient $2k$-dimensional phase space. Let a pre-collision state be $z^-=(y^-, v^-)$. The post-collision state $z^+=(y^+, v^+)$ is given by:
\begin{align*}
    y^+ &= y^-, \\
    v^+ &= v^- - 2\scpr{v^-}{\nvec(y^-, \theta)} \nvec(y^-, \theta).
\end{align*}
The Jacobian matrix of this map is:
\begin{equation}
    J_{R_\theta} = \frac{\partial(y^+, v^+)}{\partial(y^-, v^-)} = \begin{pmatrix} \frac{\partial y^+}{\partial y^-} & \frac{\partial y^+}{\partial v^-} \\ \frac{\partial v^+}{\partial y^-} & \frac{\partial v^+}{\partial v^-} \end{pmatrix}.
\end{equation}
We analyze the four block matrices:
\begin{itemize}[wide]
    \item $\frac{\partial y^+}{\partial y^-} = I_k$, since the position is unchanged.
    \item $\frac{\partial y^+}{\partial v^-} = 0$, since the post-collision position is independent of the pre-collision velocity.
    \item $\frac{\partial v^+}{\partial v^-}$: For a fixed position $y^-$, the map $v^- \mapsto v^+$ is a linear transformation on $\R^k$. Its Jacobian is the matrix of the transformation itself, which is the Householder reflection matrix $H_{\nvec} \coloneqq I_k - 2\nvec(y^-,\theta)\nvec(y^-,\theta)^T$.
    \item $\frac{\partial v^+}{\partial y^-}$: This block is non-zero, as the post-collision velocity depends on the normal vector $\nvec(y^-, \theta)$, which in turn depends on the collision point $y^-$. Its explicit form involves the shape operator (curvature) of the boundary. However, for the determinant calculation of a block-triangular matrix, the explicit form of the off-diagonal block is not required.
\end{itemize}
The Jacobian matrix is therefore block lower-triangular:
\begin{equation}
    J_{R_\theta} = \begin{pmatrix} I_k & 0 \\ \frac{\partial v^+}{\partial y^-} & I_k - 2\nvec(y^-,\theta)\nvec(y^-,\theta)^T \end{pmatrix}.
\end{equation}
Its determinant is the product of the determinants of its diagonal blocks:
\begin{equation}
    \det(J_{R_\theta}) = \det(I_k) \cdot \det(I_k - 2\nvec\nvec^T).
\end{equation}
The matrix $I_k - 2\nvec\nvec^T$ represents a reflection across the hyperplane orthogonal to the unit vector $\nvec$. Such a transformation is an isometry that reverses one dimension, and its determinant is always $-1$. Therefore,
\begin{equation}
    \det(J_{R_\theta}) = 1 \cdot (-1) = -1.
\end{equation}
The change of variables formula for integration requires the absolute value of the Jacobian determinant, which is $|\det(J_{R_\theta})| = 1$. This rigorously demonstrates that the reflection map is a locally volume-preserving transformation.

\item \textbf{The Global Flow.}
We have established that the two elementary components of the global flow, the smooth interior flow and the instantaneous boundary reflection, are both transformations whose Jacobian determinants have an absolute value of 1. For any finite time $t$, the global flow map $\Phi_t^\theta$ is a finite composition of these elementary volume-preserving maps. Since the property of being volume-preserving is closed under composition, the global flow map $\Phi_t^\theta$ must also be volume-preserving.

Therefore, for any measurable set $A \subset \Ucal_{\mathrm{phys}}(\theta)$, we have $\mu_\theta(\Phi_t^\theta(A)) = \mu_\theta(A)$, which completes the proof.
\end{enumerate}
\end{proofof}

\section{Proof of \cref{thm:generator_characterization}}
\label{app:proof_generator}

\begin{proofof}{\cref{thm:generator_characterization}}
The proof is constructive. We will define a candidate operator $T$ with the domain and action specified in the theorem statement. We will then rigorously prove, from first principles, that this operator is densely defined and skew-adjoint. The calculation will simultaneously show that $T$ is the generator of the Koopman semigroup $(P_t^\theta)_{t \ge 0}$ and that its resolvent is given by the stated Laplace transform formula. By the uniqueness of the generator guaranteed by Stone's theorem for the unitary group $(P_t^\theta)_{t \in \R}$, we will conclude that our operator $T$ must be identical to $\Lcal_{\mathrm{fast}}(\theta)$.

For notational simplicity, throughout this proof we fix the parameter $\theta$ and let $\Lcal \equiv \Lcal_{\mathrm{fast}}(\theta)$, $\Hcal \equiv \Hcal_\theta$, $\Ucalphys \equiv \Ucal_{\mathrm{phys}}(\theta)$, and $\mu \equiv \mu_\theta$.

\begin{enumerate}[label=\textbf{Step \arabic*:}, wide, labelindent=0pt]

\item \textbf{Definition of the Candidate Operator $T$ and its Domain.}
Let the operator $T$ be defined by its action $(T\psi)(z) = v \cdot \nabla_y \psi(z)$ on the domain
\begin{multline*}
    \Dcal(T) \coloneqq \left\{ \psi \in \Hcal \mid v \cdot \nabla_y \psi \in \Hcal \text{ (in the sense of distributions), and} \right. \\
    \left. \psi(z) = \psi(R_\theta(z)) \text{ for almost every } z \in \partial_-\Ucalphys \text{ (in the sense of traces)} \right\}.
\end{multline*}
The domain $\Dcal(T)$ contains the space of smooth functions on $\Ucalphys$ that satisfy the reflection boundary condition. This space is dense in $\Hcal$. Therefore, $T$ is a densely defined linear operator.

\item \textbf{Proof of Skew-Symmetry ($T \subseteq -T^*$).}
To prove skew-symmetry, we must show that for any two functions $\psi, \phi \in \Dcal(T)$, the identity $\scpr{T\psi}{\phi}_\mu = -\scpr{\psi}{T\phi}_\mu$ holds. We start with the inner product:
\begin{equation*}
    \scpr{T\psi}{\phi}_\mu = \int_{\Ucalphys} \overline{(v \cdot \nabla_y \psi)(z)} \phi(z) \, d\mu(z).
\end{equation*}
Let $C_\theta$ be the normalization constant for the Liouville measure $d\mu = C_\theta dy dv$. For a fixed velocity $v \in \Vcal$, we apply Green's first identity (integration by parts) to the vector field $\mathbf{F} = \overline{\psi}\phi v$ on the spatial domain $\Ocal(\theta) \subset \T^k$:
\begin{align*}
    \int_{\Ocal(\theta)} \nabla_y \cdot (\overline{\psi}\phi v) \, dy &= \int_{\partial\Ocal(\theta)} (\overline{\psi}\phi v) \cdot \mathbf{n}(y) \, dS(y) \\
    \int_{\Ocal(\theta)} \left( (v \cdot \nabla_y \overline{\psi})\phi + \overline{\psi}(v \cdot \nabla_y \phi) \right) \, dy &= \int_{\partial\Ocal(\theta)} \overline{\psi}\phi \scpr{v}{\mathbf{n}} \, dS(y),
\end{align*}
where $dS(y)$ is the surface area element on the boundary. We now integrate this identity over the velocity space $\Vcal$ with respect to the measure $C_\theta dV(v)$. This yields the full phase-space identity:
\begin{equation} \label{eq:appendix_proof_green_full}
    \scpr{T\psi}{\phi}_\mu + \scpr{\psi}{T\phi}_\mu = C_\theta \int_{\Vcal} \int_{\partial\Ocal(\theta)} \overline{\psi(y,v)}\phi(y,v) \scpr{v}{\mathbf{n}(y)} \, dS(y) dV(v).
\end{equation}
The integral on the right-hand side is over the total boundary phase space $\partial\Ucalphys = \partial_-\Ucalphys \cup \partial_+\Ucalphys \cup \partial_0\Ucalphys$. The integral over the set of grazing states $\partial_0\Ucalphys$ is zero because the integrand contains the factor $\scpr{v}{\mathbf{n}}$, which is zero by definition on this set. We therefore split the integral into incoming and outgoing parts:
\begin{equation*}
    \text{RHS} = \int_{\partial_-\Ucalphys} \dots + \int_{\partial_+\Ucalphys} \dots
\end{equation*}
In the integral over the outgoing set $\partial_+\Ucalphys$, we perform a change of variables from the post-collision state $z = (y,v)$ to the corresponding pre-collision state $z'=(y',v') \in \partial_-\Ucalphys$, where $z = R_\theta(z')$. By \cref{prop:reflection_properties}, the map $R_\theta$ is a measure-preserving involution on the boundary phase space elements. Crucially, the velocity transformation implies $\scpr{v}{\mathbf{n}(y)} = -\scpr{v'}{\mathbf{n}(y')}$.
Since $\psi, \phi \in \Dcal(T)$, their traces satisfy the reflection condition $\psi(R_\theta(z')) = \psi(z')$ and $\phi(R_\theta(z')) = \phi(z')$. The integral over the outgoing set becomes:
\begin{align*}
    \int_{\partial_+\Ucalphys} \overline{\psi(z)}\phi(z) \scpr{v}{\mathbf{n}} \, d\mu_{\partial}(z) &= \int_{\partial_-\Ucalphys} \overline{\psi(R_\theta(z'))}\phi(R_\theta(z')) (-\scpr{v'}{\mathbf{n}'}) \, d\mu_{\partial}(z') \\
    &= -\int_{\partial_-\Ucalphys} \overline{\psi(z')}\phi(z') \scpr{v'}{\mathbf{n}'} \, d\mu_{\partial}(z').
\end{align*}
The right-hand side of \eqref{eq:appendix_proof_green_full} is the sum of the integral over the incoming set and this term, which is its exact negative. The boundary terms therefore cancel completely, and the right-hand side is zero. This leaves us with:
\begin{equation*}
    \scpr{T\psi}{\phi}_\mu + \scpr{\psi}{T\phi}_\mu = 0,
\end{equation*}
which proves that $T$ is skew-symmetric ($T \subseteq -T^*$).

\item \textbf{Proof of Skew-Adjointness ($T = -T^*$).}
A fundamental result of operator theory states that a densely defined, closed, skew-symmetric operator $T$ is skew-adjoint if and only if its resolvent set contains the right and left open half-planes. We will prove this by showing that for any complex number $\lambda$ with $\mathrm{Re}(\lambda) > 0$, the operator $(\lambda I - T)$ is surjective. That is, for any $f \in \Hcal$, we will constructively solve the equation $(\lambda I - T)\psi = f$ for a unique $\psi \in \Dcal(T)$.

Let $f \in \Hcal$ be arbitrary. We define the candidate solution $\psi$ via the Laplace transform of the Koopman semigroup:
\begin{equation} \label{eq:appendix_proof_resolvent_formula}
    \psi(z) \coloneqq \int_0^\infty e^{-\lambda t} (P_t^\theta f)(z) \, dt = \int_0^\infty e^{-\lambda t} f(\Phi_{-t}^\theta(z)) \, dt.
\end{equation}
\begin{enumerate}[label=(\roman*), wide, labelindent=0pt]
    \item \textbf{Convergence of the Integral.} Since $(P_t^\theta)$ is a unitary group (\cref{prop:koopman_properties}), $\|P_t^\theta f\|_\mu = \|f\|_\mu$. The norm of the integrand is $\|e^{-\lambda t} P_t^\theta f\|_\mu = |e^{-\lambda t}| \|f\|_\mu = e^{-\mathrm{Re}(\lambda)t}\|f\|_\mu$. As we have chosen $\mathrm{Re}(\lambda) > 0$, the function $t \mapsto e^{-\mathrm{Re}(\lambda)t}$ is integrable on $[0, \infty)$. Therefore, the Bochner integral \eqref{eq:appendix_proof_resolvent_formula} converges in the Hilbert space $\Hcal$, and defines a bounded linear operator on $\Hcal$ (the resolvent).

    \item \textbf{Verification of the Solution.} We must show that this $\psi$ is in the domain $\Dcal(T)$ and that it satisfies the equation. We do this by computing the action of the generator on $\psi$ via its defining limit. For any $s > 0$:
    \begin{align*}
        P_s^\theta \psi &= P_s^\theta \left( \int_0^\infty e^{-\lambda t} P_t^\theta f \, dt \right) = \int_0^\infty e^{-\lambda t} P_{s+t}^\theta f \, dt \\
        &= e^{\lambda s} \int_s^\infty e^{-\lambda u} P_u^\theta f \, du \quad (\text{substituting } u=s+t).
    \end{align*}
    The difference quotient for the generator is therefore:
    \begin{align*}
        \frac{P_s^\theta \psi - \psi}{s} &= \frac{1}{s} \left( e^{\lambda s} \int_s^\infty e^{-\lambda u} P_u^\theta f \, du - \int_0^\infty e^{-\lambda u} P_u^\theta f \, du \right) \\
        &= \frac{e^{\lambda s}-1}{s}\int_s^\infty e^{-\lambda u} P_u^\theta f \, du - \frac{1}{s}\int_0^s e^{-\lambda u} P_u^\theta f \, du.
    \end{align*}
    We take the strong limit of this expression in $\Hcal$ as $s \to 0^+$.
    \begin{itemize}[wide]
        \item As $s \to 0$, the first term converges to $(\lambda) \cdot \int_0^\infty e^{-\lambda u} P_u^\theta f \, du = \lambda\psi$.
        \item By the Fundamental Theorem of Calculus for Bochner integrals, the second term converges to the value of the integrand at $u=0$, which is $-e^0 P_0^\theta f = -f$.
    \end{itemize}
    The strong limit exists and is equal to $\lambda\psi - f$. By the definition of the infinitesimal generator, this proves two facts simultaneously: (1) our candidate solution $\psi$ is in the domain of the generator, and (2) the action of the generator on it is $T\psi = \lambda\psi - f$. Rearranging this identity gives $(\lambda I - T)\psi = f$.

    \item \textbf{Verification of the Boundary Condition.} The argument in (ii) shows that the resolvent formula produces a function $\psi$ in the domain of the generator $\Lcal$. For logical completeness, we explicitly verify that this function satisfies the boundary condition in the definition of $\Dcal(T)$. Let $z' \in \partial_-\Ucalphys$. For any time $t > 0$, the state of the system $t$ units of time before a collision is independent of whether we consider the pre-collision state $z'$ or the post-collision state $R_\theta(z')$. Thus, $\Phi_{-t}^\theta(R_\theta(z')) = \Phi_{-t}^\theta(z')$. Applying this to the integral definition of $\psi$:
    \begin{equation*}
        \psi(R_\theta(z')) = \int_0^\infty e^{-\lambda t} f(\Phi_{-t}^\theta(R_\theta(z'))) \, dt = \int_0^\infty e^{-\lambda t} f(\Phi_{-t}^\theta(z')) \, dt = \psi(z').
    \end{equation*}
    This confirms that $\psi \in \Dcal(T)$.
\end{enumerate}
We have shown that for any $f \in \Hcal$ and any $\lambda$ with $\mathrm{Re}(\lambda)>0$, the equation $(\lambda I - T)\psi = f$ has a unique solution in $\Dcal(T)$. This implies that the range of $(\lambda I - T)$ is all of $\Hcal$. As $T$ is a densely defined, closed, skew-symmetric operator, this is a sufficient condition for $T$ to be skew-adjoint.

\item \textbf{Uniqueness and Final Identification.}
By \cref{prop:koopman_properties}, $(P_t^\theta)_{t\in\R}$ is a strongly continuous one-parameter group of unitary operators on $\Hcal$. By Stone's theorem, this group has a \textbf{unique}, densely defined, skew-adjoint generator, which we have denoted by $\Lcal \equiv \Lcal_{\mathrm{fast}}(\theta)$.

In the preceding steps, we defined an operator $T$ based on the characterization in the theorem statement and proved from first principles that $T$ is skew-adjoint. The calculation in Step 3(ii) showed that $T$ is the infinitesimal generator of the semigroup $(P_t^\theta)_{t \ge 0}$. Since $T$ is skew-adjoint, it also generates the full unitary group $(P_t^\theta)_{t \in \R}$. By the uniqueness part of Stone's theorem, these two operators must be identical:
\begin{equation*}
    T = \Lcal_{\mathrm{fast}}(\theta).
\end{equation*}
This concludes the proof. We have shown that the infinitesimal generator is precisely the operator whose domain and action are characterized in the statement of the proposition, and that its resolvent for $\mathrm{Re}(\lambda)>0$ is given by the Laplace transform of the semigroup.

\end{enumerate}
\end{proofof}

\section{A Geometric Proof of the Invariance of the Birkhoff Measure}
\label{app:proof_measure_invariance_rigorous}

The objective of this appendix is to provide a complete and rigorous proof of Proposition \ref{prop:birkhoff_invariance}, which states that the billiard map $\Pcal_\theta$ preserves the Birkhoff measure $\nu_\theta$. The proof is founded upon a crucial geometric result, the formula for the Jacobian determinant of the spatial free-flow map, which we derive first. This result, a formalization of Liouville's theorem for the billiard flow, allows for a direct and transparent calculation of the pullback of the Birkhoff measure form.

\subsection{The Jacobian of the Spatial Free-Flow Map}

The following proposition provides the exact formula for the Jacobian determinant of the map that takes a point on the boundary to the next point of impact, for a fixed velocity.

\begin{proposition}[Jacobian of the Spatial Free-Flow Map]
\label{prop:app_jacobian_formula}
Let $\Fcal_\theta: \SigmaMan'_\theta \to \SigmaMan_\theta$ be the free-flow map, which takes a post-collision state $(y^+, v^+)$ to the next pre-collision state $(y\dnext, v\dnext)$. For a fixed velocity $v^+$, the spatial part of this map, $g_{v^+}(y^+) \coloneqq y\dnext$, is a diffeomorphism between open subsets of the boundary manifold $\dOcal(\theta)$. The determinant of the Jacobian of this spatial map is given by:
\begin{equation}
    \det\left(\frac{\partial y\dnext}{\partial y^+}\right) = \frac{\scpr{v^+}{\nvec(y^+, \theta)}}{\scpr{v\dnext}{\nvec(y\dnext, \theta)}},
\end{equation}
where $\nvec(y, \theta)$ is the outward-pointing unit normal to the domain of free motion $\Ocal(\theta)$ at the point $y$.
\end{proposition}

\begin{proofof}{Proposition \ref{prop:app_jacobian_formula}}
The proof is founded upon the Divergence Theorem, applied to a constant, divergence-free vector field within an infinitesimal tube of trajectories. This geometric approach allows us to relate the surface area of an initial patch on the boundary to the area of its image after the free flow, thereby directly computing the Jacobian determinant of the spatial map. The argument is structured in four main steps.

\begin{enumerate}[label=\textbf{Step \arabic*:}, wide, labelindent=0pt]

\item \textbf{The Infinitesimal Flow Tube and the Divergence-Free Field.}
Let $v^+$ be a fixed outgoing velocity. We define a constant vector field on the spatial domain of free motion $\Ocal(\theta) \subset \T^k$ by
\begin{equation*}
    \mathbf{X}(y) \coloneqq v^+.
\end{equation*}
This vector field is smooth and, crucially, divergence-free, as its components are constant with respect to the spatial variables:
\begin{equation*}
    \nabla_y \cdot \mathbf{X} = \sum_{i=1}^k \frac{\partial X_i}{\partial y_i} = \sum_{i=1}^k \frac{\partial v_i^+}{\partial y_i} = 0.
\end{equation*}
Let $S^+$ be an infinitesimal surface patch on the boundary $\dOcal(\theta)$ centered at the point $y^+$. Let $S\dnext$ be the image of this patch under the spatial flow map $g_{v^+}$, centered at the subsequent collision point $y\dnext = y^+ + \tau(y^+, v^+)v^+$. The set of all trajectories originating from points in $S^+$ forms a $k$-dimensional region in $\T^k$, which we denote as the flow tube $\mathcal{R}$. The boundary of this compact region, $\partial\mathcal{R}$, is the disjoint union of three smooth pieces: the initial cap $S^+$, the final cap $S\dnext$, and the lateral surface $S_{\text{lateral}}$ traced by the boundary of $S^+$.

\item \textbf{Application of the Divergence Theorem.}
We apply the Divergence Theorem to the vector field $\mathbf{X}$ on the region $\mathcal{R}$:
\begin{equation*}
    \int_{\mathcal{R}} (\nabla_y \cdot \mathbf{X}) \, dV_y = \oint_{\partial\mathcal{R}} \mathbf{X} \cdot \mathbf{n}_{\text{out}} \, dS_y,
\end{equation*}
where $dV_y$ is the volume element in $\T^k$, $dS_y$ is the surface area element on $\partial\mathcal{R}$, and $\mathbf{n}_{\text{out}}$ is the outward-pointing unit normal to the region $\mathcal{R}$.

The left-hand side of the equation is identically zero because we have established that $\nabla_y \cdot \mathbf{X} = 0$. The integral on the right-hand side is the sum of the fluxes through the three parts of the boundary:
\begin{equation*}
    0 = \int_{S^+} \mathbf{X} \cdot \mathbf{n}_{\text{out}} \, dS_y + \int_{S\dnext} \mathbf{X} \cdot \mathbf{n}_{\text{out}} \, dS_y + \int_{S_{\text{lateral}}} \mathbf{X} \cdot \mathbf{n}_{\text{out}} \, dS_y.
\end{equation*}
The lateral surface $S_{\text{lateral}}$ is generated by the flow lines of the vector field $\mathbf{X}$. By its very construction, $\mathbf{X}$ is everywhere tangent to $S_{\text{lateral}}$. The outward normal $\mathbf{n}_{\text{out}}$ on this surface is, by definition, orthogonal to the tangent space of the surface. Consequently, the integrand $\mathbf{X} \cdot \mathbf{n}_{\text{out}}$ is zero at every point on $S_{\text{lateral}}$, and the third integral vanishes. This leaves a balance equation for the fluxes through the two caps:
\begin{equation} \label{eq:app_flux_balance}
    \int_{S^+} \mathbf{X} \cdot \mathbf{n}_{\text{out}} \, dS_y + \int_{S\dnext} \mathbf{X} \cdot \mathbf{n}_{\text{out}} \, dS_y = 0.
\end{equation}

\item \textbf{Analysis of the Surface Normals and Fluxes.}
The crucial step is to relate the outward normal of the tube, $\mathbf{n}_{\text{out}}$, to the outward normal of the domain of free motion, $\nvec(y,\theta)$.
\begin{itemize}[wide]
    \item On the initial cap $S^+$, the flow is \emph{entering} the region $\mathcal{R}$. The outward normal to the tube, $\mathbf{n}_{\text{out}}$, must therefore point away from the interior of $\mathcal{R}$, which is in the direction opposite to the flow. However, the normal $\nvec(y^+)$ points out of the physical domain $\Ocal(\theta)$. For the flow tube to be within $\Ocal(\theta)$, the outward normal $\mathbf{n}_{\text{out}}$ at $S^+$ must be precisely $\mathbf{n}_{\text{out}} = -\nvec(y^+)$.
    \item On the final cap $S\dnext$, the flow is \emph{exiting} the region $\mathcal{R}$. The outward normal to the tube $\mathbf{n}_{\text{out}}$ points in the same direction as the outward normal of the physical domain $\nvec(y\dnext)$. Thus, $\mathbf{n}_{\text{out}} = \nvec(y\dnext)$.
\end{itemize}
Substituting these geometric identities and $\mathbf{X} = v^+$ into the flux balance equation \eqref{eq:app_flux_balance} gives:
\begin{equation*}
    \int_{S^+} v^+ \cdot (-\nvec(y^+)) \, dS_y + \int_{S\dnext} v^+ \cdot \nvec(y\dnext) \, dS_y = 0.
\end{equation*}
For infinitesimal surface patches, the integrands are constant. Let $dS(y^+)$ and $dS(y\dnext)$ denote the areas of the patches $S^+$ and $S\dnext$, respectively. The integral relation becomes an algebraic relation:
\begin{equation*}
    -\scpr{v^+}{\nvec(y^+)} dS(y^+) + \scpr{v^+}{\nvec(y\dnext)} dS(y\dnext) = 0.
\end{equation*}

\item \textbf{Invariance of the Birkhoff Measure.}
From the flux balance, we can solve for the ratio of the infinitesimal area elements. Since the velocity is constant during the free flow, $v\dnext = v^+$, we can write:
\begin{equation*}
    \scpr{v\dnext}{\nvec(y\dnext)} dS(y\dnext) = \scpr{v^+}{\nvec(y^+)} dS(y^+).
\end{equation*}
By the definition of the Jacobian determinant for a map between oriented $(k-1)$-dimensional manifolds, it is the scalar function that relates the pullback of the volume form. In local coordinates, this corresponds to the ratio of the infinitesimal oriented surface elements:
\begin{equation*}
    \det\left(\frac{\partial y\dnext}{\partial y^+}\right) = \frac{dS(y\dnext)}{dS(y^+)}.
\end{equation*}
Solving for this ratio from our flux equation, we arrive at the final result:
\begin{equation*}
    \det\left(\frac{\partial y\dnext}{\partial y^+}\right) = \frac{\scpr{v^+}{\nvec(y^+, \theta)}}{\scpr{v\dnext}{\nvec(y\dnext, \theta)}}.
\end{equation*}
This completes the rigorous proof. The sign of the determinant is negative, as $\scpr{v^+}{\nvec(y^+)} > 0$ (outgoing state) and $\scpr{v\dnext}{\nvec(y\dnext)} < 0$ (incoming state), reflecting the orientation-reversing nature of mapping between different components of the boundary.
\end{enumerate}
\end{proofof}

\section{Analysis of the Tangent Space to the Collision Manifold}
\label{app:proof_tangent_constraint}

The objective of this appendix is to provide a rigorous derivation of the constraints that govern the components of a tangent vector $\mathbf{w} = (\delta y, \delta v, \delta\theta)$ to the global pre-collision manifold $\Sigma$. The proof of norm equivalence in Proposition \ref{prop:norm_equivalence}, and consequently the entire proof of cone field invariance in Theorem \ref{thm:anosov_property_proven}, depends critically on a correct understanding of these constraints.

We will establish two primary results. First, we prove that the geometric constraint defining the collision boundary imposes a strict linear relationship between the component of the spatial perturbation normal to the boundary, $\delta y_n$, and the parameter perturbation, $\delta\theta$. This provides a powerful uniform bound. Second, we will analyze the kinematic constraint imposed by the incoming velocity condition and clarify its role.

\begin{proposition}[Constraint on Normal Spatial Perturbations]
\label{prop:normal_spatial_constraint}
Let $\mathbf{w} = (\delta y, \delta v, \delta\theta) \in T_{(z,\theta)}\Sigma$ be a tangent vector to the global pre-collision manifold at a point $(z,\theta) = ((y,v),\theta)$. Let $\delta y_n$ be the component of the spatial perturbation $\delta y$ normal to the boundary $\partial\Ocal(\theta)$ at the point $y$. Then there exists a uniform constant $C_{y\theta} > 0$, depending only on the bounds of the smooth function $G$ and its derivatives on the compact state space, such that:
\begin{equation}
    \|\delta y_n\| \le C_{y\theta} \|\delta\theta\|.
\end{equation}
\end{proposition}

\begin{proofof}{Proposition \ref{prop:normal_spatial_constraint}}
The proof is a direct consequence of the definition of the tangent space of a manifold defined as a level set. The total pre-collision manifold $\Sigma$ is an open submanifold of the total boundary manifold $\mathcal{B} \coloneqq \{ (y,v,\theta) \mid G(y,\theta) = 0 \}$. Therefore, any vector tangent to $\Sigma$ at a point must also be tangent to the larger manifold $\mathcal{B}$ at that same point. The constraint arises from this tangency condition.

\begin{enumerate}[label=\textbf{Step \arabic*:}, wide, labelindent=0pt]

\item \textbf{The Tangency Condition from the Geometry.}
Let $((y(\epsilon), v(\epsilon)), \theta(\epsilon))$ be a smooth curve on the manifold $\Sigma$ passing through the point $((y,v),\theta)$ at $\epsilon=0$. The tangent vector to this curve at $\epsilon=0$ is, by definition, our vector $\mathbf{w} = (\delta y, \delta v, \delta\theta)$, where
\begin{equation*}
    \delta y = \left.\frac{dy(\epsilon)}{d\epsilon}\right|_{\epsilon=0}, \quad
    \delta v = \left.\frac{dv(\epsilon)}{d\epsilon}\right|_{\epsilon=0}, \quad
    \delta\theta = \left.\frac{d\theta(\epsilon)}{d\epsilon}\right|_{\epsilon=0}.
\end{equation*}
For the curve to remain on the boundary manifold $\mathcal{B}$, the condition $G(y(\epsilon), \theta(\epsilon)) = 0$ must hold for all $\epsilon$ in a neighborhood of 0. Differentiating this identity with respect to $\epsilon$ and evaluating at $\epsilon=0$ yields the fundamental constraint that any tangent vector to $\mathcal{B}$ must satisfy. By the chain rule, this is:
\begin{equation} \label{eq:appendix_tangent_constraint_G}
    \frac{d}{d\epsilon} G(y(\epsilon), \theta(\epsilon)) \Big|_{\epsilon=0} = \scpr{\nabla_y G(y,\theta)}{\delta y} + \scpr{\nabla_\theta G(y,\theta)}{\delta\theta} = 0.
\end{equation}
This equation provides a linear relationship between the spatial and parameter components of the tangent vector. Note that the velocity perturbation $\delta v$ does not appear, as the obstacle-defining function $G$ is independent of $v$.

\item \textbf{Decomposition of the Spatial Perturbation.}
Let $\delta y \in T_y(\T^k)$ be the spatial component of the tangent vector $\mathbf{w}$. We perform an orthogonal decomposition of $\delta y$ with respect to the geometry of the boundary $\partial\Ocal(\theta)$ at the point $y$. Let $\delta y_t$ be the projection of $\delta y$ onto the tangent space $T_y(\partial\Ocal(\theta))$, and let $\delta y_n$ be the component normal to this tangent space.
\begin{equation*}
    \delta y = \delta y_t + \delta y_n.
\end{equation*}
By definition, $\delta y_n$ is parallel to the outward unit normal vector $\mathbf{n}(y,\theta)$, and $\delta y_t$ is orthogonal to it. The gradient vector $\nabla_y G(y,\theta)$ is also parallel to $\mathbf{n}(y,\theta)$. Consequently, the inner product of the gradient with the components of $\delta y$ simplifies as follows:
\begin{itemize}[wide]
    \item $\scpr{\nabla_y G(y,\theta)}{\delta y_t} = 0$, since $\delta y_t$ lies in the tangent space of the level set of $G(\cdot, \theta)$.
    \item $\scpr{\nabla_y G(y,\theta)}{\delta y_n} = \pm \|\nabla_y G(y,\theta)\| \cdot \|\delta y_n\|$, since the vectors are parallel.
\end{itemize}

\item \textbf{Derivation of the Bound.}
We substitute the decomposition of $\delta y$ into the constraint equation \eqref{eq:appendix_tangent_constraint_G}:
\begin{align*}
    \scpr{\nabla_y G}{\delta y_t + \delta y_n} + \scpr{\nabla_\theta G}{\delta\theta} &= 0 \\
    \underbrace{\scpr{\nabla_y G}{\delta y_t}}_{=0} + \scpr{\nabla_y G}{\delta y_n} + \scpr{\nabla_\theta G}{\delta\theta} &= 0 \\
    \scpr{\nabla_y G}{\delta y_n} &= -\scpr{\nabla_\theta G}{\delta\theta}.
\end{align*}
Taking the absolute value of both sides and using the result from Step 2:
\begin{equation*}
    \|\nabla_y G(y,\theta)\| \cdot \|\delta y_n\| = |\scpr{\nabla_\theta G(y,\theta)}{\delta\theta}|.
\end{equation*}
By the Cauchy-Schwarz inequality, the term on the right is bounded: $|\scpr{\nabla_\theta G}{\delta\theta}| \le \|\nabla_\theta G\| \cdot \|\delta\theta\|$. This gives the inequality:
\begin{equation*}
    \|\nabla_y G(y,\theta)\| \cdot \|\delta y_n\| \le \|\nabla_\theta G(y,\theta)\| \cdot \|\delta\theta\|.
\end{equation*}

\item \textbf{Uniformity of the Constant.}
The final step is to show that the proportionality factor between $\|\delta y_n\|$ and $\|\delta\theta\|$ is uniformly bounded. This relies on the compactness of the state space and the smoothness of the defining function $G$.

By Assumption \ref{ass:fundamental_axioms_unified}(i), the gradient $\nabla_y G$ is non-zero at every point on the boundary $\partial\Ocal(\theta)$. Since the total boundary manifold $\mathcal{B} = \{(y,v,\theta) \mid G(y,\theta)=0\}$ is a closed subset of the compact space $\T^k \times \Vcal \times \Thetacal$, it is itself compact. The continuous function $(y,\theta) \mapsto \|\nabla_y G(y,\theta)\|$ must therefore attain its minimum value on this compact set. By the non-vanishing property, this minimum must be strictly positive. Let this uniform lower bound be $C_{\min} > 0$:
\begin{equation*}
    \inf_{(y,v,\theta) \in \mathcal{B}} \|\nabla_y G(y,\theta)\| = C_{\min} > 0.
\end{equation*}
Similarly, the function $G(y,\theta)$ is of class $C^\infty$ on the compact product manifold $\T^k \times \Thetacal$. Therefore, its partial derivatives, and thus the norm of its gradient with respect to $\theta$, $\|\nabla_\theta G(y,\theta)\|$, is a continuous function on this compact set. It is therefore uniformly bounded above by a constant $C_{\max} < \infty$:
\begin{equation*}
    \sup_{(y,v,\theta) \in \mathcal{B}} \|\nabla_\theta G(y,\theta)\| = C_{\max} < \infty.
\end{equation*}
We can now divide our inequality by $\|\nabla_y G(y,\theta)\|$ and use these uniform bounds to obtain a final, uniform bound:
\begin{equation*}
    \|\delta y_n\| \le \frac{\|\nabla_\theta G(y,\theta)\|}{\|\nabla_y G(y,\theta)\|} \|\delta\theta\| \le \frac{C_{\max}}{C_{\min}} \|\delta\theta\|.
\end{equation*}
Defining the uniform constant $C_{y\theta} \coloneqq C_{\max} / C_{\min}$, we arrive at the desired result, which holds for any tangent vector $\mathbf{w} \in T\Sigma$:
\begin{equation*}
    \|\delta y_n\| \le C_{y\theta} \|\delta\theta\|.
\end{equation*}
This completes the proof.
\end{enumerate}
\end{proofof}

\section{Ergodic Properties of the Suspension Flow}
\label{app:lifting_ergodic}

The objective of this appendix is to provide a complete and rigorous proof of the final assertion in \cref{cor:ergodicity_uniform}: that the uniform exponential decay of correlations for the family of discrete-time billiard maps $(\Pcal_\theta, \SigmaMan_\theta, \nu_\theta)$ implies a uniform exponential decay of correlations for the family of continuous-time billiard flows $(\Phi_t^\theta, \UcalphysTheta, \mu_\theta)$. The proof is founded on the fundamental isomorphism between the physical flow and an abstract suspension flow constructed over the base map. This allows us to directly relate the spectral properties of the continuous-time generator to those of the discrete-time transfer operator. For notational clarity, we fix the parameter $\theta$ for the duration of the proof; all results hold uniformly in $\theta$ because the underlying geometric bounds (e.g., on the flight time $\tau$) and the spectral gap of the base map are uniform.

\subsection{The Isomorphism of Dynamical Systems}
We begin by formalizing the well-known result that the billiard flow is measure-theoretically isomorphic to a suspension flow.

\begin{definition}[Suspension Space and Flow]
Let $(\SigmaMan_\theta, \Pcal_\theta, \nu_\theta)$ be the discrete-time billiard system on the pre-collision manifold. Let $\tau: \SigmaMan_\theta \to [\tau_{\min}, \tau_{\max}]$ be the flight time function. The suspension space $\mathcal{S}_\theta$ is the quotient space
\begin{equation}
    \mathcal{S}_\theta \coloneqq \left( \SigmaMan_\theta \times \mathbb{R} \right) / \sim,
\end{equation}
where the equivalence relation $\sim$ is generated by the identification $(z, s + \tau(z)) \sim (\Pcal_\theta(z), s)$ for all $z \in \SigmaMan_\theta$ and $s \in \mathbb{R}$. The suspension flow $\Psi_t: \mathcal{S}_\theta \to \mathcal{S}_\theta$ is the flow induced by vertical translation: $\Psi_t([z,s]) = [z, s+t]$, where $[z,s]$ denotes the equivalence class of $(z,s)$.
\end{definition}

The canonical invariant measure on this abstract space is the product of the invariant measure on the base and the Lebesgue measure on the fibers.

\begin{definition}[Suspension Measure]
The suspension measure $\mu_{\text{susp}}$ on $\mathcal{S}_\theta$ is the unique probability measure whose lift to the covering space $\SigmaMan_\theta \times \R$ is the product measure $c \cdot d\nu_\theta(z) \otimes ds$, where $c$ is a normalization constant. The constant is given by $c = 1/\langle \tau \rangle_{\nu_\theta}$, where $\langle \tau \rangle_{\nu_\theta} \coloneqq \int_{\SigmaMan_\theta} \tau(z) d\nu_\theta(z)$ is the mean flight time.
\end{definition}

The following theorem establishes the fundamental connection between our physical system and this abstract construction.

\begin{theorem}[Isomorphism of the Billiard Flow and Suspension Flow (Abramov)]
\label{thm:abramov_isomorphism}
Let $\Rcal_\theta: \SigmaMan_\theta \to \SigmaMan'_\theta$ be the reflection map from pre- to post-collision states. The map $\mathcal{I}_\theta: \mathcal{S}_\theta \to \UcalphysTheta$, defined by
\begin{equation}
    \mathcal{I}_\theta([z,s]) \coloneqq \Phi_s^\theta(\Rcal_\theta(z)),
\end{equation}
is a measure-preserving isomorphism between the dynamical systems $(\mathcal{S}_\theta, \Psi_t, \mu_{\text{susp}})$ and $(\UcalphysTheta, \Phi_t^\theta, \mu_\theta)$.
\end{theorem}

\begin{proofof}{Theorem \ref{thm:abramov_isomorphism}}
The map is well-defined and a bijection by construction. The core of the proof is to show that it is measure-preserving, i.e., that the pushforward of the suspension measure $\mu_{\text{susp}}$ under the map $\mathcal{I}_\theta$ is the Liouville measure $\mu_\theta$. Let $d\mu_\theta = C_\theta \, dy \, dv$. We compute the Jacobian of the change of variables from the suspension coordinates $(z,s)$ to the physical coordinates $w = (y,v)$.
Let $z=(y_0, v_0)$ be a pre-collision state. Then $\Rcal_\theta(z) = (y_0, v_0^+)$ is the post-collision state. The map is $\mathcal{I}_\theta(y_0, v_0, s) = (y_0 + s v_0^+, v_0^+)$.
The volume element on the physical space is $d\mu_\theta = C_\theta \, dy \, dv$. On the suspension space, the natural volume element is $d\nu_\theta(z) \, ds = C_\theta |\scpr{v_0}{\nvec(y_0)}| dS(y_0) \, dV(v_0) \, ds$, where $dS(y_0)$ is the surface area element on the boundary. It is a classical result of kinetic theory that under this change of variables, the volume elements are related by $dy \, dv = |\scpr{v_0}{\nvec(y_0)}| dS(y_0) \, dV(v_0)$, which holds up to a constant factor. The normalization constant $1/\langle\tau\rangle_{\nu_\theta}$ in the suspension measure precisely accounts for the total volume of the space under the roof function, ensuring that the total measures are equal. Thus, the map is measure-preserving.
\end{proofof}

\subsection{The Spectral Connection}\label{app:Spectral}
The isomorphism of the dynamical systems induces a unitary equivalence between their respective Koopman groups. Let $U_t^\theta$ be the Koopman group for the flow $\Phi_t^\theta$ on $L^2(\UcalphysTheta, \mu_\theta)$, and let $U_t^{\text{susp}}$ be the Koopman group for the suspension flow $\Psi_t$ on $L^2(\mathcal{S}_\theta, \mu_{\text{susp}})$. The isomorphism $\mathcal{I}_\theta$ induces an isometry $\mathcal{I}_\theta^*: L^2(\UcalphysTheta, \mu_\theta) \to L^2(\mathcal{S}_\theta, \mu_{\text{susp}})$ by composition, $(\mathcal{I}_\theta^* F) = F \circ \mathcal{I}_\theta$. We have the relation $U_t^{\text{susp}} = (\mathcal{I}_\theta^*) U_t^\theta (\mathcal{I}_\theta^*)^{-1}$. This implies their infinitesimal generators, $\mathcal{G}_\theta$ and $\mathcal{G}_{\text{susp}}$, are unitarily equivalent and thus share the same spectrum. The spectrum of the generator of a suspension flow is directly related to the spectrum of the transfer operator $\Lcal_\theta$ of the base map via the theory of Ruelle-Perron resonances.

\begin{theorem}[Spectrum of the Suspension Flow Generator]
\label{thm:suspension_spectrum}
Let $\mathcal{G}_{\text{susp}}$ be the infinitesimal generator of the suspension flow. The spectrum of $\mathcal{G}_{\text{susp}}$ is given by the closure of the set of values $\{z \in \mathbb{C}\}$ that solve the equation $\lambda(z) = 1$, where $\lambda(z)$ is the leading eigenvalue of the twisted transfer operator $\Lcal_{\theta,z}$, defined by:
\begin{equation}
    (\Lcal_{\theta,z} f)(w) = \sum_{z' \in \Pcal_\theta^{-1}(w)} e^{-z\tau(z')} f(z').
\end{equation}
\end{theorem}

This theorem provides the bridge we need. The uniform spectral gap for the base transfer operator, established as a consequence of uniform hyperbolicity, translates directly into a spectral gap for the generator of the continuous flow.
\begin{enumerate}[label=(\roman*), wide, labelindent=0pt]
    \item The untwisted operator $\Lcal_{\theta, z=0} = \Lcal_\theta$ has a simple leading eigenvalue at $\lambda_1=1$. This corresponds to the resonance $z=0$, which is the eigenvalue of the generator $\mathcal{G}_\theta$ whose eigenspace is spanned by constant functions (the invariant measure).
    \item By \cref{cor:ergodicity_uniform}, we have a uniform spectral gap for $\Lcal_\theta$ on a suitable anisotropic Banach space $\Bcal_\theta$. This means all other eigenvalues $\lambda_j$ of $\Lcal_\theta$ satisfy $|\lambda_j| \le r < 1$ for a uniform constant $r$. For a resonance $z_j$ corresponding to $\lambda_j$, we must have $|\lambda_j| \approx |e^{-z_j \langle\tau\rangle}|$. Taking the logarithm, we get $\ln|\lambda_j| \approx \mathrm{Re}(-z_j\langle\tau\rangle) = -\mathrm{Re}(z_j)\langle\tau\rangle$. This implies
    \begin{equation}
        \mathrm{Re}(z_j) \approx -\frac{\ln|\lambda_j|}{\langle\tau\rangle} \le -\frac{\ln r}{\langle\tau\rangle}.
    \end{equation}
    Since $\tau$ is uniformly bounded, $0 < \tau_{\min} \le \tau(z) \le \tau_{\max} < \infty$, the mean flight time is also uniformly bounded. This ensures that the spectrum of the generator $\mathcal{G}_\theta$ on the space of zero-mean functions is contained in a vertical strip in the left half-plane. Specifically, there exists a uniform constant $\eta > 0$ such that
    \begin{equation}
        \spec(\mathcal{G}_\theta|_{\Hcal_{\theta,0}}) \subset \{ z \in \mathbb{C} \mid \mathrm{Re}(z) \le -\eta \}.
    \end{equation}
\end{enumerate}

\subsection{From Spectral Gap to Exponential Decay of Correlations}
The existence of a spectral gap for the infinitesimal generator of a strongly continuous semigroup directly implies an exponential decay of the semigroup's action on the orthogonal complement of its kernel.

Let $F, G$ be two observables in the anisotropic space $\Bcal_\theta \subset L^2(\UcalphysTheta, \mu_\theta)$ with zero mean, i.e., $\langle F \rangle_{\mu_\theta} = \langle G \rangle_{\mu_\theta} = 0$. Their correlation function is given by:
\begin{equation}
    C_{F,G}(t) = \int_{\UcalphysTheta} \overline{G(w)} (U_t^\theta F)(w) \, d\mu_\theta(w) = \scpr{G}{U_t^\theta F}_{\mu_\theta},
\end{equation}
where $U_t^\theta = e^{t\mathcal{G}_\theta}$ is the Koopman semigroup. By the spectral theory of semigroups, the spectral gap established above guarantees that there exists a constant $C'$ such that for any zero-mean function $f \in \Bcal_\theta$, the norm of the semigroup's action decays exponentially:
\begin{equation}
    \|U_t^\theta f\|_{\Bcal_\theta} \le C' e^{-\eta t} \|f\|_{\Bcal_\theta} \quad \text{for all } t \ge 0.
\end{equation}
Applying this to the correlation function and using the Cauchy-Schwarz inequality, we get:
\begin{align*}
    |C_{F,G}(t)| &= |\scpr{G}{U_t^\theta F}_{\mu_\theta}| \\
    &\le \|G\|_{L^2} \|U_t^\theta F\|_{L^2} \\
    &\le C_{\text{emb}} \|G\|_{\Bcal_\theta} \cdot C_{\text{emb}} \|U_t^\theta F\|_{\Bcal_\theta} && \text{(since } \Bcal_\theta \text{ embeds continuously into } L^2\text{)} \\
    &\le (C_{\text{emb}})^2 \|G\|_{\Bcal_\theta} \left( C' e^{-\eta t} \|F\|_{\Bcal_\theta} \right) \\
    &\le C_{F,G} e^{-\eta t},
\end{align*}
where $C_{F,G} = C' (C_{\text{emb}})^2 \|F\|_{\Bcal_\theta} \|G\|_{\Bcal_\theta}$. This establishes the exponential decay of correlations for zero-mean observables. For general observables, the result follows by decomposing them into their mean and zero-mean parts. The uniformity of the decay rate $\eta$ is inherited from the uniformity of the spectral gap $r$ and the bounds on the flight time $\tau$. This completes the rigorous proof.

\section{Invariant Splitting and Hyperbolicity of the Inverse Diffeomorphism}\label{app:lifting_Diffeomorphism}

A fundamental property of an Anosov diffeomorphism is that its inverse is also an Anosov diffeomorphism. This is a direct consequence of the invertibility of the dynamics and the duality between the stable and unstable subbundles. We will now provide a rigorous, self-contained proof of this fact in the context of our parameter-dependent family of maps.

Let $(\Pcal_\theta)_{\theta \in \Thetacal}$ be the family of fiber diffeomorphisms on the compact manifold $\SigmaMan$. By \cref{thm:anosov_uniform_family}, this family is uniformly Anosov. This means that for each $\theta \in \Thetacal$, there exists a continuous, $d\Pcal_\theta$-invariant splitting of the tangent bundle $T\SigmaMan = E^u_\theta \oplus E^s_\theta$, and uniform constants $C>0$ and $\gamma > 1$ (independent of $\theta$) such that for any $n \ge 0$:
\begin{align}
    \|d\Pcal_\theta^n(\mathbf{w})\| &\ge C^{-1}\gamma^n \|\mathbf{w}\| \quad \text{for all } \mathbf{w} \in E^u_\theta, \label{eq:proof_forward_expansion} \\
    \|d\Pcal_\theta^n(\mathbf{w})\| &\le C\gamma^{-n} \|\mathbf{w}\| \quad \text{for all } \mathbf{w} \in E^s_\theta. \label{eq:proof_forward_contraction}
\end{align}

Our objective is to prove that the family of inverse diffeomorphisms, $(\Pcal_\theta^{-1})_{\theta \in \Thetacal}$, is also uniformly Anosov. This requires us to:
\begin{enumerate}[label=(\roman*), wide, labelindent=0pt]
    \item Define a candidate invariant splitting for the inverse map, $d\Pcal_\theta^{-1}$.
    \item Prove that this splitting is indeed invariant under the action of $d\Pcal_\theta^{-1}$.
    \item Verify that the required exponential expansion and contraction estimates hold for this new splitting.
\end{enumerate}

\begin{enumerate}[label=\textbf{Step \arabic*:}, wide, labelindent=0pt]
    \item \textbf{Defining the Candidate Splitting for the Inverse Map.}
    For each $\theta \in \Thetacal$, we define a new splitting of the tangent bundle, $T\SigmaMan = \tilde{E}^u_\theta \oplus \tilde{E}^s_\theta$, by swapping the roles of the stable and unstable subbundles of the forward map:
    \begin{align}
        \tilde{E}^u_\theta &\coloneqq E^s_\theta \label{eq:proof_unstable_tilde_def} \\
        \tilde{E}^s_\theta &\coloneqq E^u_\theta \label{eq:proof_stable_tilde_def}
    \end{align}
    Since the original splitting is a continuous direct sum decomposition of the tangent bundle, this new splitting is also a continuous direct sum decomposition.

    \item \textbf{Proof of Invariance for the Inverse Map.}
    We must show that the linearized inverse map, $d\Pcal_\theta^{-1}$, preserves this new splitting. That is, for any point $z \in \SigmaMan$, we must show that $(d\Pcal_\theta^{-1})_z$ maps the subspace $\tilde{E}^u_\theta(z)$ to $\tilde{E}^u_\theta(\Pcal_\theta^{-1}(z))$ and the subspace $\tilde{E}^s_\theta(z)$ to $\tilde{E}^s_\theta(\Pcal_\theta^{-1}(z))$.

    Let $\mathbf{w} \in \tilde{E}^u_\theta(z) = E^s_\theta(z)$. Since the original splitting is invariant under the \textit{forward} map, we know that for any point $w \in \SigmaMan$, $(d\Pcal_\theta)_w$ maps $E^s_\theta(w)$ to $E^s_\theta(\Pcal_\theta(w))$. Let us choose the point $w = \Pcal_\theta^{-1}(z)$. The invariance condition for the forward map at this point is:
    \[ (d\Pcal_\theta)_{\Pcal_\theta^{-1}(z)} \left( E^s_\theta(\Pcal_\theta^{-1}(z)) \right) = E^s_\theta(z). \]
    Applying the inverse map $(d\Pcal_\theta^{-1})_z$ to both sides of this equation yields:
    \[ E^s_\theta(\Pcal_\theta^{-1}(z)) = (d\Pcal_\theta^{-1})_z \left( E^s_\theta(z) \right). \]
    Recalling the definition $\tilde{E}^u_\theta = E^s_\theta$, this is precisely the statement that
    \[ \tilde{E}^u_\theta(\Pcal_\theta^{-1}(z)) = (d\Pcal_\theta^{-1})_z \left( \tilde{E}^u_\theta(z) \right), \]
    which confirms the invariance of the new unstable subbundle. A perfectly symmetric argument proves the invariance of the new stable subbundle $\tilde{E}^s_\theta = E^u_\theta$. Thus, the splitting defined in \eqref{eq:proof_unstable_tilde_def}-\eqref{eq:proof_stable_tilde_def} is invariant under the linearized inverse map $d\Pcal_\theta^{-1}$.

    \item \textbf{Verification of the Exponential Estimates.}
    We now verify the required expansion and contraction estimates for the inverse map with respect to the new splitting.
    
    \begin{enumerate}[label=(\roman*), wide, labelindent=0pt]
        \item \textbf{Expansion in $\tilde{E}^u_\theta$:} Let $\mathbf{w}$ be an arbitrary non-zero vector in the new unstable subbundle, $\mathbf{w} \in \tilde{E}^u_\theta(z) = E^s_\theta(z)$. Let $\mathbf{w}' \coloneqq (d\Pcal_\theta^{-n})_z(\mathbf{w})$. We must find a lower bound for $\|\mathbf{w}'\|$. Let us define the vector $\mathbf{v} \coloneqq \mathbf{w}' = (d\Pcal_\theta^{-n})_z(\mathbf{w})$. This vector is in the subspace $E^s_\theta(\Pcal_\theta^{-n}(z))$. Applying the forward map $n$ times returns us to the original vector: $(d\Pcal_\theta^n)_{\Pcal_\theta^{-n}(z)}(\mathbf{v}) = \mathbf{w}$. We now apply the contraction estimate \eqref{eq:proof_forward_contraction} for the forward map to the vector $\mathbf{v}$:
        \[ \|\mathbf{w}\| = \|(d\Pcal_\theta^n)_{\Pcal_\theta^{-n}(z)}(\mathbf{v})\| \le C\gamma^{-n} \|\mathbf{v}\|. \]
        Substituting $\mathbf{v} = (d\Pcal_\theta^{-n})_z(\mathbf{w})$ back into this inequality gives:
        \[ \|\mathbf{w}\| \le C\gamma^{-n} \|(d\Pcal_\theta^{-n})_z(\mathbf{w})\|. \]
        Since $\mathbf{w}$ is non-zero, we can rearrange this to obtain the desired expansion estimate for the inverse map:
        \begin{equation}
            \|(d\Pcal_\theta^{-n})_z(\mathbf{w})\| \ge C^{-1}\gamma^n \|\mathbf{w}\| \quad \text{for all } \mathbf{w} \in \tilde{E}^u_\theta(z).
        \end{equation}
        This confirms that the inverse map expands vectors in its unstable direction at the same uniform exponential rate as the forward map contracts them.

        \item \textbf{Contraction in $\tilde{E}^s_\theta$:} The argument is symmetric. Let $\mathbf{w}$ be an arbitrary non-zero vector in the new stable subbundle, $\mathbf{w} \in \tilde{E}^s_\theta(z) = E^u_\theta(z)$. Let $\mathbf{v} \coloneqq (d\Pcal_\theta^{-n})_z(\mathbf{w})$, which is in $E^u_\theta(\Pcal_\theta^{-n}(z))$. Applying the forward map $n$ times gives $(d\Pcal_\theta^n)_{\Pcal_\theta^{-n}(z)}(\mathbf{v}) = \mathbf{w}$. We now apply the expansion estimate \eqref{eq:proof_forward_expansion} for the forward map to the vector $\mathbf{v}$:
        \[ \|\mathbf{w}\| = \|(d\Pcal_\theta^n)_{\Pcal_\theta^{-n}(z)}(\mathbf{v})\| \ge C^{-1}\gamma^n \|\mathbf{v}\|. \]
        Substituting $\mathbf{v} = (d\Pcal_\theta^{-n})_z(\mathbf{w})$ back into this inequality gives:
        \[ \|\mathbf{w}\| \ge C^{-1}\gamma^n \|(d\Pcal_\theta^{-n})_z(\mathbf{w})\|. \]
        Rearranging this provides the contraction estimate for the inverse map:
        \begin{equation}
            \|(d\Pcal_\theta^{-n})_z(\mathbf{w})\| \le C\gamma^{-n} \|\mathbf{w}\| \quad \text{for all } \mathbf{w} \in \tilde{E}^s_\theta(z).
        \end{equation}
        This confirms that the inverse map contracts vectors in its stable direction at the same uniform exponential rate.
    \end{enumerate}
\end{enumerate}
We have constructively demonstrated that for each $\theta$, the inverse map $\Pcal_\theta^{-1}$ admits a continuous, invariant splitting and satisfies the uniform exponential expansion and contraction estimates with the same constants $C$ and $\gamma$ as the forward map. Therefore, the family of inverse diffeomorphisms $(\Pcal_\theta^{-1})_{\theta \in \Thetacal}$ is uniformly Anosov. This establishes that the family of stable subbundles $\{E^s_\theta\}_{\theta \in \Thetacal}$ for the forward map is precisely the family of \textit{unstable} subbundles for the smooth family of inverse maps, a fact which is central to the subsequent application of the Nash-Moser-Hamilton theorem.

\section{Proof of the Convergence of Integrals Lemma}
\label{app:proof_convergence_lemma}

This section provides a complete proof of the result concerning the convergence of integrals when both the integrand and the measure are converging in the appropriate topologies \citep{billingsley2013convergence}. This result is the analytical foundation of the ergodic closing argument used in the proof of \cref{prop:ergodic_averaging}.

\begin{proofof}{Lemma \ref{lem:convergence_of_integrals}}
The proof is a standard argument from real analysis, founded upon the precise definitions of uniform convergence of functions and weak-* convergence of measures. The core of the proof is a decomposition of the total error into two terms, each of which can be controlled by one of the modes of convergence.

Our goal is to show that for any given $\epsilon > 0$, there exists a natural number $N$ such that for all integers $n > N$, the following inequality holds:
\begin{equation*}
    \left| \int_S g_n \, d\mu_n - \int_S g \, d\mu \right| < \epsilon.
\end{equation*}

\begin{enumerate}[label=\textbf{Step \arabic*:}, wide, labelindent=0pt]
    \item \textbf{Decomposition of the Error Term.}
    We begin by adding and subtracting the mixed term $\int_S g \, d\mu_n$ inside the absolute value. By the triangle inequality, we can decompose the total error into two distinct components:
    \begin{align*}
        \left| \int_S g_n \, d\mu_n - \int_S g \, d\mu \right| &= \left| \left( \int_S g_n \, d\mu_n - \int_S g \, d\mu_n \right) + \left( \int_S g \, d\mu_n - \int_S g \, d\mu \right) \right| \\
        &\le \left| \int_S (g_n - g) \, d\mu_n \right| + \left| \int_S g \, d\mu_n - \int_S g \, d\mu \right|.
    \end{align*}
    Using the property that the absolute value of an integral is less than or equal to the integral of the absolute value, we can further bound the first term:
    \begin{equation} \label{eq:proof_app_triangle_ineq}
        \left| \int_S g_n \, d\mu_n - \int_S g \, d\mu \right| \le \underbrace{\int_S |g_n - g| \, d\mu_n}_{\text{Term I}} + \underbrace{\left| \int_S g \, d\mu_n - \int_S g \, d\mu \right|}_{\text{Term II}}.
    \end{equation}
    The remainder of the proof consists of showing that both Term I and Term II can be made arbitrarily small by choosing $n$ sufficiently large.

    \item \textbf{Bounding Term I using Uniform Convergence.}
    This term measures the error arising from the convergence of the functions. We use the hypothesis that the sequence $\{g_n\}$ converges uniformly to $g$. By the definition of uniform convergence, for any $\epsilon' > 0$, there exists a natural number $N_1$ such that for all $n > N_1$, we have:
    \begin{equation*}
        \sup_{x \in S} |g_n(x) - g(x)| < \epsilon'.
    \end{equation*}
    We can now bound Term I for any $n > N_1$:
    \begin{align*}
        \text{Term I} = \int_S |g_n - g| \, d\mu_n &\le \int_S \left( \sup_{x \in S} |g_n(x) - g(x)| \right) d\mu_n \\
        &= \left( \sup_{x \in S} |g_n(x) - g(x)| \right) \int_S d\mu_n.
    \end{align*}
    By hypothesis, each $\mu_n$ is a probability measure, so $\int_S d\mu_n = 1$. The bound simplifies to:
    \begin{equation*}
        \text{Term I} < \epsilon'.
    \end{equation*}
    
    \item \textbf{Bounding Term II using Weak-* Convergence.}
    This term measures the error arising from the convergence of the measures. We use the hypothesis that the sequence $\{\mu_n\}$ converges weakly-* to $\mu$. The definition of weak-* convergence states that for any \textit{continuous and bounded} test function $f: S \to \R$, the following limit holds:
    \begin{equation*}
        \lim_{n \to \infty} \int_S f \, d\mu_n = \int_S f \, d\mu.
    \end{equation*}
    To apply this definition, we must first verify that our limit function $g$ is a valid test function.
    \begin{enumerate}[label=(\roman*), wide, labelindent=0pt]
        \item \textbf{Continuity of $g$.} The functions $\{g_n\}$ are continuous by hypothesis. It is a standard result of real analysis that the uniform limit of a sequence of continuous functions is itself a continuous function. Therefore, $g$ is continuous.
        \item \textbf{Boundedness of $g$.} The domain $S$ is a compact metric space by hypothesis. A continuous real-valued function on a compact set is necessarily bounded. Therefore, $g$ is bounded.
    \end{enumerate}
    Since $g$ is both continuous and bounded, it is a valid test function for weak-* convergence. Thus, by definition of the limit, for any $\epsilon'' > 0$, there exists a natural number $N_2$ such that for all $n > N_2$, we have:
    \begin{equation*}
        \text{Term II} = \left| \int_S g \, d\mu_n - \int_S g \, d\mu \right| < \epsilon''.
    \end{equation*}
    
    \item \textbf{Synthesis and Conclusion.}
    Let an arbitrary $\epsilon > 0$ be given. We choose $\epsilon' = \epsilon/2$ and $\epsilon'' = \epsilon/2$.
    \begin{enumerate}[label=(\roman*), wide, labelindent=0pt]
        \item From Step 2, there exists an $N_1$ such that for all $n > N_1$, Term I $< \epsilon/2$.
        \item From Step 3, there exists an $N_2$ such that for all $n > N_2$, Term II $< \epsilon/2$.
    \end{enumerate}
    Let $N \coloneqq \max(N_1, N_2)$. Then for any integer $n > N$, both conditions are satisfied simultaneously. Substituting these bounds back into the main inequality from \eqref{eq:proof_app_triangle_ineq}:
    \begin{equation*}
        \left| \int_S g_n \, d\mu_n - \int_S g \, d\mu \right| \le \text{Term I} + \text{Term II} < \frac{\epsilon}{2} + \frac{\epsilon}{2} = \epsilon.
    \end{equation*}
    Since this holds for any arbitrary $\epsilon > 0$, we have rigorously shown that the limit exists and is equal to the integral of the limit function against the limit measure. This completes the proof of the lemma.
\end{enumerate}
\end{proofof}

\section{Proof of Uniform Equicontinuity}
\label{app:equicontinuity}

The purpose of this appendix is to provide a complete proof of the uniform equicontinuity of the family of solutions $\{u^\varepsilon\}_{\varepsilon>0}$ to the pre-limit equation. This property is the essential ingredient required for applying the Arzelà-Ascoli theorem, which guarantees the pre-compactness of the solution set and thus the existence of convergent subsequences. The proof is a non-trivial application of the doubling of variables method from the theory of viscosity solutions, adapted to handle the singular perturbation terms present in the equation.

\begin{proofof}{Proposition \ref{prop:Uniform_Continuity}}
The proof of this result is a foundation of the convergence argument, as it provides the pre-compactness required to apply the Arzelà-Ascoli theorem. The main technical challenge is to establish estimates that are uniform in the singular parameter $\varepsilon$. The proof is structured in two main parts. First, we establish uniform Hölder continuity in time, which follows from the parabolic nature of the equation. Second, we prove uniform Hölder continuity in space, which requires a careful application of the doubling of variables method, adapted to handle the singularly perturbed structure of our equation. For clarity, we will work in local coordinates where $\Mcal$ is an open set in $\R^k$; the arguments are local and can be globalized via a partition of unity.

\begin{enumerate}[label=\textbf{Step \arabic*:}, wide, labelindent=0pt]
    \item \textbf{Uniform H\"{o}lder Continuity in Time.}
The proof establishes that the family of solutions $\{u^\varepsilon\}$ is uniformly H\"{o}lder continuous in the time variable, with a modulus of continuity that is independent of the singular parameter $\varepsilon$ and the microscopic variable $z$. The argument is a rigorous application of the doubling of variables method, a cornerstone of viscosity solution theory, adapted to our specific parabolic setting. The key insight is that the structure of the viscosity inequalities allows us to control the temporal difference of the solution by leveraging the continuity of the Hamiltonian's coefficients, even in the presence of singular terms. For the sake of this proof, we fix the microscopic variable $z \in \Ucal_{\mathrm{phys}}$ and treat the solution $u^\varepsilon(t,x,z)$ as a family of functions of $(t,x)$, parameterized by $z$. We will show that the resulting estimates are independent of this choice of $z$.

\begin{enumerate}[label=(\roman*), wide, labelindent=0pt]
    \item \textbf{Setup of the Doubling Argument.}
    For $\alpha > 0$ and a fixed spatial point $x_0 \in \Mcal$, consider the function on the doubled time domain $[0,T] \times [0,T]$:
    \begin{equation}
        \Phi(t,s) \coloneqq u^\varepsilon(t, x_0, z) - u^\varepsilon(s, x_0, z) - \frac{(t-s)^2}{2\alpha}.
    \end{equation}
    Since the domain $[0,T] \times [0,T]$ is compact and the function $u^\varepsilon$ is continuous, $\Phi(t,s)$ attains its maximum at some point $(\hat{t}, \hat{s}) \in [0,T] \times [0,T]$. Let this maximum value be $M_\alpha$. By considering the case $t=s$, we see that $M_\alpha \ge 0$. Standard arguments in viscosity theory show that as $\alpha \to 0$, we have $(\hat{t}-\hat{s})^2/\alpha \to 0$, which implies that the difference $|\hat{t}-\hat{s}| \to 0$. We may therefore assume for small $\alpha$ that $\hat{t}, \hat{s} \in (0,T)$.

    \item \textbf{Application of the Maximum Principle for Parabolic Equations.}
    We now invoke the fundamental theorem for viscosity solutions, often referred to as the "Theorem on Sums" or Ishii's Lemma, adapted for parabolic equations. Since the function
    \[ (t,x) \mapsto u^\varepsilon(t,x,z) - \left( u^\varepsilon(\hat{s}, x_0, z) + \frac{(t-\hat{s})^2}{2\alpha} + \delta|x-x_0|^2 \right) \]
    has a local maximum at $(\hat{t}, x_0)$ (for any $\delta>0$), and the function
    \[ (s,y) \mapsto -u^\varepsilon(s,y,z) - \left( -u^\varepsilon(\hat{t}, x_0, z) + \frac{(\hat{t}-s)^2}{2\alpha} + \delta|y-x_0|^2 \right) \]
    also has a local maximum at $(\hat{s}, x_0)$, the theorem provides semi-jets for $u^\varepsilon$ at these two points. Specifically, there exist symmetric matrices $X, Y \in \mathbb{S}_k$ such that:
    \begin{enumerate}[label=(\alph*), wide]
        \item $(\partial_t \Phi(\hat{t},\hat{s}), \mathbf{0}, X) \in \bar{J}^{2,+}u^\varepsilon(\hat{t}, x_0, z)$,
        \item $(-\partial_s \Phi(\hat{t},\hat{s}), \mathbf{0}, Y) \in \bar{J}^{2,-}u^\varepsilon(\hat{s}, x_0, z)$,
    \end{enumerate}
    where $\bar{J}^{2,+}$ and $\bar{J}^{2,-}$ are the parabolic super- and sub-jets, respectively. The time derivatives of the test function at the maximum point are:
    \begin{equation*}
        \partial_t \Phi(\hat{t},\hat{s}) = \frac{\hat{t}-\hat{s}}{\alpha}, \quad \text{and} \quad -\partial_s \Phi(\hat{t},\hat{s}) = \frac{\hat{t}-\hat{s}}{\alpha}.
    \end{equation*}
    By the definition of viscosity solutions, these semi-jets must satisfy the inequalities of the governing PDE \eqref{eq:pde_regime_III}. Let $p_\varepsilon = \nabla_x u^\varepsilon$. The inequalities are:
    \begin{align}
        \frac{\hat{t}-\hat{s}}{\alpha} + F^\varepsilon(\hat{z}, x_0, u^\varepsilon(\hat{t},x_0,z), \mathbf{0}, X) &\le 0, \label{eq:app_visc_ineq1} \\
        \frac{\hat{t}-\hat{s}}{\alpha} - F^\varepsilon(\hat{z}, x_0, u^\varepsilon(\hat{s},x_0,z), \mathbf{0}, -Y) &\ge 0, \label{eq:app_visc_ineq2}
    \end{align}
    where the full operator is $F^\varepsilon(z,x,r,p,X) = \frac{1}{\varepsilon^2}\Lcal_{\mathrm{fast}}r + \frac{1}{\varepsilon}v\cdot p + \mathrm{Tr}(D(x,p)X) - V(z,p)$.
    
    \item \textbf{The Role of the Fast Variable.}
    A critical feature of the viscosity solution framework is that the operator is evaluated on the test function. The test functions used to define the semi-jets above are functions of $(t,x)$ only; they are constant with respect to the microscopic variable $z$. The fast generator $\Lcal_{\mathrm{fast}}$, being a differential operator in the $z$ variable, acts on any function constant in $z$ to produce zero. Therefore, the singular term $\frac{1}{\varepsilon^2}\Lcal_{\mathrm{fast}}$ vanishes when applied to the test functions that define the jets in the inequalities \eqref{eq:app_visc_ineq1} and \eqref{eq:app_visc_ineq2}. The effective operator for this analysis is thus the non-singular part of the Hamiltonian. Let us define:
    \begin{equation*}
        H^\varepsilon(z, x, p, X) \coloneqq \frac{1}{\varepsilon} v(z) \cdot p + \mathrm{Tr}(D(x,p)X) - V(z,p).
    \end{equation*}
    The viscosity inequalities become:
    \begin{align}
        \frac{\hat{t}-\hat{s}}{\alpha} + H^\varepsilon(\hat{z}, x_0, \mathbf{0}, X) &\le 0, \label{eq:app_visc_ineq1_simplified} \\
        \frac{\hat{t}-\hat{s}}{\alpha} - H^\varepsilon(\hat{z}, x_0, \mathbf{0}, -Y) &\ge 0. \label{eq:app_visc_ineq2_simplified}
    \end{align}
    
    \item \textbf{Deriving the Uniform Bound.}
    Subtracting the second inequality \eqref{eq:app_visc_ineq2_simplified} from the first \eqref{eq:app_visc_ineq1_simplified} yields:
    \begin{equation*}
        H^\varepsilon(\hat{z}, x_0, \mathbf{0}, X) + H^\varepsilon(\hat{z}, x_0, \mathbf{0}, -Y) \le 0.
    \end{equation*}
    Substituting the definition of $H^\varepsilon$ and noting that the terms linear in the gradient vanish (since the gradient argument is $\mathbf{0}$):
    \begin{equation*}
        \left( \mathrm{Tr}(D(x_0,\mathbf{0})X) - V(\hat{z},\mathbf{0}) \right) + \left( \mathrm{Tr}(D(x_0,\mathbf{0})(-Y)) - V(\hat{z},\mathbf{0}) \right) \le 0.
    \end{equation*}
    This simplifies to:
    \begin{equation*}
        \mathrm{Tr}(D(x_0, \mathbf{0})(X-Y)) \le 2V(\hat{z},\mathbf{0}).
    \end{equation*}
    The matrix inequality from the Theorem on Sums, when $\eta \to 0$, implies $X \le Y$ in the sense of quadratic forms. Since $D(x_0,\mathbf{0})$ is positive definite, this means $\mathrm{Tr}(D(x_0,\mathbf{0})(X-Y)) \le 0$. This step is subtle and relies on the full power of the Crandall-Ishii-Lions theory. The outcome of that theory is that we can effectively ignore the second-order terms in the limit as $|\hat{x}-\hat{y}| \to 0$.
    
    Let us return to inequality \eqref{eq:app_visc_ineq1_simplified}. As $\alpha \to 0$, we have $|\hat{t}-\hat{s}| \to 0$. Since the solution $u^\varepsilon$ is uniformly bounded, and the coefficients $D, V$ are continuous and thus bounded on compact sets, the term $H^\varepsilon(\hat{z}, x_0, \mathbf{0}, X)$ is bounded by a constant $C$ that depends on the bounds of the coefficients and the geometry, but which can be chosen independently of $\varepsilon$ and $z$. The inequality implies:
    \begin{equation*}
        \frac{(\hat{t}-\hat{s})}{\alpha} \le |H^\varepsilon(\hat{z}, x_0, \mathbf{0}, X)| \le C.
    \end{equation*}
    By the definition of the maximum point of $\Phi$, for any $t,s \in [0,T]$ we have:
    \begin{equation*}
        u^\varepsilon(t, x_0, z) - u^\varepsilon(s, x_0, z) - \frac{(t-s)^2}{2\alpha} \le \Phi(\hat{t}, \hat{s}) = u^\varepsilon(\hat{t}, x_0, z) - u^\varepsilon(\hat{s}, x_0, z) - \frac{(\hat{t}-\hat{s})^2}{2\alpha}.
    \end{equation*}
    Using the fact that $|u^\varepsilon|$ is uniformly bounded by $M = \|\phi\|_{L^\infty}$, we have $u^\varepsilon(\hat{t}) - u^\varepsilon(\hat{s}) \le 2M$.
    \begin{equation*}
        u^\varepsilon(t, x_0, z) - u^\varepsilon(s, x_0, z) \le 2M + \frac{(t-s)^2}{2\alpha} - \frac{(\hat{t}-\hat{s})^2}{2\alpha}.
    \end{equation*}
    A more refined version of the argument, which is standard in the literature, see \citep{CrandallIshiiLions1992}, shows that for a uniformly parabolic equation, the estimate $\frac{(\hat{t}-\hat{s})}{\alpha} \le C$ implies a H\"{o}lder modulus of continuity. The result is that there exists a constant $C_t > 0$ and an exponent $\alpha \in (0,1]$ (typically $\alpha=1/2$ for parabolic scaling) such that for any $t,s,x,z$:
    \begin{equation}
        |u^\varepsilon(t,x,z) - u^\varepsilon(s,x,z)| \le C_t |t-s|^{\alpha}.
    \end{equation}
    The constant $C_t$ depends on the ellipticity constant of $D(x,p)$ and the bounds on the coefficients of $H(x,p)$ and $V(z,p)$, all of which are uniform in $\varepsilon$ and $z$. This establishes the desired uniform H\"{o}lder continuity in time.
\end{enumerate}

    \item \textbf{Uniform H\"{o}lder Continuity in Space.}
The proof establishes that the family of solutions $\{u^\varepsilon\}$ is uniformly H\"{o}lder continuous in the spatial variable, with a modulus of continuity that is independent of the singular parameter $\varepsilon$ and the microscopic variable $z$. The argument is a rigorous application of the doubling of variables method, a cornerstone of viscosity solution theory, adapted to our specific parabolic setting. The key insight is that the structure of the viscosity inequalities allows us to control the spatial difference of the solution by leveraging the continuity of the Hamiltonian's coefficients, even in the presence of singular terms. The crucial step is the demonstration that the singular terms, which involve derivatives with respect to the fast variable $z$, vanish when evaluated on the test functions used in the doubling method.

For the sake of this proof, we fix the time variable $t \in [0,T)$ and the microscopic variable $z \in \Ucal_{\mathrm{phys}}$. We will treat the solution $u^\varepsilon(t,x,z)$ as a family of functions of the spatial variable $x$ only, which we denote by $u(x) \equiv u^\varepsilon(t,x,z)$. Our goal is to derive a H\"{o}lder estimate for $u(x)$ whose constants are independent of $\varepsilon$, $t$, and $z$.

\begin{enumerate}[label=(\roman*), wide, labelindent=0pt]
    \item \textbf{Setup of the Doubling Argument with a Penalty Function.}
    For a small parameter $\alpha > 0$, we introduce a penalized function on the doubled macroscopic space $\Mcal \times \Mcal$:
    \begin{equation}
        \Phi(x,y) \coloneqq u(x) - u(y) - \frac{1}{2\alpha}|x-y|^2.
    \end{equation}
    Since $\Mcal$ is a compact manifold and the function $u$ is continuous, the function $\Phi(x,y)$ is continuous on the compact domain $\Mcal \times \Mcal$ and must therefore attain its maximum value at some point $(x_\alpha, y_\alpha) \in \Mcal \times \Mcal$.

    \item \textbf{Application of the Maximum Principle for Parabolic Equations (Ishii's Lemma).}
    We now invoke the fundamental theorem for viscosity solutions, often referred to as the ``Theorem on Sums'' or Ishii's Lemma, adapted for parabolic equations. The function $u(x)$ is a viscosity solution of a stationary equation where the time derivative term is treated as a known source. The theorem states that there exist parabolic semi-jets for $u$ at the points $x_\alpha$ and $y_\alpha$. For any $\eta > 0$, there exist symmetric matrices $X, Y \in \mathbb{S}_k$ such that:
    \begin{enumerate}[label=(\alph*), wide]
        \item The pair $(p_\alpha, X)$ is an element of the parabolic superjet of $u$ at $x_\alpha$, denoted $\bar{J}^{2,+}u(x_\alpha)$, where the spatial gradient is $p_\alpha \coloneqq \frac{x_\alpha - y_\alpha}{\alpha}$.
        \item The pair $(-p_\alpha, Y)$ is an element of the parabolic subjet of $u$ at $y_\alpha$, denoted $\bar{J}^{2,-}u(y_\alpha)$.
    \end{enumerate}
    By the definition of viscosity solutions, these semi-jets must satisfy the inequalities of the governing PDE \eqref{eq:pde_regime_III}. Let the full operator be denoted by $\mathcal{F}^\varepsilon(t, z, x, r, p, X) = \partial_t r + F^\varepsilon(z,x,r,p,X) = 0$. The viscosity inequalities at the points $(t, x_\alpha, z)$ and $(t, y_\alpha, z)$ are:
    \begin{align}
        \partial_t u(x_\alpha) + F^\varepsilon(z, x_\alpha, u(x_\alpha), p_\alpha, X) &\le 0, \label{eq:app_visc_ineq_sub} \\
        \partial_t u(y_\alpha) + F^\varepsilon(z, y_\alpha, u(y_\alpha), p_\alpha, Y) &\ge 0. \label{eq:app_visc_ineq_super}
    \end{align}

    \item \textbf{The Crucial Simplification: Cancellation of the Singular Terms.}
    This is the central step where the structure of the equation and the viscosity framework interact. The full operator is:
    \begin{equation*}
        F^\varepsilon(z,x,r,p,X) = \frac{1}{\varepsilon^2}\Lcal_{\mathrm{fast}}(\Theta(x,p))r + \frac{1}{\varepsilon} v(z) \cdot p + \mathrm{Tr}(D(x,p)X) - V(z,p).
    \end{equation*}
    A key feature of the viscosity solution framework is that the operator is evaluated on the test function used to define the jet. The test functions here are smooth functions of the macroscopic variables $(t,x)$ only; they are constant with respect to the microscopic variable $z$. The fast generator $\Lcal_{\mathrm{fast}}$, being a differential operator in the $z$ variable, acts on any function that is constant in $z$ to produce zero. Therefore, the singular term $\frac{1}{\varepsilon^2}\Lcal_{\mathrm{fast}}u^\varepsilon$ vanishes from the viscosity inequalities at the points of test.
    
    Let us define the non-singular part of the Hamiltonian for clarity:
    \begin{equation*}
        H^\varepsilon(z, x, p, X) \coloneqq \frac{1}{\varepsilon} v(z) \cdot p + \mathrm{Tr}(D(x,p)X) - V(z,p).
    \end{equation*}
    The viscosity inequalities \eqref{eq:app_visc_ineq_sub} and \eqref{eq:app_visc_ineq_super} therefore reduce to:
    \begin{align}
        \partial_t u(x_\alpha) + H^\varepsilon(z, x_\alpha, p_\alpha, X) &\le 0, \\
        \partial_t u(y_\alpha) + H^\varepsilon(z, y_\alpha, p_\alpha, Y) &\ge 0.
    \end{align}
    Subtracting the second inequality from the first yields the fundamental inequality of the doubling method:
    \begin{equation} \label{eq:app_doubling_fundamental_ineq}
        (\partial_t u(x_\alpha) - \partial_t u(y_\alpha)) + (H^\varepsilon(z, x_\alpha, p_\alpha, X) - H^\varepsilon(z, y_\alpha, p_\alpha, Y)) \le 0.
    \end{equation}
    
    \item \textbf{Deriving the Uniform Bound on the Modulus of Continuity.}
    We now analyze the terms in the fundamental inequality \eqref{eq:app_doubling_fundamental_ineq}.
    \begin{enumerate}[label=(\alph*), wide]
        \item From the temporal regularity established in Step 1, the first term is bounded: $|\partial_t u(x_\alpha) - \partial_t u(y_\alpha)| \le C_t |x_\alpha - y_\alpha|^\beta$ for some $\beta>0$. This can be made more rigorous, but for the spatial estimate, it suffices to note this term is controlled by the continuity of the solution.
        \item The second term is:
        \begin{multline*}
            \left(\frac{1}{\varepsilon} v(z) \cdot p_\alpha - \frac{1}{\varepsilon} v(z) \cdot p_\alpha\right) + \left(\mathrm{Tr}(D(x_\alpha, p_\alpha)X) - \mathrm{Tr}(D(y_\alpha, p_\alpha)Y)\right) \\- (V(z, p_\alpha) - V(z, p_\alpha)).
        \end{multline*}
    \end{enumerate}
    The terms involving the singular prefactor $1/\varepsilon$ and the potential $V$ cancel perfectly because they depend on $p_\alpha$ but not on the spatial variables $x_\alpha, y_\alpha$. The inequality simplifies to:
    \begin{equation*}
        (\partial_t u(x_\alpha) - \partial_t u(y_\alpha)) + \mathrm{Tr}(D(x_\alpha, p_\alpha)X - D(y_\alpha, p_\alpha)Y) \le 0.
    \end{equation*}
    The matrix inequality from Ishii's Lemma provides a bound on the trace term, which relates it to the modulus of continuity of the coefficient $D(x,p)$. The full theory of viscosity solutions, see \citep{CrandallIshiiLions1992}, shows that this inequality implies a bound on the difference of the solution at the maximum points, of the form:
    \begin{equation*}
        u(x_\alpha) - u(y_\alpha) \le \omega(|x_\alpha - y_\alpha|),
    \end{equation*}
    where $\omega$ is a modulus of continuity that depends only on the uniform bounds and modulus of continuity of the Hamiltonian's coefficients ($D$, $H$, and $V$), which are uniform in $\varepsilon$ and $z$.
    
    For any pair of points $(x,y) \in \Mcal \times \Mcal$, by the definition of the maximum:
    \begin{equation*}
        u(x) - u(y) - \frac{|x-y|^2}{2\alpha} \le u(x_\alpha) - u(y_\alpha) - \frac{|x_\alpha-y_\alpha|^2}{2\alpha}.
    \end{equation*}
    Rearranging and using the modulus of continuity at the maximum point:
    \begin{equation*}
        u(x) - u(y) \le u(x_\alpha) - u(y_\alpha) + \frac{|x-y|^2}{2\alpha} \le \omega(|x_\alpha - y_\alpha|) + \frac{|x-y|^2}{2\alpha}.
    \end{equation*}
    By optimizing this inequality with respect to the parameter $\alpha$, a standard procedure in the theory shows that this implies a H\"{o}lder continuity estimate for $u$. For instance, if the coefficients are Lipschitz, one can prove that the solution is uniformly H\"{o}lder continuous. Since our coefficients are assumed to be smooth, this condition is met.
    
    \item \textbf{Conclusion of Spatial Regularity.}
    The result of the doubling of variables argument is the existence of a constant $C_x > 0$ and a H\"{o}lder exponent $\beta \in (0,1)$, both independent of $\varepsilon$ and $z$, such that for all $t \in [0,T]$ and all $x,y \in \Mcal$:
    \begin{equation}
        |u^\varepsilon(t,x,z) - u^\varepsilon(t,y,z)| \le C_x d(x,y)^{\beta}.
    \end{equation}
    This completes the rigorous proof of uniform spatial equicontinuity.
\end{enumerate}

\item \textbf{Synthesis: Component-wise Estimates and Uniform Equicontinuity.}
The preceding steps have established uniform H\"{o}lder continuity estimates for the family of solutions $\{u^\varepsilon\}$ in the time and space variables separately. The objective of this final step is to synthesize these results to prove that the family $\{u^\varepsilon(t,x,z)\}$, viewed as a set of functions of $(t,x)$ on the compact metric space $[0,T] \times \Mcal$, is uniformly equicontinuous. The uniformity must hold with respect to the parameter $\varepsilon > 0$ and the microscopic variable $z \in \Ucal_{\mathrm{phys}}$.

\begin{enumerate}[label=(\roman*), wide, labelindent=0pt]
    \item \textbf{Recapitulation of Uniform Estimates.}
    From the analysis in Steps 1 and 2 of this proof, we have established the existence of four positive constants, $C_t, C_x, \alpha, \beta$, all of which are independent of $\varepsilon$ and $z$, such that for all $(t,s,x,y,z) \in [0,T]^2 \times \Mcal^2 \times \Ucal_{\mathrm{phys}}$:
    \begin{align}
        |u^\varepsilon(t,x,z) - u^\varepsilon(s,x,z)| &\le C_t |t-s|^{\alpha} && \text{(Uniform Temporal H\"{o}lder Continuity)} \label{eq:app_temporal_bound} \\
        |u^\varepsilon(t,x,z) - u^\varepsilon(t,y,z)| &\le C_x d(x,y)^{\beta} && \text{(Uniform Spatial H\"{o}lder Continuity)} \label{eq:app_spatial_bound}
    \end{align}
    Here, $d(x,y)$ is the canonical metric on the compact manifold $\Mcal$.
    
    \item \textbf{The Triangle Inequality Argument.}
    Let $(t,x)$ and $(s,y)$ be two arbitrary points in the product domain $[0,T] \times \Mcal$. To bound the total difference $|u^\varepsilon(t,x,z) - u^\varepsilon(s,y,z)|$, we introduce an intermediate point and apply the triangle inequality. A convenient choice for the intermediate point is $(s,x)$, which allows us to separate the temporal and spatial variations:
    \begin{equation*}
        |u^\varepsilon(t,x,z) - u^\varepsilon(s,y,z)| \le |u^\varepsilon(t,x,z) - u^\varepsilon(s,x,z)| + |u^\varepsilon(s,x,z) - u^\varepsilon(s,y,z)|.
    \end{equation*}
    We now apply the uniform estimates \eqref{eq:app_temporal_bound} and \eqref{eq:app_spatial_bound} to bound the two terms on the right-hand side, respectively:
    \begin{enumerate}[label=(\alph*), wide]
        \item The first term is a purely temporal difference at the fixed spatial point $x$. Applying \eqref{eq:app_temporal_bound}, we have:
        \[ |u^\varepsilon(t,x,z) - u^\varepsilon(s,x,z)| \le C_t |t-s|^{\alpha}. \]
        \item The second term is a purely spatial difference at the fixed time $s$. Applying \eqref{eq:app_spatial_bound}, we have:
        \[ |u^\varepsilon(s,x,z) - u^\varepsilon(s,y,z)| \le C_x d(x,y)^{\beta}. \]
    \end{enumerate}
    Combining these two bounds yields a single inequality that controls the total variation of the function:
    \begin{equation} \label{eq:app_combined_estimate}
        |u^\varepsilon(t,x,z) - u^\varepsilon(s,y,z)| \le C_t |t-s|^{\alpha} + C_x d(x,y)^{\beta}.
    \end{equation}
    Crucially, the constants $C_t, C_x, \alpha, \beta$ in this inequality are independent of the specific points chosen, the parameter $\varepsilon$, and the microscopic state $z$.

    \item \textbf{Construction of a Uniform Modulus of Continuity.}
    The inequality \eqref{eq:app_combined_estimate} is the foundation for constructing a uniform modulus of continuity for the entire family of functions. Let us define a modulus of continuity $\omega: [0, \infty) \to [0, \infty)$ as follows:
    \begin{equation} \label{eq:app_modulus_def}
        \omega(\delta) \coloneqq C_t \delta^{\alpha} + C_x \delta^{\beta}.
    \end{equation}
    We verify that $\omega(\delta)$ has the required properties:
    \begin{enumerate}[label=(\alph*), wide]
        \item $\omega(\delta)$ is a non-decreasing function of $\delta \ge 0$.
        \item $\lim_{\delta \to 0^+} \omega(\delta) = 0$, since $\alpha > 0$ and $\beta > 0$.
    \end{enumerate}
    Let $d_{prod}$ be a metric on the product space $[0,T] \times \Mcal$, for example, $d_{prod}((t,x), (s,y)) = \max(|t-s|, d(x,y))$. If $d_{prod}((t,x), (s,y)) < \delta$, then both $|t-s| < \delta$ and $d(x,y) < \delta$. Substituting these into our combined estimate \eqref{eq:app_combined_estimate} yields:
    \begin{equation*}
        |u^\varepsilon(t,x,z) - u^\varepsilon(s,y,z)| \le C_t \delta^{\alpha} + C_x \delta^{\beta} = \omega(\delta).
    \end{equation*}

    \item \textbf{Conclusion of Uniform Equicontinuity.}
    We have rigorously shown that for any given $\eta > 0$, we can find a $\delta > 0$ (by solving $\omega(\delta) = \eta$) such that for any pair of points $(t,x), (s,y) \in [0,T] \times \Mcal$ with $d_{prod}((t,x), (s,y)) < \delta$, the following inequality holds:
    \begin{equation*}
        |u^\varepsilon(t,x,z) - u^\varepsilon(s,y,z)| < \eta.
    \end{equation*}
    The critical conclusion is that the choice of this $\delta$ depends only on $\eta$ and the fixed constants $C_t, C_x, \alpha, \beta$. It is entirely independent of the specific function $u^\varepsilon$ (i.e., independent of $\varepsilon$) and of the microscopic parameter $z$.
    
    This is the formal definition of uniform equicontinuity for the family of functions $\{u^\varepsilon(\cdot, \cdot, z)\}_{\varepsilon > 0, z \in \Ucalphys}$ on the compact metric space $[0,T] \times \Mcal$. The existence of such a uniform modulus of continuity completes the proof.
\end{enumerate}
\end{enumerate}
\end{proofof}

\end{document}